\documentclass[openany]{amsbook}

\usepackage{mathtools} \usepackage{unicode-math} 
\usepackage{fullpage}
\usepackage{fontspec}
\usepackage[shortlabels]{enumitem}
\usepackage{graphicx}
\usepackage{amsmath} \usepackage{amsthm} \usepackage{tikz} \usepackage{comment} \usepackage{algpseudocode} \usepackage{tikz-cd} \usepackage[pdfpagelabels,bookmarks,hyperindex,hyperfigures,hidelinks]{hyperref} \usepackage{csquotes}
\usepackage[french,english]{babel}

\usepackage[
backend=biber,
style=alphabetic
]{biblatex} \usepackage{rotating} 
\usepackage{indentfirst} 
\makeatletter
\renewcommand{\@maketitle}{%
  \thispagestyle{empty}%
  \begin{center}
    {\Large\bfseries\@title\par}
    \medskip
    {\large\authors\par}
  \ifx\@empty\addresses\else
    {\footnotesize\@setaddresses}%
  \fi
  \end{center}
  
  \ifx\@empty\@subjclass\else\@setsubjclass\fi
  \ifx\@empty\@keywords\else\@setkeywords\fi
  \ifx\@empty\thankses\else\@setthanks\fi
  \@setabstract
  \bigskip
}
\makeatother

\usetikzlibrary{babel}
\usetikzlibrary{arrows,shapes,positioning}
\usetikzlibrary{decorations.markings}

\newtheorem{theorem}{Theorem}[chapter]

\newtheorem{lemma}[theorem]{Lemma}
\newtheorem{proposition}[theorem]{Proposition}

\newtheorem{corollary}[theorem]{Corollary}
\newtheorem*{corollary*}{Corollary}

\theoremstyle{definition}

\newtheorem{definition}[theorem]{Definition}
\newtheorem{example}[theorem]{Example}

\theoremstyle{remark}
\newtheorem{notation}[theorem]{Notation}
\newtheorem{remark}[theorem]{Remark}

\numberwithin{equation}{section}

\newcommand{\R}{\mathbb{R}}
\newcommand{\Q}{\mathbb{Q}}
\newcommand{\Z}{\mathbb{Z}}

\newcommand{\C}{\mathbb{C}}
\renewcommand{\P}{\mathbb{P}}
\newcommand{\F}{\mathbb{F}}
\newcommand{\G}{\mathbb{G}}
\newcommand{\B}{\mathbb{B}}
\newcommand{\T}{\mathbb{T}}

\newcommand{\W}{\mathbb{W}}
\newcommand{\D}{\mathbb{D}}

\renewcommand{\O}{\mathcal{O}}
\newcommand{\cF}{\mathcal{F}}

\newcommand{\cH}{\mathcal{H}}
\newcommand{\cA}{\mathcal{A}}
\newcommand{\cL}{\mathcal{L}}
\newcommand{\cM}{\mathcal{M}}
\newcommand{\cE}{\mathcal{E}}
\newcommand{\cR}{\mathcal{R}}
\newcommand{\cB}{\mathcal{B}}
\newcommand{\cX}{\mathcal{X}}
\newcommand{\cY}{\mathcal{Y}}
\newcommand{\cP}{\mathcal{P}}
\newcommand{\red}{\mathrm{red}}
\newcommand{\llbracket}{[\![}
\newcommand{\rrbracket}{]\!]}
\newcommand{\eM}{\mathsf{M}}
\newcommand{\eG}{\mathsf{G}}

\newcommand{\A}{\mathbb{A}}

\newcommand{\im}{\operatorname{im}}

\newcommand{\ol}[1]{\overline{#1}}
\newcommand{\ul}[1]{\underline{#1}}
\newcommand{\Tr}{\operatorname{Tr }}

\newcommand{\Gal}{\operatorname{Gal}}

\newcommand{\coker}{\operatorname{coker}}
\newcommand{\val}{\on{val}}

\newcommand{\m}{\mathfrak{m}}

\newcommand{\Sha}{\ker^1}

\newcommand{\normsg}{\trianglelefteq}
\newcommand{\vol}{\operatorname{vol}}

\renewcommand{\epsilon}{\varepsilon}
\renewcommand{\phi}{\varphi}

\newcommand{\Hom}{\operatorname{Hom}}
\newcommand{\Supp}{\operatorname{Supp}}

\newcommand{\GL}{\operatorname{GL}}
\newcommand{\SL}{\operatorname{SL}}
\newcommand{\PGL}{\operatorname{PGL}}

\renewcommand{\P}{\mathbb{P}}

\newcommand{\Ab}{\mathsf{Ab}}
\newcommand{\Sh}{\mathsf{Sh}}

\newcommand{\Mor}{\operatorname{Mor}}
\newcommand{\sh}{\operatorname{sh}}
\newcommand{\Pic}{\operatorname{Pic}}
\newcommand{\Cvx}{\operatorname{Cvx}}

\newcommand{\sHom}{\mathcal{H}\kern -.5pt om}
\newcommand{\RsHom}{R\sHom}

\newcommand{\res}{\operatorname{res}}
\newcommand{\Spec}{\operatorname{Spec}}

\newcommand{\Aut}{\operatorname{Aut}}
\newcommand{\Out}{\operatorname{Out}}
\newcommand{\Inn}{\operatorname{Inn}}
\newcommand{\Map}{\ul{\operatorname{Map}}}
\newcommand{\wt}{\operatorname{wt}}
\newcommand{\Amp}{\operatorname{Amp}}

\renewcommand{\a}{\mathfrak{a}}

\renewcommand{\m}{\mathfrak{m}}
\renewcommand{\t}{\mathfrak{t}}
\newcommand{\g}{\mathfrak{g}}
\newcommand{\fr}{\mathfrak{r}}

\newcommand{\on}[1]{\operatorname{#1}}
\newcommand{\rank}{\on{rank}}

\newcommand{\Sp}{\on{Sp}}
\newcommand{\Sym}{\on{Sym}}

\newcommand{\Res}{\on{Res}}

\newcommand{\Frac}{\on{Frac}}

\newcommand{\et}{\acute{e}t}

\newcommand{\pHi}{{}^p \mathcal{H}^i}
\newcommand{\pHishift}{{}^p \mathcal{H}^{i - 2r_s}}
\newcommand{\pH}{{}^p \mathcal{H}}
\newcommand{\good}{\mathit{good}}

\newcommand{\str}{\mathrm{str}}

\newcommand{\wtilde}[1]{\widetilde{#1}}

\def\adots{\mathinner{\mkern2mu\raise\p@\hbox{.}\mkern2mu\raise4\p@\hbox{.}\mkern1mu\raise7\p@\vbox{\kern7\p@\hbox{.}}\mkern1mu}}

\AtBeginDocument{  \spaceskip=1\fontdimen2\font plus 3\fontdimen3\font minus 1.1\fontdimen4\font
}

\DeclareFieldFormat{postnote}{\mknormrange{#1}}
\DeclareFieldFormat{volcitepages}{\mknormrange{#1}}
\DeclareFieldFormat{multipostnote}{\mknormrange{#1}}

\newcommand{\M}{\mathbb{M}}
\renewcommand{\A}{\mathbb{A}}
\newcommand{\fm}{\mathfrak{m}}

\newcommand{\reg}{\on{reg}}
\newcommand{\rs}{\on{rs}}
\newcommand{\el}{\text{-ell}}

\renewcommand{\sp}{\on{sp}}
\newcommand{\ctM}{\widetilde{\cM}}
\newcommand{\ctA}{\widetilde{\cA}}
\newcommand{\ctP}{\widetilde{\cP}}
\newcommand{\ctB}{\widetilde{\cB}}
\newcommand{\ctAt}[2]{\ctA_{(#1,#2)}}
\newcommand{\ctAe}[2]{\ctA_{(#1, #2)\el}}

\newcommand{\bM}{\mathsf{M}}
\newcommand{\bA}{\mathsf{A}}
\newcommand{\bX}{\mathbb{X}}
\newcommand{\bAe}[1]{\bA_{#1\el}}
\newcommand{\bP}{\mathsf{P}}

\newcommand{\ft}{\mathfrak{t}}
\newcommand{\fc}{\mathfrak{c}}
\newcommand{\fp}{\mathfrak{p}}
\newcommand{\fq}{\mathfrak{q}}
\newcommand{\fb}{\mathfrak{b}}
\newcommand{\fh}{\mathfrak{h}}
\newcommand{\fn}{\mathfrak{n}}
\newcommand{\fl}{\mathfrak{l}}
\renewcommand{\fm}{\mathfrak{m}}
\renewcommand{\fg}{\mathfrak{g}}
\newcommand{\fj}{\mathfrak{j}}
\newcommand{\fu}{\mathfrak{u}}
\newcommand{\fD}{\mathfrak{D}}

\newcommand{\fR}{\mathfrak{R}}
\newcommand{\fM}{\mathfrak{M}}

\newcommand{\fG}{\mathfrak{G}}
\newcommand{\Greg}{G\text{-reg}}
\newcommand{\Lie}{\on{Lie}}
\newcommand{\sep}{\text{s}}
\newcommand{\obtcone}{\prescript{+}{}{\a}}
\newcommand{\narcone}{\a^+}
\newcommand{\Ad}{\on{Ad}}
\newcommand{\ad}{\on{ad}}
\newcommand{\spl}{\mathrm{spl}}

\newcommand{\Gr}{\on{Gr}}
\newcommand{\gitq}{/\!/}
\newcommand{\comp}{\mathrm{comp}}

\newcommand{\aff}{\mathrm{aff}}
\newcommand{\inv}{\mathrm{inv}}
\newcommand{\Gell}{G\text{-ell}}
\newcommand{\Mpar}{M\text{-par}}
\newcommand{\tors}{\mathit{tors}}
\newcommand{\tf}{\mathit{tf}}
\newcommand{\ab}{\mathrm{ab}}
\renewcommand{\ss}{\mathrm{ss}}
\renewcommand{\sc}{\mathrm{sc}}

\newcommand{\Bun}{\on{Bun}}
\newcommand{\ga}{\ul{\a}}
\newcommand{\codim}{\on{codim}}

\newcommand{\gobtcone}{{}^+ \ul{\a}}
\newcommand{\gnarcone}{\ul{\a}^+}
\newcommand{\Tate}{V_l \;}
\newcommand{\loc}{\mathrm{loc}}

\makeatletter
\providecommand{\leftsquigarrow}{  \mathrel{\mathpalette\reflect@squig\relax}}
\newcommand{\reflect@squig}[2]{  \reflectbox{$\m@th#1\rightsquigarrow$}}
\makeatother

\title{The Weighted Fundamental Lemma for Non-Split Groups}
\author{Cass Halleck-Dub\'e}
\address{Max Planck Institute for Mathematics, Vivatsgasse 7, 53111 Bonn, Germany}

\begin{document}


\newcommand{\newabstract}[1]{%
  \par\bigskip
  \csname otherlanguage*\endcsname{#1}%
  \csname captions#1\endcsname
  \item[\hskip\labelsep\scshape\abstractname.]
}

\begin{abstract}


In this work, we prove the Weighted Fundamental Lemma for Lie algebras for general (quasi-split) unramified $G$. The strategy is to adapt the proof of Chaudouard--Laumon from the split case. In particular, we define the notion of a $\xi$-stable Higgs bundle for a quasi-split reductive group scheme $G$ on a curve, and prove that the stack of such form a proper, finite-type Deligne--Mumford stack over the Hitchin base. We then adapt the global cohomological theorem of Chaudouard--Laumon to the quasi-split setting, and perform the required Frobenius-point count to obtain an identity of adelic orbital integrals. A local-to-global argument deduces the WFL for Lie algebras. In the last chapter, we perform a parallel reduction to deduce Waldspurger's nonstandard WFL from a global version which follows from the same geometric theory. Our work implies the Twisted Weighted Fundamental Lemma in good characteristic, which completes Moeglin--Waldspurger's stabilization of the twisted Arthur--Selberg trace formula. We therefore render unconditional the endoscopic classification for classical groups and an array of downstream consequences in automorphic theory.

\newabstract{french}
Dans ce travail, nous démontrons le lemme fondamental pondéré pour les algèbres de Lie pour un groupe $G$ quelconque. Notre stratégie consiste à adapter au cas général la preuve de Chaudouard--Laumon dans le cas déployé. Nous définissons en particulier la notion de fibré de Higgs $\xi$-stable pour un schéma en groupes réductif quasi-déployé $G$ sur une courbe, et démontrons que le champ de tels fibrés est un champ de Deligne--Mumford propre et de type fini sur la base de Hitchin. Nous adaptons ensuite au cadre quasi-déployé le théorème cohomologique global de Chaudouard--Laumon, puis effectuons le comptage de points fixes de Frobenius nécessaire afin d’obtenir une identité entre intégrales orbitales adéliques. Un argument local-global permet alors d’en déduire le lemme fondamental pondéré pour les algèbres de Lie. Dans le dernier chapitre, nous effectuons une réduction parallèle afin de déduire le lemme fondamental pondéré non standard de Waldspurger d’une version globale, laquelle découle de la même théorie géométrique. Notre travail implique le lemme fondamental pondéré tordu en bonne caractéristique, achevant ainsi la stabilisation, due à Mœglin--Waldspurger, de la formule des traces tordue d’Arthur--Selberg. Nous rendons par conséquent inconditionnelles la classification endoscopique des groupes classiques ainsi qu’une série de conséquences qui en découlent en théorie automorphe.

\end{abstract}

\maketitle
\tableofcontents

\chapter{Introduction}\label{sect.intro}

\section{What is this document?}

This manuscript proves the variant for Lie algebras of the Weighted Fundamental Lemma (WFL), conjectured by Arthur \cite{ArthurStableI}, and its ``nonstandard'' variant, conjectured by Waldspurger \cite{WaldTordu}, in equal characteristic. By \cite[Theorem 3.8]{WaldTordu}, the WFL for Lie algebras and nonstandard WFL are together known to imply the Twisted Weighted Fundamental Lemma (TWFL), an open conjecture. Many results in automorphic theory, including the stabilization of the twisted (and untwisted) Arthur--Selberg trace formula (cf.\ \cite{ArthurStableI, ArthurStableII, ArthurStableIII}, \cite{MWStabI, MWStabII}) and hence the endoscopic classification of automorphic representations for classical groups (see \cite{ArthurEndoscopic}, \cite{IntertwinersAPackets}), were until this point conditional upon the TWFL and are therefore resolved unconditionally by this work. See \S \ref{subsect.consequences} for further discussion of these implications. 

Two special cases of the WFL were previously known. The WFL depends on a local unramified (in particular, quasisplit) reductive $G$, a Levi subgroup $M \subset G$, and an elliptic endoscopic group $M'$ for $M$. The special case $M = G$ is the (Lie algebra version of the) celebrated Langlands--Shelstad Fundamental Lemma introduced in \cite{langlandsfl}, whose standard and non-standard versions were both proven in \cite{ngofl} for arbitrary (unramified reductive) $G$. On the other hand, for reductive groups $G$ which are defined and split over $\F_q$, \textit{and} for endoscopic groups $M'$ which are split, the WFL (though not its non-standard variant) was proved in \cite{CLI, CLII}. This work, then, can be thought of either as (1) the generalization of the work of Ng\^{o} from $M = G$ to arbitrary Levis, or as (2) an extension of Chaudouard--Laumon from the case of $G$ and $M'$ constant split to an arbitrary (quasi-split) unramified reductive $G$ and endoscopic group. 

We emphasize that this work resolves a critical gap in the literature even for those concerned primarily with automorphic theory for split semisimple groups. First, because the work of Chaudouard and Laumon required $M'$ as well as $G$ to be split constant, while many applications require inductive arguments that pass through all endoscopic groups. Second, because proofs of the nonstandard (and therefore twisted) WFL have not appeared in the literature even in that generality. This work provides both at the required level of generality for applications.

Methodologically, our philosophy is to avoid re-inventing the wheel whenever possible. Many results can be directly used, and many arguments can be transferred over \textit{mutatis mutandis}, from the previous work. In the interest of a self-contained and accessible presentation, and because it was necessary to go through those works with a fine-toothed comb in order to understand where additional arguments were necessary, we reproduce many of the arguments of \textit{loc.\ cit.}, as well as referencing the corresponding results. We have made every effort to indicate clearly the difference between results which are fully contained in the prior work (cited in the format ``(\cite[Corollaire 3.9.2]{CLI})''), new results whose proofs proceed closely along the prior work (cited in the format ``(cf.\ \cite[Corollaire 3.9.2]{CLI})''), and those places where new inputs or ideas are needed. Our approach should be contrasted with the complementary work of Wang (cf.\ \cite{griffinthesis}), which seeks to prove various Fundamental Lemmas (in particular, relative ones) directly at the group-level by replacing the usual Lie algebra Hitchin fibration with multiplicative Hitchin fibrations.

In the remainder of this introduction, we begin with a gentle recollection of the basic theory of local orbital integrals as they arise in applications to the Arthur--Selberg trace formula, and state the WFL as conjectured by Arthur. We then outline the content of this document over several sections. For the benefit of the reader familiar with the work of Ng\^{o} and Chaudouard--Laumon, we emphasize some of the new ideas and arguments required for our level of generality.

\section{The Weighted Fundamental Lemma.}

Following the ideas of Arthur, study of the geometric side of the trace formula can often be reduced to very simple kinds of orbital integrals. Let $G$ be a connected quasisplit unramified reductive group over a nonarchimedean local field $F_v$ with absolute value $|\cdot|$. Assume that the residue characteristic $p$ is sufficiently large and consider the hyperspecial, Lie algebra, local orbital integral
$$
    O_{\gamma}^G := |D^G(\gamma)|^{1/2} \int_{Z_G(\gamma)(F_v)\backslash G(F_v)} 1_{\fg(\O)}(\Ad(g^{-1}) \gamma) \; dg,
$$
where $\gamma \in \fg(F_v)$ is a semisimple regular element and $D^G(\gamma)$ is the Weyl discriminant. For a good choice of Haar measures, $O_{\gamma}^G$ depends only on the $G(F_v)$-conjugacy class of $\gamma$. The centralizer $J_\gamma = Z_G(\gamma)$ is a torus; when this torus is not elliptic, that is, anisotropic mod $Z_G$, the kinds of global products of such orbital integrals which appear in the trace formula diverge.  

In order to handle divergence issues arising in the global trace formula, Arthur introduced regularizing factors, local \textit{weight functions} $v_M^{G, \loc}(g)$, on the orbit. These weight functions depend on a choice of unramified Levi $M \subseteq G$ over $F_v$; they are the volume of a certain convex hull determined by the Iwasawa decompositions of $g$ with respect to parabolics with Levi $M$ (cf.\ Definition \ref{def.arthur.wt.ftn}). The function is thus identically 1 when $M = G$, and the weights become harsher truncations as we pass to smaller Levi subgroups. A \textit{weighted orbital integral}, then, is an integral
$$
    O_{M,\gamma_M}^G = |D^G(\gamma_M)|^{1/2} \int_{J_{\gamma_M}(F_v)\backslash G(F_v)} 1_{\fg(\O)}(\Ad(g^{-1}) \gamma_M) v^{G, \loc}_M(g)\; dg,
$$
where $\gamma_M \in \fm(F_v) \subseteq \fg(F_v)$ is a semisimple $G$-regular element, so that $Z_G(\gamma_M)$ is a subtorus of $M$. Calling an element $\gamma \in \fg(F_v)$ \textit{$M$-elliptic} if the maximal split subtorus $J_{\gamma}^{\spl} \subset J_{\gamma}$ is equal to the split connected center $A_M^{\spl} \subset M$, this gives us a well-behaved notion of orbital integral for all regular semisimple $\gamma$.

Orbital integrals and weighted orbital integrals are poorly adapted to the kinds of transfer one expects from Langlands functoriality: say, between pure inner forms of a group, or from an endoscopic group for $G$ to $G$ itself. Since transfer is built atop transfer of geometric data (data defined over an algebraically closed field), this failure is partially reflected in the fact that weighted orbital integrals are constant only on $F_v$-conjugacy classes rather than on \emph{stable conjugacy classes}. The latter, here equal to the set of $F_v$-points in the $G(\ol{F_v})$-conjugacy class, is the truly algebro-geometric object, and correspondingly it is natural to seek ``stabilized'' versions of orbital integrals. For $G$-elliptic orbital integrals, the stable orbital integral is simply the sum of orbital integrals for the $G(F_v)$-conjugacy classes lying in a given stable class,
$$
    SO_{\gamma}^G := \sum_{[\gamma'] \sim_{st} [\gamma]} O_{\gamma'}^G.
$$
The set of rational classes in a stable class is finite, and indeed after choice of a base-point can be equipped with the structure of a finite group. In applications to harmonic analysis, one seeks to rewrite expressions in terms of orbital integrals with \textit{stabilized} versions involving only stable orbital integrals. Of course, symmetrizing over a stable conjugacy class naturally creates ``error terms,'' Gauss sums for nontrivial characters $s$ of the group of rational classes in a stable class:
$$
    O_{\gamma}^{s, G} = \sum_{[\gamma'] \sim_{st} [\gamma]} s(\gamma') O_{\gamma'}^G.
$$
The landmark result conjectured by Langlands and proved by Ng\^{o} is the \textit{Fundamental Lemma}, which with our discriminant normalization becomes \cite[\S 1.11.3]{ngofl}
$$
    O_{\gamma}^{s, G} =  SO_{\gamma_{H_s}}^{H_s},
$$
where $H_s$ is an \textit{elliptic endoscopic group} for $G$. Much has been written about the Fundamental Lemma, see, for example, \cite{halesfl} for a general exposition or \cite[\S 1.11]{ngofl} for precise constructions in notation similar to ours.

When $M \neq G$, on the other hand, the combinatorics required to stabilize a weighted orbital integral is more involved. Motivated by the inductive structure that arose when attempting to stabilize the entire trace formula, Arthur conjectured in \cite{ArthurStableI} an identity in the same spirit, relating twisted sums of $M$-elliptic weighted orbital integrals to sums of stable orbital integrals for endoscopic groups. In the remainder of this section, we introduce the requisite notation to state the Lie algebra version of the Weighted Fundamental Lemma as it originally appeared, in terms of endoscopic notation equivalent to, but slightly different than, the formalism of this paper.

Let $M'$ be an unramified elliptic endoscopic datum for $M$, in the sense of \cite[\S 2]{arthurenddatadef}. For the purposes of this introduction, this means $M'$ is shorthand for a triple consisting of 
\begin{itemize} 
\item a quasisplit reductive $F_v$-group $M'$,
\item a subgroup of the $L$-group $\eM \subseteq {}^L M$ which is elliptic, i.e., does not lie within any proper parabolic subgroup, and
\item a semisimple element $s_M \in \widehat{M}$.
\end{itemize}
Here $\widehat{M}$ is the Langlands dual of $\M = M_{\ol{F}_v}$, taken over some algebraically closed field of characteristic 0 (say, $\C$), and the $L$-group is defined relative to its action by $\Gamma = \Gal(F_v^{\sep}/F_v)$ through some finite quotient (indeed through $\Out(\M)$, since $M$ is quasisplit). These are required to be compatible in the sense that $\eM$ is a split extension of $W_{F_v}$ by $\widehat{M'}$ and $\widehat{M'}$ is the connected centralizer of $s_M$ inside $\widehat{M}$, as well as satisfying a certain vanishing cocycle condition, cf. \cite[\S 1.1.5]{MWStabI}. The foundation of transfer from an endoscopic group is then the elementary fact that a fixed maximal torus $T \subset M$ identifies with a maximal geometric torus of $M'$ in such a way that we have a containment of Weyl groups $W^{M'} \subset W^M$. In particular, at the level of \textit{invariants}, i.e.,\ Chevalley spaces, see \S \ref{subsect.twist.chevalley} below, we have a morphism
$$
    \fm'\gitq M' \cong \ft\gitq W^{M'} \to \ft \gitq W^M \cong \fm \gitq M.
$$
We write $$\Delta_{M',M} : \fm'(F_v) \times \fm(F_v) \to \C$$ for the Langlands--Shelstad transfer factor for Lie algebras as in \cite{waldFLimpliquetransfert}, defined on pairs of semisimple $G$-regular elements (as defined by their images in $\fg\gitq G \cong \ft\gitq W$). The value of $\Delta_{M',M}(\gamma', \gamma)$ is zero unless its inputs have the same image in $\fm\gitq M$, and when $M' = M$ it is exactly the indicator function for the relation $[\gamma'] \sim_{st} [\gamma]$.

Given a triple representing an elliptic endoscopic datum for $M$, \cite[\S 4]{arthurenddatalevi} defines an \textit{induced} family of representatives for elliptic endoscopic data for $G$. One begins with $s \in s_M Z(\widehat{M})^\Gamma/Z(\widehat{G})^\Gamma$, and then defines
\begin{itemize}
    \item $\widehat{G}'_s := Z_{\widehat{G}}(s)^0$,
    \item $\eG_s := \widehat{G}'_s \eM$, and
    \item $G'_s$ to be a certain quasisplit Langlands dual $F_v$-group of $\widehat{G}'_s$ determined by $\eG_s$.\footnote{To be somewhat more precise, $G'_s$ is the quasisplit $F_v$-form of the split Langlands dual of $\widehat{G}'_s$ twisted by the $\Out$-torsor defined by the conjugation action in $\eG_s$, which is a split extension of $\widehat{G}'_s$ by $W_{F_v}$. Compare with \textit{loc.\ cit.} and \S \ref{subsect.coh.wfl} below.}
\end{itemize}
Then $(G'_s, \eG_s, s)$ define an endoscopic datum for $G$, and moreover $M'$ can be identified with a Levi subgroup of $G'_s$. Denote $\cE_{M'}(G)$ the set of such $s$ for which the resulting endoscopic datum is elliptic, though it depends on the full endoscopic datum for $M'$, not just $M'$ itself.

Now let $\gamma' \in \fm'(F_v)$ be $G$-regular semisimple. Define the endoscopic sum of orbital integrals
$$
    O_{M',M}^G(\gamma') := \sum_{\gamma \in \fm^{G\text{-reg}}(F_v)} \Delta_{M',M}(\gamma', \gamma) O_{M, \gamma}^G
$$
with the sum, taken over representatives for the conjugacy classes, finite. One also then has a recursive definition of \textit{stable} weighted orbital integrals for $\gamma \in \fm(F_v)$ semisimple $G$-regular:
$$
    SO_{M}^{G}(\gamma) := O^G_{M,M}(\gamma) - \sum_{s \neq 1}   |Z(\widehat{G}'_s)^\Gamma /Z(\widehat{G})^\Gamma|^{-1} SO_{M}^{G'_s}(\gamma),
$$
the sum taken over $s \in \cE_{M}(G)$, intuitively, over endoscopic groups $G'_s$ for $G$ which are compatible with $M$.

The bulk of this paper is dedicated to a proof of the following (under our characteristic assumptions, cf.\ Theorem \ref{thm.wfl} for a precise statement).

\begin{theorem}[{{The Weighted Fundamental Lemma for Lie Algebras, cf.\ \cite[Conjecture 5.1]{ArthurStableI}}}]
    Assume $s_M$ is not central in $\widehat{M}$, so $M' \neq M$, and pick $\gamma' \in \fm'(F_v)$ semisimple $G$-regular. Then we have the identity
    $$
        O_{M',M}^G(\gamma') = | Z(\widehat{M}')^\Gamma/ Z(\widehat{M})^\Gamma| \sum_{s \in \cE_{M'}(G)} |Z(\widehat{G}'_s)^\Gamma/Z(\widehat{G})^\Gamma|^{-1} SO_{M'}^{G'_s}(\gamma')
    $$
\end{theorem}
In the case that $G$ and $M'$ are both split groups, this is \cite[Th\'eor\'eme 19.1.1]{CLII}. We defer to \S \ref{subsect.consequences} a discussion of the applications of this theorem.

\section{The shape of prior work.}
We recall now, impressionistically, the framework of the previous works of Ng\^{o} and Chaudouard--Laumon on the FL and WFL. For all precise definitions and results we defer to the body of the paper. 

Geometric approaches have played a central role in approaches to versions of the Fundamental Lemma since Goresky--Kottwitz--Macpherson \cite{GKM} observed that the values of local orbital integrals were computable in terms of the $\F_q$-points (and so by Grothendieck--Lefschetz, in terms of the cohomology) of affine Springer fibers. A crucial insight of Ng\^{o} was that although the WFL is a local identity of orbital integrals, it can be better approached geometrically in the global function field setting. In that setting, Arthur's stabilization of the global trace formula should be incarnate within the relative cohomology of the Hitchin fibration $f_G : \cM \to \cA$ associated to a quasi-split reductive group scheme $G$ over a curve $X/\F_q$ and a divisor $D > 0$. The space is the moduli stack of pairs $(\cE, \theta)$ where $\cE$ is a $G$-torsor on $X$ and $\theta \in H^0(X, \Ad_D(\cE))$ is the Higgs field. For classical groups, the Higgs field is a $D$-twisted endomorphism $\cE \to \cE(D)$ modeled on the Lie algebra $\fg$. The Hitchin morphism is a global geometric upgrade of the Chevalley map $\fg \to \fc$.

One reduces first to groups $G$ with anisotropic center. Over the ``$G$-elliptic'' locus $\cA^{\Gell} \subset \cA$, the Hitchin fibration for groups with anistropic center is then a proper morphism of Deligne--Mumford stacks, and so its relative cohomology with constant coefficients is geometrically semisimple, i.e., admits a Beilinson--Bernstein--Deligne--Gabber decomposition after \cite{bbdg}. The Hitchin fibration and a certain rigidification admit a faithful fiberwise action of a large commutative group stack, the Picard stack $\cP \to \cA$, which acts on the relative cohomology through its fiberwise groups of components $\pi_0(\cP)$. Over the $G$-elliptic locus, $\pi_0(\cP)$ is fiberwise finite, and so after a rigidification (denoted in this work with $\wtilde{(\;)}$, e.g.\ $f : \ctM \to \ctA$) we obtain a well-behaved isotypic decomposition.

Components of this isotypic decomposition on perverse relative cohomology $\pi_0(\ctP) \curvearrowright \pHi(R f_\star^{\Gell} \ol{\Q}_l)$ can be identified with endoscopic data $s$ for the group $G$, i.e.,\ (ignoring all shifts and twists) we have a geometric decomposition over $\ctA_{\ol{\F}_q}$
$$
    \pHi(R f_{G,\star}^{\Gell} \ol{\Q}_l) = \bigoplus_{s \in \widehat{T}} \pHi(R f_{G,\star}^{\Gell} \ol{\Q}_l)_s.
$$ 
where the decomposition ranges over elements of a torus $\widehat{T}$ in the split dual group $\widehat{G}$ and is supported on finitely many. 

The (unweighted) Fundamental Lemma is an identity equating a $s$-twisted orbital integral for the group $G$ with a stable orbital integral for a corresponding endoscopic group $H_{s}$. Given a geometric endoscopic datum $(s, \rho)$, there is a map $i_{(s, \rho)} : \ctA_{H_{s}} \longrightarrow \ctA_G$, and Ng\^{o} proves a corresponding cohomological version of the Fundamental Lemma: up to semisimplification and an explicit shift and twist,
$$
    \pHi(R f_{G,\star}^{\Gell} \ol{\Q}_l)_s \cong \bigoplus_{\rho} i_{(s, \rho), \star} \pHi(R f_{H_{s},\star}^{H_{s}\el} \ol{\Q}_l)_1,
$$
where $f_{H_s} : \ctM_{H_s} \to \ctA_{H_s}$ is the corresponding Hitchin morphism for the endoscopic group, and the sum ranges over endoscopic data associated with $s$. Building an appropriate point $a \in \ctA^{\Gell}$ which is interesting at only one place $v \in X$ and applying a product formula, Ng\^{o} can deduce the Fundamental Lemma.

The \emph{Weighted} Fundamental Lemma, then, and the non-elliptic parts of Arthur's stabilization, should live geometrically within the cohomology of the part of the Hitchin fibration which is generically regular semisimple (denoted $\cA^{\heartsuit}$ in \cite{ngofl}), but not $G$-elliptic. This is the project which, entirely in the case of split constant reductive $G/\F_q$, Chaudouard--Laumon tackle. On the non-$G$-elliptic complement, which is indexed by proper Levi subgroups $M \subsetneq G$, the Hitchin fibration is worse-behaved in several ways which fundamentally stem from the existence of split central tori in $M$ (the original group, one can reduce to be semisimple). It is not quasi-compact, with the fiber $\cP_a$ having a lattice of components of rank equal to that of the split central torus. It is also wildly non-separated. The majority of \cite{CLI, CLII} builds the foundations necessary to bypass these issues and extend Ng\^{o}'s cohomological results to this locus. 

These same infinitude problems induced by split central tori had already been visible in the harmonic analysis work of Arthur, where they appeared as the nonconvergence of orbital integrals resolved by introducing the weight functions. Chaudouard and Laumon had the brilliant realization that in the geometric framework of Ng\^{o}, in which adelic orbital integrals geometrize to Hitchin fibers, these weight factors were not a purely analytic phenomenon, but instead should \emph{come from geometry}!
In their first paper, they define an open substack $\ctM^{\xi} \subset \ctM$ of bundles satisfying a certain stability condition. Point counts on these truncated Hitchin fibers, they show, correspond to weighted orbital integrals closely related to those of Arthur. They prove fundamental geometric properties of the restricted Hitchin map $f^{\xi} : \ctM^{\xi} \to \ctA$ which ensure its cohomological formalism is well-behaved over all of $\cA^{\heartsuit}$. 

In the second paper, Chaudouard and Laumon define a corresponding truncation $\ctP^1 \subset \ctP$ of the (rigidified) Picard stack which preserves $\ctM^{\xi}$. The group of components $\pi_0(\ctP^1)$ is always finite, so they obtain isotypic decompositions over the $M$-elliptic locus for all $M$. These isotypic decompositions can be identified with endoscopic data $s_M$ for $M$. The result is a cohomological incarnation of the non-elliptic stabilization of Arthur, generalized to the nonsplit case by Theorem \ref{thm.coh.wfl} below. Computing $\F_q$-points of the $s$-components, Chaudouard--Laumon show that this induces a global identity of adelic orbital integrals. A systematic procedure of \textit{scindage} breaks this identity down into one consisting of local orbital integrals at only finitely many places, of a type which can be engineered to match that in \cite{MWStabII}, and a globalization construction for their local parameter deduces the WFL.

In the remainder of the introduction, we outline the content of this manuscript in an informal way, targeted at the reader familiar with the previous work in this area and emphasizing places of difference.

\section{Geometry and cocharacters of relative reductive group schemes.}

In Chapter \ref{sect.rel.reduct.gps}, we collect the required facts about parabolics in relative reductive group schemes. Mostly extracted from \cite{conradrgs}, the theory and combinatorics of semistandard parabolics $P$ and their semistandard Levi factors $M$ adapts well, and one can work either in the function field or globally on $X$. The globalization, we note, is contingent upon the existence of a maximal torus over $X$, which we assume (and will suffice for our local applications). We assume, as in \cite{ngofl}, that the torus $T$ is split at a rational point $\infty \in X(k)$, and use this splitting to rigidify our Higgs bundles.

As in the theory of Borel and Tits \cite{BorelTits}, one can compare semistandard Levis via the centralizers of their \textit{split} central tori. In many cases, we therefore need to supplement the existing arguments to track the behavior of split tori. See for example the proof of Proposition \ref{prop.base.disjoint}.

Stability convexes naturally live in cocharacter spaces $\a_M = \Hom(X^\star(M), \R)$ associated to semistandard Levis $M$. For such a nonsplit group, there are two notions of this space: the \textit{rational} space $\a_M$ corresponding only to cocharacters defined over $X$ and the \textit{geometric} space $\ul{\a}_M$ corresponding to cocharacters over a splitting cover. Given a $\Gamma$-Galois \'etale cover splitting $G$, the former can be realized as either the Galois invariants or co-invariants of the latter. The Weyl group of the split form $\W = W(\infty)$ acts on the geometric space, while the global Weyl group $W(X) \subset W(\infty)$ acts on the rational space. 

In contrast to the discrepancy in splitting behavior between $G$ and the local $G_v$, the assumption that $T$ be split at $\infty$ implies that the splitting behavior of $T$ and the cocharacter spaces is independent of extension of the constant field of $X$; see Proposition \ref{prop.splitrankextensions}. This is useful when passing between $k_0 = \F_q$ and $k = \ol{\F}_q$ and counting points.

\section{Stratification of the Hitchin base and relative position.}

The reader familiar with the work of Chaudouard and Laumon will recall that in order to remove a Weyl-group ambiguity from the stability convex, it is necessary to work over the \'etale neighborhood $\mathcal{W}_\infty : \ctA \to \cA$ consisting of choosing an ordering $t$ of the roots of the characteristic polynomial $a$ at the point $\infty$, or a point of the cameral cover (see \S \ref{subsect.Ngo.monod.comp}). This ordering allows one to distinguish the ``relative position'' of different parabolic reductions of a $G$-torsor by comparing $t$ to the order in which the $\theta(\infty)$-eigenspaces are included in the parabolic reduction. Correspondingly, the stability conditions in Chaudouard--Laumon's convex are cones which point in the directions of different Weyl-translates of the obtuse cone spanned by the coroots, or equivalently defined by positive conditions on the fundamental weights.

The previous work defines a stratification $\bigcup_M \ctA_{M\el}$ over semistandard Levis based on the ``Levi factorization type'' of the polynomial $a$ (that is, element of the $F$-valued points of the Chevalley space). This stratification characterizes how many and which parabolic reductions the Higgs bundles $(\cE, \theta)$ in the fiber admit. Chaudouard and Laumon interpret the reductions in different relative positions to $t$ as reductions to the different semistandard parabolics which contain the Levi $M$. Crucially, their stratification of the base is the coarsest such that the set of parabolic reductions and their relative positions remains constant within a stratum, so is the appropriate stratification for studying weighted stability.

When adapting this, we run into a technical obstruction: the group $G$, by necessity, is more split at the point $\infty$ than generically. However, it is the generic behavior of $G$ which determines the set of parabolic reductions that a torsor admits. For example, if $m = (\cE, \theta)$ is a Higgs bundle living in the $T$-elliptic locus (i.e., $\chi(\theta)$ lifts to $\ft$ in the function field), then there will be a $W(X)$-torsor of Borel reductions for $m$. If we choose $t$ arbitrarily from among the $W(\infty)$-torsor of possible orderings, then the relative position of the parabolic reductions may lie in some nontrivial $W(F)$-coset, and so cannot be naturally bijected with semistandard parabolics containing $M$. 

We evade this obstruction in Chapter \ref{sect.Hitchin.base} by refining our stratification of the Hitchin base beyond that in \cite{CLI} to one indexed by pairs $$(w, M) \in W(\infty) \times \{\text{Levis $M \supseteq T$}\},$$ up to a natural equivalence relation, in which $M$ reflects the elliptic type of the Higgs field $\theta$ and $w$ measures the relative position between $t$ and $\theta$. As a consequence of our more general analysis, we can then define a \textit{compatible locus}, an open $\bA$ inside $\ctA$, which surjects onto $\cA^{\infty\text{-}r.s.}$. On this locus these questions of relative position can be ignored, and the parabolic reductions of an $M$-elliptic Hitchin triple can be identified bijectively with semistandard parabolics containing $M$ as in \cite{CLI}. In \S \ref{subsect.Ngo.monod.comp}, we show that our stratification is still subsumed by the stratification by ``monodromic invariants'' of Ng\^{o}, and recall the basic functoriality of the Hitchin base in endoscopic groups.

In Chapter \ref{sect.Hitchin.fib}, we recall the definition and basic properties of the Hitchin fibration $f : \cM \to \cA$ and its rigidification $f : \ctM \to \ctA$ over the locus of Higgs bundles $(\cE, \theta)$ for which $\theta$ is generically regular-semisimple. We define the relative position of a parabolic reduction for a rigidified Higgs bundle (which we call ``Hitchin triples'' after Chaudouard--Laumon). Our primary playground for the rest of the manuscript is the part of the Hitchin fibration over the compatible locus, which we denote $f : \bM \to \bA$.

\section{Weighted stability.}

In Chapter \ref{sect.weighted.stab}, we define degree vectors associated to the parabolic reductions. These degrees fit together into stability convexes, certain polyhedral convexes in the style of Behrend \cite{behrendthesis} inside the geometric and relative \textit{cocharacter spaces} $\ul{\a}_T$ and $\a_T$. Using them, we define the notion of $\xi$-semistability. The moduli of $\xi$-semistable Higgs bundles $\ctM^{\xi} \subset \ctM$ defines an open substack.

Our formalism for weighted stability is ultimately a mixture of rational (i.e.,\ relative) techniques which work directly on $X$ with the relative root datum $\Phi(G, S)$ and geometric (i.e.,\ absolute) techniques that work over a splitting cover $X'$ with the absolute root datum $\Phi(G_{X'}, T_{X'})$.

We define first a ``geometric'' stability convex $\ul{C}_m$, which lives within the geometric cocharacter space $\ga_T$ but imposes only the stability conditions corresponding to parabolic reductions which descend to $X$ (in particular, is \textit{not} equal to the stability convex for the pulled-back Hitchin triple, so our definition does not reduce to the split case).

Over the compatible locus, the relative positions can all be taken to be trivial, and so we can relate our convex to a ``rational'' convex $C_m$ defined within the relative space $\a_T$. This comparison is cleanest over the compatible locus. Ultimately, we use the flexibility of both perspectives: Arthur's weight functions are defined and analyzed entirely in the relative setting, and so comparison to the rational stability convex allows us to apply the point counting arguments of Chaudouard--Laumon directly, without change, and allows comparison of our stabilization formalism to Arthur's.
On the other hand, many of the geometric arguments of Chaudouard and Laumon reduce to computations on line bundles corresponding to geometric roots---computations which do not easily adapt to relative weight spaces, and so naturally occur over a splitting cover. 

We emphasize that all definitions and constructions related to $\xi$-semistability given in this paper are novel outside the split case. The previous work on stability convexes in Chaudouard--Laumon, as well as the complementary polyhedra of Behrend, are defined \textit{only} for torsors for split groups (equivalently, for group schemes over $X$ which are inner forms).\footnote{Behrend's work defines some related concepts (such as the canonical/Harder--Narasimhan parabolic reduction) outside the split case, but the way in which he does so is specific to $\xi = 0$, which has the distinguishing property of being fixed by the Galois action on the geometric cocharacter space.}

A crucial warning is that although the notion of degree of an \textit{individual} parabolic reduction is well-behaved under pullback, the stability convexes we define (i.e.,\ the set of \textit{all} parabolic reductions and their degrees) are \textit{not} compatible with descent along a splitting cover. That is, the stability convex for a $G$-Hitchin triple is \textit{not} equivalent to one defined by descending the stability convex associated to the pulled-back Hitchin triple along a splitting cover. 

This failure to behave well under descent is a reflection of the arithmetic fact that elliptic endoscopic groups can become Levi subgroups upon pullback, or that centralizing tori for regular semisimple elements can become more split. Consider, for example, a $U_{X'/X}(n)$-Higgs bundle $(\cE, \theta)$ whose characteristic polynomial lies in an elliptic endoscopic locus $\ctA_{H_s}$ corresponding to endoscopic group $H_s = U(k) \times U(n-k)$. A generic such Higgs bundle is $G$-elliptic (that is, anisotropic mod split center) in Ng\^{o}'s formalism, and the condition of being $\xi$-stable should therefore be vacuous. However, the pulled-back Hitchin bundle for $\GL_{n,X'}$ lies inside the locus corresponding to a proper Levi $\GL_k \times \GL_{n-k}$, so we obtain a nontrivial stability condition upstairs which, it turns out, would descend to cut out much of the endoscopic locus entirely.

\section{Uniformization and adelic characterization.}

In order to relate Hitchin fibers with adelic objects, we replace the uniformization results for $G$-torsors of \cite{drinfeldsimpson} with corresponding relative results, e.g.\ Steinberg's theorem and \cite{HeinlothRelUnif}. We prove in Chapter \ref{sect.adelic.fibers} an adelic description of the truncated Hitchin fiber on the compatible locus, making use of the rational stability convex.

\section{\texorpdfstring{Geometric properties of the $\xi$-stable Hitchin morphism.}{Geometric properties of the xi-stable Hitchin morphism.}}

Outside the $G$-elliptic locus, the Hitchin fibration acquires some geometric pathologies which obstruct direct adaptation of the cohomological methods of Ng\^{o}. The truncation to $\xi$-semistable bundles resolves some of these (the total space is now quasicompact in the semisimple case, for example), while introducing new concerns. In particular, the Hitchin fibration $f : \ctM \to \ctA$ satisfies the existence half of the valuative criterion (though it is very far from separated), and it is not clear that the restriction to the open substack $\ctM^{\xi}$ of $\xi$-semistable bundles still contains a lift for all specializations valued in a DVR. Indeed, for generic $\xi$ the $\xi$-semistable locus turns out to truncate the fibration exactly the right amount, as Chaudouard and Laumon realized: the restricted Hitchin map $f^\xi : \ctM^{\xi} \to \ctA$ for a semisimple split group $G$ satisfies \emph{both} halves of the valuative criterion. 

In the nonsplit setting, since our stability conditions have been defined in the relative space, the real obstruction is the existence of \emph{split} central tori, and so the natural restriction is to assume $G$ has \emph{anisotropic center}. This is also preserved under passage to elliptic endoscopic groups, while semisimplicity may not be -- crucial for the recursive construction of the stable orbital integrals. We carry this hypothesis through the geometric and cohomological parts of the paper, until we can recover our desired generality on the arithmetic side.

In order to prove the valuative criterion, we require a more systematic analysis of not only the $\xi$-stable (properly, the $\xi$-semistable) locus, but a stratification of its complement by a version of Harder--Narasimhan type. Much of the geometry of $\Bun_G$, $\cM$, and similar moduli stacks lies in not only its semistable locus, but in its full stratification by the type of the Harder--Narasimhan filtration of $G$-bundles (see for example \cite{lurieweil}, or \cite{schiederHN} for a nice exposition in the split case). A systematic formalism for such stratifications has been axiomatized by Halpern--Leistner in \cite{HalpernLeistnerInstability}, and these \emph{$\Theta$-stratifications} allow uniform treatment of many geometric results across GIT and the geometry of bundles.

In Chapter \ref{sect.can.red}, we define the $\xi$-canonical reduction, a particular parabolic reduction of a Hitchin triple which is in an appropriate sense ``maximally $\xi$-destabilizing.'' 
We prove that the $\xi$-canonical reductions can be used to define a $\Theta$-stratification of the Hitchin fibration whose minimal stratum is $\bM^\xi$. This allows us, in Chapter \ref{sect.val.exist}, to avail ourselves of a general semistable reduction theorem for $\Theta$-stratifications \cite{AHLHModuliForStacks}. We thereby streamline and avoid generalizing the corresponding proof from \cite{CLI}, a rather technical adaptation of Langton's original semistable reduction algorithm for vector bundles which we found difficult to generalize to the non-split setting.

In Chapter \ref{sect.finite.type} we show that when $G$ has anisotropic center, the substack of $\xi$-semistable Hitchin triples $\ctM^\xi$ is of finite type. We require some results generalizing classical theorems of Harder on the existence of Borel reductions of bounded stability. Our results in this direction are by reduction to the absolutely almost simple case of split rank 1. They may also be able to be extracted from \cite{lurieweil}, with a different proof, but the proof we give here is self-contained. Some parts of the finite-type argument (but not the construction of Borel reductions) follow by descent from the split case. Other parts of the argument work directly with the relative root datum. 

In Chapter \ref{sect.val.unique} we prove the uniqueness half of the valuative criterion of properness (i.e., separatedness) holds for the stack of $\xi$-\textit{stable} (not $\xi$-semistable) bundles as long as $G$ has anisotropic center. For $\xi \in \a_T$ in general position (avoiding a countable set), $\xi$-stability and $\xi$-semistability coincide, so we can obtain both halves of the valuative criterion for the same stack. In Chapter \ref{sect.Deligne.Mumford} we prove that for $\xi$ in general position the moduli of $\xi$-(semi)stable Higgs bundles is a Deligne--Mumford stack. The methods of these two sections are more or less direct adaptations of the proofs of \cite{CLI}. 

The geometric properties as we will apply them can thus be summarized as follows.
\begin{theorem}[{Theorem \ref{thm.bdd}, Theorem \ref{thm.sep}, Theorem \ref{thm.val.crit.exist}, Theorem \ref{thm.xi.locus.dm}}]
    If $G/X$ is a quasisplit reductive group with anisotropic center, $k$ satisfies the characteristic assumptions of \S \ref{subsect.char.assumptions}, and $\xi$ is in \emph{general position} in the sense of Definition \ref{def.gen.pos}, then the compatible $\xi$-(semi)stable Hitchin fibration
    $$
        f : \bM^\xi \to \bA
    $$
    is a proper morphism from a smooth Deligne--Mumford $k$-stack of finite type to a smooth $k$-scheme of finite type.
\end{theorem}

\section{Picard stacks and degree maps.}

In Chapter \ref{sect.reg.cent.Pic}, we recall the Picard stack $\cP$ of Ng\^{o}, whose action on Hitchin fibers controls their cohomology intimately. We define a degree map out of the restriction $\bP \to \bA$ of the Picard stack over the compatible locus, which on $M$-elliptic fibers takes values in $(X^\star(M))^\vee$, the dual lattice to the \textit{rational} characters on $M$ (see Proposition \ref{prop.dual.lattice.global.invs.coinvs.mod.tors} and \S \ref{subsect.Pic.deg} below). The rank of this invariant, then, is controlled by the rank of the maximal split central torus in $M$. The fiberwise kernel $\bP^{1}$ of this map preserves the $\xi$-stable truncation $\bM^{\xi} \to \bA$. 

In the second part of this section, we then recall Ng\^{o}'s construction of a map $\bX_\star = X_\star(T_\infty) \cong X^\star(\widehat{T}) \to \pi_0(\cP)$, the fiberwise-component group of the Picard stack. Using the induced action of an appropriate subgroup of $\bX_\star$ that preserves truncation, we extend the isotypic decomposition of Ng\^{o} on the relative cohomology beyond the $G$-elliptic locus. Here several features are visible that were not in the split case, and our description of the Picard stack components $\pi_0(\ctP)$ is less complete than in the split case. Nevertheless, we obtain in Theorem \ref{thm.s.decomp} \emph{two} well-behaved isotypic $s$-decompositions of the perverse cohomology over the compatible locus, one \emph{fine} and one \emph{coarse}. The discrepancy between the two, invisible in the case of split groups, measures the failure of taking $\Gamma$-invariants and taking the identity component to commute on $Z_{\widehat{G}}$, and ultimately provides a necessary scalar factor in the WFL.

\section{Cohomological theorems.}

We make use of the machinery we have generalized from the split case to carry over Chaudouard--Laumon's global cohomological theorems from which they deduce the WFL. We also adapt in Chapter \ref{sect.good.locus} the deformation-theoretic proof of a general discriminant criterion for a point of the Hitchin base to lie in Ng\^{o}'s ``good locus,'' see Theorem \ref{thm.suff.criterion.for.goodness}. This computation, which carries over with little change, will be essential in the final step of the proof of the WFL, when we globalize our local parameter and seek to apply our global theorems. 

In Chapter \ref{sect.coh.thms}, we prove the core global cohomological theorems. These are Theorems \ref{thm.no.new.supports}, \ref{thm.coh.stab.geo}, and \ref{thm.coh.wfl} below. The first half of the chapter, up to the second theorem, is essentially an exposition of certain parts of \cite[\S 7]{ngofl} and \cite[\S 10]{CLII}, with a few modifications. In order to keep $G$ and $M$ within the category of split groups, Chaudouard--Laumon restrict to an open which avoids the endoscopic loci in the Hitchin base corresponding to nonsplit endoscopic groups. Our theorem more naturally ranges over all endoscopic data for $G$ extending the coarse parameter $\ol{s}_M$ for $M$. Our cohomological theorem also exhibits a crucial asymmetry between coarse and fine $s$-isotypic components, which ultimately will allow the correct arithmetic scalars to arise. 

The restriction to split groups also allows \cite{CLII} to avoid geometrizing Arthur's general formalism for inducing endoscopic data from $M$ to $G$. In the last sections of the Chapter, we provide such a geometrization in Definition \ref{def.induced.end.datum} and analyze the corresponding support behavior. The resulting Theorem \ref{thm.coh.wfl} directly geometrizes Arthur's WFL, giving a cohomological decomposition ranging over exactly those endoscopic data $(s, \rho_s)$ which are induced from $(s_M, \rho_{s_M})$.

\section{Point counts and the global WFL.}

In Chapter \ref{sect.global.point.counts}, we perform our version of the Frobenius point counts which convert the global cohomological Theorem \ref{thm.coh.wfl} into an adelic identity which matches the form of Arthur's WFL.

\begin{theorem}[{Theorem \ref{thm.global.wfl}}]
    Let $X_0$ be a smooth geometrically integral projective curve over finite field $k_0$, and let $G_0/X_0$ a quasisplit reductive group, split at a chosen rational point $\infty$ and satisfying our running characteristic hypotheses. Fix $D$ an even divisor on $X_0$ with $\deg(D) > 2g$, and let $\bA_G$ be the associated compatible Hitchin base (Definition \ref{def.comp.loc}).
    
    Let $M_0 \subseteq G_0$ be a semistandard Levi and $\fM' = (s_M, \rho_{s_M})$ be an elliptic endoscopic datum for $M_0$ with endoscopic group $M_0'$. Assume moreover that $M_0'$ is split at $\infty$. For any $a_{M'} \in \bA_{M'}^{M'\el}(k_0)$ with images $a_M \in \bA_M$ and $a \in \bA_G^{\good}$, we have an identity of adelic orbital integrals
    $$
        J^{G, s_M}_M(a_M) = |Z^\Gamma_{\widehat{M}'}/ Z_{\widehat{M}}^\Gamma| \sum_{s \in \cE_{\fM'}(G)} |Z_{\widehat{G'_s}}^\Gamma/ Z_{\widehat{G}}^\Gamma |^{-1} S_{M'}^{G'_s}(a_{M'}),
    $$
    where $J^{G, s_M}_M(a_M)$ is the global $s$-orbital integral of Definition \ref{def.s.orbital.integral.adelic} and $S_{M'}^{G'_s}(a_{M'})$ are the stabilized avatars of Definition \ref{def.stable.orb.integral.global}.
\end{theorem}
This theorem removes the assumption that $G_0$ have anisotropic center which was required through the geometric part of the theory. Here, the main new ingredients consist of the following. 
\begin{enumerate} 
\item We carefully confirm that our rational weighted stability conditions match Arthur's formalism for $(G,M)$-families and weight functions, which he defines always in the relative cocharacter space,
\item We track the interplay between the various $\Gamma$-fixed centers which appear and the various refinements of Arthur's induced parameters $\cE_{\fM'}(G)$ which intervene in the cohomological theorems.
\item We ensure that our formalism for all prior theorems is sufficiently robust to allow them to apply to iterated elliptic endoscopic groups.
\end{enumerate}

\section{Global-to-local descent and the WFL.}

Finally, in \S \ref{sect.local.global.WFL} we deduce the local WFL for Lie algebras from its global adelic version. As in \cite{CLII}, this is more difficult than for the ordinary Fundamental Lemma because, although the Hitchin fibration satisfies a product formula over places on the curve, the condition of $\xi$-semistability is a global condition. If one knew the existence of a global parameter sufficiently regular at all places except a fixed $v$, then the local stability condition at $v$ could coincide with the global condition, but such a global parameter need not exist. Instead, one requires a formalism of scindage to iteratively deduce the WFL on smaller sets of places. 

The formalism we provide in \S \ref{subsect.splitting.lemmas} - \S \ref{subsect.transverse.places} is, we hope, more minimalist than \cite[\S 16-18]{CLII}. We specialize the general notion of a ``stabilizing family'' to the condition that the WFL holds on a subset of closed points $V \subseteq |X_0|$, for a fixed $G,M, M'$, global $\Gamma$, etc. Our global theorem then becomes the statement that the WFL holds on $V = |X_0|$, and we prove a general peeling result which allows us to remove sets of places at which we know that the WFL already holds. In accordance with Ng\^{o}'s principle of ``perverse continuation,'' the proof reduces to showing that the WFL holds for sufficiently regular parameters. Some of these are simple or reduce to the usual Fundamental Lemma of Ng\^{o}, but one is genuinely nontrivial (Lemma \ref{lem.tranv.result.2}). That result, which is an explicit computation in the case that $\dim \a_M^G = 1$, allows us to bypass substantial difficulty adapting Chaudouard--Laumon's descent argument.

Finally, we need to globalize an arbitrary local parameter in a way which is sufficiently regular at all other places, and such that the WFL at the place $v$ recovers our desired local statement. With an eye towards whatever remains of brevity, we give a minimal proof that our local $s$-orbital integrals are locally constant in the characteristic parameter, rather than developing the full geometric theory of truncated affine Springer fibers. We then appeal to Poonen's formalism for Bertini theorems on finite fields \cite{poonenbertini} in order to show that we can find a desired sufficiently-regular globalization of the local parameter. Our peeling formalism then deduces the WFL from its global version.

\section{The nonstandard WFL.}

For arithmetic applications, one often requires not only the usual WFL but also its ``nonstandard'' partner, which concerns a pair of groups with based root data which are isogeneous in a precise sense.

\begin{theorem}[{Theorem \ref{thm.nonstd.WFL}, cf.\ \cite[Conjecture 3.7]{WaldTordu}}]
    Let $F_{0,v}$ be an equal characteristic local field with ring of integers $\O_{0,v}$. Let $(M_1, G_1, M_2, G_2, \psi)$ be a nonstandard endoscopic datum over $\Spec \O_{0,v}$ in the sense of Definition \ref{def.nonstd.datum}, assumed in good characteristic. Let $$b \in \fc_{M_1}(\O_{0,v}) \cap \fc_{M_1}^{G_1\text{-reg}}(F_{0,v}) = \fc_{M_2}(\O_{0,v}) \cap \fc_{M_2}^{G_2\text{-reg}}(F_{0,v})$$ be a generically regular parameter, and write
    $$
        Y_{1,b} = \epsilon_{M_1}(b), \qquad Y_{2,b} = \epsilon_{M_2}(b)
    $$
    for the Kostant lifts. Then the local stable orbital integrals, taken with the compatible choices of measure as in Remark \ref{rmk.nonstd.measures}, satisfy
    $$
        S_{M_1, v}^{G_1}(Y_{1,b}) = c_{M_1, M_2}^{G_1, G_2} S_{M_2, v}^{G_2}(Y_{2,b}),
    $$
    where 
    $$
        c_{M_1, M_2}^{G_1, G_2} = \frac{[(\Z \langle \Phi_{G_2}^\vee \rangle)^{\Gamma} : \psi_\star((\Z \langle \Phi_{G_1}^\vee \rangle)^{\Gamma})]}{[(\Z \langle \Phi_{M_2}^\vee \rangle)^{\Gamma} : \psi_\star((\Z \langle \Phi_{M_1}^\vee \rangle)^{\Gamma})]}
    $$
    is the coefficient of Waldspurger.\footnote{Indices follow the convention of Notation \ref{not.commens.rel.index}.}
\end{theorem}

A proof of the nonstandard WFL for general (even split) groups has, to the best of the author's knowledge, not appeared in print, though \cite{Chaudouardnonstandardexample} discusses the result for the nonstandard pair $\Sp_{2n}$ and $\on{SO}_{2n+1}$, and sketches a framework for deducing it by a parallel argument from the geometric results of \cite{CLII}. In \S \ref{sect.nonstd}, we complete the proposed proof, giving first a global version of the Theorem (Corollary \ref{cor.global.nonstd.WFL}) from our point-counting theorems and the work of Ng\^{o}, and then adapting the descent and scindage arguments of \cite[\S 16-18]{CLII} (and the prior chapter) to deduce the local statement. Again, the global-to-local argument requires a new (and nontrivial) special case for parameters which intersect the resultant divisor transversally (Lemma \ref{lem.nonstd.transv.result}); the result reduces to a slightly ungainly calculation in relative semisimple rank 1.

\section{Consequences.}\label{subsect.consequences}

In this work, we have proven the Lie algebra WFL and nonstandard WFL in equal characteristic. The comparisons of \cite[\S 7.3, 8.7]{waldchangeofchar} (cf.\ also \cite{cluckerstransferwfl}) imply that these results in positive characteristic suffice to prove them in characteristic 0 for the same residue characteristic. 

By \cite[Th\'eor\`eme 3.8]{WaldTordu}, the group-level twisted WFL (TWFL) follows from this pair of Lie algebra results. For the statement of the latter we refer the reader to \cite[Conjecture 3.5]{WaldTordu} and \cite[\S II.4.4]{MWStabI}. The TWFL in sufficiently large residue characteristics for general groups in turn implies the stabilization of the twisted Arthur--Selberg trace formula of M{\oe}glin--Walspurger \cite{MWStabI}, \cite{MWStabII}. As a special case, our work completes the stabilization of Arthur's untwisted trace formula (cf.\ \cite{ArthurStableI, ArthurStableII, ArthurStableIII}). 

We emphasize that our characteristic hypotheses present no obstruction to their application to the stabilization of the twisted trace formula: in the global setting of \cite[\S X]{MWStabII}, the authors specify a finite set of places $V$ (cf.\ \cite[p.\ 592-593, p.\ 1145]{MWStabII}), in particular containing the ramified places, and identify test functions on $G(F_V)$ with adelic test functions by adjoining hyperspecial test functions outside $V$. The TWFL is then applied only to the places outside $V$, so the results of that book require only the TWFL for characteristic sufficiently large (relative to the global twisted datum, in particular $F, G$). For any such global twisted datum, our residue characteristic hypotheses (Definition \ref{def.good.characteristic}), translated into mixed characteristic by \cite{waldchangeofchar} and applied locally at each place, preclude only finitely many places of $F$ and can be absorbed into $V$. It is unknown to the author whether the results in sufficiently large characteristic imply the result in all characteristics, analogous to the reduction of \cite{halesflredgoodchar} for the unweighted FL.

The stabilization of the twisted trace formula renders unconditional a range of downstream arithmetic applications. In particular, the endoscopic classification of automorphic representations for quasisplit classical groups initiated by \cite{ArthurEndoscopic} and \cite{mokendclassunitary}, and otherwise completed in the preprint \cite{IntertwinersAPackets}, is now unconditional; this removes the conditional TWFL hypothesis from many applications resting on these classifications. See the Introduction of \cite{IntertwinersAPackets} for a collection of such applications, as well as a thorough citation of dependencies. We emphasize also the recent work on modularity of surfaces, \cite{boxeretalmodforsurfaces}, which depends upon the TWFL. Likewise, any results which depend upon generic parameters for pure inner forms of quasi-split unitary groups, as in the preprint \cite{kalethaswetalgenericparams}, are freed of that condition.

\section{Acknowledgments.}

Any cursory survey of this work makes self-evident the depth of its dependence on the ideas of G\'erard Laumon. Without his insights this work would not exist, and the author wishes to acknowledge profoundly his memory. I hope to have affirmed the elegance of his ideas. 

An earlier version of Parts 1-3 of this manuscript appeared as the author's PhD thesis at UC Berkeley, supervised by Sug Woo Shin. The author wishes to thank him for unwavering support and advice over the course of several years.

Parts of this paper were completed while visiting the Institut de Math\'{e}matiques de Jussieu---Paris Rive Gauche in Fall 2023. The author wishes to profoundly thank Pierre-Henri Chaudouard for the invitation, and for his patience in talking through many of the subtler points of his work.

Parts of this paper were completed while a visiting postdoctoral fellow at the Max Planck Institute for Mathematics in Bonn. The author is grateful to the Institute for its hospitality and financial support.

Many conversations inspired the content of this paper, including enlightening discussions with James Arthur, Alexis Bouthier, Ngô Bảo Châu, Pierre-Henri Chaudouard, Carlos Esparza, Tony Feng, Thomas Hameister, Ruofan Jiang, Tasho Kaletha, Rose Lopez, Mark Macerato, Stefan Reppen, Aaron Slipper, Griffin Wang, Wenqing Wei, and Zhiwei Yun.

Discussion with AI agents was used for the sake of literature search, proofreading, and consistency-checking during the preparation of this manuscript. No part of the document was written by AI, and the author takes full responsibility for the completeness and correctness of all of its content, mathematical and otherwise.

\section{Notation and running assumptions.}\label{subsect.notation.running.assumptions}

For the majority of this work and unless explicitly stated otherwise, 
\begin{itemize} 
\item $k$ is an arbitrary field with fixed algebraic closure $\ol{k}$, 
\item $X/k$ is a smooth, proper, geometrically integral curve, 
\item $\G$ is a split reductive group over $\Z$ with Borel-torus pair $\B \supseteq \T$,
\item $\infty_0 \in X(k)$ is a rational point at which $G$ is split, with associated geometric point $\infty \in X(\ol{k})$,
\item $\rho_G : \pi_1(X, \infty) \to \Out(\G)$ is a \emph{pointed} $\Out(\G)$-torsor, so that $G = \rho_G \times^{\Out(\G)} \G$ is a pinned quasi-split reductive group scheme over $X$ (in the sense of \S \ref{subsect.quasisplit}), equipped with an isomorphism $G_\infty \cong \G$,
\item we denote $B \supseteq T$ the Borel-torus pair in $G$ over $X$ defined by twisting $\B \supseteq \T$,
\item $(G, T)$ is of good characteristic in the sense of \ref{def.good.characteristic},
\item $D$ is a fixed effective divisor with $\deg(D) > 2g$ disjoint from $\infty$,
\item $f : \cM \to \cA$ denotes the $D$-twisted Hitchin fibration for $G, X$ (cf.\ Definition \ref{def.Hitchin.fib}),
\item $f : \ctM \to \ctA$ denotes the cameral rigidification over $\infty$ of the $\infty$-regular semisimple Hitchin fibration $\ctA \to \cA^{\heartsuit} \subseteq \cA$ (cf.\ Definition \ref{def.Hitchin.base}),
\item $\xi \in \a_T$ and $f^{\xi} : \ctM^{\xi} \to \ctA$ is the corresponding $\xi$-semistable truncation (cf.\ Definitions \ref{def.geo.stab.convex}, \ref{def.rel.stab.convex}),
\item $f^{\xi} : \bM^{\xi} \to \bA$ denotes the restriction of the $\xi$-semistable truncation to the compatible locus of Definition \ref{def.comp.loc}.
\end{itemize}
In the arithmetic applications that constitute the final part of this work, the geometric constructions are all taken over a finite field $k_0$ with algebraic closure $k = \ol{k}_0$.

\part{Background and the definition of weighted stability}

\chapter{Generalities on reductive group schemes}\label{sect.rel.reduct.gps}

One of the primary goals of this work is to adapt the constructions and arguments of \cite{CLI, CLII} to the case of quasisplit reductive group schemes $G \to X$ over a curve which are neither split nor arise by base change from the scalar field. Thus, it is useful to collect in one place the various results which allow us to work with such schemes combinatorially. In this section, $S$ is an arbitrary scheme except where stated otherwise, while $X$ is always a smooth geometrically integral curve over a field $k$. A \textit{reductive group scheme} $G \to S$ or $G \to X$ is a smooth group scheme, affine over the base, whose fibers over algebraically closed geometric points are connected groups with trivial unipotent radical. We assume always that the base scheme is of good characteristic in the sense defined below. 

\section{Characteristic assumptions.}\label{subsect.char.assumptions}

In this work, \emph{good characteristic} relative to $G/S$ is a stronger condition than the condition by that name elsewhere in the literature, closer to ``very good characteristic.''  For an arbitrary scheme $S$, say that a torus $T/S$ has \emph{prime-to-$p$ monodromy} if, for each component of $S$, there is a finite group $\Gamma$ whose order is invertible on $S$ and a finite \'etale $\Gamma$-torsor $S' \to S$ which splits the torus $T$. 
\begin{definition}\label{def.good.characteristic}
We say that a pair $(G,T)$ consisting of a reductive group scheme $G/S$ and a (fiberwise) maximal torus $T/S$ is of \emph{good characteristic} if for every $s \in S$ of residue characteristic $p$,
\begin{enumerate}
    \item $p$ does not divide $|W(G_{\ol{s}}, T_{\ol{s}})|$, and
    \item $T$ has prime-to-$p$ monodromy. 
\end{enumerate}
\end{definition}

In our cases of primary interest, where $S = X$ is a smooth curve or a local equicharacteristic trait $\Spec k[\![t]\!]$ over $k$ a field of characteristic $p$, the second condition is equivalent to
\begin{enumerate}[start=2, label=({\arabic*}$'$)]
\item the finite monodromy action $\rho_T : \pi_1(S, \infty)\to \Aut_{\Z}(\ul{X}^\star(T))$ has order prime to $p$,
\end{enumerate}
and so our combined condition could be stated as 
$$
    p \nmid |W_G| |\im(\rho_T)|.
$$

In reality more mild assumptions than (1) would perhaps suffice; our primary requirements are 
\begin{enumerate}
    \item That the Chevalley--Todd--Shepard theorem should hold, 
    \item that a Kostant section should exist,
    \item that $p$ does not divide $\pi_1(G)$, and
    \item that centralizers of semisimple elements of the Lie algebra should be (connected) Levi subgroups.
\end{enumerate}
In \S \ref{sect.Deligne.Mumford} and after, we will also need the tameness assumption for the torus. For our intended arithmetic applications, it suffices to show the results for sufficiently large characteristic, so we are free to avoid any finite number of characteristics. 

\begin{proposition}[{{{\cite[Theorem 3.14]{steinbergtors}, \cite[Proposition 14.1.1]{gkmcodims}}}}]\label{prop.goodcharcentralizers}
    Let $G$ a reductive group over a field $k$ of characteristic $p$. Then $Z_G(Y)$ is (geometrically) connected for every semisimple $Y \subseteq \fg(k)$ if and only if $p$ is not a torsion prime for $G$. Moreover, if $p$ does not divide the order of the Weyl group, then $Z_G(Y)$ is a Levi $L \supset T$ for every semisimple $Y \in \ft(k)$.
\end{proposition}
In fact, avoiding torsion primes and $p = 2,3$ suffices for this result.

\begin{remark}
    In \cite{ngofl}, the author takes $p > 2h+2$, where $h$ is the maximal height of a root, solely to ensure the existence of a Kostant section. By \cite[Theorem 3.2.2]{richekostant}, our characteristic assumptions are sufficient, see also \cite{bouthierprodform}, and indeed it is clear that the results of \cite{ngofl} hold under the weaker characteristic assumptions.
\end{remark}

\section{Maximal Tori.} For tori, we retain the usual (strong) definition of split: a torus $T \to S$ is \emph{split} if it is isomorphic to $D_S(M)$ for $M$ a $\Z$-lattice. In particular, the rank of a split torus is required to be constant over the base.

\begin{proposition}
    Assume $S$ is connected, Noetherian and normal, and $G \to S$ is a reductive group scheme. Every maximal torus $T \subseteq G$ is split by a \emph{finite}, \emph{connected} \'etale cover. 
    \begin{proof}
        The assumptions on $S$ allow us to work at the generic point, where every \'etale cover is refined by a finite \'etale one. Let $H \subseteq \pi_1(\eta, \ol{\eta})$ be the corresponding open subgroup of finite index. Its image under the surjection $\pi_1(\eta, \ol{\eta}) \to \pi_1(S, \ol{\eta})$ corresponds to a finite connected \'etale cover $S' \to S$ which splits $T$ at the generic point. But this implies that $T$ is split over all of $S'$. 
    \end{proof}
\end{proposition}

\section{Quasi-split groups, pinnings, and forms.}\label{subsect.quasisplit}
For relative reductive group schemes, we generally follow the notation and conventions of \cite{conradrgs}, with the following major exception: following \cite{ngofl}, we use the term \textit{quasi-split} to indicate that a reductive group scheme $G \to S$ arises as the $S$-form of a reductive $\Z$-group $\G$ associated to an \'etale $\Out(\G)$-torsor $\rho_G$. 

\begin{definition}
    Let $G \to S$ a reductive group scheme. A \textit{pinning} of $G$ is a triple $(T,B,s)$, where $T$ is a maximal torus over $S$, $B \supseteq T$ is a Borel (i.e.,\ all its geometric fibers are Borels), and $s \in H^0(S, \Lie B)$, such that over an \'etale-local cover $S' \to S$ the torus splits and $s = \sum_{\alpha \in \Delta} s_{\alpha}$ (a sum over the simple roots for $B$) for some nowhere-vanishing sections $s_{\alpha}$ in the root spaces $\fu_{\alpha}$. A pinning is called \textit{split} if the torus $T$ is split. 
\end{definition}

Any pinning $(\T, \B, \mathbb{s})$ of $\G$ defines a section $\varepsilon$ of the natural surjection $\Aut(\G) \to \Out(\G)$ consisting of those automorphisms which preserve the pinning; twisting $\G$ along the $\Aut(\G)$-torsor $\varepsilon_\star(\rho_G)$ gives rise to a reductive group scheme $G \to S$ which necessarily has a corresponding pinning, in particular $G$ admits a Borel-torus pair.

In other words, we have shown that our notion of quasi-split is equivalent to ``pinned quasi-split.''
\begin{proposition}
    If $G \to S$ is a form of split group $\G$ defined by an $\Out(\G)$-torsor $\rho$ on $S$, then each pinning $(\T, \B, s)$ induces one for $G$ by splitting. Conversely, if a form $G \to S$ of $\G$ admits a pinning then it is defined by an $\Out(\G)$-torsor.
\end{proposition}

This definition is strictly stronger than defining quasi-split to mean only the existence of a Borel, as in \cite{conradrgs}: reductive group schemes with Borel subgroups, even over curves, may only admit maximal tori Zariski-locally, and even if they admit a maximal torus, may require further Zariski-localization in order to obtain a full pinning. 

Given a geometric point $\infty \in X$, such a torsor corresponds to a conjugacy class of \emph{pointed twisting data}
$$
    \rho_{G}^\bullet : \pi_1(X, \infty) \to \Out(\G),
$$
with a particular such datum determined uniquely by picking a point $\infty_{\rho_G}$ in the fiber of $\rho_G$ over $\infty$. To alleviate notation, we will denote a pointed twisting datum also by $\rho_G$ when no confusion can arise.  

Over fields, the forms of a reductive group $G$ are partitioned into inner equivalence classes, and within each inner class is a unique quasi-split form. Over a broader base, uniqueness of the quasi-split inner form fails easily, but existence holds even with our stronger definition of quasi-split. 

\begin{proposition}[{{\cite[Prop. 7.2.11]{conradrgs}}}]
    Given a reductive group scheme $G \to S$, there is a (not necessarily unique) inner form $G^{qs}$ of $G$ which is quasi-split in our sense, and in particular admits a Borel-maximal torus pair $(B, T \subseteq B)$. 
\end{proposition}

Note that when $k$ is algebraically closed or finite (our two primary cases of interest), the condition that $G \to X$ admit a Borel is automatic, as in \cite[Prop. 2.13]{PappasRapoportGBundles}. This does not, however, imply that inner classes over $X$ are necessarily trivial. To be precise, they are classified by the following result.

\begin{proposition}[{{\cite[Prop. 7.2.2, Ex. 7.3.8]{conradrgs}}}]
    Let $\G$ be a split reductive group over $S$. The fibers of the map $H^1(S, \Aut_{\G/S}) \to H^1(S, \Out_{\G/S})$ are precisely the $S$-forms of $\G$ which are inner forms of each other.
\end{proposition}

\section{Character and cocharacter lattices.}
\begin{definition}
    For any algebraic group scheme $G \to S$, one can define the \'etale sheaf $\ul{X}^\star(G) := \ul{Hom}_S(G, \G_{m,S})$. We denote $X^\star(G) = \ul{X}^\star(G)(S)$ for its \textbf{global sections}, which can be identified with characters of the maximal quotient split constant $S$-torus.
\end{definition}

For general algebraic groups $G$, this sheaf is neither representable nor particularly well-behaved (consider $G = \G_a$). The situation is much better for parabolics $P \subseteq G$ inside reductive group schemes. In \cite[Thm. 5.3.1]{conradrgs}, the structure theory of reductive $G/S$ is used to construct a derived subgroup $D(G) \subseteq G$ with expected properties: 
\begin{itemize}
    \item $D(G)$ is semisimple (i.e., each geometric fiber is semisimple, or equivalently the center is finite over $S$), and its construction is compatible with arbitrary base change.
    \item $D(G) \subseteq G$ is initial among subgroups with the property that $G/D(G)$ is a torus.
    \item For $A_G \subseteq Z(G)$ the maximal central torus, there is a central isogeny $A_G \times D(G) \to G$.
\end{itemize}

\begin{notation}
    Since both will later be useful, we fix the following convention:
    \begin{itemize}
        \item the \textit{maximal central torus} (fiberwise connected center in good characteristic) of a group $H \to S$ will be denoted $A_H$, and
        \item the \textit{maximal $S$-\emph{split} central torus} of a group $H \to S$ will be denoted $A_H^{\spl}$.
    \end{itemize}
\end{notation}

\begin{corollary}
    For $P \subseteq G$ parabolic with Levi quotient $L$, $\ul{X}^\star(P)$ is representable by a separated \'etale group $S$-scheme, also denoted $\ul{X}^\star(P)$, \'etale-locally isomorphic to the group scheme $\Z^r$ for $r = \rank A_L$ (a locally constant function on $S$). 
    \begin{proof}
        By \cite[Cor. 5.2.5]{conradrgs}, parabolic subgroups have a satisfactory notion of unipotent radical $R_u(P)$. In particular, $P \to L := P/R_u(P)$ is initial among maps from $P$ to reductive groups. Compatibility with base change then implies
        \[
            \ul{\Hom}_S(P, \G_{m,S}) = \ul{\Hom}_S(L, \G_{m,S}) = \ul{\Hom}_S(L/D(L), \G_{m,S}).
        \]
        The sheaf of characters of a torus is representable: in the split case by $\Z^r$, for locally constant rank $r = \rank L/D(L)$. The central torus inclusion $A_L\to  L$ descends to an isogeny $A_L \to L/D(L)$, so $\rank L/D(L) = \rank A_L$. For the general case \'etale descent is effective for separated \'etale schemes: indeed, the descent exists as an algebraic space, and any algebraic space which is separated and locally quasi-finite is representable by a scheme, by \cite[\href{https://stacks.math.columbia.edu/tag/03XX}{Tag 03XX}]{stacks-project}.
    \end{proof}
\end{corollary}

The fiberwise dual $\ul{\Hom}(\G_{m,S}, L/D(L))$ of $\ul{X}^\star(P)$ is likewise representable, denoted $\ul{X}_\star(P)$ with global sections $X_\star(P)$. Its monodromy representation as a twist of $\Z^r$ is the contragredient of that for $\ul{X}^\star(P)$, and the pairing of the two is perfect locally, but possibly not on global sections.\footnote{For $M$ a lattice with action of finite group $\Gamma$ and dual $M^\vee$, there is a natural map $(M^\vee)^\Gamma \to \Hom(M^\Gamma, \Z)$ which becomes an isomorphism over $\Q$ but may have finite cokernel.} 

Instead the dual of the global sections admits a characterization in terms of coinvariants, which we state now. We denote $\cF \mapsto H_0(S^{\et}, \cF)$ for the functor of 0-th homology, that is, the \emph{left} adjoint to the constant sheaf functor $\Ab \to \Sh(S, \Ab)$ (in contrast to global sections, which is the right adjoint).

\begin{proposition}\label{prop.dual.lattice.global.invs.coinvs.mod.tors}
    Let $T$ a torus over connected Noetherian normal $S$. Then the dual of the global sections $X^\star(T) = H^0(S, \ul{X}^\star(T))$ is given by
    $$
        \cX_\star(T) := H_0(S, \ul{X}_\star(T))_{\tf}
    $$
    where the subscript $\tf$ indicates the quotient by torsion. If $\ol{s} \to S$ is a geometric point, then $H_0$ is given by the Galois-coinvariants, and we have
$$\Hom(X^\star(T), \Z) = H_0(S, \ul{X}_\star(T))_{\tf} = \ul{X}_\star(T_{\ol{s}})_{\pi_1(S, \ol{s}), \tf}.$$
    
    \begin{proof}
        By the adjunction and double-duality, we have
        \begin{align*}
            \Hom_{\Ab}(H_0(S, \ul{X}_\star(T)), \Z) &= \Hom_{\Sh(S^{\et}, \Ab)}(\ul{X}_\star(T), \Z_S) \\
            &= H^0(S, \sHom(\ul{X}_\star(T), \Z_S)) \\
            &= H^0(S, \ul{X}^\star(T)). 
        \end{align*}
        Of course, $\Hom_{\Ab}(-, \Z)$ kills torsion, at which point both are finite-rank torsion-free and so dualizing gives the Proposition.
        
    \end{proof}
\end{proposition}

We will sometimes need to consider cocharacters up to conjugacy.

\begin{lemma}\label{lem.conj.cochars}
    Let $G \to S$ a reductive group scheme. The \'etale sheafification of the presheaf 
    $$
        S' \mapsto \Hom_{S'}(\G_{m,S'}, G)/\{G\text{-conj.}\}
    $$
    is representable by an \'etale $S$-scheme, \'etale-locally isomorphic to the constant scheme with fiber $X_\star(T)/W_G(T)$.
\end{lemma}

\section{Parabolics and Levis.} 

Under our assumption that the quasi-split group admits a Borel pair, the usual dynamic theory of parabolics and Levis goes through unchanged.
\begin{proposition}\label{prop.parab.levi.basics}
    Assume that the reductive group scheme $G \to S$ admits a maximal torus $T$. Call a parabolic or Levi \textit{semi-standard} if it contains $T$.
    \begin{enumerate}
        \item Each semi-standard parabolic $P$ has a unique standard $L = L_P$ splitting the natural surjection $P \to P/R_u(P)$. 
        \item The pair $(P, L)$ is constructed dynamically (that is, $P = P_G(\lambda), L = Z_G(\lambda)$ arise from the limit construction for a cocharacter $\lambda : \G_{m, S} \to T$.
        \item If $\pi : S' \to S$ is a finite Galois cover with Galois group $\Gamma$ such that $(G_{S'}, T_{S'} = D_{S'}(M))$ is split, then the semi-standard parabolics of $G$ are in bijection with the $\Gamma$-stable parabolics of $G_{S'}$. That is, they correspond to Galois-stable parabolic sets of roots for $G_{S'}$.
        \item If $S$ is a field, the map $L \mapsto A_L^{\spl}$ is a containment-reversing injection from the set of standard Levi subgroups $L \supseteq T$ to the set of split sub-tori of $T$, with one-sided inverse $S \subseteq T \mapsto Z_G(S)$.
    \end{enumerate}
    \begin{proof}
        The claim (1) is \cite[Prop. 5.4.5]{conradrgs} (the uniqueness allows reduction to the split, strictly Henselian case). The claim (2) is \cite[Thm. 7.3.1.(1)]{gillergs}, which holds without the connectedness assumption therein: simply build the desired cocharacter on each connected component separately (it will not be uniform over the base in the sense expected of (co)root data, but still defines a valid cocharacter). The claim (3) follows from the previous parts and \'etale descent. For the fourth claim, if $L \subseteq M$ are standard Levis, then the maximal tori in the centers of each $A_M, A_L$ are contained in $T \subseteq L \subseteq M$. Therefore they satisfy the opposite containment $A_L \supseteq A_M$, and likewise their maximal split subtori. On the other hand, $Z_G(A_L^{\spl}) = L$ over a field by \cite[Thm. 4.15]{BorelTits}. 
    \end{proof}
\end{proposition}

\begin{remark}
    The semistandard Levi subgroups within a quasi-split reductive group $G \to S$ are also quasi-split.
\end{remark}

\begin{proposition}\label{prop.Levis.rel.standard.containment}
    Given two semistandard Levis $L \subset M$ containing $T$, there exist parabolics $P \subset Q$ such that $L$ is a Levi for $P$ and $M$ is a Levi for $Q$.
    \begin{proof}
        By (2) above, all parabolics and Levis arise from the dynamic method. Moreover, just as in the case over a field, the assignments $\lambda \mapsto Z_G(\lambda)$ and $\lambda \mapsto P_G(\lambda)$ are constant on the facets of the Weyl fan (see \S \ref{subsect.cochar.cones.measures} for a definition). They define inclusion-reversing bijections from the set of faces
        defined
        $$C_L = \{v \in X_\star(S)_{\R} : \langle \alpha, v \rangle =0 \; \forall \alpha \in \Phi(L, S)\}$$
        (resp. closures of cones) of the fan to semistandard Levis (resp. parabolics). Then it is clear that $Z_G(\lambda) \supset L$ if and only if $\lambda \in C_L$, and $Z_G(\lambda) = L$ for any $\lambda \in C_L$ avoiding a finite number of hyperplanes. Take $\mu \in C_L$ an integral cocharacter with $Z_G(\mu) = M$. Choose a full-dimensional cone in $C_L$ containing $\mu$ on its boundary, and let $\lambda$ lie in the interior of that cone. Then $P = P_G(\lambda)$ is a parabolic with semistandard Levi $L$, $Q = P_G(\mu)$ is a parabolic with Levi $M$, and $P \subset Q$. 
    \end{proof}
\end{proposition}

\section{Character and cocharacter spaces and their relations.}\label{subsect.char.cochar.spaces}

Until the end of this section, let $G$ a reductive group over connected $S$ with maximal torus $T$. 

\begin{definition}\label{def.char.space}
    For any semi-standard parabolic $P \subseteq G$ with Levi $L$, define the \textit{character space} as the \'etale sheaf
    \[
        \ul{\a}^\star_P = \ul{X}^\star(P) \otimes_{\ul{\Z}} \ul{\R},
    \]
    the tensor product taken as abelian sheaves with the constant sheaf $\ul{\R}$.  Define the \textit{cocharacter space} $\ul{\a}_P$ to be the sheaf of dual spaces, i.e., $\ul{\a}_P = \ul{\Hom}(\ul{X}^\star(P), \ul{\R})$. Write $\a_P^\star = \ul{\a}_P^\star(S)$ and $\a_P = \ul{\a}_P(S)$ for the real vector space of global sections. We will frequently abuse notation and write $\a_L^\star = \a_P^\star$ and $\a_L = \a_P$, reasonable since every character of $P$ factors through $L$.
\end{definition}
These are \'etale-locally isomorphic to the constant sheaf $\ul{\R}^r$ for $r$ the rank of the maximal central torus $A_L$, and monodromy acts by elements of $\GL_r(\Z)$ (its action on $\a^\star_P$ and on $\a_P$ differ by a transpose). These spaces are determined by the maximal central tori.
\begin{proposition}\label{prop.char.space.centr.torus}
    For each semi-standard $(P, L)$ in $G$, we have natural isomorphisms
    \[
        \ul{\a}^\star_P \cong \ul{X}^\star(A_L) \otimes_{\ul{\Z}} \ul{\R}, \qquad \text{ and } \qquad \ul{\a}_P \cong \ul{\Hom}_{\ul{\Z}}(\ul{X}^\star(A_L), \ul{\R})
    \]
    where $A_L$ is the maximal central torus subgroup in $L$. In particular, 
    $$
    \ul{\a}_P^\star(S) \cong \ul{X}^\star(A_L^{\spl}) \otimes_{\Z} \R \qquad \text{ and } \qquad \ul{\a}_P(S) \cong \ul{\Hom}_{\Z}(\ul{X}^\star(A_L^{\spl}), \R).
    $$
    \begin{proof}
        As noted above, the central isogeny $A_L \times D(L) \to L$ gives rise to an isogeny of tori $A_L \to L/D(L)$. The construction of $A_L$ and the isogeny is compatible with base change, so by functoriality of the presheaf tensor product we have a map of presheaf tensor products
        \[
            \ul{\Hom}(L/D(L), \G_{m,S}) \otimes^{pre}_{\ul{\Z}} \ul{\R} \longrightarrow \ul{\Hom}(A_L, \G_{m,S}) \otimes^{pre}_{\ul{\Z}} \ul{\R}.
        \]
        For any $S' \to S$, this is a bijection on the $S'$-sections (indeed, it suffices to work with connected $S'$, where this is obvious). This sheafifies to the desired isomorphism. The isomorphism with $\ul{\a}_P$ follows by dualizing section-by-section.
    \end{proof}
\end{proposition}

These sheaves are compatible with containments of semi-standard parabolics in the following sense, well-known in the split case over a field.
\begin{proposition}\label{prop.charshvs.decomp}
Let $S$ and $(G, T)$ as above, and $P \subseteq Q$ a containment of semi-standard parabolics. Then $\ga_Q^\star$ has a natural structure as a summand in $\ga_P^\star$, and likewise $\ga_Q$ inside $\ga_P$. That is, there are canonical decompositions of sheaves 
\begin{equation}\label{eqn.parabcontaindecomp}
    \ul{\a}_P = \ul{\a}_Q \oplus \ul{\a}_P^Q, \qquad \text{ and } \qquad \ul{\a}_P^\star = \ul{\a}_Q^\star \oplus \left(\ul{\a}_P^Q\right)^\star
\end{equation}
where the decompositions are as real vector spaces for each $S'$. These containments are proper whenever $P \subsetneq Q$.
    \begin{proof}
        The inclusion of spaces comes from the natural restriction map $X^\star(Q) \to X^\star(P)$. This is injective on sections over any $S'$: let $\chi : Q \to \G_m$ a character which sends $P$ to the identity. Such a character corresponds bijectively to $\ol{\chi} : L_Q/D(L_Q) \to \G_m$, as before, which is a morphism between group schemes of multiplicative type. In particular, whether it is trivial can be checked on geometric fibers. So we reduce to the case of an algebraically closed field, where projectivity of $Q/P \subseteq G/P$ and connectedness gives the result. The induced morphism $\a^\star_Q(S') \to \a^\star_P(S')$ is then an inclusion of subspaces for each $S' \to S$. 

        In the other direction, note that if $P \subseteq Q$ with Levis $L_P \subseteq L_Q$, then the corresponding maximal central tori inside the Levis naturally satisfy $A_Q \subseteq A_P$: we can check containment inside of $G$ on geometric fibers, hence reduce to the case of an algebraically closed field. This corresponds to a linear surjection of character lattices $X^\star(A_P) \twoheadrightarrow X^\star(A_Q)$. The previous proposition then shows there is an induced surjective map $\ul{\a}^\star_P \to \ul{\a}^\star_Q$, which is $\R$-linear on sections over any $S' \to S$. 

        The composition $\a^\star_Q \hookrightarrow \a^\star_P \twoheadrightarrow \a^\star_Q$ corresponds over any $S'$ to restricting characters along the chain $Q \supseteq P \supseteq A_P \supseteq A_Q$, which induces exactly the isomorphism of the previous proposition. Hence, we have canonical decompositions
        \[
            \ul{\a}_P^\star(S') = \ul{\a}_Q^\star(S') \oplus \left(\ul{\a}_P^Q\right)^\star(S')
        \]
        for $\left(\ul{\a}_P^Q\right)^\star = \ker(\ul{\a}_P^\star \to \ul{\a}_Q^\star)$ the complement. These decompositions are compatible with change in $S'$, hence give a decomposition of sheaves. Since for any $S'$, the space $\ul{\a}_P(S')$ is dual to $\ul{\a}_P^\star(S')$, we get a dual decomposition for $\ul{\a}_P$.
    \end{proof}
\end{proposition}

Let $G$ be a reductive $S$-group with maximal torus $T$ and let $S' \to S$ be a $\Gamma$-cover which splits $G$. By construction, we have $\ul{\a}_P(S) = \ul{\a}_P(S')^{\Gamma}$, that is, that the rational cocharacter (resp. character) space is Galois-invariants inside the geometric cocharacter (resp. character) space. Dualizing these inclusions, we can also realize the rational spaces as Galois-coinvariants, that is, there exist canonical projections $\ul{\a}_P(S') \to \ul{\a}_P(S)$.

\begin{proposition}
    Let $P$ be a parabolic with Levi factor $L$. The rational cocharacter space $\a_P(S) = X_\star(A_L^{\spl})_{\R}$ is equal to the Galois coinvariants $(X_\star(A_{L,S'})_{\R})_{\Gamma}$.
    \begin{proof}
        Recall that there is a natural $\Gamma$-invariant surjective map $X^\star(A_{L,S'})_{\R} \to X^\star(A_L^{\spl})_{\R}$ given by restriction. The equivalence of categories between character lattices and tori shows the induced map on coinvariants is an isomorphism. 
    \end{proof}
\end{proposition}

\section{Insensitivity of splitting to base field.}

Recall that the \emph{split rank} of a torus $T$ over $S$ is the maximal $n$ such that there exists a subtorus $\G^n_m \subseteq T$ over the base. When $S$ is Noetherian irreducible and normal (hence connected), the anti-equivalence of tori on $S$ and $\pi_1(S, \ol{s})$-lattices gives
\[
    \rank_{\sp} T = \dim_{\Z} L^{\pi_1(S, \ol{s})}
\]
where $L$ is the $\pi_1$-lattice corresponding to $T$. In our primary case of interest, the split rank is insensitive to extensions of scalars.
\begin{proposition}\label{prop.splitrankextensions}
     Let $G$ be a (pinned) quasi-split group scheme over $X$, as in \S \ref{subsect.quasisplit}, with maximal $X$-torus $T$, and assume that there is a point $\infty \in X(k)$ such that $T_{\infty}$ is split. Then there exists a geometrically connected \'etale cover $X' \to X$ splitting $T$. Moreover, the split rank of $T$ is insensitive to any extension of scalars, in the sense that for \textit{any} $L/k$ field extension, the split rank of $T_{X_L}$ is equal to the split rank of $T$ on $X$.
    \begin{proof}
        The torus $T$ corresponds to a discrete $\Gamma = \pi_1(X, \ol{\infty})$-lattice; denote the action $\rho$. The choice of the point $\infty \in X(k)$ gives rise to a splitting of the sequence 
        \[
            \begin{tikzcd}[ampersand replacement=\&]
                1 \ar[r] \& \pi_1(X_{k^s}, \ol{\infty}) \ar[r] \& \pi_1(X, \ol{\infty}) \ar[r] \& \Gal(k^s/k) \ar[r] \& 1
            \end{tikzcd}.
        \]
        Since $T_\infty$ is split, the image of the splitting acts trivially. The kernel of the monodromy $K \subset \pi_1(X, \ol{\infty})$ defines a connected \'etale cover splitting the group. Since $\rho(\pi_1(X, \ol{\infty})) = \rho(\pi_1(X_{k^s}, \ol{\infty}))$, the cover is geometrically connected.
        
        For the latter claim, the split rank can only increase upon base change, so it suffices to show the reverse inequality when $L$ is algebraically closed.  First consider the case $L = \ol{k}$. We will pass through the separable closure $k^s \subseteq L$ of $k$.  In particular, writing $N = \pi_1(X_{k^s}, \ol{\infty})$, we have
        \[
            N/K \cap N \xrightarrow{\sim} NK/K = \pi_1(X, \ol{\infty})/K.
        \]
        In particular, the pullback torus $T_{X_{k^s}}$ corresponds to an isomorphic lattice with an isomorphic action of a finite group, and so the split ranks are equal when $L = k^s$. The passage to $L = \ol{k}$ and then to an arbitrary algebraically closed field is a standard spreading-out argument as in \cite[Prop. 3.2.2]{conradrgs}. Alternatively when $X$ is proper, the \'etale fundamental group is canonically insensitive to choice of algebraically closed field, see \cite[Exp. X, 1.8]{SGAI}.
    \end{proof}
\end{proposition}

\begin{corollary}\label{cor.parab.levi.scalar.extensions}
    Assume that there is a point $\infty \in X(k)$ such that $T_\infty$ is split. Then for any field extension $L/k$, base change along $X_L \to X$ defines a bijection between semistandard parabolics of $G$ and those of $G_{L}$, and likewise for semistandard Levis. For any semistandard parabolic $P$ or Levi $M$, the spaces $\ul{\a}_{P} = \ul{\a}_M$ and $\a_P = \a_M$, their duals, and the maps between them are insensitive to extensions $L/k$ of the base field.
    \begin{proof}
        By the same reasoning as in the previous proposition, the monodromy image of $\pi_1(X, \ol{\infty})$ in the geometric root datum of $(G_{\ol{F}}, T_{\ol{F}})$ remains unchanged under extension of scalars. The Galois-stable parabolic subsets of the root datum are therefore in bijection, so the sets of parabolics are canonically identified at the generic point and thus globally on $X$. Since all semistandard Levi factors are uniquely determined by their parabolic i.e.,\ according to the dynamic description, the same holds for Levis.
    
        For the character spaces, every character of $P$ or $M$ factors through the maximal torus quotient, and every $X$-character through its maximal split quotient. The abelianization is always compatible with base change, and by the proposition above the maximal split quotient is also, so the split rank of $M^{ab}$ and the maps between the character and cocharacter spaces for different $(P,M)$ are canonically identified.
    \end{proof}
\end{corollary}

\section{Dual groups.}\label{subsect.dual.gps}
Fix our quasisplit $G$ and Borel-torus $(B,T)$ over the curve $X/k$. The geometric root datum of $G$ determines a split dual group and torus $(\widehat{G}, \widehat{T})$ over algebraically closed field $E$ (which we will typically take to be $\C \cong \ol{\Q}_l$). Choosing a geometric point $\infty : \Spec K \to X$ equips the dual group with an action of $\pi_1(X, \infty)$, which factors through a finite quotient $\Gamma$. Identifying parabolic subgroups with parabolic subsets of roots (and coroots) defines a bijection
$$
    P \mapsto \widehat{P}
$$
between semistandard $X$-parabolics for $G$ and $\Gamma$-stable semistandard parabolics for $\widehat{G}$. Passing to their semistandard Levi factors, we have a natural bijection 
$$
    T \subset M \subset G \mapsto \widehat{T} \subset \widehat{M} \subset \widehat{G}
$$
between semistandard $X$-Levis (meaning, Levis of semistandard $X$-parabolics) $M$ in $G$ and ($\Gamma$-stable) semistandard Levi factors $\widehat{M}$ of $\Gamma$-stable parabolics in $\widehat{G}$. Each semistandard $X$-Levi $M \subset G$ (that is, Levi factors of $X$-parabolics) defines a $\Gamma$-stable sub-geometric root datum, and under the correspondence above 
$$
    \Phi(\widehat{M}, \widehat{T}) = \Phi(M, T)^\vee;
$$
cf.\ \cite[\S 1]{arthurtransferwoi}. 

Equivalently, one could consider the bijection onto ${}^L M = \widehat{M} \rtimes \Gamma$ Levi factors of ${}^L P = \widehat{P} \rtimes \Gamma$ inside the $L$-group.

While it may seem unnatural to consider only Levi factors of $\Gamma$-stable parabolics rather than the larger class of $\Gamma$-stable Levis which may not admit a $\Gamma$-stable parabolic, this class admits the following simple characterizations.
\begin{proposition}\label{prop.characterizing.dual.Levis}
    The following are equivalent, for $\widehat{M} \supseteq \widehat{T}$ a $\Gamma$-stable semistandard Levi in $\widehat{G}$:
    \begin{enumerate}
        \item $\widehat{M}$ admits a $\Gamma$-stable parabolic (and so is dual to a semistandard $X$-Levi $M$ in $G$),
        \item $\widehat{M} = Z_{\widehat{G}}(\widehat{\lambda})$ for some $\widehat{\lambda} \in X_\star(\widehat{T})^\Gamma$, and
        \item $\widehat{M} = Z_{\widehat{G}}((Z_{\widehat{M}}^{\Gamma})^0)$.
    \end{enumerate}
    \begin{proof}
        Any $\Gamma$-stable parabolic $\widehat{P}$ for $\widehat{M}$ can be expressed geometrically as the attracting locus of some cocharacter $\widehat{\lambda}$, as in \ref{prop.parab.levi.basics}. Summing $\widehat{\lambda}$ over $\Gamma$ gives a $\Gamma$-fixed cocharacter defining $\widehat{P}$, so $(1) \implies (2)$. Given a $\Gamma$-fixed cocharacter with centralizer $\widehat{M}$ its attracting parabolic is $\Gamma$-fixed and a parabolic for $\widehat{M}$, so $(2) \implies (1)$. 

        Given $\widehat{\lambda}$ as in (2), its image lies in $(Z_{\widehat{M}}^{\Gamma})^0$. So
        $$
            \widehat{M} \subseteq Z_{\widehat{G}}((Z_{\widehat{M}}^{\Gamma})^0) \subseteq Z_{\widehat{G}}(\widehat{\lambda}) = \widehat{M}
        $$
        and $(2) \implies (3)$. Conversely, inside the $\Gamma$-fixed cocharacter lattice for $(Z_{\widehat{M}}^{\Gamma})^0$, a sufficiently generic cocharacter (avoiding finitely many root hyperplanes) satisfies $Z_{\widehat{G}}(\widehat{\lambda}) = Z_{\widehat{G}}((Z_{\widehat{M}}^{\Gamma})^0) = \widehat{M}$.
    \end{proof}
\end{proposition}

\begin{proposition}\label{prop.dual.levi.chars}
    For any semistandard Levi $M$ in $G$ and point $\infty \in X$ such that $T_\infty$ is split, there is a natural isomorphism
    $$
        \widehat{T}/Z_{\widehat{M}}^0 = \Hom_{\Z}(X_{\star}(T\cap D(M))_\infty, \C^\times).
    $$
    \begin{proof}
        The claim is dual to 
        $$
            X^\star(\widehat{T}/Z^0_{\widehat{M}}) = X_\star(T \cap D(M))_\infty.
        $$
        The claim can be checked after base change of $T,M,G$ to an algebraically closed field $K$, so we alleviate notation by omitting subscripts and assuming $(G,T)$ are already over $K$. We have
        $$X^\star(\widehat{T}/Z_{\widehat{M}}^0) = \ker(X^\star(\widehat{T}) \to X^\star(Z_{\widehat{M}}^0))$$
        but $Z^0_{\widehat{M}}$ is the dual torus to $M/D(M)= M^{\ab}$, so 
        $$
        X^\star(\widehat{T}/Z_{\widehat{M}}^0) = \ker(X_\star(T)_\infty \to X_\star(M/D(M))_\infty) = X_\star(T \cap D(M))_\infty.
        $$
    \end{proof}
\end{proposition}

\section{Endoscopic groups.}

The fundamental lemma relates orbital integrals for $G$ (over a local field) to orbital integrals for endoscopic groups of $G$. As we will recall below, the Hitchin fibration for endoscopic groups is compatible with that for $G$. Let $(\G, \B, \T, \mathbb{s})$ a pinned split reductive group scheme and $(\widehat{G}, \widehat{B}, \widehat{T}, \widehat{\mathbb{s}})$ the dual pinned quasi-split reductive group. We have therefore a corresponding isomorphism $\Out(\G) \cong \Out(\widehat{G})$.

Write $(G, B, T, s)$ for the pinned quasi-split reductive group over $X$ defined by a torsor $\rho_G \in H^1(X, \Out(\G))$.

\begin{definition}\label{def.end.datum}
    An \textit{endoscopic datum} for $G$ consists of a triple $\fG = (s, \rho_s, \rho_{s} \to \rho_G)$ defined as follows.
    \begin{itemize}
        \item $s \in \widehat{T}$.\footnote{In \cite{geostabpadic}, they prefer the formalism of homomorphisms $\varprojlim_{\gcd(N,p) = 1} \mu_N \to \widehat{T}$, each of which determines a finite-order \'etale subgroup, whose generator would give a \emph{torsion} element. We prefer to follow the conventions of \cite{ngofl, CLII}, and in particular to impose no finite-order requirement in the definition. Finite-order characters will enter in the cohomological $s$-decomposition of Theorem \ref{thm.s.decomp}.} Given such an element, consider its centralizer $Z_{\widehat{G} \rtimes \Out(\widehat{G})}(s)$ inside the semidirect product $\widehat{G} \rtimes \Out(\widehat{G})$. We write its connected-etale short exact sequence as
        \begin{equation}
        \begin{tikzcd}[ampersand replacement=\&]
                1 \ar[r] \& \widehat{G}_{s}' \ar[r] \& Z_{\widehat{G} \rtimes \Out(\widehat{G})}(s) \ar[r] \& \pi_0(s) \ar[r] \& 1
        \end{tikzcd}.\label{eqn.end.ses}
        \end{equation}
        The pinning of $\widehat{G}$ induces a pinning of $\widehat{G}_{s}'$, cf.\ \cite[Lemma 4.4]{geostabpadic}. 
        The natural quotient map $Z_{\widehat{G} \rtimes \Out(\widehat{G})}(s) \to \Out(\widehat{G})$ and conjugation map  $Z_{\widehat{G} \rtimes \Out(\widehat{G})}(s) \to \Out(\widehat{G}_{s}')$ are both maps to discrete constant \'etale groups, hence descend to $\pi_{0}(s)$:
        \[
        \begin{tikzcd}[ampersand replacement=\&]
            \& \pi_0(s) \ar[dl, "o_{\widehat{G}_s'}(s)" above left] \ar[dr, "o_{\widehat{G}}(s)" above right] \& \\
            \Out(\widehat{G}_{s}') \& \& \Out(\widehat{G})
        \end{tikzcd}.
        \]
        \item $\rho_{s}$ is a $\pi_0(s)$-torsor, and 
        \item $\rho_{s} \to \rho_G$ is an $o_{\widehat{G}}(s)$-equivariant map. Equivalently, there is a canonical isomorphism $\rho_{s} \times^{\pi_0(s)} \Out(\widehat{G}) \to \rho_G$. 
    \end{itemize}
\end{definition}

In other words, $\rho_s \to \rho_G$ is a reduction of scalars of $\rho_G$ to $\pi_0(s)$. When no ambiguity can arise we will abbreviate an endoscopic datum as $(s, \rho_s)$, with the structure of $\rho_s$ as a reduction of $\rho_G$ implicit.

\begin{definition}
    Given $s$, the \emph{split endoscopic group} $\G_{s}'$ is the pinned Langlands dual to $\widehat{G}_{s}'$, and in particular we have an isomorphism $\Out(\G_s') \cong \Out(\widehat{G}_{s}')$. The \emph{endoscopic group} for $G$ corresponding to the datum $(s, \rho_s, \rho_s \to \rho_G)$ is the quasisplit $X$-form $G_{\rho_s}'$ of $\G_{s}'$ defined by twisting along $$\rho_{G_{s}'} := o_{\widehat{G}_s'}(s)_\star(\rho_{s}),$$ viewed as a torsor for $\Out(\G_{s}') \cong \Out(\widehat{G}_{s}')$. We will denote it as $G_s'$ when no confusion can arise, but it depends on the torsor $\rho_s$. 
\end{definition}

We will frequently work with a pointed version of this construction, which renders the third datum redundant. Assume now that we have chosen $\infty : \Spec \ol{k} \to X$ and a pointed twisting datum $\rho_G : \pi_1(X, \infty) \to \Out(\G)$ defining $G$.

\begin{definition}\label{def.pointed.end.datum}
    A \emph{pointed endoscopic datum} for $\rho_G$ is a pair $(s, \rho_{s}^{\bullet})$, where
    $s \in \widehat{T}$ as above, and 
    $$
        \rho_s^\bullet : \pi_1(X, \infty) \to \pi_0(s)
    $$
    is a homomorphism refining $\rho_G^\bullet$, i.e.,
    $$
        \begin{tikzcd}[ampersand replacement=\&]
            \pi_1(X, \infty) \ar[r, "\rho_{s}^\bullet"] \ar[dr, "\rho_{G}^\bullet" below left] \& \pi_0(s) \ar[d, "o_{\widehat{G}}"] \\
            \& \Out(\G)
        \end{tikzcd}
    $$
    Given such a datum, the endoscopic group $G_{s}'$ is the form of $\G_{s}'$ defined by the pointed twisting datum $o_{\widehat{G}'_s} \circ \rho_{s}^\bullet : \pi_1(X, \infty) \to \Out(\G_{s}')$. 

    We will often write $\rho_s = \rho_s^{\bullet}$ for a pointed endoscopic torsor.
\end{definition}

Since $s \in \widehat{T}$, the dual endoscopic groups $Z_{\widehat{G}}(s)^0$ associated to an endoscopic datum naturally contain $\widehat{T}$. We have the following standard characterization.

\begin{proposition}\label{prop.end.gp.sub.root.system}
    Let $(s, \rho_s)$ be a pointed endoscopic datum for $G/X$. Then the split root system for $\widehat{G}'_s$ with respect to $\widehat{T}$ is 
    $$
        \Phi(\widehat{G}'_s, \widehat{T}) = \{\alpha \in \Phi(\widehat{G}, \widehat{T}) : \alpha(s) = 1\}.
    $$
    In particular, any proper endoscopic group $G'_s$ satisfies $\dim G'_s \le \dim G - 2$.
    \begin{proof}
        By construction, $\widehat{G}'_s = Z_{\widehat{G}}(s)^0$ contains $\widehat{T}$. The Lie algebra of the centralizer is a direct sum of weight spaces, and by definition contains $\widehat{\mathfrak{t}}$ and exactly those weight spaces on which $s$ acts trivially, i.e.\ $\alpha(s) = 1$. But $Z_{\widehat{G}}(s)^0 \neq \widehat{G}$ (equivalently, $G'_s \neq G$) if and only if $s$ is noncentral and we are missing at least one $\alpha$. But if $\alpha(s) \neq 1$ then the same is true for $-\alpha$, so 
        $$\dim G'_s = \dim \Lie \widehat{G}'_s \le \dim \Lie \widehat{G} - 2 = \dim G - 2$$
    \end{proof}
\end{proposition}

Given pointed datum $\rho_G$ and pointed endoscopic datum $\rho_s$ defining $X$-endoscopic group $H$, our pinned representatives of outer automorphisms for $G$ and $H$ give rise to two $\pi_1(X, \infty)$ actions on the cocharacter lattice of the identified geometric fibers $T_\infty := T_{H, \infty} \cong T_{G, \infty}$,
$$
    X_\star(T_{H, \infty}) \cong X^\star(\widehat{T}) \cong X_\star(T_{G, \infty}).
$$
These monodromy actions need not correspond directly, but they do correspond up to a Weyl twist, in the following sense.

\begin{proposition}[{\cite[\S 1.9, 1.10]{ngofl}}]\label{prop.Weyl.functor.endo.gps}
    Let $(s, \rho_s)$ be a pointed endoscopic datum for $G \to X$, with endoscopic group $H \to X$.  There is a natural injective homomorphism $$j_s : \W_H \rtimes \pi_1(X, \infty) \to \W_G \rtimes \pi_1(X, \infty)$$ inducing the identity on $\pi_1$ and extending the natural inclusion $$\W_{H} \cong \W_{\widehat{H}} \hookrightarrow \W_{\widehat{G}} \cong \W_G.$$ 
    This homomorphism is equivariant for the actions of both semidirect products on $T_\infty$ and so on $X_{\star}(T_\infty)$. To be precise, if $j_s$ corresponds to cocycle
    $$
        u_s \in Z^1(\pi_1(X, \infty), \W_G),
    $$
    then $\rho_{G'_s} = u_s \cdot \rho_G$ on $X_{\star}(T_\infty)$.

    After passage to any finite $\Gamma$ through which both actions factor, we can realize this homomorphism
    $$
        \W_H \rtimes \Gamma \to \W_G \rtimes \Gamma.
    $$
\end{proposition}

The cocharacter spaces satisfy the same decompositions for general endoscopic groups that they do for Levis. Indeed, for the sub-root system of Proposition \ref{prop.end.gp.sub.root.system} above associated to endoscopic datum $(s, \rho_s)$ with endoscopic group $H$, identifying the split maximal torus inside the endoscopic group gives direct sum decompositions
$$
    \ga_H = \ga_G \oplus \ga_H^G,
$$
just as for Levi subgroups: the inclusion $i : \ga_G \hookrightarrow \ga_H$ comes from the inclusion 
$$
    \{x \in X_\star(T_\infty) : \alpha(x) = 0 \; \forall \alpha \in \ul{\Phi}_G\} \hookrightarrow \{x \in X_\star(T_\infty) : \alpha(x) = 0 \; \forall \alpha \in \ul{\Phi}_H\}
$$
and the surjection $q : \ga_H \to \ga_G$ comes from the map
$$
    \ga_T/\R \langle \ul{\Phi}_H^\vee \rangle \twoheadrightarrow \ga_T/\R \langle \ul{\Phi}_G^\vee \rangle.
$$

\begin{proposition}\label{prop.cochar.space.end.grps}
    Both the inclusion $\ga_G \hookrightarrow \ga_H$ and the surjection $\ga_H \twoheadrightarrow \ga_G$ are equivariant for the homomorphism of Weyl groups $\W_H \rtimes \Gamma \to \W_G \rtimes \Gamma$. In particular, we have a corresponding canonical decomposition of rational cocharacter spaces
    $$
        \a_H = \a_G \oplus \a_H^{G}.
    $$
    \begin{proof}
        All the equivariance is clear from the construction of Ng\^{o}'s homomorphism except to observe that under restriction to the centers, although the $H$-action of $\Gamma$ and the $G$-action on $\T \cong T_\infty$ differ by a Weyl-translate, the split Weyl group acts trivially on any central geometric cocharacter.
    \end{proof}
\end{proposition}

\section{Geometric and arithmetic endoscopic data.}\label{subsect.arith.geo.endo.data} The definitions of the previous section are ambivalent about what field $k$ the curve $X$ is defined over. If $k$ is algebraically closed, we call the above data \textit{geometric} endoscopic data. Taking $k_0$ an arbitrary perfect field with algebraic closure $k$ and $X_0$ a curve over $k_0$ with base change $X/k$, one has the canonical short exact sequence
$$
    1 \to \pi_1(X, \infty) \to \pi_1(X_0, \infty) \to \Gal(k/k_0) \to 1,
$$
split if $\infty \in X(k_0)$. In particular, a pointed endoscopic datum is a pair consisting of $s \in \widehat{T}$ and a homomorphism
$$
    \pi_1(X, \infty) \rtimes \Gal(k/k_0) \to \pi_0(s);
$$
such a datum is the same as a pointed geometric endoscopic datum $\pi_1(X, \infty) \to \pi_0(s)$ equipped with $\Gal(k/k_0)$-descent data. For example, one could have a geometric endoscopic datum which is fixed on the nose by $\Gal(k/k_0)$, though not all arithmetic endoscopic data arise in this way.

\section{Twisting and Chevalley spaces.}\label{subsect.twist.chevalley} We recall the elementary theory of the space of ``characteristic polynomials'' $\fc$. In this section let $G_0$ a split reductive group scheme over a field with pinning $(B_0,T_0,s_0)$ and Weyl group $W_0 = \W$. Then the theorem of Chevalley--Todd--Shepard implies an isomorphism of affine quotients
\[
    \fc_0 := \ft_0\gitq W_0 \cong \fg_0\gitq G_0,
\]
and that these are isomorphic to affine space. We call $\chi : \fg_0 \to \fc_0$ and various derived maps \textit{characteristic polynomial maps}, and somewhat abusively call an element of $\fc_0$ a \textit{characteristic polynomial}, though the invariant space may have additional generators (i.e.,\ Pfaffians). The scaling $\G_m$-action on $\fg_0$ commutes with the adjoint action, so descends to an action on $\fc_0$ making $\chi$ equivariant. Explicitly, the action is given by the exponents of the group, the (canonically defined) degrees of homogeneous polynomial generators for $k[\ft_0]^W$.

After the theorem of Kostant (cf.\ \cite{bouthierprodform}) whenever the characteristic $p$ is sufficiently large\footnote{In particular, good in our sense suffices.} there exists a section $\varepsilon : \fc_0 \to \fg_0^{\Greg} \subseteq \fg_0$. The restriction $\chi = \chi^T_G : \ft_0 \to \fc_0$ is finite-flat, and is an \'etale $W_0$-torsor over a nonempty open subset $\fc_0^{\rs} \subseteq \fc_0$. 

Let $(G, T, B, s)$ a pinned reductive group over $X$ defined by twisting $(G_0, T_0, B_0, s_0)$ with $\rho_G \in H^1(X, \Out(G_0))$. As in \cite[\textsection 1.3]{ngofl}, the constructions above can all be twisted by $\rho_G$ to obtain corresponding constructions for $G$. In particular, we have $X$-schemes $\fg$ and $\fc$ (both vector bundles) along with a characteristic polynomial map $\chi : \fg \to \fc$ over $X$, which is $\G_m$-equivariant for an appropriate fiberwise action on $\fc/X$. In particular, given any $\G_m$-torsor $\O(D)$ on $X$, one can twist to obtain
    \[
        \ft_{D} := \ft \times^{\G_m}_X \O(D) \qquad \text{ and } \qquad \fc_{D} := \fc \times^{\G_m}_X \O(D),
    \]
    along with a \emph{characteristic polynomial map} $\chi : \ft_{D} \to \fc_{D}$ which is the finite affine quotient by $W$.

    The scaling $\G_m$-action on $\fg$ likewise defines a twist of the entire Lie algebra $\fg_D := \fg \times^{\G_m}_X \O(D)$, and Chevalley--Todd--Shepard twists to an $X$-isomorphism
    $$
        \fc_D = \ft_D\gitq W \cong \fg_D\gitq G.
    $$

If we have more strongly a \textit{pointed} twisting datum $\rho_G^\bullet : \pi_1(X, x) \to \Gamma \subseteq \Out(\G)$ for $x \in X$, it determines a $\Gamma$-Galois cover $\pi : X' \to X$ as well as a point $x' \in X'$ over $x$ and an isomorphism $G_{X'} \cong \G \times_{\Z} X'$.

\begin{definition}[The Discriminant Divisor]\label{def.discr.divisor}
    For the split group $G_0$, define the \emph{discriminant function} $D_{G_0}$ on $\ft_0$ to be the polynomial function
    $$
        \prod_{\alpha \in \Phi^{G_0}} d \alpha.
    $$
    This is clearly $W_0$-invariant, so descends to a function on $\fc_0$, also denoted $D_{G_0}$. The \emph{discriminant divisor} $\fD_{G_0} \subset \fc_0$ is the vanishing of this function.
\end{definition}
By \cite[Lemme 1.10.1]{ngofl} the discriminant divisor is reduced, with complement $\fc_0^{\rs}$, and stable under the action of $\Out(G_0)$. Thus, the divisor twists to give a corresponding discriminant divisor $\fD_G$ on $\fc_D$ satisfying the same properties.

The Chevalley spaces satisfy a basic functoriality in endoscopic groups. Let $(s, \rho_s)$ be an endoscopic datum for $G \to X$, with endoscopic group $H \to X$. 

\begin{proposition}[{\cite[\S 1.9, 1.10]{ngofl}}]\label{prop.chevalley.functor.endo.gps}
    Under the homomorphism $$\W_H \rtimes \pi_1(X, \infty) \to \W_G \rtimes \pi_1(X, \infty)$$ of Proposition \ref{prop.Weyl.functor.endo.gps}, identifying split maximal tori in $H$ and $G$, twisting defines a natural $X$-map
    $$
        \chi^H_G : \fc_{H, D} \to \fc_{D}.
    $$
    \begin{proof}
        As in the previous proposition the homomorphism is equivariant for the actions on $T_\infty$ and hence on its Lie algebra. After quotienting by $\W_G$ the correction factor disappears, and so the natural map
        $$
            \Lie(T_\infty)/\W_H \to \Lie(T_\infty)/\W_G
        $$
        is equivariant for the $\pi_1(X, \infty)$ actions on each. Twisting gives the claimed map of Chevalley spaces.
    \end{proof}
\end{proposition}

\begin{definition}\label{def.resultant.divisor}
    The \emph{resultant divisor} $\fR_H^G$ for $(s, \rho_s)$ on $\fc_{H,D}$ is defined by 
    $$
        (\chi^H_G)^\star \fD_G = \fD_H + 2\fR_H^G.
    $$
\end{definition}

This satisfies a simple transitivity property:
\begin{proposition}\label{prop.resultant.transitive}
    Let $(s, \rho_s)$ be an endoscopic datum for $G$ defining a group $H$, and let $(s_H, \rho_{s_H})$ be an endoscopic datum for $H$ defining a group $K$. Then under the maps of Proposition \ref{prop.chevalley.functor.endo.gps},
    $$
        \fR_K^G = \fR_K^H + (\chi^K_H)^\star \fR_H^G.
    $$
    \begin{proof}
        Immediate from the definition and functoriality of pullback for the discriminant.
    \end{proof}
\end{proposition}

After passing to a splitting cover and identifying maximal tori, the resultant divisor is given by
$$
    \fR_H^G = \on{div} \left(\prod_{\alpha \in \Psi} d\alpha \right),
$$
where $\Psi \subset \Phi^{G} - \Phi^{H}$ contains one of $\alpha, -\alpha$ for each pair of roots, cf.\ \cite[Lemme 1.10.2]{ngofl}. In particular, when $H = M$ happens to be a Levi subgroup of $G$, such a choice is given by a parabolic $P = MN$ for $M$; then
$$
    R_{M,P}^G := \prod_{\alpha \in \Phi(N,T)} d \alpha
$$
is locally an equation defining $\fR^G_M$. It is independent of $P$ up to sign.

We obtain a simple determinantal formula for the resultant relative to a Levi, whose proof follows immediately by descent from the split case and we leave as an exercise.
\begin{proposition}\label{prop.unip.formula.for.resultant}
    Let $Y \in \fm$ be any section of the Lie algebra of $M$ and $P = MN$ a Levi for $M$. Then the function $\det(\ad(Y)|_{\fn})$ is $M$-invariant, so descends to a corresponding section $R_{M,P}^G$ of $\fc_{M,D}$. We have
    $$
        \on{div}(R_{M,P}^G) = \fR_M^G.
    $$
\end{proposition}

\chapter{The Hitchin base and ellipticity}\label{sect.Hitchin.base}

In this section, let $k$ be an arbitrary field. Let $X$ be a smooth, proper, geometrically integral curve over $k$, and let $G \to X$ be a quasi-split reductive group scheme (cf.\ \S \ref{subsect.quasisplit}). We assume $k$ is of good characteristic for $G$. We assume there exists (and fix) $\infty \in X(k)$ such that $G_\infty$ is split---this will ensure that our combinatorics in terms of Levi subgroups remains unchanged when extending $k$, cf.\ Corollary \ref{cor.parab.levi.scalar.extensions}. 

We denote by $W$ the Weyl group, viewed as a finite \'etale group scheme on $X$. So $W(X) = W(F)$ is the Weyl group of $G_F$, while $W(\infty) = \W$ is the split Weyl group.

\section{The Hitchin base.}

Fix an effective divisor $D$ disjoint from $\infty$ on $X$. We will often assume $\deg(D) \gg 0$; until the arithmetic globalization arguments in \S \ref{subsect.globalize.param} and \S \ref{subsect.nonstd.globalization}, $\deg(D) > 2g$ is sufficient, and should be taken as the running assumption unless stated otherwise.

\begin{definition}\label{def.Hitchin.base}
    For each semistandard Levi $M \subseteq G$, define the \emph{characteristic scheme} or \emph{Hitchin base} $\cA_M$ as the $k$-scheme corresponding to the vector space $H^0(X, \fc_{M,D})$. Define the \emph{cameral base} over $\infty$, abusively denoted $\ctA_M$, as the fiber product
    \[
        \begin{tikzcd}[ampersand replacement=\&]
            \ctA_M \ar[r] \ar[d, "\mathcal{W}_{M,\infty}"] \& \ft_{D}^{\Greg} \ar[d, "\chi_M"] \\
            \cA_M \ar[r, "\on{eval}_\infty"] \& \fc_{M,D}
        \end{tikzcd}
    \]
    That is, an $S$-point of $\ctA_M$ consists of a pair $(a,t)$, where $a \in H^0(X_S, \fc_D(\cE))$ and $t \in \ft_D^{\Greg}(S)$ (these are points in the \emph{total space} of the vector bundle $\ft_D$), such that
    \[
        a(\infty_S) = \chi_M(t).
    \]
    Note this condition implies that $t$ lives over the constant $S$-point $\infty_S$. Note also that the regularity condition implies that $\ctA_M$ depends upon $G$, not just $M$.\footnote{For example, if $M = T$ is the standard torus inside $G = \GL_n$, then the condition of $T$-regularity would be trivial, while the condition of $G$-regularity means the components of $\ft_D$ are pairwise distinct.} When no confusion is possible, we write $\cA := \cA_G$ and $\ctA := \ctA_G$.

    Equivalently, defining the \emph{cameral cover} $X_a \to X \times S$ by the Cartesian square
    $$
        \begin{tikzcd}[ampersand replacement=\&]
            X_a \ar[r] \ar[d] \& \ft_{D} \ar[d, "\chi_G^T"] \\
            X \times S \ar[r, "a"] \& \fc_{G,D}
        \end{tikzcd},
    $$
    the data $t$ consists of a geometric point $\wtilde{\infty} \in X_a$ lying over $\infty$. 
    
    The map $\mathcal{W}_{M,\infty}: \ctA_{M} \to \cA_M$ is a $W_M(\infty)$-Galois finite \'etale cover of the open locus $\cA^{\infty}_M \subset\cA_M$ consisting of points which are $G$-regular semisimple at $\infty$ (and hence generically). In particular, it is a smooth $k$-scheme. We record also $\cA_M^{\heartsuit} \subset \cA_M$, the larger open subset of points which are generically $G$-regular semisimple.
\end{definition}

\begin{proposition}[{{{\cite[Proposition 3.4.1]{CLI}, \cite[Lemme 4.13.1]{ngofl}}}}]\label{prop.Hitchin.base.dim}
    If $\deg(D) > 2g-1$, the schemes $\cA_M$ and $\ctA_M$ are irreducible of dimension
    \[
        \dim(\cA_M) = \dim(\ctA_M) = \deg(D)(\# \Phi^M/2 + \on{rk} M) + (1-g) \on{rk} M
    \]
    where $\Phi^M$ is the set of roots of the split form of $M$.
    \begin{proof}
        We recall the proof. For irreducibility one remarks that the map forgetting $a$ is a surjective open map onto $\t^{\Greg}$ with fibers isomorphic to affine spaces by the vanishing of $H^1(X, \fc_{M,D}(-\infty))$, which follows from the degree hypothesis. To compute the dimension of $\cA_M$, one can use a splitting Galois cover $\pi : X' \to X$. After pullback one has $\pi^\star \fc_D = \bigoplus_{i = 1}^r \O(D)^{\otimes e_i}$, where the $e_i$ are the exponents for the split group $\G$. The degree condition implies the vanishing of $H^1(X', \pi^\star \fc_{M,D})$ and thus its summand $H^1(X, \fc_{M,D})$. The simple form of Riemann--Roch on $X$ then gives the desired formula, noting that 
        \[
            e_1 + \cdots + e_r - r = \sum_{i = 1}^r (e_i - 1) = \# \Phi^M/2
        \]
        Since $\ctA_M$ is an \'etale cover of a nonempty open in affine space, its dimension is the same. 
    \end{proof}
\end{proposition}

\section{Example: the quasi-split unitary group.}\label{subsect.ell.conv.fails}

Our goal is to adapt the stratification of the Hitchin base discussed in \cite[\S 3]{CLI} to the relative quasisplit setting. Many of the arguments translate directly. In particular, the ``split at $\infty$" assumption ensures that $T_{\Spec \O_\infty}$ is split (see \cite[Tag 0BQC]{stacks-project} and \cite[B.3.6]{conradrgs}), so when working at $F_{\infty}$ nothing changes from \cite{CLI}. 

The biggest new ingredient in this part of the theory comes from the fact that the rigidifying data $t$ is fundamentally overdetermined. In order to explain, we illustrate the failure of the converse of \cite[Corollaire 3.9.2]{CLI} in the case of the unitary group. Let $G = U(n)$ be the quasi-split unitary group associated with \'etale double cover $\pi : X' \to X$ (Galois element $\sigma$) of curves over a field $k$. This can be presented as the subgroup of $\Res_{X'/X} \GL_n$ consisting of matrices $A$ satisfying 
$$
    \ol{A}^T \Phi A = \Phi, \qquad \Phi = \begin{pmatrix}
        & & 1 \\
        & \adots & \\
        1 & &
    \end{pmatrix},
$$
and its Lie algebra with the matrices which are skew-Hermitian for the same form.

Denote the extension of function fields corresponding to by $F'/F$. With respect to the presentation above, $T$ denotes the standard diagonal maximally split maximal torus. Pick $\infty \in X(k)$ a point which splits in $X'$, denoting the fiber $\pi^{-1}(\infty) = \infty' \sqcup \infty''$. So $G_\infty$ is a split unitary group over $k$, isomorphic to $\GL_{n,k}$ by a choice of basepoint $\infty'$ or $\infty''$ over $\infty$, and the maximal torus pulls back to the diagonal maximal torus $T'$. The isomorphism is given by choosing one of the points in the fiber, without loss of generality $\infty'$, and restricting sections to it.

\begin{example}
    Let $n = 3$ and $M = T$, and let $(a,t)$ be a pair such that $a$ admits a reduction to $M = T$. That is, $a_M$ consists of an ordered factorization $(X - r_1)(X - r_2)(X - r_3)$ with $r_i \in H^0(X', \O(D))$ which satisfy $\sigma^\star(r_i) = - r_{4-i}$. This skew-symmetry condition implies that the values taken by the $r_i$ at $\infty''$ are determined by those at $\infty'$, and the same is true for $t$. we must have
    \[
        \prod_i (X - r_i|_{\infty'}) = \prod_i (X - t_i|_{\infty'})
    \]
    so there exists a permutation $\pi' \in W_G(\infty) = S_3$ such that $r_i|_{\infty'} = t_{\pi'(i)}|_{\infty'}$. However, this is data for the split group $G|_{\infty} \cong \GL_3|_{\infty'}$, so there is no requirement that this permutation lie in the \textit{global} Weyl group $W_G(X) = W_G(F) = \langle (1 \; 3) \rangle \le W_G(\infty)$.
    
    For a counterexample, choose an etale cover $\pi$ and rational functions $x \in H^0(X', \O(D))^{\sigma = 1} \subset F, y \in H^0(X', \O(D))^{\sigma = (-1)} \subset F'$ such that 
    \[
        x(\infty') = 1 = y(\infty') \quad \text{ and so } \quad x(\infty'') = 1 = - y(\infty'').
    \]
    Consider the elements
    \[
        a(\lambda) = (\lambda + x - y)(\lambda - y)(\lambda - x - y), \qquad t(\infty') = (1,2,0), \quad t(\infty'') = (0, -2,-1).
    \]
    We have a compatible pair $(a,t) \in \ctA_G$ because 
    \[
        \chi^T_G(t)(\infty') = (X - 1)(X - 2)(X - 0) = a(\infty')
    \]
    and likewise over $\infty''$. The element $a$ is clearly $T$-elliptic on its own (with the factorization given above or in reverse order), but all $T$-factorizations are incompatible with $t$, differing by a permutation $(1 \; 2 \; 3) \in W_G(\infty) \setminus W_G(X)$. Thus, $(a,t) \in \ctA_{G} \setminus \ctA_T$. We claim that the conclusion of the converse of \cite[Corollaire 3.9.2]{CLI} is false for $M = G$: no (regular semisimple) element $Y \in \fg(F)$ with characteristic polynomial $a$ is $G$-elliptic. If it were, then considering the element
    \[
        Z = \begin{pmatrix}
            y - x & & \\
            & y & \\
            &  & y + x
        \end{pmatrix},
    \]
    we proceed as in the proof of \cite[Corollaire 3.9.2]{CLI}. The two elements $Y,Z \in \g(F)$ are conjugate by $g \in G(F^{\sep})$, so their split centralizing tori are isomorphic: $A_Y^{\spl} \cong A_Z^{\spl}$. But $Z \in \ft(F)$, so its split centralizer is $A_T^{\spl}$ and $A_Y^{\spl} \neq A_G^{\spl} = 1$.  
    
    We note that this is still compatible with \cite[Proposition 3.9.1]{CLI}: the element $(a,t)$ is not in $\ctA_T$, and although we have found an element $Z = \on{diag}(y-x, y, y+x)$ such that $\chi^T_G(Z) = a_\eta$, the element $t_a$ agrees with $\on{diag}(y, x+y, y-x)$ at $\infty'$, and these are certainly not conjugate in $T(F_{\infty}) = T'(E_{\infty'})$. 
\end{example}

The essential feature of the above example is that there are $\# W_G(\infty)$ possible orderings of the eigenvalues $t$, only $\# W_G(F)$ of which are \textit{globally compatible} in the sense that they render $(a,t)$ in the image of $\ctA_T$. For the other values of $t$, the pair $(a,t)$ should still be considered as $T$-elliptic (and indeed the image of a Kostant section in the function field will be $T$-elliptic), but we must introduce a twist living in the split Weyl group $W_G(\infty)$. Indeed, any Higgs bundle $m$ with $f(m) = a$ will have a $W_G(F)$-torsor of Borel reductions, regardless of the ordering $t$. This motivates the formalism of the next section.

\section{Weyl-twisted stratification and the compatible locus.}

In this section, we construct the finest stratification of the Hitchin base over which $\xi$-semistability is well-behaved. The combinatorics of this stratification is governed by semistandard Levis $X' \to X$. We warn the reader that while the stratification is stable under extending the scalar field $k$, it is \textit{not} stable under a general \'etale base change $X' \to X$. This can be seen from the intervention of the endoscopic loci defined below: if $H_s$ is an elliptic endoscopic group whose pullback to a splitting cover becomes a Levi subgroup (for example, $H_{s} = U(1) \times U(1)$ inside $U(2)$), then points in $\ctA_{H_s}$ should still be considered $G$-elliptic, although they will not remain so upon pullback.

For each semistandard Levi $M$ and each element $w \in W_G(\infty)$, there is a twisted closed embedding
\begin{align*}
    \iota_w : \ctA_M &\hookrightarrow \ctA_G \\
        (a_M, t) &\mapsto (\chi^M_G(a_M), w^{-1}(t)),
\end{align*}
defined because $a_M|_{\infty} = \chi_M(t)$ implies 
\[
    \chi^M_G(a_M)|_{\infty} = \chi_G(t) = \chi_G(w^{-1}(t)).
\]
We define $\ctAt{w}{M}$ the image under this closed embedding, characterized as the set of $(a,t)$ such that there exists a Chevalley lift $a_M$ with $a_M|_\infty = \chi_M(w(t))$. 

\begin{proposition}\label{prop.wM.levi.equivs}
The loci $\ctAt{w}{M}$ satisfy the following redundancies:
    \begin{enumerate}
        \item $\ctAt{w}{M} = \ctAt{vw}{M}$ for all $v \in W_M(\infty)$, and
        \item $\ctAt{w}{M} = \ctAt{vw}{vM}$ for all $v \in W_G(F)$.
    \end{enumerate}
    \begin{proof}
        For (1), if $(a,t) = \iota_w(a_M, t')$ such that $a_M(\infty) = \chi_T^M(t')$, then the pair $(a_M, v(t'))$ also lies in $\ctA_M$, and 
        \[
            \iota_{vw}(a_M, v(t')) = (a, w^{-1} v^{-1} v (t')) = (a,t).
        \]

        For (2), let $(a,t) = \iota_w(a_M, t')$. Given $v \in W_G(F) = W_G(X)$, we have a well-defined conjugate semi-standard Levi $L = vMv^{-1}$. The map $\ft_D \to \ft_D$ given by the action of $v$ is equivariant for the map $W_M \to W_L$ given by conjugation by $v$, so descends to a map $ \ol{v} : \fc_{M,D} \to \fc_{L,D}$. Writing $a_L = \ol{v}(a_M)$, we can see
        \[
            a_L(\infty) = \ol{v}(a_M(\infty)) = v(\chi_T^M(t')) = \chi_T^L(v(t')),
        \]
        so $(a_L, v(t'))$ is an element of $\ctA_L$. We have
        \[
            \iota_{vw}(a_L, v(t')) = (\chi_L^G(a_L), w^{-1}(t')) = (a, t)
        \]
        because the following diagram commutes.
        $$
        \begin{tikzcd}[ampersand replacement=\&]
            \ft_D \ar[r] \ar[d, "v"] \& \ft_D\gitq W_M \ar[r] \ar[d, "\ol{v}"] \& \ft_D\gitq W_G \ar[d, equals]\\
            \ft_D \ar[r] \& \ft_D\gitq W_L \ar[r] \& \ft_D\gitq W_G
        \end{tikzcd}
        $$
    \end{proof}
\end{proposition}

\begin{definition}\label{def.wMpair.rels} Define an equivalence relation $\sim$ on the set of pairs $(w,M) \in W_G(\infty) \times \{M \supseteq T\}$ by the two conditions in Proposition \ref{prop.wM.levi.equivs}. Denote $[(w,M)]$ the equivalence class.

Define also a partial order by $[(w,M)] \le [(w, L)]$ for $L \supseteq M$. Note that this is well-defined on equivalence classes for $\sim$.

\end{definition}

We have shown the locus $\ctAt{w}{M}$ depends only on the class $[(w,M)]$. We will show a converse below, that on the \textit{elliptic loci} distinct classes have disjoint images.

\begin{proposition}
    If $[(w,M)] \le [(v,L)]$, then $\ctAt{w}{M} \subseteq \ctAt{v}{L}$.
    \begin{proof}
        Since $\ctAt{v}{L}$ depends only on the class of $(v,L)$, we may assume that $(v,L) = (w, L)$ for $L \supset M$. Then the desired factorization is clear.
    \end{proof}
\end{proposition}

\begin{definition}\label{def.wM.loci}
    Define the $(w,M)$-elliptic locus to be the open in $\ctAt{w}{M}$ which is not contained in any smaller locus. That is, if we first define the $G$-elliptic locus
    $$
        \ctA_{G\el} = \ctA_G \setminus \bigcup_{\substack{(w,M) \\ M \subsetneq G}} \ctAt{w}{M},
    $$
    and $\ctA_{M\el}$ by likewise removing all twisted images of smaller Levis, then 
    $$
        \ctA_{(w,M)\el} := \iota_w(\ctA_{M\el}).
    $$
\end{definition}

It is clear that each $\ctA_{(w,M)}$ is exactly the union of lower elliptic strata. With this definition we claim the elliptic loci stratify the base. 
The proof will appear at the end of the next section, after some preliminaries.
\begin{proposition}
    If $[(w_1,M_1)] \not\sim [(w_2, M_2)]$, then $\ctAe{w_1}{M_1} \cap \ctAe{w_2}{M_2} = \emptyset$. In particular, the $(w,M)$-loci give a stratification of $\ctA_G$.
\end{proposition}

After introducing the complexity of Weyl-twists, we observe that ultimately, it will be unnecessary.
\begin{definition}\label{def.comp.loc}
    Let $e \in W_G(\infty)$ be the identity. Define the \textit{compatible locus} to be
    $$\bA := \ctA^{\comp} = \bigcup_{M} \ctAe{e}{M}.$$
\end{definition}
That is, over the compatible locus we have a well-behaved stratification indexed purely by Levis
$$
    \bA = \bigsqcup_{M} \bAe{M}.
$$

The set of compatible classes is upwards-closed, by the following argument. If $\alpha = [(e,M)]$ is a compatible class and $\beta \ge \alpha$, then we can choose representatives
$$
    \alpha = [(w, N)], \qquad \beta = [(w, L)]
$$
with $N \subseteq L$ and the same element $w \in W_G(\infty)$. Then from the construction of the equivalence relation we have $N = u M u^{-1}$, where $w = b u$ with $b \in W_N(\infty)$ and $u \in W_G(F)$. But since $W_N(\infty) \subseteq W_L(\infty)$, we have equivalence
$$
    [(w, L)] = [(bu, L)] = [(u, L)] = [(e, u^{-1} L u)]
$$
so $\beta$ is compatible. Thus, once we have shown that we indeed have a stratification (Corollary \ref{cor.ell.strat.tw} below), it follows immediately that the union of the compatible elliptic loci is open, stratified by loci $\bAe{M}$ for semistandard Levis $M$. 

We lose no part of the geometry of $\cA^{\infty}$ in restricting to the compatible locus, as every $a \in \cA^{\infty}$ has a $t$ such that $(a, t)$ is compatible, a i.e., every $\infty$-regular Hitchin fiber is visible over the compatible locus.

Compatibility can be characterized as follows, though we won't use this characterization. For any semistandard Levi $M$ and $a \in \cA(k)$, we have a \emph{partial cameral cover} defined by the Cartesian square
$$
    \begin{tikzcd}[ampersand replacement=\&]
        X_{a,M} \ar[r] \ar[d] \& \fc_{M,D} \ar[d, "\chi^M_G"] \\
        X \ar[r, "a"] \& \fc_{G,D}
    \end{tikzcd}.
$$
The data $\chi^T_M(t) \in \fc_M(\infty)$ determines a point $\wtilde{\infty}_M$ in $X_{a,M}$ over $\infty$ (and under $\wtilde{\infty} \in X_{a,T}$). An element $a$ lies in the untwisted image $\iota_e(\ctA_M)$ if and only if $X_{a,M} \to X$ has an $F$-section, and $t$ is compatible with $a$ if and only if $\wtilde{\infty}_M$ lies in the image of such a section, or equivalently, lies on an irreducible component of $X_{a,M}$ which maps birationally onto $X$.

\section{Generic characterization of twisted loci.} 
Using the existence of the Kostant section, we can characterize the elements in the twisted loci at the generic point in the style of \cite[\S 3.9]{CLI}. For the reader familiar with the arguments in \textit{loc.\ cit.}, the new ingredient in this Weyl-twisted setting is carefully tracking \textit{split} central tori.

We recall the basic definitions over a field. 
\begin{definition}
    Let $H$ a reductive group over a field $F$ of good characteristic for $H$, and $Y \in \fh(F)$ an element of the Lie algebra. 
    \begin{itemize}
        \item We say that $Y$ is \emph{$H$-regular} if its centralizer $Z_H(Y)$ has minimal dimension.
        \item We say that $Y$ is \emph{semisimple} if $Y$ is equal to its semisimple part in the additive Jordan decomposition over the algebraic closure $\ol{F}$.
        \item If $Y$ is both $H$-regular and semisimple, then $A_Y := Z_H(Y)^0$ is a maximal torus in $H$ (defined over $F$). Let $A_Y^{\spl} \subseteq Z_H(Y)$ be the maximal $F$-split subtorus. We say that $Y$ is \emph{$H$-elliptic} if $A_Y^{\spl} = A^{\spl}_H$, i.e., if $Z_H(Y)$ is anisotropic mod center.
    \end{itemize}
\end{definition}

\begin{proposition}[{{{cf.\ \cite[Proposition 3.9.1]{CLI}}}}]\label{prop.CLI.3.8.twisted}
    Let $K/k$ a field extension. An element $(a,t) \in \ctA_G(K)$ is in the locus $\ctAt{w}{M}(K)$ if and only if there exist 
    \begin{itemize}
        \item $Y \in \fm(F)$ semisimple $G$-regular, and
        \item $m \in M(F_\infty)$, 
    \end{itemize}
    such that 
    \begin{enumerate}[label=(\alph*)]
        \item $\chi_G(Y) = a_\eta$,
        \item $\on{Ad}(m)Y = w(t_a)$.
    \end{enumerate}
    Here, $t_a \in \ft^{\Greg}(\O_\infty)$ is the unique lift of $t$ such that $\chi_G(t_a) = a_{\O_\infty}$.
    \begin{proof}
        The proof in \textit{loc. cit.} can be directly adapted, using the fact that $T|_{\O_\infty}$ is split.\footnote{This follows from e.g. \cite[Tag 0BQC]{stacks-project} applied with $X = A$ and \cite[B.3.6]{conradrgs}).}
        
        Note that a $t_a$ as in the statement exists by formal \'etaleness of $\chi$ over the regular locus. Assume that $(a,t)$ is in the image of $\cA_M^w(K)$ for some $w$, i.e., there exists $a_M$ a section of $(\m\gitq M)^{G-reg}$ with $\chi^M_G(a_M) = a$ and $a_M(\infty_K) = \chi_M(w(t))$. So we have 
\[
    \begin{tikzcd}[ampersand replacement=\&]
        \Spec(\infty) \ar[dd] \ar[rr] \& \& (\m\gitq M)^{reg} \ar[dd, "\chi^M_G"] \\
        \& \& \\
        \Spec(\O_{\infty})  \ar[rr, "a_\infty"] \ar[uurr, shift left=.75ex,"a_{M,\infty}" above left] \ar[urur, shift right=.75ex, "\chi_M(w(t_a))" below right] \& \& (\g\gitq G)^{reg}
    \end{tikzcd}.
\]
Since the right map is \'etale, repeatedly applying the lifting criterion for unramifiedness implies that $a_{M,\infty} = \chi_M(w(t_a))$. 

Now look over the generic point, where $a_{M, \eta} \in (\m\gitq M)^{reg}(F)$. The Kostant section for $M$ implies that we can find $Y \in \m(F)$, (necessarily regular semisimple because the Kostant section sends $(\m\gitq M)^{reg}$ into the regular semisimple locus), such that $\chi_M(Y)|_{F_\infty} = a_{M,\eta}|_{F_\infty} = a_{M, \infty} = \chi_M(w(t_a))$. So by regular semisimplicity there exists $L/F_{\infty}$ finite Galois and $m \in M(L)$ such that
\[
    \on{Ad}(m) Y = w(t_a).
\]
This descends to $m \in M(F_\infty)$ because, since $T|_{F_{\infty}}$ is split,
\[
    H^1(L/F_\infty, T_{F_\infty}(L)) = 1
\]
by Hilbert 90. 

Conversely, let $(a,t) \in \cA_G(K)$ and suppose there exist $Y \in \m(F)$ $G$-regular semisimple, $m \in M(F_{\infty})$, and $w \in W_G(\infty)$ such that $\chi_G(Y) = a_{\eta}$ in $(\g\gitq G)^{reg}$ and $\on{Ad}(m)Y = w(t_a)$ in $\m(F_{\infty})$. Then the element $Y$ defines a generic factorization of $a : X_K \to (\g\gitq G)_D$ through $(\m\gitq M)_D$: the solid top line of the square
\[
        \begin{tikzcd}[ampersand replacement=\&]
            \Spec(F) \ar[r] \ar[d] \& (\m\gitq M)_D \ar[d, "\chi_G^M"]\\
            X_K \ar[r] \ar[ur, dashed] \& (\g\gitq G)_D
        \end{tikzcd}.
\]
Now we consider the question of extending this morphism to $X_K$, which we do in two steps. First, spreading out (e.g. by finite-typeness of the morphism $\chi_G^M$) implies that the morphism extends to some affine open of $X_K$. Then projectivity of $\chi_G^M$ and nonsingularity of $X_K$ implies that the morphism extends to $a_M : X_K \to (\m\gitq M)_D$. 

So $(a_M, t)$ is our candidate for a preimage under the $w$-twisted embedding. We need to check that $a_M$ and $t$ are $w$-compatible, i.e., that $a_M(\infty_K) = \chi_M(w(t))$ in $(\m\gitq M)^{G-\reg}(K)$. First note that $a_{M,\eta} = \chi_M(Y)$ in $(\m\gitq M)(F) \subseteq (\m\gitq M)(F_\infty)$. The second assumption ensures that $\chi_M(Y) = \chi_M(w(t_a))$, so $a_M$ and $\chi_M(w(t_a))$ agree on the generic point of $\Spec(\O_{\infty})$. By the uniqueness in the valuative criterion of properness, they agree on all of $\Spec(\O_{\infty})$, hence on the special fiber. 
\end{proof}
\end{proposition}

Likewise we have a characterization of the elliptic loci.

\begin{proposition}[{{{cf.\ \cite[Corollaire 3.9.2]{CLI}. }}}]\label{prop.CLI.3.9.twisted}
    Let $K/k$ a field. An element $(a,t) \in \ctA_G(K)$ is in the locus $\ctAe{w}{M}(K)$ if and only if there exist 
    \begin{itemize}
        \item $Y \in \fm(F)$ semisimple $G$-regular $M$-elliptic, and
        \item $m \in M(F_\infty)$,
    \end{itemize}
    such that 
    \begin{enumerate}[label=(\alph*)]
        \item $\chi_G(Y) = a_\eta$,
        \item $\on{Ad}(m)Y = w(t_a)$.
    \end{enumerate}
    \begin{proof}
    Assume there exists $Y,m$ as in the proposition. Assume moreover that $(a,t)$ is in the image of $\ctAt{v}{L}(K)$ for some $L \subseteq M$ and $v \in W_G(\infty)$. We will show $L = M$. By the forward direction of Proposition \ref{prop.CLI.3.8.twisted}, there exist $Z \in \fl(F)$, $l \in L(F_\infty)$ satisfying $\chi_G(Z) = a_\eta = \chi_G(Y)$ and $\on{Ad}(l) Z = v(t_a)$. Since both $Y,Z \in \fg(F)$ are $G$-regular semisimple and have the same characteristic polynomial, they are conjugate by some element $g \in G(F^s)$. 

    The centralizers of $Y,Z$ in $G_F$ are maximal $F$-tori $A_Y, A_Z$. Each element $\sigma \in \Gal(F^{\sep}/F)$ satisfies
        \[
            \sigma(g) g^{-1}Z g \sigma(g)^{-1} = \sigma(g) Y \sigma(g)^{-1} =\sigma(gYg^{-1}) = \sigma(Z) = Z,
        \]
        so $\sigma(g) g^{-1}$ lies in $A_Z(F^{\sep})$. In particular, $\sigma(g)g^{-1}$ commutes with anything in $A_Z(F^{\sep})$, and so the isomorphism $\phi = \Inn(g) : A_{Y, F^{\sep}} \to A_{Z, F^{\sep}}$ is Galois-equivariant:
        \[
            \sigma(\phi(a)) = \sigma(g a g^{-1}) = (\sigma(g) g^{-1})(g \sigma(a) g^{-1})(g \sigma(g^{-1})) = g \sigma(a) g^{-1}.
        \]
        (Note $g \sigma(a) g^{-1} \in A_Z(F^{\sep})$ since $A_Y$ is defined over $F$ and so its set of $F^{\sep}$-points is Galois-stable.) So the isomorphism $\phi$ descends to one defined over $F$. In particular, it induces an isomorphism on the maximal split sub-tori, from $A_Y^{\spl} = A_M^{\spl}$ to $A_Z^{\spl}$. Since $Z \in \fl(F)$, we have $A_Z^{\spl} \supseteq A_L^{\spl}$, so the dimension of $A_M^{\spl}$ is at least that of $A_L^{\spl}$. 

        On the other hand, since $L \subseteq M$ we have $A_M^{\spl} \subseteq A_L^{\spl}$ by Proposition \ref{prop.parab.levi.basics}.(4). We conclude equality holds, and so $M = Z_G(A_M^{\spl}) = Z_G(A_L^{\spl}) = L$. 

        Conversely, assume that $(a,t)$ is in $\ctAe{w}{M}(K) \subseteq \ctAt{w}{M}(K)$, so in particular there exist $Y,m$ as in Proposition \ref{prop.CLI.3.8.twisted}. We show that $Y$ is necessarily $M$-elliptic. Assume for the sake of contradiction $Y$ is not $M$-elliptic, so that the maximal split torus in the centralizer in $G_F$ of $Y$ satisfies $A_Y^{\spl} \supsetneq A_M^{\spl}$. Taking centralizers in $G$ of both sides and applying Proposition \ref{prop.parab.levi.basics}.(4) shows $A_Y^{\spl} \subseteq M$, so up to conjugating $Y$ and $m$ in $M$ we may assume that $A_Y^{\spl} \subseteq T_F$ and the Levi $L_F = Z_G(A_Y^{\spl}) \subsetneq M_F$ is semistandard. 

        Since $\on{Ad}(m) Y = w(t_a)$, we have
        \[
            m A_Y m^{-1} = T_{F_\infty}
        \]
        in $G_{F_{\infty}}$. But on the other hand $T_{F_\infty} = Z_{G_{F_\infty}}(w(t_a))$ and $A_{Y, F_\infty}$ are both maximal (split over $F_\infty$!) tori in $L_{F_\infty}$, so are conjugate by an element $l \in L(F_\infty)$. Since $l m^{-1}$ lies in $M$ and normalizes the split torus $T_{F_\infty}$, it acts on the torus by an element $v \in \mathbb{W}_M = W_M(\infty)$. In particular, we have
        \[
            \on{Ad}(l)Y = \on{Ad}(lm^{-1}) \on{Ad}(m) Y = v w(t_a).
        \]

        Note moreover that since $L = Z_G(A_Y^{\spl}) \supseteq Z_G(Z_G(Y))$ and $\Lie Z_G(Z_G(Y)) = (\Lie G)^{Z_G(Y)}$ (which would hold even in bad characteristic), we have $Y \in \Lie L$. So the pair $(Y, l)$ satisfy the hypotheses of Proposition \ref{prop.CLI.3.8.twisted} with twist $vw$, and $(a,t)$ is in $\ctAt{vw}{L}(K)$.
    \end{proof}
\end{proposition}
\begin{remark}
    Note that it was only necessary to change the twist $w$ by an element of $W_M(\infty)$. 
\end{remark}

Now we can prove the elliptic loci give a stratification.
\begin{proposition}[{{cf.\ \cite[\S 3.10]{CLI}}}]\label{prop.base.disjoint}
    If $(w,M) \not\sim (v,L)$, then 
    \[
        \ctAe{w}{M} \cap \ctAe{v}{L} = \emptyset
    \]
    inside $\ctA_{G}$.
    \begin{proof}
        It suffices to show the intersection is empty on $K$-points for algebraically closed $K$. Moreover, we may replace $K$ with an algebrically closed extension containing $\bar{k}$. Then we apply the generic characterization above. Let $F$ denote the function field of $X_{K}$, so $G_F$ is a (possibly non-split) group over $\Spec F$, and $F^{\sep}/F$ a separable closure. Let $(a,t) \in (\ctAe{w}{M} \cap \ctAe{v}{L})(K)$. So we have regular semisimple elliptic elements $Y \in \fl(F), Z \in \fm(F)$ as well as $l \in L(F_{\infty}), m \in M(F_{\infty})$, satisfying
        \[
        \chi_G(Y) = a_{\eta} = \chi_G(Z), \quad \text{ and } \quad l.Y = v(t_a), \quad m.Z = w(t_a).
        \]
        Viewing both $Y,Z$ as elements of $\fg(F)$, they have the same characteristic polynomial $a_{\eta} \in H^0(F_K, \fc)$. Their transporter
        $$
            \mathcal{T} = \{g \in G_F : \Ad(g) Y = Z\}
        $$
        is a torsor for the centralizing torus $Z_{G_F}(Y)$. Since $F$ is of cohomological dimension $\le 1$, cohomological vanishing as in \cite[III.2.3, Thm 1' and Remarks]{serregalcoh} implies there exists an element $g \in G(F)$ such that $g.Y = Z$ (cf.\ \cite[Lemme 3.8.1]{CLI} and \S \ref{subsect.loc.triv} below). 

        Since $Y,Z$ are elliptic elements, their split centralizers are the maximal split tori in their Levis: $A_Y^{\spl} = A_L^{\spl}$, $A_Z^{\spl} = A_M^{\spl}$. Since $g$ is rational the conjugation map $\phi = \Inn(g) : A_Y \to A_Z$ identifies the maximal $F$-split sub-tori, so
        \[
            \phi(A_L^{\spl}) = A_M^{\spl}.
        \]
        Taking centralizers implies that $\phi(L) = M$, i.e., that $L$ and $M$ are conjugate by $G(F)$. Since maximal split tori are all conjugate in $M$, we may post-compose with conjugation in $M$ in order to assume that the conjugating element normalizes the fixed maximal split torus $S$ (and thus also $T_F = Z_{G_F}(S)$). Write $g' = hg \in G(F)$ with $h \in M(F)$ for such a modification, and write $u \in W_G(F) \subseteq W_G(\infty)$ for the Weyl image of $g'$. We have then $u L u^{-1} = M$.
        
        The above does not, however, imply that $w$ and $v$ differ by an element of $W_M(\infty)$. We must use the diagonalizations around $\infty$: the element
        $$
            q : = m g l^{-1} \in G(F_\infty)
        $$
        sends $v(t_a)$ to $w(t_a)$ by construction, so lies in $N_G(T)(F_\infty)$. The composition
        $$
            q (g')^{-1} = m g l^{-1} (g')^{-1} = m h^{-1} (g' l^{-1} (g')^{-1})
        $$
        lies in $M(F_\infty)$ since $g' L (g')^{-1} = M$. But both factors normalize $T$, so $q (g')^{-1}$ lies in $N_M(T)(F_\infty)$, with Weyl image $u' \in W_M(\infty)$. In the Weyl group we have
        $$
            u' u v(t_a) = w(t_a)
        $$
        so by regularity we have $u' u v = w$. Since $u' \in W_M(\infty)$ we conclude
        $$
            (w, M) = (u' u v, M) \sim (u v, M) \sim (v, L)
        $$
        which contradicts the assumption.

                                                            \end{proof}
\end{proposition}

We have therefore proven that the twisted loci give a stratification of the Hitchin base, one which will be well-adapted for weighted stability. From the dimension formulae we have density of the elliptic open in each locus, and we deduce the following.
\begin{corollary}\label{cor.ell.strat.tw}
    The loci $\ctAe{w}{M}$, ranging over classes $[(w,M)]$, define a stratification of $\ctA_{G}$ with closure relations given by $\le$ in Definition \ref{def.wMpair.rels}. In particular, the compatible locus is a nonempty open subset of $\ctA_G$.
\end{corollary}

\section{\texorpdfstring{Monodromic invariants.}{Monodromic invariants.}}\label{subsect.Ngo.monod.comp}
For later geometric arguments, we compare our twisted elliptic stratification to Ng\^o's much finer stratification of the Hitchin base by ``monodromic invariants'' (cf.\ \cite[\S 5.6]{ngofl}). Let $\wtilde{a} = (a,t) \in \ctA(\kappa)$ be a point of the Hitchin base over a field. We usually suppress the subscript $\kappa$, but all the following constructions should be understood after pullback to $\kappa$. Associated to $a$ one has the \textit{cameral cover} $X_a \to X_{\kappa}$, defined as the pullback
\begin{equation} \label{eqn.spec.curve.I}
    \begin{tikzcd}[ampersand replacement=\&]
        X_a \ar[r] \ar[d, "\pi_a"] \& \ft_D \ar[d] \\
        X_{\kappa} \ar[r, "h_a"] \& \fc_D,
    \end{tikzcd}
\end{equation}
where $a$ determines a morphism $h_a \in \Mor(X_{\kappa}, \fc_D)$. The data of $t$ gives a canonical point $\wtilde{\infty}$ in the fiber over $\infty$. 
Let $U_a = a^{-1}(\fc_{G,D,\kappa}^{\Greg}) \subset X$ denote the open on which $a$ is regular semisimple and let $V_a =(\pi_a)^{-1}(U_a)$. Write
$$
    r_a : \pi_1(U_a, \infty) \to \pi_1(X, \infty)
$$
for the everywhere-unramified quotient map.

\begin{remark} In the case of a classical group, sections of $X_a \to X$ correspond to orderings of the factorizations of the characteristic polynomial $a$. Contrast this with the theory of spectral curves, in which components correspond to the factors themselves (and which are not necessarily Galois covers of $X$). 
\end{remark}

Following \cite[\S 5.6]{ngofl}, we wish to describe monodromic invariants associated with the cameral cover $X_a \to X$. The \'etale cover $X_a \times_X U_a \to U_a$ is a $W|_{U_a}$-torsor, which defines a class in $H^1(U_a, W|_{U_a})$. Since $t$ determines a point $\wtilde{\infty} \in X_a$ above $\infty$, it corresponds to a representing cocycle
$$
    c_{\wtilde{a}} \in Z^1(\pi_1(U_a, \infty), W_{\infty}),
$$
where the action on $W_\infty$ is through the composition 
$$
    \rho_G \circ r_a : \pi_1(U_a, \infty) \to \pi_1(X, \infty) \to \Out(G_\infty).
$$
Such a cocycle is equivalent data to its graph, which induces a homomorphism
$$
    \pi_1(U_a, \infty) \to W_{\infty} \rtimes \Out(G_\infty).
$$
Our running choice of \textit{pointed} splitting datum for $G$ (cf.\ \S \ref{subsect.notation.running.assumptions}) determines in particular a choice of isomorphism $W_\infty \cong \W$, equivariant for the isomorphism $G_\infty \cong \G$, so we will interchangably consider this homomorphism to land in $\W \rtimes \Out(\G)$.

\begin{definition}\label{def.monodr.inv.Out}
    Define the \emph{monodromic homomorphism} associated to $\wtilde{a}$ to be the map above,
    \begin{align*}
        \mu_{\wtilde{a}} : \pi_1(U_a, \infty) &\to W_\infty \rtimes \Out(G_\infty) \cong \W \rtimes \Out(\G) \\
        \gamma &\mapsto (c_{\wtilde{a}}(\gamma), \rho_G(r_a(\gamma))).
    \end{align*}
    Define the \emph{monodromic subgroup} to be the image of this homomorphism, 
     $$W_{\wtilde{a}} := \im(\mu_{\wtilde{a}}).$$ Define also the \emph{(cameral) inertia subgroup},
     $$
        I_{\wtilde{a}} := \mu_{\wtilde{a}}(\ker r_a).
     $$
\end{definition}

These are both finite groups; indeed, they land within $W_\infty \rtimes \im(\rho_G)$, which is finite because $\pi_1$ is profinite and $\Out(G_\infty)$ (automorphisms of the based root datum for $G_\infty$) is discrete.

\begin{remark}
    It is immediate from the definitions that
    $$
        I_{\wtilde{a}} \normsg W_{\wtilde{a}}, \qquad I_{\wtilde{a}} \subseteq W_{\wtilde{a}} \cap \W, \quad \text{ and } \quad W_{\wtilde{a}} \twoheadrightarrow \im \rho_G.
    $$
    We will sometimes abuse notation and write $a$ as shorthand for $\wtilde{a} = (a,t)$ when no confusion can arise, and so likewise $W_a = W_{\wtilde{a}}$, $I_a = I_{\wtilde{a}}$. Of course, the rigidifying datum $t$ is necessary to obtain well-defined subgroups rather than conjugacy classes.
\end{remark}

\begin{remark}\label{rmk.monodr.inv.spl}
Given $G/X$ which is split by a (geometrically connected) $\Gamma$-Galois cover $\pi: X' \to X$, we can give a presentation of the monodromy invariants relative to $\pi$. The monodromic homomorphism $\mu_{\wtilde{a}}$ factors along the natural map $\W \rtimes \Gamma \to \W \rtimes \Out(\G)$ which is the identity on $\W$. We will sometimes therefore pick a splitting cover and work with the images of $W_{\wtilde{a}}$ and $I_{\wtilde{a}}$ in $\W \rtimes \Gamma$. Of course, when $\pi$ is a minimal pointed splitting cover this recovers the definition above. 

This presentation admits a direct geometric interpretation. Consider the further pullback
\begin{equation} \label{eqn.spec.curve.II}
    \begin{tikzcd}[ampersand replacement=\&]
        X_a' \ar[r] \ar[d, "\pi_a'"] \& X_a \ar[d] \\
        X'_{\kappa} \ar[r] \& X_{\kappa}
    \end{tikzcd}.
\end{equation}
A choice of basepoint $\infty'$ over $\infty$ determines an isomorphism $G_{X'} \cong G_0 \times X'$ (as well as compatible isomorphisms $T_{X'} \cong T_0 \times X'$ and $W_{X'} \cong \W \times X'$) by restriction to $\infty'$. In this case $X_a'$ is equipped with a point $\wtilde{\infty}' = (\infty', \wtilde{\infty})$ over $\infty$. Now $X_a' \to X_\kappa'$ has an action of $\W \cong W_{\infty}$, the split Weyl group, and $X_\kappa' \to X_\kappa$ and $X'_a \to X_a$ have compatible actions of $\Gamma$. So the composite $X_a' \to X_\kappa$ admits an action of the constant group $\W \rtimes \Gamma$, depending on the pointed splitting cover $(X', \infty')$. Define $U_a' := \pi^{-1}(U_a) \subseteq X'$ and $V_a' = (\pi_a')^{-1}(U_a')$.

Write 
$$
    j_\Gamma : \W \rtimes \Gamma \to \W \rtimes \Out(\G)
$$
for the natural map, and let $C_{\wtilde{\infty}'}'$ be the connected component of the normalization of $X_a'$ containing $\wtilde{\infty}'$. By Galois theory (cf.\ {{\cite[1.3.6, 5.4.1, 5.4.2]{ngofl}}}), we can see that the monodromic subgroup is given by $$W_{\wtilde{a}} = j_\Gamma(\on{Stab}_{\W \rtimes \Gamma}(C_{\wtilde{\infty}'})).$$ 
If $\pi$ is taken to be a minimal pointed splitting cover for $G$, then $j_{\Gamma}$ is injective and we can identify this stabilizer with $W_{\wtilde{a}}$ directly.

The inertia subgroup is the subgroup generated by the stabilizers of geometric points of $C_{\wtilde{\infty}'}$. 
\end{remark}

\begin{remark}\label{rmk.monodr.inv.reg.centr}
    The regular centralizer $J_a \to X$ (cf.\ \S \ref{sect.reg.cent.Pic}, \cite[\S 2]{ngofl}) gives an alternate perspective on these invariants, which it is convenient to record now. The restriction $J_a|_{U_a}$ is a torus, and the datum $t$ determines an isomorphism $J_a|_\infty \cong T|_\infty$. Under this isomorphism, $\mu_{\wtilde{a}}$ is precisely the monodromy of the cocharacter local system $\ul{X}_\star(J_a|_{U_a})$ through the containment
    $$
        W_\infty \rtimes \Out(G_\infty) \subseteq \Aut(X_\star(T_\infty)).
    $$
\end{remark}

\begin{example}
    Let $a \in \cA_{T}$ inside of $G = U_{X'/X}(3)$. The $\Gamma=\Z/2$ cover $X' \to X$ splits the group, so we can identify the monodromic data inside $\W \rtimes \Gamma$. The generic fiber of $X_a$ is $\Spec F \sqcup \Spec F \sqcup \Spec F' \sqcup \Spec F'$, where $F'/F$ is the extension of function fields corresponding to $X'/X$. If we choose a basepoint $\wtilde{\infty}$ in one of the two rational components, then the pair $(a,t)$ will be compatible, and the other rational basepoint consists of $(1 \; 3)\wtilde{\infty}$. 

    Pulling back along $X'/X$, the cover splits into 6 irreducible components with function fields $F'$. Galois fixes the two rational components (interchanging their fibers over $\infty$) and exchanges the other four in pairs. On the other hand, the cover $X_a' \to X'$ is on the regular semisimple locus a trivial $S_3$-torsor, so $S_3$ acts simply transitively on the set of components. 
    
    Choosing an $\infty' \in X'$ lying over $\infty$ so that $(\infty', \wtilde{\infty})$ gives a based $S_3 \rtimes \Gamma$-cover $X'_a \to X$, we can see that 
    \begin{itemize}
        \item If $\wtilde{\infty}$ is chosen to be a basepoint on an $F$-component, then anything in the Galois group preserves that irreducible component, but no nontrivial element of $S_3$ does, so $W_{\wtilde{a}} = 1 \rtimes \Gamma$.
        
        \item If $\wtilde{\infty}$ is chosen to lie on one of the $F'$-components, then a non-$\Gamma$-stable, order-2 subgroup $W_{\wtilde{a}} = \langle (w, \sigma)\rangle \cong \Z/2$, a conjugate of $1 \rtimes \Gamma$, normalizes that component. In fact $w = (1 \; 2 \; 3)$ or $(1 \; 3 \; 2)$.
    \end{itemize}
    In both cases, $I_{\wtilde{a}} = 1$. 
\end{example}

\begin{remark}\label{rmk.Ngo.counterex.untwisted}
    The above example demonstrates the necessity of deleting all the \emph{twisted} images of $\ctA_M$ when defining $\ctA_{\Gell}$. Indeed, if we choose $\wtilde{\infty}$ above, equivalent data to $t$, non-compatibility to lie on one of the $F'$-components of the cameral cover, then $\wtilde{a} = (a,\wtilde{\infty})$ is evidently not in the image of the untwisted closed embedding from the only Levi subgroup $\ctA_T \hookrightarrow \ctA_G$. However, we argue that such point $(a, \wtilde{\infty})$ should not be considered $G$-elliptic -- even ignoring the center it is not anisotropic in the sense of Ng\^{o} \cite[5.4.7]{ngofl}. Computing the action of $\W \rtimes \Gamma$ on $\T$, we see that $\T^{W_{\wtilde{a}}}$ contains a nontrivial rank 1 torus, though not one that descends to $F$. Alternatively, this $W_{\wtilde{a}}$ is conjugate to $1 \rtimes \Gamma$, so its connected fixed locus inside $\T$ is isomorphic to $(\T^{1 \rtimes \Gamma})^0 = \{\on{diag}(z, 1, z^{-1})\}$. This is not surprising: anisotropicity corresponds to quasicompactness of the Picard stack for the Hitchin fiber mod center (through finiteness of $\pi_0(\cP)$, see \S \ref{sect.reg.cent.Pic}), which is \emph{independent of our choice of point $\wtilde{\infty}$} because the Hitchin fibration is defined over $\ctA$ by pullback from $\cA$. 

    This demonstrates a minor ambiguity in the statement of \cite[Proposition 6.3.5]{ngofl}: even for groups with anisotropic center, the anisotropic locus of Ng\^{o} ($G$-elliptic, in our terminology) is the complement of the images of all \emph{twisted} embeddings of Levis as in Definition \ref{def.wM.loci}; it is not enough to consider only the untwisted embeddings $\ctA_M \to \ctA_G$. Examining the proof of \textit{loc.\ cit.}, we see that for a non-anisotropic point, the nontrivial torus $\mathbb{S}$ contained in $\T^{\W_{\wtilde{a}}}$ and the Levi $\mathbb{M}$ defined there need not be stable under the original pinned $\Gamma$-action except up to conjugacy. That is, we can find a preimage $a_M$ for $a$, but not necessarily compatible with our given data $\wtilde{\infty}$. 

    Of course, since we have only deleted a finite number of closed subvarieties isomorphic to $\ctA_M$, the codimension estimate deduced from this proposition, which is the only way it is used in \textit{loc.\ cit.}, remains intact.
\end{remark}

The monodromic invariant subgroup is sufficient to detect Levi-reductions of polynomials, in the following sense. Note that $\Out(\G)$ need not preserve $\M$, but $\im(\rho_G)$ does, so we have a natural map 
$$
    \W_M \rtimes \im(\rho_G) \to \W \rtimes \im(\rho_G).
$$

\begin{proposition}[{{cf.\ \cite[Prop. 8.2.1]{CLII}}}]\label{prop.ngo.comparison.wa}
    For a semistandard $X$-Levi $M$, an element $\wtilde{a} = (a,t)$ lies in the image of $\ctAt{e}{M}$ if and only if $W_{\wtilde{a}} \subseteq \W_M \rtimes \im(\rho_G)$, where $\W_M \subseteq \W$ corresponds to the Levi $M_\infty \subset G_{\infty}$ (given by twisting $\M \subset \G$).

    Moreover, if $\wtilde{a}_M = (a_M, t)$ is a preimage in $\ctA_M$ for $\wtilde{a} = (a,t)$, then $\mu_{\wtilde{a}_M} = \mu_{\wtilde{a}}$ as representations of $\pi_1(U_a, \infty)$ on $X_\star(T_\infty)$. 
\end{proposition}
\begin{proof}
    For the first equivalence, it suffices to check the desired containment relative to a minimal pointed splitting $\Gamma$-cover $\pi : X' \to X$, as in Remark \ref{rmk.monodr.inv.spl}. Assume first that $W_{\wtilde{a}} \subseteq \W_M \rtimes \Gamma$. We will then construct a preimage $a_M$ for $a$.

    Combining (\ref{eqn.spec.curve.I}) and (\ref{eqn.spec.curve.II}) gives a commutative diagram of $X_{\kappa}$-schemes with Cartesian squares
    \[
    \begin{tikzcd}[ampersand replacement=\&]
        \Spec(\kappa) \ar[r, "\wtilde{\infty}'"] \ar[dr, "\infty'"] \ar[ddr, "\infty" below left] \& V_a' \ar[d] \ar[r, hookrightarrow, "i_a'"] \& {\ft'}^{\reg}_{D, \kappa} \ar[d, "\chi'"] \\
        \& U_a' \ar[r] \ar[d] \& {\fc'}^{\reg}_{G', D',\kappa} \ar[d] \\
        \& U_a \ar[r, hookrightarrow, "h_a"] \& \fc_{G,D, \kappa}^{\reg}
    \end{tikzcd}.
    \]
    
    Let $C := (V_a')^0$ be the connected component of $V_a'$ containing $\wtilde{\infty}'$. Since $\pi$ is a minimal splitting cover, we can identify its stabilizer with $W_{\wtilde{a}}$ as in Remark \ref{rmk.monodr.inv.spl} above. By construction, $C \to U_a$ is Galois \'etale of group $W_{\wtilde{a}}$, so the composition of $i_a'$ with the two-step quotient 
    \[
        \begin{tikzcd}[ampersand replacement=\&]
            C \ar[r, "i_a'"] \& {\ft'}^{\reg}_{D, \kappa} \ar[r] \& {\ft'}^{\reg}_{D', \kappa}\gitq \W_M = {\fc'}^{\reg}_{M', D', \kappa}  \ar[r, "\pi"] \& {\fc'}^{\reg}_{M', D', \kappa}\gitq \Gamma = {\fc}^{\reg}_{M,D, \kappa}
        \end{tikzcd}
    \]
    descends to a lift $h_{a_M}: C \gitq W_{\wtilde{a}} = U_a \to {\fc}^{\reg}_{M,D, \kappa}$. Note that the the Chevalley spaces for $G/X$ are built by descent from those for $G'/X'$, so the corresponding Chevalley maps are compatible with pullback.

    The valuative criterion of properness along the finite map $\fc_{M,D} \to \fc_{G,D}$ now guarantees the extension of $h_{a_M}$ to a lift on all of $X_{\kappa}$. This lift $h_{a_M}$ is by construction compatible with the datum $\wtilde{\infty}'$, so by considering the two-step quotient in the other order $V_a' \to V_a \to U_a$ we see that $h_{a_M}$ is compatible with $\wtilde{\infty}$ and so $(a_M, t) \in \ctA_M$ is a preimage of $(a,t)$. This completes the first implication.

    Now assume that $\wtilde{a} = (a,t)$ has an untwisted preimage $\wtilde{a}_M = (a_M, t)$. We work over the $G$-regular semisimple locus $U_a \subseteq U_{a_M}$. The restriction $\pi_1(U_a, \infty) \to \pi_1(U_{a_M}, \infty)$ is surjective and maps the kernel of $r_a$ onto the kernel of $r_{a_M}$, so this restriction changes neither the monodromy image nor the inertia image.
    
    The cameral map $V_{a_M}|_{U_a} \to V_a|_{U_a}$ is equivariant for $W_M|_{U_a} \to W_G|_{U_a}$ and by assumption maps the distinguished point over $\infty$ to the distinguished point. So not only the cohomology classes but the representing cocycles are compatible: $c_{\wtilde{a}} = c_{\wtilde{a}_M}$ under the inclusion $W_{M,\infty} \hookrightarrow W_{G, \infty}$. The actions of $\im(\rho_G)$ on $W_{M,\infty} \cong \W_M$ and $W_{G,\infty} \cong \W_G$ are compatible, so we deduce
    $$
        \mu_{\wtilde{a}_M} = \mu_{\wtilde{a}}
    $$
    as representations of $X_\star(T_\infty)$. In particular, the monodromy subgroups and inertia subgroups for $\wtilde{a}_M$ and $\wtilde{a}$ coincide, and in particular $W_{\wtilde{a}} \subseteq \W_M \rtimes \im(\rho_G)$, giving the converse direction of the first claim.

                    \end{proof}

The twisted loci too can be characterized in terms of monodromic invariants, though we won't use this.

\begin{corollary}\label{prop.ngo.comparison.tw}
    An element $\wtilde{a} = (a, t) \in \ctA_G$ lies in $\ctAt{w}{M}$ if and only if $W_{\wtilde{a}} \subseteq w^{-1}(W_{M, \infty} \rtimes \im(\rho_G))w$.
\begin{proof}
    Change of base-point changes the representing cocycle $c_{\wtilde{a}}$ by a coboundary, which corresponds to conjugating $\mu_{\wtilde{a}}$ by the corresponding element $w \in W_{G, \infty}$. That is,
    $$
        \mu_{(a, w(t))} = \Ad((w,1)) \circ \mu_{(a,t)}
    $$
    and so $W_{(a,w(t))} = w W_{(a,t)} w^{-1}$. The previous proposition applied to $(a,w(t))$ gives the result.

                \end{proof}
\end{corollary}

Ng\^{o} has shown that the monodromic invariants are semicontinuous on $\ctA$.

\begin{proposition}[{{\cite[Lemme 5.4.5]{ngofl}}}]\label{prop.monodr.inv.semicts}
    Let $\wtilde{a} \in \bAe{M}$ and $\wtilde{a}' \in \bAe{L}$ be a specialization of $a$, so that $L \subseteq M$. Then $W_{\wtilde{a}'} \subseteq W_{\wtilde{a}}$.
\end{proposition}

We have a simple characterization of the $W_{\wtilde{a}}$-invariants, and a crucial finiteness result we will need in \S \ref{sect.reg.cent.Pic}. We have a short exact sequence of $X$-tori 
$$
    \begin{tikzcd}[ampersand replacement=\&]
        1 \ar[r] \& T \cap D(M) \ar[r] \& T \ar[r] \& M^{\ab} \ar[r] \& 1
    \end{tikzcd}.
$$
Taking $\ul{X}_\star(-)(\infty)$ gives a short exact sequence of $(W_{M,\infty} \rtimes \im(\rho_G))$-modules, and so of $W_{\wtilde{a}}$-modules. We denote the result
$$
    \begin{tikzcd}[ampersand replacement=\&]
        1 \ar[r] \& \bX_{\star, D(M)} \ar[r] \& \bX_\star \ar[r] \& \bX_{\star, M^{\ab}} \ar[r] \& 1
    \end{tikzcd}.
$$
Tensoring with $\R$ gives a decomposition of the cocharacter space,
$$
    \begin{tikzcd}[ampersand replacement=\&]
        1 \ar[r] \& \ga_{D(M)} \ar[r] \& \ga_T \ar[r] \& \ga_{M^{\ab}} \ar[r] \& 1
    \end{tikzcd}.
$$

\begin{lemma}\label{lem.monodr.inv.spl.torus}
For every Levi $M$ and $a \in \bAe{M}$, 
\begin{enumerate}[label=(\alph*)]
    \item $$
    \ga_{D(M)}^{W_{\wtilde{a}}} = 0,
$$
and moreover integrally
$$
    (\bX_{\star,D(M)})^{W_{\wtilde{a}}} = 0.
$$
\item The coinvariants $(\bX_{\star,D(M)})_{W_{\wtilde{a}}}$ is finite.
\item We have an equality
$$
    ((T_\infty^{W_{\wtilde{a}}})_{\red})^0 = (A_M^\spl)_\infty.
$$
Indeed, in good characteristic $T_{\infty}^{W_{\wtilde{a}}}$ is smooth and
$$
    (T_\infty^{W_{\wtilde{a}}})^0 = (A_M^\spl)_\infty.
$$
\end{enumerate}
Moreover, for every $b \in \ctA_G$,
\begin{enumerate}[label=(\alph*), start=4]
    \item $\a_G \subseteq \ga_T^{W_b}$, with equality if and only if $b \in \ctA_{\Gell} = \bAe{G}$.
\end{enumerate}

\begin{proof}
    For (a), since the $W_{\wtilde{a}}$-invariants is torsion-free, its vanishing can be checked after tensoring with $\R$ ($\Q$ is enough), so the statement about $\ga_{D(M)}$ suffices. Now let $v \neq 0 \in \ga_{D(M)}^{W_{\wtilde{a}}}$. Up to conjugating $v$ and $W_{\wtilde{a}}$ by an element of $W_{M, \infty}$, we may assume that $v$ lies in the closed dominant geometric chamber with respect to a fixed rational semistandard Borel of $M$. This preserves $M$-ellipticity of the corresponding twisted $\wtilde{a}$, since the elliptic loci are defined by deleting all \textit{twisted} loci for proper Levis.
    
    In particular, this chamber is stable under $\im(\rho_G)$. Note that $\ga_{D(M)}$ is the usual geometric reflection representation of $W_{M, \infty}$, and is stratified by the geometric Weyl fan. Since $v \neq 0$, it lives in a face corresponding to a proper geometric parabolic subgroup $Q_\infty \subsetneq M_\infty$ with Levi $L_\infty \subsetneq M_\infty$. 

    We claim now that $v$ and the face it lies in are $\im(\rho_G)$-stable. Since $W_{\wtilde{a}}$ surjects onto $\im(\rho_G)$, for each $\gamma \in \im(\rho_G)$ we have some $(x_\gamma, \gamma) \in W_{\wtilde{a}}$, which fixes $v$, mapping to it. Then
    $$
        x_\gamma \gamma(v) = v
    $$
    so $\gamma(v), v$ are $W_M(\infty)$-conjugate but both in the closed dominant chamber, hence \textit{equal}. The face it lies in corresponds to the set of geometric simple roots satisfying $\langle \alpha, v\rangle = 0$, hence is also $\im(\rho_G)$-stable.
    
    The above implies that $Q_\infty$ and $L_\infty$ descend to a rational parabolic and Levi. But this implies $W_a$ lies in the stabilizer of the face, which is $W_{L, \infty} \rtimes \im(\rho_G)$. By Proposition \ref{prop.ngo.comparison.wa}, this contradicts $M$-ellipticity of $\wtilde{a}$.

        For (b), we know $(\bX_{\star, D(M)})_{W_{\wtilde{a}}}$ is finite type, so it is finite if and only if it vanishes when tensored to $\R$. Averaging over $W_{\wtilde{a}}$ induces an isomorphism
        $$
            (\bX_{\star, D(M)})_{W_{\wtilde{a}}} \otimes_{\Z} \R \xrightarrow{\sim} (\bX_{\star, D(M)})^{W_{\wtilde{a}}} \otimes_{\Z} \R,
        $$
        so we are done by (a).

        For (c), note that taking $W_{\wtilde{a}}$-invariants of the (real) cocharacter spaces is exact. On $M^{\ab}$, $W_{\wtilde{a}}$ acts through $\im(\rho_G)$ (and surjects onto it), so
        $$
            \ga_{M^{\ab}}^{W_{\wtilde{a}}} = \ga_{M^{\ab}}^{\im(\rho_G)}.
        $$
        So in the short exact sequence
        $$
        \begin{tikzcd}[ampersand replacement=\&]
        1 \ar[r] \& \ga_{D(M)}^{W_{\wtilde{a}}} \ar[r] \& \ga_T^{W_{\wtilde{a}}} \ar[r] \& \ga_{M^{\ab}}^{W_{\wtilde{a}}} \ar[r] \& 1
        \end{tikzcd}
        $$
        the left side vanishes and the right side has dimension equal to the split rank of $M^{\ab}$, which is equal to $\dim A_M^{\spl}$. So $\dim({(T_\infty^{W_{\wtilde{a}}})}_{\red}) = \dim A_M^{\spl}$. By definition, $A_M^{\spl}$ is reduced and fixed pointwise by ${W_{\wtilde{a}}}$, so we have
        $$
            (A_M^{\spl})_\infty \subseteq {(T_{\infty}^{W_{\wtilde{a}}})}_{\red}
        $$
        and indeed within its identity component. But we have shown they have the same dimension, so equality holds. In good characteristic (Definition \ref{def.good.characteristic}), $p \nmid |W_G| |\im(\rho_G)|$, so $p$ does not divide the order of the monodromy subgroup and the last claim is clear.

        Then (d) is essentially by the same argument as above. We have
        $$
            \ga_T^{W_b} = \ga_{D(G)}^{W_b} \oplus \ga_{G}^{W_b}
        $$
        and $W_b$ acts through and surjects onto $\im(\rho_G)$, so the right factor is $\a_G$. That gives the inclusion. If moreover $b$ is $G$-elliptic, then it is compatible, so the left factor vanishes. If it is not $G$-elliptic, then it lies in some twisted locus for a proper $X$-Levi $L$, and so by Proposition \ref{prop.ngo.comparison.tw}
        $$
            W_b \subseteq w^{-1}(W_{L, \infty} \rtimes \im(\rho_G)) w.
        $$
        This latter group fixes $w^{-1} \a_L$ pointwise, which strictly contains $\a_G$.
\end{proof}
\end{lemma}

\section{Endoscopic loci and ellipticity.}\label{subsect.hitchin.base.endoscopic}

We here recall various standard results relating ellipticity and the Hitchin base for endoscopic groups. Various notions of transfer in the Langlands program begin with transfers at the level of coarse quotients. In our setting, the fundamental such transfer is the functoriality of the spaces of invariants, and thus the Hitchin base, in endoscopic groups. We fix a pointed twisting datum $\rho_G$ for $G$, and all endoscopic data are pointed endoscopic data relative to this choice. Until the final remark, everything in this subsection occurs over algebraically closed $k$.

\begin{definition}\label{def.end.loci}
    Let $(s, \rho_s)$ a (pointed) endoscopic datum for $G$ with endoscopic group $H = G'_s$. The natural map of Chevalley spaces from Proposition \ref{prop.chevalley.functor.endo.gps} induces on global sections a finite map, depending on $\rho_s$,
    $$
        \cA_{G'_s} \to \cA_G.
    $$
    Define the rigidification $\ctA_{G'_s}$ as for Levi subgroups, enforcing \textit{$G$-regular} (rather than just $G'_s$-regular) semisimplicity at $\infty$:
    $$
        \ctA_{G'_s} = \{(a_{G'_s}, t) : t \in \ft_\infty^{G\rs}, \; \chi_{G'_s}(t) = a_{G'_s}(\infty)\}.
    $$
\end{definition}

\begin{proposition}[{{\cite[\S 1.9, Proposition 6.3.2]{ngofl} or \cite[{\S 5.1-5.3}]{geostabpadic}}}]\label{prop.hitchin.base.functorial.endoscopic.groups}
The rigidified map 
$$
    i_{\rho_s} = i_{(s, \rho_s)}: \ctA_{G'_s} \to \ctA_G
$$
is a closed embedding.
\end{proposition}

We can give a monodromy-theoretic characterization of the images of these maps. Recall also that given $s \in \widehat{T}$, one has the coendoscopic connected-\'etale sequence
\[
\begin{tikzcd}[ampersand replacement=\&]
        1 \ar[r] \& \widehat{G}'_s \ar[r] \& Z_{\widehat{G} \rtimes \Out(\widehat{G})}(s) \ar[r] \& \pi_0^G(s) \ar[r] \& 1
\end{tikzcd}.
\]
There is a corresponding short exact sequence of Weyl groups, noting that $\W \cong \widehat{\W}$ and $\Out(\G) \cong \Out(\widehat{G})$:
$$
\begin{tikzcd}
    1 \ar[r] & \W_{\widehat{G}_s} \ar[r] & Z_{\W \rtimes \Out(\G)}(s) \ar[r, "p_s"] & \pi_0^G(s) \ar[r] & 1.
\end{tikzcd}
$$
Denote $W_s$ for the middle term in this sequence and $I_s$ the left term. By \cite[Lemme 10.1]{ngoendo}, this short exact sequence is canonically split by those elements that leave fixed the set of positive roots corresponding to the Borel, so we have a canonical isomorphism $W_s \cong I_s \rtimes \pi_0^G(s)$. 

\begin{proposition}[{\cite[\S 6.3]{ngofl}}]\label{prop.endoscopic.locus.decomposition} For any $s \in \widehat{T}$ the locus 
$$
    \ctA_{s} := \{\wtilde{a} : W_{\wtilde{a}} \subseteq W_s, I_{\wtilde{a}} \subseteq I_s \} \subset \ctA_G
$$
is exactly the finite disjoint union of the images of the reduced endoscopic closed embeddings
$$
    i_{\rho_s} : \ctA_{G_s'} \to \ctA_G
$$
of Proposition \ref{prop.hitchin.base.functorial.endoscopic.groups} over all pointed endoscopic data $\rho_s : \pi_1(X, \infty) \to \pi_0^G(s)$ for that particular $s$. 

More precisely, $\wtilde{a} \in i_{\rho_s}(\ctA_{G'_s})$ if and only if $W_{\wtilde{a}} \subseteq W_s$ and
$$
    p_s \circ \mu_{\wtilde{a}} = \rho_{s} \circ r_a : \pi_1(U_a, \infty) \longrightarrow \pi_0^G(s).
$$
In particular, if such a lift $\wtilde{a}_{G'_s}$ exists, then 
$$
    W_{\wtilde{a}_{G'_s}} = W_{\wtilde{a}}, \qquad I_{\wtilde{a}_{G'_s}} = I_{\wtilde{a}}.
$$
\begin{proof}[Proof Outline] First, the set of $\rho_s$ for a given $s$ is finite, despite the possible infinitude of $\pi_0^G(s)$: any such $\rho_s$ lands in the finite group 
$$
    \frac{W_s \cap (\W_G \rtimes \im(\rho_G))}{I_s},
$$
and there are finitely many possible homomorphisms $\rho_s$ from the (topologically finitely generated) geometric fundamental group into this finite group. 

Given any point $\wtilde{a} = (a,t) \in \ctA_s$ gives rise to a natural map 
$$
    \pi_1(U_a, \infty) \to W_{\wtilde{a}}/I_{\wtilde{a}} \to W_s/I_s = \pi_0^G(s)
$$
which descends by the inertia condition to a unique $\rho_{s} : \pi_1(X, \infty) \to \pi_0^G(s)$, i.e.\ a pointed endoscopic datum for $s$. So $\ctA_s$ splits as a disjoint union over loci $\ctA_{(s,\rho_s)}$ corresponding to a particular endoscopic datum.

One can then lift $\wtilde{a}$ and its monodromy homomorphism to a unique point and homomorphism for the Hitchin base $\ctA_{G'_s}$, by an argument analogous to Proposition \ref{prop.ngo.comparison.wa}, cf.\ \textit{loc.\ cit.} The only additional subtlety not visible in the Levi case is that the endoscopic monodromy cocycle $c_{G'_s}$ must be defined to differ from $c_{\wtilde{a}}$ by the Weyl twist in the embedding of Proposition \ref{prop.Weyl.functor.endo.gps},
$$
    c_{G'_s}(\gamma) := c_{\wtilde{a}}(\gamma) u_{\rho_s}(r_a(\gamma))^{-1}
$$
where $u_{\rho_s} : \pi_1(X, \infty) \to \W_G$ is the cocycle corresponding to that embedding. This $c_{G'_s}$ defines a section of the endoscopic Chevalley space over $U_a$, and corresponds to a monodromy homomorphism $\mu_{G'_s}$ restricting $\mu_{\wtilde{a}}$, and equivariant for the action on $X_\star(T_\infty)$. This gives an endoscopic lift of $\wtilde{a}$ over $U_a$, which extends to all of $X$ by finiteness of the map $\fc_{G'_s, D} \to \fc_{G,D}$ and normality of $X$.
\end{proof}
\end{proposition}

\begin{proposition}[{{cf.\ \cite[10.6.6, 10.6.10]{CLII}}}]\label{prop.ell.end.data}
    Given a pointed endoscopic datum $(s, \rho_{s})$ with endoscopic group $H = G_s'$, the following are equivalent.
    \begin{enumerate}[label=(\alph*)]
        \item The closed embedding $i_{\rho_s} : \ctA_{G_s'} \to \ctA_G$ has image which intersects the $G$-elliptic locus,
        \item The canonical surjection $\a_{G'_s} \to \a_G$ is an isomorphism.
        \item The endoscopic group satisfies 
        $$
            A_{G'_s}^\spl = A_{G}^\spl.
        $$
        \item $(Z_{\widehat{G}}^{\pi_1(X, \infty)})^0 = (Z_{\widehat{G}_s'}^{\pi_1(X, \infty)})^0$.
    \end{enumerate}
    Here $\pi_1(X, \infty)$ acts on $Z_{\widehat{G}_s'}$ through the endoscopic action, the canonical composition
    $$
        \pi_1(X, \infty) \xrightarrow{\rho_s} \pi_0(s) \to \Out(\widehat{G}'_s),
    $$
    noting that inner automorphisms act trivially on the center.
    
    In the case that these conditions hold, $i_{\rho_s}^{-1}(\ctA_{G\el}) = \ctA_{G_s'\el}$, where the elliptic locus in $\ctA_{G_s'}$ is defined the same way as for $G$. 
    
    \begin{proof}
        The equivalence of (b) and (c) is clear, and (c) and (d) correspond under duality. Assume $\wtilde{a}_{G'_s} \in \ctA_{G'_s}$ is a point with image $\wtilde{a} \in \ctA_G$. Then by Lemma \ref{lem.monodr.inv.spl.torus}.(d) applied both to $G$ and to $G'_s$, and the identification of the monodromy invariants $W_{\wtilde{a}} = W_{\wtilde{a}_{G'_s}}$ in Proposition \ref{prop.endoscopic.locus.decomposition}, we have canonical inclusions
        $$
            \a_G \subseteq \a_{G'_s} \subseteq \ga_T^{W_{\wtilde{a}}} =\ga_T^{W_{\wtilde{a}_{G'_s}}}.
        $$

        If $\wtilde{a}_{G'_s}$ is $G'_s$-elliptic (such a point which is $G$-regular semisimple at $\infty$ exists under our running degree assumptions, by \cite[4.7.1, 4.7.2, 4.7.5]{ngofl}), then the second containment is an equality. Condition (c) is equality in the first containment, and condition (a) at the point $\wtilde{a}$ is the composite equality, so (a) implies (c). On the other hand, (c) implies that \emph{every} $G'_s$-elliptic point has $G$-elliptic image, so (c) implies (a) and the final claim of the Proposition.
    \end{proof}
\end{proposition}

\begin{definition}\label{def.ell.end.datum}
    We say an endoscopic datum $(s, \rho_s)$ is \emph{elliptic} if it satisfies the equivalent conditions of Proposition \ref{prop.ell.end.data}. 
\end{definition}

In order to obtain nice isotypic decompositions on cohomology, we will need to know that endoscopic parameters $s$ can be taken to be torsion. The crucial result is the following.

\begin{proposition}\label{prop.ell.end.torsion}
    Let $(s, \rho_s)$ be an elliptic endoscopic datum for $G$, with endoscopic group $G'_s$. Then the image of $s$ in
    $$
        \widehat{T}/(Z_{\widehat{G}}^{\pi_1(X, \infty)})^0
    $$
    has finite order. In particular, if $G$ has anisotropic center, then $s$ itself has finite order. 
    \begin{proof}
        By construction, $s$ lies in $Z_{\widehat{G'_s}}^{\pi_1(X, \infty)}$ relative to the endoscopic action. Characterization (d) of ellipticity identifies the identity component of this group with that for $\widehat{G}$. The component group is finite, so $s$ has finite order modulo this identity component. 
    \end{proof}
\end{proposition}

We will ultimately need to perform a recursion through the elliptic endoscopic groups for $G$. This is possible by the following observation.

\begin{proposition}[{{cf.\ \cite[Remarque 10.6.7, Proposition 10.6.10, Remarque 10.7.2]{CLII}}}]\label{prop.ell.end.finite}
    If $G$ has anisotropic center, then there are only finitely many pointed geometric elliptic endoscopic data $(s, \rho_s)$ for $G$, and all its elliptic endoscopic groups also have anisotropic center.
    \begin{proof}
        Proof similar to \textit{loc.\ cit.} For the finiteness, we have shown already that there are finitely many possible $\rho_s$ for a given $s$. For an elliptic datum in a group with anisotropic center, $s$ must lie in 
        $\widehat{T}^{K_{s, \rho_s}}$, where $$K_{s, \rho_s} = p_s^{-1}(\im \rho_s) \subseteq \W_G \rtimes \im(\rho_G)$$
        contains $I_s = \W_{G'_s}$ and we have
        $$
            \left(\widehat{T}^{K_{s, \rho_s}} \right)^0 = \left(Z_{\widehat{G}'_s}^{\pi_1(X, \infty)} \right)^0 = 1.
        $$
        For a given $s$ the groups $K_{s, \rho_s}$ vary, but all are contained in the finite group $\W_G \rtimes \im(\rho_G)$, so the set of possible $s$ lie in a finite union of finite fixed-point groups, and is finite. 
    \end{proof}
\end{proposition}

\begin{remark}\label{rmk.arith.geo.data.correspond}
The above definitions of monodromic invariants are valid arithmetically (over some $k_0$ rather than algebraically closed $k$), if $X, D, \infty, \ldots$ are all defined over $k_0$ and one works with arithmetic endoscopic data in the sense of \S \ref{subsect.arith.geo.endo.data}. The only exception is that the $\Gal(k/k_0)$-stable loci $\ctA_s$ may consist of several components for endoscopic parameters defined over extensions of $k_0$, which are permuted by this group. The geometric pointed parameters are given by restricting the arithmetic ones to $\pi_1(X, \infty) \subseteq \pi_1(X_0, \infty)$. 

If one assumes, as we will in applications, that one is given an arithmetic pointed endoscopic datum $(s, \rho_s)$ such that the endoscopic group $G'_s$ is split in the $k_0$-fiber $\infty$ (as well as $G$), then acts on the root datum at $\infty$ faithfully through the splitting $\pi_0^G(s) \hookrightarrow W_s$. Under the splitting $\pi_1(X_0, \infty) \cong \pi_1(X, \infty) \rtimes \Gal(k/k_0)$ corresponding to $\infty$, the Galois group factor acts on the rational fiber of the torus over $\infty$ trivially. The isomorphism of geometric fibers $$T_{H, \infty} \cong T_{G, \infty}$$ is then defined over $k_0$, and $i_{\rho_s}$ descends to $k_0$.

Since $\im(\rho_s) = \im(\rho_s|_{\pi_1(X, \infty)})$, the condition of ellipticity is stable under extension of scalars from $k_0$ to $k$. 
\end{remark}

\chapter{The Hitchin fibration and parabolic reductions}\label{sect.Hitchin.fib}

As in the previous chapter, we denote $X$ a smooth, proper, geometrically integral curve over field $k$, and let $G \to X$ be a quasi-split reductive group scheme (cf.\ \S \ref{subsect.quasisplit}) with split form $\G$, associated to $\Out(\G)$-torsor $\rho_G$. We assume $k$ is of good characteristic for $G$ in the sense of \S \ref{subsect.char.assumptions}. We fix $\infty \in X(k)$ such that $G_\infty$ is split.

In this chapter, we recall the definition of the moduli of Higgs bundles $\cM$ and the Hitchin fibration $f : \cM \to \cA$ as well as the pullback to the \'etale open $f : \ctM \to \ctA$. We analyze the set of parabolic reductions of $m \in \ctM$ in terms of $f(m)$.

\section{The Hitchin fibration.}

Recall $\Bun_{G} := \Bun_{G,X}$, a smooth Artin stack whose $S$-points consists of the groupoid of $G$-torsors on $X \times S$. 

\begin{definition}\label{def.Higgs.bundles}
    The \emph{moduli of Higgs bundles} $\cM_G$, denoted $\cM$ when no confusion is possible, is defined to be the $k$-stack whose groupoid of $S$-points has objects
    $$
        \{(\cE \in \Bun_G(S), \theta \in H^0(X \times S, \Ad_D(\cE))\}
    $$
    where $\Ad_D(\cE) = \cE \times^G \fg_D$, the $D$-twisted Lie algebra bundle of $\cE$. We call a point $(\cE, \theta)$ a \emph{Higgs bundle}, and $\theta$ its \emph{Higgs field}.
    
    The isomorphisms in $\cM(S)$ are isomorphisms $\cE \cong \cE'$ over $X \times S$ such that under the induced isomorphism $\Ad_D(\cE) \cong \Ad_D(\cE')$, the Higgs fields $\theta$ and $\theta'$ correspond. 
\end{definition}
In other words, $\cM = \Map_{X/k}(X, [\fg_D/G])$, the $k$-stack parameterizing $X$-morphisms $X \to [\fg_D/G]$. The map of forgetting $\theta$ induces a map $\cM \to \Bun_G$, which can be thought of as post-composition
$$
    \cM = \Map(X, [\fg_D/G]) \to \Map(X, [\star/G]) = \Bun_G.
$$

\begin{remark}
Classically, one chooses $D = K_X$, after which Serre duality identifies $\cM$ with $\mathcal{H}^{-1}$ of the cotangent complex of $\Bun_G$. For that reason, over the complex numbers it sits on one side of the ``non-abelian Hodge correspondence.'' For a general $D$, the map $\cM \to \Bun_G$ is often referred to as a ``vector bundle.'' We make this precise, in a derived sense which will play no further role in this work. Associated to $D$ is a perfect complex 
$$
    C_D = Rp_\star(\Ad_D(\cE_{univ})) \in D^b(\Bun_G),
$$
where $p : X \times \Bun_G \to \Bun_G$. Since $X$ is a curve, $C_D$ is perfect of Tor-amplitude $[0,1]$. Taking $$V_D :=\ul{\Spec}_{\Bun_G}(\Sym(C^\vee_D))$$ defines a derived linear stack on $\Bun_G$, whose fiber over $\cE$ is the derived vector space $R \Gamma(X, \Ad_D(\cE))$, which has locally constant Euler characteristic. The moduli of Higgs bundles is then the underlying classical stack $\pi_0 V_D$. Of course, this now has only semicontinuous fiber dimension.
\end{remark}

The moduli of Higgs bundles is also fibered over the Hitchin base $\cA = \cA_G = \Map_{X/k}(X, \fc_{D})$, as follows. 

\begin{definition}[The Hitchin fibration]\label{def.Hitchin.fib}
    The $D$-twisted Chevalley map $\chi : \fg_D \to \fc_D$ defined in \S \ref{subsect.twist.chevalley} induces a map $\chi : \Ad_D(\cE) \to \fc_D$ by descent: over each trivializing cover there is a map $\fg_D \cong \Ad_D(\cE) \to \fc_D$, and two such differ by $G$-conjugation, so descend to a well-defined global map $\Ad_D(\cE) \to \fc_D$. There is an induced morphism on global sections, which induces the \emph{Hitchin morphism} for $G$
    \begin{align*}
        f : \cM &\to \cA \\
            (\cE, \theta) &\mapsto \chi(\theta).
    \end{align*}
    In other words, the Chevalley map is $G$-invariant, so descends to a map $[\fg_D/G] \to \fc_D$. Pushforward along this morphism defines the map 
    $$
        f : \cM = \Map(X, [\fg_D/G]) \to \Map(X, \fc_D) = \cA.
    $$
    When $\deg(D) > 2g-2$, a deformation argument due to Biswas and Ramanan deduces from the simple case of Riemann--Roch that the generically regular-semisimple Higgs space $\cM \times_{\cA} \cA^{\heartsuit}$ is smooth over $k$ (\cite[Th\'eor\`eme 4.14.1]{ngofl}). 
    
    We define the \emph{rigidified Hitchin fibration}, also denoted $f$, to be the pullback $\ctM := \cM \times_{\cA} \ctA \to \ctA$, where $\ctA$ is the rigidified Hitchin base of Definition \ref{def.Hitchin.base}. That is, $\ctM$ is the moduli space of triples $(\cE, \theta, t)$ such that $f(\cE, \theta) = a$ agrees with $\chi^T_G(t)$ at infinity, and the rigidified Hitchin map is 
    $$
        (\cE, \theta, t) \mapsto (\chi(\theta), t).
    $$
\end{definition}

\begin{remark}
The utility of the Hitchin fibration in arithmetic can be understood as follows. Restricting to the generic point and assuming the bundle $\cE$ trivial there (which it is up to finite extension of constants), a Higgs field is a conjugacy class of elements of the Lie algebra over the function field $F$ whose matrix entries have zeroes and poles bounded by $D$ over the curve. So this moduli space provides a geometrically well-behaved global setting within which to probe $\fg(F)$ down to a finite depth. As we allow $D$ to range over all sufficiently positive divisors, we can access all of $\fg(F)$. The Hitchin fibration, then, is a globalized version of the map $\fg(F) \to \fc(F)$, and so provides information about the structure of global semisimple $F$-conjugacy classes in the Lie algebra.
\end{remark}

\section{Parabolic reductions of Higgs bundles.}

We recall that reductions of a $G$-torsor $\cE \to X$ to a parabolic subgroup are in bijection with sections of a smooth projective $X$-variety.
\begin{proposition}
     Given $\cE \to X$ a $G$-torsor and $P \subseteq G$ a parabolic subgroup, the relative quotient sheaf $\cE/P \to X$ is representable by a smooth projective $X$-variety, and represents the functor of $P$-reductions of $\cE$. 
    \begin{proof}
        By twisting, the functor of $P$-reductions of $\cE$ is isomorphic to that of parabolic subgroups of $G^/{\cE} := \cE \times^{G, ad} G$ of the same type as $P$ (\cite[Exp. XXVI, 3.20]{SGAIII}). So the proposition follows from the fact that $G^{\cE}/P^{\cE} \cong \ul{\on{Par}}_{t(P)}(G^{\cE})$ is representable by a smooth projective $X$-scheme \cite[Exp. XXVI, 3.20]{SGAIII} (whose connected components define the type).
    \end{proof}
\end{proposition}

\begin{definition}\label{def.P.red}
    Given $m = (\cE, \theta) \in \cM_G(S)$ and a semistandard parabolic $P$, we define a \textit{$P$-reduction} of $m$ to be a pair $(\cE_P, \theta_P)$ consisting of 
    \begin{itemize} 
    \item a $P$-torsor $\cE_P$ on $X_S$,
    \item $\theta_P \in H^0(X_S, \on{Ad}_D(\cE_P))$, and
    \item an isomorphism $\phi : \cE_P \times^P G \cong \cE$,
    \end{itemize}
        such that under the corresponding injection $\on{Ad}(\phi) : \on{Ad}_D(\cE_P) \to \on{Ad}_D(\cE)$, we have $\theta_P \mapsto \theta$.
\end{definition}

Since $\theta$ is generically regular semisimple, there are at most finitely many parabolic reductions up to the natural notion of equivalence. The rigidifying data $t$ allows us to define a relative position for each parabolic reduction, and therefore to organize corresponding stability conditions on each parabolic reduction.

\begin{definition}\label{def.rel.pos}
    Given $m = (\cE, \theta, t) \in \ctM_G(S)$ and a semistandard parabolic $P$, we say that a $P$-reduction $m_P = (\cE_P, \theta_P)$ of $(\cE, \theta)$ is in \textit{relative position }$w \in W_G(\infty_S)$ if it satisfies 
    $$\chi_P(\theta_P)(\infty) = \chi_P(w(t)).$$
    We will sometimes call such a reduction a \textit{$(w, P)$-reduction. } 
\end{definition}

\begin{remark}
    Note that the pair $(w, P)$ hold the same ambiguity as the pairs $[(w, M)]$ of the Hitchin base: $w$ is defined only up to $W_M(\infty)$, where $M$ is the semistandard Levi factor of $P$, and for any $v \in W_G(F)$, a $(w, P)$-reduction is equivalent to a $(vw, vP)$-reduction. In the split case, one can always reduce to considering $w = e$, so our definitions reduce to those of \cite{CLI} in that case.
\end{remark}

\begin{remark}
    For classical groups the relative position $w$ of a $P$-reduction is given by comparing the order of the $\theta$-eigenspaces in $t$ with the partial order induced by the partial flag $\cE_P$ (or its associated graded).
\end{remark}

The essential property of the Hitchin base we require is that the $(w,M)$-stratification we have defined detects parabolic reductions of Higgs bundles. 

\begin{proposition}\label{prop.base.detects.p.reds}
    Let $m \in \ctM_G(K)$ such that $f(m) = (a_M, t) \in \ctAe{w}{M}(K)$. For each semistandard parabolic $P \subset G$ and $v \in W_G(\infty)$, there exists a unique $(v,P)$-reduction of $m$ if and only if $[(w,M)] \le [(v, M_P)]$ (where $M_P$ is the semistandard Levi factor of $P$).
    \begin{proof}
        If there exists a $(v,P)$-reduction $m_P = (\cE_P, \theta_P, t)$ of $m$ for any parabolic with Levi $M_P$, then its characteristic polynomial $(\chi_P(\theta_P), t)$ is an element of $\ctAt{v}{M_P}(K)$ whose image in $\ctA_{G}(K)$ is $f(m)$. Since $f(m)$ is in the image of $\ctAt{v}{M_P}$, the stratification of Corollary \ref{cor.ell.strat.tw} implies that $[(w, M)] \le [(v, M_P)]$.

        Conversely, let $(w,M) \le (w, M_P)$, so $M \subseteq P$. We will construct a $P$-reduction of $m = (\cE, \theta, t)$. Choose a preimage $(a_M, w(t)) \in \ctA_M$ for $(a,t)$. First, there exists a nonempty open $U \subseteq X_K$ (containing $\infty_K$) on which $a_M$ is $G$-regular. By the Kostant section, there exists a regular semisimple lift $Y \in \fm(U)$. Consider the centralizer $Z_Y = Z_G(Y)$, a $U$-torus subgroup of $M_U$ (which agrees with $T$ over $\infty$). 
        
        Let $\cE' \to U$ be the scheme whose sections on $V \to U$ are sections $\phi$ of $\cE_V$ such that the isomorphism $\on{Ad}_D(\cE_V) \to \fg_V$ induced by the section sends $\theta_V$ to $Y_V$. This is a $Z_Y$-torsor, which has sections \'etale-locally because $Z_Y$ is \'etale-locally conjugate to $T_U$, and because $\theta_U$ and $Y$ are both regular semisimple with $G$-characteristic polynomial $\chi_G^M(a_M)$.
        
        The torsor $\cE'$ is a reduction of $\cE_U$: there is a map
        \begin{align*}
            (\cE' \times G_U)(V) &\to \cE(V) \\
            (\phi, g) &\mapsto g.\phi
        \end{align*}
        which descends to an isomorphism on the quotient by the $Z_Y$-action. We claim now that $\theta_U$ defines a tautological section of $\on{Ad}_D(\cE'_P)$ over $U$. We have
        $$
            \on{Ad}_D(\cE'_P) = \cE' \times^{Z_Y} P \times^P \fp_D \cong \cE' \times^{Z_Y} \fp_D.
        $$
        All the local sections of this associated product $(\phi, Y)$ are equivalent, so by uniqueness they descend to a $U$-section $\theta_{P, U}$ whose image under the natural map $\on{Ad}_D(\cE'_P) \to \on{Ad}_D(\cE_U)$ is $\theta_U$. 

        Over $\infty$, the section $\theta_{P,U}$ has the same $M$-characteristic polynomial as $Y$, namely, $a_M(\infty) = \chi_M^T(w(t))$. So the triple $(\cE'_P, \theta_{P,U}, t)$ is a $(w,P)$-reduction over $U$ of $m$. By the valuative criterion of properness for $\cE/P$, the torsor extends to all of $X_K$ and $\theta_{P,U}$ likewise.
        
        Finally, we must prove uniqueness. Let $m_P = (\cE_P, \theta_P, t)$ and $m_P' = (\cE_P', \theta_P', t)$ be two $(w,P)$-reductions. The torsors $\cE_P, \cE'_P$ correspond to two different sections \[
            \sigma_P, \sigma_P' : X_K \to \cE/P,
        \] and
        by uniqueness in the valuative criterion it suffices to check they are isomorphic (and that $\theta_P \leftrightarrow \theta_P'$ under an isomorphism) at a geometric generic point $\ol{\eta} = \Spec F^s$ of the curve. But now the group $G|_{F_s}$ is split constant and the proof in \cite[Proposition 4.10.1]{CLI} applies, replacing the reference to \cite[Proposition 3.9.1]{CLI} with Proposition \ref{prop.CLI.3.8.twisted}.
        \end{proof}
\end{proposition}

\section{The compatible Hitchin fibration.}

As noted in the previous section, the locus $\bA \subset \ctA$ of pairs $(a,t)$ which are compatible (in the sense of Definition \ref{def.comp.loc}) is an open subset which surjects onto $\cA$.

\begin{definition}
    Define the \emph{compatible Hitchin fibration} to be the restriction of $f$ to the compatible locus, also denoted $f$:
    $$
        \bM := f^{-1}(\bA) \to \bA.
    $$
\end{definition}
Proposition \ref{prop.base.detects.p.reds} implies that for $f(m) \in \bAe{M}$ \emph{$M$-elliptic} and \emph{compatible} in the sense of Definition \ref{def.comp.loc}, all the parabolic reductions can be taken to be compatible (i.e., in relative position $w = e$) and to range over exactly the semistandard Levis $P \supset M$. 

Since $\bA$ surjects onto $\cA^\infty$ for sufficiently positive divisors, every generically regular-semisimple Hitchin fiber is visible in the compatible locus for some choice of $\infty$. In particular, since the parameter $t$ (and the stability parameter $\xi$ defined later) are ultimately auxiliary, all relevant geometry for our desired applications can be seen over the compatible locus.

\chapter{Weighted semistability conditions}\label{sect.weighted.stab}

In this section, we define a $\xi$-weighted analog of Harder--Narasimhan stability for Hitchin triples. Such a description is possible because, since our Higgs fields are generically regular semisimple, a Hitchin triple has at most finitely many parabolic reductions. On the compatible locus, we showed above we can identify these with distinct semistandard parabolics, and thus it is reasonable to impose positivity conditions pointing in different directions, which can be combined into a polyhedral convex $C_m$. Generalizing the usual notion of semistability, we say that a Hitchin triple $m$ is $\xi$-semistable if $\xi$ lies within this convex.

After briefly defining the polyhedral convex, $\xi$-stability, and $\xi$-semistability over the $\infty$-rigidified Hitchin base $\ctA$, we will quickly restrict to the compatible locus $\bA$, where, for rational parameters $\xi \in \a_T$, the stability can be expressed purely in terms of a convex in the relative root system and the combinatorics is thereby streamlined. The restriction to the compatible locus will then be in effect for much of the remainder of this work. Since $\bA$ surjects onto the unrigidified ($\infty$-regular semisimple) Hitchin base $\cA^{\infty}$, this restriction loses no essential part of the geometry. 

In order to obtain a well-behaved truncation of the Hitchin fibration for our purposes, it is necessary to engage at least briefly both with the notions of $\xi$-stability and $\xi$-semistability. In particular, uniqueness in the valuative criterion holds only for the $\xi$-stable stack, but existence holds only for the $\xi$-semistable stack. An insight of Chaudouard--Laumon was that for a sufficiently general parameter $\xi$, the conditions coincide. In particular, this disallows the classical notation of Harder--Narasimhan stability $\xi = 0$. In this work, we find it more natural to focus on semistability except when proving uniqueness in Chapter \ref{sect.val.unique}, and restrict to a parameter in general position when we need to avail ourselves of properties of the $\xi$-stable stack.

In the first section, we collect notation around root and coroot data in both the geometric (i.e., absolute, $\Phi(\G, \T)$) and rational (i.e., relative, $\Phi(G, S)$) character and cocharacter spaces. We also fix once and for all our choice of Euclidean measures on all such spaces. The reader is encouraged to skip that section on first reading, absorbing only the convention that geometric objects are defined via a geometric point or splitting cover and are denoted with an underline, while the relative objects are denoted without.

\section{Cocharacter spaces, cones, and measures.}\label{subsect.cochar.cones.measures}

Let $G/X$ a reductive group scheme which has a Borus pair over $X$ and choose a Galois cover $X' \to X$ splitting $G$ with Galois group $\Gamma$. Recall the cocharacter space sheaf $\ga_P$ associated to each semistandard parabolic $P \supset T$. We have then the \textit{geometric} cocharacter space $\ga_P(X')$, which we abusively denote $\ga_P$, with Galois action, and the \textit{rational} cocharacter subspace $\a_P = \ga_P^\Gamma = \ga_P(X)$. Both are real vector spaces. Note that by Proposition \ref{prop.splitrankextensions}, these spaces remain canonically isomorphic upon base change to $X_L$ for any extension $L/k$. If $M \supset T$ is the semistandard Levi factor of $P$, then we will sometimes write $\ga_M = \ga_P$ and $\a_M = \a_P$; indeed, all characters factor through the Levi quotient and its dual space consists of cocharacters landing in $A_M$, the central torus of the Levi factor.

As in \S \ref{sect.rel.reduct.gps}, we have for any containment $P \subseteq Q$ of parabolics canonical decompositions $\ga_P = \ga_Q \oplus \ga_P^Q$ and $\a_P = \a_Q \oplus \a_P^Q$. We will have need to consider various cones and roots in these spaces. Geometric data is generally denoted with an underline, and relative/rational data is generally denoted without. For example, the geometric coroots $\ul{\Phi}^\vee$ live inside $\ga_T = \ga_B$, and the geometric roots $\ul{\Phi}$ live inside $\ga_T^\star = \ga_B^\star$.

Write $\ul{\Delta}_B$ (resp. $\ul{\Delta}_B^\vee$) for the set of simple positive geometric roots (resp. coroots) associated to $B$. For any $P \supset B$ one has corresponding subsets $\ul{\Delta}_{P,B} \subset \ul{\Delta}_B$ (resp. $\ul{\Delta}_{P,B}^\vee \subset \ul{\Delta}_B^\vee$) of those simple roots (resp. coroots) corresponding to $R_u(P)$. The image of $\ul{\Delta}_{P,B}$ under the restriction $\ga_B^\star \twoheadrightarrow \ga_P^\star$ defines a basis $\ul{\Delta}_P \subset (\ga_P^G)^\star$ of \emph{simple roots} which is independent of $B$, and likewise $\ga_B \twoheadrightarrow \ga_P$ defines a basis of \emph{simple coroots} $\ul{\Delta}_P^\vee \subset \ga^G_P$ which is independent of $B$. Write $\widehat{\ul{\Delta}}_P$ for the basis of \emph{fundamental weights} in $(\ga_P^G)^\star$, dual to the coroots $\ul{\Delta}_P^\vee$. These form the extremal rays, or corners, of the dominant chambers for the various geometric Borels over $X'$ containing $T_{X'}$. 

We will make much use of certain cones in the cocharacter spaces. Denote the \emph{geometric obtuse cone} associated to $P$
$$
    \gobtcone_P := \{x \in \ga_P : \omega(x) > 0 \; \forall \omega \in \widehat{\ul{\Delta}}_P\} \subseteq \ga_P
$$
and the \emph{geometric narrow cone} associated to $P$
$$
    \gnarcone_P := \{x \in \ga_P : \alpha(x) > 0 \; \forall \alpha \in \ul{\Delta}_P\} \subseteq \ga_P.
$$
It is most natural to extend the obtuse cones from $\ga_P$ to $\ga_T$ by inverse image. Indeed as in \cite[Lemme 5.3.1]{CLI}, $\gobtcone_P$ is the intersection of $\gobtcone_Q + \ga_P^Q$ ranging over all the maximal proper geometric parabolics $Q \supset P$.

On the other hand, the canonical embedding $\ga_P \hookrightarrow \ga_T$, we can identify the narrow cone $\gnarcone_P$ with its image, the degenerate (not full-dimensional) cone consisting of $x \in \ga_T$ which satisfy $\alpha(x) > 0$ for $\alpha \in \ul{\Delta}_P$ and $\alpha(x) = 0$ for all $\alpha \in (\ga_T^P)^\star$ (that is, the roots in the Levi factor of $P$; if we choose a $B \subseteq P$, then it suffices to check $\alpha \in \ul{\Delta}_B \setminus \ul{\Delta}_P$). The narrow cones are then disjoint, and together encode the \emph{geometric Weyl fan}
$$
    \ga_T = \bigsqcup_{P \supset T_{X'}} \gnarcone_P
$$
where the union ranges over semistandard $X'$-parabolics. Its maximal facets are the Weyl chambers for the semistandard Borels. The closure of $\gnarcone_P$ is the union over cones for parabolics containing $P$. In other words, this is the simplicial structure on the apartment corresponding to $T_{X'}$ in the spherical Tits building for $G_{X'}$.

The constructions of this section are valid for any root system, and so can be performed identically with respect to the relative root system of a maximal split torus $S \subseteq T$. We call such objects \emph{relative} or \emph{rational} to distinguish them from the geometric objects, and denote them without an underline (so $\Delta_B$, $\Delta_P$, $\widehat{\Delta}_P$, $\narcone_P$, $\obtcone_P$, etc.). All the claims made above hold verbatim for the relative versions, with $P$ ranging over the set of $X$-parabolics.

\begin{definition}\label{def.cochar.space.measure.choice}
    For the sake of defining orbital integrals, we will need to define a measure on each of the $\ga_P$ compatibly. For the sake of convenience in our relative arguments, we will normalize our choice in a slightly contorted way so that the restriction to the \emph{relative} spaces has a particularly simple form. We do this once and for all as follows: on the geometric cocharacter space
    $$\ga_T = \ga_B,$$
    consider the positive-definite inner product on $\ga_T^G$, defined by
    $$
        (v, w) := \sum_{\ul{\Phi}_i} c_{i}^{-1} \sum_{\alpha \in \ul{\Phi}_i} \langle v, \alpha \rangle \langle w, \alpha \rangle,
    $$
    where $\ul{\Phi}_i$ ranges over $\Gamma$-orbits of irreducible components of the absolute root system
    $$
        \ul{\Phi} = \bigsqcup_i \ul{\Phi}_i, \qquad W = \prod_i W_i,
    $$
    (equivalently, $\Phi = \bigsqcup_i \Phi_i$ relatively) and the scalar $c_i$ is engineered on each component to convert the absolute root-sum into a relative root-sum: for any $W_i(F)$-invariant norm on the relative space $\a_T$,
    $$
        c_i := \frac{\sum_{\beta \in \Phi_i} m_\beta\| \beta\|^2}{\sum_{\beta \in \Phi_i}\| \beta\|^2}
    $$
    where $m_\beta$ is the root multiplicity of a relative root $\beta \in \Phi_i$. Fix any $W_G(\infty) \rtimes \Gamma$-invariant positive-definite form on $\ga_G$, and define the inner product on 
    $$
        \ga_T = \ga_G \oplus \ga_T^G
    $$
    to be the orthogonal inner product. Define the measure on $\ga_T$ to be the Lebesgue measure with respect to any orthonormal basis for this inner product. 
    
    It is immediate that our choice of inner product is invariant under $W_G(\infty) \rtimes \Gamma$, and the scalars are engineered such that
    $$
        (-,-)|_{\a_T^G} = \sum_{\alpha \in \Phi} \langle -, \alpha \rangle \langle -, \alpha \rangle,
    $$
    in particular, we get a $W_G(F)$-invariant inner product on $\a_T^G$. 
    
    By restriction, we then get an induced choice of $W_M(\infty) \rtimes \Gamma$-invariant inner product and measure on each $\ga_M \hookrightarrow \ga_T$, and a $W_M(F)$-invariant inner product and measure on each $\a_M$, for each semistandard Levi. The latter will even retain a simple root-sum form, cf.\ Proposition \ref{prop.root.sum.formula.inner.prod} below. Even more, for any endoscopic datum $(s, \rho_s)$ for $G$ with endoscopic group $H$, each space $\ga_H$ and $\a_H$ are equipped with the restricted inner product and induced measure. Equivariance under Ng\^{o}'s homomorphism of Weyl groups implies in particular that the geometric (resp.\ relative) endoscopic cocharacter spaces are all equipped with measures which are $W_H(\infty) \rtimes \Gamma$-invariant (resp.\ $W_H(F)$-invariant). The decompositions
    $$
        \ga_H \cong \ga_G \oplus \ga_H^G
    $$
    are moreover all orthogonal with respect to the inner product.
\end{definition}

\section{The geometric stability convex.}

\begin{definition}\label{def.geo.stab.convex}
    Let $P \subseteq G$ a parabolic subgroup over $X$ with Levi $M$ and $m = (\cE, \theta) \in \cM(L)$ for $L/k$ an arbitrary field extension. The \textit{geometric degree} $\deg(m_P)$ of a $P$-reduction $m_P = (\cE_P, \theta_P)$ is defined to be the element of the geometric cocharacter space $\ga_P$ which pairs with geometric characters $\chi \in \ga_M^\star = \ga^\star_P$ (defined over splitting $\Gamma$-cover $X'/X$) by the degree of the associated line bundle
    \[
        \deg(m_P)(\chi) = \frac{1}{|\Gamma|} \deg_{X'}(\cE_P \times^{\chi}_{X'} \G_{m,X'}).
    \]
    This is normalized to be independent of choice of splitting cover. It naturally lands in $\a_M = \ga_M^\Gamma \subseteq \ga_M$: given two $X'$-characters of $M' = M \times_X X'$ in the same Galois-orbit $\chi' = \sigma^\star \chi$, we have
    $$
        \sigma^\star(\cE_P \times_{X'}^{\chi} \G_{m,X'}) = \cE_P \times_{X'}^{\chi'} \G_{m,X'}.
    $$
    Now let $m = (\cE, \theta, t) \in \ctM(L)$ be a rigidified Higgs bundle, and assume that we are given a $P$-reduction $m_P$ of $(\cE, \theta)$ in relative position $w$. 
    The \textit{closed twisted cone} of the reduction $m_P$ is defined to be a shift of the closed obtuse cone centered at the (negative) degree of the reduction:
    \[
    \ul{C}_{m_P} = w^{-1}\left(- \deg(m_P) - \ol{\gobtcone_P} + \ga_T^P\right) \subseteq \ga_T.
    \]
    The \textit{open twisted cone} of the reduction $m_P$ is likewise defined by
    \[
        \ul{C}_{m_P}^{\str} = w^{-1}\left(- \deg(m_P) - \gobtcone_P + \ga_T^P\right) \subseteq \ga_T,
    \]
    so that $\ul{C}_{m_P} = \ol{\ul{C}_{m_P}^{\str}}$. Note these are well-defined up to the ambiguity in $[(w,P)]$: for example, a $(w,P)$-reduction can also be viewed as a $(vw, vP)$-reduction for $v \in W_G(F)$, and these give the same cone.

    The \textit{weak (geometric) stability convex} of $m$ is defined to be the intersection of the closed twisted cones
    \[
        \ul{C}_{m} = \bigcap_{m_P} \ul{C}_{m_P} \subseteq \ga_T,
    \]
    where $m_P$ ranges all parabolic reductions of $m$. The \textit{strict (geometric) convex} is defined to be the corresponding intersection of the open twisted cones 
    $$
        \ul{C}_m^{\str} = \bigcap_{m_P} \ul{C}_{m_P}^{\str}.
    $$ 
    Given $\xi \in \ga_T$ and a scheme $S/k$, we say that $m \in \ctM(S)$ is \emph{$\xi$-semistable} if $\xi \in \ul{C}_{m_s}$ for all geometric points $s \in S$, and likewise \emph{$\xi$-stable} if always $\xi \in \ul{C}_{m_s}^{\str}$. Denote $\ctM^\xi$ for the substack of $\ctM$ consisting of $\xi$-\textit{semistable} bundles. We will show below it is an open substack. Denote also
    $$
        f^\xi : \ctM^\xi \to \ctA
    $$
    for the restriction of the Hitchin morphism to the $\xi$-semistable substack.
\end{definition}

\begin{remark}
    As in \cite[Lemme 5.5.2]{CLI}, the stability convexes are fully determined by the reductions to maximal proper parabolics $P \subset G$. Note that because $P \subset G$ may not remain a maximal proper parabolic upon passage to a splitting cover, their geometric obtuse cones need \emph{not} be half-spaces.
\end{remark}

\begin{remark}
    Unlike in \cite{CLI}, even for parameters which are elliptic for $M = T$ the geometric semistability convex is \emph{not} a complementary polyhedron in the sense of Behrend: by Proposition \ref{prop.base.detects.p.reds} it will have only as many cones as there are \emph{rational} Borels containing $T$, which is in bijection with chambers in the relative root system $\Phi(G, S)$, not $\Phi(G_{\ol{s}}, T_{\ol{s}})$. We will define a more satisfactory rational convex below.
\end{remark}

By Proposition \ref{prop.base.detects.p.reds}, the stratification $\ctA_G = \bigsqcup_{[(w,M)]} \ctAe{w}{M}$ detects \textit{which} $(w,P)$-reductions $m \in \ctM(L)$ admits. 
\begin{corollary}
    For $m$ such that $f(m) \in \ctAe{w}{M}$, the intersection defining the stability convex $\ul{C}_m$ ranges over $(v,P)$ such that $[(v,M_P)] \ge [(w,M)]$, where $M_P$ is the semistandard Levi factor of $P$.
\end{corollary}

\begin{remark}\label{rmk.stab.convex.stable.extensions}
    Together with Corollary \ref{cor.parab.levi.scalar.extensions}, this implies that the semistability and stability convexes $\ul{C}_m$, $\ul{C}_m^{\str}$ associated to $m \in \ctM(L)$ is insensitive to extensions of $L/k$. 
\end{remark}

\begin{remark}\label{rmk.xi.stab.parabolic.bundles}
    One can to define $\xi$-stability beyond the regular-semisimple locus. Of course, one cannot there avail oneself of the rigidification $\ctM \to \cM$, necessary to organize the various (finitely many!) parabolic reductions into a single convex. Instead the natural object on which one can define $\xi$-stability is $\cM^{\infty\text{-par}}$, the moduli of Higgs bundles with full parabolic structure at $\infty$. On the $\infty$-regular semisimple locus $\cA^{\infty\text{-reg}}$, we have an isomorphism $\cM^{\infty\text{-par}} \cong \ctM$ over $\cM$, since any full flag at $\infty$ must be an eigenflag for some ordering of the eigenvalues. 
    
    Our notion of $\xi$-stability, translated to parabolic bundles over $\cA^{\infty\text{-reg}}$ under this isomorphism, admits a natural generalization to \emph{all} parabolic bundles $\cM^{\infty\text{-par}}$, as in, for example, \cite{heinlothschmitt}. In the split case, \textit{loc.\ cit.} also proves some of the same geometric results as \cite{CLI}, in particular the semistable reduction theorem. It is likely that our geometric results likewise extend to $\cM^{\infty\text{-par}}$ over the entire Hitchin base. This formalism seems to move us away from the relative characterizations (as in the next section) which allow us to interface cleanly with Arthur's framework for stabilization, so we do not pursue this direction in this work. It would be interesting to investigate whether such a truncation allows some part of the cohomological story from \cite{zhiweithesis} and \cite{ngofl} to extend beyond the generically regular semisimple locus. 
\end{remark}

\section{The rational stability convex.} For the sake of comparison to Arthur's formulation of weighted orbital integrals, which works entirely with relative data $A_P^{\spl}$, $\Delta_P$, etc., it will be essential to understand the interaction of the geometric stability convex over a splitting cover with the $X$-rational cocharacter spaces $\a_P$.

Note that if $P/X$ is a parabolic with Levi $M$ and $\chi \in X^\star(M)$ is a character defined over $X$, then the associated line bundle in the above definition of degree descends from $X'$ to $X$, so that the degree $\deg(m_P) \in \a_P = \a_M$ can be defined purely over $X$ via the pairing relation
$$
    \deg(m_P) = \deg_X(m_P) := \left( \chi \mapsto  \deg(\cE_P \times_X^\chi \G_{m,X}) = \frac{1}{|\Gamma|} \deg_{X'}(\cE_P \times_{X'}^{\chi} \G_{m,X'})\right).
$$

If we work over the compatible locus, then all parabolic reductions are in trivial relative position, and we can encode $\xi$-stability over the compatible locus purely in terms of rational data.
\begin{definition}\label{def.rel.stab.convex}
    For $m$ such that $f(m) \in \ctAe{e}{M}$ is compatible, we define the \textit{(weak)} \textit{rational stability convex} $C_m \subseteq \a_T$ (resp.\ the \textit{strict rational stability convex} $C_m^{\str} \subseteq \a_T$) to be the convex given by intersecting the closed cones
    $$
        C_{m_P} = -\deg_X(m_P) - \ol{\obtcone_P} + \a_T^P
    $$
    (resp.\ the open cones $C_{m_P}^{\str}$) over all reductions to semistandard parabolics $m_P$ which $m$ admits. The resulting convex is visibly independent of any splitting cover $X', \Gamma$.
\end{definition}

\begin{remark}
    Since the degree of a $P$-reduction is defined by its pairing with rational characters, it naturally lands inside the lattice $\cX_\star(M) = \cX_\star(P) \subset \a_M = \a_P$ dual to the global sections $X^\star(M) = X^\star(P)$, cf.\ Proposition \ref{prop.dual.lattice.global.invs.coinvs.mod.tors}. For the geometric part of this work, it will suffice to work always in the cocharacter space, but we will later make use of this lattice, cf.\ Definition \ref{def.gen.pos}, Definition \ref{def.cochar.sheaf}, \S \ref{sect.global.point.counts}.
\end{remark}

\begin{definition}
For $\xi \in \a_T$ and a scheme $S/k$, we say that a compatible family $m \in \bM(S)$ is \emph{$\xi$-semistable} if $\xi \in C_{m_s}$ for all $s \in S$, and \emph{$\xi$-stable} if $\xi \in C_{m_s}^{\str}$ for all $s \in S$. Denote $\bM^\xi$ the substack of compatible $\xi$-semistable bundles. This agrees with the geometric definition above on the compatible locus, by Proposition \ref{prop.rat.geo.convex.compat} below.

    We write
    $$
        f^\xi : \bM^\xi \to \bA
    $$
    for the restriction of the Hitchin map to $\bM^\xi$, the ``$\xi$-truncated Hitchin fibration.''
\end{definition}

To obtain a well-behaved stack outside the anisotropic locus, \cite{CLI} observed that it is best to work in a setting where $\xi$ is not a lattice point, so that $\xi$ never lies on the boundary of the stability convex.

\begin{definition}\label{def.gen.pos}
We say that $\xi \in \a_T$ is in \textit{general position} if for every proper $X$-parabolic $P$, the projection of $\xi$ along every $\a_T^G \twoheadrightarrow \a_{P}^G$ avoids the image of the lattice $\cX_\star(P)$.
\end{definition}

\begin{remark}
    For $\xi$ in general position, $\xi$-semistability and $\xi$-stability coincide. Indeed, equality in any of the weak inequalities of semistability would give a reduction to some maximal proper parabolic $Q$ such that the projection of $\xi$ into $\a_Q^G$ equals $-\deg(m_Q)$, contradicting general position.
\end{remark}

So, over the compatible locus, we recover exactly the picture of the stability convex as in \cite{CLI} or \cite{behrendsemistable}, in the \textit{relative root system}, i.e., a polyhedral convex with facets corresponding to semistandard $P \supseteq M$ cut out by inequality conditions pointing in the direction of different relative Weyl chambers. It is also compatible with the geometric convex defined above in the following senses.

\begin{proposition}\label{prop.rat.geo.convex.compat}
    Given an element $m = (\cE, \theta, t) \in \bM(K)$ (i.e., such that $f(m) \in \ctAe{e}{M}$ lies in the \textit{compatible} locus) and $\xi \in \a_T$, we have for any splitting $\Gamma$-cover $\pi : X' \to X$ with orthogonal projection $q : \ga_T \twoheadrightarrow \ga_{T, \Gamma} \cong \ga_T^\Gamma = \a_T$, 
    $$
        q(\ul{C}_m) = C_m = \ul{C}_m \cap \a_T.
    $$
    In particular, a rational parameter $\xi \in \a_T$ lies in $C_m$ if and only if it lies in $\ul{C}_m$. The same claims hold for the strict convexes.
    \begin{proof} 
        For any $\Gamma$-stable convex subset $\ul{C} \subseteq \ga_T$, the fact that the projection map $q : \ga_T \to \a_T$ is given by averaging implies
        $$
            q(\ul{C}) = \ul{C} \cap \a_T.
        $$
        Since $m$ is assumed compatible, all of its parabolic reductions $m_P$ can be taken to be in trivial relative position, and the set of positivity conditions defining the corresponding cone $\ul{C}_{m_P}$ is $\Gamma$-stable. In particular, on the compatible locus each weak and strict geometric cone appearing in the definition of $\ul{C}$ is a $\Gamma$-stable convex subset. 
        
        Fix a parabolic reduction $m_P$ in trivial relative position $e$. We will first check the projection compatibility for a single cone,
        $$
            q(\ul{C}_{m_P}) = C_{m_P}. 
        $$
        This then implies, by the observation above, the intersection comparison
        $$
            C_{m_P} = \ul{C}_{m_P} \cap \a_T.
        $$
        The vertex $(-\deg_{X}(m_P))$ naturally lives in $\a_M$, and is used as the vertex for both cones. The decompositions $\ga_T = \ga_P \oplus \ga_T^P$ are compatible with $\Gamma$-invariants (and likewise coinvariants) by Proposition \ref{prop.charshvs.decomp}.         By construction, the relative coroots $\Delta_P^\vee$ differs from the image under $q$ of $\Gamma$-orbit averages of geometric coroots $\ul{\Delta}_P^\vee$ by positive scalars (cf.\ \cite[\S 2.6]{KPBruhatTits}); in particular, they generate the same positive rays. It follows that $\obtcone_P$ is the image of $\gobtcone_P$ under the projection $q$, and the same for $C_{m_P}, \ul{C}_{m_P}$, as desired. The same argument with strict inequalities gives the same equalities for the open cones.

        For the full convex, then, compatibility with intersections and the cone-by-cone result gives the intersection comparison:
        $$
            \ul{C}_m \cap \a_T = \bigcap_P (\ul{C}_{m_P} \cap \a_T) = \bigcap_P C_{m_P} = C_m.
        $$
        This implies the projection claim because $\ul{C}_m$ over the compatible locus is a $\Gamma$-stable convex subset.
    \end{proof}
\end{proposition}

The projection of $C_m$ into $\a_M$ is a version of Behrend's complementary polyhedron in the relative root system \cite{behrendthesis}, by the following result. This is in some sense the primary geometric rationale for restricting to the compatible locus. 

\begin{proposition}[{{cf.\ \cite[Lemme 5.7.1]{CLI}, \cite[Lemma 3.15]{HeinlothGIT}}}] \label{prop.behrend.conds.rat}
    Let $m \in \bM(K)$ such that $f(m) \in \bAe{M}(K)$ is compatible $M$-elliptic, and let $P,Q$ be two adjacent parabolics (meaning there is a unique relative coroot $\alpha^\vee$ in $(\Delta_P)^\vee \cap (-\Delta_Q)^\vee$) with Levi $M$.\footnote{Even if the relative root system is nonreduced, the relative coroot system will only see the underlying reduced root system, so there is no ambiguity about multipliable roots.} Then
    $$
        \deg(m_Q) - \deg(m_P) \in \R_{\ge 0} \alpha^\vee.
    $$
    \begin{proof}
        Note that this is a purely rational condition, since adjacent parabolics may pull back to non-adjacent parabolics over a splitting cover. Let $R$ be the minimal parabolic over $X$ containing both $P$ and $Q$. $\deg(m_Q)$ and $\deg(m_P)$ agree when evaluated on any $X$-characters $\mu \in X^\star(R)$, so $\deg(m_Q) - \deg(m_P) \in \R \alpha^{\vee}$. We reduce to the case where $\dim \a_M^G = 1$ by replacing $G$ with the Levi quotient of $R$ and replacing $\cE_P, \cE_Q$ with their pushforwards. We may now assume that $P$ and $Q$ are opposite parabolics. Writing $\mathfrak{u} = \on{Lie}(R_u(P))$ for the unipotent radical and
        $$
            V_1 = \cE_P \times^P \on{Lie}(R_u(P)), \qquad V_2 = \cE_Q \times^Q \mathfrak{q},
        $$
        for two subbundles of $V = \cE \times^G \fg$, complementary at the generic point, we have a morphism of sheaves of the same rank
        $$
            V_1 \to V/V_2
        $$
        which is generically an isomorphism, say by using the generic torus-reduction as in Proposition \ref{prop.base.detects.p.reds}, and so $\deg(V_1) \le \deg(V/V_2)$. If $\alpha$ is a multipliable root, then this homomorphism sends the filtration $(V_1)_{2\alpha} \subset V_1$ into the corresponding filtration of $V/V_2$ (generically, hence everywhere), so we conclude
        $$
            \deg((V_1)_\alpha) \le \deg((V/V_2)_\alpha)
        $$
        and likewise for $2\alpha$ if it is also a root.
        
        The degrees of these vector bundles compute the numerical invariants $\deg(\cE_P)(\mu)$ and $\deg(\cE_Q)(\mu)$ on the common $M$-character $\mu = \det(\on{Ad}|_{\mathfrak{u}})|_M = \det(\on{Ad}|_{\fg/\mathfrak{q}})|_M \in X^\star(M)$, so we compute that 
        $$
            0 \le \deg(m_Q)(\mu) - \deg(m_P)(\mu) = x_\alpha \mu(\alpha^\vee).
        $$
        This determinant character $\mu$ is a positive multiple of roots in $\mathfrak{u}$, namely $\alpha$ and possibly $2 \alpha$, so it pairs positively with $\alpha^\vee$. It follows that $x_\alpha \ge 0$. 
    \end{proof}
\end{proposition}

\begin{corollary}[{{cf.\ \cite[Proposition 5.8.1]{CLI}, \cite[Lemma 5.2.5]{behrendthesis}}}]\label{cor.cvx.hull} 
For $m \in \bM(K)$ compatible $M$-elliptic as above, the closed convex
    $C_m$ is equal to the sum
    $$
        \Cvx\left(\{-\deg(m_Q)\}_{Q \supset M}\right) + \a_T^M + \a_G,
    $$
    with the convex hull taken over parabolics $Q$ with Levi factor equal to $M$.
    \begin{proof}
        This follows from the Proposition above exactly as in \textit{loc.\ cit.}
    \end{proof}
\end{corollary}

\begin{example}
    Consider the case of quasisplit $G = U_{X'/X}(4)$. With respect to the standard presentation of the quasi-split unitary group, as in \S \ref{subsect.ell.conv.fails}, Galois-orbits of geometric coroots for the upper triangular Borel in $G$ are 
$$
\alpha_{1} = (1, -1, 0, 0) \longleftrightarrow \alpha_{3} = (0,0, 1, -1)
$$
and
$$
\alpha_{2} = (0, 1, -1, 0)
$$
The projection of $\ga_T = \mathbb{R}^4$ onto the Galois invariants is given by 
$$
q : (a,b,c,d) \mapsto \frac{1}{2}(a-d, b-c, c-b, d-a),
$$
i.e., the matrix
$$
\frac{1}{2}\begin{pmatrix}
1 & 0 & 0 & -1\\
0 & 1 & -1 & 0\\
0 & -1 & 1 & 0\\
-1 & 0 & 0 & 1
\end{pmatrix}
$$
We use the basis for the invariant subspace $\a_T$ which is the image of the $\alpha_i$, i.e.,
$$
\beta_{2} = \alpha_{2}, \qquad \beta_{13} = \frac{1}{2} (\alpha_{1} + \alpha_{3})
$$
The rational obtuse cone $\obtcone_B$ corresponding to $B$ is the positive span of these. Its inverse image is the set of $(a,b,c,d)$ such that 
$$
   a - d > 0, \quad a+b - c - d > 0.
$$
On the other hand, the geometric cone $\gobtcone_B$ consists of the set of $(a,b,c,d)$ such that 
$$
    3a - b - c -d > 0, \quad a + b - c - d > 0, \quad a + b + c - 3d > 0.
$$
So the rational cone is indeed the image of the geometric cone under the projection, and also its intersection with the invariant subspace. Note that $\gobtcone_B$ is \textit{not} the inverse image $q^{-1}(\obtcone_B)$.
\end{example}

The following observation is elementary and natural to state, but we won't use it later.

\begin{proposition}\label{prop.tw.Hitchin.fiber.compat.comparison}
    Let $K/k$ a field extension. If $(a,t) \in \ctA(K)$ and $\xi \in \ga_T$ (with respect to some splitting cover), there is for each $w \in W_G(\infty)$ an isomorphism of groupoids between the truncated Hitchin fibers
    $$
        (f^\xi)^{-1}(a,t) \cong (f^{w(\xi)})^{-1}(a, w(t)).
    $$
    \begin{proof}
        Clearly, $m =(\cE, \theta, t) \mapsto m_w = (\cE, \theta, w(t))$ defines an isomorphism of groupoids between the non-truncated Hitchin fibers. It remains to check the relationship between $\ul{C}_m$ and $\ul{C}_{m_w}$. The set of parabolic reductions is intrinsic to $(\cE, \theta)$; if a reduction $(\cE_P, \theta_P)$ is a $(v,P)$-reduction of $m$, then
        $$
            \chi_P(\theta_P)(\infty) = \chi_P(v(t)),
        $$
        as in Definition \ref{def.rel.pos}, so the same parabolic reduction will be a $(vw^{-1},P)$-reduction of $m_w$. 

        Comparing the cones $\ul{C}_{m_P}$ and $\ul{C}_{m_{P,w}}$, we see that 
        $$
            \ul{C}_{m_{P,w}} = wv^{-1}(-\deg(m_P) - \ol{\gobtcone_P} + \ga_T^P) = w \ul{C}_{m_P}.
        $$
        So $\xi \in \ul{C}_{m_P}$ if and only if $w(\xi) \in \ul{C}_{m_{P,w}}$. Thus, the isomorphism induced by $t \mapsto w(t)$ sends $\ctM^\xi(K)$ to $\ctM^{w(\xi)}(K)$.
    \end{proof}
\end{proposition}
In particular, every truncated Hitchin fiber corresponds to a fiber over a compatible point truncated with respect to a different parameter $\xi$. We cannot simultaneously choose $(a,t)$ compatible and $\xi \in \a_T$ a rational parameter, however, so there is a mild loss of generality. For our purposes the role of $\xi$ will be auxiliary, so there is no harm in the loss.

\section{Semicontinuity.}

Now let $S$ an arbitrary affine $k$-scheme. For each $s \in S$, we obtain $G_s,B_s,T_s$ on the curve $X_{s} = X_{\kappa(s)}$. By Proposition \ref{prop.splitrankextensions}, the spaces $\a_{P_s}$ and the morphisms between them are constant over $S$. It is thus reasonable to ask, for an $m \in \ctM(S)$, how $\ul{C}_{m_s}$ varies in $s$. We prove a semi-continuity result on the behavior of these convexes in families. 
\begin{proposition}[{{cf.\ \cite[Prop. 5.6.1]{CLI}}}] \label{prop.semicts}
    Assume $S$ is Noetherian of finite type. Then the function $s \mapsto \ul{C}_{m_s}$ is \textit{lower-semicontinuous}. That is, for any $\Xi \subseteq \ga_T$, the set of points $s \in S$ such that $\Xi \subseteq \ul{C}_{m_s}$ is \textit{open} in $S$. 

    The same holds for rational stability convexes $C_{m_s}$ inside $\a_T$ over the compatible locus.
    \begin{proof}
        Fix such a $\Xi$ and $m = (\cE, \theta, t) \in \ctM(S)$. We show first that $\ul{C}_m$ is constructible and takes only finitely many values. This can be checked on irreducible components, so we reduce to $S$ irreducible. Let $\eta$ its generic point, $X_{\eta}$ the generic fiber. Assume $m_{\eta}$ is $(w,M)$-elliptic. If we show that there exists $U \subset S$ open such that $m_u$ is $(w,M)$-elliptic for all $u \in U$ and such that the degree vectors of all parabolic reductions are constant, then by Noetherian induction the convexes $\ul{C}_{m_s}$ are constructible. By Proposition \ref{prop.base.detects.p.reds} and the definition of the convex $\ul{C}_{m}$, it suffices to check that for any $(v,P)$ such that $[(v, M_P)] \ge [(w, M)]$, the $(v,P)$-reduction of $m_{\eta}$ extends to an open $U \subset S$ with constant degrees. But this is clear because $P$-reductions correspond to sections of $\cE/P$, which extend from the generic point to an open, and moreover since the degree of a line bundle is locally constant $\deg(m_{P,u})$ is constant. The relative position $v$ of a particular parabolic reduction is discretely valued and so extends from the generic point to an open. There are only finitely many semi-standard parabolics, so we can intersect to obtain a nonempty open $U$ which satisfies the desired conditions for all relevant $(v,P)$

        Since we have shown constructibility of $\{s \in S : \Xi \subset \ul{C}_{m_s}\}$, it suffices now to show it is stable under generization. Our hypotheses imply that $S$ is catenary, so we can reduce to a height 1 generization $\{s, \eta\}$. By normalization there always exists a map from a DVR with image these two points, and compatibility with pullback implies that we can check the desired property on DVRs, so we may assume that $S$ is a DVR. We wish to show that $\ul{C}_{m_s} \subset \ul{C}_{m_\eta}$. By Proposition \ref{prop.base.detects.p.reds}, the set of classes $[(v,P)]$ to which $m_\eta$ admits a reduction is naturally a subset of those for $m_s$, so it suffices to show 
        $$
            \ul{C}_{m_{P,s}} \subseteq \ul{C}_{m_{P, \eta}}.
        $$
        for each such pair $(v,P)$ to which both fibers admit reduction.
        Since the two reductions have the same relative position, the positivity conditions of these cones point in the same direction (cancelling the twist by $v$) and so the condition is equivalent to
        $$
            \omega(\deg(m_{P,s}) - \deg(m_{P, \eta})) \ge 0
        $$
        for each geometric fundamental weight $\omega \in \hat{\ul{\Delta}}_P$. 

        Since our inequality is now a purely geometric condition on the degree of an individual specialization of parabolic reductions, it can be checked after pulling all the data in question back along a splitting cover. So we may now assume that $G$ is split. We can then conclude using the Tannakian argument of \cite[Prop. 5.6.1]{CLI}, as follows. The generic parabolic reduction extends away from a set of codimension 2. Choose a representation $\rho : G \to V$ of highest weight $N\omega$ for some positive integer $N$, with $P$-stable highest-weight line $V_\omega$. We can then compute $\omega(\deg(\cE_{P,s}))$ and $\omega(\deg(\cE_{P, \eta}))$ as the degrees of the associated sub-line-bundles 
        $$
            \mathcal{V}_{\omega, \eta} := \cE_{P, \eta} \times^{P, \rho} V_\omega \subset \cE_{P, \eta} \times^{P, \rho} V
        $$
        and
        $$
            \mathcal{V}_{\omega, s} := \cE_{P,s} \times^{P, \rho} V_\omega \subset \cE_{P, s} \times^{P, \rho} V.
        $$
        The special line bundle is the saturation of the extension of the generic bundle on the special fiber, so we deduce the corresponding inequality for degrees of line bundles, from which follows the result.

        For the rational stability convexes, the characterization $C_m = \ul{C}_m \cap \a_T$ implies the same result: if $\Xi \subset \a_T$ is a subset, then $\Xi \subset C_{m_s} = \ul{C}_{m_s} \cap \a_T$ if and only if $\Xi \subset \ul{C}_{m_s}$. The set of $s$ satisfying the latter is open, so the former is likewise.
        \end{proof}
\end{proposition}

The same proof, of course, would give the same result for the strict inequalities defining $\xi$-stability.

\begin{corollary}[{{cf.\ \cite[Proposition 6.1.4]{CLI}}}]
    The stack of $\xi$-semistable bundles $\ctM^\xi$ is an open substack of $\ctM$ for any $\xi \in \ga_T$. Likewise, $\bM^\xi$ is an open substack of $\bM$ for any $\xi \in \a_T$.
    \begin{proof}
        This follows from the above by Noetherian approximation. Since $\ctM$ is locally of finite type, we may check functor-of-points criteria on affine Noetherian $S$. Moreover, any such $S$ can be written as a limit of schemes $\varprojlim_i S_i$ finite-type over $k$, and for $\ctM$ locally of finite type, we have 
        $$
        \ctM(S) = \varinjlim_i \ctM(S_i)
        $$
        by \cite[Tag. 01ZC]{stacks-project}. (To be precise, we could replace $\ctM$ with a smooth cover by a scheme). So every $m \in \ctM(S)$ descends to some $m_i \in \ctM(S_i)$, in the sense that $\ctM^{\xi} \times_{\ctM, m} S \to S$ arises from base change along $\ctM^{\xi} \times_{\ctM, m_i} S_i \to S_i$. The latter is an open immersion by the previous proposition, so the latter is also. 
    \end{proof}
\end{corollary}

\chapter{Adelic characterization of truncated Hitchin fibers}\label{sect.adelic.fibers}

In this chapter, $k$ always denotes an algebraically closed field.

\section{Preliminaries on local triviality and uniformization.}\label{subsect.loc.triv}

In order to pass between torsors and Beauville--Laszlo descent data, Chaudouard--Laumon make frequent use of the following fundamental local triviality results.
\begin{theorem}[\cite{drinfeldsimpson}]
    Let $X \to S$ a curve over an arbitrary $\Z[n^{-1}]$-scheme, and let $G \to S$ a reductive group scheme over the base. Let $\cE$ a $G$-bundle on $X$.
    \begin{enumerate}
        \item (Zariski-triviality) After an \'etale base change $S' \to S$, the torsor $\cE|_{S'}$ is Zariski-locally trivial on $X_{S'}$.
        \item (One-Point Uniformization) Assume $G$ is semisimple. Let $D \subseteq X$ any closed subscheme of $X$ which maps isomorphically to $S$, and write $U = X \setminus D$. After an \'etale base change $S' \to S$, $\cE|_{U_{S'}}$ is trivial. 
    \end{enumerate}
\end{theorem}
The proof strategy of Drinfeld--Simpson is to reduce to the case of split $G$ by an \'etale change of the base, then to show $B$-reductions exist. Since for split groups $B$ admits a filtration by $\G_{a,S}$ and $\G_{m,S}$,\footnote{See the remarks before Thm. 5.3.5 in \cite{conradrgs}, for example.} we may induct on dimension to show that every $B$-torsor is Zariski-locally trivial, with base cases given by Hilbert 90 and the comparison isomorphism $H^1_{Zar}(X, \O_X) \cong H^1_{\et}(X, \G_a)$. 

We require analogous local triviality results in the case when $G$ is defined over $X$, not constant. The simplest such result, over an algebraically closed field of constants, follows from Steinberg's theorem.
\begin{theorem}[Steinberg, cf.\ \protect{\cite[III.2.3, Thm 1' and Remarks]{serregalcoh}}]\label{thm.steinberg}
    Let $X \to \Spec k$ be a curve over an algebraically closed field $k = \ol{k}$ with generic point $\eta$ and let $G_\eta \to \eta$ be a reductive group scheme over the generic point. Then $H^1(\eta, G) = 0$. In particular, if $G \to X$ is a reductive group scheme with generic fiber $G_\eta$, then all $G$-torsors on $X$ are trivial at the generic point (hence on some Zariski open).
\end{theorem}

This implies that adelic uniformization holds for $G$-torsors on $X$, in the sense of Weil and Beauville--Laszlo. We state the result more generally for a class of generically reductive group schemes, for application in Chapter \ref{sect.coh.thms}.
\begin{theorem}[Weil]\label{thm.weil.unif}
Let $H \to X$ be a smooth separated group scheme, locally of finite type, on $X$, a curve over an algebraically closed field $k$ with function field $F$. Denote the complete local rings at $v \in |X|$ by $\O_v$, with fraction fields $F_v$. Assume the generic fiber $H_F$ is connected reductive. Then we have an isomorphism of groupoids
$$
\on{Bun}_H(X)(k) = \left[H(F) \Big\backslash \left(\prod'_{v \in |X|} H(F_v) \Big/\prod_{v\in |X|} H(\O_{v})\right)\right],
$$
where the right quotient is a set, the restricted product is with respect to the $H(\O_v)$, and the left quotient is the action groupoid.
\begin{proof}
    This is well known in the case that the group is split constant. The Beauville--Laszlo patching argument extends directly to our level of generality, using Steinberg vanishing.

    To be precise, given $\cE \in \Bun_H(X)(k)$, generic reductivity implies generic triviality of torsors by Steinberg's theorem above. This triviality extends to a Zariski-open $U$ because the torsor is locally of finite presentation over $X$. On the complement, torsors on the formal disks $\O_v$ for $v \not\in U$ are smooth, so their trivializations lift from the special fiber, essentially by Hensel's lemma, which does not require finite-typeness (cf.\ \cite[Th\'eor\`eme 11.7, Remarque 11.8(3)]{giraudcohofsch}). The discrepancies between these two trivializations on their overlaps give elements $h_v \in H(F_v)$. Changing the disk trivializations multiplies the $h_v$ by $H(\O_v)$, and changing the generic trivialization changes the family on the left by $H(F)$, so we obtain a well-defined map to the action groupoid on the right.

    Conversely, an adelic family patches by Beauville--Laszlo, which holds for separated flat group schemes locally of finite type. This level of generality follows, for example, from \cite[Proposition 3.20]{achingeryoucisBL}. The action by $h \in H(F)$ extends to an isomorphism of the torsors if and only if $(h_v')^{-1} h h_v \in H(\O_v)$ for every $v$, as required by the action groupoid. 
\end{proof}
\end{theorem}

\begin{notation}\label{not.BL.conventions}
    We fix here our conventions for Beauville--Laszlo gluing. All torsors are right torsors. If
    $$
        \beta : \cE|_{\eta} \xrightarrow{\sim} G|_{\eta}, \qquad \alpha_v : \cE|_{\Spec \O_v} \xrightarrow{\sim} G|_{\Spec \O_v}
    $$
    are the generic and formal-disk trivializations, respectively, then the transition elements on the punctured formal disks are given by
    $$
        \beta \circ \alpha_v^{-1} = g_v.
    $$
    In other words, the modifications are measured \emph{from the disk trivializations to the generic trivialization}. 

    Note that in particular, when $G = \G_m$ that taking the modification at a single point $g_v = z_v^n$ (compared to local uniformizer $z_v$) gives the $\G_m$-torsor corresponding to the line bundle $\O(-nv)$.
\end{notation}

\begin{example}
    Consider the case of $G = U_{X'/X}(1)$ associated to \'etale double cover $X' \to X$ with involution $\sigma$ and extension of function fields $F'/F$. A $G$-torsor is a line bundle $\cL$ on $X'$ equipped with an anti-linear structure $\sigma^\star \cL \cong \cL^{\vee}$ (that is, a trivialization of the norm). Beauville--Laszlo data, that is, an adelic uniformization, for this $G$ reduces to a divisor for the line bundle. The theorem above asserts that we can find a divisor $D$ on $X'$ representing $\cL$ which is anti-linear on the nose, i.e., with $D + \sigma^\star D = 0$, rather than up to linear equivalence. Concretely, given any rational section $s$ for $\cL$, then the skew-symmetry isomorphism identifies $s \sigma(s)$ with a rational function $c \in (F')^\times)^\sigma = F^\times$. Steinberg's theorem in this setting is the observation that the norm map $N : (F')^\times \to F^\times$ is surjective (cf.\ \cite[X.\S 7, Prop.\ 11]{serrelocfields}), so we can write $c = d \sigma(d)$. Then $D = \on{div}(sd^{-1})$ is a anti-linear divisor for $\cL$. 
\end{example}

The corresponding claim fails to hold for curves over finite fields because Steinberg's theorem requires algebraically closed constants, but the obstruction lies purely in the field of constants. The following is an easy consequence of effective Galois descent.

\begin{corollary}
    Let $k_0$ a finite field with algebraic closure $k$. If $G_0 \to X_0$ is a reductive group scheme over smooth proper geometrically connected curve $X_0 \to \Spec k_0$, then writing $G \to X$ for its pullback to $k = \ol{k}_0$ and $\tau = 1 \times \tau$ the arithmetic Frobenius on $X = X_0 \times_{k_0} k$, we have
    $$
        \Bun_{G_0}(X_0)(k_0) = \left[G(F) \Big\backslash G(\mathbb{A}_{X})\Big/G\left(\prod_{v\in X} \O_{v}\right)\right]^\tau,
    $$
    where the fixed locus is taken in the stacky sense.
\end{corollary}

\section{One-point uniformization.} In the semisimple case, adelic uniformization is vastly strengthened by the following one-point uniformization theorem of Heinloth. Strictly speaking, this is stronger than needed in the rest of the paper, but it feels thematically consistent to include, as the natural generalization of the result of Drinfeld--Simpson to the relative setting. 

\begin{theorem}[{cf.\ {\cite{HeinlothRelUnif}}}] \label{thm.heinlothunif}
    Let $G \to X$ be a semisimple group scheme whose algebraic fundamental group $\pi_1(G)$ is not a multiple of the characteristic $p$.\footnote{By an examination of types, this follows from assuming $p$ does not divide the order of the Weyl group.} Fix a closed point $x \in X$ and $X^\circ = X \setminus \{x\}$. Let $S$ Noetherian and $\cE \in \on{Bun}_G(S)$ a $G$-torsor on $X \times S$. Then there exists an \'etale cover $S' \to S$ such that $\cE|_{X^\circ \times S'}$ is trivial. 
    \begin{proof}
        \cite[Thm. 4]{HeinlothRelUnif} is exactly this result in the case that $G \to X$ is also simply connected. Writing 
        \[
            \begin{tikzcd}
            0 \ar[r] & Z \ar[r] & G^{sc} \ar[r] & G \ar[r] & 0
        \end{tikzcd}
        \]
        for the simply connected cover of $G$, we descend the result from $G^{sc}$ to $G$ in 3 steps of increasing generality.

        First, consider the case that $S = \Spec k$ is a field. Taking \'etale covers and spreading out reduces immediately to the case where $k$ is algebraically closed. The obstruction to lifting $\cE|_{X^\circ \times S}$ to a $G^{sc}$-torsor is class in $H^2(X^{\circ}, Z)$, which vanishes since $Z \to X$ is finite \'etale and $X^\circ$ is an affine curve (cf. \cite[Tag 03SC]{stacks-project}). So there exists a $G^{sc}$-torsor $\cE^{\circ, sc}$ on $X^\circ$ such that $$\cE^{sc, \circ} \times^{G^{sc}}_{X^\circ} G \cong \cE|_{X^\circ}.$$ By Steinberg's theorem, this $G^{sc}$-torsor is trivial on an open. Choose such a trivialization and use it to glue $\cE'$ with the trivial $G^{sc}$-bundle on a formal neighborhood of $x$. The result is a $G^{sc}$-torsor $\cE^{sc}$ on $X$ extending $\cE^{\circ, sc}$: 
        \[
            \cE^{sc}|_{X^\circ} \cong \cE^{\circ, sc}.
        \]
        By the simply connected case, $\cE^{sc}$ is trivial on $X^{\circ}$. We deduce the desired triviality of
        \[
            (\cE^{sc}|_{X^\circ}) \times_{X^\circ}^{G^{sc}} G \cong (\cE^{\circ, sc} \times^{G^{sc}}_{X^{\circ}} G)|_{X^\circ} \cong \cE|_{X^\circ}.
        \]

        Second, allow $S$ to be an arbitrary finite type $k$-scheme and take $\cE \in \Bun_G(S)$. Again we may assume $k$ algebraically closed. Choose a $k$-point $s \in S$. We will twist $\cE$ by its fiber $\cE_s \in \Bun_G(k)$ to obtain a torsor such that the lifting obstruction vanishes. To be precise, let $H := \ul{\Aut}_G(\cE_s)$, a (semisimple) inner form of $G$ over $X$. Its simply connected cover fits into the diagram
        \[
        \begin{tikzcd}
            0 \ar[r] & Z \ar[r] & H^{sc} \ar[r] & H \ar[r] & 0
        \end{tikzcd};
        \]
        here the kernel is canonically the same $Z$ because inner automorphisms extend canonically to act trivially on the center of the simply connected cover. Consider the right $H$-torsor 
        \[
            \mathcal{Q} := \ul{\on{Isom}}((\cE_s)_S, \cE)
        \]
        on $X \times S$. This is canonically trivial at $s$, so in particular $\mathcal{Q}_s$ lifts to (the trivial) $H^{sc}$-torsor on $X_s$. 
        
        Proper base change (\cite[Tag 095T]{stacks-project}) implies that restriction is an isomorphism
        $$
            H^2(X_{S^{sh}}, Z) \xrightarrow{\sim} H^2(X_s, Z),
        $$
        where $S^{sh}$ is the strict Henselization of $S$ at $s$. The obstruction to lifting $\mathcal{Q}$ to a $H^{sc}$-torsor lies in these cohomology groups and vanishes on the special fiber by the previous paragraph, so we can lift $\mathcal{Q}_{S^{sh}}$ to $\mathcal{Q}^{sc}_{S^{sh}} \in \Bun_{H^{sc}}(S^{sh})$. By spreading out, this lift extends to some \'etale neighborhood $S_1 \to S$ of $s$.  The simply connected special case then implies that $\mathcal{Q}^{sc}_{S_1}$ is trivial on $X^{\circ}$ after passage to some \'etale cover $S_2 \to S_1$, so the same holds for $\mathcal{Q}_{S_2}$. This triviality means, by definition of $\mathcal{Q}$, that
        $$
            \cE|_{X^\circ \times S_2} \cong (\cE_s)_S|_{X^{\circ} \times S_2},
        $$
        and the latter is trivial by the first case, so we deduce that the result holds on some \'etale neighborhood of $s \in S$. Repeating the argument for every closed point in $S$ and taking the disjoint union proves the theorem. 

        Finally, a standard finite-presentation argument removes the requirement that $S$ be of finite type over $k$: affine-locally, $S = \Spec A$ is the spectrum of a filtered union of finitely generated $k$-subalgebras, so any torsor $\cE$ on $X \times S$ descends to a torsor on some finite-type $k$-scheme $S_0$. The previous case shows the desired uniformization on $S_0$, then pulling back its cover along $S \to S_0$ gives the desired uniformizing \'etale cover for $S$.
    \end{proof}
\end{theorem}

This result applied to a DVR is very useful in checking valuative criteria.
\begin{corollary}\label{cor.dvr.Zar.loc.triv}
    Let $G \to X$ be a semisimple group scheme whose algebraic fundamental group $\pi_1(G)$ is not a multiple of the characteristic, and $R$ be a DVR $k$-algebra. Then for each torsor $\cE$ on $X_R$, there is an \'etale cover $R \to R'$ such that $\cE|_{R'}$ is Zariski-locally trivial on $X_R$.
\end{corollary}

\section{Truncated Hitchin fibers and Beauville--Laszlo descent.}\label{subsect.adelic.trunc.Hitchin}
In the remainder of this section we describe the points of a truncated Hitchin fiber adelically. With our eye on point counting, these notions will be defined only $X$-\textit{rationally} and only on the compatible locus of Definition \ref{def.comp.loc}. Recall that $k$ is algebraically closed, and denote $\Gr_{G,v}$ the affine Grassmannian associated to $G|_{\O_v}$.

\begin{definition}\label{def.adelic.hitchin.deg}
    For each semistandard parabolic $P \supseteq T$ with Levi $M$ and any place $v \in X(k)$, we have a map $H_{P,v} : \Gr_{G,v}(k) \to \a_M$ defined by its pairing with \textit{globally rational} characters: for all $\mu \in X^\star(M) = \ul{X}^\star(M)^\Gamma$ and any $g \in G(F_v)$,
    \[
        \mu(H_{P,v}(g)) = -\val_v(\mu(p))
    \]
    where $g = pk \in P(F_v) G(\O_v)$ is a local Iwasawa decomposition. Here, $\val_v$ is normalized by $\val_v(\pi_v) = 1$.

    Given an adelic family $(g_v)_v$, we define 
    $$
        H_P((g_v)_v) = \sum_v H_{P,v}(g_v).
    $$
    For semistandard Levi $M$, we write $H_M : M(\A) \to \a_M$ for the restriction of $H_P$ to Levi factor $M$, where $P$ is any semistandard parabolic for $M$. This is independent of $P$, by construction. 
\end{definition}

\begin{remark}      
    Note that $H_P(g_v)$ is defined by virtue of the Iwasawa decomposition in $G_{v} = G_{F_v}$ relative to $P_v$. The local groups $G_{v}$ can be more split than $G$, so note that our definitions of $H_{P,v}$ and $H_P$ include an implicit projection onto $\a_M = \ul{\a}_M^\Gamma = \ul{\a}_{M,\Gamma}$: we have evaluated only on \textit{globally} defined characters.
\end{remark}

\begin{proposition}\label{prop.deg.p.red.BL}
    Suppose that $\cE$ is a $G$-torsor on $X$ which corresponds to an adelic family $(g_v)_v \in \prod'_v \Gr_{G,v}(k)$ as in Theorem \ref{thm.weil.unif}. Fix a parabolic $P$ and suppose that we have Iwasawa decompositions $g_v = p_v G(\O_v)$. Then the family $(p_v)_v$ patch together to give a $P$-reduction $\cE_P$ of $\cE$, and its degree is given by
    \[
        \deg(\cE_P) = H_P((g_v)_v) \in \a_M.
    \]
    \begin{proof}
        First, Beauville--Laszlo patching is compatible with pushforward in the obvious sense. That is, given $\phi : G \to H$ an $X$-homomorphism of affine, smooth, finite type group schemes on $X$ and a $G$-torsor $\cE$ on $X$ with B--L gluing data $(g_v)_v$, the pushforward torsor $\cE \times^G H$ corresponds to B--L data $(\phi(g_v))_v$. This follows from the construction. 
        
        Now Beauville--Laszlo gluing for $P$ (a smooth, affine group scheme over $X$) implies that the family $(p_v)_v$ do combine to give a $P$-torsor $\cE_P$ on $X$ extending our fixed generic trivialization. Since $p_v \equiv g_v \in \Gr_{G,v}(k)$, the same family of elements viewed as transition functions for a $G$-torsor gives a torsor isomorphic to $\cE$, so by the compatibility observed above, $\cE_P \times^P G \cong \cE$. 

        To confirm the claim about degrees, let $\mu \in X^\star(P)$ be an $X$-rational character. The compatibilities above imply that $\cE_P \times^{P, \mu} \G_m$ is the $\G_m$-torsor corresponding to B--L data $(\mu(p_v))_v$, and by the elementary case of $\G_m$ we have
        $$
            \mu(\deg(\cE_P)) = \deg(\cE_P \times^{P, \mu}_X \G_m) = -\sum_{v} \val_v(\mu(p_v)) = \mu(H_P((g_v)_v)).
        $$
    \end{proof}
\end{proposition}

\begin{definition}
    Let $(a,t) \in \bAe{M}(k)$ be a point in the compatible base, with unique lift $(a_M, t) \in \bA_M(k)$. Then $\cX_{(a,t)}^\xi$ is the set of pairs $(Y, (g_v)_{v \in |X|})$ as follows:
    \begin{enumerate}
        \item $Y \in \fm(F)$ is a $G$-regular semisimple $M$-elliptic element \textit{which is conjugate by an element of $M(F_{\infty})$ to $t_a$.} (In particular, $\chi_G(Y) = a|_{\eta}$).
        \item $(g_v)_{v \in |X|}$ is a family of affine Grassmannian elements $g_v \in \Gr_{G,v}(k)$ such that 
        \begin{enumerate}
            \item $g_v$ is trivial cofinitely often,
            \item ($D$-integrality) $\on{Ad}(g_v^{-1})Y \in z_v^{-d_v} L^+_{\fg,v}(k)$ (in particular, at $\infty$ we have $\on{Ad}(g_{\infty}^{-1})Y \in L^+_{\fg,v}(k)$ and so $\chi_G(Y) \in \fc(\O_{\infty})$), and            
            \item (weighted semistability) $\xi$ lands within $C_g$, the convex defined by the intersection of the cones
            \[
                - (H_{P}((g_v)_v)) - \ol{\obtcone_P} + \a_T^P
            \]
            over all the rational semistandard parabolics $P$ in $G$ containing $M$. 
        \end{enumerate}
    \end{enumerate}

    Note that there is a natural action of $M(F)$ on this set given by
    \[
        m.(Y, (g_v)_v) = (\on{Ad}(m) Y, (mg_v)_v);
    \]
    it preserves the weighted stability condition: normality of $N_P$ implies that if $g_v = p_v k_v$ is an Iwasawa decomposition of $g_v$ then $m g_v = (m p_v) k_v$ is an Iwasawa decomposition of $m g_v$, and so for any $P \supseteq M$,
    \[
        \mu(H_P((mg_v)_{v})) = -\sum_v \Big(\val_v(\mu(m)) + \val_v(\mu(p_v))\Big) =\mu(H_P((g_v)_v)) - \sum_{v} \val_v(\mu(m)),
    \]
    and the product formula ensures cancellation. 
    
    We denote $[M(F) \backslash \cX_{(a,t)}^\xi]$ the action groupoid.
\end{definition}

\begin{proposition}[{{cf.\ \cite[\S 7]{CLI}}}]\label{prop.BLdescent.weighted}
    Let $(a,t) \in \bAe{M}(k)$ (recall $k = \ol{k}$). There is an equivalence of categories between the $k$-points of the truncated Hitchin fiber ${f^\xi}^{-1}(a,t)(k)$ and $[M(F) \backslash \cX_{(a,t)}^\xi]$.
    \begin{proof}
    Given a Higgs bundle $(\cE, \theta, t)$ lying over $(a,t)$, the bundle $\cE$ is generically trivial by Steinberg's Theorem \ref{thm.steinberg}. Picking such a trivialization, the generic fiber $\theta_\eta$ corresponds to an element $Y \in \fg(F)$ with $\chi_G(Y) = a_\eta$. By Proposition \ref{prop.CLI.3.9.twisted}, there exists $Z \in \fm(F)$ such that $\chi_G(Z) = a_\eta$ and $Z$ is conjugate by $M(F_\infty)$ to $w(t_a)$. Since $Y$ and $Z$ have the same characteristic polynomial, they are conjugate in $G(F)$, so we may assume the trivialization of $\cE$ sends $\theta_\eta$ to $Z$ and we have condition (1).

    Choosing trivializations of $\cE$ also on the formal neighborhoods $\O_v$, the transition data $g_v \in \Gr_{G,v}(k)$ between these neighborhoods and the generic one gives a well-defined family of affine Grassmannian elements $(g_v)_v$. The condition that $\theta \in H^0(X, \on{Ad}_D(\cE))$ corresponds exactly to the $D$-integrality condition (2)(b). 

    It remains to prove that the weighted stability condition translates to (2).(c). We claim that $\deg(m_P) = \sum_{v \in |X|} H_P(g_v)$. That is, if we write $g_v \in p_v L^+ G(k)$ for the Iwasawa decomposition of $g_v$, then $\chi_P(Z) = \chi_P(w(t_a))$.
    
    We claim that the family $(Z, (p_v)_v, t)$ define a $(w,P)$-reduction of $m$. By Beauville--Laszlo patching we can glue the $(p_v)_v$ with the trivial generic torsor to obtain a $P$-reduction $\cE_P$ of $\cE$. The $D$-integrality conditions for $(g_v)_v$ apply equally to the $(p_v)_v$, so $Z \in \fm(F)$ corresponds to a section of $\on{Ad}_D(\cE_P)$. So we have found a $P$-reduction of $m$ which is in relative position $w$. By Proposition \ref{prop.deg.p.red.BL}, the degree of this $P$-reduction is given by the formula $\deg(m_P) = H_P((g_v)_v)$, so (2).(c) corresponds to $\xi$-stability.
    \end{proof}
\end{proposition}

We will find a streamlined version of this description useful in our arithmetic calculations. Let $(a,t) \in \bAe{M}(k)$ with preimage $a_M \in \bA_M(k)$. Fix $Y \in \fm(F)$ a $G$-regular semisimple element with characteristic polynomial $a_{M,\eta}$, for example, by taking a Kostant lift, and write $J_{Y}$ for its centralizer in $M$. Consider the set $\cX_{(a,t)}^{\xi, Y}$ of families $(g_v)_v$ satisfying condition (2) above with respect to $Y$. This set has an action of $J_Y(F)$ by restriction from the action of $M(F)$ as above.

\begin{proposition}\label{prop.adelic.hitchin.fiber.reduced.kostant}
    Under the notations as above, there is an equivalence of categories between the $k$-points of the truncated Hitchin fiber ${f^{\xi}}^{-1}(a,t)(k)$ and the action groupoid $[J_{Y}(F) \backslash \cX_{(a,t)}^{\xi, Y}]$.
    \begin{proof}
        The natural functor 
        $$
            [J_{Y}(F) \backslash \cX_{(a,t)}^{\xi, Y}] \to [M(F) \backslash \cX^{\xi}_{(a,t)}]
        $$
        is clearly fully faithful. Any two choices of $Y$ are rationally conjugate by \cite[Lemme 3.8.1]{CLI}, so it is also essentially surjective.
    \end{proof}
\end{proposition}

\part{\texorpdfstring{Geometric properties of $\xi$-stability}{Geometric properties of xi-stability}}

\chapter{\texorpdfstring{$\xi$-instability, canonical reductions, and the $\Theta$-stratification}{xi-instability, canonical reductions, and the theta-stratification}}\label{sect.can.red}
\chaptermark{\texorpdfstring{$\xi$-instability and canonical reductions}{xi-instability and canonical reductions}}

In this chapter, we return to our running conventions from \S \ref{subsect.notation.running.assumptions}, in particular, that $k$ be an arbitrary field.

The primary goal of this chapter is to show that $\bM^\xi$ is the open part of a certain type of stratification of $\bM$ which generalizes, in the case $\xi = 0$, the Harder--Narasimhan--Shatz stratification of $\on{Bun}_G$ and $\cM$. We begin in the first section by defining canonical ``Harder--Narasimhan'' $P$-reductions and a degree of instability for compatible Hitchin triples relative to $\xi$, in the style of \cite{CLI} and Behrend's thesis \cite{behrendthesis}. We then prove the basic semicontinuity and extension results about these invariants. In the second section, we give a general characterization of ``filtered objects'' for the Hitchin fibration, relating them to parabolic reductions in families.

In the third section, we use this characterization to show that the $\xi$-canonical reductions we have defined can be packaged into a $\Theta$-stratification of the Hitchin fibration in the sense of Halpern--Leistner \cite{HalpernLeistnerInstability}, whose minimal stratum is $\bM^{\xi}$. In Chapter \ref{sect.val.exist}, this will allow us to apply a general semistable reduction theorem of Alper--Halpern-Leistner--Heinloth \cite{AHLHModuliForStacks} for such stratifications: since the untruncated Hitchin fibration satisfies the existence half of the valuative criterion, its open $\xi$-semistable part does likewise. While the introduction of the $\Theta$-filtration machinery is a minor detour from the formalism for $\xi$-stability introduced thus far, it packages the essential geometric and combinatorial properties of $\xi$-stability we will require. It also bypasses the need to adapt by hand the rather technical inductive semistable reduction argument which Chaudouard--Laumon use (based on the method of Langton \cite{langton}, which is also the heart of \cite{AHLHModuliForStacks}) to prove existence in the valuative criterion.

\section{\texorpdfstring{The $\xi$-canonical reduction.}{The xi-canonical reduction.}}\label{sect.xi.canonical.red}

Given $\xi \in \a_T$ (relative cocharacter space), the rational stability convexes allow us to define measures of $\xi$-instability and $\xi$-canonical parabolic reductions for each compatible Hitchin triple. 

Let $d$ be the distance in $\a_{T}$ defined in Definition \ref{def.cochar.space.measure.choice}. In particular, it consists of an arbitrary positive-definite inner product in the central direction and on $\a_T^G$ is induced from the canonical Weyl-invariant bilinear form
$$
    (v,w) := \sum_{\alpha \in \Phi} \langle v, \alpha \rangle \langle w, \alpha \rangle.
$$
This is a nondegenerate inner product on $\a_T^G$, but sees only the projection into this space. In fact, this defines by restriction a family of compatible inner products on the various $\a_M^G$:
$$
    (-,-)|_{\a_M^G \times \a_M^G} = \sum_{\alpha \in \Phi \setminus \Phi_M} \langle -, \alpha \rangle \langle -, \alpha \rangle.
$$
We write $d$ to denote the infimum distance between closed subsets of these spaces. Recall that we have a family of direct sum decompositions of these spaces corresponding to intermediate Levis (or parabolics).

\begin{proposition}\label{prop.root.sum.formula.inner.prod}
    Given an inclusion of semistandard Levis $M \subset L$, the decomposition $\a_M = \a_L \oplus \a_M^L$ is orthogonal with respect to the inner product above.
    \begin{proof}
        Let $x \in \a_L$ and $y \in \a_M^L = \ker(\a_M \to \a_L)$. Both pair trivially with all roots in $\Phi_M$ and $y$ has trivial central projection, so the inner product simplifies as we claimed above. We have (rational) Weyl groups $W_M \subset W_L$. The element $x$ is fixed by $W_L$, and both are fixed by $W_M$. Decompose the set of roots $\Phi$ into $W_L$-orbits. For each such orbit $\O$, the contribution to $(x,y)$ is
        $$
            \sum_{\alpha \in \O} \langle x, \alpha \rangle \langle y, \alpha \rangle.
        $$
        Since $x$ is $W_L$-invariant, its pairing with each $\alpha \in \O$ is constant, and we can write this sum as 
        $$
            \langle x, \alpha_0 \rangle \sum_{\alpha \in \O} \langle y, \alpha \rangle = \langle x, \alpha_0 \rangle \left\langle y, \sum_{\alpha \in \O} \alpha \right\rangle
        $$
        for any representative of the orbit. By definition of $\a_M^L$, the element $y$ vanishes when paired with any $W_L$-invariant character (which would be an element of $\a_L^\star$), so this expression is 0 for any orbit of roots. We conclude $(x,y) = 0$. 
    \end{proof}
\end{proposition}

\begin{notation}\label{not.proj.aPQ}
    Given vector $\xi \in \a_T$ and semistandard parabolics $P \subseteq Q$, we write $\xi_P^Q$ for the image of $\xi$ under the natural projection $\a_T \twoheadrightarrow \a_P^Q$. 
\end{notation}

We can now define a measure of the $\xi$-instability of a compatible Hitchin triple $m$. 

\begin{definition}\label{def.xi.deg.HN}
    Given $\xi \in \a_T$ and $m \in \bM(K)$ for a field $K/k$, define the 
    \emph{$\xi$-degree of instability} $\deg_{i,\xi}(m)$ (or $\deg_i(m)$ if $\xi$ is understood) to be the distance $d(\xi, C_m)$ from $\xi$ to $C_m$. Define the \emph{$\xi$-Harder--Narasimhan point} of $m$ to be the unique closest point in $C_m$ to $\xi$.
\end{definition}

\begin{remark}
    Note that since $C_m$ is invariant under central translation in $\a_G$, we would obtain the same definition if we worked in $\a_T^G$. 
\end{remark}

\begin{remark}
    One could define such a distance for the geometric stability convex and arbitrary $\xi \in \ga_T$ outside the compatible locus, at the cost of having to track additional $w$-twists. This would extend the definition above.
\end{remark}

\begin{remark}
    Note that even when $\xi = 0$, this is superficially different from the standard definition due to Behrend, which is
    $$\max_{\cE_P} \{\deg(\cE_P)(\det_P) \}$$
    for $\cE_P$ ranging over all parabolic reductions ($\theta$-stable ones for a Higgs bundle). They encode basically the same information. Though ours is a real number, it takes discrete values (since $\deg(\cE_P) \in \cX_\star(P)$ always lies at a lattice point and there are only finitely many possible face directions, cf.\ Proposition \ref{prop.can.parab} below) and jumps semicontinuously by Proposition \ref{prop.semicts}. 
\end{remark}

Next, we wish to define a $\xi$-canonical semistandard parabolic and $\xi$-canonical reduction for a compatible Higgs triple. In terms of our stability convex, the condition of the canonical reduction being ``maximally $\xi$-destabilizing'' translates to the following condition. Compare also \cite[Prop. 3.13]{behrendsemistable}.

\begin{proposition}[{{cf.\ \cite[Prop. 6.3.2]{CLI}}}]\label{prop.can.parab}
    Let $m \in \bM(K)$ have $f(m) \in \bAe{M}$, and write $\varrho \in C_m$ for the $\xi$-HN point of $m$. Then there exists a unique parabolic $P \supseteq M$ such that, writing $m_P$ for the associated $P$-reduction of $m$,
    \begin{enumerate}
        \item $\varrho$ lies in the supporting affine space
        $$
            -\deg(m_P) + {\a_T^P} + \a_G
        $$
        and
        \item $\xi - \varrho \in {\narcone_P} \cap \a_T^G$, where
        $$
            {\narcone_P} := \{p \in \a_P : \alpha(p) > 0 \; \forall \alpha \in \Delta_P\}
        $$
        is the narrow cone, i.e., the Weyl facet corresponding to $P$.
    \end{enumerate}
    The reduction $m_P$ to this parabolic satisfies
    \begin{equation}
        \xi - \varrho = \xi_P^G + \deg(m_P)^G,\label{eqn.xi.dist.canonical}
    \end{equation}
    where we have used the notational conventions of Notation \ref{not.proj.aPQ}.
\begin{proof}
    The first claim follows identically to the proof of \emph{loc.\ cit.}\ or \cite[Lemma 5.3.12]{behrendthesis}, performed in the relative root system. Note in particular that the rays in $\a_T$ from $\varrho$ are stratified by semistandard parabolics, since
    $$
        \a_T = \bigcup_{P \supseteq T} {\narcone_P}.
    $$
    So the vector $\xi - \varrho \in {\narcone_P}$ for some semistandard $P$. Because of the $M$-ellipticity of $m$, $C_m$ is invariant under translation in the ${\a_T^M}$-direction, and so the nearest point $\varrho$ will satisfy $\xi - \varrho \in \a_M$. In particular, the direction $\xi - \varrho$ will necessarily correspond to a parabolic $P$ which contains $M$.

    For the last claim, note that (1) is equivalent to 
    $$
        \varrho_P^G = - \deg(m_P)^G,
    $$
    and (2) implies
    $$
        (\xi - \varrho)_T^P = 0.
    $$
    Adding gives the result.
\end{proof}
\end{proposition}

\begin{figure}
\centering
\includegraphics[scale=0.2]{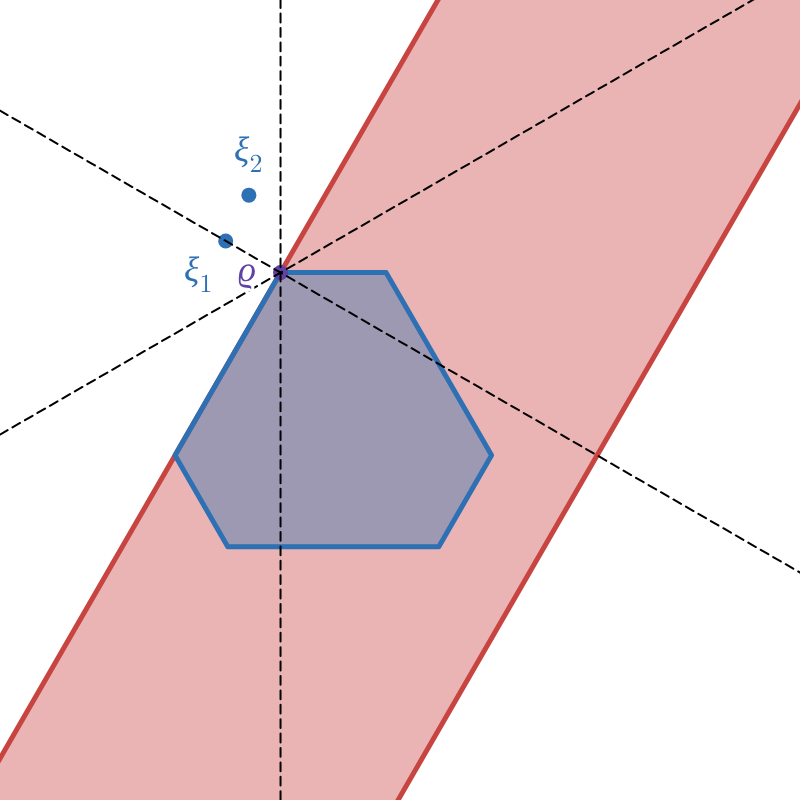}
\caption{Two stability convexes and $\xi$-HN points for $\SL_3$}
\label{fig.stab.convexes}
\end{figure}

The stratification of rays from $\varrho$ is illustrated for $\SL_3$ in Figure \ref{fig.stab.convexes}, which depicts two stability convexes, a bounded one associated with a $T$-elliptic Hitchin triple and an unbounded one associated with an intermediate Levi $T \subsetneq M \subsetneq G$. The decomposition of the rays from $\varrho$ into narrow cones is denoted with dashed lines. $\varrho$ is the $\xi$-HN point in the bounded convex for both $\xi_1$ and $\xi_2$, while the $\xi_2$-HN point of the unbounded convex is different (unmarked). 

\begin{definition}\label{def.can.parab}
    For $\xi \in \a_T$ and $m \in \bM(K)$ with $f(m) \in \bAe{M}$, define the \emph{$\xi$-canonical parabolic $P$} (resp.\ the \emph{$\xi$-canonical reduction $m_P$}) to be the unique parabolic $P \supseteq M$ (resp.\ the reduction $m_P$ of $m$ to that parabolic) satisfying conditions (1), (2) of Proposition \ref{prop.can.parab}.
\end{definition}

\begin{example}\label{ex.weak.ss.SL3}
    The strict positivity in the second condition tells us which, among the set of $Q$-reductions corresponding to closed facets containing the $\xi$-HN-point, actually has $\xi$-positive numerical invariants. In particular, the canonical parabolic is \emph{not} necessarily the parabolic such that $\varrho$ lies in the interior of the face of $C_m$ corresponding to $P$; instead the important condition is which Weyl facet $\xi - \rho$ lies interior to. This is comparable to Behrend's characterization of the canonical parabolic as the maximal parabolic among the set of those which maximize the instability. Compare the direction condition in \cite[Prop. 5.3.14]{behrendthesis}. 
    
    For example, consider in Figure \ref{fig.stab.convexes} the bounded convex as a stability convex $C_m$ for $\SL_3$. The Higgs bundle $m$ has the same $\xi$-HN-point for  $\xi = \xi_1, \xi_2$ and the same degree of instability. However, the $\xi_2$-canonical reduction is a Borel reduction, while the $\xi_1$-canonical reduction is to an intermediate parabolic: writing $\xi_1 = (\xi_1^1, \xi_1^2, \xi_1^3)$ with respect to the standard basis of cocharacters, the arrangement in the Figure corresponds to a vector bundle which admits a full flag $0 \subset \cE_1 \subset \cE_2 \subset \cE_3 = \cE$ such that
    $$
        \deg(\cE_1) + \xi_1^1 = \deg(\cE_2/\cE_1) + \xi_1^2 > \deg(\cE_3/\cE_2) + \xi_1^3;
    $$
    the first equality implies that the numerical invariant corresponding to root $\alpha = e^1 - e^2$ is equal to zero, not positive, and so the $\xi_1$-canonical reduction is the partial flag $0 \subset \cE_2 \subset \cE_3$. 
\end{example}

In what follows, we will need the following purely distance-theoretic characterization of the canonical reduction, a converse to the final claim of Proposition \ref{prop.can.parab}.

\begin{proposition} \label{prop.can.red.dist.characterization}
Let $P \supset T$ be a semistandard parabolic.
    \begin{enumerate}
    \item  For any $x \in \narcone_P \cap \a_P^G$, the point $0$ is the point in the closure of the negative obtuse cone $(- \ol{\obtcone_P})$ closest to $x$.
    \item Let $m \in \bM$ with $P$-reduction $m_P$. Suppose that 
    $$
        \nu := \xi_P^G + \deg(m_P)^G \in \narcone_P \cap \a_T^G.
    $$
    Then
    $$
        \deg_i(m) \ge \|\nu\|,
    $$
    with equality if and only if $m_P$ is the canonical reduction.
    \end{enumerate}
    \begin{proof}
        For (1), write $\{\omega_i\} = \widehat{\Delta}_P^\vee$ for the image of the fundamental coweights, the basis for $\a_P$ dual to $\Delta_P$. On the other hand we have $\{\alpha_j^\vee\} = \Delta_P^\vee$, the basis of simple positive coroots for $P$, dual to the fundamental weights. By definition the narrow cone is
        $$
            \narcone_P = \R_{> 0} \langle \widehat{\Delta}_P^\vee \rangle + \a_G,
        $$
        and the closure of the negative obtuse cone is
        $$
            - \ol{\obtcone_P} = \R_{\le 0} \langle \Delta_P^\vee\rangle + \a_G. 
        $$
        Since the closest point to $x$ necessarily also lives in $\a_P^G$, we can quotient this central direction. It follows from our definition of the inner product that $(\omega_i, \alpha_j^\vee) = 0$ if $i \neq j$, and $(\omega_i, \alpha_i^\vee) > 0$. In particular, we can see from this that for any $x \in \narcone_P \cap \a_P^G$ and any $y \in - \ol{\obtcone_P} \cap \a_P^G$, the inner product satisfies $(x,y) \le 0$, with equality only if $y = 0$. In particular, then
        $$
            \|x - y\|^2 = \|x\|^2 + \|y\|^2 - 2(x,y) \ge \|x\|^2,
        $$
        with equality only if $y = 0$. So $y = 0$ uniquely minimizes the distance. This proves (1).

        Now projecting (1) into $\a_P^G$ and translating it by $\deg(m_P)^G$, we see that under the hypothesis of (2) the closest point to $\xi^G$ in the single cone $C_{m_P}$ is 
        $$
            \xi^P - \deg(m_P)^G = \xi^G - \nu.
        $$
        The picture is invariant under central translation, so $d(\xi, C_{m_P}) = \|\nu\|$. Since $C_m \subseteq C_{m_P}$, this gives the inequality, with equality precisely when $\xi - \nu \in C_m$, i.e., $\xi - \nu = \varrho$ is the $\xi$-H--N point. By assumption $\nu$ lies in the narrow cone for $P$, and by the defining formula $\xi - \nu$ lies in
        $$
            -\deg(m_P) + \a_T^P + \a_G.
        $$
        Thus, $m_P$ satisfies the two conditions of Proposition \ref{prop.can.parab} and is the $\xi$-canonical reduction.
        
        Conversely, for the canonical reduction, we have (\ref{eqn.xi.dist.canonical}), so
        $$
            \xi - \varrho = \nu
        $$
        and equality holds. 
    \end{proof}
\end{proposition}

It follows from Proposition \ref{prop.semicts} that the degree of $\xi$-instability is upper-semicontinuous in families. We have a partial converse (whose proof follows \cite[7.1.3]{behrendthesis}):

\begin{proposition}[{cf.\ \cite[7.1.3]{behrendthesis}}]\label{prop.can.P.loc.const}
    Let $m \in \bM(R)$ be a compatible Hitchin triple over a DVR $R$ with generic point $\eta$ and special point $s$. Assume $\deg_i(m_s) = \deg_i(m_\eta)$. Then the $\xi$-canonical parabolics $P_\eta, P_s$ satisfy $P_\eta = P_s$, and the canonical reductions glue to a corresponding reduction of $m$ over $R$.
    \begin{proof}
        Let $\varrho \in \a_T$ denote the $\xi$-HN point of $m_\eta$. Lower-semcontinuity implies $C_{m_\eta} \supseteq C_{m_s}$. The only way we can have $\deg_i(m_s) = \deg_i(m_\eta)$ is if $\varrho \in C_{m_s}$ is the $\xi$-HN point for $m_s$ also. Conditions (1) and (2) of Proposition \ref{prop.can.parab} are intrinsic to the point $\varrho$, so the Weyl fan face in which $\xi - \varrho$ lands determines the canonical parabolics $P_\eta = P_s$.
        
        Now writing $P = P_\eta$, it suffices to show that the canonical reduction of $m_\eta$ and the unique $P$-reduction of $m_s$ glue to a reduction of $m$ over $R$. Write $m_\eta = (\cE_\eta, \theta_\eta, t)$ and likewise for $m_s$. It suffices to show the gluing condition for bare $G$-bundles (cf.\ \cite[Prop. 7.1.3]{behrendthesis} applied to the group schemes $G^{\cE}$, or \cite[Lemma 4.4.2]{heinlothschmitt}), as the condition of $\theta$-stability can be checked fiberwise. As in \emph{loc.\ cit.}, consider the short exact sequence of Lie algebras associated to the canonical reduction $\cE_{P, \eta}$:
        $$
            \begin{tikzcd}[ampersand replacement=\&]
                0 \ar[r] \& \Ad(\cE_{P, \eta}) \ar[r] \& \Ad(\cE_\eta) \ar[r] \& \fq_\eta \ar[r] \& 0
            \end{tikzcd}.
        $$
        We can by properness of relative Quot schemes extend this across the special fiber to a short exact sequence of coherent sheaves
        $$
            \begin{tikzcd}[ampersand replacement=\&]
                0 \ar[r] \& \fr \ar[r] \& \Ad(\cE) \ar[r] \& \fq \ar[r] \& 0
            \end{tikzcd}
        $$
        with $\fq$ coherent, flat over $R$ and $\fr$ defined to be the kernel. Behrend shows that $\fq$ is locally free at the generic point of the special fiber, so, writing $U \subset X_R$ for the maximal open on which $\fq$ is locally free (i.e.\ on which $\fr$ is a subbundle), it contains $X_\eta$ and $V := U \cap X_s$ is nonempty. 

        By \cite[Lemma 7.1.2]{behrendthesis}, there is a parabolic reduction $\cE_{P,U} \subset \cE$ whose Lie algebra (adjoint subbundle) is $\fr_U$. Properness implies that we can extend $\cE_{P,V} = \cE_{P,U}|_{V}$ to a $P$-reduction $\wtilde{\cE}_{P,s}$ of $\cE_s$ on the entire special fiber, whose adjoint bundle is the saturation of the subsheaf $\fr_s \subseteq \Ad(\cE)$. This reduction need not glue with the reduction on the generic fiber in general, but note that it is $\theta$-stable generically, hence everywhere: since $\Ad_D(\wtilde{\cE}_{P,s}) \subset \Ad_D(\cE)$ is a subbundle, the condition of a global section $\theta$ factoring through the subbundle is closed.  We conclude that $(\cE_{P,s}, \theta, t)$ is a $P$-reduction of $m_s$. Indeed, over a regular-semisimple neighborhood of $\infty_R$, the family of compatible reductions of fixed type are parameterized by a finite \'etale scheme, so relative position is locally constant and this extension is in the same relative position, i.e.\ still an $(e,P)$-reduction. Thus, it must be equal to the unique $(e,P)$-reduction of $m_s$, which is the canonical reduction.

        Degree increases under saturation, so writing $\det_P = \det(\Ad_P | \mathfrak{p})$,
        $$
        \deg(\cE_{P,\eta})(\det_P) = \deg(\fr_\eta) = \deg(\fr_s) \le \deg(\Ad(\wtilde{\cE}_{P,s})) = \deg(\cE_{P,s})(\det_P)
        $$
        We assumed degree of instability was constant, so the HN-point $\varrho$ lies in $-\deg(\cE_{P,\eta}) + \a_T^P + \a_G$ \emph{and} in $-\deg(\cE_{P, s}) + \a_T^P + \a_G$. So the two degrees are equal as elements of $\a_P^G$. In particular, $\deg(\cE_{P,\eta})(\det_P) = \deg(\cE_{P,s})(\det_P)$, so $\fr_s$ has the same degree as its saturation $\Ad(\wtilde{\cE}_{P,s})$. This implies they are equal, so $\fq$ was locally free everywhere and $\cE_{P,s}, \cE_{P, \eta}$ glue to a $P$-reduction of $\cE$ over $R$ which is $\theta$-stable everywhere. 
    \end{proof}
\end{proposition}

We summarize the behavior of our invariants in families across $\bM$.

\begin{corollary}\label{cor.semicts.invars}
    Given a family $m \in \bM(S)$,
    \begin{itemize}
        \item The Levi $M_s$ for which $f(m_s)$ is $M$-elliptic is lower-semicontinuous on $S$ with respect to containment,
        \item the degree of instability $\deg_i(m_s)$ is upper-semicontinuous on $S$, and
        \item the $\xi$-canonical parabolic $P_s \supseteq M_s$ is continuous (that is, locally constant) on $S$ for families $m \in \bM(S)$ which have constant degree of instability.
    \end{itemize}
    These are all discretely valued. 

    In particular, the condition $\deg_i(m_s) \le c$ defines for each $c \ge 0$ an open substack
    $$
        \bM_{\le c} \subset \bM
    $$
for every $c \in \R_{\ge 0}$, with $\bM_{\le 0} = \bM^{\xi}$. 
    \begin{proof}
    These have all been established (\ref{prop.semicts}, \ref{def.comp.loc}, \ref{cor.ell.strat.tw}) except the last claim. Constructibility of the canonical parabolic is clear from Proposition \ref{prop.semicts}, and equality over a DVR is the previous proposition. 

    These imply the result as in \cite[Theorem 7.2.5]{behrendthesis}. Note that as discussed there, since we are in good characteristic there is no concern about purely inseparable extensions. Concretely, as remarked above, the stability convex $C_m$ and $\xi$-canonical reduction are insensitive to extensions of $L$.

    That $\deg_i(m_s)$ takes discrete values follows from (\ref{eqn.xi.dist.canonical}) and that $\deg(m_P)$ is valued in the lattice $\cX_\star(M)$. In particular, 
    $$
        \bM_{\le c} = \bigcup_{\|\zeta - \xi\| \le c} \bM^{\zeta}
    $$
    is certainly open.
    \end{proof}
\end{corollary}

\section{Rees constructions and equivariant degenerations.}

\newcommand{\cZ}{\mathcal{Z}}
\newcommand{\Filt}{\on{Filt}}
\newcommand{\ev}{\on{ev}}
\newcommand{\gr}{\on{gr}}

In the previous section, we defined the degree of instability $\deg_{i, \xi}(m)$ in the spirit of \cite{behrendthesis} or \cite{CLI}. In this subsection, we connect this with the formalism of Halpern-Leistner, inspired by Kempf's results in GIT \cite{kempf}, in which one ranges over a space of \emph{filtered objects}. In this section, we characterize filtered objects first for $\bM \to \bA$ and then for $\bM_{\le c} \to \bA$. A filtered Hitchin triple will include canonical reduction, as well as the data of a global cocharacter $\lambda \in X_\star(S)$ such that $P_\lambda$ is the canonical parabolic (in particular picking out the Weyl facet containing $\xi - \varrho$). This cocharacter acts on the associated graded $M$-bundle as well as the reduction $\cE_{P_\lambda}$. 

Write $\Theta = [\A^1/\G_m]$, a $k$-stack. Recall that restriction to 1 defines a morphism
    \[
        \ev_1 : \Filt(\cX) \to \cX.
    \]
    Likewise, restriction along $B\G_m = [0/\G_m] \to \Theta$ defines a map to ``graded objects''
    \[
        \gr : \Filt(\cX) \to \Gr(\cX) := \Map_S(B\G_m, \cX).
    \]Our goal is to describe the stack of filtered objects with respect to $\bM \to \bA$:
\[
    \Filt_{\bA}(\bM) := \Map_{\bA}(\Theta_{\bA}, \bM).
\]
Note that all relevant filtration stacks are algebraic by \cite[Prop. 1.1.2]{HalpernLeistnerInstability}.

We begin with recalling the structure of filtered objects in $\Bun_G$. We can describe maps $\Theta \to \Bun_G$ using a version of the Rees construction, cf.\ \cite[\S 1.F.b]{HeinlothGIT}. We begin with the construction for bare torsors, which uses the following data:
\begin{itemize} 
\item a base scheme $U/k$,
\item a reductive group scheme $G/U$,
\item a $G$-torsor $\cE$ on $U$,
\item a cocharacter $\lambda : \G_{m, U} \to G$, and \item a reduction $\cE_P$ of $\cE$ to $P = P_G(\lambda) \subset G$.
\end{itemize}
By the dynamic definition of $P := P_G(\lambda)$, the $\lambda$-conjugation map
\begin{align*}
    c_\lambda : P \times \G_m &\to P \times \G_m \\
        (p, t) &\mapsto (\lambda(t) p \lambda(t)^{-1}, t)
\end{align*}
 extends uniquely to $P \times \A^1 \to P \times \A^1$. We denote this extension also by $c_\lambda$. We can then define the Rees bundle, a $P$-torsor on $U \times \A^1$
as the associated product
$$
    \on{Rees}(\cE_P, \lambda) := (\cE_P \times \A^1) \times^{P \times \A^1, c_\lambda}_{\A^1} (P \times \A^1).
$$
This has $\G_m$-equivariant structure given by acting on the equivalence class of a triple $(e, p, t)$ as 
$$
    \alpha.(e,p,t) = (e, \lambda(\alpha) p, \alpha t),
$$
so descends to a $P$-torsor on $U \times [\A^1/\G_m]$. Extending scalars to $G$ gives a $G$-torsor on $U \times [\A^1/\G_m]$. In the fiber over $0$, it has value induced from the Levi quotient of $\cE_P$, with grading given by $\lambda$. This construction commutes with every base change $V \to U$, so it gives a morphism of stacks over $U$
$$
    \on{Rees}_\lambda(U) : B_U P \to \Filt_U(B_U G).
$$
Indeed, these define open and closed embeddings, with disjoint images for distinct conjugacy classes of cocharacters. When $G$ is obtained by base change from a split reductive $k$-group, these morphisms over all $\lambda$ (without loss of generality taken dominant) are jointly surjective, and so one has the following theorem.

\begin{theorem}[{\cite[Theorem 1.4.8]{HalpernLeistnerInstability}}]\label{thm.Filt.BG.HL}
    Let $\G/k$ a split reductive group with split maximal torus $\T$ and Weyl group $\W$. Then the Rees construction gives an isomorphism
    $$
        \bigsqcup_{\lambda \in X_\star(\T)/\W} B\P_{\lambda} \xrightarrow{\sim} \Filt_k(B\G).
    $$
\end{theorem}

Specialize now to the case $U = X_R$, $G/X$ as usual. The construction in this case gives maps 
$$
    \on{Rees}_\lambda(R) : \Bun_{P_{\lambda},X}(R) \to \Filt_k(\Bun_{G,X})(R).
$$
compatible with pullback on $\Spec R$, i.e., a morphism
$$
    \on{Rees}_\lambda : \Bun_{P_{\lambda},X} \to \Filt_k(\Bun_{G,X}).
$$
Here, allowing $\lambda$ to range over just the pullbacks of globally-defined cocharacters on $X$ suffices to give a jointly surjective family, as the following result shows.

\begin{lemma}[{cf.\ \cite[Prop.\ 1.13, Lemma 3.9]{HeinlothGIT}}]\label{lem.Rees.BunG}
    The Rees construction induces an isomorphism
    $$
        \bigsqcup_{\lambda \in X_\star(S)/W_G(F)} \Bun_{P_\lambda} \xrightarrow{\sim} \Filt_k(\Bun_G),
    $$
    where $S \subseteq T$ is a maximal split $X$-torus in $G$.
    \begin{proof}
        We will descend the isomorphism from Theorem \ref{thm.Filt.BG.HL}. Let $R$ a ring and consider an arbitrary $R$-filtered object on $\Bun_G$. The restriction to $[0/\G_m]$ gives a $G_R$-torsor $\cE_0$ on $X \times R$ with a $\G_m$-action over $X \times R$:
    $$
        \lambda : \cE_0 \times \G_m \to \cE_0.
    $$
    Such an action induces a section $\lambda' \in \Hom_{X \times R}(\G_{m, X\times R}, \ul{\Aut}_G(\cE_0))$. For each local trivialization of $\cE_0$, this automorphism group is identified with $G$, and the change of frame conjugates the corresponding cocharacter of $G$. So $\lambda'$ descends to a well-defined section on $X \times R$ of the \'etale sheaf of cocharacters up to conjugacy of $G$, cf.\ Lemma \ref{lem.conj.cochars} above. 

    We identify such conjugacy classes with globally-defined dominant cocharacters, as follows. Write $Y \to X$ for a geometrically connected finite \'etale splitting $\Gamma$-cover, as in Proposition \ref{prop.splitrankextensions}. Then $X_\star(T_Y)/W_G(T_Y) \cong X_\star(\T)/W_G(\T)$. By \'etale descent, a conjugacy class of cocharacter on $X_R$ is isomorphic to a section of 
    $$
        \Map_{\text{loc.\ const.}}(Y_R, X_\star(\T)/W_G(\T))^{\Gamma}.
    $$
    The fibers of $Y_R \to \Spec R$ are geometrically connected, so such a map from $Y_R$ into this discrete set descends to a function $\Spec R \to X_\star(\T)/W_G(\T)$, locally constant because $q$ is smooth and surjective. But the original map was $\Gamma$-equivariant and $Y_R \to \Spec R$ is surjective, so the descended map must land within the $\Gamma$-invariants:
    $$
        \Map_{\text{loc.\ const.}}(Y_R, X_\star(\T)/W_G(\T))^{\Gamma} \cong \Map_{\text{loc.\ const.}}(\Spec R, (X_\star(\T)/W_G(\T))^{\Gamma}).
    $$
    
    Now recall our choice $S$ of a maximal split torus in $G$ over the entire curve $X$. Since $G$ is quasi-split, we have by \cite[Lemma 1.1.3]{KottwitzTransfer} an isomorphism
    $$X_\star(S)/W(S) \cong (X_\star(\T)/W(\T))^{\Gamma}.$$ 
    The idea is that for any $\Gamma$-fixed conjugacy class the pinned $\Gamma$-action can transport our cocharacter to some Weyl-translate, but every geometric cocharacter admits a dominant representative, and the dominant chamber corresponding to our Borel is stable under $\Gamma$, so that representative must be itself $\Gamma$-fixed.

    In summary, we have above associated to an arbitrary $R$-filtered object on $\Bun_G$ a \textit{type}, a locally constant function valued in conjugacy classes of global cocharacters on $X$
    $$
        [\lambda] : \Spec R \to X_\star(S)/W_G(S),
    $$
    which we can represent with $\lambda$ in the rational dominant chamber. For any such $\lambda$, the Rees construction lands within the full substack $$\on{Rees}_\lambda : B_X P_\lambda \longrightarrow  \Filt_X(B_X G)_{[\lambda]} \subset Filt_X(B_X G)$$ 
    of bundles with type $[\lambda]$. After pullback to $Y \to X$, this becomes the inclusion of the $\lambda$-summand, in the Theorem above. In particular, it is \'etale-locally an equivalence, hence an equivalence over $X$ by \'etale descent.\footnote{We can check full-faithfulness by checking the map of isomorphism sheaves is an isomorphism after the cover. Essential surjectivity then follows from \cite[Lemma 046N]{stacks-project}.} 

    Finally we pass to mapping stacks. For $R$ arbitrary and locally constant dominant cocharacter $\Spec R \to X_\star(S)$, decompose $R = \bigsqcup_{\lambda} R_\lambda$. The Rees map on each $R_\lambda$ becomes the composition
    $$
        \Bun_{P_\lambda}(R_\lambda) \cong B_X P_\lambda (X_{R_\lambda}) \longrightarrow \Filt_X(B_X G)(X_{R_\lambda}) \cong \Filt_k(\Bun_G)(R_\lambda),
    $$
    an equivalence onto the part with type $\lambda$. The disjoint union of these equivalences gives the lemma.
    \end{proof}
\end{lemma} 

\begin{remark}
    Here is another perspective on the argument above, closer to \cite{HeinlothGIT} and perhaps clearer, in the case that $R = \Spec \ol{k}$. 
    
    The cocharacter $\lambda'$ was valued in the group $G' = \Aut_G(\cE_0)$, which is an inner form of $G$ over $X$. The cocharacter then determines an $X$-parabolic subgroup $P'_{\lambda} \subset G'$, and $\cE' = \on{Isom}_G(\cE, (\cE_0)_{X \times [\A^1/\G_m]})$ is a $G'$-torsor on $X \times [\A^1/\G_m]$ with canonical trivialization (and thus canonical $P_{\lambda}'$-reduction) over $X \times [0/\G_m]$. As in \textit{loc.\ cit.}, a deformation-theoretic argument (which works even on $(X \times R)$-families) lifts the parabolic reduction iteratively to Artinian neighborhoods of $[0/\G_m]$, and it then algebraizes by the results of \cite{equivalgebraiz} to a parabolic reduction on $X \times [\A^1/\G_m]$. This provides a construction of the inverse isomorphism, except that it produces data not for $G$ but for some inner form.
    
    Now we must really assume that $R = \Spec \ol{k}$. Since $G$ is quasi-split there is a well-defined injective transfer of conjugacy classes of cocharacters
$$
    X_\star(G^{\cE_0})/G^{\cE_0}(F) \hookrightarrow X_\star(G)/G(F)
$$
constructed as follows. Lifting the conjugacy class to $F^s$ and transferring along an $F^s$-isomorphism $G^{\cE_0}|_{F^s} \cong G|_{F^s}$. The result is a Galois-stable conjugacy class of cocharacters. Conjugating them to lie within fixed maximal torus $T$ with maximal split subtorus $S$, the isomorphism from \cite[Lemma 1.1.3]{KottwitzTransfer} used above furnishes a lift to an $F$-conjugacy class in $X_\star(G)$, which we also denote by $\lambda$.

Since we are working over an algebraically closed constant field, Steinberg's theorem \ref{thm.steinberg} implies $G$ and $G'$ are isomorphic over a Zariski-open in $X$, and so by the valuative criterion we can transfer the parabolic $P_\lambda' \subset G'$ to a parabolic $P_\lambda \subset G$ over all of $X$, and the reduction $\cE_\lambda'$ transfers to a $P_\lambda$-reduction $\cE_\lambda$ of $\cE$.\footnote{Recall we have assumed a maximal torus $T$ exists globally, so a maximal split torus does likewise and the cocharacters extend from the generic point to all of $X$. Note the contrast with the case of $G$ a parahoric group scheme, as in \textit{loc.\ cit.}, where a maximal split torus may only exist generically.}
\end{remark}

\begin{proposition}\label{prop.filt.Hitchin}
The stack $\Filt_{\bA}(\bM)$ classifies tuples $(m, \lambda, \cE_{\lambda})$, where
\begin{itemize}
    \item $m = (\cE, \theta, t) \in \bM$,
    \item $\lambda \in X_\star(S)$, defining a parabolic $P_\lambda \supset T$, and
    \item $\cE_\lambda$ a $P_\lambda$-reduction of $m$ (in trivial relative position, i.e., an $(e, P_\lambda)$-reduction).
\end{itemize}
\begin{proof}
     Since the map $\cM \to \Bun_G$ is affine and representable,  \cite[Prop. 1.3.2]{HalpernLeistnerInstability} implies that $\Filt_k(\cM) \hookrightarrow \Filt_k(\Bun_G) \times_{\Bun_G} \cM$ is a closed immersion. Its image consists of the set of parabolic reductions compatible with the Higgs field, as in \cite[\S6.3]{AHLHModuliForStacks}, which splits as a disjoint union across $\lambda \in X_\star(S)/W(S)$ by the previous Lemma. We have a formal 2-Cartesian diagram of mapping stacks
        $$
            \begin{tikzcd}[ampersand replacement=\&]
                \Filt_{\cA}(\cM) \ar[r] \ar[d] \& \cA \ar[d] \\
                \Filt_{k}(\cM) \ar[r] \& \Filt_k(\cA)
            \end{tikzcd},
        $$
        where the right vertical map corresponds to the trivial $\cA$-family $\pi_2 : \Theta_{\cA} \to \cA$. Since $\cA \to k$ is representable, it is easy to see $\Filt_k(\cA) \cong \cA$ along the ``evaluation at 1'' map, a one-sided inverse to $\pi_2$, so we see that $\Filt_{\cA}(\cM) \cong \Filt_k(\cM)$, where the latter is an $\cA$-scheme along the map $m \mapsto f(m|_1)$. Compatibility with base change \cite[Lemma 1.1.1]{HalpernLeistnerInstability} implies $\Filt_{\bA}(\bM) \cong \Filt_{\cA}(\cM) \times_{\cA} \bA$. In particular, the disjoint decomposition across $\lambda \in X_\star(S)/W(S)$ rigidifies over $\bA$ into a decomposition over the relative position $(w,P_\lambda)$ between $t$ and $\theta|_{\cE_{P_\lambda}}$ as in Definition \ref{def.rel.pos}: the class $[(w,P_\lambda)]$ is a locally constant class on $\Filt_{\bA}(\bM)$, which we can always present as $(e, P')$ since we are over the compatible locus. That is, $\Filt_{\bA}(\bM)$ decomposes first over the conjugacy class $[\lambda]$, then over semistandard $P$ conjugate to $P_\lambda$. Choosing the standard $P_\lambda \supset B$ as base-point, the set of semistandard conjugates $P' \sim P_\lambda$ is canonically in bijection with the Weyl conjugacy class of $\lambda$, so 
        $$
            \Filt_{\bA}(\bM) = \bigsqcup_{\lambda \in X_\star(S)} \Filt^{\lambda}_{\bA}(\bM),
        $$
        where $\Filt^\lambda_{\bA}(\bM)$ parameterizes $(m, \cE_\lambda)$ consisting of a Hitchin triple and a $P_\lambda$-reduction. This is the stated decomposition. 
\end{proof}
\end{proposition}

We now wish to characterize those filtrations lying over a stratum $\bM_{\le c}$, i.e., to understand
$$
    \Filt_{\bA}(\bM_{\le c}) \subset \Filt_{\bA}(\bM) \times_{\bM} \bM_{\le c}.
$$
Since all nonzero fibers of a Rees family are isomorphic to the fiber at 1, it suffices to relate the instability of the fiber at 1 to the instability of the associated graded fiber at 0. The fundamental geometric relationship is the following compatibility among faces of the convex $C_m$. We write $$F_P(m) := C_m \cap (-\deg(m_P) + \a_T^P + \a_G),$$ the face of $C_m$ corresponding to $P$.

\begin{lemma}\label{lemma.face.compat.assoc.graded}
    Let $m_P$ be a compatible $P$-reduction of $m \in \bM$, and let $\lambda$ be a cocharacter such that $P_\lambda = P$, so that $(m, \lambda, m_P)$ define a filtered object in $\bM$ over $\A$ by Proposition \ref{prop.filt.Hitchin}. Write $m_0$ be the underlying $G$-Hitchin triple associated to $\gr(m, \lambda, m_P)$ by the Proposition above. Then
    \begin{equation}
        C_{m_0} = F_P(m) \subseteq C_m.\label{eqn.face.convex.Levi}
    \end{equation}
    It follows that if $m_P$ is the canonical reduction of $m$, then $m_0$ has the same $\xi$-H--N point, degree of instability, and canonical parabolic.
    \begin{proof}
        Pulling back our filtered object $\Theta \to \bM$ along the map $\A^1 \to \Theta$ gives a family with special fiber $m_0$ (grading forgotten) and generic fiber $m$ over $\G_m$. Semicontinuity (Proposition \ref{prop.semicts}) then gives
        $$
            C_{m_0} \subseteq C_m.
        $$
        The $\lambda$-grading identifies $m_0$ as the pushforward from a Hitchin triple for $L = Z_G(\lambda)$, the Levi quotient of $P$. In particular, its $P$-reduction and $\ol{P}$-reduction have the same degree, so its convex degenerates into the affine space
        $$
            C_{m_0} \subseteq (-\deg(m_P) + \a_T^P + \a_G).
        $$
        Combining gives $C_{m_0} \subseteq F_P(m)$.

        Conversely, write $M$ for the Levi for which $m$ (and thus also $m_0$, since they have the same characteristic polynomial) is elliptic. We therefore have $M \subseteq L$. For every $Q \subseteq P$ with Levi factor $M$, we now relate the $Q$-reductions of $m$ and of $m_0$. The compatible reduction $m_Q$ is a reduction of scalars of $m_P$ (by uniqueness). Write $\cE_Q$ for its underlying bundle. The pushforward along $P \twoheadrightarrow L$, which is $\cE_Q/U_P$, is a reduction of $\cE_L$ to $Q \cap L \cong Q/U_P$, and stable under the associated graded Higgs field $\theta_L$. Pushing this datum forward along $Q \cap L \hookrightarrow Q$ gives the unique $Q$-reduction of $m_0$,
        $$
            \cE_{0, Q} := (\cE_Q/U_P) \times^{Q \cap L} Q.
        $$
        Taking its pushforward along $Q \twoheadrightarrow M$ and using $Q/U_Q \cong M$, we can compute
        $$
            \cE_{0, Q} \times^Q M = (\cE_Q/U_P) \times^{Q/U_P} M \cong \cE_Q \times^Q M;
        $$
        that is, that $\cE_{0,Q}$ and $\cE_Q$ have the same $M$-associated graded. In particular, if we write $v_Q(m) = -\deg(m_Q)$ and $v_Q(m_0)$ likewise, then
        $$
            v_Q := v_Q(m) = v_Q(m_0)
        $$
        for those parabolics which satisfy $Q \subseteq P$.

        We will now use the convex hull characterization of Corollary \ref{cor.cvx.hull} to conclude. Applying it to $m_0$ gives
        $$
                \Cvx(\{v_Q\}_{Q \subseteq P}) + \a_T^M + \a_G \subseteq \Cvx(\{v_Q(m_0)\}_{Q}) + \a_T^M + \a_G = C_{m_0}
        $$
        where both convex hulls are taken over semistandard parabolics with Levi factor exactly $M$ and live in $\a_M^G$. We have shown above that the vertices for the smaller convex hull coincide with the corresponding vertices of $m$, so it remains to show that 
        $$F_P(m) \subseteq \Cvx(\{v_Q\}_{Q \subseteq P}) + \a_T^M + \a_G.$$ Working modulo $\a_T^M + \a_G$ and using the corresponding convex hull characterization of $C_m$, this becomes the statement that taking the convex hull over \textit{all} the minimal parabolic reductions adds no new points in the slice $A := v_P^G + \a_M^L$.
        This follows from the negative inner product pairing between $\narcone_P$ and $(- \ol{\obtcone_P})$ of Proposition \ref{prop.can.red.dist.characterization}.(1), as we argue now. For any points $u \in \narcone_P \cap \a_P^G$ and $x \in C_m$,
        $$
            (u, x) \le (u, v_P^G),
        $$
        with equality if and only if $x \in A$. Indeed, adding a positive component from any vertex not in the slice would contribute a smaller value of the functional. But if $x \in A \cap C_m$ is a point not in the convex hull $ \Cvx(\{v_Q^G\}_{Q \subseteq P})$, then there exists a functional $(v, -)$ such that
        $$
            (v,x) > \max_{Q \subseteq P} (v, v_Q^G).
        $$
        Again, $Q$ ranges over those with Levis equal to $M$. For $\epsilon > 0$ sufficiently small and $v \in \a_M^G$ generic among those which separate $x$ from the convex hull as above, $u + \epsilon v$ lies in some chamber $\narcone_Q$ for $Q \subseteq P$ with Levi $M$. The defining condition of the cone corresponding to $Q$ implies
        $$
            (u + \epsilon v, x) \le (u + \epsilon v, v_Q^G).
        $$
        Since $x$ and $v_Q^G$ lie in $A$, their inner products with $u$ coincide (indeed
        $$
            (u,x) = (u, v_P^G) = (u, v_Q^G)
        $$
        by compatibility of the reductions under pushforward $Q \subseteq P$). Cancellation gives $(v,x) \le (v, v_Q^G)$, a contradiction.

        We have shown $C_{m_0} = F_P(m)$. For the final claim, by Proposition \ref{prop.can.parab} the $\xi$-H--N point of $m$ lies in $F_P(m)$, so is also the $\xi$-H--N point of $m_0$. The canonical parabolics coincide under the degeneration as in Proposition \ref{prop.can.P.loc.const}.
    \end{proof}
\end{lemma}

\begin{proposition}\label{prop.filt.Hitchin.bdd}
    For any $c \ge 0$, the correspondence of Proposition \ref{prop.filt.Hitchin} identifies the stack $\Filt_{\bA}(\bM_{\le c})$ with the substack of triples $(m, \lambda, m_P) \in \Filt_{\bA}(\bM)$ such that 
    $$
        d(\xi, F_P(m)) \le c.
    $$
    \begin{proof}
        On each geometric fiber, a morphism $\Theta \to \bM$ lands inside $\bM_{\le c}$ if and only if its fiber over $0$ does, since $C_{m_0} \subseteq C_m$ by semicontinuity. Lemma \ref{lemma.face.compat.assoc.graded} then shows this is equivalent to $d(\xi, F_P(m)) \le c$.

        Since $\bM_{\le c} \subset \bM$ is an open substack, checking the condition on geometric fibers implies it in arbitrary families.
    \end{proof}
\end{proposition}

\section{\texorpdfstring{$\Theta$-stratification and semistable reduction.}{Theta-stratification and semistable reduction.}}

Recall the notion of $\Theta$-stratification introduced by Halpern-Leistner (\cite[Definition 2.1.1, Definition 2.1.2]{HalpernLeistnerInstability}).
\begin{definition}\label{def.theta.stratification}
    Let $\cX$ an algebraic stack locally of finite type over Noetherian algebraic space $B$.
    \begin{enumerate}
        \item A \emph{$\Theta$-stratum} in $\cX$ consists of a union of connected components $\cZ^+ \subset \Filt(\cX)$ such that $\ev_1: \cZ^+ \to \cX$ is a closed immersion.
        \item A \emph{well-ordered $\Theta$-stratification} indexed by a totally ordered set $I$ with minimal element $0$ is a cover of $\cX$ by open substacks $\cX_{\le i}$ for $i \in I$ such that $\cX_{\le i'} \subset \cX_{\le i}$ for $i' < i$, along with a $\Theta$-stratum $\cZ_i^+ \subset \Filt(\cX_{\le i})$ whose complement is the union of $\cX_{\le i'}$ for $i' < i$. We require that for every $x \in |\cX|$,
        \begin{itemize}
            \item the set $\{i \in I : x \in \cX_{\le i}\}$ has a minimal element, and
            \item the (totally ordered) set $\{i \in I : \ev_1(\cZ_{i}^+) \cap \ol{\{ x\}} \neq \emptyset \}$ is well-ordered.
        \end{itemize}
    \end{enumerate}
\end{definition}

In \cite[\S 8]{AHLHModuliForStacks}, the authors show that canonical (that is, HN) reductions of $G$-bundles define a $\Theta$-stratification on $\Bun_G$. We now collect the properties proven in this section to show that the $\xi$-canonical reduction lifts to define a $\Theta$-stratification of $\bM$ over $\bA$.

Of course, in order to define a $\Theta$-stratification, the characterization of filtered Hitchin triples in Proposition \ref{prop.filt.Hitchin} suggests that we will need to associate to every canonical reduction a cocharacter $\lambda$ defining the associated parabolic. This choice of weights on the associated graded can sometimes be done more canonically, but for our purposes the cocharacter can be chosen completely \textit{arbitrarily}, as long as they lie in the correct Weyl facet such that $P_\lambda$ is the canonical parabolic and the choice is locally constant in families of reductions. See Remark \ref{rmk.choice.of.cochars} below for further discussion.

\begin{notation}\label{not.cochar.for.parabolic}
    Fix once and for all an integral cocharacter 
    $$
        \lambda_P \in X_\star(S) \cap \a_P^G
    $$
    for each semistandard parabolic $P$ such that $P_{\lambda_P} = P$.
\end{notation}
Of course, when $P = G$ we must choose $\lambda_G = 0$. Such a $\lambda_P$ exists for each parabolic because $\narcone_P \cap \a_P^G$ is a nonempty rational cone. While the actual components in the $\Theta$-stratification will depend upon our choices above, the image stratification in $\bM$ will not.

\begin{definition}\label{def.theta.strata}
    For each $c \ge 0$, define the set of \emph{admissible pairs} $T_c$ as the finite set of $(P,d)$ such that
    \begin{itemize}
        \item $P$ is a semistandard parabolic,
        \item $d \in \cX_\star(P)^G$,
        \item $\nu_{P,d} := \xi_P^G + d \in \narcone_P \cap \a_T^G$,
        \item $\| \nu_{P,d}\| = c$.
    \end{itemize}
    Define then
    $$
        \cZ^+_c := \bigsqcup_{(P, d) \in T_c} \Filt_{\bA}^{\lambda_P, d}(\bM_{\le c}),
    $$
    where the superscripts indicate the substack of triples in Proposition \ref{prop.filt.Hitchin.bdd} with $\lambda = \lambda_P$ and specified degree $d$. 
\end{definition}

\begin{remark}
    The superscript $G$ denotes projection into $\a_P^G$, as in Notation \ref{not.proj.aPQ}. We have taken $d \in \a_P^G$ rather than $\a_P$ because the central component of the degree of $m_P$ is already determined by $m$, so any other choice would make the stack empty. Likewise, the stack would also be empty if $d$ did not lie in the projected lattice $\cX_\star(P)^G$. 
\end{remark}

\begin{remark}
Note that $T_0 = \{(G, 0)\}$ and so
$$
    \ev_1 : \cZ^+_0 = \Filt_{\bA}^{0,0}(\bM_{\le 0}) \xrightarrow{\sim} \bM_{\le 0}.
$$
\end{remark}

Since the labels $\lambda_P, d$ are locally constant on $\Filt_{\bA}(\bM)$ and so on $\Filt_{\bA}(\bM_{\le c})$, the strata $\cZ^+_c$ are indeed each a union of connected components in the latter. It remains to check its behavior under evaluation at $1$.

\begin{proposition}
    For every $c \ge 0$, evaluation at 1 induces a closed immersion
    $$
        \ev_1 : \cZ_{c}^+ \hookrightarrow \bM_{\le c}.
    $$
    Its underlying image consists precisely of the triples of $\xi$-instability equal to $c$.
    \begin{proof}
        For each $(P,d) \in T_c$, the fiber in $\Filt^{\lambda_P, d}_{\bA}(\bM_{\le c})$ over a family $m \in \bM_{\le c}(S)$ parameterizes all compatible $P$-reductions of degree $d$, which become filtrations when equipped with the constant cocharacter $\lambda_P$. The conditions on $\nu_{P,d}$ imply by Proposition \ref{prop.can.red.dist.characterization} that the reductions which appear are canonical in each geometric fiber. The associated graded of each such reduction has the same instability by the last statement of Lemma \ref{lemma.face.compat.assoc.graded}, so its corresponding Rees filtration factors through $\Filt_{\bA}(\bM_{\le c})$, on each geometric fiber and then in families, by Proposition \ref{prop.filt.Hitchin.bdd}.

        In the case of a family over a geometric point $m \in \bM(K)$, $K$ algebraically closed, uniqueness of the canonical reduction and the fixed choice of $\lambda_P$ implies that $\ev_1$ is universally injective. Moreover, a triple of instability equal to $c$ has canonical reduction such that $(P,d)$ is an admissible pair by Proposition \ref{prop.can.parab}, so every point of $\bM$ with instability exactly $c$ lies in the image (and no other, since every other point has canonical reduction with instability $< c$ and thus is missed by $\cZ_c^+$). This gives the last claim.

        Now let $(m, \theta, t) : \Spec R \to \bM_{\le c}$ a finite-type affine family. The summand of parabolic reductions $m_P$ corresponding to $(P, d)$ is the closed subspace of sections of $\cE/P \to X_R$ which have the correct (locally constant) degree, are Higgs-compatible, and have prescribed relative position at $\infty$. 

        The morphism $\cE/P \to X_R$ is projective, with anticanonical bundle relatively ample. After twisting by a sufficiently positive pullback from $X_R$, the graphs of these sections have fixed Hilbert polynomials over $B$, and so their parameter space inside the Hilbert scheme is separated and finite type. Specifying the degree defines a locus which is closed and open within. Taking the disjoint union over the finite set $T_c$ shows that $\ev_1 : \cZ_c^+ \to \bM_{\le c}$ is representable, separated, and of finite type.
        
        Now we show that $\ev_1$ is unramified. Infinitesimally, the relative tangent space at a point $z = m_P \in \cZ^+_c(K)$ over $m = (\cE, \theta, t) \in \bM(K)$ consists of a deformation of $m_P$ over the constant family $m \in \bM(K[\epsilon]/\epsilon^2)$. Thus, it is only the parabolic reduction, viewed as a section $$\sigma : X_K \to \cE/P,$$ which can vary. The first order deformations of such a section consist of the vector space
        $$
            H^0(X_K, \sigma^\star T_{(\cE/P)/X_K}) = H^0(X_K, \cE_P \times^P (\fg/\fp)).
        $$
        The conditions of degree and relative position are locally constant, but we should restrict to those deformations which preserve $\theta \in H^0(X_K, \Ad_D(\cE))$. This linearizes to the condition that $[\theta, v] = 0$, where
        $$
            [\theta, -] : H^0(X_K, \cE_P \times^P (\fg/\fp)) \to H^0(X_K, \cE_P \times^P (\fg/\fp)(D))
        $$
        is induced from $\fg$ down to $\fg/\fp$. At a geometric generic point on $X_K$, we can trivialize the twist and diagonalize $\theta$ in a maximal torus contained in $P$, after which its bracket acts by nonzero scalar $\alpha(\theta)$ on  each root space of $\fg/\fp$. Since $[\theta, -]$ is a generically invertible map of locally free sheaves, its induced map on global sections is injective, so at every geometric point the relative tangent space
        $$
            \ker\left(H^0(X_K, \cE_P \times^P (\fg/\fp)) \to H^0(X_K, \cE_P \times^P (\fg/\fp)(D)) \right)
        $$
        vanishes. For a representable morphism locally of finite type, this implies unramified (cf.\ \cite[Tag 0B2D, Tag 02GF]{stacks-project}, applied on charts by \cite[Tag 03ZK]{stacks-project}).

        Finally, we apply the valuative criterion to show $\ev_1$ is proper. Given some family $m = (m_s, m_\eta) \in \bM_{\le c}(R)$ over DVR $R$ whose generic member lifts along $\ev_1$ to $\cZ_c^+$, the computation of geometric fibers above implies $$c = \deg_i(m_\eta) \le \deg_i(m_s) \le c.$$
        By Proposition \ref{prop.can.P.loc.const} the canonical reduction defined by the lift of $m_\eta$ extends, with constant degree $d$ (mod center). Choosing the constant cocharacter $\lambda_P$ gives a Rees filtration which lifts the family over $R$, so the existence part of the valuative criterion holds. In summary, $\ev_1$ is a proper, universally injective, unramified and representable morphism, hence a closed immersion (cf.\ \cite[Tag 05W8]{stacks-project}, applied locally on scheme charts).
    \end{proof}
\end{proposition}

We can now show directly from the definition that Definition \ref{def.theta.strata} defines a $\Theta$-stratification of $\bM$ with minimal open $\bM^{\xi}$. Compare also \cite[\S 8]{AHLHModuliForStacks}.
\begin{theorem}\label{thm.theta.strat.Hitchin}
    For arbitrary $\xi$, the data of Definition \ref{def.theta.strata} define a $\Theta$-stratification for $\bM \to \bA$.
    \begin{proof}
        We take $I = \{ c \ge 0 : T_c \neq \emptyset\}$. For any $C$, the set $I \cap [0, C]$ is finite: for each of the finitely many $P$, a bounded ball contains finitely many points of the translated projected image of $\cX_\star(P)$. By Proposition \ref{prop.can.parab} the corresponding opens $\bM_{\le c}$ cover $M$. By the characterization of the image of $\cZ_c^{+} \hookrightarrow \bM_{\le c}$ above, 
        $$
            \bM_{\le c} \backslash \ev_1(\cZ_c^{+}) = \bigcup_{c' < c} \bM_{\le c'}.
        $$
        We have shown all the criteria of Definition \ref{def.theta.stratification}.
    \end{proof}
\end{theorem}

In Chapter \ref{sect.val.exist}, this will allow us to use a powerful semistable reduction theorem for stacks equipped with a $\Theta$-stratification in order to prove existence in the valuative criterion. To match their hypotheses we will need to record the following.

\begin{lemma}\label{lem.aff.diag}
    The Hitchin fibration $\cM \to \cA$ has affine diagonal. The same holds for $\bM \to \bA$.
    \begin{proof}
        This is standard. Let $S \to \cM \times_{\cA} \cM$ a test family, corresponding to the pair of Higgs bundles $(\cE_1, \theta_1), (\cE_2, \theta_2)$ on $X_{S}$. The isomorphism scheme between the bare torsors
        $$
            \ul{\on{Isom}}_{G_{X_S}}(\cE_1, \cE_2) \to X_S
        $$
        is fppf-locally isomorphic to $G_{X_S}$, hence affine over $X_S$. The condition that $\theta_1$ and $\theta_2$ correspond along such an isomorphism defines a closed subscheme
        $$
            I := \ul{\on{Isom}}_{G_{X_S}}^{(\theta_1, \theta_2)}(\cE_1, \cE_2) \subseteq \ul{\on{Isom}}_{G_{X_S}}(\cE_1, \cE_2),
        $$
        hence the structure map $\pi :\ul{\on{Isom}}^{(\theta_1, \theta_2)}(\cE_1, \cE_2) \to X_S$ is also affine. The pullback of the diagonal
        $$
            \Delta \times_{\cM \times_{\cA} \cM} S
        $$
        is the functor of sections of $\pi$. This space of sections is affine over $S$ as a consequence of \cite[Tag 08K6]{stacks-project}: \textit{loc. cit.} implies (since $\O_{X_S}$ is flat over $S$, of finite presentation, and has proper support; and $\pi$ is affine) that the functor of $\O_{X_S}$-\textit{module} sections $\pi_\star \O_{I} \to \O_{X_S}$ is affine over $S$. The condition of being a map of algebras defines a closed subscheme of such sections, hence is also affine.

        Since $\bM \to \bA$ is the pullback of $\cM \to \cA$ along $\bA \to \cA$, its diagonal is also affine.
    \end{proof}
\end{lemma}

\begin{remark}\label{rmk.choice.of.cochars}
    The careful reader may be bothered by our arbitrary choice of $\lambda_P$ for each semistandard $P$. In some situations, specific choices connect the theory with versions that appear elsewhere. 
    \begin{enumerate}
        \item If one chooses $\lambda_P$ to be some positive integral scalar of the determinant character
        $$
            \det_P := \det(\Ad_P | \Lie(U_P)),
        $$
        (more precisely, its dual with respect to our choice of metric $(-, -)$), then up to scalar 
        $$
            (\deg(m_P)^G, \lambda_P) = \deg(m_P)(\det_P).
        $$
        This allows us to recover Behrend's formalism as in \cite{behrendsemistable}, where the ($\xi = 0$) canonical parabolic can be characterized as the largest parabolic under containment among those which maximize this measure of instability. 
        \item More generally, if $\xi^G$ is a \textit{rational} parameter in the sense that it lies in $X_\star(S)_{\Q} \subset \a_T$, then we can make a systematic choice of $\lambda$ which varies with $m$, not just with the $\xi$-canonical parabolic $P$ of $m$, as follows.
        
        In this case, for any $m$ the vector $\xi - \varrho$ is a rational vector whose primitive integral multiple defines a cocharacter $\lambda$ for which $m_{P_\lambda}$ is the $\xi$-canonical reduction. The triple $(m, \lambda, m_{P_\lambda})$ is then characterized as the unique filtered object over $m$, up to scaling on $\lambda$, which maximizes the \emph{destabilizing weight} 
        $$
        \wt_{\cE_{P_\lambda}}^\xi(\lambda) := \frac{(\deg(\cE_{P_\lambda}) - \xi, \lambda)}{\|\lambda\|}.
        $$
        Optimizing this function can be used to define a more natural $\Theta$-stratification (with the same set of $\xi$-canonical parabolics), as in \cite[Theorem 2.2.2]{HalpernLeistnerInstability} and the simplifications which follow. This gives a formulation much closer to the GIT-style formalism of \cite{kempf}. 
        
        Unfortunately, this construction works only when $\xi^G$ is a $\Q$-multiple of a genuine cocharacter, otherwise the weight will be optimized at some irrational vector, and one would need to arbitrarily choose $\lambda$ to be a rational approximation, depending on $m$.
    \end{enumerate}
\end{remark}

\chapter{\texorpdfstring{The $\xi$-semistable Higgs moduli stack is of finite type}{The xi-semistable Higgs moduli stack is of finite type}}\label{sect.finite.type}

In this chapter, we retain our running hypothesis, cf.\ \S \ref{subsect.notation.running.assumptions}, except that we assume $k$ is \emph{algebraically closed} until the final statement. In this section, we prove finite-typeness of $\ctM^\xi$ for $G/X$ a quasisplit reductive group with anisotropic center, $A_G^\spl = 1$. The key new result outside the split case is a the following lemma guaranteeing the existence of Borel reductions of bounded stability for $\cE$-torsors. In the split case it is due to Harder (cf. \cite{harderborels}, \cite{harderminkowski}). 

\begin{lemma}[{{cf.\ \cite[Prop. 6.1.6]{CLI}}}]
    Let $B \supset T$ be a maximal Borel-torus pair inside quasisplit reductive $G$. There exists a constant $c > 0$ depending only on $G$ such that every $G$-torsor $\cE$ on $X$ admits a Borel reduction $\cE_B$ whose degree lies in the corresponding shift of the relative dominant chamber, i.e.
    $$
        \mathcal{C}(c) = \{H \in \a_T : \alpha(H) \ge -c \; \forall \alpha \in \Delta_B \}.
    $$
    \begin{proof}
        Assume first that $G$ is semisimple adjoint. The condition on the degree states that the ``first numerical invariants'' $n_i := \deg(\cE_B)(\alpha_i) \ge -c$ for all simple roots $\alpha_i$. We imitate the proof in Harder \cite[Satz 2.2.6]{harderborels}. In particular, we need as base case the result for groups of split rank 1. It is immediate that the inequalities $\deg(\cE_B)(\alpha_i) > -c$ are independent of central isogeny, so we assume $G$ is adjoint. Indeed this suffices: given a central isogeny $G_0 \to G$ sending Borel $B_0 \mapsto B$ and $\cE_0$ a $G_0$-torsor, we have $\cE_0/B_0 \cong (\cE_0 \times^{G_0} G)/B$, so the Borel reductions found in the adjoint case imply the existence of corresponding Borel reductions for isogeneous $G_0$.
        
        By \cite[Prop. 6.4.4]{conradrgs}, such a $G$ is necessarily of the form $\on{Res}_{X'/X} G'$ for some $X' \to X$ finite \'etale and some $G' \to X'$ adjoint absolutely simple and quasisplit. By \cite{titsclass}, there are only two such possibilities for such a $G'$ with $G$ of split rank 1: adjoint type $\prescript{1}{}{A}_1$ (that is, split $\on{PGL}_2$) or $\prescript{2}{}{A}_2$ (that is, quasisplit $\on{PU}_3$). Moreover, in order that $G$ still have split rank 1, the \'etale cover $X'/X$ must be connected. Moreover, since a $\on{Res}_{X'/X} G'$ torsor on $X$ is simply the same data as a $G'$-torsor on $X'$, Borel reductions correspond, and the degrees are equal up to a factor of $\deg X'/X$, the problem reduces to the absolutely simple case. 

        That is, we have reduced to 2 base cases: $\on{PGL}_{2,X}$, and $\on{PU}_{3, X'/X}$ associated to a double cover. The former is done by Harder, and we recall the argument. Since $H^2(X, \G_m) = 0$ (algebraically closed constants), every $\PGL_2$-torsor lifts to a rank 2 vector bundle $\cE$, considered up to twists by line bundles. Any such can be twisted by a line bundle $\O(p)$ until it has degree $2g-1$ or $2g$ (since $k$ is algebraically closed, points of degree 1 exist on $X$, and each twist increases the degree of $\cE$ by 2). The simple form of Riemann--Roch then applies, and $\cE$ has a nonzero global section. The closure of the image of this section provides a sub-line bundle $\cL \subset \cE$ which admits a global section and therefore is of nonnegative degree. This flag is a $B$-reduction such that
        $$
            \deg(\cE_B) = \deg(\cL) - \deg(\cE/\cL) = 2\deg(\cL) - \deg(\cE) \ge -\deg(\cE) \ge -2g,
        $$
        so we have found that a bound of $c = 2g(X')$ suffices for the degree of the pullback of $\cE_B$ to $X'$, where $X' \to X$ is any cover that realizes $G$ as $\on{Res}_{X'/X} G'$ for absolutely simple $G'$. By Riemann--Hurwitz, 
        $$
            g(X') \le g(X)\deg \pi,
        $$
        but the same factor of $\deg \pi$ appears when relating the degree of a $\on{Res}_{X'/X} G'$-torsor to the degree of its underlying $G'$-torsor on $X'$, so we obtain a uniform bound of $c = 2g(X)$.

        We also need to handle the other base case of quasisplit $\on{PU}_3$. By considering the long exact sequence on cohomology associated to the norm sequence
        $$
            \begin{tikzcd}
                1 \ar[r] & U_1 \ar[r] & \Res_{X'/X} \G_{m, X'} \ar[r] & \G_{m, X} \ar[r] & 1
            \end{tikzcd},
        $$
        and using that the norm map is surjective on divisors, we can see that $H^2(X, U_1) = 0$. It follows that any $\on{PU}_3$-torsor lifts to a $U_3$-torsor, that is, a rank 3 vector bundle $\cE$ on $X'$ with an isomorphism $h : \cE \cong \sigma^\star(\cE^\vee)$. A Borel reduction corresponds to an isotropic line bundle $\cL \subset \cE$, in the sense that its complement $$\cP = \ker(\cE \to \sigma^\star \cE^{\vee} \to \sigma^\star \cL^{\vee})$$ should satisfy $\cL \subset \cP$; its degree is necessarily a triple of the form $\left(\frac{d}{2}, 0, -\frac{d}{2}\right)$ for $d = \deg(\cL)$ (relative to $X$-normalization). We seek a uniform bound $c$ such that every $\cE$ has a Borel reduction with $\frac{d}{2} > - c$. Write $g' = g(X')$, and fix a divisor $A$ on $X'$ of degree $2g' - 1$. Then $W := H^0(X', \cE(A))$ admits a quadratic form 
        $$
            q_A : W \to H^0(X', \O_{X'}(A + \sigma(A))).
        $$
        By Riemann--Roch, we have $h^0(X', \cE(A)) \ge 3g'$, and $h^0(X', \O_{X'}(A + \sigma(A))) = 3g' - 1$. The map $q_A$ is then defined by $3g' - 1$ quadratic forms on a vector space of dimension at least $3g'$. Projectivizing, the isotropic locus is the intersection of $3g' - 1$ quadric hypersurfaces in $\mathbb{P}(W)$ and so is nonempty. Let $s \in W$ be a nonzero isotropic section. The closure of its image is a line subbundle $\cF \subset \cE(A)$ which is still isotropic, and twisting down $\cL = \cF(-A) \subset \cE$ is an isotropic line subbundle. Since $\deg(\cF) \ge 0$, we have $\deg(\cL) > - 2g'$. As above, we have $g(X') \le 2g(X)$, and the same normalization factor of 2 appears in the degree pairing on $\on{PU}_3$, so the uniform bound $c = 2g(X)$ works here as well. 
        
        Now let $G$ arbitrary semisimple adjoint. Given a $G$-torsor $\cE$, define a partial order on the (nonempty) set of Borel reductions $\cE_B$ via the partial order on $\Z^n$ applied to the ``second numerical invariants'' 
        $$
        (p_1(\cE_B),\ldots, p_n(\cE_B)) :=(\deg(\cE_B)(\omega_1), \ldots, \deg(\cE_B)(\omega_n)) \in \Q^n.
        $$
        Here the $\omega_i$ are the fundamental weights, dual to the (nondivisible) coroots.
        By a simple computation in terms of highest weight representations over a splitting cover, each coordinate $p_i(\cE_B)$ has bounded denominator and is bounded above as one ranges over the (nonempty set of) Borel reductions of a fixed $\cE$, so there exists a $B$-reduction which is maximal for this partial order. Let $\cE_B$ be such a maximal reduction. Then we claim that $n_i \ge -c$ for all simple roots $\alpha_i$ for $c = 2g(X)$.          
        Assume it did not, so $\deg(\cE_B)(\alpha_i) < - c$. We reduce to a split rank 1 computation as follows: let $P_i$ be the minimal standard parabolic whose Lie algebra contains $\fg_{-\alpha_i}$, and write $M_i$ for its associated Levi quotient. Write $\cE_{P_i} = \cE_B \times^B P_i$, and $\cE_{M_i} = \cE_{P_i} \times^{P_i} P_i/R_u(P_i)$ for the pushed-forward $M_i$-torsor. Write $\ol{B}_i$ for the image of $B$ inside the Levi, which is a Borel subgroup.

        By the third isomorphism theorem for groups, we have a natural isomorphism
        $$
        P_i/B \xrightarrow{\sim} (P_i/R_u(P_i))/(B/R_u(P_i)) = M_i/\ol{B}_i.$$
        This globalizes to a corresponding isomorphism for $\cE_{P_i}$:
        $$
        \cE_{P_i}/B \xrightarrow{\sim} (\cE_{P_i}/R_u(P_i))/(B/R_u(P_i)) = \cE_{M_i}/\ol{B}_i.$$
        Thus, $B$-reductions $\cF_B$ of $\cE_{P_i}$ are in bijection with $\ol{B}_i$-reductions $\cF_{B,M}$ of $\cE_{M_i}$. From the construction, it is immediate that $n_i$ corresponds to $\deg(\cF_{B,M})(\alpha)$, where $\alpha$ is the unique nondivisible relative root of $M$. In fact, we can further descend everything to the adjoint semisimple form. 
        
        Since the result is already proven for semisimple groups of split rank 1 and $\cE_{B,M} = \cE_B \times^B B/R_u(P_i)$ by assumption fails to meet the bound in this case, we can find a different reduction $\cF_{B, M}$ which does satisfy the bound. We can lift this to a different $B$-reduction $\cF_B$ of $\cE_{P_i}$ (and therefore of $\cE$). The ``second numerical invariants" $p_j = \deg(\cE_B)(\omega_j)$ can be computed in terms of $\cE_B \times^B P^j$, where $P^j$ is the \textit{maximal} proper parabolic whose Lie algebra contains weight spaces for the negation of all simple roots \textit{except} $\alpha_j$ (cf. \cite[\S 2.1]{harderborels}). Since $\cE_B$ and $\cF_B$ agree upon pushforward to $P_i \subset P^j$, it is clear then that $\deg(\cF_B)(\omega_j) = \deg(\cE_B)(\omega_j)$ for all $j \neq i$. 

        Since $p_j(\cE_B) = p_j(\cF_B)$ for all $j \neq i$ but $n_i(\cE_B) < n_i(\cF_B)$, we show it follows that $p_i(\cE_B) < p_i(\cF_B)$. The fundamental weights form a basis and $\omega_j$ is orthogonal to $\alpha_i$, so one can write
        \begin{equation}\label{eqn.combopjni}
            \omega_{i} = \sum_{{j \neq i}} c_j\omega_j + c_i \alpha_i.
        \end{equation}
        Pairing with $\alpha_i^\vee$ gives
        $$
            1 = c_i \langle \alpha_i^\vee, \alpha_i \rangle.
        $$
        In particular, $c_i > 0$ and so, pairing the identity (\ref{eqn.combopjni}) with $\deg(\cE_B)$ gives
        $$
            p_i(\cE_B) = \sum_{j \neq i} c_j p_j(\cE_B) + c_i n_i(\cE_B)
        $$
        and likewise for $\cF_B$. All the terms are equal except that $n_i(\cE_B) < n_i(\cF_B)$, so $p_i(\cE_B) < p_i(\cF_B)$, violating the maximality assumption on $\cE_B$ and giving the desired contradiction. 

        Finally, let $G$ be arbitrary reductive. Writing $\ol{G} = G/Z_G$ for its adjoint quotient and likewise for $\ol{B}, \ol{T}$, the root systems are canonically identified. Every $G$-torsor $\cE$ pushes forward to a $\ol{G}$-torsor $\ol{\cE}$, and one has $G/B \cong \ol{G}/\ol{B}$ so $B$-reductions can be canonically identified,
        $$ \cE_B \longleftrightarrow (\cE_B \times^B \ol{B}) =: \ol{\cE}_{\ol{B}}.$$ 
        For any $\alpha \in \Phi$ one has
        $$
            \deg(\cE_B)(\alpha) = \deg(\ol{\cE}_{\ol{B}})(\alpha),
        $$
        so the conclusions lift.
    \end{proof}
\end{lemma}

In order to apply this result about bare $G$-torsors to Hitchin triples $(\cE, \theta, t)$, we need the following purely degree-theoretic characterization of when reductions of the bundle automatically induce reductions of $\theta$, and thus $P$-reductions of the Hitchin pair $(\cE, \theta)$.

\begin{lemma}[{{cf.\ \cite[Lemme 6.1.7]{CLI}}}]
    Let $c >0$ be the constant of the previous lemma. Let $P$ a standard parabolic and $\Delta_P \subset \Delta_B$ the set of simple relative roots in $N_P$ (relative to $B$). There exists a $d > 0$ such that for any $(\cE, \theta, t) \in \ctM^\xi(k)$, every $B$-reduction $\cE_B$ of $\cE$ whose degree satisfies
    $$
    \deg(\cE_B)(\alpha) \ge -c \; \;\forall \alpha \in \Delta_B 
    $$
    and \begin{equation}\label{eqn.deg.cp.def}
        \begin{cases}
            \deg(\cE_B)(\alpha) \ge d &\; \forall \alpha \in \Delta_P \\
            \deg(\cE_B)(\alpha) \le d &\; \forall \alpha \in \Delta_B \setminus \Delta_P
        \end{cases}
    \end{equation}
    necessarily admits a lift of $\theta$ along $\on{Ad}_D(\cE_P) = \cE_B \times^{B,\on{Ad}} \fp_D \hookrightarrow \on{Ad}_D(\cE)$.
    \begin{proof}
        As at the end of the previous proof, replacing $G$ with $\ol{G} = G/Z_G$, it is easy to reduce to the semisimple (adjoint) case. Now let $X' \to X$ be a splitting $\Gamma$-cover with split pullback $G' \to X'$. Likewise denote all objects arising over $X'$ by pullback with an apostrophe, e.g. $P', \fg'$, etc. We wish to reduce to the split case. Note that since degree $\deg(\cE_B)$ is compatible with pullback up to a $|\Gamma|$ scalar and the rational roots are the image of the geometric roots under the canonical projection onto the rational space, the conditions of Equation (\ref{eqn.deg.cp.def}) are equivalent after pullback, i.e., whether considered in the relative space $\ul{\a}_T(X)$ or the geometric space $\ul{\a}_T(X')$. Consider now the following diagram
        $$
        \begin{tikzcd}[ampersand replacement=\&]
            H^0(X, \on{Ad}_D(\cE_P)) \ar[r] \ar[d] \& H^0(X, \on{Ad}_D(\cE)) \ar[d] \\
            H^0(X', \on{Ad}_D(\cE_P')) \ar[r] \& H^0(X',\on{Ad}_D(\cE')) 
        \end{tikzcd}$$
        in which the top row is the $\Gamma$-invariants of the bottom row by descent, and both horizontal maps are injective because they are induced by the inclusion of a subbundle. The diagram is therefore Cartesian, i.e., if the pullback $\theta'$ lifts to some $\theta'_P$, then that $\theta'_P$ descends and so $\theta$ likewise lifts to $\theta_P$. So we reduce to the split semisimple case, with the corresponding $d$, which follows from the degree argument in the proof of \textit{loc.\ cit.}\footnote{We reference here the proof and not the statement of \cite[Lemma 6.1.7]{CLI}, because the pullback Higgs bundle need not be $\xi$-stable! But their argument constructs the lift from the degree bounds alone.}
    \end{proof}
\end{lemma}

\begin{lemma}[{{cf.\ \cite[{Lemme 6.1.6}]{CLI}}}]
Assume now that $G$ has anisotropic center. Let $B \supset T$ be an $X$-maximal Borel-torus pair. There exists a finite subset $S \subset X_\star(T)$ such that for every $(\cE, \theta, t) \in \ctM^\xi(k)$, the torsor $\cE$ admits a $B$-reduction $\cE_B$ with $\deg(\cE_B) \in S$.
\begin{proof}
By the above, there exists a $c > 0$ such that every $\cE$ admits a Borel reduction $\cE_B$ whose degree lies in
\[
    \mathcal{C}(c) = \{H \in \a_T : \alpha(H) \ge - c \;\;\forall \alpha \in \Delta_B\}.
\]
We can cover the relative space $\a_T$ with cones $\mathcal{C}_P(d)$ defined by (\ref{eqn.deg.cp.def}) above for different choices of $P$, so it suffices to find a finite set $S_P$ such that every Borel reduction whose degree lies in $\mathcal{C}(c) \cap \mathcal{C}_P(d)$ actually lies in $S_P$. Then $S = \bigcup_P S_P$ suffices. Since $A_G^\spl = 1$, we have $\a_G = 0$, so the simple relative roots span $\a_T^\star$. In particular, a subset of $\a_T$ on which every $\alpha \in \Delta_B$ is bounded above and below is compact. 

If $P = G$, then $\mathcal{C}(c) \cap \mathcal{C}_P(d)$ is already compact by definition, so we assume $P \neq G$. This set is bounded in both directions for those $\alpha \in \Delta_B \setminus \Delta_P$, but unbounded in the positive direction for $\alpha \in \Delta_P$. By the previous lemma, $\deg(\cE_B) \in \mathcal{C}(c) \cap \mathcal{C}_P(d)$ implies that $\cE_P = \cE_B \times^B P$ is compatible with $\theta$, so we obtain a reduction of the Hitchin pair $(\cE_P, \theta_P)$. This may not be compatible with the original $t$, but it certainly gives a reduction of $(\cE, \theta, w(t)) \in \ctM^{w(\xi)}$ for some $w \in W_G(\infty)$. 

We use the semistability condition of this triple to bound above $\det(\cE_B)(\alpha)$ for $\alpha \in \Delta_P$, as follows. The definition of $\xi$-semistability implies that $\deg(\cE_P) \in -w(\xi) - \gobtcone_P + \a_T^P$. Consider now the auxiliary character $\delta_P = \det(\Ad_P | \Lie R_u(P))$. It pairs positively with the simple coroots generating the obtuse cone and zero with the Levi coroots, so we have
$$
    \deg(\cE_B)(\delta_P) = \deg(\cE_P)(\delta_P) \le - \langle w(\xi), \delta_P \rangle.
$$
Taking the max over $w \in W_G(\infty)$ gives a uniform upper bound on this pairing value.

Now, $\delta_P$ is a nonnegative linear combination of the positive roots in $R_u(P)$, with strictly positive coefficients on the simple roots $\Delta_P$. Applying $\deg(\cE_B)$ to such a decomposition $\delta_P = \sum_\beta m_\beta \beta$, our upper bound for $\delta_P$ and the lower bounds for all $\beta \neq \alpha$ imply a global upper bound on $\deg(\cE_B)(\alpha)$. So
$$
    S_P = \mathcal{C}(c) \cap \mathcal{C}_P(d) \cap 
    \bigcup_{w \in W_G(\infty)} (-w(\xi) - \gobtcone_P + \a_T^P)
$$
bounds the desired finite set. 
\end{proof}
\end{lemma}

\begin{theorem}[{{cf. \cite[{Proposition 6.1.5}]{CLI}}}]\label{thm.bdd}
    If $G/X$ is quasisplit reductive with anisotropic center, then $\ctM^\xi$ and $\bM^\xi$ are of finite type over $k$.
    \begin{proof}
        Let $(B,T)$ a Borel-torus pair in $G/X$. The map $\on{Bun}_B \to \on{Bun}_T$ is a bijection on components. By \cite[Lemma 16(1)]{HeinlothRelUnif}, we have $\pi_0(\Bun_T) \cong (\ul{X}_\star(T)_\infty)_\Gamma$, and a well-defined degree map 
        $$
            (\ul{X}_\star(T)_\infty)_\Gamma \to ((\ul{X}_\star(T)_\infty)_\Gamma)_{\tf}.
        $$
        Its kernel is the torsion subgroup, necessarily finite; we write $\on{Bun}_B^\delta$ and $\on{Bun}_T^\delta$ for the finite union of components of fixed numerical degree $\delta \in ((\ul{X}_\star(T)_\infty)_\Gamma)_{\tf}$. 
        
        Each component $\on{Bun}_T^\delta$ is of finite type, as follows: When $T$ is split this moduli stack is a product of Jacobians hence clearly of finite type. For general $T$, splitting it over an \'etale $\Gamma$-cover $X' \to X$ gives a pullback map
        $$
        \on{Bun}_T \to \on{Bun}_{T_{X'}},
        $$
        affine of finite presentation by \cite[Proposition 4.4.3]{behrendthesis}. The image has constant degree in each $\G_m$ component of $T$, so maps into a finite product of Jacobians. The map $\Bun_B \to \Bun_T$ is smooth of finite type by \cite[Proposition 8.1.7]{behrendthesis}, so the inverse images $\on{Bun}_B^\delta$ are likewise of finite type. We have shown above the image of $\ctM^{\xi}$ under the finite-type map $\ctM \to \on{Bun}_G$ is covered (say, on geometric points) by the images of finitely many $\on{Bun}_B^{\delta}$. Hence $\ctM^{\xi}$ is of finite type (namely, quasicompact), and its open $\bM^\xi$ likewise.
    \end{proof}
\end{theorem}

\chapter{The valuative criterion of properness: uniqueness}\label{sect.val.unique}

\newcommand{\rel}{\mathrm{rel}}

We maintain our notation and running assumptions from \S \ref{subsect.notation.running.assumptions}. In this section, we confirm the following theorem.

\begin{theorem}[Separation of the $\xi$-stable Hitchin morphism] \label{thm.sep}
    Let $G/X$ be a quasisplit reductive group which has anisotropic center, $A_G^{\spl} = 1$. For any $\xi \in \a_T$, the restriction of $f : \ctM \to \ctA$ to the $\xi$-stable locus is separated. That is, let $R$ a complete DVR $k$-algebra with field of fractions $K$, residue field $\kappa$, and uniformizer $\pi$. Let $m = (\cE, \theta, t), m' = (\cE', \theta', t)$ be two $\xi$-stable elements of $\ctM(R)$ such that (1) $f(m) = f(m')$ and (2) there exists an isomorphism $\phi_K : m'_K \xrightarrow{\sim} m_K$ of their generic fibers. Then there exists a \textit{unique} isomorphism $\phi : m' \to m$ extending $\phi_K$.

    In particular, if $\xi$ is in general position, then $\xi$-stability and semistability coincide, so the $\xi$-semistable Hitchin fibration $f^\xi : \bM^{\xi} \to \bA$ is separated.  
\end{theorem}

Note that since $\ctA$ is a scheme, no 2-morphism compatibilities appear in the valuative criterion \cite[Tag. 0CLD]{stacks-project}. 

\begin{remark} The anisotropic central condition on $G$ is essential. For example, a non-separated Hitchin fiber for $\GL_2$ (and the $\xi$-truncation) can be found in \cite[\S 5,6]{CGeoofFL}.
\end{remark}

The proof will occupy the remainder of this chapter. The structure of the proof adapts \cite[\S 9]{CLI}, although some new arguments are needed to handle the generality of anisotropic center, rather than semisimple, groups.

\section{Preliminary reductions and constructions.}\label{subsect.sep.prelim.reduction}
\begin{proof}
When checking separatedness via the valuative criterion, we can reduce to the case that $R$ is a complete DVR and $\kappa$ is algebraically closed because $\ctA$ is locally Noetherian and the morphism $f^\xi$ is locally of finite type (\cite[Prop. 7.8]{lmb}).\footnote{All of our stacks are indeed algebraic stacks in the sense of \cite{lmb}, in particular, their diagonals are separated and quasi-compact, so this version of the valuative criterion is sufficient to check separatedness. This follows for $\ctM$ from the corresponding result for $\on{Bun}_G$, which is well-known.}

A crucial role in the proof is played by the \emph{generic trait}: let $B$ be the local ring in $X_R$ of the \textit{generic point} of the \textit{special fiber} $X_{\kappa}$. So $\Spec B \to X_R$ is a DVR supported on the generic points, with residue field $\kappa(X)$ and fraction field $F = K(X)$. We write $\eta$ for the generic point of the generic fiber and $s$ for the generic point of the special fiber.

The isomorphism $\phi_K$ consists of the data of an isomorphism $\cE'_K \to \cE_K$, also denoted $\phi_K$, such that the induced map $\on{Ad}_D(\phi_K)$ sends $\theta'_K$ to $\theta_K$. We observe that it is necessary and sufficient to extend $\phi_K$ over $\Spec B$.

\begin{lemma}
    The isomorphism $\phi_K$ extends to $X_R$ if and only if $(\phi_K)_{\eta}$ extends to $\Spec B$. These extensions are unique if they exist.
    \begin{proof}
        If $(\phi_K)_{\eta}$ extends to $\Spec B$, then it extends to an open $U \subset X_R$ which contains all the points of codimension $\le 1$ in $X_R$. Such an isomorphism is the same as a section of $\on{Isom}_{G}(\cE', \cE)$, which is an affine scheme over $X_R$ (indeed, a torsor for the affine group $\on{Aut}_G(\cE)$). By algebraic Hartogs' lemma (\cite[IV, Cor. 20.4.12]{EGA}), the section over $U$ extends uniquely to all of $X_R$. The induced isomorphism $\on{Ad}_G(\phi_K)$ sends $\theta'$ to a section of the adjoint bundle which agrees with $\theta$ on a dense, so must be equal to $\theta$ everywhere (because the base is reduced and the total space is separated over the base, so morphisms are determined by their behavior on dense opens). The same reasoning implies that the isomorphism is determined by its restriction to $X_K$, and so the extension is unique if it exists.
    \end{proof}
\end{lemma}

Moreover, we will reduce from $B$ to its completion, $\widehat{B}$, a complete DVR with residue field $\kappa(X)$ and fraction field $\widehat{F}$. 

\begin{lemma}\label{lem.complete.ext.sep.reduct}
    Write $\widehat{B}$ for the completion of $B$, and $\widehat{F} = \Frac(\widehat{B})$ for its fraction field, with point $\widehat{\eta}$. The isomorphism $(\phi_K)_\eta$ extends to $\Spec B$ if and only if its base change $(\phi_K)_{\widehat{\eta}}$ extends to $\Spec \widehat{B}$.
    \begin{proof}
        Since these isomorphisms are morphisms from their domain into the affine scheme $\on{Isom}(\cE', \cE)$, this is immediate from the fact that $B = \widehat{B} \cap F$ inside $\widehat{F}$ induces a pushout in the category of affine schemes.
    \end{proof}
\end{lemma}

Assume that $m$ is $(w,L)$-elliptic on the generic fiber and $(w,M)$-elliptic on the special fiber, $M \subseteq L$ (see Definition \ref{def.wMpair.rels}). In particular, $f(m) = (a,t) \in \ctAt{v}{L}(R)$, i.e., there is a preimage $(a_L, w(t)) \in \ctA_L$, and $[(w,L)] \ge [(w,M)]$.

By a Kostant section over $\widehat{B}$, we can find $Y \in \fl(\widehat{B})$ whose $L$-characteristic polynomial is $a|_{\Spec \widehat{B}}$. Moreover, Proposition \ref{prop.CLI.3.9.twisted} applied with ambient group $L$ and subgroup $(e, M)$ implies that we can find an element $\ol{Z} \in \fm(\kappa(X))$ and $m \in M(\kappa(X)_\infty)$ such that 
\[
    \chi_L(\ol{Z}) = a_s \quad \text{and} \quad \on{Ad}(m)\ol{Z} = w(t_a)|_{X_\kappa}.
\]
Since $\ol{Y}, \ol{Z} \in \fl(\kappa(X))$ have the same $L$-characteristic polynomial, they are $L(\kappa(X))$-conjugate by \cite[Lemme 3.8.1]{CLI} (the proof of which does not require $L$ to be defined over $k$, only $\kappa(X)$). Since $L(\widehat{B}) \to L(\kappa(X))$ is surjective by smoothness, we have shown we can replace $Y$ by a conjugate such that $\ol{Y}$ is $M(\kappa(X)_\infty)$-conjugate to $w(t_a)|_{X_\kappa}$. 

By Steinberg's theorem, we know that the torsors $\cE$ and $\cE'$ are trivial on the special fiber of $\widehat{B}$. These are smooth torsors on a Henselian DVR, so the trivializations lift to $\widehat{B}$. Under the induced trivializations of $\on{Ad}_D(\cE), \on{Ad}_D(\cE')$, the sections $\theta, \theta'$ correspond to regular elements of $\fg(\widehat{B})$ with the same characteristic polynomial as $Y \in \fl(\widehat{B})$. 

The transporter from $\theta$ to $Y$ is a smooth torsor for the centralizing torus $T_Y$, and has a point on $\kappa(X)$ by Steinberg's theorem, so by Henselian lifting it has a section on $\widehat{B}$. The same argument holds for $\cE'$. In other words, we may and do fix trivializations of $\cE, \cE'$ over $\widehat{B}$ such that $\theta$ and $\theta'$ both correspond to $Y$.

Under these trivializations, $(\phi_K)|_{\widehat{\eta}}$ corresponds to left multiplication by an element $\delta \in T_Y(\widehat{F}) \subset G(\widehat{F})$ such that 
\[
    \on{Ad}(\delta^{-1})Y = Y.
\]
The isomorphism extends across $\widehat{B}$, and thus across all of $X_R$ by the Lemmas above, if and only if $\delta \in G(\widehat{B})$. 

Since our Hitchin base consists only of polynomials which are regular semisimple, $f(m) = a$ is generically regular on both the generic and special fiber and so $Y$ is regular semisimple over $\widehat{B}$, with centralizer given by the $\widehat{B}$-torus $T_Y$. A torus over an irreducible normal scheme is isotrivial, and so splits after a finite \'etale cover of $\widehat{B}$.\footnote{Note that this cover is purely a cover of $B$, and in particular does not necessarily come from an extension of the scalar ring $R$. So we cannot absorb this extension using a formal base change argument.}

By construction $\delta \in T_Y(\widehat{F})$. We claim that there exists a unique $\lambda_{\widehat{B}} \in X_\star(T_Y)(\widehat{B})$ and $h \in T_Y(\widehat{F}) \cap G(\widehat{B}) = T_Y(\widehat{B})$ such that
\[
    \delta = \lambda_{\widehat{B}}(\pi) h,
\]
where $\pi$ is a uniformizer for $B$ and so for $\widehat{B}$. Certainly, such a decomposition exists in an unramified extension splitting the torus, since $\pi$ remains a uniformizer. But by uniqueness, $\lambda_{\widehat{B}}(\pi)$ and $h$ are necessarily Galois-invariant, so in fact both descend to be defined over $\widehat{B}$.

Since $h \in G(\widehat{B})$ centralizes $Y$ we can shift one of our trivializations, say $\cE$, by $h$ without affecting our setup. So we may assume $h = 1$. It remains to prove that $\lambda_{\widehat{B}} = 0$, so that $\delta = 1$ extends across $\widehat{B}$.

Now, the completion process may have affected ellipticity on the completed generic fiber, so we return to the special fiber, obtaining a $\kappa(X)$-rational cocharacter:
$$
    \ol{\lambda} := (\lambda_{\widehat{B}})_{\kappa(X)} \in X_\star(T_{\ol{Y}}) = X_\star(T_Y)(\kappa(X)).
$$
Since we know $\ol{Y}$ is $M$-elliptic,
$$
    X_\star(T_{\ol{Y}}) = X_\star(A_M^{\spl}).
$$
The specialization map on cocharacters for this torus over a Henselian DVR is an isomorphism, so $\lambda = 0$ if and only if $\ol{\lambda} = 0$. Of course, this global cocharacter extends constantly to a cocharacter of $A_M^{\spl}$ defined on all of $X_R$ (not necessarily equal after base change to $\lambda_{\widehat{B}}$: the latter is valued in $T_Y$ which may not contain $A_M^{\spl}$ on the generic fiber). We write $\lambda$ for this constant extension.

\section{\texorpdfstring{Proof that $\lambda = 0$: construction of parabolic reductions.}{Proof that lambda = 0: construction of parabolic reductions.}}

Assume for the sake of contradiction that $\ol{\lambda} \neq 0$ and define semistandard parabolic $P = P_G(\ol{\lambda})$ (resp. $\ol{P} = P_G(-\ol{\lambda})$). These are proper parabolics because $\lambda \neq 0$ and the split center is trivial, and they contain $M$ because $\ol{\lambda} \in X_\star(A_M^{\spl})$.

Consider the Lie algebra $\fg/X$, and decompose it as a sum of subbundles according to the weights of the action of $\G_m$ through $\ol{\lambda}$: 
\[
    \fg = \bigoplus_{i} \fg_i.
\]
Of course, each summand is a sum of relative weight spaces
\[
    \fg_i = \bigoplus_{\substack{\alpha \in X^\star(S) \\ \alpha(\lambda) = i}} \fg_{\alpha}.
\]
Write $r = \max\{i : \fg_i \neq 0\}$, so that by symmetry $-r$ is the smallest nonzero summand. We take the complement of the lowest part:
\[
    \fg_+ = \bigoplus_{-r < i \le r} \fg_i.
\]
We have $\fg = \fg_+ \oplus \fg_{-r}$. Note that $\fg_+$ is stable under $P$ and $\fg_{-r}$ is stable under $\ol{P}$. Moreover, the weights $i$ which appear for these subbundles on $X$ (or $X_{\kappa}$) coincide exactly with the weights for the corresponding action of $\lambda_{\widehat{B}} \curvearrowright \fg_{\widehat{B}}$, since the latter subbundles have locally constant rank and the two coincide on the special fiber. The subbundles themselves may well differ on the completed generic point.

Consider the composite isomorphism on the generic fiber $X_K$ (the whole generic curve, not just the generic point)
\begin{equation}
    \phi'_K := \pi^r \Ad_D(\phi_K) : \Ad_D(\cE') \to \Ad_D(\cE).\label{eqn.sep.gen.isom.cEs}
\end{equation}
We now show that since $r$ bounds how negative the weights of $\delta = \lambda(\pi)$ can be, this isomorphism extends to a morphism over all of $X_R$. Certainly, after base change to $\widehat{B}$, trivializing $\O(D)$, and applying our given trivializations, this isomorphism corresponds to the map
$$
    \pi^r \Ad(\lambda_{\widehat{B}}(\pi)) : \fg(\widehat{F}) \to \fg(\widehat{F}).
$$
Since $\lambda_{\widehat{B}}$ has the same weights as $\lambda$, this acts integrally, so this map extends across $\widehat{B}$. By the same reasoning as in Lemma \ref{lem.complete.ext.sep.reduct}, this descends from $\widehat{B}$ to $B$. So the isomorphism (\ref{eqn.sep.gen.isom.cEs}) extends to an open containing every codimension 1 point and by Hartog's lemma extends to a morphism $$\phi' : \Ad_D(\cE') \to \Ad_D(\cE)$$ on all of $X_R$. Consider its special fiber $\phi'_{\kappa}$, and write $W'$ for its kernel and $Q$ for its image as a map of coherent sheaves. First $Q$, and therefore $W'$, are vector bundles, and in fact $W' \subset \Ad_D(\cE'_{\kappa})$ is a subbundle. Since the untwisted adjoint bundles have degree 0, after twisting, the first isomorphism theorem for $\phi'$ gives\footnote{In the proof of \cite[\S 9.4]{CLI}, the degree numerics forget to account for the $D$-twist. Our following bookkeeping modifications complete their argument. Alternatively, one could twist all the vector bundles in question down by $(-D)$ and then use their numerics.}
$$
    \deg(W') + \deg(Q) = \dim(G) \deg(D).
$$

At the generic point $s$ of $X_{\kappa}$, which can be recovered before or after completion $B \to \widehat{B}$, the map $\phi'_{\kappa}$ acts by zero except on $\fg_{-r}$, i.e.,
$$
    (W')_s = \fg_+ \otimes \kappa(X), \qquad Q_s = \fg_{-r} \otimes \kappa(X).
$$

Recall that we fixed trivializations of $\cE'$ and $\cE$ on $\widehat{B}$, which specialize to trivializations over $\kappa(X)$. We can push these forward to give generic parabolic reductions for any parabolic, which extend by the valuative criterion of properness to all of $X_{\kappa}$. Consider two such,
$$
    \cE'_{P, \kappa} \quad \text{ and } \quad \qquad \cE_{\ol{P}, \kappa}.
$$
We can relate $W'$ and $Q$ to these parabolic reductions:
$$
    W' = (\cE'_{P, \kappa} \times^P \fg_+)(D)
$$
because both are subbundles (i.e., saturated) inside $\Ad_D(\cE'_{\kappa})$ and they have the same generic fiber. Similarly, 
$$
    Q \hookrightarrow (\cE_{\ol{P}, \kappa} \times^{\ol{P}} \fg_{-r})(D)
$$
as subsheaves of $\Ad_D(\cE_{\kappa})$, the codomain visibly the saturation of $Q$. In the next section, we will use these parabolic presentations to relate the degrees of these vector bundles to the degrees of the parabolic reductions. The weak inequality obtained from the inclusion of subbundles above will then contradict a strict inequality from the strict $\xi$-stability of $\cE'_{\kappa}$ and $\cE_{\kappa}$, completing the proof by contradiction.

\section{\texorpdfstring{Conclusion from $\xi$-stability.}{Conclusion from xi-stability.}}

In order to apply $\xi$-stability, we need to know that $\cE'_{P, \kappa}$ and $\cE_{\ol{P},\kappa}$ extend to reductions of the Higgs bundles $m_{\kappa}, m_{\kappa}'$, in the same relative position. By construction, these parabolic reductions are compatible with $\ol{Y} \in \fm(\kappa(X))$, which is the Higgs field under our chosen generic trivializations. The condition of being in relative position $w$ follows from the condition above that $\ol{Y}$ is $M(\kappa(X)_\infty)$-conjugate to $w(t_a)$.

So we have shown that $\cE'_{P, \kappa}$ extends to a $(w,P)$-reduction of $m'_{\kappa}$ and that $\cE_{\ol{P},\kappa}$ extends to a $(w,\ol{P})$-reduction of $m_\kappa$. Since these are both $\xi$-stable by assumption, we obtain inequalities in opposite directions from the twisted geometric cones:
\[
    w(\xi) \in - \deg(m'_{P,\kappa}) -\gobtcone_P + \ga_T^P \qquad \text{ and } \qquad  w(\xi) \in - \deg(m_{\ol{P},\kappa}) - \gobtcone_{\ol{P}} + \ga_T^{\ol{P}}.
\]

In order to contract these positivity conditions in the geometric cocharacter space $\ga_T$ into a single linear inequality, we make use of the character
$$
    \mu := \deg(P : \fg_+).
$$
Viewed as a geometric character, it has gradient in a direction positive for $P$ and negative for $\ol{P}$. Indeed, $\mu$ is a nonzero nonnegative linear combination of the geometric fundamental weights $\hat{\ul{\Delta}}_P$. In particular, then, it pairs strictly positively with $\gobtcone_P$, and likewise for $-\mu$ on $\gobtcone_{\ol{P}}$. Applying these to the stability conditions thus induces strict opposite inequalities:
\[
    \deg(\mu(\cE_{P,\kappa}')) < - \mu(w(\xi)) < \deg(\mu(\cE_{\ol{P},\kappa})).
\]
On the other hand, since $W' = (\cE'_{P, \kappa} \times^{P}_{X_\kappa} \fg_+)(D)$ we have $$\deg(W') = \deg(\mu(\cE'_{P, \kappa}))+ \dim(\fg_+)\deg(D).$$ 
On the other hand, the injection $Q \hookrightarrow (\cE_{\ol{P},\kappa} \times_{X_\kappa}^{\ol{P}, \on{Ad}} \fg_{-r})(D)$ above gives
$$
    \deg(Q) \le \deg((\cE_{\ol{P},\kappa} \times_{X_\kappa}^{\ol{P}, \on{Ad}} \fg_{-r})(D)) = - \deg(\mu(\cE_{\ol{P}, \kappa})) + \dim(\fg_{-r})\deg(D).
$$
Adding these two degree conditions and cancelling the common factor of $\dim(\fg) \deg(D)$ gives
$$
    0 \le \deg(\mu(\cE'_{P,\kappa})) - \deg(\mu(\cE_{\ol{P},\kappa})),
$$
contradicting the opposite inequality from $\xi$-stability. Thus $\ol{\lambda} = 0$, hence $\lambda_{\widehat{B}} = 0$, and so $\delta$ (i.e., $\phi_K$) extends across $\widehat{B}$ and so all of $X_R$.
\end{proof}

\chapter{The valuative criterion of properness: existence}\label{sect.val.exist}

We continue to work over the compatible locus. Fix an arbitrary element $\xi \in \a_T$. We prove the following $\xi$-semistable reduction theorem for $\bM^\xi \to \bA$.

\begin{theorem}\label{thm.val.crit.exist}
    Let $G$ a quasisplit reductive group over $X$, $\kappa$ an algebraically closed extension of $k$ and $R$ a complete DVR $k$-algebra with residue field $\kappa$ (in particular, $R$ is strictly Henselian). Let $K = \on{Frac} R$ and pick an algebraic closure $\ol{K}$. Fix $(a, t) \in \bA(R)$ and $m_{K} \in \bM^\xi(K)$ such that $f(m_{K}) = (a,t)$. There exists a finite extension $K'/K$ and $m \in \bM^{\xi}(R')$, where $R'$ is the integral closure of $R$ in $K'$, such that
    \begin{enumerate}
        \item The image of $m$ in $\bM^\xi(K')$ is isomorphic to the pullback of $m_K$.
        \item $f(m) = (a,t)$.
\end{enumerate}\end{theorem}
    \[
        \begin{tikzcd}[ampersand replacement=\&]
            \Spec K' \ar[r] \ar[d] \& \Spec K \ar[r] \ar[d] \& \bM^\xi \ar[d, "f"] \\ 
            \Spec R' \ar[rru, dashed] \ar[r] \& \Spec R \ar[r] \& \bA
        \end{tikzcd}
    \]

This theorem is exactly \cite[Th\'{e}or\`{e}m 8.1.1]{CLI} when $G$ is defined (and split) over $k$. They proceed by a conceptually simple but technically involved variant of Langton's inductive argument: since the untruncated Hitchin fibration satisfies the same existence criterion, up to finite extension they can choose an integral $R$-model, then iteratively perform local modifications to make its $\xi$-instability $\deg_{i,\xi}(m)$ decrease until it becomes $\xi$-semistable. The following theorem of Alper--Halpern-Leistner--Heinloth greatly generalizes this algorithm to apply to any reasonable stack with a $\Theta$-stratification. 

\begin{theorem}[{{\cite[Theorem 6.5]{AHLHModuliForStacks}}}]
    Let $\cX$ be a quasi-separated algebraic stack with affine diagonal which is locally of finite type over Noetherian algebraic space $S$. Assume that $\cX$ is equipped with a well-ordered $\Theta$-stratification. Then for any morphism $g : \Spec(R) \to \cX$ from a DVR whose generic fiber $g|_K : \Spec(K) \to \cX$ lands inside a stratum, there is (up to a finite extension of $R$ and $K$) a modification of $g$ which lands entirely within that stratum. 
    
    \end{theorem}

\begin{proof}[Proof of Theorem \ref{thm.val.crit.exist}]
Fix $R$ as in the statement. In the next subsection, Lemma \ref{lem.untrunc.val.exist}, we recall that the untruncated Hitchin fibration satisfies the same valuative existence. So up to finite extensions $R'/R, K'/K$, we have a lift $m' \in \bM(R')$. 

In Theorem \ref{thm.theta.strat.Hitchin}, we proved that the data constructed in Definition \ref{def.theta.strata} define a $\Theta$-stratification whose minimal stratum is the $\xi$-semistable locus $\bM^{\xi} \subseteq \bM$. Applying \cite[Theorem 6.5]{AHLHModuliForStacks} to this stratification gives, up to a further finite extension $R''/R, K''/K'$, a modified lift $m \in \bM^{\xi}(R'')$. The composition of finite extensions is finite, so this completes the proof.
\end{proof}

\begin{remark}
    We expect that the framework of \S \ref{sect.can.red} can also be used to show the hypotheses of ``$\Theta$-reductivity'' and ``S-completeness'' of \cite[Theorem 5.4]{AHLHModuliForStacks}. These imply that the stack $\bM^{\xi}$ admits an adequate moduli space. On the other hand, by \cite[Corollaire 4.11.3]{ngofl} we also know that the stabilizers over $\cA^{\heartsuit}$ are subgroups of a split torus, hence linearly reductive, so in fact $\bM^{\xi}$ should admit a good moduli space. No such considerations are necessary for our applications, so we omit a proof.
\end{remark}

\begin{remark}
    As in Remark \ref{rmk.xi.stab.parabolic.bundles}, the definition of $\xi$-stability can be extended beyond the regular semisimple locus if one works with bundles with parabolic structure at $\infty$. We expect that one can identify the $\xi$-semistable locus in such a stack as the minimal stratum in a $\Theta$-stratification as in \S \ref{sect.xi.canonical.red}, and thereby deduce from \cite{AHLHModuliForStacks} the $\xi$-semistable reduction and the existence of moduli spaces for such parabolic bundles. For $G$ semisimple split, this is done on $\Bun_G$ in \cite{heinlothschmitt}.
\end{remark}

Combining with the finiteness results of Chapter \ref{sect.finite.type} gives the following corollary.
\begin{corollary}\label{cor.univ.closed}
    Suppose $G$ has anisotropic center, $A_G^{\spl} = 1$. For every $\xi \in \a_T$, the compatible semistable Hitchin morphism $f^\xi : \bM^\xi \to \bA$ is universally closed.
    \begin{proof}
        By Theorem \ref{thm.bdd}, $f^\xi$ is of finite type over Noetherian scheme $\bA$. By \cite[Th\'eor\`eme 7.10]{lmb}, it suffices for universal closure to verify the existence of lifts up to finite extensions of $R$, for complete DVRs with algebraically closed residue field, which is Theorem \ref{thm.val.crit.exist}.
    \end{proof}
\end{corollary}

\section{The valuative criterion for \texorpdfstring{$\ctM$.}{ctM}}

It remains only to recall the standard result that the untruncated Hitchin fibration satisfies existence in the valuative criterion. See for example \cite[Lemma 6.20]{AHLHModuliForStacks}, \cite[\S 8.4]{CLI} in the split case. Here, the restriction to the compatible locus plays no role, so we prove the result for $\ctM \to \ctA$.

\begin{lemma}\label{lem.untrunc.val.exist}
Under the same hypotheses on $R$, $\ctM \to \ctA$ satisfies the existence part of the valuative criterion as in Theorem \ref{thm.val.crit.exist}.
\begin{proof}
    Writing $m_K = (\cE_K, \theta_K, t)$, we consider the restriction $(\cE_K)_{\eta}$ to the generic point $\eta$ of $X_K$. The group $G_{\eta}$ is a connected reductive group scheme over the function field $F_K$ of $X_K$. Writing $G_{\ol{\eta}}$ for the pullback to $X_{\ol{K}}$, its function field $F_{\ol{K}}$ is a function field of a curve over an algebraically closed field, so in particular is of cohomological dimension $\le 1$. 
    It follows from \cite[III.2.3, Thm 1' and Remarks]{serregalcoh} that $H^1(F_{\ol{K}}, G_{\ol{\eta}}) = 0$, and so every $G$-bundle on $X_{\ol{K}}$ is trivial at the generic point. The isomorphism $\cE_{\ol{\eta}} \cong G_{\ol{\eta}}$ descends to a finite extension $K'/K$, so it follows that, up to replacing $K'$ by a finite extension, we may assume that $\cE_K$ is trivial at $\eta$. 
    
    Fix such an isomorphism $\cE_{\eta} \cong G_{\eta}$. Under this isomorphism, the Higgs field generically corresponds to an element $Y \in \fg(F_K)$, whose characteristic polynomial $\chi(Y)$ is equal to $a|_{\eta}$. Now, inside of $X_R$, consider the local ring of the generic point of $X_{\kappa}$. This is a DVR supported on the two generic points, which we denote $B$. Its field of fractions is $F_K$. By the Kostant section, we can choose a $Z \in \fg(B)$ whose characteristic polynomial is $a|_{\Spec B}$. 

    We have now two elements $Y,Z \in \fg(F_K)$ whose characteristic polynomials are equal. Viewing the same elements in $\fg(F_{\ol{K}})$, they have the same characteristic polynomial (which is $G$-regular semisimple at $\infty$, and thus at the generic point). The transporter functor $T_G(Y,Z) = \{g \in G : g.Y = Z\}$ is nonempty (it has points over $\ol{F_{\ol{K}}}$), and is a torsor for the torus $Z_G(Y)$. In particular $T_G(Y,Z)$ is smooth scheme, hence splits over $F_{\ol{K}}^s$. Since $Z_G(Y)$ is connected reductive, the same vanishing $H^1(F_{\ol{K}}, Z_G(Y)) = 0$ as above implies that $T_G(Y,Z)$ is trivial, so there exists an element $g \in G(F_{\ol{K}})$ taking $Y$ to $Z$. Since tensor product commutes with colimits and $G$ is locally of finite presentation over $F_K$, the element $g$ lies in $F_{L}$ for some $L/K$ finite. So, up to replacing $K$ by a finite extension, we may assume $g \in G(F_K)$. 
    
    The element $g$ allows us to glue $(\cE_K, \theta)$ on $X_K$ with $(G_R \times_{X_R} \Spec B, Z)$ (a trivial bundle with Higgs field $Z$) on $\Spec B$: the latter spreads out to a bundle on some Zariski neighborhood $V$, and the isomorphism between the two pairs spreads from $\Spec F_K$ to $V \cap X_K$, after possibly shrinking $V$. The result is an extension $(\cE_U, \theta_U)$ of $(\cE_K, \theta_K)$ to an open $U \subseteq X_R$ which contains all points of codimension $\le 1$. By \cite[Thm. 6.13]{collsantors}, every $G$-torsor on $U$ extends uniquely to $X_R$. By \cite[IV, Cor. 20.4.12]{EGA}, the section $\theta_U$, defined away from a set of codimension 2, extends to $X_R$. Thus, we have found an extension $(\cE, \theta)$ of $(\cE_K, \theta_K)$ to $X_R$. Since the elements $\chi(\theta)(\infty_R)$ and $\chi(t)$ in $\fc(R)$ agree in $\fc(K)$, they are equal, so $(\cE, \theta, t) \in \ctM(R)$ is the desired extension.
\end{proof}
\end{lemma}

\chapter{\texorpdfstring{The $\xi$-semistable Higgs moduli stack is Deligne--Mumford}{The xi-semistable Higgs moduli stack is Deligne--Mumford}}\label{sect.Deligne.Mumford}

Here we must assume $\xi$ is in general position, in the sense of Definition \ref{def.gen.pos}, and restrict to the compatible locus (the latter mainly for convenience).

\begin{theorem}[{{cf.\ \cite[Th\'{e}or\`{e}me 10.1.1]{CLI}}}]\label{thm.xi.locus.dm}
    Let $G/X$ be quasisplit reductive with anisotropic center and $\xi \in \a_T$ a parameter in general position. Then $\bM^\xi$ is a Deligne--Mumford stack.
    \begin{proof}
        It suffices to show that for every point $m = (\cE, \theta, t) \in \bM^\xi(k)$, the stabilizer group is unramified. This can be checked after faithfully flat base change, so we may assume that $k$ is algebraically closed. Since the stabilizer is representable by a group scheme locally of finite type, it suffices to show its Lie algebra vanishes. Let $M$ be the Levi such that $a=f(m)$ is $M$-elliptic. The stabilizer group of $m$ is $\Aut_G(\cE, \theta)$, the centralizer of $\theta$ in $\Aut_G(\cE)$. Its tangent space consists of the subspace of $H^0(X, \Ad(\cE))$ annihilated by $[\theta, -]$. Note the lack of a twist by $D$. Let $\phi \in \Lie\Aut_G(\cE, \theta)$ be such a global section. 
        
        By Lemma \ref{lem.monodr.inv.spl.torus}.(c) under our characteristic assumptions, and \cite[Corollaire 4.11.3, Lemma 4.14.2]{ngofl}, we have
        $$
            \Lie \Aut_G(\cE, \theta) \subseteq \Lie(T^{W_a}_\infty) = \Lie(A_M^{\spl})_\infty.
        $$
        Write $Z \in \Lie(A_M^{\spl})_\infty$ for the image of $\phi$ under this inclusion. Explicitly, the infinitesimal automorphism $\phi$ gives a global section of the Lie algebra of the smoothened centralizer $I^{\mathrm{lis}}_{(\cE, \theta)}$, which pushes forward to a global section of the Lie algebra of the N\'eron model $J_a^{\flat}$ determined by its value $Z$ at $\infty$. Since $A_M^{\spl}$ is split over proper $X$ (i.e., constant), there is a canonical constant global section, also denoted $Z$, of $\Lie A_M^{\spl}$ restricting to $Z$ at infinity. Since $I^{\mathrm{lis}}_{(\cE, \theta)} \cong J_a^{\flat}$ on the regular open and Ng\^{o}'s inclusion is given by pushing forward, it follows that $Z_\eta \in \fg(F)$ is the generic coordinate of $\phi$ in any generic trivialization in which $\theta_\eta$ corresponds to a $G$-regular, $M$-elliptic element $Y \in \fm(F)$. 
        
        Now consider $L := Z_G(Z)$. By \ref{prop.goodcharcentralizers}, our characteristic hypotheses ensure $L$ is a \textit{geometric} Levi. Clearly $M$ centralizes $Z$, so $L \supseteq M \supseteq T$ is semistandard. Moreover, $L$ is defined over $X$ because $Z$ is. It remains to show that $L$ is the Levi factor of an $X$-rational parabolic, which we argue as follows. Working over a splitting cover, $Z$ is fixed by $W_L$, hence is equal to its own $W_L$-average ($p \nmid |W_L|$!). Under the canonical identification $\Lie T = X_\star(T) \otimes_{\Z} k$, we can write $Z$ as a $k$-linear combination of cocharacters. Averaging puts these characters into $X_\star(T)^{W_L} = X_\star(A_L)$, so
        $$
            Z \in \Lie(A_L)(X') \subseteq \ft(X').
        $$
        Containment in such a subbundle can be checked \'etale-locally, so in fact $Z \in \Lie(A_L)(X)$, which by Lemma \ref{lem.good.char.torus.lie.alg} below is equal to $\Lie(A_L^{\spl})$. It follows that $L = Z_G(A_L^{\spl})$, so a generic cocharacter of $A_L^{\spl}$ defines $L$ as its centralizing Levi and defines an $X$-parabolic $Q$ for $L$ by its attracting locus, \ref{prop.parab.levi.basics}.

        Denote the unipotent radical of $Q$ by $N$. As in Prop. \ref{prop.BLdescent.weighted}, the bundle $\cE$ is determined by the generic trivialization above and a family $(g_v)_{v \in X(k)} \in \prod'_v \Gr_{G,v}(k)$ of modifications satisfying conditions as in Definition \ref{def.adelic.hitchin.deg}. In particular,
        $$
            \Ad(g_v^{-1})Z = \Ad(g_v^{-1})(\phi_\eta) \in L^+_{\fg, v}(k)
        $$
        is \emph{integral} at each point because $\phi$ is a global section of $\Ad(\cE)$ without $D$-twist. 
        
        We claim that all $g_v$ have representatives in $\Gr_{L,v}(k)$. By the Iwasawa decomposition, we can assume up to disk reparameterization that each $g_v$ can be written $g_v = l_v n_v$ with $l_v \in L(F_v)$ and $n_v \in N(F_v)$. The above argument shows $\Ad(g_v^{-1})Z$ is integral at each place, but $Z$ is constant, so $Z \in L_{\fg, v}^+(k)$ also, and so
        $$
            \Ad(n_v^{-1})Z - Z \in L_{\fg, v}^+(k) \cap \mathfrak{n}(F_v) = \mathfrak{n}(\O_v).
        $$
        By the Lemma below, this implies the $n_v$ all lie in $N(\O_v)$ and $g_v \sim l_v \in \Gr_{G,v}(k)$. This gives a collection of Beauville--Laszlo data gluing together to an $L$-reduction $\cE_L$ of $\cE$. Moreover, writing $Y = \theta_{\eta}$ for an $M$-elliptic element in a generic trivialization, we have
        $$
            \Ad(l_v^{-1}) Y = \Ad(n_v) \Ad(g_v^{-1}) Y \in z_v^{-d_v} \fg(\O_v),
        $$
        but also $\Ad(l_v^{-1})Y \in \fl(F_v)$, so
        $$
            \Ad(l_v^{-1})Y \in z_v^{-d_v} \fl(\O_v).
        $$
        The same gluing arguments of \S \ref{sect.adelic.fibers}, applied to the group $L$, show that $Y$ glues to an $L$-reduction
        $$
            \theta_L \in H^0(X, \Ad_D(\cE_L)) \to H^0(X, \Ad_D(\cE)).
        $$

        In summary, we  have shown that $m = (\cE, \theta, t)$ admits a ``reduction'' to $L$, in the sense that $(\cE, \theta)$ is the extension of scalars of $m_L = (\cE_L, \theta_L)$. For each parabolic $P$ having Levi $L$, (in particular, $Q$ and $\ol{Q}$) the extension of $m_L$ to $P$ is by uniqueness equal to the $P$-reduction of $m$. In particular, this implies that all the $P$-reductions of $m$ for such $P$ have the same degree $\deg(\cE_L)$. Both $m_Q$ and $m_{\ol{Q}}$ are compatible with $t$, since $Y$ is $M(F_\infty)$-conjugate to $t_a$ and $M \subseteq L$. The rational $\xi$-semistability inequalities imply that $\xi_L^G + \deg(\cE_L)^G \in \a_L^G$ (the projections into $\a_L^G$) lies in two opposite pointed cones. In particular, the projection of $\xi$
        $$
            \xi_L^G = - \deg(\cE_L)^G
        $$
        is a lattice point in $\cX_\star(L)$, contradicting $\xi$ being in general position unless $L = G$. 
        
        In the latter case $Z \in H^0(X, \Lie(A_G))$ is central. But our characteristic hypotheses imply by Lemma \ref{lem.good.char.torus.lie.alg} below that
        $$
            H^0(X, \Lie(A_G)) = \Lie((A_G^{\spl})_\infty) = 0,
        $$
        using the anisotropic center hypothesis. By injectivity, $\phi = 0$.

\begin{lemma}\label{lem.unip.adjoint.isom}
    Fix any base scheme $S$ and any reductive group scheme $G \to S$ with fiberwise maximal torus $T$, parabolic $Q$ with unipotent radical $N$ and semistandard Levi factor $L \supset T$. Let $Z \in \fg(S)$ such that its scheme-theoretic centralizer is exactly $L$. Then the map
    \begin{align*}
        \iota : N &\to \mathfrak{n} \\
        n &\mapsto \Ad(n^{-1})Z - Z
    \end{align*}
    is an $S$-isomorphism.
    \begin{proof}
        The map is globally defined, and so can be checked to be an isomorphism on geometric fibers or after a splitting cover. This reduces to the classical statement: since the centralizer of $Z$ is exactly $L$, the adjoint $\on{ad}(Z)$ acts by an $S$-unit scalar $d\alpha(Z)$ on each weight space $\fg_\alpha \subset \mathfrak{n}$. Both the unipotent radical and its Lie algebra are affine spaces, and filtering both by height of the roots, $\iota$ preserves the filtration. The Chevalley commutation relations imply that the morphism $\iota$ is a triangular polynomial $S$-morphism with invertible diagonal entries, hence is an isomorphism. 
    \end{proof}
\end{lemma}

    \begin{lemma}\label{lem.good.char.torus.lie.alg}
        Let $A/X$ be a torus with prime-to-$p$ monodromy and split at $\infty \in X(k)$. Then
        $$
            H^0(X, \Lie A) = \Lie(A^{\spl}).
        $$
        In particular, if $A^\spl = 1$, then $H^0(X, \Lie A) = 0$.
        \begin{proof}
            Since $\infty$ defines a section $\Gal(\ol{k}/k) \to \pi_1(X, \ol{\infty})$, the image of the geometric and arithmetic monodromy coincide. In particular, a minimal splitting cover $X'$ is \textit{geometrically} connected (and projective) and so $H^0(X', \O_{X'}) = k$. By descent, we have
            $$
                H^0(X, \Lie A) = (\Lie A_{X'})^\Gamma = (\ul{X}_\star(A)(\infty) \otimes_{\Z} k)^\Gamma.
            $$
            The characteristic hypothesis implies taking $\Gamma$-invariants is exact, so commutes with the tensor, and
            $$
                H^0(X, \Lie A) = \ul{X}_\star(A)(\infty)^\Gamma \otimes_{\Z} k = \Lie(A^{\spl}).
            $$
        \end{proof}
    \end{lemma}
    
    \end{proof}
\end{theorem}

\part{Cohomological stabilization}

\chapter{\texorpdfstring{Regular centralizers, Picard stacks, and the $s$-decomposition}{Regular centralizers, Picard stacks, and the s-decomposition}}\label{sect.reg.cent.Pic}
\chaptermark{\texorpdfstring{Regular centralizers and the $s$-decomposition}{Regular centralizers and the s-decomposition}}

In this section, we recall the necessary ingredients required to define the Picard stack $\cP \to \cA$ and its action on Hitchin fibers, which is essential to analysis of the fibration (cf.\ \cite[\S 4]{ngofl}). In parallel, we adapt the constructions of \cite[\S 3, 5]{CLII} to our generality and define a variant of the degree map on the restriction to the compatible locus $\bP \to \bA$, whose kernel $\bP^{1}$ preserves $\xi$-(semi)stability for any $\xi$.

\section{Regular Centralizers.}
Fix a parabolic $P$ over $X$ with semistandard Levi $M$. Let $I^P$ denote the universal centralizer subgroup over $\fp$ (recall the latter is a vector bundle over $X$):
\[
    I^P = \{(p, Y) \in P \times_k \fp : \on{Ad}(p) Y = Y\}.
\]

\begin{proposition}
    The natural closed immersion $I^P \to I^G$ is an isomorphism on the $G$-regular locus $\fp^{\Greg}$.
    \begin{proof}
        It suffices by the fibral isomorphism criterion to check surjectivity on geometric fibers over $X$ (cf. \cite{conradrgs}, Lemma B.3.1 and the proof of Corollary 5.2.8). This reduces us to the constant geometric setting of \cite[Lemme 3.1.1]{CLII}.
    \end{proof}
\end{proposition}

One can then define the \textit{regular centralizer} group scheme $J^M$ over $\fc_M$ either by pullback along a Kostant section or by faithfully flat descent (cf. \cite[2.1.1]{ngofl}). Recall that $\fc_G$ is equipped with a $\G_m$-action making the Chevalley map equivariant. The group scheme $J^G$ over $\fc_G$ can be equipped with a $\G_m$-equivariant structure for this action (cf.\ \cite[\S 2.2.3]{ngofl}, \cite[\S 5.1]{CLII}).

\begin{proposition}[{{\cite[3.2.1]{CLII}}}]\label{prop.reg.to.centralizer}
    The natural isomorphism
    $$
    (\chi_G^\star J^G)|_{\fp^{\Greg}} \to I^G|_{\fp^{\Greg}} = I^P|_{\fp^{\Greg}}$$
    extends to a $(\G_m \times P)$-equivariant homomorphism over $\fp$:
    $$
    (\chi_G^\star J^G)|_{\fp} \to I^P
    $$
    whose composition with $i_P : I^P \to I^G$ is the homomorphism $(\chi_G^\star J^G)|_{\fp} \to I^G|_{\fp}$ from \cite[2.1.1]{ngofl}. Twisting by $D$ gives a factorization of morphisms over $[\fp_D/P]$:
    $$
        \begin{tikzcd}[ampersand replacement=\&]
        {[(\chi_G^\star J_D^G)|_{\fp_D}/P]} \ar[r] \ar[dr] \& {[I_D^P/P]} \ar[d, "i_P"] \\
        \& {i_P^\star [I_D^G/G]}.
        \end{tikzcd}
    $$
    \begin{proof}
        The complement of the regular locus has codimension at least 2 (as can be counted on geometric fibers), the domain is normal, and the codomain is affine, so we are done by Algebraic Hartog's Lemma. The equivariance is clear on the regular loci, so holds everywhere, and the twisted statement follows immediately.
    \end{proof}
\end{proposition}

\section{Degree of a regular centralizer.}

There is a natural map over the curve $X$
$$
    P \to \ul{X}_\star(P) \otimes_{\Z_X} \G_{m, X} = \ul{X}_\star(M) \otimes_{\Z_X} \G_{m,X}
$$
(the tensor product taken in abelian sheaves on $X$) which is fiber-wise dual to the evaluation map $P \times_X \ul{X}^\star(P) \to \G_{m,X}$. The codomain is a torus; in fact,
\begin{proposition}
    The torus $T_M := \ul{X}_\star(M) \otimes_{\Z_X} \G_{m,X}$ satisfies $\ul{X}^\star(T_M) \cong \ul{X}^\star(M)$ canonically. Equivalently, $T_M$ is canonically isomorphic to the maximal torus quotient of $M$
    $$
        T_M \cong M^{\ab}.
    $$
    In particular, $X^\star(T_M) = X^\star((M^{\ab})^\spl)$, the maximal split torus quotient.
    \begin{proof}
        We compute by internal tensor-Hom adjunction
        \begin{align*}
            \ul{X}^\star(T_M) &= \sHom_{X}(\ul{X}_\star(M) \otimes_{\Z_X} \G_{m,X}, \G_{m,X}) \\
            &\cong \sHom_X(\ul{X}_\star(M), \sHom_{X}(\G_{m,X}, \G_{m,X})) \\
            &= \sHom_X(\ul{X}_\star(M), \Z_X) \\
            &= \ul{X}^\star(M). 
        \end{align*}
    But every character of $M$ factors through $M/D(M) = M^{\ab}$.
    \end{proof}
\end{proposition}

Pulling the coevaluation map back along $\fp \to X$ and restricting to $I^P$ gives a map 
$$
    I^P \to \ul{X}_\star(P) \otimes_{\Z_X} \G_{m, \fp}.
$$
At the level of \emph{regular} centralizers, we can descend this to the Chevalley space:

\begin{proposition}[{{\cite[Proposition 3.3.1]{CLII}}}]\label{prop.deg.map.J}
    There exists a unique morphism of group schemes over $\fc_M$ (itself over $X$)
    \[
        (\chi^M_G)^\star J^G \to \ul{X}_\star(M) \otimes_{\Z_X} \G_{m, \fc_M}
    \]
    whose pullback along $\chi_P : \fp \to \fc_M$ is the natural composition
    \[
        \begin{tikzcd}[ampersand replacement=\&]
            (\chi_P)^\star (\chi_G^M)^\star J^G = (\chi_G^\star J^G)|_{\fp} \ar[r] \ar[dr] \& I^P \ar[d] \\
            \& \ul{X}_\star(M) \times_{\Z_X} \G_{m, \fp}
        \end{tikzcd}.
    \]
    The correseponding results hold after twisting by $D$, that is, for $J^G_D$, $\fc_{M,D}$, etc.
    \begin{proof}
        The argument is by faithfully flat descent along $\chi_P$. The required gluing condition can be checked geometric fiberwise over $X$, where the proof reduces to the one given in \cite{CLII}. All maps in the picture are $\G_m$-equivariant, 
    \end{proof}
\end{proposition}

\begin{remark}
Let us emphasize some intuition about this construction when $M = T$. We have shown that there is a $\ft$-morphism $(\chi_G^T)^\star J^G \to \ul{X}_\star(T) \otimes_{\Z_X} \G_{m, \ft} = T \times_X \ft$ which, given $Y \in \fb$, extends the natural projection 
$$
I_Y^B = Z_B(Y) \subseteq B \to T.
$$
The above map says that this projection, extended to the regular centralizer over $\fb$, is well-defined at the level of characteristic polynomials. That is, given any $Y \in \fb^{\Greg}$ with that characteristic polynomial and any $b$ centralizing $Y$ (which is necessarily in $B$, see \cite[Lemme 2.4.3]{ngofl}), the image of $b$ in the torus is \textit{independent of $Y$} under the canonical isomorphisms $Z_G(Y) \cong Z_G(Y')$.
\end{remark}

\section{Description Galoisienne.}\label{subsect.descr.galois}

Making use of the degree map, we can characterize $J^G$ in terms of data along $\chi^T_G$. Recall that $W/X$ is the finite \'etale Weyl group scheme. The morphism constructed above when $M = T$,
$$
    (\chi_G^T)^\star J^G \to X_\star(T) \otimes_{\Z_X} \G_{m, \ft} = T \times_X \ft,
$$
is $W$-equivariant for the diagonal action on $T \times_X \ft$, (this can be checked fiberwise over $X$), and so the adjoint morphism
$$
    J^G \to (\chi_G^T)_\star (T \times_X \ft),
$$ 
where $(\chi_G^T)_\star$ denotes the Weil restriction, factors through the $W$-fixed points
$$
    ((\chi_G^T)_\star (T \times_X \ft))^W.
$$
It is clear that this is a smooth group scheme, isomorphic to a torus over the regular locus, and in particular commutative everywhere. It follows from the definitions that there is a natural map 
$$
    X_\star(T) \otimes_{\Z_X} (\chi_G^T)_\star \G_{m, \ft}  \to (\chi_G^T)_\star (X_\star(T) \otimes_{\Z_X} \G_{m, \ft})
$$
which can be checked to be an isomorphism \'etale-locally, so we have an alternative characterization
$$
    ((\chi^T_G)_\star(T \times_X \ft))^W = (X_\star(T) \otimes_{\Z_X} (\chi_G^T)_\star \G_{m, \ft})^W.
$$

\begin{proposition}[{{\cite[Corollaire 2.4.8]{ngofl}}}]
    The morphism $J^G \to ((\chi^T_G)_\star(T \times_X \ft))^W$ is an isomorphism over the connected component containing the identity. In particular, we have a infinitesimal characterization
    $$
        \on{Lie}(J^G) = (X_\star(T) \otimes_{\Z} (\chi_G^T)_\star \O_{\ft})^W.
    $$
\end{proposition}

\section{\texorpdfstring{$D$-twisting and the degree map on the Picard stack.}{D-twisting and the degree map on the Picard stack.}}\label{subsect.Pic.deg}

We now define the Picard stack associated to our Hitchin fibration and a degree map on its elements. 

\begin{definition}
    The \emph{Picard stack} for $G \to X$ is the relative stack of $J_D$-bundles on $X$ over $\cA$:
    $$\cP := \on{Bun}_{J_D/\cA}(X) = \Map_{X \times_k \cA/\cA}(X \times_k \cA, B_{\cA} J_D),$$
    a stack over $\cA$. Explicitly, it has functor of points associating to $S \to \cA$ the groupoid of $(X \times \cA)$-morphisms
    $$
        X \times_k S \to B_{\cA}(J_D) 
    $$
    Since $J$ is commutative, the contracted product indeed equips $\cP$ with the structure of a Picard stack over $\cA$ in the sense of \cite[Exp.\ XVIII, Def. 1.4.2]{SGAIV}.
    
    Specializing, for any $k$-scheme $R$ and $a \in \cA(R)$ we have by formal base change relations (cf. \cite[Lemma 1.1.1]{HalpernLeistnerInstability}) a Picard $R$-stack
    $$
        \cP_a := a^\star \cP = \Map_{X_R/R}(X_R, B_{X_R} J_a),
    $$
    where $J_a := a^\star J_D^G$. We write $\bP \to \bA$ for the restriction of $\cP$ to the compatible locus of Definition \ref{def.comp.loc}. 
\end{definition}

Now we proceed to define the degree map. The following sheaf, constructed as in \cite[\S 6]{CLII}, is its natural target -- as in Proposition \ref{prop.dual.lattice.global.invs.coinvs.mod.tors}, the dual is a space of coinvariants mod torsion, or more canonically, $H_0$ mod torsion.

\begin{definition}\label{def.cochar.sheaf}
Denote by $\cX_\star$ the constructible sheaf on $\bA$ whose restriction to $\bAe{M}$ is the constant sheaf of value $H_0(X, \ul{X}_\star(M))_{\tf}$ and such that the specialization maps corresponding to $M \subset L$ are those induced from the natural map $$\ul{X}_\star(M) \to \ul{X}_\star(L).$$ Abusively, we write also $\cX_\star(M) = H_0(X, \ul{X}_\star(M))_{\tf}$, the dual of the relative character lattice, as in Proposition \ref{prop.dual.lattice.global.invs.coinvs.mod.tors}.
\end{definition}

Given $a \in \bA_{M}(R)$ (compatible) with a preimage $a_M \in \bA_M(R)$, the $D$-twisted map over $\fc_{M,D}$ constructed in Prop. \ref{prop.deg.map.J} pulls back along $a_M$ to a morphism over $X_R$
\begin{align*}
    J_{a} = a^\star J^G_D = a_M^\star (\chi^M_G)^\star J^G_D \longrightarrow \; & a_M^\star\left[(\ul{X}_\star(M) \otimes_{\Z_X} \G_{m,X}) \times_X \fc_{M,D}\right] \\
        &= \ul{X}_\star(M) \otimes_{\Z_X} \G_{m,X} \\
        &\cong M^{\ab}_R.
\end{align*}
Applying $B_{X_R}(-)$ gives a map
$$
    B_{X_R} J_{a} \to B_{X_R}(M^{\ab}_R),
$$
and taking mapping stacks from $X_R$ this becomes
$$
    \bP_a \to \Map_{X_R/R}(X_R, B_{X_R} M^{\ab}_R). 
$$

Taking degrees of $M^{\ab}$-torsors then defines the degree, which is naturally valued in the dual space $\cX_\star(M)$ of the global characters $X^\star(Z)$. 

\begin{lemma}\label{lem.deg.B.torus}
    Given a torus $Z$ over proper $X$, there is a homomorphism
    $$
        \deg : (B_X Z)(X) \to \Hom_{\Z}(X^\star(Z), \Z) \cong H_0(X, \ul{X}_\star(Z))_{\tf}
    $$
    natural in homomorphisms of tori and compatible with \'etale pullback.
    \begin{proof}
        Given a $Z$-torsor $\cE$ on $X$ and an $X$-character $\chi : Z \to \G_m$, one has a line bundle $\cE \times^{T, \chi} \G_m$. This defines the first map, whose compatibilities are left to the reader. The isomorphism is Proposition \ref{prop.dual.lattice.global.invs.coinvs.mod.tors}.
    \end{proof}
\end{lemma}

\begin{remark}\label{rmk.deg.reg.centr.def.integrally}
    We want to emphasize this subtlety which has first arisen here: the degree of a regular centralizer element or torsor is naturally dual to rational characters, and so lives in the dual $\cX_\star(M)$ of $X^\star(M)$ -- this is \emph{not} the integral space $X_\star(M)$, but rather Galois coinvariants for the geometric space of cocharacters. In the geometric part of this work, it was more convenient to tensor to $\R$, after which the distinction between $\ul{X}_\star(T)^\Gamma$ and $\ul{X}_\star(T)_{\Gamma, \tf}$ vanishes. Here and in the forthcoming arithmetic calculations, it is necessary to work integrally.
\end{remark}

Since degree is locally constant in families, we can compose with this lemma point-by-point over $R$, giving a composition $\deg_M : \bP_a \to \Map_{X_R/R}(X_R, B_{X_R} J_a) \to H_0(X, \ul{X}_\star(M))_{\tf, R} = \cX_\star(M)_R$ for each $R$-point $a$ which lands inside the image of $\bA_M$ inside $\bA_G$.

In particular, applying the map over the universal point inside $\bA_{M}$ gives
\begin{equation}
    \deg_M : \bP|_{\bA_M} \to \cX_\star(M) \times \bA_{M}. \label{eqn.deg.M.Pic}
\end{equation}
The degree maps over each $M$-elliptic locus will be assembled into one map below, using the following compatibility between $\deg_M$ and $\deg_L$.

\begin{proposition}[{{\cite[Prop.\ 5.3.1]{CLII}}}]\label{prop.deg.Pic.compat}
    Let $M \subseteq L$ be an inclusion of semistandard Levis, and let $a \in \bA_M \subseteq \bA_L$. Then we have a commutative triangle
    $$
        \begin{tikzcd}[ampersand replacement=\&]
            \bP_a \ar[r, "\deg_M" above] \ar[dr, "\deg_L" below left] \& \cX_\star(M) \ar[d] \\
            \& \cX_\star(L)
        \end{tikzcd}
    $$
    where the vertical map is the specialization map of \ref{def.cochar.sheaf}.
\begin{proof}
    The proof is exactly as in \textit{loc.\ cit.\ } By Prop. \ref{prop.Levis.rel.standard.containment} we can choose parabolics $P$ for $L$ and $Q$ for $M$ such that $Q \subseteq P$. The compatibilities between the various Chevalley maps imply that the homomorphisms $\chi_P^\star (\chi^L_G)^\star J^G|_{\fp^{\reg}} \to I^P|_{\fp^{\reg}} \to \ul{X}_\star(L) \otimes_{\Z_X} \G_{m, \fc_L}$ of Prop. \ref{prop.deg.map.J} are compatible with that for $M,Q$. Applying $B$ preserves that commutativity. Likewise, naturality of the degree evaluation of Lem. \ref{lem.deg.B.torus} and compatibility of the maximal torus quotients for $M \subset L$ implies that the maps $(B_X L^{\ab})(X) \to \cX_\star(L)$ and $(B_X M^{\ab})(X) \to \cX_\star(M)$ commute with the natural projections. 
\end{proof}
\end{proposition}

\begin{proposition}[{{\cite[Prop. 6.4.1]{CLII}}}] \label{prop.deg.Pic} The degree maps of (\ref{eqn.deg.M.Pic}) glue to a unique map 
$$
    \deg : \bP \to \cX_\star
$$
over $\bA$.
\begin{proof}
    Uniqueness is clear. For existence, note that the universal $\deg_M$ from (\ref{eqn.deg.M.Pic}) pulls back along $\chi_G^M$ to give
    $$
        (\chi_G^M)^\star \bP \to \cX_\star(M)_{\bA_M},
    $$
    or by adjunction
    $$
        \bP \to (\chi_G^M)_\star (\cX_\star(M)_{\bA_M}).
    $$
    Assembling these over $M \supset T$ gives a map
    $$
        \deg : \bP \to \bigoplus_{M} (\chi_G^M)_\star (\cX_\star(M)_{\bA_M}).
    $$
    The sheaf $\cX_\star$ can be identified with the subsheaf of this direct sum consisting of compatible elements, i.e., the sheaf whose stalk at a point $a \in \bAe{M}$ consists of elements $\bigoplus_{L \supseteq M} p_L$ such that $p_L \mapsto p_{L'}$ under $\cX_\star(L) \to \cX_\star(L')$. The compatibility of Proposition \ref{prop.deg.Pic.compat} implies that $\deg$ lands within this subsheaf.
    
\end{proof}
\end{proposition}

\begin{definition}\label{def.P1}
    Denote by $\bP^{1} = \ker_{\bA}(\deg)$ the fiberwise kernel of the degree map. That is, over $a \in \bAe{M}$ it is the substack of $BJ_a$-torsors whose degree when paired with each character in $X^\star(M)$ is \emph{zero}. 
\end{definition}

\section{Action on Hitchin fibers.}\label{subsect.pic.action.hitchin.fibers}

There is an action of $\cP$ on $\cM$ over $\cA$, which restricts to an action of $\bP$ on $\bM$ over $\bA$.

\begin{proposition}[{{\cite[\S4.3]{ngofl}, \cite[\S5.5]{CLII}}}]
    Let $a \in \cA$ and $m \in \cM_a$. Then there is a map $J_a \to m^\star[I_D^G/G] = \Aut(m)$, where we have interpreted $m$ as a map $X \to [\fg_D/G]$.
    \begin{proof}
        The map $\chi^\star J \to I$ is $\G_m \times G$-equivariant, so descends to a homomorphism $[(\chi^\star_{G,D} J_D^G)/G] \to [I_D^G/G]$ over $[\fg_D/G]$. Pulling back along $m : X \to [\fg_D/G]$ gives a map
        $$
            J_a = m^\star [(\chi_{G,D}^\star J^G_D)/G] \to m^\star [I_D^G/G].
        $$
        The latter is the group of automorphisms of $m = (\cE, \theta)$: the automorphisms of the torsor $\cE$ are given by $[\cE]^\star [G/G]$, where $[\cE] : X \to BG$ is the classifying map. The condition of commuting with $\theta$ restricts to $m^\star [I_D^G/G]$. 
    \end{proof}
\end{proposition}

So we have a map $J_a \to \Aut(m)$. This induces an action of $\cP_a$ on $\cM_a$ by associated product of torsors:
$$
    \cF \in BJ_a : m \mapsto m' := (\cF \times^{J_a} m).
$$
The new Higgs bundle $m'$ lies in the same Hitchin fiber $m$ because they are isomorphic locally over the curve, so their characteristic polynomials agree everywhere.

\section{\texorpdfstring{The $\bP$-action and $\xi$-stability.}{The P-action and xi-stability.}}

Let $(m,t) \in \bM$ and $P$ a semistandard parabolic subgroup. A $P$-reduction of $(m,t)$ is exactly the data of a factorization
$$
    \begin{tikzcd}[ampersand replacement=\&]
        X \ar[r, "m_P"] \ar[dr, "m" below left] \& {[\fp_{D}/P]} \ar[d, "i_P"] \\
        \; \& {[\fg_{D}/G]}
    \end{tikzcd}
$$
such that the composition $a_P : X \xrightarrow{m_P} [\fp_{G,D}/P] \to \fc_P = \fc_M$ satisfies $a_P(\infty) = \chi_M^T(t)$, i.e., lands inside $\bA_M$.  

\begin{proposition}[{{cf.\ \cite[5.5.1]{CLII}}}] \label{prop.Pic.action.shift.convex}
Let $(a,t) \in \bAe{M}$ and $(m,t) \in \bM_a$. Then for any $p \in \bP$, we have
\begin{equation*}
    C_{p\cdot m} = C_m - \deg_M(p). 
\end{equation*}
In particular, the degree 0 substack $\bP^{1}$ preserves the $\xi$-semistable substack $\bM^\xi$ for any $\xi$.
    \begin{proof}
        Recall that $J_a := a^\star J_D^G$ along $a: C \to \fc_{G,D}$. For any parabolic $P \supset M$ let $m_P : X \to [\fp_D/P]$ the unique $P$-reduction (that is, $(e,P)$-reduction) of $(m,t)$. 

        From \ref{prop.reg.to.centralizer}, we have a factorization
        \begin{equation}
            {[(\chi_G^\star J_D^G)|_{\fp_D}/P]} \to {[I_D^P/P]} \to {i_P^\star [I_D^G/G]}. \label{eqn.reg.to.centralizer}
        \end{equation}
        Note that ${[(\chi_G^\star J_D^G)|_{\fp_D}/P]} \cong i_P^\star [\chi_{G,D}]^\star J_D^G$ by the commutativity of 
        $$
            \begin{tikzcd}[ampersand replacement=\&]
                \fp_D \ar[r, "i_P"] \ar[d] \& \fg_D \ar[d] \ar[dd, "\chi_G", bend left=55] \\
                {[\fp_D/P]} \ar[r, "i_P"] \ar[d, "{[\chi_{P,D}]}" left] \& {[\fg_D/G]} \ar[d, "{[\chi_{G,D}]}" left] \\
                \fc_{P,D} \ar[r] \& \fc_{G,D}
            \end{tikzcd}.
        $$
        Pulling back (\ref{eqn.reg.to.centralizer}) along $m_P$ gives a factorization 
        $$
            a^\star J_D^G = a_P^\star [\chi_{G,D}]^\star J_D^G = m_P^\star i_P^\star [\chi_{G,D}]^\star J_D^G \to m_P^\star [I_D^P/P] \to m_P^\star i_P^\star [I_D^G/G] = \Aut(m).
        $$
        The middle term $m_P^\star[I^P_D/P]$ is canonically $\Aut(m_P)$. We have factored $J_a \to \Aut(m)$ through $J_a \to \Aut(m_P)$, so the first map defines an action of $J_a$-torsors on $P$-reductions of Hitchin triples compatible with the action on the triples. From this compatibility, the twisted $P$-Hitchin triple $p\cdot m_P := p \times^{J_a} m_P$ defines a $P$-reduction of $p \cdot m$. Indeed, it is the unique $(e,P)$-reduction because of the compatibility of the $\bP$-action with the characteristic polynomial map: that is, $p \times^{J_a} m_P$ looks locally like $m_P$, hence has characteristic polynomial which is equal to $a_P$ over a cover, and so $a_P(\infty) = \chi_M^T(t)$ is preserved. 
        
        Finally, we can compute for any $\mu \in X^\star(P)$
        $$
            \deg_P(p \times^{J_a} m_P)(\mu) = \deg(m_P)(\mu) + \deg_M(p)(\mu|_M).
        $$
        This is equivalent to commutativity of the following diagram, in which all mapping stacks are taken in the sense above:
        $$
            \begin{tikzcd}[ampersand replacement=\&]
                \Map(X, B_X J_a) \times \Map(X, [\fp_D/P]) \ar[rr, bend left = 15, "{(\deg_M, \deg_P)}"] \ar[r]\ar[d] \& \Map(X, B_X M^{\ab})^2 \ar[r] \ar[d] \& \cX_\star(M)^2 \ar[d, "\sum"] \\
                \Map(X, [\fp_D/P])  \ar[rr, bend right=15, "\deg_P" below] \ar[r] \& \Map(X, B_X M^{\ab}) \ar[r] \& \cX_\star(M)
            \end{tikzcd}.
        $$
        Here $M^{\ab} = M/D(M)$, the horizontal maps are the compositions used to define the degree maps on Hitchin triples and $J_a$-torsors, and we have used the fact that the degree of an $M$-torsor factors through the associated $M^{\ab}$-torsor.

        The square on the left commutes before taking $\Map(X, -)$ by functoriality of the associated product. But now the middle vertical map is the associated product, which \textit{is} the group operation on the Picard stack, so the commutativity of the right square is just the statement that the map of Lemma \ref{lem.deg.B.torus} is a homomorphism. 

        Since the cones defining $C_m$ are centered at the negative of these degrees, their intersection is translated by $-\deg_M(p)$, as desired.
    \end{proof}
\end{proposition}

\section{\texorpdfstring{The morphism $\bX_\star \to \pi_0(\ctP)$.}{The morphism cocharacters to Picard components.}}

The starting point for Ng\^{o}'s endoscopic decomposition of the Hitchin fibration is an analysis of the fiberwise component group of $\cP$, or more explicitly, its pullback $\ctP$ to $\ctA$. 

\begin{proposition}[{{\cite[6.2]{ngoendo}}}]
    There exists a constructible sheaf of groups $\pi_0(\cP)$ over $\cA$ whose fiber over $a \in \cA$ is $\pi_0(\cP_a)$.
\end{proposition}
A crucial (and defining) property of $\cA^{G\el}$ in Ng\^{o}'s work is that it is exactly the locus on which $\pi_0(\cP_a)$ is finite, so the $\cP_a$-actions on cohomology of Hitchin fibers factor though finite groups. The restriction to $\bP^{1}$ will allow this to extend to $\bA$. Recall the definition of the intrinsic monodromy invariant $W_{\wtilde{a}} \subset \W_G \rtimes \Out(\G)$ from Definition \ref{def.monodr.inv.Out}. The essential result on $\pi_0(\cP)$ is the following. 

\begin{definition}\label{def.Wa.coinvar.sheaf}
    Define $\bX_\star := \ul{X}_\star(T)(\infty)_{\bA}$, the constant group scheme over $\bA$. Define $\cY$ to be the quotient of $\bX_\star$ over $\bA$ whose fiber over $\wtilde{a} = (a,t)$ is $(\bX_\star)_{W_{\wtilde{a}}}$.
\end{definition}

\begin{proposition}[{{\cite[\S 5.5]{ngofl}, \cite[6.3-6.7]{ngoendo}}}]\label{prop.ngo.cpt.map}
     There exists a canonical commutative diagram of fiberwise-surjective maps
    $$
        \begin{tikzcd}[ampersand replacement=\&]
        \bX_\star \ar[r] \ar[dr]\& \cY \ar[d] \\
        \& \pi_0(\ctP)
        \end{tikzcd}.
    $$
    \begin{proof}[Proof sketch]
        We recall the construction of the map $(\bX_\star)_{W_{\wtilde{a}}} \to \pi_0(\ctP_a)$ from various work of Ng\^{o}. Assume $a$ is a geometric point with residue field $K$. Define $\ctP_a'$ to be the stack of $J_a^0$-torsors, where $J_a^0$ is the connected component of the identity. Generically, $a$ is regular semisimple, so the cokernel $\pi_0(J_a) = \coker(J_a^0 \hookrightarrow J_a)$ is supported on a proper closed subset. We therefore have an exact sequence
        $$
            H^0(X_K, \pi_0(J_a)) \to H^1(X_K, J_a^0) \to H^1(X_K, J_a) \to H^1(X_K, \pi_0(J_a)) = 0.
        $$
        The middle $H^1$ terms are the groups of isomorphism classes of $J_a^0$-torsors (resp. $J_a$-torsors), so the map $\ctP_a' \to \ctP_a$ is essentially surjective on geometric points, hence surjective on $\pi_0$. 
        
        We then make a second reduction to the N\'{e}ron model. Explicitly, write $U_a \subset X_K$ for the open on which $a$ is regular semisimple and $V_a$ its preimage in the cameral cover $\wtilde{X}_a$, so that $J_a^0|_{U_a} = J_a|_{U_a}$ is a torus, canonically isomorphic to $T$ on a complete formal neighborhood $\widehat{\O}_{\infty_K}$ of $\infty_K$. The theory of connected N\'{e}ron models implies that any torus $Z$ on $U_a$, i.e., $J_a^0|_{U_a}$, extends to a universal connected group scheme $N^0(Z)$ over $X_K$.
        We write $N_a^0 := N^0(J_a|_{U_a})$. In the notation of \S \ref{subsect.neron.models} $N_a^0 = (J^\flat_a)^0$. The N\'{e}ron mapping property gives a map
        $$
            J_a^0 \to N_a^0
        $$
        extending the identity on $U_a$. Again, this is a generically an isomorphism, so the long exact sequence in cohomology implies a surjection $H^1(X_K, J_a^0) \to H^1(X_K, N_a^0)$, inducing a corresponding surjection on connected components of Picard stacks. But, unlike above, the defect to injectivity $H^0(X_K, N_a^0/J_a^0)$ is an affine part coming from the normalization, in particular, is connected. It therefore maps into $(\ctP_a')^0$, so the map on connected components is also injective.

        Finally, Ng\^{o} defines a morphism $(\bX_\star)_{W_{\wtilde{a}}} \to \pi_0(\Bun_{N_a^0})$ using the formalism of Kottwitz invariants. We consider additive functors $F : \mathsf{Tori}_{U_a} \to \Ab$ which are
        \begin{enumerate}
            \item right-exact and
            \item satisfy $F(Z) \cong \Z$ naturally for induced tori $Z$.
        \end{enumerate}
        It is easy to see that $$A : Z \mapsto [\ul{X}_\star(Z)(\infty_K)]_{\pi_1(U_a, \infty_K)}$$
        satisfies these conditions. Kottwitz's lemma says that these properties determine the functor uniquely, in the following sense. The proof proceeds by writing an arbitrary torus as a quotient of induced tori.

        \begin{lemma}[{\cite[2.2]{kottwitzisocrystals}, \cite[6.5]{ngoendo}}]\label{lem.kottwitz.lemma}
            Let $B : \mathsf{Tori}_{U_a} \to \Ab$ be any additive functor satisfying (1) and (2) above. Then $A \cong B$, with natural isomorphisms
            $$
            \Hom(\ul{X}_\star(-)_{\infty_K}, B) \cong \Hom(A, B)
            $$
            and 
            $$
                \Hom(\ul{X}_\star(-)_{\infty_K}, B) \cong B(\G_{m, U_a}).
            $$
            Moreover, $B(\G_{m, U_a})$ is free of rank 1, and $\phi : A \to B$ is an isomorphism if and only if it corresponds to a generator of this free group. 
        \end{lemma}
        
        Ng\^{o} shows that $B : Z \mapsto \pi_0(\Bun_{N^0(Z)})$ satisfies the necessary conditions, so a choice of generators for $\pi_0(\Bun_{N^0(\G_{m, U_a})}) = \pi_0(\Pic_X)$ (which we normalize by the positive generator $\O_X(\infty)$, cf.\ Notation \ref{not.BL.conventions}) induces a morphism $$\ul{X}_\star(J_a)(\infty_K) \to \pi_0(\Bun_{N^0_a}) = \pi_0(\ctP_a') \to \pi_0(\ctP_a)$$ which factors through the $\pi_1(U_a, \infty)$-coinvariants. 
        
        Now, the restricted cameral map $V_a \to U_a$ is the pullback of $\ft^{\reg}_D \to \fc^{\reg}_D$, and the isomorphism \ref{prop.reg.to.centralizer} implies that 
        $$J_a \times_{U_a} V_a \cong J_{D}^G \times_{\fc^{\reg}_D} V_a \cong T \times_X V_a.$$

        This isomorphism identifies the $\pi_1(U_a, \infty)$-action on $J_a$ with the monodromy action defined by $\mu_{\wtilde{a}}$ (Definition \ref{def.monodr.inv.Out}), cf.\ Remark \ref{rmk.monodr.inv.reg.centr}. Our choice of point $\wtilde{\infty}$ in the cameral cover (via $t$) therefore determines an isomorphism $\ul{X}_\star(J_a)|_\infty \cong \ul{X}_\star(T)|_\infty$ such that $\pi_1(U_a, \infty)$ acts on the composition through $W_a$. So Ng\^{o}'s map is the same data as a map $(\bX_\star)_{W_a} \to \pi_0(\ctP_a)$, as desired.
    \end{proof}
\end{proposition}

The morphism $\bX_\star \to \pi_0(\ctP)$ constructed above is not very explicit, passing through Kottwitz's yoga of induced tori and also through the surjection $\pi_0(\ctP') \to \pi_0(\ctP)$, whose kernel is inexplicit in general. 

Drinfeld proposed a functorial Beauville--Laszlo construction of a refinement $\bX_\star \to \ctP'$ (or to $\ctP$, more weakly) which induces Ng\^{o}'s map above upon passing to components. Functoriality of this construction also justifies why the maps constructed above geometric-fiber-by-geometric-fiber interpolate over $\bA$.

\begin{proposition}[{{\cite[\S 7.3]{CLII}, \cite[6.8]{ngoendo}}}]\label{prop.inv.map.BL}
    There exists an $\ctA$-morphism $\bX_\star \to \ctP$ (in fact, to $\ctP'$) refining the map above. 
    \begin{proof}
        Let $(a,t) \in \ctA(S)$. The data of $t$ determines an isomorphism $J_a|_{\infty_S} \cong T|_{\infty_S}$ by \cite[2.1.1]{ngofl}. Since $\chi^T_G$ is \'etale on the regular locus and $a(\infty) = \chi^T_G(t)$ is regular, formally-\'etale lifting extends $t$ to a section of $\ft^{\reg}$ on the complete neighborhood $\Spec(\widehat{\O}_{\infty_S})$ of $\infty_{\ctA}$ (where $\widehat{\O}_{\infty_S} := \widehat{\O}_{X_S, \infty_S}$), and thus a canonical isomorphism
        $$
            J_a \times_X \widehat{\O}_{X_S, \infty_S} \cong T \times_X \widehat{\O}_{X_S, \infty_S}.
        $$
        
        Since $\infty$ is a geometric point (though we may choose it to live over a $k$-point), the restrictions $T|_{\infty}$ and $T|_{\widehat{\O}_{\infty}}$ are split. Writing $z_{\infty_S}$ for a uniformizer at $\infty_S$, we can thus build for each $\lambda \in X_\star(T)(\infty_{\bA})$ a $J_a$-torsor (or indeed a $J_a^0$-torsor) by modifying the trivial torsor by $z_{\infty_S}^{-\lambda} \in J_a(F_{\infty_S})/J_a(\widehat{\O}_{\infty_S})$, that is, by gluing the trivial torsors on 
        $$
        \Spec(\widehat{\O}_{\infty_S}) \quad \text{ and the complement } \quad X_S^\circ = X_S \setminus \{\infty_S\}
        $$
        along the automorphism given by multiplication by $z_{\infty_S}^{-\lambda} \in T(F_{\infty_S})$. As in \cite[6.8]{ngoendo}, the induced map into $\pi_0(\ctP)$ corresponds to the one above, and in particular is independent of the uniformizer. 
    \end{proof}
\end{proposition}

\begin{remark}
    We used above that the data of $t$ (a point in the cameral cover over $\infty$) determined the isomorphism $J_a|_{\infty} \cong T|_{\infty}$. Likewise, given a $\Gamma$-splitting cover $\pi : X' \to X$, a choice of $\infty'$ lying over $\infty$ (i.e., a \textit{pointed} twisting datum) determines an isomorphism $T_{\infty} \cong T_0$, the maximal torus in the split form. We obtain a homomorphism $X_\star(T_0) \to \ctP$ which depends upon $\infty'$, but the argument of \cite[6.8]{ngoendo} shows that at the level of $\pi_0(\ctP')$, the induced map is well-defined. Anyways, it corresponds to the canonically constructed morphism above. 
\end{remark}

\begin{remark}\label{rmk.Kottwitz.signs.geom}
    Note that we needed to chose the normalization $\lambda \mapsto [z_\infty^{-\lambda}]$ in the construction of the map above. This in particular fixes the choice of normalization in the geometric component map. We have chosen the one consistent with our Beauville--Laszlo conventions. In \cite{CLII}, the opposite convention is chosen for Kottwitz, with the net effect that the statement of \cite[Lemma 8.3.2]{CLII} (generalized by Lemma \ref{lem.compos.deg.pic} below) has an error by a negative sign. In Chapter \ref{sect.global.point.counts}, they (and we) choose the opposite normalization for the arithmetic Kottwitz map, cf.\ Remark \ref{rmk.Kottwitz.signs} below. The statements of the main arithmetic theorems in \textit{loc.\ cit.} are ultimately independent of the choice of signs.
\end{remark}

\section{\texorpdfstring{The component group $\pi_0(\bP^{1})$ and lattice actions.}{The component group of the truncated Picard stack.}}\label{subsect.comp.grp.action}

The previous section defined a canonical surjective morphism $\bX_\star = \ul{X}_\star(T)(\infty)_{\ctA} \to \pi_0(\ctP)$, which we can of course restrict to the compatible locus. On the other hand, on that locus we have previously defined a degree morphism $\bP \to \cX_\star$, which factors through $\pi_0(\bP)$ because $\cX_\star$ is fiberwise discrete. Note the asymmetry, not visible in the split case, that the former map has as domain the \textbf{split} or \textbf{geometric} cocharacter group, while the latter has codomain valued in the \textbf{rational} data $\cX_\star(M) = H_0(X, \ul{X}_\star(M))_{\tf}$. The composition $\ul{X}_\star(T)(\infty)_{\bA} \to \cX_\star$ is easy to understand. 

\begin{lemma}\label{lem.compos.deg.pic}
    The composition $\bX_\star \to \pi_0(\bP) \to \cX_\star$ over $\wtilde{a} = (a,t) \in \bAe{M}$ is the natural map
    $$
        \ul{X}_\star(T)(\infty) \to \ul{X}_\star(M)(\infty) \to \ul{X}_\star(M)(\infty)_{\pi_1(X, \infty), \,\tf} = \cX_\star(M)
    $$
    dual to $X^\star(M) \to \ul{X}^\star(M)(\infty) \to \ul{X}^\star(T)(\infty)$.
    \begin{proof}
        Fix $(a,t)$ compatible, $M$-elliptic, and work fiberwise. We can compute the composition as $\bX_\star \to \bP \to \cX_\star$ making use of the explicit Beauville--Laszlo patching construction. This reduces the problem to computing the $M$-degree of the $J_a$-torsor $j_\lambda$ defined by Beauville--Laszlo data $z_{\infty}^{-\lambda}$ at $\infty$ (and trivial data everywhere else on $X$). That is, for each character $\mu \in X^\star(M)$, $\deg(j_\lambda)(\mu) := \deg(j_\lambda \times^{J_a} \G_{m,X})$ is the pushforward along
        $$
            J_a \to \ul{X}_\star(M) \otimes_{\Z_X} \G_{m,X} \to \G_{m,X}.
        $$
        On the disk $\widehat{\O}_{\infty}$, $J_a \cong T$ is determined by $t$. Combining Definition \ref{def.adelic.hitchin.deg}, Prop \ref{prop.deg.p.red.BL}, and our conventions for Beauville--Laszlo gluing (Notation \ref{not.BL.conventions}), we calculate 
        $$\deg(j_\lambda)(\mu) = -\val_{\infty}(z_\infty^{-\mu(\lambda)}) = \mu(\lambda).$$ The composition $\lambda \mapsto \mu(\lambda)$ is the map in the statement of the Lemma because the coinvariants mod torsion is the space dual to $X^\star(M)$, as in Proposition \ref{prop.dual.lattice.global.invs.coinvs.mod.tors}.
    \end{proof}
\end{lemma}

We can now introduce a lattice whose action will control the fine $s$-decomposition below.
\begin{definition}
    Define 
    $$
        \cY^1 := \ker(\cY \to \cX_\star)
    $$
    and 
    $$
        \bX_\star^1 := \ker(\bX_\star \to \cX_\star).
    $$
    Then we have a commutative triangle of surjective maps restricting those previously
    $$
        \begin{tikzcd}[ampersand replacement=\&]
            \bX_\star^1 \ar[r] \ar[dr] \& \cY^1 \ar[d] \\
            \& \pi_0(\bP^{1})
        \end{tikzcd}.
    $$
    For each semistandard Levi $M$, define the \emph{degree-zero lattice} $$\bX_{\star,M}^1 := \ker(\bX_\star = \ul{X}_\star(T)(\infty) \to \ul{X}_\star(M)(\infty)_{\pi_1(X, \infty), \tf} = \cX_\star(M)).$$ This is the fiber of $\bX_\star^1$ over each $M$-elliptic point $a$.
\end{definition}
It is immediate from the definitions that for any pair $M \subset L$, one has 
$$
    \bX_{\star,M}^1 \subseteq \bX_{\star,L}^1.
$$

\begin{proposition}\label{prop.dual.levi.chars.fine}
    For semistandard Levi $M$, the dual group is
    $$
        \widehat{\bX_{\star,M}^1} = \Hom_{\Z}(\bX_{\star,M}^1, \ol{\Q}_l^\times) \cong \widehat{T}/Z_{\widehat{M}}^{\Gamma,0},
    $$
    where
    $$
        Z_{\widehat{M}}^{\Gamma, 0} := (Z_{\widehat{M}}^{\Gamma})^0 = ((Z_{\widehat{M}}^0)^{\Gamma})^0.
    $$
    In particular, each finite-order character of $\bX_{\star,M}^1$ is a torsion element of this dual group. These identifications are natural in the Levi: if $M \subseteq L$ then the map
    $$
        \bX_{\star,M}^1 \hookrightarrow \bX_{\star,L}^1
    $$
    is dual to the natural quotient
    $$\widehat{T}/Z_{\widehat{L}}^{\Gamma,0} \twoheadrightarrow \widehat{T}/Z_{\widehat{M}}^{\Gamma,0}.$$
    \begin{proof}
        One has a short-exact sequence
        $$
            \begin{tikzcd}[ampersand replacement=\&]
                0 \ar[r] \& \bX_{\star,M}^1 \ar[r] \& \bX_{\star} \ar[r] \& \cX_\star(M) \ar[r] \& 0
            \end{tikzcd}.
        $$
        By definition, $\cX_\star(M)$ is dual to $X^\star(M)$. The corresponding subtorus is $(Z_{\widehat{M}}^\Gamma)^0 =((Z_{\widehat{M}}^0)^{\Gamma})^0$. The naturality is clear.
    \end{proof}
\end{proposition}
Note the order of operations! $Z_{\widehat{M}}^{\Gamma, 0} \neq (Z_{\widehat{M}}^0)^\Gamma$.

A sublattice of $\bX_{\star,M}^1$ will induce the \emph{coarse} $s$-decomposition below, which will be central in what follows. We have a short exact sequence of tori
$$
\begin{tikzcd}[ampersand replacement=\&]
    1 \ar[r] \& T \cap D(M) \ar[r] \& T \ar[r] \& M^{\ab} \ar[r] \& 1
\end{tikzcd}
$$
and taking $\ul{X}_\star(-)(\infty)$ gives a short exact sequence of $(W_{M,\infty} \rtimes \im(\rho_G))$-modules, hence of $W_a$-modules, which we denote
\begin{equation}\label{eqn.def.XDM.lattices}
    \begin{tikzcd}[ampersand replacement=\&]
    0 \ar[r] \& \bX_{\star,D(M)} \ar[r] \& \bX_\star \ar[r] \& \bX_{\star,M^{\ab}} \ar[r] \& 0
\end{tikzcd}.
\end{equation}
Since the $M$-degree map factors through $M^{\ab}$, it follows that $\bX_{\star, D(M)} \subseteq \bX_{\star,M}^1$.

The distinction between $\bX_{\star, D(M)}$ and $\bX^1_{\star,M}$ is only visible in the non-split setting, as follows. Since $W_{M,\infty}$ acts trivially on $\bX_{\star, M^{\ab}} = \ul{X}_\star(M^{\ab})(\infty)$, 
the monodromy subgroup $W_a \subseteq \W_M \rtimes \im(\rho_G)$ of Definition \ref{def.monodr.inv.Out} acts on it through $\im(\rho_G)$ (onto which $W_a$ surjects). So the coinvariants by $W_a$ and $\pi_1(X, \infty)$ coincide, and the map whose kernel defines $\bX_\star^1$ is the diagonal map in the diagram
$$
\begin{tikzcd}[ampersand replacement=\&]
    0 \ar[r] \& \bX_{\star, D(M)} \ar[r] \ar[d, two heads] \& \bX_{\star} \ar[r] \ar[d, two heads, "q_T"] \ar[dr, "\psi_a"] \& \bX_{\star, M^{\ab}} \ar[r] \ar[d, two heads] \& 0 \\
     0 \ar[r] \& (\bX_{\star, D(M)})_{W_a, \, \tf} \ar[r] \& (\bX_{\star})_{W_a, \, \tf} \ar[r] \& (\bX_{\star, M^{\ab}})_{W_a, \, \tf} \ar[r] \& 0
\end{tikzcd}.
$$
But by Lemma \ref{lem.monodr.inv.spl.torus}, $(\bX_{\star, D(M)})_{W_a}$ is finite, and so the bottom-left term vanishes. So
$$
    (\bX_{\star}^1)_a = \ker(\psi_a) = \ker(q_T),
$$
and we deduce a short-exact sequence
$$
    \begin{tikzcd}[ampersand replacement=\&]
        0 \ar[r] \& \bX_{\star,D(M)} \ar[r] \& \bX_{\star, M}^1 \ar[r] \& \ker(\bX_{\star, M^{\ab}} \to (\bX_{\star, M^{\ab}})_{\Gamma, \tf}) \ar[r] \& 0
    \end{tikzcd}.
$$
In the split case, the $W_a$-action (equivalently, $\Gamma$-action) on $\ul{X}_\star(M^{\ab})(\infty) = X_\star(M^{\ab})$ is trivial, and the last term vanishes. Dualizing this sequence gives easily the following computation.

\begin{proposition}
    The inclusion $\bX_{\star, D(M)} \subseteq \bX_{\star,M}^1$ is dual to the surjection
    $$
        \widehat{T}/(Z_{\widehat{M}}^\Gamma)^0 \to \widehat{T}/Z_{\widehat{M}}^0.
    $$
\end{proposition}

We can describe more fully the fiber of $\cY^1$ at $a \in \bAe{M}$. By definition and the computation above, it is
$$
    \cY^1_a = \ker((\bX_\star)_{W_a} \to (\bX_{\star, M^{\ab}})_{\pi_1(X, \infty), \,\tf} = (\bX_{\star, M^{\ab}})_{W_a, \,\tf} \cong (\bX_{\star})_{W_a, \, \tf}) = ((\bX_\star)_{W_a})_{\tors}.
$$
We deduce immediately finiteness of the fibers.
\begin{lemma}[{{cf.\ \cite[Lemme 8.3.3]{CLI}}}]\label{lem.coinv.finite}
    For every Levi $M$ and $a \in \bAe{M}$, the fiber of $\cY^1$ over $a$ is finite.
\end{lemma}

This group has a useful decomposition that will be conceptually meaningful for us. 
\begin{definition}\label{def.end.symmetry.and.residual.torsion}
    For $a \in \bAe{M}$, call $\cY^1_a$ the \emph{fine symmetry group}. Writing
    $$
        i_{M, W_a} : (\bX_{\star, D(M)})_{W_a} \longrightarrow (\bX_{\star})_{W_a},
    $$
    which may not be injective, the domain is finite and so lands within $\cY^1_a$. Define the \emph{coarse}, or \emph{endoscopic}, \emph{symmetry group} to be $\Sigma_a := \im(i_{M, W_a}) \subseteq \cY^1_a$, and the \emph{residual torsion group} to be $\Xi_M := (\bX_{\star, M^{\ab}})_{\Gamma, \tors}$. 
\end{definition}

We have a short-exact sequence
\begin{equation}\label{eqn.Y1.fiber.descr}
    \begin{tikzcd}[ampersand replacement=\&]
        0 \ar[r] \& \Sigma_a \ar[r] \& \cY^1_a \ar[r] \& \Xi_M \ar[r] \& 0
    \end{tikzcd}.
\end{equation}
In the split case, $W_a$ acts trivially, so $(\bX_{\star, M^{\ab}})_{W_a}$ is torsion-free and the right term vanishes, so $\cY^1_a = \im(i_{M, W_a})= \Sigma_a$.

We record the dual short exact sequence of (\ref{eqn.Y1.fiber.descr}).
\begin{proposition}\label{prop.dual.Y1.fibers}
    Taking $\Hom(-, \ol{\Q}_l^\times)$, the dual of (\ref{eqn.Y1.fiber.descr}) is given by
    $$
        \begin{tikzcd}[ampersand replacement=\&]
            1 \ar[r] \& \widehat{\Xi}_M \ar[r] \ar[d, equals] \& \widehat{\cY^1_a} \ar[r] \ar[d, equals] \& \widehat{\Sigma}_a \ar[r] \ar[d, equals] \& 1 \\
            1  \ar[r] \& \pi_0((Z_{\widehat{M}}^0)^\Gamma) \ar[r] \& \pi_0(\widehat{T}^{W_a}) \ar[r] \& \widehat{T}^{W_a}/(Z^0_{\widehat{M}})^\Gamma \ar[r] \& 1
        \end{tikzcd}.
    $$
    The dual fine symmetry group is equal to
    $$
        \widehat{\cY^1_a} = \widehat{T}^{W_a}/ Z^{\Gamma,0}_{\widehat{M}} ,
    $$
    and the dual residual torsion group is equal to
    $$
        \widehat{\Xi}_M = \frac{(Z_{\widehat{M}}^0)^\Gamma}{Z_{\widehat{M}}^{\Gamma,0}}
    $$
    measuring the failure of the $\Gamma$-invariants to be connected.
    \begin{proof}
        This is an elementary computation. Note that since $W_a \subseteq \W_M \rtimes \im(\rho_G)$ acts on the center through its surjection onto $\im(\rho_G)$, 
        $$
            Z_{\widehat{M}}^{W_a} = Z_{\widehat{M}}^{\Gamma}.
        $$
        Note also that by the finiteness of $(\bX_{\star, D(M)})_{W_a}$, pushforward of cocharacters gives an isomorphism
        $$
            (\bX_{\star})_{W_a, \tf} \xrightarrow{\sim} (\bX_{\star, M^{\ab}})_{W_a, \tf}.
        $$
        The left side is the character lattice of $(\widehat{T}^{W_a})^0$, and the right side is the character lattice of 
        $$
            Z_{\widehat{M}}^{\Gamma,0} = ((Z_{\widehat{M}}^0)^\Gamma)^0.
        $$
        Thus the two descriptions of $\widehat{\cY^1_a}$ agree.
    \end{proof}
\end{proposition}

\begin{remark}\label{rmk.torus.ex.P.deg.rat}
    In order to justify the choice that the degree map should be valued in \textit{rational} cocharacters while Ng\^{o}'s map should be defined on the \textit{geometric} cocharacters, as well as to illustrate the discrepancy between the lattices $\bX_{\star, D(M)} \subseteq \bX_{\star,M}^1$, it may be helpful to recall the case when $G = T$ is an $X$-torus. Assume $G$ splits over a geometrically irreducible $\Gamma$-Galois \'etale cover $X'/X$. For example, we could take the non-split rank 1 torus $G = U_{X'/X}(1) = \Res^1_{X'/X} \G_m$ along a quadratic cover. By commutativity, the Chevalley map $\chi : \fg_D \to \fc_D$ is an isomorphism and 
    $$
        \cM \cong \Map(X, [\fg_D / G]) = \Map(X, \fg_D \times BG) = \Map(X, \fg_D) \times \Bun_G,
    $$
    with the Hitchin map the projection onto the first factor. We then have $\cM_a \cong \cP_a \cong \Bun_G$. The components of $\Bun_G$ are easy to compute, e.g.\ by the sledgehammer of \cite[Theorem 6]{HeinlothRelUnif}:
    $$
        \pi_0(\Bun_G) \cong \pi_1(G_{\ol{\eta}})_\Gamma \cong (\ul{X}_\star(T)(\infty))_\Gamma.
    $$
    Note that this can have torsion: for $G = U_{X'/X}(1)$ it is $\Z/2$. By the Beauville--Laszlo construction of Proposition \ref{prop.inv.map.BL}, Ng\^{o}'s map is easily seen to induce an isomorphism as above. On the other hand, the natural degree map lands inside the rational cocharacter space and necessarily cannot see this torsion. For $G = U(1)$, the only possible degree map $\Z/2 \to \Z$ is the zero map, reflecting the fact that a purely imaginary line bundle has numerical degree 0. In terms of our notation above,
    $$
        \bX_{\star, D(M)} = 0 \subset \bX^1_{\star,M} = \Z,
    $$
    and so all of the fine symmetry is residual torsion:
    $$
        \cY_a^1 = \Z/2, \qquad \Sigma_a = 0, \qquad \Xi_M = \Z/2.
    $$
    
    At the end of the day, this finite number of components causes us no issue; indeed, the residual torsion is essential to obtain the correct constants when we compare to Arthur's endoscopic parameterization. It is only split tori which give rise to a lattice of components, just as in automorphic theory it is only split tori which cause noncompactness convergence issues.
\end{remark}

The structural feature we need from the full $\cY^1$ is that the map $\bX_\star^1 \to \pi_0(\bP^{1})$ is surjective and factors through a finite quotient---thus, the restriction from $\bP$ to $\bP^{1}$ is enough to ensure finite groups of components, so that we can obtain a well-behaved isotypic decomposition for the action of $\pi_0(\bP^{1})$ beyond the $G$-elliptic locus.

The specialization maps between the $\cY^1_a$ above are also natural. 

\begin{proposition}\label{prop.nat.Y1.obvious}
    Let $a \in \bAe{L}$ and $a' \in \bAe{M}$ be a specialization of $a$. Then $M \subset L$ and $W_{a'} \subset W_a$. Moreover, the specialization map $\cY^1_{a'} \to \cY^1_a$ induces a morphism of short exact sequences
    $$
        \begin{tikzcd}[ampersand replacement=\&]
            0 \ar[r] \& \Sigma_{a'} \ar[r] \ar[d] \& \cY^1_{a'} \ar[r] \ar[d] \& \Xi_{M} \ar[r] \ar[d] \& 0 \\
            0 \ar[r] \& \Sigma_a \ar[r] \& \cY^1_{a} \ar[r] \& \Xi_L\ar[r] \& 0
        \end{tikzcd},
    $$
    where the morphism $\Sigma_{a'} \to \Sigma_a$ is induced from the natural map 
    $$
        (\bX_{\star, D(M)})_{W_{a'}} \to (\bX_{\star, D(L)})_{W_a}.
    $$
    \begin{proof}
        The fact that $M \subset L$ is the content of \S \ref{sect.Hitchin.base}, with $W_{a'} \subset W_a$ given by Proposition \ref{prop.monodr.inv.semicts}. The last claim is immediate from the definitions and the computation of $\cY^1_a$ above.
    \end{proof}
\end{proposition}

\section{\texorpdfstring{The coarse and fine $s$-decompositions.}{The s-decomposition.}}

\begin{definition}\label{def.Mpar.open}
    For each semistandard $M$, define the $M$-parabolic open $$U_M = \bigcup_{L \supset M} \bAe{L} \subset \bA_G.$$ 
    Denote $f^{\xi, \Mpar} : \bM^{\xi, \Mpar} \to U_M$ the restriction of the Hitchin map to $U_M$. Of course, if $M = G$ then $\Mpar = \Gell$.
\end{definition}

Collating the properties of $\cY^1$ from the previous section gives the following.

\begin{lemma}\label{lem.glblsect.Pcmpts.finite}
    Over $U_M$, the homomorphisms of sheaves $(\bX_{\star, D(M)})_{U_M} \subseteq \bX_{\star}^1 \to \cY^1 \to \pi_0(\bP^{1})$ induce group homomorphisms
    $$
        \bX_{\star, D(M)} \subseteq \bX_{\star,M}^1 \to \Gamma(U_M, \cY^1) \to \Gamma(U_M, \pi_0(\bP^{1})).
    $$
    These maps are natural in the obvious sense given a containment of Levis $M \subseteq L$ and the corresponding open immersion $U_L \subseteq U_M$. Moreover, $\Gamma(U_M, \pi_0(\bP^{1}))$ is a finite abelian group.
    \begin{proof}
        Given $a \in \bAe{L}$ for $M \subseteq L$ we get a composition $$\bX_{\star, M}^1 \subseteq \bX_{\star, L}^1 \to \cY^1_a.$$ It is clear from the definition of $\cY^1$ as a quotient that these combine to a map to the global sections of $\cY^1$, which we can then compose with $\Gamma(U_M, \cY^1) \to \Gamma(U_M, \pi_0(\bP^{1}))$. The finiteness on global sections of the sheaf of components follows from the finiteness on fibers of Lemma \ref{lem.coinv.finite}. The naturality follows from that in Proposition \ref{prop.nat.Y1.obvious}.
    \end{proof}
\end{lemma}

We have assembled all necessary ingredients to generalize the isotypic decomposition of \cite[8.5]{CLII} in the nonsplit setting. In fact, we will obtain two isotypic decompositions corresponding to the two lattices $\bX_{\star, D(M)} \subseteq \bX_{\star,M}^1$. The fine decomposition provides the most direct generalization of \cite[\S 6.2]{ngofl} outside the elliptic locus, while the coarse decomposition will appear naturally for the purposes of point-counting.

We have constructed an action of $\bP^{1}$ on $\bM^{\xi}$ by \ref{prop.Pic.action.shift.convex}, which we restrict to $f^{\xi, \Mpar} : \bM^{\xi, \Mpar} \to U_M$. Since the group is commutative, pullback induces an action on $\pHi(R f_\star^{\xi, \Mpar} \ol{Q}_l)$, which factors through its fiberwise component group $\pi_0(\bP^{1})$ by the \textit{lemme d'homotopie} (\cite[Lemme 3.2.3]{ngolaumonunitaire}). Thus we have an isotypic decomposition for the action of $\Gamma(U_M, \pi_0(\bP^{1}))$ on $\pHi(R f_\star^{\xi, \Mpar} \ol{\Q}_l)$:
$$
    \pHi(R f_\star^{\xi, \Mpar} \ol{\Q}_l) = \bigoplus_{s' \in \Hom_{\Z}(\Gamma(U_M, \pi_0(\bP^{1})), \ol{\Q}_l^\times)} \pHi(R f_\star^{\xi, \Mpar} \ol{\Q}_l)_{s'}.
$$
But we have maps
$$
    \bX_{\star,D(M)} \subseteq \bX_{\star,M}^1 \to \Gamma(U_M, \pi_0(\bP^1)).
$$
The fine $s$-decomposition is the decomposition which remembers only the action of $\bX_{\star, M}^1$, and the coarse $s$-decomposition is the decomposition which remembers only the action of $\bX_{\star,D(M)}$. 

\begin{theorem}[{{cf.\ \cite[Thm.\ 8.5.1, Cor.\ 8.5.3]{CLII}}}]\label{thm.s.decomp}
    For every semistandard $M$, denote 
    $$
        S^f_M := (\widehat{T}/Z_{\widehat{M}}^{\Gamma,0})_{\tors}, \qquad S^c_M := (\widehat{T}/Z_{\widehat{M}}^0)_{\tors}.
    $$
    For any $i \in \mathbb{Z}$ and $s \in S_M^f$, let 
    $$
        \pHi(Rf^{\xi, \Mpar} \ol{\Q}_l)_{s}
    $$
    denote the sum of all $s'$-isotypic components where $s' \in \Hom_{\Z}(\Gamma(U_M, \pi_0(\bP^{1})), \ol{\Q}_l^\times)$ restricts to $s$ on $\bX_{\star,M}^1$. 

    Then one has the \emph{fine $s$-decomposition}, a direct-sum decomposition
    $$
        \pHi(R f_\star^{\xi, \Mpar} \ol{\Q}_l) = \bigoplus_{s \in S_M^f} \pHi(R f_\star^{\xi, \Mpar} \ol{\Q}_l)_s
    $$
    with only finitely many nonzero summands. Grouping the summands along $S_M^f \to S^c_M$ gives the \emph{coarse $s$-decomposition}
    $$
        \pHi(R f_\star^{\xi, \Mpar} \ol{\Q}_l) = \bigoplus_{\ol{s} \in S_M^c} \pHi(R f_\star^{\xi, \Mpar} \ol{\Q}_l)_{\ol{s}}, \qquad \pHi(R f_\star^{\xi, \Mpar} \ol{\Q}_l)_{\ol{s}} = \bigoplus_{s \mapsto \ol{s}} \pHi(R f_\star^{\xi, \Mpar} \ol{\Q}_l)_{s}.
    $$

    If $M \subseteq L$, write $j_M^L  : U_L \to U_M$ for the inclusion. Then the fine decomposition is functorial in the sense that for any $s \in S^f_M$,
    $$
        (j_M^L)^\star(\pHi(Rf_\star^{\xi, \Mpar} \ol{\Q}_l)_s) = \bigoplus_{s'} \pHi(Rf_\star^{\xi, \Mpar} \ol{\Q}_l)_{s'},
    $$
    where the sum runs over the preimage of $s$ under the natural map 
        $$
        S^f_L \to S^f_M.
        $$
    The coarse decomposition satisfies the analogous functoriality for the fibers over 
    $$
        S^c_L \to S^c_M.
        $$
    
                                                                                                    \begin{proof}
        The proof of \textit{loc.\ cit.} transfers to our setting. We know by Lemma \ref{lem.glblsect.Pcmpts.finite} that $\Gamma(U_M, \pi_0(\bP^{1}))$ is a finite abelian group, so we have an action of $\bX_{\star, M}^1$ on the cohomology which factors through a finite quotient. In particular, the characters of $\Gamma(U_M, \pi_0(\bP^{1}))$ pull back to a finite number of finite-order characters on $\bX_{\star, M}^1$, which correspond under Proposition \ref{prop.dual.levi.chars.fine} to a finite number of elements of $S^f_M$. This gives the existence and finiteness of the decomposition, and the coarse decomposition follows by restriction.

        For the compatibility along $M \subseteq L$, functoriality of $(j_M^L)^\star$ implies that $\Gamma(U_M, \pi_0(\bP^{1}))$ acts on $(j_M^L)^\star(\pHi(Rf_\star^{\xi, \Mpar} \ol{\Q}_l)_s)$ through the restriction map $\Gamma(U_M, \pi_0(\bP^{1})) \to \Gamma(U_L, \pi_0(\bP^{1}))$. Naturality in \ref{lem.glblsect.Pcmpts.finite} says this corresponds to the natural map $\bX_{\star, M}^1 \to \bX_{\star, L}^1$. Naturality in Proposition \ref{prop.dual.levi.chars.fine} implies this dualizes to the claimed decomposition.
    \end{proof}
\end{theorem}

Note that the coarse decomposition and the fine decomposition differ on the $M$-elliptic locus essentially in whether they capture the residual torsion group $\Xi_M$. The coarse decomposition is the decomposition appropriate for weighted point-counting, but in order to compare with the results of Ng\^{o}, we will need to pass the coarse and the fine decomposition on the elliptic locus. This will ensure our global version of the WFL matches the expected form, as in \cite[\S VII.2.2]{MWStabII}. Over a fiber $a \in \bAe{M}$, of course, the set of isotypic characters $s \in \widehat{\cY^1_a}$ which restrict to a given $\ol{s}$ can be described as a torsor under the dual $\widehat{\Xi}_M$.

Our descriptions of $\widehat{\Sigma}_a$ and $\widehat{\cY^1_a}$ imply a strong bound on the support of the $s$-isotypic component.

\begin{proposition}[{cf.\ \cite[Th\'eor\`eme 8.5.1(1)]{CLII}}]\label{prop.s.isotypic.cpt.support}
    Let $a \in \bAe{M}$ have monodromic invariant $W_a$. Then the stalk of the fine $s$-isotypic component, for $s \in (\widehat{T}/Z_{\widehat{M}}^{\Gamma, 0})_{\tors}$, vanishes unless $s$ admits a lift to $\widehat{T}^{W_a}$, i.e., unless
    $$
        s \in \widehat{\cY_a^1} \cong \pi_0(\widehat{T}^{W_a}) \subset (\widehat{T}/Z_{\widehat{M}}^{\Gamma, 0})_{\tors}.
    $$
    Likewise, the stalk of the coarse $\ol{s}$-component vanishes unless 
    $$
        \ol{s} \in \widehat{\Sigma}_a \cong \widehat{T}^{W_a}/(Z_{\widehat{M}}^0)^{\Gamma}.
    $$
    \begin{proof}
        If the fine $s$-component is nonzero at $a$, then by definition there is a character $s'$ of $\Gamma(U_M, \pi_0(\bP^{1}))$ with nonzero isotypic component at $a$. The action of $\bX_{\star, M}^1$ on the stalk $\pHi(R f_\star^{\xi, \Mpar} \ol{\Q}_l)_{a}$ is defined through the composition
        $$
            \bX_{\star, M}^1 \to \Gamma(U_M, \pi_0(\bP^{1})) \to \pi_0(\bP^{1}_a),
        $$
        so considering the (finer) isotypic decomposition by $\pi_0(\bP^{1}_a)$, we conclude at least one character $s''$ of this last group has nonzero component restricting to $s$. But the composition $\bX_{\star, M}^1 \to \pi_0(\bP^{1}_a)$ factors through $\cY^1_a$; the image of $s''$ in $\widehat{\cY^1_a}$ is a character which maps to $\widehat{\bX_{\star, M}^1}$. It follows from the dual computations above that the dual map to 
        $$\bX_{\star, M}^1 \twoheadrightarrow \cY^1_a$$ is
        $$ \widehat{T}^{W_a}/Z^{\Gamma,0}_{\widehat{M}} \hookrightarrow (\widehat{T}/Z^{\Gamma,0}_{\widehat{M}})_{\tors}.
        $$
        The statement for the coarse decomposition follows from that for the fine decomposition, or by a parallel argument. 
    \end{proof}
\end{proposition}

\chapter{\texorpdfstring{The good locus and $\delta$-invariants}{The good locus and delta-invariants}}\label{sect.good.locus}

In this section, we recall the definition of the locus $\cA^{\good}$ over which Ng\^{o}'s support theorem applies. Applying this to the terms in the $s$-decomposition of Theorem \ref{thm.s.decomp}, we will in the next section deduce the cohomological version of Arthur's stabilization over $\bA \cap \ctA^{\good}$, see Theorem \ref{thm.coh.wfl} below. In \cite[\S 9]{CLII}, the authors prove by a deformation argument a general sufficient criterion (Th\'eor\`eme 9.1.3 in \textit{loc.\ cit.}) for a point to lie in $\cA^{\good}$. This chapter is a direct adaptation of their argument to our generality. This characterization will be most crucially used in the final proof of the WFL: a local parameter is globalized, and we must do so in a way which lands in the good locus in order to apply our geometric theorems.

\section{N\'eron models.}\label{subsect.neron.models}

Recall that the Picard stack $\cP \to \cA$ parameterizes $J_a$-torsors on the curve. In \S \ref{sect.reg.cent.Pic}, we have described its group of connected components. The structure of its identity component is also analyzed over $\cA^{\heartsuit}$ by Ng\^{o}, as follows.\footnote{In the commutative case compare also with geometric class field theory as in \cite{serregeocft}.}

Ng\^{o} in \cite[Proposition 4.8.2]{ngofl} shows the fiber $\cP_a^0$ of the identity component of the Picard stack over a geometric point $a$ admits a canonical decomposition
\begin{equation}\label{eqn.aff.ab.decomp.Pic.stack}
    1 \to \cR_a \to \cP_a^0 \to (\cP_a^{\flat})^0 \to 1,
\end{equation}
with $\cR_a$ affine and $(\cP_a^{\flat})^0$ an abelian stack.\footnote{Note that this $\cR_a$ is appropriate for our purposes, but can be smaller than Ng\^{o}'s $\cR_a = \ker(\cP_a \to \cP_a^{\flat})$. The two have the same dimension.} Here, $\cP_a^{\flat}$ is the moduli space of torsors for the global finite-type N\'eron model $J_a^{\flat}$ of the regular centralizer $J_a$. 
Recall that by the latter we mean the N\'eron model whose set of $\O_v$-points at each geometric point on the curve is the maximal bounded subgroup $J_a(F_v)^{\flat}$ of $J_a(F_v)$; its component group in the special fiber is the torsion subgroup of that for the lft N\'eron model of \cite[Ch.\ 10]{blrneron} (cf.\ \cite[\S 3.b]{pappasrapoport} and the references therein). This N\'eron model $J_a \to J_a^{\flat}$ is the terminal object in the category of smooth finite-type $X$-group schemes $Z$ equipped with maps $J_a \to Z$ extending a specified identification over $U_a = a^{-1}(\fc^{\Greg}_D)$. It is easier to describe than $J_a$. 

\begin{proposition}[{{\cite[Corollaire 4.8.1]{ngofl}}}]\label{prop.Jflat.Weil.res.description}
    We have a natural isomorphism
    $$
        J_a^\flat \cong (\Res_{X_a^{\flat}/X} (T \times_X X_a^{\flat}))^W,
    $$
    where $X_a^{\flat} \to X_a$ is the normalization of the cameral cover associated to $a$.
    \begin{proof}
        The unit of the adjunction induces a natural map 
        $$J_a \to \Res_{X_a^{\flat}/X} (J_a \times_X X_a^{\flat}),$$ 
        which lands inside the $W$-fixed points (since $W$ acts on $X_a$ over $X$, by functoriality of normalization it also acts on $X_a^{\flat}$). By the Galois description of $J$ in \S \ref{subsect.descr.galois}, $J_a|_{U_a}$ and $T|_{U_a}$ are canonically isomorphic after pullback to $X_a^{\flat}$, and we know that $T$ is the global Neron model of $T|_{U_a}$, so by its universal property we get a $W$-equivariant $X_a^{\flat}$-map
        $$
            J_a \times_X X_a^{\flat} \to T \times_X X_a^{\flat}.
        $$
        Applying Weil restriction and composing gives the natural map $\varphi: J_a \to (\Res_{X_a^{\flat}/X} (T \times_X X_a^{\flat}))^W$. But the same argument applies replacing $J_a$ with any smooth finite-type group $X$-scheme isomorphic to $J_a$ on $U_a$, so we find that $\varphi$ is terminal in the category of above and is the structure map of the global N\'eron model. 
        
        Alternatively Beauville--Laszlo gluing and the universal property of the N\'eron model, we can see that $J_a^{\flat}$ is defined by gluing local N\'eron models for $J_{a}|_{\widehat{\O}_{v}}$ at each $v \in X \setminus U_a$ with $J_a|_{U_a}$. The global map above is an isomorphism on $U$, so it suffices to show the corresponding natural map in the local setting is an isomorphism at each place $v \in X \setminus U_a$. There, a local version of the same adjunction argument above shows that $(\Res_{X_a^{\flat}/X} T \times_X X_a^{\flat})^W$ satisfies the universal property of the local finite-type N\'eron model.
    \end{proof}
\end{proposition}

For the Picard stacks under consideration, we denote $(-)^{\#}$ for the functor which constructs the fppf sheaf of isomorphism classes, i.e., the coarse quotient. From the description above, it is easy to see that $(\cP_a^{\flat})^0$ is an \textit{abelian stack}, that is to say, it is of the form $A \times BD$ for $A$ an abelian variety and $D$ a diagonalizable group. In particular, its coarse quotient is an abelian variety, whose properness is the core of the later cohomological arguments.

Since the universal $X$-map $J_a \to J_a^\flat$ is an isomorphism on the open $U$ where $a$ is regular semisimple, $\cR_a$ is a product of affine, finite-type group schemes $\cR_{a,v}$ supported on $v \in X \setminus U$. 

\section{\texorpdfstring{The $\delta$-invariant.}{The delta-invariant.}}
Following Ng\^{o}'s support theorem, we study the cohomology only over a locus $\cA^{\good}$ where the abelian variety part of $\cP_a$ controls its action on cohomology.
\begin{definition}
    For $a \in \cA^{\heartsuit}(k)$ (and thereby for $(a,t) \in \ctA$, define the $\delta$-invariant
    $$
        \delta_a := \dim(\cR_a).
    $$
\end{definition}

The above characterization of $J_a^\flat$ implies a basic formula for the $\delta$-invariant.

\begin{proposition}[{{\cite[Corollaire 4.9.4]{ngofl}}}]\label{prop.norm.form.delta}
    For any $a \in \cA^{\heartsuit}(k)$, one has
    $$
        \delta_a = \dim H^0(X, (\ft \otimes_{\O_X} (\pi_{a \star}^{\flat} \O_{X_a^\flat}/\pi_{a\star} \O_{X_a}))^W),
    $$
    where $\pi_a^\flat : X_a^\flat \to X$ is the composite of the normalization map with $\pi_a : X_a \to X$.
\end{proposition}

On the other hand, there is also the global version of Bezrukavnikov's dimension formula \cite{bezrukavnikovdimension}---note that locally at each place over $k = \ol{k}$ our groups become split, which is the generality of \textit{loc.\ cit.}

\begin{proposition}[{{\cite[Proposition 4.9.5]{ngofl}}}]\label{prop.bezru.form.delta}
    For any $a \in \cA^{\heartsuit}(k)$, one has
    $$
        \delta_a = \frac{1}{2} \left(\|\Phi\| \deg(D) - \sum_{v \in X \setminus U} c_i(v) \right),
    $$
    where $c_i(v)$ is the difference between $\rank T$ and the toric rank of the special fiber of $J_a^\flat$ at $v$. 
\end{proposition}

By pullback, we have a corresponding $\delta$-invariant $(a, t) \mapsto \delta(a)$ on $\ctA$. 

\begin{proposition}[{{\cite[\S 5.6, Lemme 5.6.3]{ngofl}}}]
    As a function on the base $\ctA$, the $\delta$-invariant is upper-semicontinuous.
                                                                                        \end{proposition}

Note that the proof in \textit{loc.\ cit.} does not work, as, in their notation,
$$
    \# P_s[l] = l^{\mu_s + 2\tau_s}
$$
rather than the claimed
$$
    \# P_s[l] = l \dim(T_x) + l^2 \dim(A_x).
$$
The result is nevertheless true, cf.\ \cite[Exp.\ X, Remarques 8.7]{SGAIII}.

\section{\texorpdfstring{The good locus.}{The good locus.}}

With these preliminaries out of the way, we can define the good locus as the largest open satisfying a semismallness condition for the affine part of the Picard stack.

\begin{definition}\label{def.good.locus}
    Define $\ctA^{\good} \subset \ctA$ to be the largest open subset consisting of Zariski points such that 
    \begin{equation}\label{eqn.good.loc}
        \codim_{\ctA}(a) \ge \delta_a. 
    \end{equation}
    That is, $\ctA^{\good}$ is the topological interior of the subset of such points. Define $\bA^{\good}$ to be the intersection of $\ctA^{\good}$ with the compatible locus.
\end{definition}

\begin{remark}
    Note a subtlety of the definition above: $\ctA^{\good}$ does not necessarily consist of \emph{all} points satisfying (\ref{eqn.good.loc}), which could be topologically poorly behaved. In Ng\^{o}'s work, he works over an even (\textit{a priori}) smaller open, defined as follows. Consider the stratification of $\ctA$ by $\delta$-value:
    $$
        \ctA = \bigsqcup_\delta \ctA_{\delta}.
    $$
    Then (\ref{eqn.good.loc}) is satisfied on a given stratum $\ctA_{\delta}$ if and only if it is satisfied at the generic points of the irreducible components of that stratum. If $\ctA^{\good} \neq \ctA$, write $\delta_0$ for the minimal value such that $\codim_{\ctA}(\ctA_{\delta_0}) < \delta_0$, that is, that the generic point $\eta$ of some irreducible component of $\ctA_{\delta_0}$ satisfies $\codim_{\ctA}(\eta) < \delta_0$. Then the subset
    $$
        \bigsqcup_{\delta \ge \delta_0} \ctA_{\delta}
    $$
    is clearly closed, and its complement is open. This is the good locus as defined in Ng\^{o}; it is more strict than the definition above, as we have thrown away all components of all strata that could potentially lie in the closure of a bad point, regardless of whether they actually do.
\end{remark}

Results like \cite[Th\'eor\`eme 6.1.3]{CLII} indicate that the good locus is quite large. It is conjectured by Ng\^{o} that the good locus is all of $\ctA$, and this is known in characteristic 0.

\begin{proposition}[{{\cite[Lemme 5.6.4, Proposition 5.7.2]{ngofl}}}]
    If $\deg(D) > 2g$, then the $0$-stratum $\ctA_0$ is nonempty. In particular, the good locus is nonempty. 

    Moreover, for any given $\delta$, there exists a bound $N$ such that if $\deg(D) > N$ then $\codim_{\ctA}(\ctA_{\delta}) \ge \delta$.
\end{proposition}

\section{\texorpdfstring{The affine $\delta$-invariant.}{The affine delta-invariant.}}\label{subsect.aff.delta.inv}

While the decomposition of (\ref{eqn.aff.ab.decomp.Pic.stack}) resembles a Chevalley decomposition, caution is warranted by the fact that group stacks do not have a uniquely defined Chevalley decomposition. By the below, we could rigidify our Picard stack to obtain a scheme. Applying the Chevalley decomposition to such a rigidification, we can define an invariant $\delta_a^{\aff}$ which is closely related to, but not in general equal to, $\delta_a$, measuring the kernel of the map from $\cP^0$ to a genuine abelian \textit{variety}. In the next section, we will make use of some results of Ng\^{o} for group schemes which when applied to a rigidification involve $\delta_a^{\aff}$.

\begin{proposition}[{{\cite[Proposition 4.5.7]{ngofl}}}]\label{prop.pic.stack.rigidification}
    For each $a \in \cA^\infty$, define $P^\infty_a$ to be the stack classifying $J_a$-torsors equipped with a trivialization at $\infty$. Then $P^\infty_a$ is representable by a smooth scheme locally of finite type. By modifying the trivialization we have an action of $J_{a,\infty}$ on $P^{\infty}_a$. This action induces an canonical isomorphism $[P_a^\infty / J_{a,\infty}] \cong \cP_a$.
    
    Moreover, a choice of $t$ determines an isomorphism $T_\infty \cong J_{a,\infty}$, so $(a,t) \in \ctA$ determines a canonical isomorphism $[P_a^\infty / T_{\infty}] \cong \ctP_{(a,t)}$.
\end{proposition}

Let $a \in \cA^{\heartsuit}$ and choose any presentation of $\cP_{a}^0$ as a cokernel of smooth fiberwise-connected commutative group schemes 
$$
    Q_{-1,a} \to Q_{0,a} \to \cP_{a}^0 \to 1,
$$
with $Q_{-1}$ affine. We then have $(\cP_a^0)^\# = Q_{0,a}/\im(Q_{-1,a} \to Q_{0,a})$, the quotient as fppf sheaves, and the inertia of $\cP_{a}^0$ at every point is (\cite[Proposition 4.5.5]{ngofl}
$$
    H^0(X, J_a) = Z_{\G}^{\im(\rho_G)},
$$
with dimension $\dim A_G^{\spl}$. For example, such a presentation can come from a framing at a regular point on the curve, as above. 

Consider the Chevalley decomposition of $Q_{0,\ol{a}}$ at a geometric point $\ol{a}$ over $a$:
$$
    1 \to Q_{0, \ol{a}}^{\aff} \to Q_{0, \ol{a}} \to Q^{\ab}_{0, \ol{a}} \to 1.
$$
As observed at the beginning of \cite[\S 4.12]{ngofl}, since $Q_{-1, \ol{a}}$ is affine the composition $ Q_{-1, \ol{a}} \to Q_{0, \ol{a}} \twoheadrightarrow Q_{0,\ol{a}}^\ab$ associated to any Chevalley decomposition has image which is proper, affine, smooth, and connected, hence trivial, and we have an induced surjection with kernel which we denote $K_{\ol{a}}$:
$$
    1 \to K_{\ol{a}} \to \cP_{\ol{a}}^0 \to Q_{0,\ol{a}}^\ab \to 1.
$$

\begin{definition}\label{def.aff.delta.inv}
    Given $a \in \cA^{\heartsuit}$, define the \emph{affine $\delta$-invariant} to be the kernel of the map to the abelian quotient: $$\delta_a^{\aff} := \dim K_{\ol{a}}.$$
\end{definition}

Note that this is independent of choice of geometric point $\ol{a}$ over $a$. Moreover, $\dim(Q_{0,\ol{a}}^{\ab})$ and thus $\delta_a^{\aff}$ are also independent of choice of rigidification (as long as $Q_{-1,a}$ is required to be affine). Indeed, $Q_{0, \ol{a}}^{\ab}$ is the abelian Chevalley quotient of the coarse moduli of isomorphism classes $(\cP^0_{\ol{a}})^\#$. 

The basic comparison between the two invariants $\delta_a$ and $\delta_a^{\aff}$ is the following computation.
\begin{proposition}[{{cf.\ \cite[10.5.1]{CLII}}}]\label{prop.delta.invs.comparison}
    For a geometric point $\wtilde{a} = (a,t)$, we have
    $$
        \delta_a - \delta_a^\aff = \dim(\T^{W_{\wtilde{a}}}),
    $$
    where $W_{\wtilde{a}}$ is the finite monodromy subgroup of Definition \ref{def.monodr.inv.Out}. In particular, $\delta_a^{\aff}$ is upper semicontinuous, and $$\delta_a^{\aff} \le \delta_a - \dim A_G^{\spl}$$ with equality exactly on the $G$-elliptic locus. 
    \begin{proof}
        By the short exact sequences defining them,
        $$
            \delta_a = \dim \cP_a^0 - \dim(\cP^\flat_a)^0, \quad \text{ and } \quad \delta_a^{\aff} = \dim(\cP_a^0) - \dim(Q^{\ab}_{0,a}).
        $$
        Subtracting gives
        $$
            \delta_a - \delta_a^{\aff} = \dim(Q^{\ab}_{0,a}) - \dim(\cP^\flat_a)^0.
        $$
        
        To compute this difference, the exact sequence (\ref{eqn.aff.ab.decomp.Pic.stack}) descends to the fppf sheaves of isomorphism classes:
        $$
            \R_a \to (\cP_a^0)^\# \to ((\cP_a^{\flat})^0)^\# \to 1,
        $$
        so the kernel of the right map is affine. Since the affine part of $(\cP_a^0)^\#$ maps trivially to the abelian variety $((\cP_a^{\flat})^0)^\#$, 
        the map descends to
        $Q_{0,a}^{\ab} \to ((\cP^{\flat}_a)^0)^\#$, whose kernel is a quotient of $\cR_a$, hence affine, but also proper as a closed subvariety of $Q_{0,a}^{\ab}$. So the descended map is an isogeny. It follows that the residual difference in dimensions is the stabilizer dimension of $\cP^{\flat}_a$, which by \cite[\S 4.11]{ngofl} is $$\dim H^0(X, J_a^{\flat}) = \dim \T^{W_{\wtilde{a}}}.$$ 
        
        Upper semicontinuity of $\delta_a^{\aff}$ follows from 
        $$\delta_a^{\aff} = \dim \cP_a^0 - \dim Q_{0,a}^{\ab};$$
        the first term is locally constant and the second is lower semicontinuous by \cite[Exp.\ X, Rem.\ 8.7]{SGAIII}. We have a canonical containment $A_G^{\spl} \subseteq \T^{W_{\wtilde{a}}}$ with equality of dimensions on the elliptic locus.
    \end{proof}
\end{proposition}

\section{A sufficient criterion for goodness.}

A substantial amount of \cite[\S 9]{CLII} is devoted to giving a deformation-theoretic argument for a general criterion under which a point $a \in \ctA_G$ lies in the good locus. In this section we adapt their argument.

Recall from \cite[\S 5.2]{ngofl} the moduli problem of normalized cameral curves.

\begin{definition}
    Let $\cB = \cB_G$ the functor which associates to a $k$-scheme $S$ the groupoid of triples
    $(a, X_a^\flat, \nu)$, where
    \begin{itemize}
        \item $a \in \cA^{\heartsuit}_G(S)$,
        \item $X_a^\flat$ is a proper smooth curve over $S$,
        \item $\nu : X_a^\flat \to X_a$ is a $W$-equivariant normalization in family of the cameral curve $X_a$ (cf.\ \S \ref{subsect.Ngo.monod.comp}, \cite[\S 5.1]{ngofl}). 
    \end{itemize}
    There is a natural map $\cB_G \to \cA^{\heartsuit}$ which is a bijection on geometric points. We define then
    $$
        \ctB := \cB \times_{\cA^{\heartsuit}} \ctA.
    $$
    
\end{definition}
Both $\cB$ and $\ctB$ are representable by $k$-schemes of finite type. Indeed, they admit an equivalent description which is essentially the definition of \cite{CLII}.

\begin{lemma}[{\cite[Lemme 5.2.2]{ngofl}}]\label{lem.CL.presentation.cB}
    The functor $\cB_G$ is isomorphic to the functor which associates to $k$-scheme $S$ the set of pairs $(X_a^\flat, \gamma)$, where
    \begin{itemize}
        \item $X_a^\flat$ is a smooth proper curve over $S$ equipped with a $W$-action and a morphism $\pi : X_a^\flat \to X \times S$ which is a $W$-torsor on a nonempty open which surjects onto $S$, and
        \item $\gamma : X_a^\flat \to \ft_D \times S$ is a $W$-equivariant $S$-morphism such that the preimage of the regular semisimple locus is dense in the fiber over each $s \in S$.
    \end{itemize}
\end{lemma}

Since the map $\ctB \to \ctA$ is a bijection on geometric points, as in \cite[\S 5.4-5.6]{ngofl} the $\delta$-stratification on $\ctA$ induces a \emph{locally constant} $\delta$-stratification on $\ctB$.

\begin{definition}[The cotangent complex on $\cB$]
    Fix a geometric point $b = (a, X^\flat_a, \nu) \in \cB(\ol{k})$, identified with $\pi : X_a^\flat \to X$ and $f : X_a^\flat \to \ft_D$ as in Lemma \ref{lem.CL.presentation.cB}. Consider the cotangent complex of $f$ over $X$:
    $$
        L_b := L_{X^\flat_a/\ft_D} = [f^\star \Omega^1_{\ft_D/k} \to \Omega^1_{X_a^\flat/k}] \cong [f^\star \Omega^1_{\ft_D/X} \to \Omega^1_{X_a^\flat/X}]
    $$
    supported in degrees $(-1,0)$. Note that $\cH^{-1}(L_b)$ is a vector bundle of rank $r = \rank G$, and $\cH^0(L_b) = \Omega^1_{X_a^\flat/X_a}$ is torsion. The complex $L_b$ is naturally $W$-equivariant, where $W/X$ is the finite \'etale Weyl group scheme. Define the deformation complex
    $$
        \cE_b := (\pi_\star \RsHom_{X^\flat_a}(L_b, \O_{X_a^\flat}))^W,
    $$
    a complex of $\O_X$-modules. Then define $E_b := R\Gamma(X, \cE_b)$. 
\end{definition}
This is the appropriate deformation complex which controls the infinitesimal structure of $\cB$ at $b$. In particular, as in \textit{loc.\ cit.},
\begin{itemize}
    \item $H^0(E_b)$ is the space of infinitesimal automorphisms of $\cB$ at $b$,
    \item $H^1(E_b) = T_b \cB$, and
    \item $H^2(E_b)$ is the obstruction to smoothness.
\end{itemize}

\begin{remark}\label{rmk.H0Eb,vanishes}
    In fact, $H^0(E_b) = 0$. Choosing the presentation
    $$
        L_b = [f^\star \Omega^1_{\ft_D/k} \to \Omega^1_{X_a^\flat/k}],
    $$
    by smoothness this perfect complex dualizes to the complex of tangent bundles
    $$
        \RsHom(L_b, \O_{X_a^\flat}) \cong [T_{X_a^\flat/k} \to f^\star T_{\ft_D/k}].
    $$
    Its 0-th cohomology is the kernel of this map, which in particular lies in the kernel of $T_{X_a^\flat/k} \to T_{X/k}$. But generically $X_a^\flat \to X$ is a $W$-torsor, hence finite \'etale, and the map on tangent bundles is generically an isomorphism. So a section of $T_{X_a^\flat/k}$ which vanishes upon pushforward is generically, hence identically, 0. We have shown
    $$
        \cH^0(\RsHom(L_b, \O_{X_a^\flat})) = 0,
    $$
    and taking $\pi_\star$ and $W$-invariants are exact, so the same holds for $\cH^0(\cE_b)$. The hypercohomology spectral sequence shows that the only possible contribution to $H^0(E_b)$ is $H^0(X, \cH^0(\cE_b)) = 0$. 
\end{remark}

\begin{remark}
When the group $W$ is constant (\'etale under our characteristic hypotheses), its invariants are exact and we have
$$
    R\Gamma(X, (\pi_\star \RsHom_{X^\flat_a}(L_b, \O_{X_a^\flat}))^W) = R\Gamma(X_a^\flat, \RsHom_{X^\flat_a}(L_b, \O_{X_a^\flat}))^W = R\Hom_{\O_{X_a^\flat}}(L_b, \O_{X_a^\flat})^W,
$$
so this generalizes the definition of \cite{CLII}. However, the right side makes no sense for $W$ a group scheme over $X$. In our generality one could interpret $E_b$ as the equivariant derived Hom
$$
    R \Hom_W(L_b, \O_{X_a^\flat})
$$
but the definition above is more explicit (and equivalent, since $W/X$ is linearly reductive).
\end{remark}

In the same way, the deformations of $\cA_G$ at $a$ are controlled by
$$
    R\Gamma(X, \pi_{X_a,\star} \RsHom_{X_a}(L_{X_a/\ft_D}, \O_{X_a})^W).
$$
This is in fact a complex with cohomology only in degree 1, as befits a smooth scheme. The deformation complex term has a simple description.

\begin{lemma}[{\cite[Lemme 9.6.1]{CLII}}]\label{lem.H-1.LX_a}
    We have
    $$
        L_{X_a/\ft_D} \cong \cH^{-1}(L_{X_a/\ft_D})[1] =\pi^\star_{X_a} a^\star \Omega^1_{\fc_D/X}[1]
    $$
    for $\pi_{X_a} : X_a \to X$ the natural map.
    \begin{proof}
        The argument given in \textit{loc.\ cit.} applies verbatim. The constructions in question are \'etale-local, anyways, so the result transfers by outer twisting from the split case.
    \end{proof}
\end{lemma}

\begin{definition}\label{def.kappa.inv}
    For any $W$-equivariant coherent sheaf $\cL$ on $X_a^\flat$, define
    $$
        F_L(\cL) := (\pi_\star(\cH^{-1}(L_b) \otimes \cL))^W.
    $$
    This is an exact functor because $W$ is linearly reductive and $\pi$ is finite. Then define for each $x \in |X|$ the \emph{local invariant obstruction}
    $$
        \kappa_x(b) := \mathrm{length}_{\O_{X,x}} F_L(\Omega^1_{X_a^\flat/X})_x
    $$
    and the \emph{global invariant obstruction}
    $$
        \kappa(b) = \sum_{x \in |X|} \kappa_x(b).
    $$
    For any geometric point $a \in \cA_G^{\heartsuit}$, we write $\kappa_x(a) = \kappa_x(b)$ and $\kappa(a) = \kappa(b)$ for $b$ the unique point of $\cB$ corresponding to $a$ under the natural bijection $\cB \to \cA^{\heartsuit}$.
\end{definition}

We collect the following basic properties of this invariant, whose proofs follow exactly as in \cite[\S 9.7]{CLII}.

\begin{lemma}[{\cite[Remarques 9.7.2, Proposition 9.7.4]{CLII}}]\label{lem.kappa.vanishing}
    The function $\kappa$ is upper-semicontinuous on $\cB$. The local invariant $\kappa_x$ vanishes if either
    \begin{itemize}
        \item $X_a^\flat \to X$ is \'etale over $x$, or
        \item $X_a$ is smooth above $x$. 
    \end{itemize}
\end{lemma}

The basic global bound is the following.

\begin{lemma}[{cf.\ \cite[Lemme 9.8.1]{CLII}}]
    If $\kappa(b) < \deg(D) - 2g + 2$, then $H^2(E_b) = 0$ and $\cB$ is smooth at $b$.
    \begin{proof}
        The proof proceeds essentially unchanged until the final calculation, where we pass to a splitting cover and then descend. By Serre duality and exactness of $W$-invariants ($=$ coinvariants for vector spaces), the obstruction space is dual to 
        $$
            H^0(X, F_L(\Omega^1_{X_a^\flat/k})).
        $$
        Applying the exact functor $F_L$ to the short exact sequence
        $$
            \begin{tikzcd}[ampersand replacement=\&]
                0 \ar[r] \& \pi^\star \Omega^1_{X/k} \ar[r] \& \Omega^1_{X_a^\flat/k} \ar[r] \& \Omega^1_{X_a^\flat/X} \ar[r] \& 0
            \end{tikzcd},
        $$
        we have an inclusion of vector bundles whose cokernel is torsion of length $\kappa(b)$, cf.\ \cite[Remarque 9.7.3]{CLII}.
        In the local ring at each such point, we have a containment of lattices with cokernel length $\kappa_x(b)$, so we need to contract by at most $\kappa_x(b)$ to reverse the containment. Summing over points, there exists an effective divisor $K$ of degree at most $\kappa(b)$ such that 
        $$
            F_L(\Omega^1_{X_a^\flat/k}) \subset F_L(\pi^\star \Omega^1_{X/k})(K).
        $$
        By the projection formula and the injection $\cH^{-1}(L_b) \hookrightarrow f^\star \Omega^1_{\ft_D/X}$, it suffices as in \textit{loc.\ cit.} to show
        $$
            H^0(X_a^\flat, f^\star \Omega^1_{\ft_D/X} \otimes \pi^\star \Omega^1_{X/k}(K)) = 0.
        $$
        The nonsplit version of the computation of this first tensor factor is
        $$
            f^\star \Omega^1_{\ft_D/X} = \pi^\star(\ft^\vee \otimes \O_X(-D)),
        $$
        where $\ft^\vee$ is a vector bundle over $X$, so the tensor bundle is 
        $$
            \pi^\star(\ft^\vee \otimes \Omega^1_{X/k}(K - D)).
        $$
        Pulling back to a finite \'etale Galois splitting cover renders $\ft^\vee$ trivial, and by our assumption on $\kappa(b)$ the right factor has negative degree. So every summand has negative degree on every connected component of the pulled-back normalization, so there are no global sections. \'Etale descent for global sections then implies the vanishing on the original $X_a^\flat$.
    \end{proof}
\end{lemma}

The second result we need to extract is the infinitesimal dimension estimate.

\begin{definition}
    Define the exact functor
    $$
        \cL \mapsto F_\Omega(\cL) := (\pi_\star(\cL \otimes \Omega^1_{X_a^\flat/k}))^W.
    $$
\end{definition}

\begin{lemma}[{cf.\ \cite[Lemme 9.9.1]{CLII}}]\label{lem.dim.est.cB}
    If $H^2(E_b) = 0$ and $\deg(D) > 2g$, then
    $$
        \dim T_b(\cB) \le \dim \cA_G - \delta_a.
    $$
    \begin{proof}
        The relative cotangent complexes (which here all have length 2) fit into a distinguished triangle
        $$
            \begin{tikzcd}[ampersand replacement=\&]
                \nu^\star L_{X_a/\ft_D} \ar[r] \& L_{X_a^\flat/\ft_D} \ar[r] \& \left[ \nu^\star \Omega^1_{X_a/X} \to \Omega^1_{X_a^\flat/X}\right] \ar[r, "+1"] \& \;
            \end{tikzcd}.
        $$
        We apply the exact functor $F_{\Omega}$ and take Euler characteristics. 
        \begin{itemize}
            \item The first term $F_{\Omega}(\nu^\star L_{X_a/\ft_D})$ reduces by Lemma \ref{lem.H-1.LX_a} to its $\cH^{-1}$, which can be computed to be
            $$
                \cH^{-1}(F_{\Omega}(\nu^\star L_{X_a/\ft_D})) = a^\star \Omega^1_{\fc_D/X} \otimes (\pi_\star \Omega^1_{X_a^\flat/k})^W.
            $$
            By tame descent as in \cite[Lemme 9.11.1(2)]{CLII}, we have $(\pi_\star \Omega^1_{X_a^\flat/k})^W \cong \Omega^1_{X/k}$, so the Euler characteristic is 
            $$
                - \chi(X, a^\star \Omega^1_{\fc_D/X} \otimes \omega_X)
            $$
            or by Serre duality
            $$
                \chi(X, a^\star T_{\fc_D/X}).
            $$
            Our degree hypothesis $\deg(D) > 2g$ ensures the vanishing of $H^1$, so that this Euler characteristic is just its $H^0$. But $\cA_G$ is the space of global sections of $\fc_D$ so its tangent space is the space of global sections of the tangent bundle. This first term therefore gives $\dim \cA_G$.

            \item The middle term is
            $$
                F_\Omega(L_{X_a^\flat/\ft_D}) = (\pi_\star(L_{X_a^\flat/\ft_D} \otimes \Omega^1_{X_a^\flat/k}))^W;
            $$
            duality commutes with proper pushforward along $\pi : X_a^{\flat} \to X$, so
            $$
                \RsHom_X(\pi_\star L_b^\vee, \Omega^1_{X/k}) \cong \pi_\star \RsHom_{X_a^\flat}(L_b^\vee, \Omega^1_{X_a^\flat/k}) \cong \pi_\star(L_b \otimes \Omega^1_{X_a^\flat/k}).
            $$
            This then passes to the $W$-invariants, in particular identifies $\cE_b = (\pi_\star L_b^\vee)^W$ with the dual of $F_\Omega(L_{X_a^\flat/\ft_D})$. Taking $R\Gamma(X, -)$ and applying Serre duality gives $R\Gamma(X, F_{\Omega}(L_{X_a^\flat/\ft_D})) = E_b^\vee[-1]$. Its Euler characteristic is therefore, by the vanishing of $H^0(E_b)$ as in Remark \ref{rmk.H0Eb,vanishes} and the vanishing of $H^2$ by assumption, equal to 
            $$
                \dim H^1(E_b) = \dim T_b \cB.
            $$

            \item The right term contributes two torsion terms
            $$
                 \chi(F_\Omega(\Omega^1_{X_a^\flat/X})) - \chi(F_\Omega(\nu^\star \Omega^1_{X_a/X})) = \deg F_\Omega(\Omega^1_{X_a^\flat/X})-\deg F_\Omega(\nu^\star \Omega^1_{X_a/X}).
            $$
            Write $p$ for the second degree. The short exact sequence of relative differentials for $\ft_D \to \fc_D \to X$ gives, after pullback along $f : X_a^\flat \to \ft_D$ and taking $F_\Omega$,
            $$
                p = \deg F_{\Omega}(f^\star \Omega^1_{\ft_D/X}) - \deg F_\Omega(f^\star \chi^\star \Omega^1_{\fc_D/X}),
            $$
            just as in \cite[Lemme 9.10.2]{CLII}. The first term splits across the decomposition $f^\star \Omega^1_{\ft_D/X} = \ft^\vee \otimes \pi^\star \O_X(-D)$ as in \cite[Lemme 9.10.3]{CLII}, and splitting $\Omega^1_{X_a^\flat}$ across the short exact sequence of relative differentials for $X_a^\flat \to X \to k$, we deduce just as in \cite[Lemme 9.10.4]{CLII} from Ng\^{o}'s formula (Proposition  \ref{prop.norm.form.delta}) that
            $$
                p = \delta + \deg((\ft^\vee \otimes \pi_\star \Omega^1_{X^\flat_a/X})^W).
            $$
            There is a mild discrepancy between the $\ft^\vee$ which appears here and the $\ft$ in the formula above. Our hypotheses on characteristic are enough to ensure the equality of these two formulae: in \textit{loc.\ cit.} they fix a global nondegenerate $W$-invariant bilinear form, but we need not do so globally: \'etale locally the group splits and we can choose on the split Cartan a Weyl-invariant bilinear form $\ft_0 \cong \ft_0^\vee$ and conclude that the local lengths agree. But of course the value of these local lengths descends from an \'etale cover, though the pairing itself may not.
            
            \item It remains to show that 
            $$
                \deg((\ft^\vee \otimes \pi_\star \Omega^1_{X^\flat_a/X})^W) - \deg ((\pi_\star(\Omega^1_{X_a^\flat/X} \otimes \Omega^1_{X^\flat_a/k}))^W) \ge 0.
            $$
            As in \cite[Lemme 9.11.1(3)-(4)]{CLII}, we have at a geometric point $\wtilde{x} \in X_a^\flat$ over $x \in X$ with inertia $W_{\wtilde{x}} \subset W$ that
            $$
                \on{length}_x((\ft^\vee \otimes \pi_\star \Omega^1_{X^\flat_a/X})^W) = \dim(\ft) - \dim(\ft^{W_{\wtilde{x}}}),
            $$
            and
            $$
                \on{length}_x((\pi_\star(\Omega^1_{X_a^\flat/X} \otimes \Omega^1_{X^\flat_a/k}))^W) = \begin{cases}
                    0 &\text{ if $X_a^\flat$ is \'etale over $x$}\\
                    1 &\text{ otherwise}
                \end{cases}.
            $$
            These are both \'etale-local computations and reduce to the split case. At \'etale points both vanish, and at ramified points, since the Weyl group acts faithfully, $\ft^{W_x} \neq \ft$, our desired inequality holds.
        \end{itemize}
        This concludes the proof.
    \end{proof}
\end{lemma}

We collect finally the handful of local calculations from \cite{CLII} that we will require in the application. All of the assertions reduce to the split case: for any geometric point $\ol{x} \to X$, after base change to the strict Henselization $\Spec \O^{\sh}_{X, \ol{x}}$ the twisting datum defining $G/X$ is trivial, and the Cartan, Chevalley map, Weyl group scheme, discriminant, and cameral cover identify with those for the corresponding split group. Hence, all assertions concerning geometric properties of the completed local cameral cover, discriminant multiplicities, stabilizers, smoothness, and the local invariants $\delta_x, \kappa_x$ may be checked after this base change.

\begin{lemma}[{\cite[Lemmes 9.13.2, 9.14.1]{CLII}}]\label{lem.kappa.local.computations}
\noindent
    \begin{enumerate} 
    \item For any $b \in \cB$ living over $a \in \cA_G$ and any closed point $x \in X$, one has the inequality
    $$
        \kappa_x(b) \le \frac{|W_G|}{2} [k(x) : k] v_x(\fD_G(a)).
    $$
    \item For any $b \in \cB$ living over $a \in \cA_G$ and any closed point $x \in X$, if the intersection multiplicity of $\fD_G$ and $a(X)$ at $x$ is $1$, then the cameral cover $X_a$ is smooth over $x$. In particular by Lemma \ref{lem.kappa.vanishing} $\kappa_x(b)$ vanishes for any $b \in \cB$ over $a$.
    \end{enumerate}
\end{lemma}

We can now state and prove Chaudouard--Laumon's general criterion for goodness of a parameter of the Hitchin base. Let $H/X$ be an arbitrary endoscopic group for $G$, including, for example, a semistandard Levi or an endoscopic group for a Levi. Recall the discriminant and resultant divisors $\fD_H$, $\fR_{H}^G$ from \S \ref{subsect.twist.chevalley}. Write $\Delta = (\fD_H + \fR_{H}^G)_{\red}$. 

\begin{theorem}[{\cite[Th\'eor\`eme 9.1.3]{CLII}}]\label{thm.suff.criterion.for.goodness}
    For $H/X$ as above, let $a_H \in \ctA_{H}(k)$ and write $a = \chi^{H}_G(a_H) \in \cA_G$ for its image under the appropriate Chevalley map. Assume that there is a finite set of points $S \subset |X|$ such that the following two conditions hold.
    \begin{enumerate}
        \item For every point $x \not\in S$, the curve $a_H(X)$ either avoids $\Delta$ or intersects exactly one of $\fD_H$, $\fR_{H}^G$ transversally at a smooth point.
        \item One has
        $$
            \frac{|W_G|}{2} \sum_{x \in S} [k(x) : k] v_x(\fD_G(a)) < \deg(D) - 2g + 2
        $$
        where $v_x(\fD_G(a))$ is the intersection multiplicity of $a(X)$ and $\fD_G$ at $x$. 
    \end{enumerate}
    Then $a$ lies in the good locus $\ctA_G^{\good}$. 
    \begin{proof}
        The proof is essentially the same as \textit{loc.\ cit.}, reducing the local computations to the split case as above. We outline the argument. The functions $\delta, \kappa$ are defined on $\ctA_G$ by pullback from the corresponding function on $\cA_G$; that is, they are independent of $t$. Since the good locus is defined as the \emph{interior} of the set of points satisfying 
        $$
            \codim_{\ctA}(a) \ge \delta_a
        $$
        we need to show that the inequality holds not just for $a$ but for any generization of it. It suffices to consider one-parameter generizations
        $$
            \beta : \Spec K[[u]] \to \cA_G
        $$
        for some field $K$ with special point $a$. Write $\eta$ for the generic point of the DVR and $s$ for its special point. We now calculate $\kappa(\beta_\eta)$. 
        
        Let $x_\eta$ be a point on the generic curve which specializes to $x \not\in S$. If $a_{H}(x) \not\in \Delta$ or $a_{H}(x)$ crosses $\fD_{H}$ transversely, then at worst the same is true for $x_\eta$. In the first case the cameral cover is \'etale and in the second the $G$-discriminant multiplicity is 1 and by Lemma \ref{lem.kappa.local.computations}.(2) the cameral cover is smooth. In either case $\kappa_{x_\eta} = 0$. 

        If, on the other hand, the special point intersects $\fR_H^G$ transversally, then (passing to a splitting \'etale cover or the strict Henselization) any chosen cameral lift $\wtilde{x}_s$ of the special point $x_s$ identifies a unique root hyperplane $\pm \alpha \subset \Phi_G \setminus \Phi_H$ such that $$\on{Stab}_{W_G}(\wtilde{x}_s) = \langle s_\alpha \rangle.$$ By generization, we have $|W_{\wtilde{x}_\eta}| \le 2$, and the ambient discriminant valuation is 0, 1, or 2. In the first two cases we are done. In the last case $v_{x_\eta}(\fD_G(\beta_\eta)) = 2$, for $\wtilde{x}_\eta'$ a point in the normalization $X^\flat_a$ over $x_\eta$, the Ng\^{o}--Bezrukavnikov formula gives
        $$
            2 = v_{x_\eta}(D_G(a_\eta)) = 2 \delta_{x_\eta} + \dim \ft - \dim \ft^{W_{\wtilde{x}_\eta'}}.
        $$
        Since $W_{\wtilde{x}_\eta'} \subseteq \langle s_\alpha \rangle$ and $\ft^{\langle s_\alpha \rangle}$ is a root hyperplane, we obtain a parity contradiction unless $W_{\wtilde{x}_\eta'}$ is trivial. In that case the normalization is \'etale over $x_\eta$, and again $\kappa_{x_\eta} = 0$ by Lemma \ref{lem.kappa.vanishing}. 

        So $\kappa(\beta_\eta)$ is supported entirely on points $x_\eta$ specializing to $x \in S$. By Lemma \ref{lem.kappa.local.computations}.(1),
        $$
            \kappa(\beta_\eta) = \sum_{x_\eta \rightsquigarrow S} \kappa_{x_\eta}(\beta_\eta) \le \frac{|W_G|}{2} \sum_{x_\eta \rightsquigarrow S} v_{x_\eta}(\fD_G(\beta_\eta)).
        $$
        Clearly we have 
        $$
            \sum_{x_\eta \rightsquigarrow S} v_{x_\eta}(\fD_G(\beta_\eta)) \le \sum_{x \in S} v_x(\fD_G(a))
        $$
        (in fact equality holds, though we don't need this), so by our hypothesis bound
        $$
            \kappa(\beta_\eta) < \deg(D) - 2g + 2.
        $$

        Since this bound holds for every generic deformation of $a$ and $\kappa$ is constructible, the locus of points satisfying this inequality on $\kappa$ contains a Zariski neighborhood $U$ of $a$. For an irreducible component $Z$ of the $\delta = r$ stratum meeting $U$, take its generic point and the corresponding normalized cameral curve $b \in \cB$. Then 
        $$
            \kappa(b) < \deg(D) - 2g + 2,
        $$
        so the deformation obstruction vanishes and $\dim Z \le \dim_b \cB \le \dim \cA_G - r$ by the dimension estimate of Lemma \ref{lem.dim.est.cB}. Every $\delta$-stratum meeting $U$ has the required codimension, so $U \subseteq \cA_G^{\good}$ and $a \in \cA_G^{\good}$.
    \end{proof}
\end{theorem}

One immediate application of this theorem is a characterization of how the good locus intersects the various endoscopic loci, which we record now for use in \S \ref{sect.coh.thms} and \S \ref{sect.global.point.counts}.

\begin{lemma}[{cf.\ \cite[Lemme 14.3.1]{CLII}}]\label{lem.good.locus.compat.end}
    Work over algebraically closed $k$ and assume $\deg(D) > 2g$. Let $M \subset G$ be a semistandard Levi, and $(s, \rho_s)$ a pointed endoscopic datum for $M$ with endoscopic group $H$. Then under the composite closed embedding as in \S \ref{subsect.hitchin.base.endoscopic}
    $$
        i_{\rho_s} : \ctA_{H} \hookrightarrow \ctA_M \hookrightarrow \ctA_G,
    $$
    (all spaces with imposed $G$-regularity at $\infty$),
    \begin{enumerate}
        \item $i_{\rho_s}^{-1}(\ctA_G^{\good}) \subseteq \ctA_{H}^{\good}$ and
        \item $i_{\rho_s}^{-1}(\ctA_G^{\good}) \cap \ctA_{H}^{H\el} \subseteq \ctA_{H}^{H\el}$ is dense, in particular nonempty. 
    \end{enumerate}
    \begin{proof}
        The first part is proven exactly as in \textit{loc.\ cit.}, i.e., follows from the elementary computation of the dimension of the Hitchin base in Proposition \ref{prop.Hitchin.base.dim} and the Ng\^{o}--Bezrukavnikov formula for the invariant $\delta_a$.

        For the second part, we apply the previous theorem when $S = \emptyset$. A parameter $a_H \in \ctA_H$ has image which is $G$-good as long as it intersects $(\fD_H + \fR_{H,D}^G)_{\red}$ transversally, and in the case that $M = G$ a nonempty open subset do by \cite[Lemme 8.5.1]{ngofl}. When $M \subsetneq G$ the same Bertini proof (as in \cite[Proposition 4.7.1]{ngofl}) applies: at each point the condition of a bad intersection imposes two conditions, and $\deg(D) > 2g$ implies surjectivity of the map restricting a section of $\fc_{H,D}$ to the first jet at a point. See also \S \ref{subsect.globalize.param} for a similar argument. Intersecting with the $H$-elliptic, $G$-regular-at-$\infty$ locus gives a dense open by irreducibility and the nonemptiness from Proposition \ref{prop.ell.end.data}.
    \end{proof}
\end{lemma}

The second claim gives a negative answer, on our (larger) good locus, to the empty intersection possibility raised in \cite[\S 7.8.4]{ngofl}.

\chapter{Support and the cohomological WFL}\label{sect.coh.thms}

\section{The Hitchin fibration.}

We recall the setup of the Hitchin fibration we will apply. In particular, we change notation slightly from the previous. Let $X_0/k_0$ a smooth geometrically integral curve over a finite field $k_0 = \F_q$, and let $G_0 \to X_0$ a reductive group scheme with maximally split fiberwise-maximal torus $T_0 \to X_0$. Write $X$, $G$, $T$, etc. for the base change of these objects to $k = \ol{k}_0$. Assume now that $G$ has anisotropic center, i.e., $A_G^{\spl} = 1$, and when necessary that $\xi 
\in \a_T$ is in general position. This is sufficient for all our previous theorems to apply. 

Write $\tau : x \mapsto x^q$ for the (arithmetic) Frobenius automorphism of $k$ and the induced Frobenius of $X, G, T$, etc. Let $D$ an effective $\tau$-stable divisor on $X$ of degree $\deg(D) > 2g$, and assume $D = 2 D'$, so that the choice of $D'$ defines a Kostant--Hitchin section $\cA \to \cM$ as in \cite[\S 4.3.3]{bouthierprodform}. Pick $\infty \in X(k)$ lying over a rational point of $X_0$. Assume that the fiber of $G_0$ is split at $\infty$.

The Hitchin base $\ctA_G$, resp. the compatible Hitchin base $\bA_G$, is equipped with an action of $\tau$, under which the closed subschemes $\ctA_M$, resp. $\bA_M$, for semistandard $X$-Levis $M$ are stable, cf.\ Corollary \ref{cor.parab.levi.scalar.extensions}. Likewise $\ctM_G$ and $\bM_G$ have an action of $\tau$, for which the Hitchin morphism is equivariant. 

For any $\xi \in \a_T$, the $\xi$-stable substack $\bM_G^\xi$ is $\tau$-stable because the stability convex $C_m$ is insensitive to the extension $k/k_0$, cf.\ Remark \ref{rmk.stab.convex.stable.extensions}. Indeed, we define the truncated Hitchin fibration $f_0^{\xi} : \bM_{0, G}^{\xi} \to \bA_{0, G}$ over $k_0$, and upon base change to $k$ we recover the corresponding constructions over $k$. Note in particular that our stability combinatorics remain unchanged under this extension: parabolic reductions of Higgs bundles $m$ for $G_0$ are determined by the $M_0$-elliptic type of $f_0(m)$ in \S \ref{sect.Hitchin.fib}, and the combinatorics of relative position remain unchanged because the group $G_0$ is already split over $\infty$. We denote with the adjective ``geometric'' those constructions and objects which are defined over $k$, and ``arithmetic'' those which are defined over $k_0$. 

Write $d_{f}$ for the relative dimension of $f : \bM_G \to \bA_G$. This is well-defined: by \cite[Corollaire 4.16.4]{ngofl} the Hitchin morphism $f$ over $\cA^{\heartsuit}$ is of constant relative dimension (it is equal to the dimension of $\bP_a$, which is locally constant). By \cite[Corollaire 4.16.3]{ngofl} each Hitchin fiber is equidimensional, so the relative dimension remains the same constant after restricting to the compatible truncated Hitchin morphism $f^\xi$. We denote $d_{f^\xi} = d_f$. Write $d_{\bM^\xi}$ for the dimension of $\bM^\xi$ (equal to the dimension of $\bM$) and $d_{\bA}$ for the dimension of the base, so $d_{\bM^\xi} = d_{\bA} + d_{f^\xi}$.

\section{Decomposition and actions of the Tate module.}

Now we recall the essential cohomological machinery in our setting---this section is entirely an exposition of \cite[\S 7]{ngofl} and \cite[\S 10.1-10.3]{CLII}, with all results and arguments transferring directly. Many of the results below apply in the general setting of a proper morphism $f_0 : M_0 \to A_0$ of smooth $k_0$-schemes of constant fiber dimension equipped with a relative action of a smooth commutative group scheme $g_0 : P_0 \to A_0$ satisfying mild hypotheses, cf. \cite[\S 7]{ngofl} which works in this generality. In particular, the Chevalley decomposition of the fibers of $P_0$ induces a $\delta$ invariant as above, and one requires the good locus hypothesis
$$
    \codim_A(a) \ge \delta_a
$$
of Definition \ref{def.good.locus}. Of course, our Hitchin fibrations are maps of Deligne--Mumford stacks, not schemes, but the results apply anyways, either by translating the proofs directly or by rigidifying our torsors with trivializations at $\infty$, with the only caveat that the arguments below will involve $\delta_a^{\aff}$ instead of the original $\delta_a$. See in particular \cite[4.5.6, 4.12]{ngofl} for this point, though note that over the elliptic locus, in groups with anisotropic center, the two agree. Since our focus in this work is on the application to the Hitchin fibration, we have not usually bothered to state the results generally.

Now assume $\xi$ is in \emph{general position} in the sense of Definition \ref{def.gen.pos}. The geometric results of \S \ref{sect.finite.type}, \S \ref{sect.val.unique}, \S \ref{sect.val.exist}, and \S \ref{sect.Deligne.Mumford} imply that $f^\xi_0$ and $f^\xi$ are proper morphisms (of constant relative dimension $d_{f^{\xi}}$) of smooth Deligne--Mumford stacks. As noted above, the results of \cite{bbdg} extend to this setting (say, by choosing appropriate rigidifications, cf.\ \cite[\S 6.2, \S 7.8]{ngofl}), and in particular Deligne purity implies that over $k$ we have
$$
R f_\star^\xi \ol{\Q}_l[d_{\bM^{\xi}}] \cong \bigoplus_i \pHi(Rf_\star^\xi \ol{\Q}_l[d_{\bM^{\xi}}])[-i],
$$
with each graded piece splitting as a sum of intermediate extensions of local systems on locally closed subsets $Z_\alpha \hookrightarrow \bA$:
$$
    \pHi(Rf_\star^\xi \ol{\Q}_l[d_{\bM^{\xi}}]) = \bigoplus_{i_\alpha : Z_\alpha \hookrightarrow \bA} i_{\alpha,!\star} \cL_\alpha[\dim(Z_\alpha)].
$$
Switching the order of the summation, we can write the full relative cohomology complex as a sum over the Zariski points of $\bA$ of intermediate extensions of graded local systems. To be precise, 
\begin{equation}
    \pHi(Rf_\star^\xi \ol{\Q}_l[d_{\bM^\xi}]) = \bigoplus_{a \in \bA} j_{a,!\star} \cL_a^i[d_a],\label{eqn.socle.decomp}
\end{equation}
where $\cL_a^{\bullet}$ is a graded local system on solely the point $\{a\}$ with the property that it extends to some open in $\ol{\{a\}}$, and we have somewhat abusively written $j_{a, !\star}$ for the well-defined functor which performs intermediate extension to $\bA$ from any such open. The pieces $\cL_a^{n}$ are moreover pure of weight $n + d_{\bM^{\xi}} - d_a$. 

Define the \emph{support} of the sheaf $Rf_\star^\xi \ol{\Q}_l[d_{\bM^\xi}]$ to be the set of irreducible closed subsets appearing in this decomposition with nonzero local systems. Define the \emph{socle} to be the finite subset of Zariski points of $\bA$ which are the generic points of the support. In the decomposition of (\ref{eqn.socle.decomp}) it is the finite set of points at which $\cL_a^\bullet$ is nonzero. For each $\cL_a^{\bullet}$, define its \emph{amplitude} to be the range of degrees in which it is supported
$$
    \Amp(\cL_a^{\bullet}) = \max \{n : \cL_a^n \neq 0\} - \min \{n : \cL_a^n \neq 0\}.
$$

A Poincar\'e duality argument due to Goresky--MacPherson gives a general bound on the cohomological amplitude and the codimension of the socle.

\begin{proposition}[{{\cite[Th\'eor\`eme 7.3.1]{ngofl} or \cite[Lemme 10.2.1]{CLII}}}]\label{prop.amp.codim.bound.basic}
    Let $a \in \bA$ a point in the socle of $Rf_\star^\xi \ol{\Q}_l[d_{\bM^\xi}]$. Then 
    $$
        \Amp(\cL_a^\bullet) \le 2(d_{f^\xi} - d_{\bA} + d_a),
    $$
    and
    $$
        \codim_{\bA}(a) \le d_{f^\xi},
    $$
    where $d_a$ is the dimension of the closure of $a$.
\end{proposition}

Ng\^{o}'s support theorem is built on alternative bounds for the amplitude and the codimension, found by analyzing an action of the Tate module associated to the $k_0$-Picard stack $g_0 : \bP_0^0 \to \bA_0$ (the fiberwise identity component)
$$
    \Tate \bP_0^0 := R^{2d_f - 1} g_{0, !} \ol{\Q}_l (d_f)
$$
on the factors appearing in the decomposition (\ref{eqn.socle.decomp}). We recall the construction of this action and its properties. Writing $K^\bullet_0 = R f_{0,\star}^\xi \ol{\Q}_l$, there is a natural cup product
$$
    \Tate \bP_0^0 \otimes K^\bullet_0 \to K^\bullet_0[-1]
$$
defined by convolution along the correspondence $\bM_0^{\xi} \xleftarrow{a} \bP_0^0 \times_{\bA_0} \bM_0^{\xi} \xrightarrow{\pi_2} \bM_0^{\xi}$, i.e.
$$
    (v \in \Tate \bP_0^0, x \in R^k f_{0,\star}^{\xi} \ol{\Q}_l) \mapsto \pi_{2,!}(\pi_1^\star v \cup a^\star x),
$$
where $a$ is the action map (which preserves $\xi$-truncation since $\bP_0^0 \subseteq \bP_0^1$). Here the trace map $\pi_{2,!}$ lowers degree by $2d_f$ and cancels the twist $(d_f)$, giving an output in $R^{k-1} f_{0, \star} \ol{\Q}_l$. Indeed, this is more naturally interpreted as the $(-1)$-part of the action of
$$
    \Lambda^\bullet := R g_{0, !}\ol{\Q}_l[2d_f](d_f)
$$
on $K^\bullet_0$ by the same formula. Extending the usual argument for families of abelian varieties, there is likewise a multiplication $\Lambda^\bullet \otimes \Lambda^\bullet \to \Lambda^\bullet$ which equips its cohomology with a graded-commutative algebra structure, and in particular a morphism
$$
    \bigwedge^i \Tate \bP_0^0 \to H^{-i}(\Lambda^\bullet),
$$
which is actually an isomorphism. In fact, $\Lambda^\bullet$ is formal, in the sense that the complex splits canonically as a sum of its cohomology, cf. \cite[\S 7.4.4]{ngofl}. We obtain an isomorphism
$$
    \Lambda^\bullet = \bigoplus_{i \ge 0} \bigwedge^i \Tate \bP_0^0[i] = \bigwedge^\bullet \Tate \bP_0^0
$$
compatible with the multiplication. 

The action of the Tate module (which we view as a sheaf supported in degree 0) induces an action on perverse cohomology which decreases degree by 1 (\cite[\S 7.4.6]{ngofl}):
$$
    \Tate \bP_0^0 \otimes \pHi(K^\bullet_0) \to  \pH^{i-1}(K^\bullet_0). 
$$
We consider just the diagonal components of the action of the Tate module for the decomposition (\ref{eqn.socle.decomp}). It follows from the projection formula that for a point $a$ in the socle, $\Tate \bP_{0, a}^0$ acts on each $\cL_a^{\bullet}$ as an endomorphism of degree -1.

For a point $a$ in the socle, one has a canonical sequence of rational Tate modules into an affine and abelian part, built from the Chevalley decomposition. Let us briefly re-emphasize that we are working with Deligne--Mumford stacks, which do not have canonical Chevalley decompositions, so this warrants some explanation. 

Let $a \in \bA$ a geometric point. Note first that under our anisotropic center hypothesis, the inertia of $\bP_{0,a}^0$ is of dimension $\dim H^0(X, J_a) = 0$, which implies that the coarse quotient $\bP_{0,a}^0 \to (\bP_{0,a}^0)^\#$ induces an isomorphism on rational Tate modules.

Now choose any presentation of $\bP_{0,a}^0$, as a cokernel of smooth fiberwise-connected commutative group schemes (as in \S \ref{subsect.aff.delta.inv})
$$
    Q_{-1,a} \to Q_{0,a} \to \bP_{0,a}^0 \to 1,
$$
with $Q_{-1}$ affine, and take the corresponding sequence of rational Tate modules
$$
    0 \to \Tate Q_{-1, a} \to \Tate Q_{0, a} \to \Tate \bP_{0,a}^0 \to 0.
$$

Consider the Chevalley decomposition of $Q_{0,a}$:
$$
    1 \to Q_{0, a}^{\aff} \to Q_{0, a} \to Q^{\ab}_{0, a} \to 1.
$$
As in the discussion before Definition \ref{def.aff.delta.inv}, $Q_{-1,a}$ has trivial image in $Q_{0, a}^{\ab}$ and we obtain a map $\bP_{0,a}^0 \to Q_{0, a}^{\ab}$. We can descend this to the coarse quotients by taking $(-)^{\#}$; the resulting kernel is $$Q_{0,a}^{\aff}/\im(Q_{-1,a} \to Q_{0,a}) \cong K_a^\#.$$ 
Taking Tate modules and identifying the modules for $\bP_{0,a}^0$ with its coarse quotient, we obtain a filtration
\begin{equation}\label{eqn.tate.module.geo}
    0 \to \Tate K_a^\# \to \Tate \bP_{a}^0 \to \Tate Q_{0, a}^{\ab} \to 0.
\end{equation}
The filtration is independent of choice of presentation; we denote the subobject $(\Tate \bP_{a}^0)^{\aff}$ and the quotient $(\Tate \bP_{a}^0)^{\ab}$ to emphasize this independence.

We can extend (\ref{eqn.tate.module.geo}) to any Zariski point of $\bA$ by choosing a geometric point lying over it and descending, in fact. This is slightly subtle because the Chevalley decomposition need only exist after a finite radiciel extension: the geometrically maximal normal affine subgroup need not be defined over the residue field, as in \cite{pseudoabvars} or \cite[Exp.\ XVII, App.\ III, Cor.\ to Prop.\ C.5.1]{SGAIII}. However, since finite radiciel extensions are universal homeomorphisms, the Tate modules and the filtration above descend along such an extension.

The weights for these components of the Tate module are easy to compute using a choice of presentation (and working over whatever extension of $a$ necessary for the objects to be defined): since $Q^{\ab}_{0,a}$ is an abelian variety, its Tate module is dual to $H^1(Q_{0,a}^{\ab}, \Q_l)$, so is pure of weight -1. The affine part $K_a^\#$, on the other hand, has $l$-torsion only in its toric part $T_a$, whose rational Tate module is isomorphic to $X_\star(T_a) \otimes_{\Z} \Q_l(1)$ with geometric Frobenius acting by finite order (with eigenvalues roots of unity) on the cocharacter lattice and by $q^{-1}$ on the $\Q_l(1)$ factor. Overall, $\Tate K_a^{\#}$ is pure of weight -2. But we know that the action of $\Tate \bP_{0,a}^0$ takes $\cL_a^n$, which is pure of weight $(n + d_{\bM^{\xi}} - d_a)$, to $\cL_{a}^{n-1}$, which is pure of weight one less. So the affine part of the Tate module acts \textit{trivially} on the factors $\cL_a^n$. Extending multiplicatively to $\Lambda_a^\bullet$, we have the following.

\begin{proposition}\label{prop.tate.action.La}
    For each $a$ in the socle of the $\xi$-truncated Hitchin fibration, the corresponding factor $\cL_a^\bullet$ in the decomposition of (\ref{eqn.socle.decomp}) has an action of $\bigwedge^\bullet (\Tate \bP_{0,a}^0)^{\ab}$, the homology of an abelian variety. 
\end{proposition}

In the case that $\ol{a}$ is a geometric point lying over a closed point $a$, the Frobenius weights on the affine and abelian parts of (\ref{eqn.tate.module.geo}) are (-2) and (-1), as above, so there is a unique Frobenius-equivariant splitting. This gives an action of $\bigwedge^\bullet (\Tate \bP_{0,a}^0)^{\ab}$ on the full truncated fiber-cohomology, whose induced action on the perverse graded pieces is the action above.

\section{Freeness of the cohomology over the Tate module.}

The most technical part of Ng\^{o}'s work on the support theorem is his analysis of this action. In the remainder of this section, we will show that this action makes $\cL_a^\bullet$ into a free $\bigwedge^\bullet (\Tate \bP_{0,a}^0)^{\ab}$-module. The crucial step is to to show a freeness result for the full cohomology of geometric closed fibers. 

\begin{proposition}[{{\cite[Lemme 10.3.3, Th\'eor\`eme 10.5.1]{CLII}}}]\label{prop.tate.action.geo.fibers.free}
    For every geometric point $\ol{a} \in \bA_0$ living over a closed point, the Frobenius-splitting action of $\bigwedge^\bullet (\Tate \bP_{\ol{a}}^0)^{\ab}$ on $H^\bullet(\bM^\xi_{ \ol{a}}, \ol{\Q}_l)$ makes it into a free module.
\end{proposition}

We will deduce the proposition from the following fundamental freeness result for \'etale-locally trivial families over an abelian variety.

\begin{lemma}[{{\cite[Lemme 10.3.3]{CLII}}}]\label{lem.free.ab.var}
    Let $Y$ a scheme or Deligne--Mumford stack over algebraically closed $k$, equipped with a proper $k$-morphism $f : Y \to B$ to an abelian variety. Assume that there exists an isogeny $\lambda : B' \to B$ and an action of $B'$ on $Y$ such that
    $$
        f(b\cdot y) = \lambda(b) + f(y).
    $$
    Then the action of $H_\bullet(B', \ol{\Q}_l) = \bigwedge^\bullet \Tate B' \cong \bigwedge^\bullet \Tate B$ on $H^\bullet(Y, \ol{\Q}_l)$ makes it into a free module.
\end{lemma}


Following Chaudouard--Laumon, we apply this lemma not to the truncated Hitchin fiber, but to a contracted product homeomorphic to it. We can use a slightly different construction (using $G$-affine Springer fibers rather than those for $M$) because the product formula of Ng\^{o} for anisotropic Hitchin fibers (\cite[4.15.1]{ngofl}) has been extended in \cite{bouthierprodform} to all of $\cA^{\heartsuit}$. We recall the theorem. Given $v \in X$ and $a \in \cA$, we denote $\cM_{a,v}$ for the \textit{affine Springer fiber} of $a$ at the place $v$, an ind-scheme defined as in \cite[\S 3.2]{ngofl}. We denote also $\cP_{a,v}$ for the local Picard stack, parameterizing $J_a$-torsors on $\widehat{\O}_v$ equipped with trivializations on the punctured disk, as in \cite[\S 3.3]{ngofl}. This acts on the affine Springer fiber with a dense orbit, just as in the global setting.

\begin{theorem}[{{\cite[Theorem 4.3.8]{bouthierprodform}}}]\label{thm.bouthier.prod.gluing}
    For $a \in \cA^{\heartsuit}(k)$, we write $j : U_a := a^{-1}(\fc_{D}^{\reg}) \hookrightarrow X$ for the open on which $a$ is regular semisimple. Then Beauville--Laszlo gluing defines a universal homeomorphism of locally finite type algebraic $k$-stacks 
    $$
        \on{BL}_a : \cM_a^{\on{prod}} := \prod_{v \in X \setminus U_a} \cM_{a,v}^{\red} \times^{\prod_{v \in X\setminus U_a} \cP_{a,v}^{\red}} \cP_a \longrightarrow \cM_a.
    $$
\end{theorem}
It is evident from the construction that this map is equivariant for the action of $\cP_a$ on the right factor and on $\cM_a$. Pulling the map back to a homeomorphism $\bM_a^{\on{prod}} \to \bM_a$ for $(a,t)$ in the compatible locus, we can denote $\bM_a^{\on{prod}, \xi} \to \bM_a^\xi$ for its restriction to the inverse image of the $\xi$-stable locus. 

In fact, under the hypotheses of \cite[Theorem 4.3.8(a)]{bouthierprodform}, which are weaker than ours, the gluing map is moreover \emph{finite} and \emph{representable by schemes} on any geometric fiber over $\cA^{\heartsuit}$. From our geometric results (Theorems \ref{thm.bdd}, \ref{thm.sep}, \ref{thm.val.crit.exist}, \ref{thm.xi.locus.dm}) $\bM_a^{\mathrm{prod},\xi}$ is therefore also proper Deligne–Mumford and canonically has the same étale cohomology as $\bM_a^\xi$, which is the fact we require.

\begin{proof}[Proof of finite representability in Theorem \ref{thm.bouthier.prod.gluing}]
    We follow the steps of the proof of \cite[Theorem 4.3.8(b)]{bouthierprodform}, with different arguments. Consider the factorization
    $$
        \cM_a^{\on{prod}} \xrightarrow{\;r\;} \prod_{v \in X \setminus U_a} \cM_{a,v} \times^{\prod_{v \in X\setminus U_a} \cP_{a,v}} \cP_a \xrightarrow{\;g\;} \cM_a,
    $$
    as in \S 4.3.6. In \textit{loc.\ cit.} they prove that (1) $g$ is fully faithful on arbitrary algebras and that (2) the reduction map $r$ is an equivalence on reduced test algebras, so that $\on{BL}_a$ is fully faithful on reduced test algebras. 
    
    We claim $r$ is in fact faithful on arbitrary test objects. A morphism $((m_v)_v, p) \to ((m_v')_v, p')$ in the (prestack) associated product consists of a reduced family of local modifications 
    $$
        h = (h_v)_v \in \prod_{v \in X\setminus U_a} \cP_{a,v}^{\red}
    $$
    such that $h_v. m_v = m_v'$ and an isomorphism $\beta : p \otimes h \cong p'$. The non-reduced version admits the corresponding characterization, and the image $r(h, \beta)$ is the same pair, where $h$ is now viewed in the product of nonreduced local Picard stacks. Since
    $$
        \prod_{v \in X\setminus U_a} \cP_{a,v}^{\red} \hookrightarrow \prod_{v \in X\setminus U_a} \cP_{a,v}
    $$
    is a monomorphism, the prestack $r$ is injective on arrows. Equality of arrows can be checked locally, so this passes to the stackification. Since the composite $\on{BL}_a$ is a \emph{faithful} map of algebraic stacks, it is representable by algebraic spaces (\cite[Tag. 04Y5]{stacks-project}). 

    Next we show the composite gluing map is separated. Since it is representable by algebraic spaces, its diagonal is a \emph{monomorphism} representable by schemes and locally of finite type (\cite[Tag. 04YQ]{stacks-project}). Choosing an arbitrary test scheme $T$, a point
    $$
        T \to \cM^{\on{prod}}_a \times_{\cM_a} \cM^{\on{prod}}_a
    $$
    consists of two product families together with an isomorphism between their glued global images. Write $\Delta_T$ for the pullback of $\Delta_{\on{BL}_a}$ along this point, a monomorphism. After further pullback to $T^{\red}$, full faithfulness of the map on reduced test objects implies that isomorphism came from a unique isomorphism in $\cM^{\on{prod}}_a$, i.e., that $\Delta_{T^{\red}} = \Delta_{T} \times_T T^{\red}$ has a section. But it is a monomorphism, hence an isomorphism. In particular, $\Delta_{T^{\red}}$ and hence $\Delta_T$ are homeomorphisms. So $\Delta_{\on{BL}_a}$ is a universal homeomorphism, in particular universally closed. A universally closed, locally of finite type, monomorphism is a closed immersion (\cite[Tag. 05W8]{stacks-project}), so $\on{BL}_a$ is separated.

    Finally, as in \textit{loc.\ cit.}, $\on{BL}_a$ is a separated universal homeomorphism of algebraic spaces, so representable by schemes and integral after \cite[Corollary 5.22]{rydhdesc}. But it is also locally of finite type, hence finite.
\end{proof}

\begin{proof}[Proof of Proposition \ref{prop.tate.action.geo.fibers.free}]

Let $S = X \backslash U_a$, and consider the full lft N\'eron extension of $J_a$ across $S$,
$$
    J_a^{(S)} := j_\star j^\star J_a.
$$
This imposes no integrality condition on $S$. Write $\mathcal P_a^{(S)}=\Bun_{J_a^{(S)}}$. Theorem \ref{thm.weil.unif} applied to this generically reductive group scheme (over algebraically closed constants $k$) gives a uniformization of the $k$-points
$$
    \Bun_{J_a^{(S)}}(k) \cong [J_a(F) \backslash J_a(\A) / \prod_{v \not\in S} J_a(\O_v) \prod_{v \in S} J_a(F_v)].
$$

One advantage of the product formula for Hitchin fibers is that there is a natural Beauville--Lazlo patching map $\cM_a^{\on{prod}} \to \Bun_{J_a^{(S)}}$, which in general may not descend to $\cM_a$. One moreover has a short exact sequence 
$$
\begin{tikzcd} 0 \ar[r] & J_a^{\flat} \ar[r] & J_a^{(S)} \ar[r] & \bigoplus_{v \in S} i_{v,\star} \ul{\Lambda_v} \ar[r] & 0
\end{tikzcd},
$$
where $\Lambda_v := J_a(F_v)/J_a^{\flat}(\O_v)$ is free of finite rank, induced by the universal property of the lft N\'eron model from their identification over $U_a$. The induced pushforward of torsors gives a surjective map $\cP_a^{\flat} \to \Bun_{J_a^{(S)}}$ which induces an equivalence of Picard stacks
$$
    \Bun_{J_a^{(S)}} \cong \left[\cP_a^{\flat}/\ul{\Lambda_S}\right],
$$
where $\ul{\Lambda_S} := \bigoplus_{v \in S} \ul{\Lambda_v}$ as a $k$-scheme. On $k$-points the quotient map corresponds to the natural quotient from the maximal bounded subgroup $J_a(F_v)^{\flat}$ to $J_a(F_v)$:
$$
[J_a(F) \backslash J_a(\A) / \prod_{v \not\in S} J_a(\O_v) \prod_{v \in S} J_a(F_v)^{\flat}] \to [J_a(F) \backslash J_a(\A) / \prod_{v \not\in S} J_a(\O_v) \prod_{v \in S} J_a(F_v)].
$$
This quotient is not an isogeny---even after pullback to $\bA$ and restriction to the degree 0 part the quotient map has additional automorphisms. However, after passing to coarse neutral components, it will be.

Denote the coarse space of a Picard stack by a superscript $\#$. The quotient presentation identifies the coarse fppf quotient $\mathcal P_a^{(S),\#}$ (cf.\ \S \ref{sect.good.locus}) with the quotient of $\mathcal P_a^{\flat,\#}$ by the image of $\ul{\Lambda_S}$. The subgroup of $\ul{\Lambda_S}$ mapping into the identity component is finitely generated, and its image there is finite because $k=\overline{\mathbf F}_q$. Consequently
\[
B_a:=(\mathcal P_a^{(S),\#})^0
\]
is an abelian variety, and $(\mathcal P_a^{\flat,\#})^0\to B_a$ is an isogeny. The existing presentation $Q_{0,a}\to\bP_a^0$ therefore induces an isogeny $Q_{0,a}^{\ab}\to B_a$, identifying their rational Tate modules. Applying \cite[Proposition 7.5.3]{ngofl} to $Q_{0,a}$ (which also uses $k = \ol{\mathbb{F}}_q$), and composing with an isogeny in the opposite direction, gives a homomorphism $s_a:B_a\to \bP_a^0$ whose composite with the natural map $\bP_a^0\to B_a$ is multiplication by some $N>0$.

Compose the product-formula map with the coarse quotient to obtain
\[
f_a: \bM_a^{\mathrm{prod},\xi}\longrightarrow
\mathcal P_a^{(S),\#}.
\]
Through $s_a$, the group $B_a$ acts on this product model, preserving the truncation since $\bP_a^0\subseteq\bP_a^1$, and
\[
f_a(s_a(b)\cdot m)=[N]b+f_a(m).
\]
Properness implies that the image meets only finitely many components of $\cP_a^{(S),\#}$. Translating each such component to $B_a$, we may apply Lemma \ref{lem.free.ab.var} with $\lambda=[N]_{B_a}$; taking their finite direct sum proves the required freeness. Under the induced identification $V_\ell B_a\simeq(V_\ell\bP_a^0)^{\mathrm{ab}}$, the normalized map $N^{-1}V_\ell(s_a)$ is the canonical Frobenius-compatible splitting,\footnote{Indeed, the splitting descends to some finite extension, so commutes with some power of the Frobenius $\tau^m$. But the weights of $\tau^m$ remain distinct, as in the discussion prior to Proposition \ref{prop.tate.action.La}, so the $\tau^m$-equivariant splitting is unique, and must equal to the $\tau$-equivariant one.} so the resulting action agrees with the one in the proposition.

\end{proof}

Ng\^{o} next transfers this freeness result to the action of the Tate module on $\cL_a^\bullet$. This relies on the crucial fact that on the good locus of the Hitchin fibration, the Tate module $\Tate \bP_0^{0}$ is \emph{polarizable}, in the sense that it admits an alternating form
$$
    \Tate \bP_0^0 \times \Tate \bP_0^0 \to \ol{\Q}_l(1)
$$
which is on each geometric fiber identically zero on the affine part and perfect on the abelian part. This is \cite[Proposition 4.12.1]{ngofl}. The argument is an elegant reduction to the case of $\G_m$ using a map $$J_a \to \pi_{a, \star}'(\T \times X'_a)$$ defined as in \S \ref{subsect.descr.galois}, where $\T$ is a split torus and we have chosen a splitting $\Gamma$-cover $X' \to X$ for the group and written $\pi_a' : X'_a := X' \times_X X_a \to X$ for the composite $\W \rtimes \Gamma$. The map induces an extension of scalars from $J_a$-torsors on $X$ to $\T$-torsors on $X'_a$. This extension is a group homomorphism with respect to the group structure, so sends $l^n$-torsion points of $\cP^{\heartsuit}$ to $l^n$-torsion points of $\Bun_T \cong \Pic_{X_a'} \otimes X_\star(\T)$. Tensoring a symmetric invariant form on the lattice part with the classical Weil pairing on the Jacobian gives a form with the desired properties. 

Ng\^{o} proves that polarizability, as well as the freeness on relative cohomology of Proposition \ref{prop.tate.action.geo.fibers.free}, imply freeness of the Tate module action on each local system summand $\cL_a^\bullet$. Its proof constitutes \cite[\S 7.5-7.7]{ngofl}, and would take us too far into the weeds to expound. 

\begin{proposition}[{{\cite[Proposition 10.3.1]{CLII}, \cite[7.4.10]{ngofl}}}]
    For any point $a$ in the socle of the $\xi$-truncated Hitchin fibration, the fiber $\cL_a^\bullet$ is a free module for $\bigwedge^\bullet (\Tate \bP^0_{0,a})^{\ab}$.
\end{proposition}

\begin{corollary}\label{cor.amp.codim.bound.hard}
    Given a point $a$ in the socle of the truncated Hitchin fibration, one has the inequalities
    $$
        \Amp(\cL_a^{\bullet}) \ge 2(d_f - \delta_a^{\aff})
    $$
    and
    $$
        \delta_a^{\aff} \ge \codim_{\bA}(a).
    $$
    \begin{proof}
        The rank of the abelian part of the Tate module is 2 times the dimension of the abelian variety (after some rigidification), which is $d_f - \delta_a^{\aff}$ by Definition \ref{def.aff.delta.inv}. Since the fiber is nonzero, choosing a homogeneous free generator and taking its span under $\bigwedge^\bullet (\Tate \bP^0_{0,a})^{\ab}$ gives nonzero elements spread across at least $2(d_f - \delta_a^{\aff})$ degrees.

        Composing this lower bound on the amplitude with the upper bound of Proposition \ref{prop.amp.codim.bound.basic} gives
        $$
            2(d_{f^\xi} - \delta_a^{\aff}) \le 2(d_{f^\xi} - d_{\bA} + d_a), 
        $$
        from which we deduce the codimension bound.
    \end{proof}
\end{corollary}

In a slogan, Ng\^{o}'s support theorem states that the supports are large relative to the dimension of their affine Picard part.

\begin{theorem}[{{cf.\ \cite[Th\'eor\`eme 10.5.1]{CLII}}}]\label{thm.no.new.supports}
    According to our running hypothesis, let $G$ be reductive with anisotropic center. The socle (the set of generic points of support components) of the restriction of $Rf_\star^\xi \ol{\Q}_{l} [d_{\bM^\xi}]$ to $\bA^{\good}$ is entirely contained within the elliptic locus $\bA_{\Gell} \cap \bA^{\good}$.
    \begin{proof}
        Let $a \in \bA^{\good}$ be a point in the socle. Since we are on the good locus of Definition \ref{def.good.locus}, we have 
        $$
            \delta_a \le \codim_{\bA}(a).
        $$
        The computation of \ref{prop.delta.invs.comparison} implies that if $a$ is not $G$-elliptic
        $$
            \delta_a^{\aff} < \delta_a \le \codim_{\bA}(a).
        $$
        But the freeness of the action of the Tate module implies the bound of Corollary \ref{cor.amp.codim.bound.hard},
        $$
            \delta_a^{\aff} \ge \codim_{\bA}(a),
        $$
        a contradiction. 
    \end{proof}
\end{theorem}

This implies that the perverse cohomology of the truncated Hitchin fibration is determined by its restriction to the elliptic locus, as below. 

\begin{corollary}\label{cor.coh.intermed.ext.from.ell.locus}
    For each $i$, we have on $\bA^{\good}$ a canonical isomorphism 
    $$
        \pHi(Rf_\star^\xi \ol{\Q}_l[d_{\bM^\xi}]) \cong j_{!\star} j^\star \pHi(Rf_\star^\xi \ol{\Q}_l[d_{\bM^\xi}]),
    $$
    where $j : \bA^{\Gell} \cap \bA^{\good} \hookrightarrow \bA^{\good}$ is the open immersion.
    \begin{proof}
        By the decomposition theorem and stability of perversity under intermediate extension and restriction to an open, both are semisimple perverse sheaves which agree on those summands whose supports intersect the elliptic locus. But the theorem implies all summands do, so the canonical map induces isomorphisms on all summands.
    \end{proof}
\end{corollary}

\section{Geometric stabilization.}

With the requisite machinery built, and in particular Theorem \ref{thm.no.new.supports} indicating that on the good locus all cohomological supports intersect the $G$-elliptic locus, one can now, transfer Chaudouard and Laumon's cohomological weighted fundamental lemma to our setting. Write $f_G^{\xi} : \bM_G^{\xi} \to \bA_{G}$ for the compatible $\xi$-stable Hitchin fibration associated to $G$. With the aid of $\xi$-truncation, we can extend Ng\^{o}'s cohomological result, which we recall now, from the elliptic locus to $\bA_G^{\good}$. 

For each elliptic endoscopic datum $(s, \rho_s, \rho_s \to \rho_G)$, write $i_{\rho_s} : \ctA_{G_s}^{G_s\el} \hookrightarrow \bA_G^{G\el}$ for the closed embedding of Proposition \ref{prop.hitchin.base.functorial.endoscopic.groups}, restricted to the elliptic loci.\footnote{For the reader who has been carefully tracking questions of compatibility in the sense of the discrepancy between $\bA$ and $\ctA$, note that compatibility is a vacuous condition on the elliptic locus, so $\ctA_{G}^{G\el} = \bA_G^{G\el} \subset \bA_G$.}

\begin{theorem}[{{The cohomological Fundamental Lemma, \cite[6.4.1, 6.4.2]{ngofl}}}]\label{thm.ngo.coh.fl}
    Let $s \in \widehat{T}_{\tors}$. On $\bAe{G}$, over $k = \ol{k}_0$, one has an isomorphism of fine isotypic parts in the sense of Theorem \ref{thm.s.decomp}:
    $$
        \pH^\bullet(R f^{G\el}_{\star} \ol{\Q}_l)_{s} \cong \bigoplus_{(s, \rho_s)} (i_{\rho_s})_\star \pH^\bullet(R f^{G_{\rho_s}'\el}_{G_{\rho_s}',\star} \ol{\Q}_l)_1[-2r_s](-r_s),
    $$
    where $r_s = \codim_{\ctA_G}(\ctA_{G_{\rho_s}'})$, the shift is in the perverse cohomological degree, and the sum ranges over (geometric) elliptic endoscopic data $(s, \rho_s)$ with fixed component $s$, for $G_{\rho_s}'$ the corresponding elliptic endoscopic group. Equivalently, for each $(s, \rho_s)$ elliptic,
    $$
        (i_{\rho_s})^\star \pH^\bullet(R f^{G\el}_{\star} \ol{\Q}_l)_{s} \cong \pH^{\bullet}(R f^{G_{\rho_s}'\el}_{G_{\rho_s}',\star} \ol{\Q}_l)_1[-2r_s](-r_s).
    $$
    In fact, if a given geometric endoscopic datum extends to an arithmetic one over finite coefficient field $k_0$, then the second isomorphism above is defined over $k_0$, up to taking semisimplifications of both sides.
\end{theorem}

\begin{remark}
By the elementary bound of \cite[Proposition 6.2.3]{ngofl}, the support of the $s$-isotypic part on the elliptic locus lies within $\ctA_{s} = \bigcup_{\rho_{s}} i_{\rho_s}(\ctA_{G_{\rho_s}})$, the disjoint union of all the endoscopic loci for geometric endoscopic data extending $s$. Pulling back or pushing forward between one of these components justifies the equivalence between the two geometric statements. 
\end{remark}

\begin{remark}
    We remark that under our running hypothesis that $A_G^{\spl} = 1$, $$Z_{\widehat{G}}^{\Gamma, 0} = 1$$ and Ng\^{o}'s decomposition is exactly the fine decomposition of \ref{thm.s.decomp}. For groups with anisotropic centers, $G$-ellipticity in our sense coincides with anisotropicity. For groups with nontrivial split central tori, Ng\^{o}'s anisotropic locus is empty, so we have lost no generality.

    Moreover, by Proposition \ref{prop.ell.end.torsion}, every geometric elliptic endoscopic datum occurring here already has $s$ of finite order. 
\end{remark}

Now pick once and for all a semistandard Levi $M$. Recall that the open $U_M \subset \bA_G$ (inside the compatible locus of Definition \ref{def.comp.loc}) is defined to be the union of the $L$-elliptic loci for $L \supseteq M$. It is Frobenius-stable by Corollary \ref{cor.parab.levi.scalar.extensions} and its consequence that the constructions of \S \ref{sect.Hitchin.base} are compatible with extending scalars. Write $f^{\xi, \Mpar}$ for the restriction of $f^{\xi}$ to this open. Since the $\delta$-invariant is likewise compatible with extending scalars, the good locus $\ctA^{\good}$ of Definition \ref{def.good.locus} is Frobenius-stable. So we can define over $k_0$ the open embedding
$$
    j_{\on{ell}} : \bAe{G}^{\good} \hookrightarrow U_M^{\good}.
$$

We can now give a description of each coarse summand in terms of fine pieces on the $G$-elliptic locus.

\begin{theorem}\label{thm.coh.stab.geo}
    Let $G$ quasisplit reductive with $A_G^{\spl} = 1$ and $\xi \in \a_T$ in general position. Let $\ol{s}_M \in (\widehat{T}/Z^0_{\widehat{M}})_{\tors}$. Then as perverse sheaves on $U_M^{\good}$, the corresponding degree $i$ coarse $\ol{s}_M$-isotypic component of Theorem \ref{thm.s.decomp} satisfies
    $$
        \pHi(R f_\star^{\xi, \Mpar} \ol{\Q}_l)_{\ol{s}_M} \cong \bigoplus_{\substack{s \in \widehat{T}_{\tors} \\ s \mapsto \ol{s}_M}} \bigoplus_{(s, \rho_{s})} (j_{\on{ell}})_{!\star} (i_{\rho_{s}})_\star \pHishift(R f_{G_{\rho_s}',\star}^{G_{\rho_s}'\el} \ol{\Q}_l)_1(-r_s),
    $$
    where
    \begin{itemize} 
    \item the sum $s$ ranges over the fiber of $\ol{s}_M$ under $\widehat{T}_{\tors} \to (\widehat{T}/Z^0_{\widehat{M}})_{\tors}$,
    \item $(s,\rho_{s})$ ranges over geometric elliptic pointed endoscopic data for $G$ with fixed $s$.
    \item $\pHi(R f_{G_{\rho_s}',\star}^{G_{\rho_s}'\el} \ol{\Q}_l)_1$ is the fine $1$-isotypic part in the decomposition of Theorem \ref{thm.s.decomp} applied to the (quasisplit reductive) $X$-endoscopic $G_{\rho_s}'$, and
    \item $r_s = \codim_{\ctA_G}(\ctA_{G_{\rho_s}'})$.
    \end{itemize}
\begin{proof}
    By Corollary \ref{cor.coh.intermed.ext.from.ell.locus}, we know that on the good locus we have 
    $$
        \pHi(R f_\star^{\xi} \ol{\Q}_l) \cong (j_{ell})_{!\star} (j_{ell})^\star \pHi(R f_\star^{\xi} \ol{\Q}_l). 
    $$
    Once restricted to $U_M^{\good}$, this splits across the coarse $s$-decomposition of Theorem \ref{thm.s.decomp}, so writing $K_{\ol{s}_M} := \pHi(R f_\star^{\xi, \Mpar} \ol{\Q}_l)_{\ol{s}_M}|_{U_M^{\good}}$, it decomposes in terms of the \emph{fine} summands on the $G$-elliptic locus:
    $$
        K_{\ol{s}_M} \cong (j_{\on{ell}})_{!\star} (j_{\on{ell}})^\star K_{\ol{s}_M} \cong (j_{\on{ell}})_{!\star} \bigoplus_{\substack{s \in \widehat{T}_{\tors} \\ s \mapsto \ol{s}_M}} \pHi(R f^{G\el}_{G, \star} \ol{\Q}_l)_{s}.
    $$
    By the geometric Fundamental Lemma \ref{thm.ngo.coh.fl}  of Ng\^{o}, the fine $s$-isotypic factor is a sum of stable factors over elliptic endoscopic data $(s, \rho_{s})$ for that fixed $s$
    $$
        K_{\ol{s}_M} \cong (j_{\on{ell}})_{!\star} \bigoplus_{\substack{s \in \widehat{T}_{\tors} \\ s \mapsto \ol{s}_M}} \bigoplus_{(s, \rho_s)} (i_{\rho_{s}})_\star \pHishift(R f_{G_{\rho_s},\star}^{G_{\rho_s}\el} \ol{\Q}_l)_1(-r_s)
    $$
    where the index sets are as desired.

                \end{proof}
\end{theorem}

In particular, Ng\^{o}'s full-support result for these summands on the elliptic locus extends.
\begin{corollary}\label{cor.full.supp}
    Under the hypotheses of Theorem \ref{thm.coh.stab.geo}, fix a geometric elliptic pointed endoscopic datum $(s, \rho_s)$ for $G$, with endoscopic group $G'_s$, indexing a summand of the coarse $\ol{s}_M$-component. Every nonzero geometric simple constituent $P$ of the $(s, \rho_s)$-summand has closed support
    $$
        i_{\rho_s}(\ctA_{G_s'}) \cap U_M^{\good}.
    $$
    \begin{proof}
        By Proposition \ref{prop.ell.end.data} and Lemma \ref{lem.good.locus.compat.end}, $V_{G'_s} := i_{\rho_s}^{-1}(\ctA_{G}^{G\el} \cap \ctA_G^{\good}) \subseteq \ctA_{G'_s}^{G'_s\el} \cap \ctA_{G'_s}^{\good}$
        as a nonempty open. Theorem \ref{thm.coh.stab.geo} implies that every such $P$ arises by closed pushforward and intermediate extension from a simple constituent of the stable piece $$\pHishift(Rf_{G'_{s}, \star}^{G'_{s}\el} \ol{\Q}_l)_1$$
        on $V_{G'_s}$. But as in the proof of \cite[Theorem 7.8.3]{ngofl}, such stable summands have full support $V_{G'_s}$: $G'_s$-goodness and \cite[Proposition 7.2.3]{ngofl} identify them as having the same support as nonzero summands of the ordinary top-degree stable cohomology, which is constant of rank 1 by \cite[Proposition 6.5.1]{ngofl}. The support of $P$ is therefore the closure of $i_{\rho_s}(V_{G'_s})$ in $U_M^{\good}$, which is the claimed intersection.
    \end{proof}
\end{corollary}

We end with a few remarks on the nature of Theorem \ref{thm.coh.stab.geo}.
\begin{remark}
    This theorem is stated geometrically, over $k$ rather than $k_0$, and indeed not all the elliptic endoscopic data $(s, \rho_{s})$ which arise are arithmetic for a given $k_0$. Correspondingly, this decomposition is not fully defined over $k_0$. However, each term is defined over an appropriate finite extension, and there are finitely many such terms by Proposition \ref{prop.ell.end.finite}, so the entire decomposition is defined over an appropriate finite extension of $k_0$. Our application will ultimately need an arithmetic version, making use of \cite[6.4.2]{ngofl} rather than \cite[6.4.1]{ngofl}.
\end{remark}

\begin{remark}\label{rmk.split.end.datum.restriction.CL}
    When $G = \G \times_k X$ is split (constant) $\rho_G = 1$ is the trivial cocycle, and so there is for each $s \in \widehat{T}$ a canonical choice $\rho_s = 1$ (and so $G_s' = \G_s'$) of split endoscopic group. In \cite{CLII}, for each $s$ the authors restrict to the open locus which contains only the endoscopic data associated to that split endoscopic group by manually removing the (finitely many by \ref{prop.ell.end.finite}) closed loci $\ctA_{\rho_s}$ associated to endoscopic data with $\rho_s \neq 1$. The result is that their statement of analogous to the theorem above contains no sum over $\rho_s$. For a quasi-split group there is no such canonical choice, and indeed applications require a WFL for which neither $G$ nor the endoscopic groups which appear need be split. 
\end{remark}

\begin{remark}\label{rmk.discrep.cWFL.vs.Arthur}
    Arthur's statement of the WFL differs from this theorem in a few ways:
\begin{enumerate}
    \item The WFL contains only a single sum over $s$, with a single endoscopic datum for $G$ induced from an endoscopic datum for $M$ and the choice of $s$, whereas the theorem above is a double sum.
    \item The WFL allows only $s \in s_M Z_{\widehat{M}}^\Gamma$, whereas the theorem has no Galois-fixed requirement.
    \item The WFL allows translates by the full (Galois-fixed) center, whereas the theorem ranges only over $s \in s_M Z^0_{\widehat{M}}$.
\end{enumerate}
\end{remark}

All of the discrepancies in these remarks, except the bookkeeping between $Z_{\widehat{M}}^\Gamma$ and $(Z_{\widehat{M}}^0)^\Gamma$, will naturally resolve when we specify a full endoscopic datum for $M$, rather than just $s_M$. We explain this in the next section.

\section{Arthur's induced endoscopy.}\label{subsect.ind.endoscopy}

The Weighted Fundamental Lemma depends on an endoscopic datum $\fM = (M', \ldots)$ for a Levi $M$, and sums over a set of (representatives of) endoscopic data $G'_s$ for $G$ which are compatible with it (and in particular contain $M'$ as a Levi) as defined in \cite[\S 4]{arthurenddatalevi}. One version of this construction was recalled in the Introduction. In order to give a statement of our cohomological theorem that more closely matches the WFL, we reinterpret this construction in terms of a variant of Ng\^{o}'s formalism for endoscopic data. Reinterpreting our cohomological theorem in terms of this formalism, we then deduce in the next chapter a global version of the WFL in the stalk over a point on the Hitchin base which is $M'$-elliptic.

For the sake of the general construction and definitions, we briefly relax our hypothesis of anisotropic center and let $G$ arbitrary quasisplit reductive. We fix a semistandard Levi $M \subset G$. Let $\fM' = (s_M, \rho_{s_M})$ be a pointed elliptic endoscopic datum for $M$, defining quasisplit endoscopic group $M' \to X$. We now construct a collection of induced endoscopic data $(s, \rho_s)$ for $G$, a task mildly complicated by the fact that there need be no direct map 
$$
    \widehat{M} \rtimes \Out(\widehat{M}) \to \widehat{G} \rtimes \Out(\widehat{G})
$$
to induce a map $\pi_0^M(s_M) \to \pi_0^G(s)$ that can be used to push forward the torsors. The cleanest workaround seems to be to pass through a variant of the construction of Definition \ref{def.end.datum} relative to a common splitting cover, as we now explain.

Fix $\Gamma$ some finite quotient of the fundamental group large enough to split both $G$ and $M'$. Then $\widehat{G}$ has a $\Gamma$-action preserving $\widehat{M}$, and so a corresponding inclusion
$$
    \widehat{M} \rtimes \Gamma \subseteq \widehat{G} \rtimes \Gamma.
$$
Consider the variant of the connected-\'etale short exact sequence from (\ref{eqn.end.ses}) which occurs in $\widehat{M} \rtimes \Gamma$ instead of $\widehat{M} \rtimes \Out(\widehat{M})$. That is, $s_M \in \widehat{T}$, defining the connected-\'etale short exact sequence
\[
\begin{tikzcd}[ampersand replacement=\&]
        1 \ar[r] \& \widehat{M}'_{s_M} \ar[r] \& Z_{\widehat{M} \rtimes \Gamma}(s_M) \ar[r] \& \pi_0^{M,\Gamma}(s_M) \ar[r] \& 1
\end{tikzcd}.
\]
This sequence is simply related to (\ref{eqn.end.ses}): the subgroup is the same, and the quotient is given by 
$$
    \pi_0^{M, \Gamma}(s_M) = \pi_0^M(s_M) \times_{\Out(\widehat{M})} \Gamma.
$$
The endoscopic torsor $\rho_{s_M} : \pi_1(X, \infty) \to \pi_0^M(s_M)$ then lifts naturally to a $\Gamma$-refinement
$$
    \rho_{s_M}^{\Gamma} : \pi_1(X, \infty) \to \pi_0^{M,\Gamma}(s_M),
$$
whose pushforward along $$\pi_0^{M,\Gamma}(s_M) \to \Gamma \to \Out(\G)$$ is $\rho_G$ and whose pushforward along $\pi_0^{M,\Gamma}(s_M) \to \pi_0^M(s_M) \to \Out(\widehat{M}')$ defines $M' \to X$ as a form of the split Langlands dual $\M'$ of $\widehat{M}'_{s_M}$.

Each $s \in s_M Z_{\widehat{M}}^{\Gamma}$ induces a corresponding pointed endoscopic datum $(s, \rho_s)$ for $G$, defining endoscopic group $G'_s \to X$, as follows. We have a $\Gamma$-SES for $G$ just as above for $M$:
\[
\begin{tikzcd}[ampersand replacement=\&]
        1 \ar[r] \& \widehat{G}'_s \ar[r] \& Z_{\widehat{G} \rtimes \Gamma}(s) \ar[r] \& \pi_0^{G,\Gamma}(s) \ar[r] \& 1
\end{tikzcd}.
\]
Choosing a representative $s = s_M z$ with $z \in Z_{\widehat{M}}^\Gamma \subseteq Z(\widehat{M} \rtimes \Gamma)$, it is clear from the centrality of $z$ that
$$
    Z_{\widehat{M} \rtimes \Gamma}(s_M) \subseteq Z_{\widehat{G} \rtimes \Gamma}(s),
$$
so we get an induced map \begin{equation} q_s^{\Gamma} : \pi_0^{M,\Gamma}(s_M) \to \pi_0^{G,\Gamma}(s),\label{eqn.endo.map.pi0}\end{equation}. In fact it is injective: any cocharacter whose centralizer in $\widehat{G}$ is $\widehat{M}$ has its centralizer in $\widehat{G}'_s$ equal to $\widehat{M}'$, as in the proof of Proposition \ref{prop.end.gp.sub.root.system}, which implies the kernel is trivial. 

The pushforward
$$
    \rho_s^{\Gamma} := (q_s^\Gamma)_\star \rho_{s_M}^{\Gamma} : \pi_1(X, \infty) \to \pi_0^{G,\Gamma}(s) 
$$
gives a $\Gamma$-endoscopic torsor for $G$ in the obvious sense. Writing $\rho_s$ for its pushforward along $\pi_0^{G,\Gamma}(s) \to \pi_0^{G}(s)$, it is immediate from the construction that $(s, \rho_s)$ is a pointed endoscopic datum for $G$. 

\begin{definition}\label{def.induced.end.datum}
    Given (pointed) elliptic endoscopic datum $\fM'=(s_M, \rho_{s_M})$ for $M$ and $s \in s_M Z_{\widehat{M}}^{\Gamma}$, define the \textit{induced (pointed) endoscopic datum for $G$} to be the pair $\fG'_s = (s, \rho_s)$ constructed above.
    
    We denote $G'_s \to X$ the corresponding endoscopic group for $G$. Correspondingly, we equip $\widehat{G'_s}$ always with the (endoscopic) outer $\Gamma$-action induced by conjugation along $\rho_s^{\Gamma}$, \emph{not} the restriction of the $\Gamma$-action on $\widehat{G}$. This then lifts to actual automorphisms by a choice of pinning.
    
    Denote \emph{Arthur's endoscopic parameter set} $\cE_{\fM'}(G)$ to be the set of such $$s \in s_M Z_{\widehat{M}}^{\Gamma}/Z_{\widehat{G}}^{\Gamma}$$ such that the induced datum is elliptic for $G$, in the sense of Definition \ref{def.ell.end.datum}. In particular, by Proposition \ref{prop.ell.end.finite} all such endoscopic groups have anisotropic center if $G$ does. 

    It is useful to define two refinements of this set. Define the \emph{coarse parameter set} $\cE_{\fM'}^c(G)$ the set of
    $$
        s \in s_M (Z_{\widehat{M}}^0)^{\Gamma}/(Z_{\widehat{G}}^0)^{\Gamma},
    $$
    respectively the \emph{fine parameter set} $\cE_{\fM'}^f(G)$ the set of
    $$
        s \in s_M (Z_{\widehat{M}}^0)^{\Gamma}/Z_{\widehat{G}}^{\Gamma,0},
    $$
    in both cases ranging only over $s$ such that the induced datum is elliptic. Note that when the group has anisotropic center, $Z_{\widehat{G}}^{\Gamma,0} = 1$. We call a class in the above sets an \emph{endoscopic parameter}.

    If $\fM' = 1_M = (1, \rho_M)$ is the trivial endoscopic datum for $M$ defining $M$ itself, then denote $\cE_{M}(G) := \cE_{\fM'}(G)$ and $\cE_{M}^f(G), \cE_M^c(G)$ likewise.
\end{definition}

\begin{remark}\label{rmk.ind.endo.compatibilities}
    The construction above satisfies a few obvious compatibilities, whose verifications are immediate:
    \begin{itemize}
        \item Given a map of Galois groups $\Gamma' \to \Gamma$, we have
        $$
            \pi_0^{G, \Gamma'}(s) \cong \pi_0^{G, \Gamma}(s) \times_{\Gamma} \Gamma'
        $$
        and likewise for $M$. In particular, the intrinsic induced datum $(s, \rho_s)$ is independent of the choice of cover. 

        \item The induced endoscopic group and endoscopic $\Gamma$-centralizer depend only on the Arthur class $s \in s_M Z_{\widehat{M}}^\Gamma/Z_{\widehat{G}}^\Gamma$. 
        \item If $M \subseteq L \subseteq G$ is a chain of semistandard Levi containments, then the construction above is functorial. 
    \end{itemize}

    One recovers the usual formulation of endoscopy by defining $\eM \subseteq Z_{\widehat{M} \rtimes \Gamma}(s_M)$ to be the inverse image of $\Gamma \to \pi_0^{M,\Gamma}(s_M)$. It then fits into a short exact sequence
    $$
        \begin{tikzcd}
            1 \ar[r] & \widehat{M}' \ar[r] & \eM \ar[r] & \Gamma \ar[r] & 1
        \end{tikzcd},
    $$
    so $(M', \eM, s_M)$ is an elliptic endoscopic datum in the same sense as the introduction (though globally, rather than locally). Performing the same construction on the induced datum will give the usual formula $\eG_s = \widehat{G}_s' \eM$.
\end{remark}

\begin{remark}\label{rmk.ind.data.arithmetic}
The constructions and definition above are equally valid whether one works geometrically over $k$, or arithmetically. In other words, if the endoscopic datum $\fM'$ was defined in terms of some finite quotient of the arithmetic $\pi_1(X_0, \infty)$, i.e., is defined over finite constant field $k_0$, then the induced endoscopic data for $G$ will all be defined over the same constant field.

As will be useful in \S \ref{sect.global.point.counts} below, we observe that if \textit{both} $G$ and the datum $\fM'$ is split at the rational point $\infty$ (i.e.\ if the endoscopic group $M'$ is), then the induced data $\fG'_s$ are likewise. In particular, by the argument of Remark \ref{rmk.arith.geo.data.correspond} they define individual Frobenius-stable loci inside the Hitchin base.
\end{remark}

\begin{remark}\label{rmk.ell.end.data.finite}
    For arbitrary reductive $G$ and elliptic endoscopic datum $\fM'$ for $M \subset G$, the set $\cE_{\fM'}(G)$ is finite: by consideration of sub-root systems, there are at most finitely many possibilities for $\widehat{G_s'} = Z_{\widehat{G}}(s)^0$, for each, the induced outer $\Gamma$-action is fixed, and the ellipticity condition implies finitude of the quotient $Z_{\widehat{G}_s'}^\Gamma/Z_{\widehat{G}}^\Gamma$ which depends only on $\widehat{G'_s}$, so at most finitely many classes of $s$ generate each $\widehat{G}_s'$. (The torsor $\rho_s$ is then uniquely determined as the pushforward.)
\end{remark}

\begin{remark}\label{rmk.EMG.vs.EMG0}
    There are natural surjective maps
    $$
        \cE_{\fM'}^f(G) \to \cE_{\fM'}^c(G) \to \cE_{\fM'}(G).
    $$
    (The latter is surjective by Arthur's lemma
    $$
        Z_{\widehat{M}}^\Gamma = Z_{\widehat{G}}^\Gamma (Z_{\widehat{M}}^\Gamma)^0
    $$
    cf.\ \cite[Lemma 1.1]{arthurtransferwoi}.) The first map has fibers which are torsors under the dual residual torsion group of Proposition \ref{prop.dual.Y1.fibers}
    $$
        \frac{(Z_{\widehat{G}}^0)^\Gamma}{Z_{\widehat{G}}^{\Gamma,0}} = \widehat{\Xi}_G.
    $$
    The second map has fibers which are torsors under the kernel
    $$
        \frac{(Z^0_{\widehat{M}})^\Gamma \cap Z_{\widehat{G}}^\Gamma}{(Z_{\widehat{G}}^0)^\Gamma}.
    $$
    These discrepancies will create the expected scalar factors between our cohomological Theorem \ref{thm.coh.wfl} indexed by $\cE^f_{\fM'}(G)$ and the global Weighted Fundamental Lemma indexed by $\cE_{\fM'}(G)$. 
\end{remark}

All the induced endoscopic $G'_s$ contain $M'$ over $X$ as a Levi.

\begin{proposition}\label{prop.ind.G'.cont.M'.Levi}
    We have a natural identification of root data
    $$
        \Phi(\widehat{M'}, \widehat{T}) = \Phi(\widehat{G'_s}, \widehat{T}) \cap \Phi(\widehat{M}, \widehat{T})
    $$
    inside $\Phi(\widehat{G}, \widehat{T})$. Moreover, fix a Galois-fixed cocharacter $\widehat{\lambda} \in X_\star(\widehat{T})^\Gamma$ (for the original action of $\Gamma$ on $\widehat{G}$) whose centralizer is $\widehat{M} \subset \widehat{G}$. Then, relative to an appropriate pinning for $\widehat{G'_s}$, the endoscopic action fixes $\widehat{\lambda}$ and extends the endoscopic action on $\widehat{M}'$. In particular, $M'$ is a semistandard $X$-Levi in each induced $G'_s$, canonical up to such a choice of adapted pinning.
    \begin{proof}
        Write $s = s_M z$ for $z \in Z_{\widehat{M}}^\Gamma$. By definition $\widehat{G'_s} = Z_{\widehat{G}}(s)^0$, so 
        $$
            \widehat{\alpha} \in \Phi(\widehat{G'_s}, \widehat{T}) \Longleftrightarrow \widehat{\alpha}(s) = 1.
        $$
        But for any roots $\widehat{\alpha}$ of $\widehat{M}$, we have $\widehat{\alpha}(z) = 1$, so this last condition is equivalent to $\widehat{\alpha}(s_M) = 1$, which defines $\Phi(\widehat{M'}, \widehat{T})$. This gives the first statement. 

        By Proposition \ref{prop.characterizing.dual.Levis}, we can find $\widehat{\lambda} \in X_\star(\widehat{T})^{\Gamma}$ such that $Z_{\widehat{G}}(\widehat{\lambda}) = \widehat{M}$. The root calculation above shows that
        $$
            Z_{\widehat{G'_s}}(\widehat{\lambda}) = \widehat{M}'.
        $$

        Now consider the positive set of roots defined by
        $$
            \Phi^+_{\widehat{G'_s}} := \{\widehat{\alpha} : \widehat{\alpha}(\widehat{\lambda}) > 0\} \sqcup \Phi^+_{\widehat{M'}},
        $$
        the latter relative to the fixed choice of pinning for $\widehat{M'}$. By the first claim of the proposition, this defines a Borel for $\widehat{G'_s}$ contained in $P_{\widehat{G'_s}}(\widehat{\lambda})$. By \cite[Lemme 10.1]{ngoendo} (cf.\ \S \ref{subsect.hitchin.base.endoscopic}), there exists a splitting 
        \begin{equation}
        \pi_0^{M,\Gamma}(s_M) \to Z_{\W_M \rtimes \Gamma}(s_M)
        \end{equation}
        whose image are those lifts preserving $\Phi^+_{\widehat{M'}}$, and likewise a splitting
        \begin{equation}
            \pi_0^{G, \Gamma}(s) \to Z_{\W_G \rtimes \Gamma}(s)\label{eqn.end.splitting.Weyl}
        \end{equation}
        preserving $\Phi^+_{\widehat{G'_s}}$. We can use this to give a refinement of the pushforward $(q_s^\Gamma)_\star$ constructed above, as follows. 
        
        The image of $\rho_{s_M}^\Gamma$ along this section preserves $\widehat{\lambda}$ because both $W_M$ and the original $\Gamma$-action did. The image of this section also centralizes $z \in Z_{\widehat{M}}^\Gamma$, so it maps into $\pi_0^{G, \Gamma}(s)$. This composition immediately induces the pushforward $(q_s^\Gamma)_\star$, but it explicitly maps to representatives in $Z_{W_G \rtimes \Gamma}(s)$ which preserve $\Phi_{\widehat{G'_s}}^+$, that is, ones in the image of the splitting (\ref{eqn.end.splitting.Weyl}) for $G$. Pushing back down the splitting, the endoscopic torsor $\Gamma \to \pi_{0}^{G,\Gamma}(s)$ also preserves $\widehat{\lambda}$, i.e., $\widehat{\lambda}$ is also fixed under the endoscopic $\Gamma$-action on $\widehat{G'_s}$. Its parabolic and Levi factor are therefore likewise stable, and 
        so by the dual correspondence at the beginning of \S \ref{subsect.dual.gps} they correspond to an $X$-parabolic and Levi factor inside $G'_s$.
    \end{proof}
\end{proposition}

\section{The cohomological WFL.}\label{subsect.coh.wfl}

Now place ourselves in the arithmetic setting from the top of the chapter. In particular, $G_0$ has anisotropic center and is split at the chosen rational point $\infty$. Assume also that $\fM'$ is an arithmetic pointed elliptic datum with group $M'_0$ which is also split at $\infty$, as in Remark \ref{rmk.ind.data.arithmetic}. By Remark \ref{rmk.arith.geo.data.correspond}, the datum $\fM'$ remains elliptic geometrically. By Proposition \ref{prop.ell.end.torsion} its coarse parameter 
$$
    \ol{s}_M \in \widehat{T}/Z_{\widehat{M}}^0
$$
defines an isotypic component of the coarse $s$-decomposition (\ref{thm.s.decomp}). Since $G$ has anisotropic center, the same remark shows that every $s \in \cE_{\fM'}^f(G)$ has finite order. 

We now seek to isolate a summand within the coarse factor $$K_{\ol{s}_M}^i := \pHi(Rf_\star^{\xi, \Mpar} \ol{\Q}_l)_{\ol{s}_M}|_{U^{\good}_M}$$ which corresponds to $\rho_{s_M}$. We would like to define it to be the sum of those pieces which have support containing $\iota_{\rho_{s_M}}(\ctA_{M'})$. Geometrically, this is the right definition, and by Remark \ref{rmk.arith.geo.data.correspond} this locus in the Hitchin base is Frobenius-stable, but in order to obtain a definition which descends to $k_0$, we need to take the semisimplification. 

\begin{definition}\label{def.fM'.cpt}
    Given an arithmetic (defined over $k_0$) elliptic pointed endoscopic datum $\fM' = (s_M, \rho_{s_M})$ for $M$ which is split at $\infty$, define the \emph{coarse semisimple $\fM'$-component} to be the maximal summand of the semisimplified coarse component
    $$
        (K_{\ol{s}_M}^i)^{\ss}_{\fM'} \subseteq (K_{\ol{s}_M}^i)^{\ss}
    $$
    consisting of those simple components whose geometric support contains the generic point $\eta_{\fM'} \in \iota_{\rho_{s_M}}(\ctA_{M'}^{M'\el})$.   
\end{definition}

\begin{remark}
    Indeed, $\eta_{\fM'} \in U_M^{\good}$, so this is a well-founded condition:
    \begin{itemize} 
    \item by Proposition \ref{prop.ell.end.data} the $M'$-elliptic locus is vacuously compatible and so embeds into the compatible $M$-stratum of $G$, i.e., into $U_M$, and then
    \item by Lemma \ref{lem.good.locus.compat.end}.(2), $\eta_{\fM'} \in U_M^{\good}$.
    \end{itemize}
\end{remark}

Moreover, the next result characterizes the summands which appear (out of those parameterized by elliptic endoscopic data $(s, \rho_s)$ for $G$ after Theorem \ref{thm.coh.stab.geo}) as exactly those represented by $s \in \cE_{\fM'}^f(G)$ and induced $\rho_s$. As in Remark \ref{rmk.ind.data.arithmetic}, each induced geometric endoscopic datum is Frobenius-stable and split at $\infty$, and so each such summand is Frobenius-stable.

\begin{lemma}\label{lem.M.G.comp.of.supports}
    Fix $i$, let $\fG' = (s, \rho_s)$ be a geometric elliptic endoscopic datum for $G$ which appears in the decomposition of Theorem \ref{thm.coh.stab.geo} for $K^i_{\ol{s}_M}$, and let $P$ be a nonzero geometric simple constituent of that summand. Then $P$ is a component of $(K_{\ol{s}_M}^i)^{\ss}_{\fM'}$ if and only if $s \in \cE_{\fM'}^f(G)$ and $\fG'$ is induced from $\fM'$ in the sense of Definition \ref{def.induced.end.datum}.

    In particular, if $a \in U_M^{\good}$ is the image of a point $a_{M'} \in \ctA_{M'}^{M'\el}$, then
    $$
        (K_{\ol{s}_M}^i)^{\ss}_{\fM',a} = (K_{\ol{s}_M}^i)^{\ss}_a.
    $$
    \begin{proof}
        If $\fG' = \fG'_s$ is induced, then $M'$ is a semistandard Levi in $G'_s$ by Proposition \ref{prop.ind.G'.cont.M'.Levi}. The Hitchin base is functorial in endoscopic groups and in particular in Levi subgroups, so we have a diagram of closed embeddings 
$$
\begin{tikzcd}[ampersand replacement=\&]
    \& \ctA_{M} \ar[dr, "\chi^M_G"] \& \\
     \ctA_{M'} \ar[ur, "\iota_{\rho_{s_M}}"] \ar[dr, "\chi^{M'}_{G'}"] \& \& \ctA_G \\
     \& \ctA_{G'_s} \ar[ur, "\iota_{\rho_s}"] \& 
\end{tikzcd}
$$
We can easily see that this square commutes, so $\eta_{\fM'}$ belongs to $\iota_{\rho_{s}}(\ctA_{G'_s})$, and is in $U_M^{\good}$. Corollary \ref{cor.full.supp} implies $P$ has full support on the intersection, so $\eta_{\fM'}$ is in the support of $P$. 

Conversely, if \textit{any} point $a$ of $\Supp(P)$ intersects $\iota_{\rho_{s_M}}(\ctA_{M'}^{M'\el})$, then we will show $\fG'$ is induced. Such an $a$ lies in particular in $\iota_{\rho_s}(\ctA_{G_s'}) \cap \iota_{\rho_{s_M}}(\ctA_{M'}^{M'\el})$. Since $s$ lies above $\ol{s}_M$, we can write
for some $s = s_M z$, with $z \in Z_{\widehat{M}}^0$. Then we must have $z$ Galois-fixed and so $z \in (Z_{\widehat{M}}^{0})^{\Gamma}$, as follows. Indeed, at any point $a$ in the intersection, we have $W_{\wtilde{a}} \subseteq W_s$, so $W_{\wtilde{a}}$ stabilizes $s = s_M z$ and $\mu_{\wtilde{a}}$ induces $\rho_s$, but also stabilizes $s_M$. For any $\gamma \in \pi_1(U_a, \infty)$, let $(c_{\wtilde{a}}(\gamma), \ol{\gamma}) \in \W_G \rtimes \Gamma$ living above $\mu_{\wtilde{a}}(\gamma) \in W_{\wtilde{a}} \subseteq \W \rtimes \Out(\G)$ in the monodromic subgroup. Since $\ctA_{M'} \subset \ctA_{M}$ and so $c_{\wtilde{a}}(\gamma) \in W_M$ acts trivially on the center of $M$, we obtain
$$
    s_M z = (c_{\wtilde{a}}(\gamma), \ol{\gamma})(s_M z) = s_M (c_{\wtilde{a}}(\gamma), \ol{\gamma})(z) = s_M \ol{\gamma}(z)
$$
and so $z = \ol{\gamma}(z)$ for all $\ol{\gamma} \in \Gamma$, as claimed. 

By that calculation, $s \in s_M Z_{\widehat{M}}^\Gamma$, so defines an induced datum $(s, \rho_s')$ by the construction above Definition \ref{def.induced.end.datum}. The commuting square above for this new datum shows $a$ also lies in $\iota_{\rho_s'}(\ctA_{G'_s})$. But by Proposition \ref{prop.endoscopic.locus.decomposition} the images of the embeddings $\iota_{\rho_s} : \ctA_{G'_s} \to \ctA_G$ are disjoint as $\rho_s$ ranges over the finite collection of torsors for a fixed $s$, so we must have $\rho_s = \rho_s'$ and $(s, \rho_s)$ is induced. Since $\fG'$ is elliptic, $s \in \cE_{\fM'}^f(G)$. Taking $a = \eta_{\fM'}$ gives the converse, but since $a$ was arbitrary we get also the stalk-wise statement.
    \end{proof}
\end{lemma}

We deduce our main cohomological theorem, under our running hypothesis of $G$ quasisplit reductive with anisotropic center.
\begin{theorem}\label{thm.coh.wfl}
    Let $\fM' = (s_M, \rho_{s_M})$ an arithmetic (defined over $k_0$) elliptic endoscopic datum for $M$ with endoscopic group $M'$ which is split at $\infty$. We have an isomorphism of semisimplified perverse Weil sheaves on $U_M^{\good}$
    $$
(\pHi(Rf_\star^{\xi, \Mpar} \ol{\Q}_l)_{\ol{s}_M})^{\ss}_{\fM'} \cong \left(\bigoplus_{s \in \cE_{\fM'}^f(G)} (j_{ell})_{!\star}(i_{\rho_s})_\star \pHishift(Rf_{G_{s}',\star}^{G'_s\el} \ol{\mathbb{Q}}_l)_1(-r_s)\right)^{\ss},
    $$
    where
    \begin{itemize}
        \item $\rho_s$ and $G_s'$ are the induced endoscopic data of Definition \ref{def.induced.end.datum} for $s \in \cE_{\fM'}^f(G)$,
        \item $\ss$ denotes the semisimplification of perverse sheaves,
        \item $(\pHi(Rf_\star^{\xi, \Mpar} \ol{\Q}_l)_{\ol{s}_M})^{\ss}_{\fM'}$ is the coarse semisimple $\fM'$-component of Definition \ref{def.fM'.cpt},
        \item $\pHi(R f_{G_{\rho_s}',\star}^{G_{\rho_s}'\el} \ol{\Q}_l)_1$ is the $1$-isotypic part in the fine decomposition of Theorem \ref{thm.s.decomp} applied to the (quasisplit reductive) $X$-endoscopic $G_{\rho_s}'$,
    \item $r_s = \codim_{\ctA_G}(\ctA_{G_{\rho_s}'})$.
    \end{itemize}
\begin{proof}
    Again denote $K_{\ol{s}_M}^i = \pHi(Rf_\star^{\xi, \Mpar} \ol{\Q}_l)_{\ol{s}_M}|_{U^{\good}_M}$, viewed as a Weil sheaf. After base change to $k = \ol{k}_0$, we have a geometric decomposition $$
        j_{ell}^\star (K^i_{\ol{s}_M})^{\ss}|_k = \bigoplus_{s} \bigoplus_{(s, \rho_s)} K^i_{G, s, \rho_s},
    $$
    which coarsens to an arithmetic decomposition when we group the endoscopic data into Frobenius orbits. By Lemma \ref{lem.M.G.comp.of.supports} above, the support condition of Definition \ref{def.fM'.cpt} includes every geometric simple constituent of the induced summands and excludes every other constituent. Each summand is Frobenius-stable by Remark \ref{rmk.ind.data.arithmetic}. We obtain an arithmetic decomposition (up to semisimplification) over induced data
    $$
        j_{ell}^\star (K_{\ol{s}_M}^i)^{\ss}_{\fM'} \cong \bigoplus_{s \in \cE_{\fM'}^f(G)} K_{s, \rho_s}^{\ss}.
    $$
    The last statement of Theorem \ref{thm.ngo.coh.fl} implies that each $K_{s, \rho_s}^{\ss}$, the semisimplified fine summand with support $\iota_{\rho_s}(\ctA_{G'_s})$, is isomorphic as a Weil sheaf to $$\left((i_{\rho_s})_\star \pHishift(Rf_{G_s',\star}^{G_s'\el} \ol{\Q}_l)_1(-r_s)\right)^{\ss}.$$ Theorem \ref{thm.no.new.supports} implies that
    $$
        K^i_{\ol{s}_M} \cong (j_{ell})_{!\star} j_{ell}^\star K^i_{\ol{s}_M}
    $$
    over $k_0$, i.e., as Weil sheaves, and likewise for its semisimplified summand $(K^i_{\ol{s}_M})^{\ss}_{\fM'}$, so intermediate extension gives the theorem.
\end{proof}
\end{theorem}

\part{Arithmetic applications}

\chapter{Point counts, $s$-orbital integrals, and the global WFL}\label{sect.global.point.counts}

In this chapter, we define the global $s$-orbital integrals for our setting and carry over roughly unchanged the Frobenius point-counting argument of Chaudouard--Laumon \cite[\S 11, 12]{CLII} for $\xi$-weighted Hitchin fibers. In the first section, we define $s$-orbital integrals, prove their invariance up to $Z_{\widehat{M}}^\Gamma$. In the second section, we express the trace of Frobenius on a \emph{coarse} component of the relative $\xi$-stable Hitchin cohomology as a global (adelic) weighted orbital integral. For $\xi$ in general position, the result is in terms of a global version of Arthur's weighted orbital integrals, and in particular is independent of $\xi$. This provides the arithmetic interpretation of the left side of the cohomological WFL \ref{thm.coh.wfl}. In the third subsection, we define recursively the stabilized weighted orbital integrals, and inductively prove that the trace of Frobenius on a \emph{fine} stable summand equals the corresponding stable orbital integral. Together, these give in Theorem \ref{thm.global.wfl} a global incarnation of the WFL that matches Arthur's expected form. 

The moral of this chapter is that the techniques applied in the split case by \cite{CLI} will apply equally well here, as long as we are working over the compatible locus and all weight conditions are defined relative to the \textit{rational} notion of $\xi$-stability from Definition \ref{def.rel.stab.convex}. This is absolutely necessary in order to interface with Arthur's weight functions, which are defined purely in terms of relative/rational data.

We recall our running situation and fix notation: let $k_0 = \F_q$ and $k = \ol{\F}_q$. Let $X_0$ be a curve over $k_0 = \F_q$ with generic point $\eta_0$, function field $F_0$, and ad\`eles $\A_0$. Let $G_0$ a quasisplit\footnote{This means pinned quasisplit, as in \S \ref{sect.rel.reduct.gps}.} reductive $X_0$-group. We will often assume $G_0$ to have anisotropic center. Fix Borel $B_0$ and maximal torus $T_0$, and assume the characteristic of $k_0$ is good for $G_0$. We fix a geometrically connected Galois $\Gamma$-cover $X'_0 \to X_0$ which splits $G_0$, cf.\ Proposition \ref{prop.splitrankextensions}. Denote without $0$ the corresponding objects over $X = X_0 \times_{k_0} k$. We denote by $\tau$ the arithmetic Frobenius automorphism $x \mapsto x^q$ of $k$ as well as the induced automorphism $1 \times \tau$ of $X$, $G$, $B$, etc. Proposition \ref{prop.splitrankextensions} and its Corollary imply that the torus $T$ is no more generically split than $T_0$, that the sets of semistandard Levis and parabolics in $G_0$ and in $G$ are in bijection, and that the rational character and cocharacter lattices remain the same when base changing from $k_0$ to $k$, i.e., the rational lattices satisfy
$$
    X^\star(M) = X^\star(M_0) \qquad \text{ and }\qquad X_\star(M) = X_\star(M_0)
$$
for each semistandard Levi $M$. The same invariance applies to the dual space $\cX_\star(M)$ and the vector spaces $\a_M$ and $\a_M^\star$.

For each closed point $v_0 \in |X_0|$, write $k_{v_0}$ for its residue field, $z_{v_0}$ for a uniformizer at $v_0$, $\O_{v_0}$ for its (complete) local ring of integers and $F_{v_0}$ for its local fraction field. Let $L^+_{G,v_0}$, $L_{G,v_0}$, and $\Gr_{G,v_0}$ denote the arc group, loop group, and affine Grassmannian of $G_0|_{\O_v}$, respectively. Likewise, write $L^+_{\fg, v}$, $L_{\fg, v}$ for the arc and loop Lie algebras, with $k_{v_0}$-points given by $\fg(\O_{v_0})$ and $\fg(F_{v_0})$. Pick an effective divisor 
\[
    D_0 = \sum_{v_0 \in |X_0|} d_{v_0}  v_0
\]
which is even (i.e.\ $d_{v_0}$ is even for all $v_0$), of degree $\deg(D_0) > 2g$, and disjoint from $\infty$. Denote its base change to $X$ by $D$. Write $1_D$ the indicator function of $(z_v^{-d_v})_v \fg(\O)$. Denote 
\[
    K_v := L^+_{G,v}(k_v)
\]
and $K := \prod_v K_v$ the compact open in $G(\A_0)$ corresponding to the integral group scheme $G_0$.

\section{Global weighted $s$-orbital integrals.}\label{subsect.weighted.s.orbital.global}

In this section, $G$ is not required to have anisotropic center. We recall Chaudouard--Laumon's construction of global weighted $s$-orbital integrals. The results of this section follow exactly as in \cite[\S11]{CLII}, with adapted definitions (in particular, the monodromic subgroup $W_{a_M} \subseteq W_M \rtimes \Out(\M)$ of Definition \ref{def.monodr.inv.Out}) and some careful tracking of $\Gamma$-invariants and coinvariants. We begin with an application of Kottwitz's Lemma \ref{lem.kottwitz.lemma} above.

\begin{proposition}[{\cite[11.2.1]{CLII}}]\label{prop.kottwitz.Ja}
    Let $U_0 \subset X_0$ an open containing $\infty$, and consider the functor $B : \mathsf{Tori}_{U_0} \to \Ab$ given by
    $$
        B(Z) \mapsto (Z(F) \backslash Z(\A))_\tau.
    $$
    This satisfies the hypotheses (1), (2) of Kottwitz's lemma, and so there is a canonical isomorphism of functors $\ul{X}_\star(-)(\infty)_{\pi_1(U_0, \infty)} \cong B(-)$. Moreover, the association $[\mu] \mapsto [z_\infty^{\mu}]$ defines a natural map $\ul{X}_\star(-)_\infty \to B$ such that the natural triangle
    $$
        \begin{tikzcd}[ampersand replacement=\&]
            \ul{X}_\star(-)_\infty \ar[d] \ar[dr] \& \\
            \ul{X}_\star(-)(\infty)_{\pi_1(U_0, \infty)} \ar[r, "{\sim}" below] \& B(-)
        \end{tikzcd}
    $$
    commutes.
\end{proposition}

Denote $\kappa_Z : (Z(F) \backslash Z(\A))_\tau \xrightarrow{\sim}  \ul{X}_\star(Z)(\infty)_{\pi_1(U_0, \infty)}$ the inverse isomorphism, the \emph{Kottwitz map} associated to the torus $Z/U_0$.

Choose now a semistandard Levi $M$ and a point $a_M = (h_{a_M}, t) \in \ctA_M(k_0) = \ctA_M(k)^{\tau}$. The $M$-regular centralizer $$J_{a_M} := h^\star_{a_M} J^M_D$$ is a torus over $U_{0,a_M} := h^{-1}_{a_M}(\fc_{M,D}^{\Greg})$, so by the previous Proposition we have an isomorphism
$$
    (J_{a_M}(F)\backslash J_{a_M}(\A))_{\tau} \cong \ul{X}_\star(J_{a_M})(\infty)_{\pi_1(U_{0, a_M}, \infty)}.
$$
Note that the action of the arithmetic and geometric monodromy on the cocharacter lattice coincide; the Frobenius image induced by the choice of the split point $\infty$ acts trivially. As in Remark \ref{rmk.monodr.inv.reg.centr}, the action of the (arithmetic) $\pi_1(U_{0, a_M}, \infty)$ on the $k_0$-fiber $\bX_\star := \ul{X}_\star(T)|_\infty \cong \ul{X}_\star(J_a)|_\infty$ is the arithmetic monodromy homomorphism of Definition \ref{def.monodr.inv.Out}, which has image equal to the image of the geometric monodromy because we have chosen rational rigidification at $\infty$ where $T$ is split, so the image of the arithmetic section induced by infinity
$$
    \Gal(\ol{k}_0/k_0) \to \pi_1(U_{0, a_M}, \infty)
$$
acts trivially on $\ul{X}_\star(J_a)|_\infty$. So we deduce an isomorphism 
$$\kappa_J : (J_{a_M}(F) \backslash J_{a_M}(\A))_\tau \cong (\bX_\star)_{W_{a_M}}.$$

\begin{remark}\label{rmk.Kottwitz.signs}
    As in all instances of Kottwitz's lemma, there is a choice of sign normalization coming from a choice of generator for $B(\G_m) \cong \Z$. Here we have chosen the positive \emph{valuation} normalization
    $$
        [\mu] \mapsto [z_\infty^{\mu}]
    $$
    to match the conventions in \cite[Proposition 11.2.1]{CLII}. Note that this differs by a sign from the positive \emph{degree} Beauville--Laszlo convention used above (cf.\ \ref{not.BL.conventions}, \ref{lem.kottwitz.lemma}), in which $\mu \in \ul{X}_\star(T)(\infty)$ maps to the modification of the trivial torsor by $z_\infty^{-\mu}$ at $\infty$. This allows our statements and definitions below to extend those in \textit{loc.\ cit.}, but the slightly awkward consequence that 
    $$
        \begin{tikzcd}
             (J_{a_M}(F) \backslash J_{a_M}(\A))_\tau \ar[r, "\kappa_J"] \ar[dr, "H_M" below left] & (\bX_\star)_{W_{a_M}} \ar[d]\\ & \cX_\star(M)
        \end{tikzcd}
    $$
    does not commute but anticommutes, where $H_M$ is the negative sum of valuations from Definition \ref{def.adelic.hitchin.deg} and the vertical arrow is the natural coinvariants-mod-torsion map, as in Lemma \ref{lem.compos.deg.pic}.
\end{remark}

We can now define the invariant map which appears in the definition of $s$-orbital integrals. Consider the sequence of maps of homogeneous spaces with Frobenius-actions
$$
    \begin{tikzcd}[ampersand replacement=\&]
        1 \ar[r] \& J_{a_M}(\A) \ar[r] \& G(\A) \ar[r] \& J_{a_M}(\A) \backslash G(\A) \ar[r] \& 1
    \end{tikzcd}.
$$
This is short-exact in the obvious sense. From its long-exact sequence in nonabelian Galois cohomology, one obtains a coboundary map
$$
    \epsilon : (J_{a_M}(\A) \backslash G(\A))^{\tau} \to \ker(J_{a_M}(\A)_\tau \to G(\A)_\tau) \subset J_{a_M}(\A)_\tau
$$
with fibers orbits under $G(\A_0)$. Explicitly, the map identifies a class represented by $g \in G(\A)$ with the $\tau$-conjugacy class $g \tau(g)^{-1} \in J_{a_M}(\A)_\tau$, and two classes have the same image if and only if they differ by right translation by $G(\A_0)$.

\begin{definition}
    Given $a_M \in \ctA_{M}(k_0) = \ctA_M(k)^\tau$, define
    $$
        \inv_{a_M} : (J_{a_M}(\A) \backslash G(\A))^{\tau} \to (\bX_\star)_{W_{a_M}}
    $$
    to be the composition
    $$
        \begin{tikzcd}[ampersand replacement=\&] (J_{a_M}(\A) \backslash G(\A))^{\tau} \ar[r, "{\epsilon}"] \& J_{a_M}(\A)_\tau \ar[r] \& (J_{a_M}(F)\backslash J_{a_M}(\A))_{\tau} \ar[r, "{\kappa_J}"] \& (\bX_\star)_{W_{a_M}}
        \end{tikzcd}.
    $$
    Write $\langle -, -\rangle$ for the pairing $\widehat{T}^{W_{a_M}} \times (\bX_\star)_{W_{a_M}} \to \C^\times$ induced by the contragredient action of the monodromic subgroup $W_{a_M}$ on $\widehat{T} = \Hom(\bX_\star, \C^\times)$.
\end{definition}

\begin{proposition}[{\cite[Lemme 11.4.2]{CLII}}]\label{prop.inv.sc.cover}
    For any $a_M \in \ctA_{M}(k_0)$ and $g \in (J_{a_M}(\A) \backslash G(\A))^{\tau} \to (\bX_\star)_{W_{a_M}}$, the invariant map takes values
    $$
    \inv_{a_M}(g) \in \im\left(\ul{X}_\star(T_{M_{\sc}})(\infty)_{W_{a_M}} \to  \ul{X}_\star(T)(\infty)_{W_{a_M}}\right)
    $$
    \begin{proof}
        Exactly as in \textit{loc.\ cit.}
    \end{proof}
\end{proposition}

\begin{definition}\label{def.Haar.measure.max.cmpt.convention} Equip $G(\A_0)$ with the Haar measure which assigns volume 1 to the compact open $K$ corresponding to global integral model $G_0$, and equip $J_{a_M}(\A_0)$ with the Haar measure which assigns volume 1 to its unique maximal compact open (recall it is commutative!). As in \cite[\S 11.5]{CLII}, we observe that $J_{a_M}(\A)_\tau$ is countable, so along the coboundary map $(J_{a_M}(\A) \backslash G(\A))^{\tau}$ is the disjoint union of countably many orbits under $G(\A_0)$, each of which can be equipped with the quotient measure from $G(\A_0)$. So we obtain a measure on $(J_{a_M}(\A) \backslash G(\A))^{\tau}$ which is right $G(\A_0)$-invariant. \end{definition}

Now we can define weighted $s$-orbital integrals. As in \S \ref{subsect.adelic.trunc.Hitchin}, Definition \ref{def.adelic.hitchin.deg}, there is the adelic characterization of the weight function of Arthur: there is for each rational semistandard parabolic $P$ with Levi $M$ a well-defined map 
$$
    H_P : G(\A) \to \a_M = \ul{\a}_M^\Gamma
$$
defined by the requirement that for each rational character $\lambda \in X^\star(P) = \ul{X}^\star(P)^\Gamma$,
$$
    \langle \lambda, H_P(g)\rangle = -\sum_{\ol{v} \in X(k)} \val_{\ol{v}}(\lambda(p_{\ol{v}}))
$$
where $g \in pG(\O)$ in a (geometric) Iwasawa decomposition. Recall that we have fixed a choice of Weyl-invariant inner product on the rational character space $\a_T$ and so on all its induced summands. 
\begin{definition}\label{def.arthur.wt.ftn}
    For $g \in G(\A)$, define \emph{Arthur's weight} $v_M^G(g)$ to be the volume in $\a_M^G$ (the \emph{rational} cocharacter space) of the convex envelope spanned by the projections of the $-H_P(g)$ for $P$ with Levi $M$, with respect to the measure on $\a_M^G$ induced from the global choice in Definition \ref{def.cochar.space.measure.choice}.
\end{definition}
In particular, the measure obtained in this way is invariant for the relative Weyl group, which is all that Arthur ever requires in his weight formalism.

\begin{remark}
    For a closed point $v \in |X_0|$, we can likewise define a closed-place contribution to the global height using the same valuation normalization $\val_v(\pi_v) = 1$. For $g_v \in G(F_{0,v})$, this is then related by
    $$
        H_{P,v}(g_v) = \sum_{\ol{v}|v} H_{P,\ol{v}}(g_v) = \deg(v) H_{P, \ol{v}}(g_v).
    $$
    The global $H_P$ then satisfies on $G(\A_0)$ the closed-point formula 
    $$
        H_P(g) = \sum_{v \in |X_0|} H_{P,v}(g_v).
    $$
\end{remark}

Recall that $H_P|_{M(\A)}$ is independent of $P$ and satisfies $H_P(mg) = H_M(m) + H_P(g)$. It follows that left multiplication by $m \in M(\A)$ translates each vertex $-H_P(g)$ by $-H_M(m)$, and so $v_M^G(mg) = v_M^G(g)$. The weight function is also clearly invariant under right multiplication by $G(\O)$. 

\begin{definition}[{\cite[D\'efinition 11.7.1]{CLII}}]\label{def.s.orbital.integral.adelic}
    Let $M$ be a semistandard Levi and $a_M \in \ctA_M(k_0) = \ctA_M(k)^\tau$. Let $s \in \widehat{T}^{W_{a_M}}$. The \emph{weighted $s$-orbital integral} associated to $a_M$ is defined by
    $$
        J_M^{G,s}(a_M) = q^{-\deg(D)|\Phi^G|/2} \int_{(J_{a_M}(\A) \backslash G(\A))^{\tau}} \langle s, \inv_{a_M}(g)\rangle 1_D(\Ad(g)^{-1} X_{a_M}) v_M^G(g) \; dg,
    $$
    where \begin{itemize} 
    \item $\Phi^G$ is the set of absolute roots of $G$, 
    \item $X_{a_M} \in \fm(F)$ is the image of $a_{M, \eta}$ under the Kostant section, and
    \item the integral is taken with respect to the measure on $(J_{a_M}(\A) \backslash G(\A))^{\tau}$ defined above.
    \end{itemize}
\end{definition}

It will be helpful (and psychologically useful) to rewrite this as a sum of integrals over the values of the invariant. Given $a_M \in \bAe{M}$, recall the (finite) coarse endoscopic group of Definition \ref{def.end.symmetry.and.residual.torsion}:
$$
    \Sigma_a := \im(\ul{X}_{\star}(T \cap D(M))(\infty)_{W_{a_M}} \to \ul{X}_{\star}(T)(\infty)_{W_{a_M}}). 
$$
As a corollary of Proposition \ref{prop.inv.sc.cover}, $\inv_{a_M}(g) \in \Sigma_a$ for every $g \in (J_{a_M}(\A) \backslash G(\A))^{\tau}$. So writing for $\lambda \in \Sigma_a$
\begin{equation}\label{eqn.sigma.cpt.adelic.pts}
    (J_{a_M}(\A) \backslash G(\A))^{\tau, \lambda} := \left\{g \in (J_{a_M}(\A) \backslash G(\A))^{\tau} : \inv_{a_M}(g) = \lambda\right\},
\end{equation}
and equipping it with the induced measure as a closed and open subset of $(J_{a_M}(\A) \backslash G(\A))^{\tau}$, we deduce the following decomposition of the $s$-orbital integral.

\begin{proposition}\label{prop.s.orbital.int.decomp}
    Under the same notations as Definition \ref{def.s.orbital.integral.adelic}, we have
    $$
        J_M^{G,s}(a_M) = q^{-\deg(D)|\Phi^G|/2} \sum_{\lambda \in \Sigma_a} \langle s, \lambda \rangle \int_{ (J_{a_M}(\A) \backslash G(\A))^{\tau,\lambda}} 1_D(\Ad(g)^{-1} X_{a_M}) v_M^G(g) \; dg.
    $$
\end{proposition}

Finally, in order to obtain statements of the Weighted Fundamental Lemma which match Arthur's (in particular, ranging over $s \in s_M Z_{\widehat{M}}^\Gamma$ rather than $s \in s_M (Z_{\widehat{M}}^0)^\Gamma$, cf.\ Remark \ref{rmk.discrep.cWFL.vs.Arthur}), it will be important to know that the $s$-orbital integrals satisfy a corresponding invariance. 

\begin{proposition}[{cf.\ \cite[Proposition 11.7.3]{CLII}}]\label{prop.s.orbital.int.invariant.under.center}
    The map $s \in \widehat{T}^{W_{a_M}} \mapsto J_M^{G,s}(a_M)$ is invariant under translation by $Z_{\widehat{M}}^\Gamma$.
    \begin{proof}
        Let $M_{\sc}$ denote the simply connected $X$-cover of the derived group of $M$ and $T_{M_{\sc}}$ the inverse image of $T$ inside $M_{\sc}$. By Proposition \ref{prop.inv.sc.cover}, one has $$\inv_{a_M}(g) \in \im\left(\ul{X}_\star(T_{M_{\sc}})(\infty)_{W_{a_M}} \to  \ul{X}_\star(T)(\infty)_{W_{a_M}}\right).$$ Unlike in the split case, $W_{a_M}$ need not fix $Z_{\widehat{M}}$ pointwise because of its outer-automorphism component. But dualizing the map of coinvariants shows that the pairing $\langle s, - \rangle$ is invariant under the annihilator 
        $$\ker\left(\widehat{T}^{W_{a_M}} \to \widehat{T}_{\sc}^{W_{a_M}}\right),$$
        which is canonically $Z_{\widehat{M}}^{W_{a_M}}$ by left-exactness of fixed points. The $W_M$ part, of course, acts trivially on the center, so $$Z_{\widehat{M}}^{W_{a_M}} = Z_{\widehat{M}}^{\pi_1(X,\infty)} = Z_{\widehat{M}}^\Gamma,$$ as desired.
    \end{proof}
\end{proposition}

\section{The trace of Frobenius on the coarse $\ol{s}$-component.}

Fix semistandard Levi $M$, and assume that $G$ has anisotropic center. In this section, we compute the alternating sum of traces of Frobenius on the stalk of the coarse $s$-decomposition of Theorem \ref{thm.s.decomp} over $a \in \bAe{M}(k_0)$. By Proposition \ref{prop.s.isotypic.cpt.support}, we know that $\ol{s}_M \in \widehat{T}/Z_{\widehat{M}}^0$ can only contribute components to the stalk over $a$ if $s_M$ admits a lift $s \in \widehat{T}^{W_a}$. We will perform the computation first in terms of such a lift. In this section we adapt Chaudouard--Laumon's point-counting argument to give the following theorem.

\begin{theorem}[{cf.\ \cite[Th\'eor\`eme 12.1.1]{CLII}}]\label{thm.s.summand.orbital.integral}
    Let $a_M \in \bAe{M}(k_0)$ with image $a \in \bA_G$, and $\xi \in \a_T$ in general position. Let $s \in \widehat{T}^{W_a}$ have image $\ol{s} \in \widehat{T}^{W_a}/(Z^0_{\widehat{M}})^\Gamma \subset (\widehat{T}/Z^0_{\widehat{M}})_{\tors}$. Then the trace of Frobenius of the corresponding coarse summand satisfies
    $$
        \sum_{i} (-1)^i \Tr(\tau^{-1}, \pHi(R f_\star^{\xi, \Mpar} \ol{\Q}_l)_{a, \ol{s}}) = \frac{c(J_{a_M})}{\vol(\a_M/\cX_\star(M)) |\widehat{T}^{W_a}/(Z^0_{\widehat{M}})^{\Gamma}|} q^{\deg(D) |\Phi^G|/2} J_M^{G,s}(a_M),
    $$
    where
    \begin{itemize}
        \item $\a_M$ is the relative cocharacter space of Definition \ref{def.char.space},
        \item $\cX_\star(M)$ is the dual of the relative character lattice $X^\star(M)$, as in Proposition \ref{prop.dual.lattice.global.invs.coinvs.mod.tors} and Definition \ref{def.cochar.sheaf},
        \item $J_{a_M}$ is the fiber of the $M$-regular centralizer of Ng\^{o},
        \item $\Phi^G$ is the set of absolute roots for $G$,
        \item $c(J_{a_M})$ is a positive constant defined below (\ref{eqn.def.of.cJ}).
    \end{itemize}
\end{theorem}

By the Theorem and Proposition \ref{prop.s.orbital.int.invariant.under.center}, it follows that the left side depends only on the image of $s$ in $\widehat{T}/Z_{\widehat{M}}^{\Gamma}$.

Before we begin the proof, a few remarks comparing this section with \cite[\S 12]{CLII}. We follow the same strategy, with the natural corrections in order to keep everything relative. Note the appearance of $\cX_\star(M)$, the natural target lattice for rational degree computations, in the statement of the theorem, and that $\ol{s}$ naturally lives in $\widehat{T}^{W_a}/(Z^0_{\widehat{M}})^\Gamma$ in our generality, where $W_a$ is the monodromic invariant of Definition \ref{def.monodr.inv.Out}. Several negative signs on lattices and $\Sigma_a$ are corrected from \textit{loc.\ cit.}; comments about this are mostly relegated to footnotes and can be safely ignored, since the final statement in Proposition \ref{prop.s.orbital.int.decomp} is symmetric under re-indexing $\lambda \mapsto -\lambda$.

We order the argument slightly differently and try to offer a complementary exposition. In particular, we streamline the computation by performing Fourier inversion directly on the endoscopic symmetry group of Definition \ref{def.end.symmetry.and.residual.torsion}:
$$
    \Sigma_a := \im(\ul{X}_{\star}(T \cap D(M))(\infty)_{W_a} \to \ul{X}_{\star}(T)(\infty)_{W_a}),
$$
rather than on $\ul{X}_{\star}(T \cap D(M))(\infty)_{W_a}$. By Proposition \ref{prop.dual.Y1.fibers}, we know that 
$$
    \widehat{\Sigma}_a \cong \widehat{T}^{W_a}/(Z^0_{\widehat{M}})^{\Gamma},
$$
finite when $a_M \in \bAe{M}$. Note the second constant factor in the denominator is therefore exactly 
$$
    |\Sigma_a| = |\widehat{\Sigma}_a|.
$$

\begin{proof}
    Recall that as in \S \ref{subsect.comp.grp.action}, we have actions on the cohomological stalks
    $$
        K_a^i := \pHi(R f_\star^{\xi, \Mpar} \ol{\Q}_l)_{a}
    $$
    through the composition
    $$
        \bX_{\star,D(M)} \twoheadrightarrow \Sigma_a \to \pi_0(\bP^{1}_a).
    $$
    We will absuively identify $\lambda \in \bX_{\star, D(M)}$ with its image in $\Sigma_a$.
    
    Splitting across the summands of the $s$-decomposition, we know that pullback by $\lambda \in \bX_{\star, D(M)}$ acts on the $\ol{s}$-component by $\langle s, \lambda \rangle$, where $\ol{s} \in (\widehat{T}/Z_{\widehat{M}}^0)_{\tors}$.\footnote{Note that our choice of signs is consistent with Chapter \ref{sect.reg.cent.Pic}: $\lambda$ acts on this Hitchin fiber via the $J$-torsor with Beauville--Laszlo data $[z_\infty^{-\lambda}]$. We are then going to consider the pullback action of cohomology from this action, also denoted with $\lambda$.} So for any $\lambda \in \bX_{\star, D(M)}$ we have an equality of alternating traces
    $$
        \Tr(\lambda \tau^{-1}, K_a^{\bullet}) = \sum_{\ol{s} \in (\widehat{T}/Z_{\widehat{M}}^0)_{\tors}} \langle \ol{s}, \lambda\rangle\Tr(\tau^{-1}, K_{a,\ol{s}}^{\bullet}).
    $$
    But by Proposition \ref{prop.s.isotypic.cpt.support} we know only $\ol{s} \in \widehat{\Sigma}_a \subset (\widehat{T}/Z_{\widehat{M}}^0)_{\tors}$ contribute, and on the other hand $\lambda$ acts on cohomology through its image in $\Sigma_a$, so we can coarsen our decomposition to
    $$
    \Tr(\lambda \tau^{-1}, K_a^{\bullet}) = \sum_{\ol{s} \in \widehat{\Sigma}_a} \langle \ol{s}, \lambda\rangle \Tr(\tau^{-1}, K_{a,\ol{s}}^{\bullet}),
    $$
    where the pairing is now that between $\Sigma_a$ and $\widehat{\Sigma}_a$. By Fourier inversion, we recover
    $$
        \Tr(\tau^{-1}, K_{a,\ol{s}}^{\bullet}) = \frac{1}{|\Sigma_a|} \sum_{\lambda \in \Sigma_a} \langle \ol{s}, \lambda \rangle^{-1} \Tr(\lambda \tau^{-1}, K_a^{\bullet}).
    $$

    By abstract nonsense, the alternating trace of an endomorphism on stalks of the perverse cohomology is equal to the alternating trace on stalks of the classical cohomology: in the Grothendieck group of the derived category of pairs $(K \in D^b_c(U_M), \phi : K \to K)$, 
    $$
        [R f_\star^{\xi,\Mpar} \ol{\Q}_l, \lambda \tau^{-1}] := \sum_i (-1)^i [\pHi(R f_\star^{\xi} \ol{\Q}_l), \lambda \tau^{-1}],
    $$
    which pulls back along $i_a^\star$ to the corresponding equality in the stalk. But on the other hand, proper base change gives $i_a^\star R f_\star^{\xi, \Mpar} \ol{\Q}_l \cong R\Gamma(\bM^\xi_a, \ol{\Q}_l)$. Taking alternating trace therefore gives 
    $$
    \Tr(\lambda \tau^{-1}, K_a^{\bullet}) = \Tr(\lambda \tau^{-1}, H^\bullet(\bM^{\xi}_a, \ol{\Q}_l)).
    $$
    Grothendieck--Lefschetz then identifies 
    $$
        \Tr(\lambda \tau^{-1}, H^\bullet(\bM^{\xi}_a, \ol{\Q}_l)) = \# \bM_a^\xi(k)^{\lambda \tau},
    $$
    where the fixed-points and the cardinality are taken in the groupoid sense, cf.\ \cite[\S 8.1]{ngofl}.

    Geometrically, we have already given a description of the fiber groupoid $\bM_a^{\xi}(k)$ for $a \in \bAe{M}(k)$ in Proposition \ref{prop.adelic.hitchin.fiber.reduced.kostant}: it is the groupoid of $g \in G(\A)/G(\O)$ such that \begin{equation}\label{eqn.wtd.stab.conds.g}
    1_D(\Ad(g)^{-1}(X_{a_M})) = 1 \qquad \text{ and } \qquad 
        \xi_M \in \Cvx((-H_P(g))_{P})
    \end{equation}
    under the left multiplication action of $J_{a_M}(F)$.
    
    The data $t$ determines a canonical isomorphism $J_{a_M}|_{\O_\infty} \cong T|_{\O_\infty}$. Write $j_{\lambda} \in J_{a_M}(F_\infty)$ as the element corresponding to $z_\infty^{-\lambda}$ for $\lambda \in \bX_{\star, D(M)}$, and identify it with an adelic element with all other components 1.We somewhat abusively denote $j_\lambda$ for $\lambda \in \Sigma_a$, which implicitly means to choose a lift to $\wtilde{\lambda} \in \bX_{\star, D(M)}$ and then perform the construction above. The action of $\lambda \in \bX_{\star, D(M)}$ on the Hitchin fiber is through the action of the Picard stack and its action on cohomology factors through $\Sigma_a$, so the resulting formulae will all be independent of the lift $\wtilde{\lambda}$. As described in \S\ref{subsect.pic.action.hitchin.fibers}, this action is by taking the associated product with $J_a$-torsors, and as described in Proposition \ref{prop.inv.map.BL}, the canonical map $\bX_\star \to \ctP$ is given by constructing a $J_a$-torsor with gluing data $j_{\lambda}$ at $\infty$ and trivial data at all other points. Arithmetically, then, the action of $\lambda$ on the Beauville--Laszlo family $\tau(g)$ is simply left multiplication by $j_{\lambda}$, viewed as an adelic element.\footnote{This is the element corresponding not to $\lambda$ but to $-\lambda$ under the Kottwitz isomorphism, cf.\ Remark \ref{rmk.Kottwitz.signs}.}

    The fixed-point groupoid $\bM_a^{\xi}(k)^{\lambda \tau}$, then, has objects
    \begin{equation}\label{eqn.fixed.gpoid}
        (g \in G(\A)/G(\O), \delta \in J_{a_M}(F)) : \delta j_\lambda \tau(g) G(\O) = g G(\O),
    \end{equation}
    such that $g$ satisfies (\ref{eqn.wtd.stab.conds.g}), with action as before on $g$ and $\tau$-conjugation on $\delta$. We now perform a reduction.

    Recall $K = G(\O_0)$. Now, as in \cite[Lemme 12.3.3]{CLII}, we reduce to the action groupoid, denoted $[J_{a_M}(F) \backslash \mathcal{W}^\xi_{\lambda}]$, of pairs $(g \in G(\A)/K, \delta \in J_{a_M}(F))$ satisfying (\ref{eqn.wtd.stab.conds.g}) and  \begin{equation}\label{eqn.fixed.gpoid.jlambda}
        \delta j_{\lambda} \tau(g) = g \in G(\A) 
    \end{equation}
    exact equality rather than as cosets. This equality is well-defined on $g K$ since every element of $K$ is $\tau$-fixed. The action on the groupoid is the same action of $J_{a_M}(F)$. There is clearly a fully faithful functor from this groupoid to the one above; to see it is essentially surjective is to show that every class $(g G(\O), \delta)$ such that $$\delta j_\lambda \tau(g G(\O)) = g G(\O)$$ is isomorphic to one satisfying the same condition with $K$. For any given representative $g$, one has
    $$
        \delta j_\lambda \tau(g) = g u
    $$
    for some $u \in G(\O)$. By Lang's theorem for $G(\O)$ (Lemma \ref{lem.Lang.GO} below), we can write $u = x \tau(x)^{-1}$ for some $x \in G(\O)$. Then taking $g' = gx$, we find a new representative that satisfies (\ref{eqn.fixed.gpoid.jlambda}) on the nose.

    Combining the above steps, we have shown
    $$
        \Tr(\tau^{-1}, K^\bullet_{a,\ol{s}}) = \frac{1}{|\Sigma_a|} \sum_{\lambda \in \Sigma_a} \langle \ol{s}, \lambda \rangle^{-1} \# [J_{a_M}(F)\backslash \mathcal{W}^\xi_\lambda].
    $$
    It remains to relate the cardinality of the groupoid $[J_{a_M}(F) \backslash \mathcal{W}^\xi_\lambda]$ to orbital integrals. We do this in Lemma \ref{lem.gpoid.orb.integral} below, which implies the formula
    $$
        \Tr(\tau^{-1}, K^\bullet_{a,\ol{s}}) = \frac{c(J_{a_M})}{\vol(\a_M/\cX_\star(M))|\Sigma_a|} \sum_{\lambda \in \Sigma_a} \langle \ol{s}, \lambda \rangle^{-1} \int_{(J_{a_M}(\A) \backslash G(\A))^{\tau,-\lambda}} 1_D(\Ad(g)^{-1} X_{a_M}) v_M^G(g) \; dg.
    $$
    Re-indexing in terms of $-\lambda$ and combining along the decomposition of Proposition \ref{prop.s.orbital.int.decomp} completes the proof of Theorem \ref{thm.s.summand.orbital.integral}.
\end{proof}

\begin{lemma}\label{lem.Lang.GO}
    The Lang morphism of pro-algebraic group schemes $L = 1 - \tau : G(\O) \to G(\O)$ is surjective.
    \begin{proof}
    It suffices, for each closed point $v_0$ of degree $d$, to show that 
    $$
        1 - \tau : \prod_{v | v_0} L^+_{G,v}(k) \to \prod_{v | v_0} L^+_{G,v}(k)
    $$
    is surjective, where $v$ ranges over the Frobenius orbit of geometric points over $v_0$. The endomorphism $\tau$ permutes the factors cyclically, so a coordinate calculation shows surjectivity reduces to surjectivity of $1 - \tau^d$ on any given factor $L^+_{G,v}(k)$. Replacing $k_0$ with its degree-$d$ extension (and $\tau$ with $\tau^d$, we have reduced to checking Lang's theorem for an arc group at a single Frobenius-fixed place.
    
    Recall that these arc groups are representable by group schemes of infinite type:
        \[
            L^+_{G,v} = \lim_{m \in \Z_{\ge 0}} L^+_{m,G,v} 
        \]
        where $L^+_{m,G,v}: R \mapsto G(R\llbracket z_v\rrbracket/(z_v^{m+1}))$ are twisted jet schemes, representable by irreducible group varieties over $k$ (cf.\ \cite[Prop.\ 2]{HeinlothRelUnif} for the representability). We construct a preimage orbit-by-orbit. For each $m$ we have by Lang's theorem a nonhomomorphic morphism
        \[
            \begin{tikzcd}[ampersand replacement=\&]
                 L^+_{m,G,v}(k_0) \ar[r, hook] \& L^+_{m,G,v}(k) \ar[r, two heads] \& L^+_{m,G,v}(k)
            \end{tikzcd}
        \]
        with finite fibers which are right torsors under $L^+_{m,G,v}(k_0)$. We wish to know that the right map remains surjective upon passage to the limit. This follows if the family $\big(L^+_{m,G,v}(k_0) \big)_m$ satisfies the Mittag--Loeffler condition (cf.\ \cite[Tag 0594]{stacks-project}). The reduction maps $L^+_{m',G,v}(k_0) \to L^+_{m,G,v}(k_0)$ are surjective by the formal smoothness criterion:
        \[
            \begin{tikzcd}[ampersand replacement=\&]
                \Spec k_0 \llbracket z_v\rrbracket/(t_v^{m-1}) \ar[r] \ar[d] \& G \ar[d] \\
                \Spec k_0 \llbracket z_v\rrbracket/(t_v^m) \ar[r] \ar[ur, dashed] \& \Spec k_0\llbracket z_v\rrbracket
            \end{tikzcd}.
        \]
        Thus we can proceed along the proof of \cite[Tag 0598]{stacks-project}, which works without assuming commutativity. Explicitly, let $x \in L^+_{G,v}(k)$ arbitrary, and for a given $m$ let $h_m \in L_{m, G,v}^+(k)$ satisfy $h_m \tau(h_m^{-1}) \equiv x \mod t^m$. We wish to show we can lift $h_m$ to arbitrary level. For any $m' \ge m$ and any $h_{m'}$ mapping to $x$ mod $t^{m'}$, the element $e_m := h_{m}^{-1} h_{m'} \in L^+_{m,G,v}(k)$ is in the kernel of the Lang map $L^+_{m,G,v}(k_0)$ at level $m$. Choosing a Frobenius-fixed lift $e_{m'}$ of $e_m$, the element $h_{m'} e_{m'}^{-1}$ is a lift of $h_m$ which maps to $x$. It follows that surjectivity of the Lang map is preserved upon passing to the limit $L^+_{G,v}(k) \to L^+_{G,v}(k)$.
    \end{proof}
\end{lemma}

Now we perform the fundamental computation of $\# [J_{a_M}(F) \backslash \mathcal{W}^\xi_\lambda]$ in terms of Arthur's orbital integrals. It is easier to first pass through an orbital integral more directly related to the definition of our stability condition. 

\begin{definition}\label{def.weights.CL}
Let $\xi_M \in \a_M$ and define \emph{Chaudouard--Laumon's weight function} on $G(\A)$ by
$$
    w^\xi_M(g) = \# \{\mu \in \cX_\star(M) : 1_{M,g}(\xi_M + \mu) = 1\};
$$
this is the count of points in the shifted lattice $\xi_M + \cX_\star(M)$ which land inside the compact convex $C_g \in \a_M$.\footnote{Recall that this is indeed compact because $G$ has anisotropic center and we are working in the rational space $\a_M$.}
\end{definition}

Since the left action of $M(\A)$ translates the rational cones $$-H_P(g) + \obtcone_P + \a_T^P$$
by $-H_P(m) = -H_M(m)$, which lands inside $\cX_\star(M)$, the weight function is left $M(\A)$-invariant. 

\begin{lemma}[{cf.\ \cite[Proposition 12.2.2]{CLII}, \cite[Corollaire 12.2.4]{CLII}}]\label{lem.gpoid.orb.integral}
    Let $a_M \in \bAe{M}(k_0)$, $\lambda \in \Sigma_a$, and consider $[J_{a_M}(F) \backslash \mathcal{W}^\xi_\lambda]$, the action groupoid of $J_{a_M}(F)$ on the set $\mathcal{W}^\xi_\lambda$ of pairs $(g, \delta) \in G(\A)/K \times J_{a_M}(F))$ satisfying (\ref{eqn.wtd.stab.conds.g}) and (\ref{eqn.fixed.gpoid.jlambda}) above. The cardinality of this groupoid satisfies the following two formulae:
    \begin{enumerate}
        \item For any $\xi \in \a_T$, $$\# [J_{a_M}(F) \backslash \mathcal{W}^\xi_\lambda] = c(J_{a_M}) \int_{(J_{a_M}(\A) \backslash G(\A))^{\tau,-\lambda}} 1_D(\Ad(g)^{-1} X_{a_M}) w_M^\xi(g) \; dg,$$
        \item For any $\xi \in \a_T$ in general position, $$\# [J_{a_M}(F) \backslash \mathcal{W}^\xi_\lambda] = \frac{c(J_{a_M})}{\vol(\a_M/\cX_\star(M))} \int_{(J_{a_M}(\A) \backslash G(\A))^{\tau,-\lambda}} 1_D(\Ad(g)^{-1} X_{a_M}) v_M^G(g) \; dg$$
    \end{enumerate}
    where $(J_{a_M}(\A) \backslash G(\A))^{\tau,-\lambda}$ is the subset, considered above, of those cosets represented by $g$ such that $\inv_{a_M}(g) = -\lambda$ and
    \begin{equation}\label{eqn.def.of.cJ}
        c(J_{a_M}) := \vol(J(F_0) \backslash J(\A_0)^1) \cdot (\#\ker(J(F)_{\tau} \to J(\A)_{\tau})).
    \end{equation}
\end{lemma}

\begin{proof}
    The proof of (1) proceeds almost unchanged from \cite[Proposition 12.2.2]{CLII}.\footnote{The sign difference from \textit{loc.\ cit.} is a consequence of the fact that our geometric modification associated to $\lambda$ is represented by $z_\infty^{-\lambda}$, while the arithmetic Kottwitz map retains the valuation normalization $\kappa_J([z_\infty^\mu]) = [\mu]$, cf.\ Remark \ref{rmk.Kottwitz.signs}.} The essential new point is that our relative degree map $H_M : J(\A_0) \to \cX_\star(M)$ is \emph{surjective}. By the split rigidification $J|_{F_{0,\infty}} \cong T|_{F_{0, \infty}}$ at $\infty$, we can build a datum with specified geometric local degree: for any $\mu \in X_\star(T_\infty)$ we obtain an element $z_\infty^{-\mu}$, whose degree when paired with any $\chi \in X^\star(T_\infty)$ is $\langle \mu, \chi\rangle$. Viewing this as an adelic element trivial except at $\infty$, it suffices to check surjectivity of 
    $$
        X_\star(T_\infty) = \Hom(X^\star(T_\infty), \Z) \to \Hom(X^\star(M), \Z).
    $$
    This follows because $X^\star(M) \subset X^\star(T_\infty)$ is saturated; after all, it is characterized as the set of integral points satisfying linear equations.
    
    For the sake of completeness, we recall the rest of the argument. Our $J = J_{a_M}$ is $M$-elliptic, so in particular anisotropic mod center, and $J(F_0)\backslash J(\A_0)^1$ is compact, so finite volume. Lemme 12.2.5 applies directly in our setting to give finiteness of $\ker(J(F)_\tau \to J(\A)_\tau)$. In particular, $c(J)$ is finite. 

    Write $j = j_\lambda \in J(\A)$ for the specific adelic representative constructed above (which has Kottwitz invariant $-\lambda$) and write $\psi = \psi_j: G(\A) \to G(\A)$ for the map $g \mapsto g \tau(g)^{-1} j^{-1}$. The equation 
    $$
        \delta j \tau(g) = g \qquad \Longleftrightarrow \qquad \delta = g \tau(g)^{-1} j^{-1} = \psi(g)
    $$
    allows us to replace the data $\delta$ with the condition that $g \in \psi^{-1}(J(F))$. The stabilizer of such a pair is $J(F_0) \cap g K g^{-1}$, so 
    $$
        \# [J(F) \backslash \mathcal{W}^\xi_\lambda] = \sum_{g \in J(F) \backslash \psi^{-1}(J(F)) / K} |J(F_0) \cap g K g^{-1}|^{-1} 1_D(\Ad(g^{-1}) X_{a_M}) 1_{M,g}(\xi_M).
    $$
    Equipping $J(F)\backslash \psi^{-1}(J(F))$ with the appropriate $G(\A_0)$-invariant quotient measure as in \textit{loc.\ cit.}, this becomes an integral
    \begin{equation}\label{eqn.orbital.integral.formula.gpoid}
        \# [J(F) \backslash \mathcal{W}^\xi_\lambda] = \int_{J(F) \backslash \psi^{-1}(J(F))} 1_D(\Ad(g^{-1}) X_{a_M}) 1_{M,g}(\xi_M) \; dg.
    \end{equation}
    One decomposes this integral: $(J(F) \backslash J(\A))^\tau$ acts without fixed points on $J(F) \backslash \psi^{-1}(J(F))$, and the orbits are exactly the fibers of the natural map\footnote{Note that the fixed point equation $g \tau(g)^{-1} = \delta j$ implies that $\inv_{a_M}(g) = -\lambda$, with our Kottwitz convention.}
    $$
        J(F) \backslash \psi^{-1}(J(F)) \twoheadrightarrow (J(\A) \backslash G(\A))^{\tau, -\lambda}.
    $$
    The associated quotient is compatible with choices of measures, so Fubini's theorem allows the integral of (\ref{eqn.orbital.integral.formula.gpoid}) to be computed as a double integral. The inner integral is
    $$
        \int_{(J(F) \backslash J(\A))^\tau} 1_D(\Ad((hg)^{-1}) X_{a_M}) 1_{M,hg}(\xi_M) \; dh,
    $$
    but since all such $h$ lie in $J(\A)$ and centralize $X_{a_M}$, the left factor in the integrand is contant and can be pulled out. The right factor, on the other hand, satisfies
    $$
        1_{M,hg}(\xi_M) = 1_{M,g}(\xi_M + H_M(h)),
    $$
    as in Definition \ref{def.weights.CL}. The domain $(J(F) \backslash J(\A))^\tau$ admits a filtration with subquotients (1) $J(F_0) \backslash J(\A_0)^1$, (2) $J(\A_0)^1 \backslash J(\A_0)$, and (3) $\ker(J(F)_\tau \to J(\A)_\tau)$, all compatible with the choices of measures. 
    \begin{enumerate} \item The first integral contributes a constant factor $\vol(J(F_0) \backslash J(\A_0)^1)$
    \item The second integral contributes, by the surjectivity of the first paragraph, the sum over integral translates of $\xi_M$ which defines $w_M^\xi$.\footnote{The plus sign in Definition \ref{def.weights.CL} is correct, since $C_{hg} = C_g - H_M(h)$ and so $1_{M, hg}(\xi_M) = 1_{M,g}(\xi_M + H_M(h))$.}
    \item The third integral contributes, by $M(\A)$-invariance of $w_M^\xi$, a constant factor of $\# \ker(J(F)_\tau \to J(\A)_\tau)$.
    \end{enumerate}
    In summary, then, (\ref{eqn.orbital.integral.formula.gpoid}) becomes
    $$
        \# [J(F) \backslash \mathcal{W}^\xi_\lambda] = \vol(J(F_0) \backslash J(\A_0)^1) \cdot (\#\ker(J(F)_{\tau} \to J(\A)_{\tau})) \int_{(J(\A) \backslash G(\A))^{\tau, -\lambda}} 1_D(\Ad(g^{-1}) X_{a_M}) w_M^\xi(g) \; dg,
    $$
    which is (1).

    To deduce (2), it remains to see that for $\xi$ in general position, the two kinds of weighted orbital integrals differ by the constant factor $\vol(\a_M/\cX_\star(M))$. This follows \emph{mutatis mutandis} as in \cite[\S 11.10-13]{CLI}; indeed they follow closely arguments and constructions of Arthur which require no splitness hypotheses. We emphasize the few modifications necessary and provide citations:
    \begin{itemize}
        \item The calculations of \cite[\S 11.10]{CLI} make use of the sublattice $X_\star(M_{\sc}) \subset X_\star(M)$, which in the split case has the convenient property that $\Delta_P^\vee$ is a $\Z$-basis for each parabolic $P$ with Levi $M$. In the relative setting, this may fail to be the case, but the role of $X_\star(M_{\sc})$ is auxiliary: in \cite[\S 6]{arthurlocalTF}, he performs the same computation with an arbitrary sufficiently divisible lattice $R_M = N \Z(\Delta_P^\vee) \subset \cX_\star(M)$ (recall $G$ has anisotropic center). One then gets by the same method a corresponding analogue of \emph{loc.\ cit.}, Proposition 11.10.1:
        $$
            w_M^\xi(g) = \lim_{\Lambda \to 0} \sum_{\mu_0 \in \cX_\star(M)/R_M} \sum_{P \supset M} c_{P,N}(\Lambda)^{-1} e^{-\Lambda(H_P(g) + [\mu_0 + \xi_M]_{P,N})},
        $$
        where $c_{P,N}(\Lambda) := \prod_{\alpha \in \Delta_P} \left(e^{\Lambda(N \alpha^\vee)} - 1\right)$ and for any $\lambda \in \a_M$ the expression $[\lambda]_{P,N} \in R_M$ is defined with regards to the basis $\{N \alpha^\vee : \alpha \in \Delta_P\}$ for $R_M$, i.e.,
        $$
            \lambda - [\lambda]_{P,N}
        $$
        is a real linear combination of $N \alpha^\vee$ with coefficients in $(0, 1]$. 
        \item The formula for Arthur's weight functions as in \S 11.11 transfers directly to the relative setting. One should normalize $d_P$ by the covolume of the appropriate auxiliary lattice $R_M$, and the formula Proposition 11.11.1 adapts directly. 
        \item The general formalism for $(G,M)$-families as in \S 11.12 is developed for general reductive groups in \cite{ArthurGMFamilies} relatively and without splitness hypotheses. Replace Chaudouard--Laumon's definition of $w_P$ with
        $$
            w_{P,N}(\mu, \Lambda) = \frac{d_P(\Lambda)}{c_{P,N}(\Lambda)} e^{-\Lambda([\mu]_{P,N})},
        $$
        and obtain by the same proof the analog of Lemme 11.12.1 with conclusion
        $$
            \Lambda([\mu]_{P,N}) = \Lambda([\mu]_{P',N}),
        $$
        so that indeed $(w_{P,N})_P$ is a $(G,M)$-family. 
        \item In \S 11.13, the same finite-index cancellation which produces a volume factor
        $$
            \frac{\vol(\a_L/X_\star(L))}{\vol(\a_M/X_\star(M))}
        $$
        cancels any terms involving the auxiliary lattices $R_M \subset \cX_\star(M)$ and $R_L \subset \cX_\star(L)$, and we deduce a version of Lemme 11.13.2 of the form
        $$
            \sum_{\mu \in \cX_\star(M)/R_M} w_L(\mu + \xi_M) = \frac{\vol(\a_L/\cX_\star(L))}{\vol(\a_M/\cX_\star(M))} w_L^\xi(1).
        $$
        One requires the fact that $\cX_\star(M) \to \cX_\star(L)$ is surjective, which follows because $X^\star(L) \hookrightarrow X^\star(M)$ is saturated, just as argued at the start of this proof.
        \item With the previous modifications made, \S 11.14 goes through directly, and requires only $\a_G = 0$ for its general-position argument, so applies in our generality of anisotropic center. One deduces the desired constant factor $\vol(\a_M/\cX_\star(M))^{-1}$ as in Corollaire 11.14.3.
    \end{itemize}
    Details are left to the reader.
\end{proof}

\section{Stabilization and the endoscopic point count.}

Now we wish to compare the stable ($s = 1$) components of the cohomology with stable orbital integrals. Our formalism for stabilization can be more simply inductive than \cite[\S 13]{CLII}. Recall the recursive definition, as in \cite{ArthurStableI}.

\begin{definition}\label{def.stable.orb.integral.global}
    For any quasisplit reductive $H/X$ satisfying the hypotheses at the top of the chapter and semistandard Levi $M \subset H$, define the \emph{stable orbital integral} $S^H_M$ for $a_M \in \bAe{M}$ (implicitly, recursively) by
    $$
        J_M^{H,1}(a_M) = \sum_{s \in \cE_{M}(H)} \left| Z_{\widehat{H}'_s}^\Gamma / Z_{\widehat{H}}^\Gamma \right|^{-1} S^{H'_s}_M(a_M),
    $$
    where $\cE_M(G) = \cE_{1_M}(G)$ is the set of induced endoscopic groups for the trivial datum for $M$. That is, since the identity datum $1_H = (1, \rho_H)$ is included in the sum, the above identity is equivalent to
    \begin{equation}
        S^{H}_M(a_M) = J_M^{H,1}(a_M) - \sum_{\substack{s \in \cE_{M}(H),\\ s \neq 1}} \left| Z_{\widehat{H}'_s}^\Gamma / Z_{\widehat{H}}^\Gamma \right|^{-1} S^{H'_s}_M(a_M).\label{eqn.stable.def.recurse}
    \end{equation}
\end{definition}

A few remarks about the well-founded nature of this recursion are in order. By Remark \ref{rmk.ell.end.data.finite}, the sum appearing in the definition is finite. If a class $[s]\in \cE_M(H)$ satisfied $H'_s = H$, then $s$ would be central in $\widehat{H}$ and the class would be trivial, so every term on the right side of (\ref{eqn.stable.def.recurse}) is a proper endoscopic group. By Proposition \ref{prop.end.gp.sub.root.system}, any proper endoscopic group has smaller dimension than $H$, so the recursion has only finite depth.

More subtly, the definition of $J_M^{H,1}$ required our range of standing assumptions on $H$, and we need to check those transfer down to endoscopic groups.
\begin{itemize}
    \item We work always with quasisplit reductive $X$-groups, which is preserved by passage to endoscopic groups.
    \item Since our original $H$ is assumed to be arithmetic, i.e.\ defined over a fixed $X_0/k_0$, then by the Remark after Definition \ref{def.induced.end.datum}, the various $H'_s$ are likewise.
    \item Let $\Gamma_T := \im(\pi_1(X, \infty) \to \Aut_{\Z}(\ul{X}_\star(T)(\infty)))$, and $\Gamma_{T_{H'_s}}$ likewise for the corresponding maximal torus in $H'_s$. By the construction of the endoscopic group, one can embed $\W_{H'_s} \rtimes \Gamma_{T_{H_s'}} \subset \W_{H} \rtimes \Gamma_T$, so our standing assumption of good characteristic transfers from $H$ to its endoscopic groups.
    \item Since $H$ is assumed split at $\infty$, then $M \subset H$ is the same, so all $H_s'$ (which contain $M$ as a Levi by \ref{prop.ind.G'.cont.M'.Levi}) contain the same split maximal torus at $\infty$.
    \item Many of our geometric results require $H$ to have anisotropic center. This is inherited by the elliptic endoscopic groups, by Proposition \ref{prop.ell.end.finite}.
\end{itemize}

Combining the above, we can see that the class of admissible groups for all the results proven so far is closed under passage to elliptic endoscopic groups.

\begin{remark}
    The auxiliary weight factor in general position $\xi$ is not naturally inherited or in general position for endoscopic groups. This is not necessary: in our cohomological Theorem \ref{thm.coh.wfl} one considers only the elliptic part of the Hitchin fibration for the endoscopic groups, and this restriction is independent of a choice of $\xi$. When applying the inductive hypothesis, one can and must choose a different auxiliary $\xi$ in general position in order to apply our cohomological theorems to each endoscopic group. But such a $\xi$ certainly exists for each endoscopic group, and our results are all ultimately independent of the choice. In particular, the Definition \ref{def.s.orbital.integral.adelic} of $s$-orbital integrals above is stated in terms of Arthur's weight function and independent of $\xi$.
\end{remark}

With these preliminaries resolved, we can now inductively compute the Frobenius trace on the stable endoscopic side of Theorem \ref{thm.coh.wfl}.

\begin{theorem}[{cf.\ \cite[Th\'eor\`eme 14.1.1]{CLII}}]\label{thm.stable.point.count}
    Let $G$ a quasisplit reductive $X$-group with anisotropic center, $M$ a semistandard Levi, and $a_M \in \bAe{M}(k_0)$ with image $a \in \bA_G$. Suppose that $a \in \bA^{\good}$. Let $j_{ell} : \bAe{G} \to \bA_G$, an open immersion. We have an equality
    \begin{equation}\label{eqn.pt.count.stable.endo}
        \Tr(\tau^{-1}, ((j_{ell})_{!\star}(\pH^\bullet(Rf_{G,\star}^{G\el} \ol{\Q}_l)_1))_a) = \frac{c(J_{a_M}) \cdot q^{\deg(D) |\Phi^G|/2} |\pi_0(Z^\Gamma_{\widehat{M}})|}{\vol(\a_M/\cX_\star(M)) \cdot |\pi_0(\widehat{T}^{W_a})| \cdot |\pi_0(Z^\Gamma_{\widehat{G}})|} S_M^G(a_M),
    \end{equation}
    where $\pH^\bullet(Rf_{G,\star}^{G\el} \ol{\Q}_l)_1$ is the \emph{fine} stable summand of Theorem \ref{thm.s.decomp} and $\pi_0(\widehat{T}^{W_a}) = \widehat{T}^{W_a}/Z_{\widehat{M}}^{\Gamma,0}$ is the dual fine symmetry group.
\begin{proof} 
    We proceed by a streamlined version of the inductive strategy from \textit{loc.\ cit.} For each $H \in \cE_{M}(G)$ (elliptic endoscopic groups induced from the trivial endoscopic datum for $M$), define $\wtilde{S}_M^H(a_M)$ so that equation (\ref{eqn.pt.count.stable.endo}) holds with $\wtilde{S}_M^H(a_M)$ in place of $S_M^H(a_M)$. Note that by the considerations above, each such $H$ is an admissible group and admits a truncated Hitchin fibration with an $s$-decomposition, well-defined $s$-orbital integrals, etc. Additionally, for each $H$ the induced point from $a_M$ is $H$-good by Lemma \ref{lem.good.locus.compat.end}.(1).
    
    By induction on dimension, we may assume
    $$
        \wtilde{S}_M^H(a_M)={S}_M^H(a_M)
    $$
    for all proper endoscopic groups $H \neq G$. 

    We apply now the cohomological WFL to a collection of trivial endoscopic data for $M$. The pointed endoscopic data for $M$ whose underlying endoscopic group is $M$ itself are parameterized, modulo the coarse character equivalence, by $Z_{\widehat{M}}^\Gamma/(Z_{\widehat{M}}^0)^\Gamma$. For each such class, choose a representative $z \in Z_{\widehat{M}}^\Gamma$. Then $z$ and $\rho_M$ canonically define a pointed endoscopic datum $\fM_z$ for $M$ whose group is $M$,\footnote{In the convention of \S \ref{subsect.ind.endoscopy}, the refined data relative to $\Gamma$ is $\rho_{\fM} \gamma \mapsto [(1, \gamma)] \in \pi_0^{M,\Gamma}(z) \cong \pi_0^{M,\Gamma}(1)$. } but may have nontrivial coarse character $\ol{z} \in \widehat{T}/Z_{\widehat{M}}^0$. Applying Theorem \ref{thm.coh.wfl} to each $\fM_z$ and restricting to the stalk over $a$ gives a decomposition in terms of fine parameters over the coarse parameter (note by the support Lemma \ref{lem.M.G.comp.of.supports} that only the summands induced from the various $\fM_z$ appear):
    \begin{equation}
        (\pHi(R f_\star^{\xi, \Mpar} \ol{\Q})_{\ol{z}})_{a}^{\ss} \cong \bigoplus_{\substack{s \in \cE_{\fM_z}^f(G)\\ s\mapsto \ol{z}}} (j_{ell})_{!\star} (i_{\rho_s})_\star \pHishift(R f_{G_s', \star}^{G'_s \el} \ol{\Q}_l)_1(-r_s)_{a}^{\ss}.\label{eqn.stalk.wfl}
    \end{equation}
    We apply alternating trace of Frobenius to both sides (which is insensitive to semisimplification). One obtains by Theorem \ref{thm.s.summand.orbital.integral} that the left-hand side corresponding to $z$ is
    $$
        C \cdot q^{\deg(D) |\Phi^G|/2} J^{G,z}_M(a_M), \qquad \text{ where } \qquad C = \frac{c(J_{a_M})}{\vol(\a_M/\cX_\star(M)) \cdot |\widehat{T}^{W_a}/(Z_{\widehat{M}}^0)^\Gamma|}
    $$
    By Proposition \ref{prop.s.orbital.int.invariant.under.center} we have $J^{G,z}_M(a_M) = J^{G,1}_M(a_M)$ and these are all equal. On the other side, the shift $-2r_s$ is even, so does not affect the alternating trace, and the Tate twist contributes $q^{r_s}$. 
    
    Consider the extended fine index set
    $$
        \wtilde{\cE}_M^f(G) := \bigsqcup_{z \in Z_{\widehat{M}}^\Gamma/(Z_{\widehat{M}}^0)^\Gamma} \cE_{\fM_z}^f(G) = \{s \in Z^\Gamma_{\widehat{M}}/Z_{\widehat{G}}^{\Gamma,0} : G'_s \text{ elliptic}\}.
    $$
    Averaging the equations (\ref{eqn.stalk.wfl}) over $z \in Z_{\widehat{M}}^\Gamma/(Z_{\widehat{M}}^0)^\Gamma$ gives the combined identity
    \begin{equation}
        C \cdot q^{\deg(D) |\Phi^G|/2} J^{G,1}_M(a_M) = \frac{1}{|Z_{\widehat{M}}^\Gamma/(Z_{\widehat{M}}^0)^\Gamma|} \sum_{\wtilde{s} \in \wtilde{\cE}_{M}^f(G)} q^{r_s} \Tr(\tau^{-1}, (j_{ell})_{!\star} \pH^\bullet(Rf_{G'_s,\star}^{G'_s\el} \ol{\Q}_l)_{a, 1}).\label{eqn.combined.stable}
    \end{equation}
    Now we collect the right hand side down to a sum indexed by endoscopic parameters. There is a natural map to Arthur's index set from Definition \ref{def.induced.end.datum}
    $$
        \mu : \wtilde{\cE}_M^f(G) \to \cE_M(G)
    $$
    with the induced endoscopic group and corresponding terms depending only on the image. Indeed, shifting by $Z_{\widehat{G}}^\Gamma$ preserves not only the connected dual centralizer $\widehat{G}_{s}'$, but the full $Z_{\widehat{G} \rtimes \Gamma}(s)$. 
    So the terms on the right side of (\ref{eqn.combined.stable}) and the ellipticity of a given $\wtilde{s}$ depend only on $s = \mu(\wtilde{s})$. The fiber size can be calculated in two steps: passing from fine to coarse parameters gives a factor of 
    $$
        |\widehat{\Xi}_G| = \left|\frac{(Z_{\widehat{G}}^0)^\Gamma}{Z_{\widehat{G}}^{\Gamma,0}}\right|,
    $$
    and passing from coarse parameters to Arthur's endoscopic parameters gives a factor of
    $$
        \left|\frac{Z_{\widehat{G}}^\Gamma}{(Z_{\widehat{G}}^0)^\Gamma}\right|.
    $$
    We can then write (\ref{eqn.combined.stable}) as
    $$
        C \cdot q^{\deg(D) |\Phi^G|/2} J^{G,1}_M(a_M) = \frac{|Z_{\widehat{G}}^\Gamma/Z^{\Gamma,0}_{\widehat{G}}|}{|Z_{\widehat{M}}^\Gamma/(Z_{\widehat{M}}^0)^\Gamma|} \sum_{s \in {\cE}_{M}(G)} q^{r_s} \Tr(\tau^{-1}, (j_{ell})_{!\star} \pH^\bullet(Rf_{G'_s,\star}^{G'_s\el} \ol{\Q}_l)_{a, 1}).
    $$

    Substitute the definition of $\wtilde{S}_M^{G'_s}(a_M)$ for each endoscopic group $G'_s$ and cancel $C$ and $|Z^\Gamma_{\widehat{M}}/(Z^0_{\widehat{M}})^\Gamma|$ to obtain
    $$
        q^{\deg(D) |\Phi^G|/2} J^{G,1}_M(a_M) =  \sum_{s \in {\cE}_{M}(G)} \frac{|Z_{\widehat{G}}^\Gamma/Z_{\widehat{G}}^{\Gamma,0}|}{|Z_{\widehat{G'_s}}^\Gamma/Z_{\widehat{G'_s}}^{\Gamma,0}|} q^{r_s} q^{\deg(D)|\Phi^{G'_s}|/2} \wtilde{S}_M^{G'_s}(a_M).  
    $$
    We have $r_s = \frac{\deg D}{2}(|\Phi^G| - |\Phi^{G'_s}|)$, so the powers of $q$ agree and can be cancelled. Ellipticity implies $Z^{\Gamma,0}_{\widehat{G'_s}} = Z^{\Gamma,0}_{\widehat{G}}$ (indeed, both are zero under our anisotropic center hypothesis), and moreover $Z^\Gamma_{\widehat{G}} \subset Z^\Gamma_{\widehat{G'_s}}$, so the remaining factor can be rewritten
    $$
        \frac{|Z_{\widehat{G}}^\Gamma/Z_{\widehat{G}}^{\Gamma,0}|}{|Z_{\widehat{G'_s}}^\Gamma/Z_{\widehat{G'_s}}^{\Gamma,0}|} = |Z_{\widehat{G'_s}}^\Gamma/Z_{\widehat{G}}^\Gamma|^{-1}.
    $$
    By the inductive hypothesis, all terms satisfy $\wtilde{S}^{G'_s}_M(a_M) = S^{G'_s}_M(a_M)$ except possibly $G'_s = G$, so rearranging gives
    $$
        \wtilde{S}^{G}_M(a_M) = J_M^{G,1}(a_M) - \sum_{\substack{s \in \cE_M(G) \\ s \neq 1}} |Z_{\widehat{G'_s}}^\Gamma/Z_{\widehat{G}}^\Gamma|^{-1} S^{G_s'}_M(a_M)
    $$
    and comparing with (\ref{eqn.stable.def.recurse}) gives $\wtilde{S}^{G}_M(a_M) = S^{G}_M(a_M)$ by induction. 
\end{proof}
\end{theorem}

\section{The global WFL.}

With the results of the previous two sections, we can immediately deduce the global analog of the WFL for all quasisplit reductive $G$. 

\begin{theorem}[{The global Weighted Fundamental Lemma, cf.\ \cite[Th\'eor\`eme 15.1.1]{CLII}}]\label{thm.global.wfl}
    Let $G_0/X_0$ be a quasisplit reductive group such that $G_0$ is split at chosen point $\infty$ and the characteristic is good. Let $M_0 \subseteq G_0$ be a semistandard Levi with base change $M \subseteq G$. Let 
    $$
        \fM' = (s_M, \rho_{s_M})
    $$
    be an arithmetic elliptic endoscopic datum for $M$ with endoscopic group $M'$. Assume that $M'$ is split in the rational fiber $\infty$, and that $\Gamma$ is chosen large enough to split $M'$ as well as $G$. Let
    $$
        a_{M'} \in \bA_{M'}^{M'\el}(k_0)
    $$
    and denote $a_M = i_{\rho_{s_M}}(a_{M'}) \in \bA_{M}(k_0)$ and $a = \chi_G^M(a_M) \in \bA_G(k_0)$ its images. Assume 
    $$
        a \in \bA_G^{\good}.
    $$
    For every $s \in \cE_{\fM'}(G) \subset s_M Z_{\widehat{M}}^\Gamma/Z_{\widehat{G}}^\Gamma$, let $\fG'_s = (s, \rho_s)$ be the induced datum with endoscopic group $G'_s$. Then $M'$ is naturally a semistandard Levi of $G'_s$ and
    $$
        J^{G,{s_M}}_M(a_M) = |Z^\Gamma_{\widehat{M}'}/ Z_{\widehat{M}}^\Gamma| \sum_{s \in \cE_{\fM'}(G)} |Z_{\widehat{G'_s}}^\Gamma/ Z_{\widehat{G}}^\Gamma |^{-1} S_{M'}^{G'_s}(a_{M'}).
    $$

\begin{remark}\label{rmk.M'.split.at.infty}
    Note the hypothesis that $M'$ should be split at infinity. This is necessary to apply our previous point-counting arguments to the pair $M' \subset G'_s$. Of course, no such hypothesis exists in the local WFL. We will later absorb this into our globalization argument, showing that we can engineer an arbitrary local endoscopic setup $M'_v, M_v, G_v$ at one place on a curve while mandating that $M'$ should be split at $\infty$.
\end{remark}

\begin{proof}
     First, assume that $G$ has anisotropic center, so that our various geometric results apply. The claim that $M' \subset G'_s$ as a Levi is Proposition \ref{prop.ind.G'.cont.M'.Levi}. The scalar quotients appearing are finite by ellipticity of the endoscopic groups. We consider the stalk of the coarse summand $\pHi(R f_\star^{\xi, \Mpar} \ol{\Q}_l)_{\ol{s}_M}$ at $a$. By the last claim of Lemma \ref{lem.M.G.comp.of.supports}, we have
    $$
        \left(\pHi(R f_\star^{\xi, \Mpar} \ol{\Q}_l)_{\ol{s}_M}^{\ss}\right)_a = \left(\pHi(R f_\star^{\xi, \Mpar} \ol{\Q}_l)_{\ol{s}_M}^{\ss}\right)_{\fM', a}.
    $$
    We proceed by the same yoga of fine, coarse, and Arthur's endoscopic parameters as in the proof of Theorem \ref{thm.stable.point.count}. We will average the results of the cohomological WFL \ref{thm.coh.wfl} over twisted data
    $$
        \fM'_z = (s_M z, \rho_{s_M}), \qquad z \in Z_{\widehat{M}}^\Gamma/(Z_{\widehat{M}}^0)^\Gamma
    $$
    defining the same endoscopic group $M'$. Again by this we mean the endoscopic datum defined using the canonical isomorphism
    $$Z_{\widehat{M} \rtimes \Gamma}(s_M z) \cong Z_{\widehat{M} \rtimes \Gamma}(s_M)$$
    to transfer the refined torsor $\rho_{s_M}^\Gamma$ as in \S \ref{subsect.ind.endoscopy} to a $\Gamma$-torsor for $s_M z$, which projects to a usual endoscopic torsor that we abusively denote by $\rho_{s_M}$.
    
    The argument of Lemma \ref{lem.M.G.comp.of.supports} moreover implies that $W_a$ fixes each $s_M z$, and each has the same point-count, so our first fixed-point computation Theorem \ref{thm.s.summand.orbital.integral} applies and gives on the average
    \begin{equation}
        C_M \cdot q^{\deg(D) |\Phi^G|/2} J_M^{G,{s_M}}(a_M) = \frac{1}{|Z_{\widehat{M}}^\Gamma/(Z^0_{\widehat{M}})^\Gamma|} \sum_{\wtilde{s} \in \wtilde{\cE}_{\fM'}^f(G)} q^{r_s} \Tr(\tau^{-1}, (j_{ell})_{!\star} \pH^\bullet(R f_{G'_s, \star}^{G'_s\el} \ol{\Q}_l)_{1,a}).\label{eqn.global.wfl.avg}
    \end{equation}
    where $\wtilde{\cE}_{\fM'}^f(G) = \{s \in s_M Z^\Gamma_{\widehat{M}}/Z_{\widehat{G}}^{\Gamma,0} : G'_s \text{ elliptic}\}$ and
    $$
        C_M = \frac{c(J_{a_M})}{\vol(\a_M/\cX_\star(M)) \cdot |\widehat{T}^{W_a}/(Z^0_{\widehat{M}})^{\Gamma}|}.
    $$
    By the good locus Lemma \ref{lem.good.locus.compat.end}, the stable point-count Theorem \ref{thm.stable.point.count} applies to the pair $M' \subset G'_s$. The right side of (\ref{eqn.global.wfl.avg}) then becomes (multiplying by a factor of $|(Z^0_{\widehat{M'}})^\Gamma/Z_{\widehat{M'}}^{\Gamma,0}|$ on the numerator and denominator)
    $$
        \frac{1}{|Z_{\widehat{M}}^\Gamma/(Z^0_{\widehat{M}})^\Gamma|} \sum_{\wtilde{s} \in \wtilde{\cE}_{\fM'}^f(G)} q^{r_s} \frac{c(J_{a_{M'}}) \cdot q^{\deg(D) |\Phi^{G'_s}|/2} |Z^\Gamma_{\widehat{M'}}/(Z^0_{\widehat{M'}})^\Gamma|}
        {\vol(\a_{M'}/\cX_\star(M')) \cdot |\widehat{T}^{W_a}/(Z_{\widehat{M'}}^0)^\Gamma| \cdot |Z^\Gamma_{\widehat{G'_s}}/Z^{\Gamma,0}_{\widehat{G'_s}}|} S_{M'}^{G'_s}(a_{M'}).
    $$
    Several factors cancel:
    \begin{itemize} 
        \item the factors of $|\widehat{T}^{W_a}/(Z^0_{\widehat{M}})^{\Gamma}|$ and $|\widehat{T}^{W_a}/(Z_{\widehat{M'}}^0)^\Gamma|$ cancel to $|(Z^0_{\widehat{M'}})^{\Gamma}/(Z^0_{\widehat{M}})^{\Gamma}|$
    \item the powers of $q$ on the two sides are equal by the codimension computation as before, 
    \item $c(J_{a_M}) = c(J_{a_{M'}})$ because the regular centralizers are generically isomorphic, see \cite[Proposition 2.5.1]{ngofl}, and
    \item $\a_M = \a_{M'}$ and $\cX_\star(M) = \cX_\star(M')$, compatibly with the measures (this is essentially by definition of ellipticity; note that $\cX_\star(M)$ is dual to $(Z_{\widehat{M}}^\Gamma)^0$, and likewise for $M'$ -- these duals are equal by ellipticity). 
    \end{itemize}
    So (\ref{eqn.global.wfl.avg}) becomes, collecting the $|Z^\Gamma_{\widehat{G}}/Z^{\Gamma,0}_{\widehat{G}}|$ fine terms lying over each endoscopic parameter,
    $$
        J_M^{G,{s_M}}(a_M) =         \frac{|Z^\Gamma_{\widehat{M'}}/(Z^0_{\widehat{M'}})^\Gamma| \cdot |(Z^0_{\widehat{M'}})^{\Gamma}/(Z^0_{\widehat{M}})^{\Gamma}|}{|Z_{\widehat{M}}^\Gamma/(Z^0_{\widehat{M}})^\Gamma|}
        \sum_{s \in \cE_{\fM'}(G)} \frac{|Z^\Gamma_{\widehat{G}}/Z^{\Gamma,0}_{\widehat{G}}|}{ |Z^\Gamma_{\widehat{G'_s}}/Z^{\Gamma,0}_{\widehat{G'_s}}|} S_{M'}^{G'_s}(a_{M'}).
    $$
    Both remaining constants simplify to those in the statement of the Theorem.

    Now assume that $G_0$ is a general quasisplit reductive group (but otherwise satisfying our running hypotheses). By the same reduction as \cite[\S 8.6.1]{ngofl}, the quotient $\ol{G} = G/A_G^{\spl}$ has anisotropic center, and likewise $\ol{M} = M/A_G^{\spl}$ and $\ol{M}'$ and each $\ol{G}_s'$ descend, and the same for all recursively defined endoscopic groups. At the level of dual tori we have
    $$
        \begin{tikzcd}[ampersand replacement=\&]
        1 \ar[r] \& \widehat{\ol{T}} \ar[r] \& \widehat{T} \ar[r] \& \widehat{A_G^\spl} \ar[r] \& 1
        \end{tikzcd},
    $$
    so a given $s_M$ may not lie within $\widehat{\ol{T}}$. However, the $s$-orbital integrals are invariant under translation by $Z_{\widehat{M}}^\Gamma$ after Proposition \ref{prop.inv.sc.cover}, and the natural map $Z_{\widehat{M}}^\Gamma \to \widehat{A_G^\spl}$ is surjective because the dual map $$X^\star(\widehat{A_G^\spl}) \to X_\star(A^\spl_G) \to \pi_1(M)_{\Gamma} = X^\star(Z_{\widehat{M}}^{\Gamma})$$ can be checked to be injective. So assume $s_M \in \widehat{\ol{T}}$. We can likewise identify their torsors: we have isomorphisms
    $$
        \pi_0^{\ol{M},\Gamma}(\ol{s_M}) \cong \pi_0^{M,\Gamma}(s_M), \qquad \pi_0^{\ol{G},\Gamma}(\ol{s}) \cong \pi_0^{G,\Gamma}(s_M)
    $$
    compatible with the induced construction of Definition \ref{def.induced.end.datum}.

    The weighted $s$-orbital integrals $J^{G,s}_M(a_M)$ and $J^{\ol{G},s}_{\ol{M}}(\ol{a_M})$ are equal with quotient-compatible Haar measures, by the following observations: 
    \begin{itemize} 
    \item the root systems and therefore the relevant cocharacter spaces $\a_M^G$ are equal, 
    \item the value of the invariant lies in the image of the cocharacter lattice of $T_{M_{\sc}}$ by Proposition \ref{prop.inv.sc.cover}, hence is insensitive to our central quotient,
    \item the parabolics and Levis are in bijection, 
    \item the $H_P$ are equal modulo the common central direction, so $v_M^G(g) = v_{\ol{M}}^{\ol{G}}(g)$,
    \item since $A_G^{\spl} \subset J_{a_M}$, the adelic quotient in the domain of integration is canonically isomorphic, and
    \item $X$ necessarily splits as $\ol{X} + Z$ for central component $Z$ which is fixed under conjugation. Since $a$ is a global section of $\fc_{G,D}$, one can check that $Z$ satisfies $D$-integrality, so $X$ is $D$-integral if and only if $\ol{X}$ is.
    \item The adjoint-equivariant Lie algebra splitting $\Lie(G) = \Lie(\ol{G}) \oplus \Lie(A_G^{\spl})$ induces a corresponding decomposition of the Hitchin base, compatible with the constructions of the regular centralizer, cameral cover, and the constant $\delta$. In particular, $\ol{a}$ is still good. 
    \end{itemize}
    The constants appearing in the recursive definition of the stable orbital integrals are likewise invariant under passage to the quotient, and the endoscopic recursion is likewise preserved, so the stable orbital integrals are equal. 
    
    Applying the theorem for $\ol{G}$, $\ol{M}$, $\ol{M}'$, $\ldots$, we therefore obtain the global WFL for $G$, $M$, $M'$, $\ldots$, except with respect to the choice of Haar measure on $J(\A_0)$ chosen to be compatible with the \textit{quotient}, which may not correspond to our conventions above. Rescaling by some factor $c > 0$ to the maximal-compact-volume-one normalization of Definition \ref{def.Haar.measure.max.cmpt.convention} will rescale every weighted orbital integral by $c^{-1}$, and the same factor carries through the stable recursion to apply to every stable orbital integral, so it cancels from the identity.
\end{proof}
\end{theorem}

\chapter{Global-to-local and the Weighted Fundamental Lemma}\label{sect.local.global.WFL}

\section{Local orbital integrals and the statement of the WFL.}\label{subsect.local.setup}
In this section, let $F_{0,v} = k_v (\!(z_v)\!)$ and $\O_{0,v} = k_v [\![z_v]\!]$ for finite field $k_v = \F_{q^d}$, and denote $F_v$, $\O_v$ the completed base changes to $k = \ol{\F}_q$. Assume always that $p$ does not divide the (split) cardinality of the Weyl group $|W_G|$. Let $G/\O_{0,v}$ a quasisplit reductive group (unramified), split by some finite quotient of $\Gal(F_v/F_{0,v})$. Fix $M \subset G$ a Levi. Fix a Kostant section over $\O_{0,v}$
$$
    \epsilon_M : \fc_M \to \fm.
$$

We now imitate the weighted $s$-orbital formalism, as in \S \ref{subsect.weighted.s.orbital.global}, locally. Let $b \in \fc_M(\O_{0,v}) \cap \fc_M^{\Greg}(F_{0,v})$ be a generically $G$-regular semisimple parameter with $Y_b = \epsilon_M(b)$ its Kostant image. Consider the local regular centralizer $J_b = Z_M(Y_b) = Z_G(Y_b)$ over $F_{0,v}$. Choosing a geometric point $\wtilde{b}$ of the local cameral cover 
$$
    \ft \to \fc_M
$$
over $b_{\eta}$ (and in particular, a geometric point over $v$), we obtain identifications
$$
    X_\star(J_{b,F_{v}}) \cong X_\star(J_{b,F^{\sep}_{v}}) \cong \bX_\star := X_{\star}(\T).
$$
The action of the full Galois group $\Gamma_{0,v} = \Gal_{F_{0,v}}$ can permute the fibers of the cameral cover; as in \S \ref{subsect.Ngo.monod.comp}, these package together into a local monodromy homomorphism
$$
    \mu_{b,v} : \Gamma_{0,v} \to \W_{M} \rtimes \Out(\M).
$$
We define the \emph{local monodromic subgroup} $W_{b,v}$ to be the image of this homomorphism. Then
$$
    X_\star(J_{b,F_v^{\sep}})_{\Gamma_{0,v}} \cong (\bX_\star)_{W_{b,v}},
$$
so dually
$$
    \widehat{J}_b^{\Gamma_{0,v}} \cong \widehat{T}^{W_{b,v}}.
$$
These constructions are compatible with the global ones upon restriction to a geometric point lying over $v$.

We can give a local version of the Kottwitz invariant construction: writing $\tau_v$ for the arithmetic Frobenius, we have an exact sequence of pointed sets
$$
    \begin{tikzcd}[ampersand replacement=\&]
        1 \ar[r] \& J_b(F_{0,v}) \ar[r] \& G(F_{0,v}) \ar[r] \& (J_b(F_v) \backslash G(F_v))^{\tau_v} \ar[r] \& \ker(J_b(F_v)_{\tau_v} \to G(F_v)_{\tau_v}) \ar[r] \& 1
    \end{tikzcd}.
$$
By \cite[\S7.6]{kottwitzisocrystalsii}, or by an argument parallel to Proposition \ref{prop.kottwitz.Ja}, we have an isomorphism
$$
    J_b(F_v)_{\tau_v} \cong X_\star(J_{b, F_v^{\sep}})_{\Gamma_{0,v}} \cong (\bX_\star)_{W_{b,v}}.
$$
As in the global setting, equipping $J(F_{v})_{\tau_v}$ with the counting measure and $J(F_{0,v})\backslash G(F_{0,v})$ with the quotient of Haar measures normalized to give the maximal compact in $J(F_{0,v})$ and to give $G(\O_{0,v})$ measure 1, the exact sequence above induces a right-$G(F_{0,v})$-invariant measure on $(J_b(F_v) \backslash G(F_v))^{{\tau_v}}$.

\begin{definition}\label{def.local.s.integral}
    Define the invariant map to be the composition
    $$
        \inv_{b,v} : (J_b(F_v)\backslash G(F_v))^{\tau_v} \to \ker(J_b(F_v)_{\tau_v} \to G(F_v)_{\tau_v}) \subseteq J_b(F_v)_{\tau_v} \cong (\bX_\star)_{W_{b,v}}.
    $$
    For $Y_b$ and $M$ as above and $s \in \widehat{J}_b^{\Gamma_{0,v}} \cong \widehat{T}^{W_{b,v}}$, the \emph{local weighted $s$-integral} is defined to be
    $$
        J_{M,v}^{G,s}(b) = J_{M,v}^{G,s}(Y_b)= |D^G(Y_b)|_{F_{0,v}}^{1/2} \int_{(J_b(F_v) \backslash G(F_v))^{{\tau_v}}} \langle s, \inv_{b,v}(g) \rangle 1_{\fg(\O_v)}(\Ad(g)^{-1} Y_b) v_M^{G, \loc}(g) \; dg.
    $$
    The pairing is that between $\widehat{T}^{W_{b,v}}$ and $(\bX_\star)_{W_{b,v}}$ (or more invariantly, between $\widehat{J}_b^{\Gamma_{0,v}}$ and $ X_\star(J_{b, F_v^{\sep}})_{\Gamma_{0,v}}$).
    
    The function $v_M^{G, \loc}$ is the intrinsic local version of Arthur's weight function from Definition \ref{def.arthur.wt.ftn}. We have local functions $H_{P,v}^{\loc}$ valued in $\a_{M,v} := \Hom_{\Z}(X^\star_{F_0,v}(M), \R)$, defined for each $F_{0,v}$-rational character $\lambda$ by
    $$
        \langle \lambda, H_{P,v}^{\loc}(g) \rangle = - \val_v(\lambda(p)).
    $$
    The function $v_M^{G, \loc}$ is then the volume in $\a_{M,v}^G$ of the convex hull of the projections of $-H_{P,v}^{\loc}(g)$, ranging over $F_{0,v}$-parabolics $P$ with Levi $M$.
\end{definition}

\begin{remark}
    Note that $H_{P,v}^{\loc}$ and $v_M^{G, \loc}$ differ from the global definition appearing in \S 6 and Definition \ref{def.arthur.wt.ftn} by a projection (there may be more local parabolics) and a scalar: for $g_v \in G(F_{0,v})$ at a closed point $v$ and global $P$,
    $$
        H_{P,v}(g_v) = \deg(k_v/k_0) \on{proj}_{\a_M} H_{P_v,v}^{\loc}(g_v).
    $$
    Relative to a geometric point $H_{P, \ol{v}}$, the scalar equals 1 and the two differ only up to projection.
\end{remark}

\begin{proposition}\label{prop.comparison.local.s}
    Choose an unramified local endoscopic datum $\fM'_v = (s_M, \rho_{v,s_M})$ with endoscopic group $M'/\O_{0,v}$. Let $Y' \in \fm'(F_{0,v})$ semisimple and $G$-regular. Then, writing $b' = \chi_{M'}(Y') \in \fc_{M'}(F_{0,v})$ with image $b \in \fc_M(F_{0,v})$, denote
    $$
        Y_0' := \epsilon_{M'}(b'), \qquad Y_0 := \epsilon_{M}(b).
    $$
    we have
    $$
        J^{G,s_M}_{M,v}(Y_0) = \sum_{[Y]} \Delta_{M',M}(Y',Y) |D^G(Y)|^{1/2}_{F_{0,v}} \int_{J_Y(F_{0,v})\backslash G(F_{0,v})} 1_{\fg(\O_{0,v})}(\Ad(g)^{-1} Y) v_M^{G,\loc}(g) \; dg,
    $$
    the sum taken over $M(F_{0,v})$-conjugacy classes of $G$-regular semisimple elements in $\fm(F_{0,v})$. Here $\Delta_{M', M}$ is Waldspurger's variant of the Langlands--Shelstad transfer factor for Lie algebras with the factor $\Delta_{IV}$ omitted, as in \cite{waldFLimpliquetransfert}.
    \begin{proof}
        This is standard and follows as at the start of the proof of \cite[Theorem 19.1.1]{CLII}. For completeness, we give more details. To compare the two formulae, we need to translate the invariant formalism above into the more traditional language of stable conjugacy. The invariant defined above is valued in $\ker(J_b(F_v)_{\tau_v} \to G(F_v)_{\tau_v})$, where the subscripts denote ${\tau_v}$-conjugacy classes. A standard argument in Galois cohomology shows that 
        $$
            \ker(J_b(F_v)_{\tau_v} \to G(F_v)_{\tau_v}) \cong \ker(H^1(F_{0,v}, J_b) \to H^1(F_{0,v}, M))
        $$
        We give an argument in Lemma \ref{lem.inv.red.stable.conj} for exposition's sake.

        Choosing the basepoint $Y_0$, we can identify $\ker(H^1(F_{0,v}, J_b) \to H^1(F_{0,v}, M))$ with rational classes in the $M$-stable class of $Y_0$. We denote $\inv(Y_0, Y)$ for the image of $\inv_{s}(Y)$ along this identification. The value of the invariant decomposes the domain of the integral across rational classes $[Y]$ in the $M$-stable class of $Y_0$:
        $$
            (J_{Y_0}(F_v)\backslash G(F_v))^{{\tau_v}} = \bigsqcup_{[Y]} (J_Y(F_{0,v}) \backslash G(F_{0,v})).
        $$
        So, using invariance of $D^G$ under stable conjugacy,
        $$
            J_{M,v}^{G,s_M}(Y_0) = \sum_{[Y]} \langle s_M, \inv(Y_0, Y) \rangle |D^G(Y)|^{1/2}_{F_{0,v}} \int_{J_Y(F_{0,v}) \backslash G(F_{0,v})} 1_{\fg(\O_{0,v})}(\Ad(g)^{-1} Y) v_M^{G, \loc}(g) \; dg.
        $$

        The comparison between our pairing and the transfer factor of Langlands--Shelstad was done in \cite[\S 1.5-1.11]{ngofl}. We outline the steps. For $s_M \in \widehat{T}^{W_{b,v}}$, the character $\langle s_M, - \rangle : J_b(F_v)_{\tau_v} \to \C^\times$ is given by passing through the Kottwitz identification 
        $$
            J_b(F_v)_{\tau_v} \cong X_\star(J_{b, F_{0,v}^{\sep}})_{\Gamma_{0,v}} \cong X_\star(\T)_{W_{b,v}}. 
        $$
        Under the isomorphism above, the corresponding character of $\ker(H^1(F_{0,v}, J_b) \to H^1(F_{0,v}, M))$ is simply the restriction to the kernel of the Newton map, and can be checked to be equal to the Tate--Nakayama dual with the same normalization as in \cite[\S 1.6, 1.9]{ngofl}.
        
        By Kottwitz's normalization of the transfer factor \cite{kottwitztransferliealg}, $\Delta_{M',M}(Y_0', Y_0) = 1$, and of course $\langle s_M, \inv(Y_0, Y_0) \rangle  = 1$. The Langlands--Shelstad variation formula \cite[cf.\ \S 3.4, \S 3.7]{lstransferfactors} shows that the transfer factor varies by a multiplicative factor of the Tate--Nakayama pairing
        $$
            \frac{\Delta_{M',M}(Y_0', Y)}{\Delta_{M',M}(Y_0', Y_0)} = \langle s_M, \inv(Y_0, Y)\rangle
        $$
        and so
        $$
            \Delta_{M',M}(Y_0', Y) = \langle s_M, \inv(Y_0, Y) \rangle
        $$
        over all classes.
        Moreover, the factor is invariant under stable conjugacy in the endoscopic variable, so we may replace $Y_0'$ with $Y'$. The result is the desired formula.
    \end{proof}
\end{proposition}

\begin{lemma}\label{lem.inv.red.stable.conj}
    For $J$ a torus inside Levi $M \subseteq G$ as above, we have a canonical bijection
    $$
        \ker(J(F_v)_{\tau_v} \to G(F_v)_{\tau_v}) \cong \ker(H^1(F_{0,v}, J) \to H^1(F_{0,v}, M)).
    $$
    \begin{proof}
        First, we claim that $M(F_v)_{\tau_v} \to G(F_v)_{\tau_v}$ has trivial kernel, so the kernel in question reduces to
        $$\ker(J(F_v)_{\tau_v} \to M(F_v)_{\tau_v}).$$ Choose a parabolic $P$ for $M$. The existence of the projection $P \to M$ and section $M \to P$ ensures the map $M(F_v)_{\tau_v} \to P(F_v)_{\tau_v}$ is injective, so it suffices to check injectivity with $M$ replaced by $P$. Assume $p \in P(F_v)$ and for some $g \in G(F_v)$ we have
        \[
            g p {\tau_v}(g^{-1}) = 1. 
        \]
        Then $gp = {\tau_v}(g)$, so $gP \in (G/P)(F_{0,v})$. Since $G \to G/P$ is surjective on points over any field (indeed, it has Zariski-local sections, see for example \cite[Prop. 7.2.1]{conrad2}) we can find a $g_0 \in G(F_{0,v})$ and a $p_0 \in P(F_v)$ such that $g_0 p_0 = g$. We have 
        \[
            p = g^{-1} {\tau_v}(g) = p_0^{-1} {\tau_v}(p_0),
        \]
        as desired. 

        Second, we claim natural isomorphisms
        $$
            \ker(J(F_v)_{\tau_v} \to M(F_v)_{\tau_v}) \cong \ker(H^1_{un}(F_{0,v}, J) \to H^1_{un}(F_{0,v}, M)) \cong \ker(H^1(F_{0,v}, J) \to H^1(F_{0,v}, M))
        $$
        so that the invariant above is constant on rational conjugacy classes within a given stable conjugacy class. The right isomorphism follows from Steinberg's theorem and the inflation-restriction sequence. The first is an exercise in Kottwitz's Newton map formalism (\cite{kottwitzisocrystals}, \cite{kottwitzisocrystalsii}): there is a natural map for any reductive $H$
        $$
            \nu_H : H(F_v)_{\tau_v} \to (\Hom_{F_v}(\D, H)/H(F_v))^{{\tau_v}},
        $$
        where $\D$ is the pro-torus with character group $\Q$. Here naturality means that given $H \to G$, the natural map on ${\tau_v}$-conjugacy classes corresponds to post-composition of homomorphisms. The kernel of the Newton map is $H^1_{un}(F_{0,v}, H)$ by, for example, \cite[Lemma 4.7]{GisoshtukasHK}. If $[j] \in \ker(J(F_v)_{\tau_v} \to M(F_v)_{\tau_v})$, then $\nu_M([j]) = 0$, which means that the postcomposed homomorphism $ i\circ\nu_{J} : \D \to J \to M$ is trivial. But $J \to M$ is injective, so $\nu_{J}([j]) = 0$. 
    \end{proof}
\end{lemma}

As in the global case, we have our recursive definition of stable orbital integrals. Recall that $\Gamma = \Gamma(F_v/F_{0,v})$ captures the monodromy of our unramified group $G$. Locally, we will only consider \emph{unramified} elliptic endoscopic data $\fM'$ for semistandard Levi $M \subset G$. As defined in Definition \ref{def.induced.end.datum}, it is clear that the induced endoscopic data $\cE_{\fM'}(\fG)$ are also split by the $\Gamma$-extension. So taking $\Gamma$ to be either $\Gal(F_v/F_{0,v})$ or a finite quotient corresponding to an extension large enough to split both $G$ and $M'$, all pinned actions on the ambient groups and endoscopic data will factor through $\Gamma$, and we can write
$$
    (-)^{\Gamma_{0,v}} = (-)^{\Gamma}
$$
unambiguously.

\begin{definition}\label{def.stable.integral.local}
    For any quasisplit reductive $H/X$ satisfying the hypotheses at the start of the chapter and semistandard $M \subset H$ with $Y_b$ as above, define the \emph{stable orbital integral} $S_{M,v}^H$ for $Y_b$ implicitly by
    $$
        J_{M,v}^{H,1}(Y_b) = \sum_{s \in \cE_{M}(H)} |Z_{\widehat{H}'_s}^\Gamma/Z_{\widehat{H}}^\Gamma |^{-1} S_{M,v}^{H_s'}(Y_b),
    $$
    where $\cE_M(G) = \cE_{1_M}(G)$ is the set of induced endoscopic groups for the trivial datum for $M$. Equivalently,
    $$
        S_{M,v}^H(Y_b) = J_{M,v}^{H,1}(Y_b) - \sum_{\substack{s \in \cE_M(H) \\ s \neq 1}} | Z_{\widehat{H}'_s}^\Gamma/Z_{\widehat{H}}^\Gamma|^{-1} S_{M,v}^{H_s'}(Y_b).
    $$
    This recursive definition is well-founded as in the discussion after Definition \ref{def.stable.orb.integral.global}.
\end{definition}

We now state the version of the Weighted Fundamental Lemma which is the primary goal of this work. 

\begin{theorem}[The Weighted Fundamental Lemma for Lie algebras]\label{thm.wfl}
    Assume that $p$ does not divide the (split) cardinality of the Weyl group $|W_G|$. Choose an unramified local endoscopic datum $\fM'_v = (s_M, \rho_{v,s_M})$ with endoscopic group $M'/\O_{0,v}$, and let $\Gamma_v$ be an unramified quotient prime to $p$ large enough to split $G$ and $M'$. Let $b' \in \fc_{M'}(\O_{0,v}) \cap \fc_{M'}(F_{0,v})^{\Greg}$ be a generically $G$-regular parameter, and write $b \in \fc_M(\O_{0,v})$ for its image under the natural map $\fc_{M'} \to \fc_M$. Denote
    $$
        Y_{b'} = \epsilon_{M'}(b') \in \fm'(\O_{0,v}), \qquad Y_{b} = \epsilon_{M}(b) \in \fm(\O_{0,v})
    $$
    for the Kostant lifts. Then
    \begin{equation}
        J_{M,v}^{G,s_M}(Y_b) = |Z_{\widehat{M'}}^{\Gamma_v}/Z_{\widehat{M}}^{\Gamma_v} |\sum_{s \in \cE_{\fM'}(G)} |Z_{\widehat{G}'_s}^{\Gamma_v}/Z_{\widehat{G}}^{\Gamma_v} |^{-1} S_{M',v}^{G_s'}(Y_{b'}).\label{eqn.local.wfl}
    \end{equation}
\end{theorem}

In the case where $G$ and the endoscopic datum $\fM'$ are split, this is \cite[Th\'eor\`eme 19.1.1]{CLII}.

\section{Splitting lemmas.}\label{subsect.splitting.lemmas}

In \cite[\S 16]{CLII}, the authors develop a full formalism for scindage and descent over finite sets of places on the curve $X_0$. Some part of the same formalism (in our relative, nonsplit setting) is unavoidable, but for the sake of brevity, and we hope, transparency of the ideas, we try for a minimalist path to the WFL. 

In this section and the next, we keep the running global geometry ($X_0, D_0, \ldots$) described at the top of \S \ref{sect.global.point.counts}, in particular $G$ split at infinity, (but not necessarily having anisotropic center). In particular, $M \subset G$ will be fixed, and we fix also $\fM' = (s_M, \rho_{s_M})$ a global elliptic endoscopic datum for $M$ also split at $\infty$. Until the end of this section, we write $\Gamma$ for a finite quotient of the global arithmetic Galois group $\pi_1(X_0, \infty)$ through which $G, \fM'$, and consequently every induced datum from $\fM'$ factor. By our running splitting hypotheses, its image in relevant data over $\infty$ coincides with the geometric Galois group $\pi_1(X, \infty)$. 

Assume now that we have a parameter $a_{M'} \in \bA_{M'}^{M'\el}(k_0)$ with images $a_M \in \bA_M$, $a \in \bA_{G}$. We will need Arthur's coefficient $d$ defined in \cite[\S 7-9]{arthurinvlocal}.

\begin{definition}\label{def.Arthur.d}
    For semistandard Levis
    $$
        R \subset L_1, L_2 \subset H,
    $$
    define $d_R^H(L_1, L_2) = 0$ unless the natural map
    $$
        \a_R^{L_1} \oplus \a_R^{L_2} \to \a_R^H
    $$
    is an isomorphism, in which case $d_M^H(L_1, L_2)$ is the absolute value of the (real) determinant of this natural map with respect to the chosen Euclidean structures on each (cf.\ Definition \ref{def.cochar.space.measure.choice}). In particular, $d_R^H(H,R) = d_R^H(R,H) = 1$. 
\end{definition}

\begin{remark}\label{rmk.end.gp.Arthur.d}
    Recall that if $G'$ is any elliptic endoscopic group for $G$, we have a natural isomorphism $\a_{G'} \cong \a_G$, and that we equipped $\a_{G'}$ with the measure which makes this map an isometry. If we likewise have elliptic endoscopic groups $R'$ for $R$, $L_1'$ for $L_1$, and $L_2'$ for $L_2$, then by construction we have
    $$
        d_{M}^H(L_1, L_2) = d_{M'}^{H'}(L_1', L_2').
    $$
\end{remark}

\begin{definition}\label{def.s.integral.V}
For any nonempty set of closed points $V \subset |X_0|$ and $s \in \widehat{T}^{W_{a_M}}$, define the \emph{weighted $s$-orbital integral on $V$} by
$$
    J_{M,V}^{G,s}(a_M) := q^{-\deg(D|_V) |\Phi^G|/2} |D^G(Y_{a_M})|_V^{1/2} \int_{(J_{a_M}(\A_V) \backslash G(\A_V))^{\tau}} \langle s, \inv_{a_M}(g) \rangle 1_D(\Ad(g^{-1}) Y_{a_M}) v_M^G(g) \; dg,
$$
where $|-|_V$ denotes the product of the local absolute values and $Y_{a_M} \in \fm(\A_V)$ is the image of $a_M|_V$ under the Kostant section. For $g \in G(\A_V)$, the weight function $v_M^G(g)$ is the global one, evaluated by adjoining the identity outside $V$. Define $S_{M,V}^{G}$ by the same recursion of Definition \ref{def.stable.orb.integral.global}.
\end{definition}

\begin{remark}\label{rmk.partial.local.comparison.J}
    The above definition is built to interpolate between the global $s$-orbital integral of Definition \ref{def.s.orbital.integral.adelic} when $V = |X_0|$ (note that all the same normalizations of measures were chosen) and the local $s$-orbital integral of Definition \ref{def.local.s.integral} when $V = \{v\}$ with residue field $k_v$. However, since the weight $v_M^G$ is always the global one rather than the intrinsic local weight $v_M^{G, \loc}$, the comparison is not exact in general: they measure volumes of convexes lying in different spaces! When $D_v = \Gamma$, the projection $\ga_M^{D_v} \to \ga_M^\Gamma = \a_M$ distinguishing the two weights is the identity, and so the degree scalar contributes a factor of $\deg(k_v/k_0)^{\dim \a^L_M}$ to the orbital integrals. That is, for $v \not\in \on{supp}(D)$ we have
    $$
        J^{L,s}_{M, \{v\}}(a_M) = \deg(k_v/k_0)^{\dim \a_M^L} J_{M,v}^{L,s}(Y_{a_M}),
    $$
    where the left side is the partial adelic integral on $\{v\}$ and the right side is the local orbital integral for the local Kostant image $Y_{a_M} \in \fm(F_v)$. The same comparison holds for the stable integrals, by induction through their defining recursion: by ellipticity, every term has the same scalar factor, and the endoscopic coefficients agree when $D_v = \Gamma$.
\end{remark}

\begin{remark}\label{rmk.local.renorm.D}
    If $v$ lies on the support of the divisor $D$, the discrepancy of indicator functions can also be absorbed: writing $Y:= z_v^{d_v} Y_{a_M}$, for $d_v$ the multiplicity of $D$ at $v$, we have
    $$
        1_{D}(\Ad(g^{-1}) Y_{a_M}) = 1_{\fg(\O_v)}(\Ad(g^{-1}) Y).
    $$
    Moreover, the scalar multiplication does not change the centralizer, the invariant map, or Arthur's weight. Our choice of normalizing factor exactly cancels the change of the discriminant valuation on replacing $Y_{a_M}$ with $Y$, so combining with the previous remark gives when $D_v = \Gamma$
    $$
        J^{L,s}_{M, \{v\}}(a_M) = \deg(k_v/k_0)^{\dim \a^L_M} J_{M,v}^{L,s}(Y),
    $$
    where the left side is the partial adelic $s$-orbital integral and the right side is the local definition. The same renormalization passes through the recursion to the stable orbital integrals.
\end{remark}

The inductive structure we will use is the following.
\begin{definition}\label{def.global.wfl.V.places}
    For arbitrary $V \subseteq |X_0|$, we say that our \emph{the $\fM'$-WFL holds on $V$ for $M \subset G$} to mean that for every semistandard Levi $M \subseteq L \subseteq G$, the WFL holds for the $V$-orbital integrals, with $L$ in place of $G$ and $\Gamma$ the \emph{global} monodromy action, i.e.,
    \begin{equation}
        J^{L,{s_M}}_{M,V} (a_M) = |Z^\Gamma_{\widehat{M}'}/ Z_{\widehat{M}}^\Gamma| \sum_{s \in \cE_{\fM'}(L)} |Z_{\widehat{L'_s}}^\Gamma/ Z_{\widehat{L}}^\Gamma |^{-1} S_{M', V}^{L'_s}(a_{M'}).\label{eqn.V.wfl}
    \end{equation} 
    Say that $V$ is \emph{peelable} if in addition $J^{M,s_M}_{M, V}(a_M) \neq 0$. See Lemma \ref{lem.peeling} below for the meaning of the name.
\end{definition}

\begin{remark}\label{rmk.monodr.local.global.discrep.wfl}
    We warn that despite appearances, when $V$ is a single place the above definition does \emph{not} automatically capture the local WFL of Theorem \ref{thm.wfl}. In addition to the discrepancies in the orbital integrals discussed in Remark \ref{rmk.partial.local.comparison.J},
    the identity is also stated using \emph{global} endoscopic parameters. The $\Gamma$ appearing in the above definition is the global monodromy of $G$ and $\fM'$, while the relevant Galois group to recover the actual WFL at $v$ is $D_v \subseteq \Gamma$, the decomposition group. Not only would this give us the wrong constants, in general, but it also gives us only data corresponding to some of the endoscopic terms: a parameter
    $$
        s = s_M z, \qquad z \in Z_{\widehat{M}}^{D_v}
    $$
    need not satisfy $z \in Z_{\widehat{M}}^{\Gamma}$, so defines a perfectly valid induced local endoscopic datum which does not extend to the global induced setting. The larger center can also identify multiple global endoscopic parameters, and can change ellipticity. Moreover, even the orbital integrals appearing are not exactly local orbital integrals: the Arthur weight function in question is the global one, formed only from global parabolics and heights $H_{P,v}$ which pair with $\Gamma$-invariant characters, rather than all $D_v$-invariant characters. 

    We will avoid these concerns in our globalization step below by engineering our global monodromy to satisfy $D_v = \Gamma$. This is sufficient to ensure the fixed-center constants are correct, the sum ranges correctly over all local endoscopic parameters, and the weight function pairs with all local characters, differing only by a scalar. Moreover, for each parameter $s$, the construction of induced endoscopic data $(s, \rho_s)$ from Definition \ref{def.induced.end.datum} is compatible with restriction to $D_v$ as in Remark \ref{rmk.ind.endo.compatibilities}. All global induced endoscopic groups will therefore correspond to the expected local induced groups upon restriction to $\Spec \O_{0,v}$, and their ellipticity conditions and index sets will naturally biject. The scalar factor $\deg(k_v/k_0)^{\dim \a_M^L}$ then cancels from both sides of the equation, or we can engineer $k_v = k_0$ so that the scalar is 1.
    
    By these considerations, when we construct an appropriate global elliptic endoscopic parameter $\fM'$ which restricts to our given unramified elliptic endoscopic local parameter $\fM_v'$ and prove that the WFL holds on $V = \{v\}$ for $M \subset G$ in the sense of Definition \ref{def.global.wfl.V.places}, the result will be Theorem \ref{thm.wfl}.
\end{remark}

\begin{lemma}[{Splitting Lemma, \cite{arthurinvlocal}, \cite{arthurtransferwoi}, \cite[Proposition 16.4.1,17.4.1]{CLII}}]\label{lem.splitting.lemma}
    If $V = V_1 \bigsqcup V_2$ is a disjoint decomposition of $V \subset |X_0|$, then one has for any $a_M$ and $s$ as above
    \begin{equation}
        J_{M,V}^{G,s}(a_M) = \sum_{L_1, L_2 \supseteq M} d_M^G(L_1, L_2) J_{M,{V_1}}^{L_1, s}(a_M) J_{M,{V_2}}^{L_2, s}(a_M).\label{eqn.splitting.formula}
    \end{equation}
    Likewise, the stable orbital integrals satisfy a splitting formula
    \begin{equation}
        S_{M,V}^{G}(a_M) = \sum_{L_1, L_2 \supseteq M} e_M^G(L_1, L_2) S_{M,V_1}^{L_1}(a_M) S_{M,V_2}^{L_2}(a_M),\label{eqn.splitting.formula.stable}
    \end{equation}
    where $$e_M^G(L_1,L_2) := d_M^G(L_1, L_2) \left|\frac{Z_{\widehat{L_1}}^\Gamma \cap Z_{\widehat{L_2}}^\Gamma}{Z_{\widehat{G}}^\Gamma} \right|^{-1}$$
    are Arthur's stable coefficients. Writing then for any intermediate Levi
    $$
        E_{M,V}^{L}(a_{M'}) := |Z^\Gamma_{\widehat{M}'}/ Z_{\widehat{M}}^\Gamma| \sum_{s \in \cE_{\fM'}(L)} |Z_{\widehat{L'_s}}^\Gamma/ Z_{\widehat{L}}^\Gamma |^{-1} S_{M',V}^{L'_s}(a_{M'})
    $$
    for the endoscopic side of the WFL, these satisfy the original splitting formula
    \begin{equation}
        E_{M,V}^{G}(a_{M'}) = \sum_{L_1, L_2 \supseteq M} d_M^G(L_1, L_2) E_{M,{V_1}}^{L_1}(a_{M'}) E_{M,{V_2}}^{L_2}(a_{M'})\label{eqn.splitting.endoscopic}
    \end{equation}
    In particular, when $G = M$ we have the product formulas
    $$
        J_{M,V}^{M,s}(a_M) = J_{M,{V_1}}^{M, s}(a_M) J_{M,{V_2}}^{M, s}(a_M), \qquad E_{M,V}^{M}(a_{M'}) = E_{M,{V_1}}^{M}(a_{M'}) E_{M,{V_2}}^{M}(a_{M'})
    $$
    \begin{proof}
        That this holds for the weighted orbital integrals on $V$ is asserted in the proof of \cite[Proposition 8.1]{ArthurStableI}. To be precise, the original splitting formula is proven for $(G,M)$-families in \cite[Corollary 7.4]{arthurinvlocal}, then applied for weighted orbital integrals in \cite[Proposition 4.3]{arthurtransferwoi}. That the stable orbital integrals satisfy the corresponding splitting identity, and then that this packages to recover the endoscopic splitting formula  (\ref{eqn.splitting.endoscopic}) is the content of \cite[Theorem 6.1]{arthurtransferwoi}. 

        In the split case, the first claim is \cite[17.4.1]{CLII}, and their proof carries over verbatim (though interpreted in the relative spaces, of course). That is the only actual analysis required -- the citations of Arthur above, and in particular the proof of \cite[Theorem 6.1]{arthurtransferwoi}, demonstrate that the combinatorics which deduce the corresponding splitting formula for stable orbital integrals, and then recombine to obtain the original splitting formula for the endoscopic side, work in the relative setting without any splitness hypotheses. Carrying over the argument into our notation is left as an exercise to the reader.
    \end{proof}
\end{lemma}

In order to inductively remove places, one uses the following criterion.

\begin{lemma}[Peeling lemma]\label{lem.peeling}
    Fix all notations from the beginning of this section, and let $V = V_1 \bigsqcup V_2$ be a nontrivial decomposition of sets of places. Suppose
    \begin{enumerate}
        \item that the $\fM'$-WFL holds on $V$ for $M \subset G$,
        \item that the $\fM'$-WFL holds on $V_2$ for $M \subset G$,
        \item that
        $$
            J_{M, V_2}^{M, s_M}(a_M) \neq 0.
        $$
        Then the $\fM'$-WFL holds on $V_1$ for $M \subset G$.
    \end{enumerate}
    In other words, if the WFL holds on $V = V_1 \sqcup V_2$ for $M \subset G$ and $V_2$ is peelable, then the WFL holds on $V_1$ for $M \subset G$. 
    \begin{proof}
        We induct upwards on the intermediate Levi $M \subset L \subset G$. Subtract the two splitting formulae on $V$, which are equal by (1), and apply the $\fM'$-WFL on $V_2$, as well as the inductive equalities on $V_1$. Every term cancels except $(L_1, L_2) = (L,M)$, whose coefficient is 1. So
        $$
            J_{M, V_1}^{L, s_M}(a_M) J_{M, V_2}^{M, s_M}(a_M) = E_{M,{V_1}}^{L}(a_{M'}) J_{M, V_2}^{M, s_M}(a_M),
        $$
        and we conclude by nonvanishing on $V_2$, (3).
    \end{proof}
\end{lemma}

\begin{remark}\label{rmk.peelable.union}
    It is easy to check that if $V_1, V_2$ are disjoint and peelable, then $V_1 \sqcup V_2$ is peelable, but we won't use this, so we don't prove it. 
\end{remark}

\section{Peeling the regular places.}\label{subsect.peel.regular}

In this section, we collect the local computations which we will use to show various places are peelable. In light of the Peeling Lemma, this will allow us to inductively reduce the global WFL of Theorem \ref{thm.global.wfl} down to one place and deduce Theorem \ref{thm.wfl}. Of course, the claim that a place is peelable includes the fact that a version of the WFL holds at that place. Thus, our proof should still be considered as an application of Ng\^{o}'s ``Perverse Continuation Principle'' -- by packaging the WFL into families and understanding the support summands in the relative cohomology, we can by hand prove the identity for sufficiently simple characteristic polynomials and then transfer its validity to arbitrarily complicated characteristic polynomials. Thus, the results of this section and the next should be interpreted as essentially affirming the validity of the WFL at good places (except for the discrepancy discussed in Remark \ref{rmk.monodr.local.global.discrep.wfl}), from which its global validity spreads along the uniformity imposed by our cohomological theorems.

Recall that associated to an endoscopic group $H$ for $G$ we have the discriminant divisor $\fD_H$ of Definition \ref{def.discr.divisor} and the resultant divisor $\fR_H^G$ of Definition \ref{def.resultant.divisor}, both on $\fc_H$. In particular, these specialize to take $H = M$ a standard Levi, or $H = M'$ an endoscopic group associated to Levi $M$. In any case, these provide a measure of complexity for $$b \in \fc_{H}(\O_{0,v}) \cap \fc_{H}^{\Greg}(F_{0,v})$$ a generically $G$-regular parameter. We write
$$
    v(\fD_{H}(b))
$$
for the multiplicity of the closed point in the Cartier divisor $(b)^\star \fD_{H}$ on $\Spec \O_{0,v}$, and 
$$
    v(\fR_{H}^G(b)) 
$$
similarly. Of course, if $f$ is any local function cutting out the divisor, then the above multiplicity is computed by
$$
    v(f(b')),
$$
so the notation is only mildly abusive. 

The basic local computations are encoded in the following Lemma.

\begin{lemma}\label{lem.unit.disc.result.vanish}
    Let $b \in \fc_M(\O_{0,v}) \cap \fc_M^{\Greg}(F_{0,v})$, and write $Y := \epsilon_M(b)$.
    \begin{enumerate} 
    \item (Unit-discriminant normalization) Suppose
    $$
        v(\fD_M(b)) = 0.
    $$
    Then for every $s \in \widehat{T}^{W_{b,v}}$ we have
    $$
        J_{M,v}^{M,s}(Y) = 1.
    $$
    More precisely, every $[m] \in (J_b(F_v) \backslash M(F_v))^{\tau}$ contributing to the integral has a representative in $M(\O_{0,v})$ and in particular has invariant 0.
    \item (Unit-resultant vanishing) Let $M \subsetneq L$ an intermediate Levi, so $b$ is generically $L$-regular semisimple. If
    $$
        v(\fR_M^L(b)) = 0
    $$
    then for every $s \in \widehat{T}^{W_{b,v}}$ we have
    $$
        J_{M,v}^{L,s}(Y) = 0.
    $$
    More precisely, every $g \in L(F_v)$ contributing to the integral lies in $M(F_v) L(\O_v)$.
    \end{enumerate}
    If moreover $M \subseteq L$ is obtained by restricting global Levi subgroups to $v$ and $b$ is the restriction of a global parameter $a_M$, the same conclusions hold for the singleton partial adelic orbital integrals $J_{M, \{v\}}^{L,s}(a_M)$ for every $s \in \widehat{T}^{W_{a_M}}$, viewed locally by restriction, after the divisor renormalization of Remark \ref{rmk.local.renorm.D}.
    \begin{proof}
        For the first claim, let $[m] \in (J_b(F_v) \backslash M(F_v))^{\tau}$ be a coset contributing to the integral, and choose representative $m \in M(F_v)$. Put $Y_1 = \Ad(m^{-1})Y \in \fm(\O_v)$. Since the coset is $\tau$-fixed, there exists $j \in J_b(F_v)$ such that $\tau(m) = jm$. Since $j$ centralizes $Y$,
        $$
            \tau(Y_1) = Y_1 \qquad \implies \qquad Y_1 \in \fm(\O_{0,v}).
        $$
        So $Y$ and $Y_1 = \Ad(m^{-1})Y$ are rational elements which have the same (integral) characteristic polynomial. Reducing mod $z_v$, we have regular semisimple elements over finite residue field with the same characteristic polynomial. These are rationally conjugate by Lang's theorem, and the transporter is smooth because it is a torsor for the centralizer, a torus, so Henselian lifting implies conjugacy by $m_1 \in M(\O_{0,v})$. Then $m_1 m^{-1} \in Z_M(Y)(F_{v}) = J_Y(F_{v})$, so 
        $$
            m \in J_Y(F_{v}) M(\O_{0,v}).
        $$
        So every such coset $[m]$ has a rational representative $m_1$, and the only elements contributing range over a single $J_Y$-orbit. The local invariant map is then identically 0 on the domain by definition, so $\langle s, \inv_{b,v}(g)\rangle = 1$ identically for each $s$. The discriminant factor of course contributes 1 and the local weight factor $v_M^{M,\loc}$ is identically 1, so we compute simply that with our normalization of measure
        $$
            J_{M,v}^{M,s}(b) = 1
        $$
        for every $s$. The singleton claim follows similarly, since $v_M^M = 1$ likewise.

        For the second claim, we will show that the resultant condition implies 
        $J_{M}^{L,s}(Y) = 0$ unless $M = L$. Let $g \in L(F_{0,v})$ be an element such that 
        $$
            \Ad(g^{-1}) Y \in \fl(\O_{0,v}). 
        $$
        Choose a parabolic $P = M \rtimes N$ with Levi $M$ and consider an Iwasawa decomposition $g = mnk$. The decomposition 
        $$
            \Ad(k) \Ad(g^{-1}) Y = (\Ad(n^{-1}) Z - Z) + Z
        $$
        where $Z = \Ad(m^{-1})Y \in \fm(F_{0,v})$ and $(\Ad(n^{-1}) Z - Z) \in \fn(F_{0,v})$ implies that both components are integral. 

        The resultant hypothesis implies, after the formula of Proposition \ref{prop.unip.formula.for.resultant}
        $$
            \det(\ad(Y)|_{\fn}) = \pm R_M^L(\chi_M(Y)),
        $$
        that $\ad(Y)$ and $\ad(Z)$, maps $\fn(\O_{0,v}) \to \fn(\O_{0,v})$, are integral isomorphisms. 

        We claim that $n \in N(\O_v)$. Over the maximal unramified extension, pick a cocharacter $\lambda$ defining $P$ as its attracting locus and $M$ as its centralizer, and filter $N$ by positive weights for this cocharacter, $N_1 \supset N_2 \supset \cdots$. The map
        $$
            \phi : n \mapsto \Ad(n^{-1})Z - Z
        $$
        preserves the filtration, and on the associated graded it becomes a homomorphism
        $$
            \phi(n_1 n_2) = \Ad(n_2^{-1}) \phi(n_1) + \phi(n_2) = \phi(n_1) + \phi(n_2),
        $$
        indeed a linear $\O_v$-isomorphism $N_i/N_{i+1} \to \fn_i/\fn_{i+1}$ which is equal to its differential $\ad(Z)$ under the identification $N_i/N_{i+1} \cong \fn_i/\fn_{i+1}$. We can then prove by induction that $n \in N(\O_v)$: we know that $\phi(n)$ is integral, and on the first subquotient it induces an $\O_v$-isomorphism, so the class of $n$ is integral. The multiplicative difference, then, lies in $N_2(F_v)$ and has integral image under $\phi$. Inductively, we can show all the differences are integral.
        
        The original element $g$ therefore lies in $M(F_{v}) L(\O_{v})$. But this implies that all $H_{P,v}^{\loc}(g)$ are equal, and Arthur's weight $v_M^{L,\loc}(g)$ vanishes unless $M = L$, so we have $J^{L,s}_{M,v}(Y) = 0$ for any $L \supsetneq M$. The last claim about the single integrals follows similarly: the global $H_P$ are also equal and the global $v_M^L$ vanishes.
    \end{proof}
\end{lemma}

These local computations quickly imply the peelability of all regular places and a reduction of our global WFL down to finitely many places. Fix global $\fM'$, $M \subset G$ and $a_{M'}$ satisfying the assumptions from the top of \S \ref{subsect.splitting.lemmas}. 

\begin{corollary}\label{cor.reg.v.peelable}
    Let $V \subseteq |X_0|$ be any set (finite or infinite) of regular places, in the sense that $a_M$ avoids $\fD_G$ on $V$. Then $V$ is peelable.

    In particular, the $\fM'$-WFL holds on the finite set $B \subset |X_0|$ for $M \subset G$, where
    $$
        B = \on{supp}(a_M^\star \fD_M + a_M^\star\fR_M^G).
    $$
    \begin{proof}
        Up to the renormalization on $\on{supp}(D)$ described in Remark \ref{rmk.local.renorm.D}, the restriction of $a_{M}$ (resp. $a_{M'}$) to a formal disk around each $v \in V$ defines a valid local parameter $b_v$ (resp. $b'_v$). These are valid local parameters in the sense of the first section. Pick an intermediate Levi $L$. By the decomposition $\fD_G|_{\fc_M} = \fD_M + 2 \fR_M^G$ and likewise for $M'$, the local parameters satisfy
        $$
            v(\fD_{M'}(b'_v)) = v(\fD_M(b_v)) = 0, \qquad v(\fR_{M'}^L(b'_v)) = v(\fR_M^L(b_v)) = 0.
        $$
        We apply the ``more precisely'' statements of Lemma \ref{lem.unit.disc.result.vanish} at each place in $V$. The weight functions, endoscopic parameter sets, and stabilization coefficients remain defined relative to the global $\Gamma$, cf.\ Remark \ref{rmk.monodr.local.global.discrep.wfl}. 
        
        Adelically on $V$, we deduce that the partial \emph{adelic} $s$-orbital integral satisfies
        $$
            J^{L, s}_{M,V}(a_M) = \delta_{L = M}.
        $$
        For the case $L = M$, the support of the $V$-orbital integral modulo $J_{a_M}(\A_V)$ is the trivial class represented by 1 by the result at each place, so 
        $$
            J^{M, s}_{M,V}(a_M) = 1.
        $$
        For the case $L \supsetneq M$, we have shown that every $g_v$ which can occur in the support of the integral can be written as $g_v = m_v k_v$ for $m_v \in M(F_v), k_v \in L(\O_v)$. In particular,
        $$
            H_{P,v}(g_v) = H_{M,v}(m_v)
        $$
        is independent of parabolic $P$ for $M$, so each adelic point $(g_v)_v$ in the support of the integral satisfies that
        $$
            H_{P,V}(g) = \sum_v H_{P,v}(g_v)
        $$
        is independent of $P$. The adelic convex hull defining $v_M^L(g)$ in $\a_M^L$ is therefore a single point, with volume zero unless $L = M$. We have shown
        $$
            J^{L, s}_{M,V}(a_M) = \delta_{L = M}.
        $$
        
        Now we collate the computation $$J_{M,V}^{L,s}(a_M) =\delta_{L = M}$$ into the form of the WFL. We have, by the same avoidance of discriminants and resultants for $M'$ that the same Lemma applies, and repeating the argument above we have
        $$
            J_{M',V}^{L', 1}(a_{M'}) = \delta_{L' = M'}
        $$
        for induced endoscopic $L'$. 

        The $\fM'$-WFL on $V$ is then the following identity for each $L$:
        $$
            J_{M,V}^{L,s}(a_M) = |Z^\Gamma_{\widehat{M}'}/ Z_{\widehat{M}}^\Gamma| \sum_{s \in \cE_{\fM'}(L)} |Z_{\widehat{L'_s}}^\Gamma/ Z_{\widehat{L}}^\Gamma |^{-1} S_{M',V}^{L'_s}(a_{M'}),
        $$
        where 
        $$
            S_{M',V}^H(a_{M'}) = J_{M',V}^{H,1}(a_{M'}) - \sum_{\substack{s \in \cE_{M'}(H) \\ s \neq 1}} | Z_{\widehat{H}'_s}^\Gamma/Z_{\widehat{H}}^\Gamma|^{-1} S_{M', V}^{H_s'}(a_{M'})
        $$
        Consider this last recursion. If $M' \subsetneq H$, then any induced elliptic group $H'_s$ is not equal to $M'$ (the split central torus changed!) so we will never get a nonzero term in the recursion defining $S_{M',V}^H(a_{M'})$, and it is identically 0. On the other hand, $S_{M',V}^{M'}(a_{M'}) = 1$. Considering now the WFL, both the left and right side vanish for $M \neq L$ and the WFL is $0 = 0$. When $M = L$, the identity reduces to
        $$
            J_{M,V}^{M,s}(a_M) = J_{M',V}^{M',1}(a_{M'})
        $$
        and both are 1. So we have checked the $\fM'$-WFL holds on the set of regular points. And of course, we have shown the nonvanishing condition for peelability, so $V$ is peelable.

        The last claim is now immediate from the Peeling Lemma \ref{lem.peeling} and the global WFL of Theorem \ref{thm.global.wfl}: the $\fM'$-WFL on $|X_0|$ is exactly Theorem \ref{thm.global.wfl} applied to $\fM', M$, and each intermediate Levi $L$, each of which inherits the required hypotheses for $G$. Then we can peel off the set of regular places by the Peeling Lemma. 
    \end{proof}
\end{corollary}

\section{Transverse places.}\label{subsect.transverse.places}

If we could arrange a globalization of our local parameter at $v$ which was regular at all other points, then we would now be done, but this is too much to hope for. Instead, as in \cite[\S 8.5]{ngofl}, consider the divisor 
$$
    \Delta = (\fD_G|_{\fc_{M',D}})_{\red} = (\fD_{M'} + \fR_{M'}^G)_{\red} = (\fD_{M'} + \fR_{M'}^M + (\chi_M^{M'})^\star\fR_{M}^G)_{\red}.
$$
We will show below that a Bertini theorem allows us to choose a globalization of our local parameter which intersects this divisor \emph{transversally} at smooth points of the divisor. If $b'$ is a local parameter given by restriction of $a_{M'}$ to a formal disk around the place $v$ and $v((\fD_G|_{\fc_{M',D}})_{\red}(b')) = 1$, then exactly one of
$$
    v(\fD_{M'}(b')), \qquad v(\fR_{M'}^M(b')), \qquad v(\fR_M^G(\chi^{M'}_M(b')))
$$
is equal to 1 and the others vanish. We need to show that such places are peelable by hand. The following three lemmas do this. We assume $v \not\in \on{Supp}(D)$, so that $J_{M, \{v\}}^{L,s}$ and $S_{M', \{v\}}^{H}$ can be written with local arguments $Y_b, Y_{b'}$ rather than $a_M, a_{M'}$.

\begin{remark}\label{rmk.transv.lemma.loc.global.discrep} As discussed in Remarks \ref{rmk.partial.local.comparison.J} and \ref{rmk.monodr.local.global.discrep.wfl}, $\Gamma$ is the global splitting Galois group of Definition \ref{def.global.wfl.V.places}, not the local decomposition group, so the weight factors and $\Gamma$-fixed scalars are not necessarily the usual local ones, and even when $D_v = \Gamma$ the orbital integrals differ by a scalar. Peelability concerns the partial orbital integrals of Definition \ref{def.s.integral.V} at $V = \{v\}$, \emph{not} the genuinely local orbital integrals. 

In the first two lemmas, we will pass to the genuine local objects in order to reduce to the usual Langlands--Shelstad fundamental lemma, applied to the restriction of that global datum along $D_v \subseteq \Gamma$. In this case, both the global and local weight functions will be identically 1 or 0 on the support of the integrand, and in the former case the required $\Gamma$-fixed scalars cancel directly, so we reduce to an equality of genuinely local unweighted orbital integrals.\end{remark}   

\begin{lemma}[{\cite[Lemme 8.5.7]{ngofl}, cf.\ \cite[Lemme 17.7.3]{CLII}}]\label{lem.transv.disc}
    Let $v \in |X_0|$ be a place disjoint from $\on{supp}(D)$, and write $b$ (resp.\ $b'$) for the restriction of the global parameter $a_M$ (resp.\ $a_{M'}$) to the completed disk around $v$. Suppose
    \begin{equation}
        v(\fD_{M'}(b')) = 1, \qquad v(\fR_{M'}^M(b')) = 0, \qquad v(\fR_M^G(b)) = 0.\label{eqn.tranv.D}
    \end{equation}
    Then $v$ is peelable.
    \begin{proof}
        The hypothesis of the second claim in Lemma \ref{lem.unit.disc.result.vanish} applies both for $M$ and $M'$, so we know for $L \supsetneq M$
        $$
            J_{M,\{v\}}^{L,s_M}(Y_b) = 0
        $$
        and likewise for $M'$-orbital integrals. On the endoscopic side, every group is $L$-elliptic, with strictly smaller global split-center dimension than $M$. In particular, $M'$ itself will never appear. The recursive definition of $S^{L'_s}_{M', \{v\}}$ will therefore never contribute any nonzero terms, so the endoscopic side also vanishes. 
        
        On the other hand, when $L = M$ the WFL identity (\ref{eqn.V.wfl}) reduces to the claim
        \begin{equation}
            J_{M,\{v\}}^{M,s_M}(Y_b) = J_{M', \{v\}}^{M', 1}(Y_{b'}).\label{eqn.M=L.fl}
        \end{equation}
        We now apply the (unweighted) Langlands--Shelstad fundamental lemma in the simple case of transverse discriminant, case (2) in the proof of \cite[Lemme 8.5.7]{ngofl}, to the local orbital integrals for the restriction of $(s_M, \rho_{s_M})$ to $D_v \subseteq \Gamma$. Indeed, this case is proven by hand in \textit{loc.\ cit.}: the affine Springer fibers in question are zero-dimensional and therefore are equal to their regular locus, which transfers from $M'$ to $M$. That local identity is exactly (\ref{eqn.M=L.fl}): the global fixed-center coefficients cancel and the local and global weights are identically 1. In fact, the restriction along $D_v \subseteq \Gamma$ discussed in Remarks \ref{rmk.partial.local.comparison.J} and \ref{rmk.transv.lemma.loc.global.discrep} preserves both orbital integrals when $L = M$: the integration measures are identified, and the local and global Kottwitz maps correspond, identifying the local and global pairings with $s_M$. The partial integrals on $\{v\}$ therefore coincide with the local integrals at $v$. Note that our normalizations absorb the standard power of $q$. The integral is nonzero: the stable class is a sum (with coefficients 1) of rational classes, and an open compact neighborhood of the identity contributes nontrivially to the distinguished rational class:
        $$
            J_{M',v}^{M', 1}(Y_{b'}) \ge |D^{M'}(Y_{b'})|_v^{1/2} \int_{J_{b'}(F_{0,v}) \backslash M'(F_{0,v})} 1_{\fm'(\O_{0,v})}(\Ad(g^{-1}) Y_{b'}) \; dg > 0.
        $$
        So $v$ is peelable.
    \end{proof}
\end{lemma}

\begin{lemma}[{cf.\ \cite[\S 8.3]{ngofl}, \cite[Lemme 17.7.2]{CLII}}]\label{lem.tranv.result.1}
    Let $v \in |X_0|$ be a place disjoint from $\on{supp}(D)$, and write $b$ (resp.\ $b'$) for the restriction of the global parameter $a_M$ (resp.\ $a_{M'}$) to the completed disk around $v$. Suppose
    \begin{equation}
        v(\fD_{M'}(b')) = 0, \qquad v(\fR_{M'}^M(b')) = 1, \qquad v(\fR_M^G(b)) = 0.\label{eqn.tranv.R}
    \end{equation}
    Then $v$ is peelable.
    \begin{proof}
        If $v(\fR_{M'}^M(b')) = 1$, then 
        $$
            v(\fD_M(b)) = 2.
        $$
        As above, the case $L = M$ restricts along $D_v \subseteq \Gamma$ to a special case of the Langlands--Shelstad Fundamental Lemma, case (3) in the proof of \cite[Lemme 8.5.7]{ngofl}, checked by hand using the one-dimensional affine Springer fiber computation of \cite[\S 8.3]{ngofl}, and the right side is 1 by unit-discriminant for $M'$. When $L \supsetneq M$, then we have
        $$
            J_{M,\{v\}}^{L, s_M}(Y_b) = 0
        $$
        by the unit-resultant part of Lemma \ref{lem.unit.disc.result.vanish}. As before, global $L$-ellipticity in the stable recursion ensures every $L'_s \supsetneq M'$ for every term which appears. We need to check that we also have unit resultant on the endoscopic side, $v(\fR_{M'}^{L_s'}(b')) = 0$. This intersection multiplicity can be checked after passage to $\O_v$, at which point we can use the root-theoretic formula for the resultant as discussed after Proposition \ref{prop.resultant.transitive}: we can choose $Y' \in \ft(\O_v)$ a \emph{diagonal} lift of $b'$ by \'etale lifting, and then
        $$
            v(\fR_{M'}^{L_s'}(b')) = \sum_{\{\pm \alpha\} \subset \Phi(L_s', T) \setminus \Phi(M', T)} v(d\alpha(Y')).
        $$
        Since $v(\fR_M^L(b)) = 0$, the corresponding product over $\Phi(L,T) \backslash \Phi(M, T)$ is equal to 0, which implies that $$
            v(d \beta(Y')) = 0 \quad \forall \beta \in \Phi(L,T) \setminus \Phi(M,T).
        $$
        So it suffices to check the containment of geometric roots
        $$
            \Phi(L_s', T) \setminus \Phi(M', T) \subseteq \Phi(L,T) \setminus \Phi(M,T).
        $$
        This can be checked in terms of coroots on the dual side. But for sub-root systems defined from connected centralizers of torus elements, a coroot being contained in the sub-coroots is equivalent to the corresponding root being contained in the sub-roots. So it suffices to check
        $$
            \Phi(\widehat{L_s'}, \widehat{T}) \setminus \Phi(\widehat{M'}, \widehat{T}) \subseteq \Phi(\widehat{L},\widehat{T}) \setminus \Phi(\widehat{M},\widehat{T}),
        $$
        which is exactly the root computation of Proposition \ref{prop.ind.G'.cont.M'.Levi}. We have shown that $v(\fR_{M'}^{L_s'}(b')) = 0$, so the unit-resultant part of Lemma \ref{lem.unit.disc.result.vanish} also applies to all the integrals appearing on the endoscopic side, and the WFL $0 = 0$ holds for $L \supsetneq M$; we have shown $v$ is peelable in this subcase.
\end{proof}
\end{lemma}

The final case is the only genuinely weighted special case, and is also the place where the most new work is required in order to handle the nonsplit case. 
\begin{lemma}[{cf. \cite[Lemme 17.6.2, Lemme 17.8.1]{CLII}}]\label{lem.tranv.result.2}
    Let $v \in |X_0|$ be a place disjoint from $\on{supp}(D)$, and write $b$ (resp.\ $b'$) for the restriction of the global parameter $a_M$ (resp.\ $a_{M'}$) to the formal disk around $v$. Suppose
    \begin{equation}
        v(\fD_{M'}(b')) = 0, \qquad v(\fR_{M'}^M(b')) = 0, \qquad v(\fR_M^G(b)) = 1.\label{eqn.tranv.R.wt}
    \end{equation}
    Then $v$ is peelable.
    \begin{proof}
        First we handle all the easy parts. In the case $L = M$, we have done all the work:
        $$
            v(\fD_M(b)) = v(\fD_{M'}(b')) + 2 v(\fR_{M'}^M(b')) = 0
        $$
        so both $b$ and $b'$ are regular and the unit-discriminant lemma \ref{lem.unit.disc.result.vanish}.(1) applied to the restriction along $D_v \subseteq \Gamma$ gives\footnote{The weight function is identically 1, so these are genuinely the local $s$-orbital integrals.}
        $$
            J_{M, \{v\}}^{M, s_M}(Y_b) = 1 = J_{M', \{v\}}^{M', 1}(Y_{b'}).
        $$ 
        This is the WFL identity, since the global fixed-center scalar for $M'$ cancels the coefficient of the unique endoscopic term, and is nonzero as we require for peelability. It remains only to check the WFL identity (\ref{eqn.V.wfl}) holds for $L \supsetneq M$. If $v(\fR_M^L(b)) = 0$, then the unit-resultant Lemma \ref{lem.unit.disc.result.vanish}.(2) implies the WFL reduces to $0 = 0$ as in the previous cases. So assume $v(\fR_M^L(b)) = 1$. 
        
        The proof in the remaining case is rather involved, and will occupy us until the end of the section. For the reader's convenience, we outline the argument. It will broadly have two steps, along the same strategy as in \cite[{Lemma 17.6.2, 17.8.1}]{CLII}. In the first step, we identify a pair of bad roots $\{\pm \alpha\}$ contributing to $v(\fR_M^L)$, build a geometric Levi $R_\alpha^{\on{geom}}$ containing only the torus and that pair of roots, and then relate our orbital integrals to integrals for that group. Following Chaudouard--Laumon, we will show that $R_\alpha^{\on{geom}}$ is stable for a twisted Frobenius endomorphism, and thus descend our orbital integral almost entirely to one for an unramified group in relative semisimple rank 1. The form will be simple enough that we can compute the $\{v\}$-orbital integrals explicitly, uniformly both for $L \supsetneq M$ and for the endoscopic pairs $H \supsetneq M'$ which appear in the stable recursion. 
        
        The value will always be either a simple nonzero constant $C$ (if $H$ still contains the bad pair of roots $\pm\alpha$ and $\dim \a_M^L = 1$) or else $0$, so in the second step we analyze the endoscopic sums in the stabilization formula (\ref{eqn.stable.def.recurse}) and in the WFL (\ref{eqn.V.wfl}) to count recursively the total multiple of $C$ contributed by the stable orbital integrals. Both sides contribute $1 \cdot C$, so the WFL holds. 
        
        We remark briefly on the comparison between our proof and that of \textit{loc.\ cit.}. In the former, they also reduce the integral to $R_\alpha^{\on{geom}}$, which they do not explicitly compute. They then show the collection of orbital integral values appearing descend to a family for the pair $T \subseteq R_\alpha^{\on{geom}}$ in semisimple rank 1, in some sense orthogonal to $M$. Stabilization for the original family is then a consequence of the (vacuous) stabilization of the descended family. Although Arthur's descent formalism is quite general and allows nonsplit groups, we have not found a way to make the same formal argument fit within it: instead, our reductions simplify the stabilization combinatorics enough to compute both sides by hand.

        Before properly beginning, we reduce our $s$-orbital integral to the distinguished rational component in the stable class as in claim (1) of Lemma \ref{lem.unit.disc.result.vanish}: given any coset contributing to the integral $[g] \in (J_b(F_v) \backslash L(F_v))^{\tau}$ with Iwasawa decomposition of a representative
        $$
            g = mnk \in M(F_v) N(F_v) L(\O_{v})
        $$
        relative to parabolic $P = MN$, the nonvanishing of $\fD_M$ on the special fiber implies by Henselian lifting that $m \in J_b(F_v) M(\O_v)$. The cocycle defining $\inv_{s_M}(g)$ can then be represented by an element of $J_b(\O_v)$. But this is a smooth unramified torus, so Lang's theorem kills the cocycle and $\inv(g) = 0$. In particular, the $s_M$-orbital integral reduces from $(J_b(F_v) \backslash L(F_v))^{\tau}$ to its distinguished fiber, $J_b(F_{0,v}) \backslash L(F_{0,v})$.

        \subsection{Step 1: Computation of the weighted orbital integrals}         After restricting to the completed algebraically closed disk $\O_v$ and choosing compatible split presentations for $M$ and $M'$, choose by Hensel-lifting from the (algebraically closed) special fiber an element $Y \in \ft(\O_v)$  with $\chi_{M'}(Y) = b'$ (so also $\chi_M(Y) = b$). By regular semisimplicity, Henselian-lifting the transporter from $Y$ to $Y_b := \epsilon_M(b) \in \fm(\O_{0,v}) \subset \fm(\O_v)$ gives $m \in M(\O_v)$ with $Y = \Ad(m) Y_b$. Let $\tau$ denote the local Frobenius action at $v$ on the geometric based root datum. The original $Y_b$ was fixed by Frobenius $\tau$, so $Y$ is fixed by the twisted Frobenius
        $$
            \theta = \on{Int}(\wtilde{w}) \circ \tau,
        $$
        where $\wtilde{w} = m \tau(m)^{-1} \in N_M(T)(\O_v)$ with image $w \in W_M(\ol{v})$. This $\theta$ preserves the integral model $L_{\O_v}$ and $Y$, hence the subgroup $T_{\O_v}$. 
        
        We now translate our orbital integral from $Y_b$ to $Y$ along conjugation by $m$. The discriminant factor clearly remains unchanged. Since $m \in M(\O_v)$, conjugation by it preserves $\fl(\O_v)$ and sends the maximal compact in $J_b$ to $T(\O_v)^{\theta}$. Quotient measures normalized to give these volume 1 transfer without a scalar. Moreover, the Arthur weight is unchanged: it is right-invariant by $L(\O_v)$ and left-invariant by $M(F_v)$, so is conjugation-invariant by $m \in M(\O_v)$:
        $$
            v_M^L(m g m^{-1}) = v_M^L(g).
        $$
        In summary,
        \begin{align*}
            J_{M,\{v\}}^{L, s_M}(Y_b) &= |D^L(Y_b)|_v^{1/2} \int_{J_b(F_{0,v}) \backslash L(F_{0,v})} 1_{\fl(\O_v)}(\Ad(g^{-1}) Y_b) v_M^L(g) \; dg \\
            &= |D^L(Y)|_v^{1/2} \int_{T(F_v)^{\theta} \backslash L(F_{v})^{\theta}} 1_{\fl(\O_v)}(\Ad(g^{-1}) Y) v_M^L(g) \; dg.
        \end{align*}

        The geometric-rank-1 construction now adapts to the nonsplit setting as follows. Recall our fixed maximal torus $T$. By the differential formula for the resultant, there exists a \emph{unique} pair of geometric roots $\{\pm \alpha\} \subset \Phi(G_{F_v},T_{F_v}) \setminus \Phi(M_{F_v},T_{F_v})$ such that 
        $$
            v(d\alpha(Y)) = 1
        $$
        and the others are units. We claim $\theta(\alpha) = \alpha$. On the one hand, we must have 
        $$
            v(d(w\tau \alpha)(Y)) = v(d\alpha(Y)) = 1
        $$
        so $w\tau \alpha = \pm \alpha$. Choose a global parabolic $P = MN$ for $M$, and define $+\alpha$ to be the geometric root contained in $N_{F_v}$. Then $N_{F_v}$ is stable under $\tau$ and $w \in W_M$, so $w\tau \alpha \in \Phi(N_{F_v})$. So $\theta(\alpha) = \alpha$. 
        
        Consider therefore the \emph{geometric} rank-one Levi cut out by the bad roots 
        $$
            R_\alpha^{\on{geom}} := \langle T_{\O_v}, U_\alpha, U_{-\alpha} \rangle \subset L_{\O_v}.
        $$
        This is a $\theta$-stable (but not necessarily $\tau$-stable!) semistandard Levi over $\O_v$ of absolute semisimple rank 1, containing $\theta$-stable Borel $T U_\alpha$. Interpreting $\theta$ as a different choice of Frobenius endomorphism on $L_{\O_v}$, the inclusion $R_\alpha^{\on{geom}} \subset L_{\O_v}$, as well as $T, U_\alpha$, descends from $\O_v$ to inclusions of abstract reductive groups over $\O_{0,v}$ (possibly incompatibly with the original Frobenius structure). Compare \cite[\S 8.3.3]{ngofl}. Likewise, since $\theta$ and $\tau$ differ by an inner twist by $M(\O_v)$, every global $X$-parabolic $P \subset L$ with Levi $M$ restricts to an $\O_v$-parabolic which is $\theta$-stable, so likewise descends to an $\O_{0,v}$-parabolic with respect to this Frobenius. We henceforth write $\O_{0,v}$-points and $F_{0,v}$-points for all these groups relative to $\theta$, i.e., to mean $\theta$-fixed points. Pick with respect to this structure an $\O_{0,v}$-integral coordinate $x_\alpha : \G_a \xrightarrow{\sim} U_\alpha$, with generating differential $X_\alpha := (dx_\alpha)_0(1) \in \Lie(U_\alpha)(\O_{0,v})$, and corresponding $x_{-\alpha}, X_{-\alpha}$.
        
        We argue now that the $\theta$-orbital integral above is supported on $R_\alpha^{\on{geom}}$. Choose geometric parabolic $Q \subset L_{F_v}$ with Levi $R_\alpha^{\on{geom}}$ and unipotent radical $N$. Every root value of $Y$ on $\fn$ is a unit. By Iwasawa and the same filtration argument as in the proof of Lemma \ref{lem.unit.disc.result.vanish}.(2), one can show that all $g \in L(F_v)$ with $\Ad(g^{-1}) Y \in \fl(\O_v)$ are of the form 
        $$
            g = lk, \qquad l \in R_{\alpha}^{\on{geom}}(F_v), \quad k \in L(\O_v).
        $$
        If moreover $g$ is $\theta$-fixed (as is necessary to contribute to the integral) then
        $$
            \theta(lk) = lk \qquad \implies \qquad l^{-1} \theta(l) \in R_{\alpha}^{\on{geom}}(\O_v).
        $$
        Lang's theorem applies to $R_{\alpha}^{\on{geom}}(\O_v)$ with Frobenius $\theta$, so we can arrange $l = \theta(l)$ (and $k = \theta(k)$). 
        The map of double cosets
        $$
            T(F_v)^{\theta} \backslash R_{\alpha}^{\on{geom}}(F_v)^{\theta}/R_{\alpha}^{\on{geom}}(\O_v)^{\theta} \to T(F_v)^{\theta} \backslash L(F_v)^{\theta} / L(\O_v)^\theta
        $$
        is therefore a bijection on those satisfying the integrality condition.
        With respect to the quotient measures which equip $L(\O_v)^\theta$ and $R_{\alpha}^{\on{geom}}(\O_v)^\theta$ with volume 1, the containment of $\O_v$-groups implies for any $r \in R_\alpha^{\on{geom}}(F_v)^{\theta}$ an equality of stabilizers
        $$
            T(F_v)^{\theta} \cap r R_\alpha^{\on{geom}}(\O_v)^\theta r^{-1} = T(F_v)^{\theta} \cap r L(\O_v)^\theta r^{-1}
        $$
        and so the cosets have the same mass whether considered in $R_{\alpha}^{\on{geom}}$ or in $L$. Our orbital integral therefore reduces to
        $$
            J_{M,\{v\}}^{L, s_M}(Y_b) = |D^L(Y)|_v^{1/2} \int_{T(F_v)^{\theta} \backslash R_{\alpha}^{\on{geom}}(F_v)^{\theta}} 1_{\Lie(R_{\alpha}^{\on{geom}})(\O_v)}(\Ad(l^{-1}) Y) v_M^L(l) \; dl.
        $$
        Since $D_M(b)$ is a unit and the bad pair is unique, all roots of $L$ other than $\{\pm \alpha\}$ contribute units to the discriminant, so
        $$
            |D^L(Y)|_v = |D^{R_\alpha^{\on{geom}}}(Y)|_v = q_v^{-2}, \quad q_v := \# k_v.
        $$
        In terms of the $\theta$-Frobenius structure, the orbital integral can therefore be written (note $l$ and $Y$ are $\theta$-stable, so $\Ad(l^{-1})Y$ is also)
        $$
            J_{M,\{v\}}^{L, s_M}(Y_b) = q_v^{-1} \int_{T(F_{0,v}) \backslash R_{\alpha}^{\on{geom}}(F_{0,v})} 1_{\Lie(R_{\alpha}^{\on{geom}})(\O_{0,v})}(\Ad(l^{-1}) Y) v_M^L(l) \; dl.
        $$
        
        The formula above is almost a standard weighted orbital integral over $\O_{0,v}$ for this group $R_{\alpha}^{\on{geom}}$ of semisimple split rank 1, except that the weight function $v_{M}^L(l)$ is both global and not necessarily equal in a simple way to the weight function for $R_{\alpha}^{\on{geom}}$. Luckily, on the domain of $l$ which actually compute to the integral, we can compute this weight function, and thus the weighted orbital integral, directly. The standard Iwasawa-and-Fubini argument relative to the Borel $R_\alpha^{\on{geom}}(F_{0,v}) = T(F_{0,v}) U_\alpha(F_{0,v}) R_\alpha^{\on{geom}}(\O_{0,v})$ and our integral coordinate $x_\alpha : \G_a \cong U_\alpha$ give
        $$
            J_{M,\{v\}}^{L, s_M}(Y_b) = q_v^{-1} \int_{F_{0,v}} 1_{\Lie(R_{\alpha}^{\on{geom}})(\O_{0,v})}(\Ad(x_\alpha(u)^{-1}) Y) v_M^L(x_\alpha(u)) \; du
        $$
        relative to a measure on $\G_a$ which normalizes $\G_a(\O_{0,v}) = \O_{0,v}$. Since $Y$ acts on $\Lie(U_\alpha)$ by $d\alpha$, we calculate
        $$
            \Ad(x_\alpha(u)^{-1}) Y = Y + u d\alpha(Y) X_\alpha.
        $$
        The bad geometric root satisfies $v(d\alpha(Y)) = 1$, so this is integral if and only if $v(u) \ge -1$, i.e.,
        $$
            J_{M,\{v\}}^{L, s_M}(Y_b) = q_v^{-1} \int_{\pi^{-1} \O_{0,v}} v_M^L(x_\alpha(u)) \; du.
        $$
        It remains to calculate $H_{P,v}(x_\alpha(u))$ for each global $X$-parabolic $P$, and thus the global weight function, on $u$. Definition \ref{def.adelic.hitchin.deg} gives the definition for a geometric point, so for a closed point $v$ the weight used in the weighted orbital integral is
        $$
            H_{P,v} = \sum_{\ol{v}|v} H_{P, \ol{v}} = \deg(v) H_{P, \ol{v}}.
        $$
        Write $h_\alpha \in \a_M^L$ for the image under the natural projection of $\alpha^\vee \in \ga_T^L$.\footnote{A globally rational character defining $M$ as a Levi of $L$ pairs nontrivially with $\alpha^\vee$, so $h_\alpha \neq 0$.} For $u \in \O_{0,v}$, the weight $H_{P,v}(x_\alpha(u))$ is zero for all parabolics, and the integrand is identically 0. 
        
        Assume otherwise $v(u) = -1$. then the global parabolics $P \subset L$ with Levi $M$ each restrict to a $\theta$-stable parabolic. There are two possibilities.
        \begin{itemize}
        \item If $\alpha \in \Phi(R_u(P)_{\ol{v}},T_{\ol{v}})$, then $x_\alpha(u) \in R_u(P)$ is its own Iwasawa decomposition relative to $P$ with Levi image 1. Then $H_{P, v}(x_\alpha(u)) = 0$.
        \item If $-\alpha \in \Phi(R_u(P)_{\ol{v}}, T_{\ol{v}})$, then we need to compute the Iwasawa decomposition of $x_\alpha(u)$ relative to the parabolic $P$. This can be done using the open cell commutator relationship between $x_\alpha$ and $x_{-\alpha}$. For example, \cite[Theorem 4.2.6.(2)]{conradrgs} or \cite[Exp.\ XX, Th\'eor\`me 2.1]{SGAIII} applied to $x_\alpha(u)$ and $x_{-\alpha}(1)$ give the formula
        $$
            x_\alpha(u) =x_{-\alpha}\left(\frac{1}{1 + u}\right) \alpha^\vee(1 + u)  \cdot x_\alpha\left( \frac{u}{1 + u} \right) x_{-\alpha}(-1).
        $$
        The left half of the expression lies in $P(F_{0,v})$ and the right half lies in $R_\alpha^{\on{geom}}(\O_{0,v})$, so this is an Iwasawa decomposition for any such $P$, with Levi component $\alpha^\vee(1+u).$ For any globally rational character $\mu \in X^\star(M)$, we compute
        $$
            \langle \mu, H_{P,v}(x_\alpha(u)) \rangle = - \deg(v) v(\mu(\alpha^\vee(1 + u))) = -\deg(v)\langle \mu, \alpha^\vee \rangle v(u) = \deg(v) \langle \mu, \alpha^\vee \rangle,
        $$
        so in $\a_M^L$ we have $H_{P,v}(x_\alpha(u)) = \deg(v) h_\alpha$. 
        \end{itemize} 
        For any parabolic, its opposite lies in the other case, so both cases occur. The weight function extends the local parameter to an adelic parameter by adjoining 1 at all other places, so $H_P = H_{P,v}$. The convex hull characterization of Corollary \ref{cor.cvx.hull} and Proposition \ref{prop.deg.p.red.BL} in ambient group $L$ therefore give
        $$
            C_m = [0, -\deg(v) h_\alpha] + \a_T^M + \a_L,
        $$
        where $[0, -\deg(v) h_\alpha]$ is the straight line segment in $\a_M^L$. In particular, it has nonzero volume only if the space is 1-dimensional:
        $$
            v_M^L(x_\alpha(u)) = \begin{cases}
                \deg(v) \|h_\alpha\|, \quad & v(u) = -1, \;\dim \a_M^L = 1,\\
                0, & \text{otherwise.}
            \end{cases}
        $$
        In the orbital integral, then, the weight function is constant on the cosets of $\pi^{-1} \O_{0,v} / \O_{0,v}$, vanishing on the trivial coset and nonvanishing on the others if $\dim \a_M^L = 1$. We conclude our fundamental calculation
        $$
            J_{M,\{v\}}^{L,s_M}(Y_b) = \begin{cases} C, \quad & \dim \a_M^L = 1,\\
            0, & \text{otherwise,}
            \end{cases}
        $$
        where $C := (1 - q_v^{-1}) \deg(v) \|h_\alpha\| \neq 0$.

        Now let $H$ be any induced global elliptic endoscopic group for $L$, or any group arising in the stabilization recursion by iteratively taking elliptic endoscopic groups induced from the trivial datum on $M'$. Every such $H$ contains $M'$ as a proper Levi, by Proposition \ref{prop.ind.G'.cont.M'.Levi}, and we have by ellipticity
        $$
            \a_{M'} = \a_M, \qquad \a_H = \a_L, \qquad \a_{M'}^H = \a_M^L
        $$
        with compatible choices of measures, as in \S \ref{subsect.cochar.cones.measures}. Dually to the argument of the previous Lemma, the geometric roots of $H$ relative to $T$ are a subset of the geometric roots of $L$, so there are two possibilities.
        \begin{itemize}
            \item If $\pm \alpha \not\in \Phi(H_\infty, T_\infty)$, then $v(\fR_{M'}^H(b')) = 0$. In this case, the unit-resultant computation gives $J_{M',\{v\}}^{H,1}(b') = 0$. 
            \item If $\pm \alpha \in \Phi(H_\infty, T_\infty)$, then $v(\fR_{M'}^H(b')) = 1$ still. In this case, we can repeat the computation above using the same diagonalization $Y$ for $M' \subsetneq H$. The constants remain exactly the same: since on the dual side the roots of $\widehat{H}$ are exactly roots of $\widehat{L}$, as argued above, the projection $h_\alpha$ of the coroot $\alpha^\vee$ is the \emph{same} vector in $\a_M^L = \a_{M'}^H$, whether considered as a coroot for $L$ or $H$. 
        \end{itemize}
        In summary,
        $$
            J_{M', \{v\}}^{H,1}(Y_{b'}) = \begin{cases} C, \quad & \dim \a_M^L = 1, \; \alpha \in \Phi(H_\infty, T_\infty),\\
            0, & \text{otherwise.}
            \end{cases}
        $$

        \subsection{Step 2: Stabilization in relative dimension 1.}
        Given the computations above, it remains to understand the global endoscopic combinatorics appearing in the stabilization recursion. Of course, if $\dim \a_M^L > 1$, then all the orbital integrals appearing in the stabilization vanish, and the WFL is verified
        $$
            J_{M,\{v\}}^{L, s_M}(Y_b) = 0 = S_{M',\{v\}}^{H}(Y_{b'}).
        $$

        Assume henceforth $\dim \a_M^L = \dim \a_{M'}^H = 1$. For any $H$ which does not contain $\pm \alpha$, its iterative descendants also lack $\pm \alpha$, so every orbital integral appearing in the recursion vanishes and
        $$
            S_{M', \{v\}}^{H}(Y_{b'}) = 0.
        $$
        The remainder of the proof reduces now to counting those $H$ which appear in the recursion and retain $\pm \alpha$; the constant $C$ will cancel from both sides.

        Fix a global $H \supset M'$ appearing in the recursion. The stable orbital integral $S^H_{M',\{v\}}(Y_{b'})$ is defined (\ref{eqn.stable.def.recurse}) as a sum over endoscopic parameters $t \in \cE_{M'}(H)$ induced from the trivial datum for $M'$. The space of such parameters is the subset of $Z_{\widehat{M}'}^\Gamma/Z_{\widehat{H}}^{\Gamma}$ which generate an elliptic datum. From Proposition \ref{prop.end.gp.sub.root.system}, we can see that the coroot $\alpha^\vee$ (equivalently, the root $\alpha$) lies in the sub-root datum for the induced group $H_t'$ if and only if $\alpha^\vee(t) = 1$, where $\alpha^\vee : \widehat{T} \to \G_m$ is a character of the dual. So the set of parameters retaining $\pm \alpha$ (not \textit{a priori elliptic}) is
        $$
            \cE_{M'}^{\alpha}(H) :=  \ker\left(\alpha^\vee|_{Z_{\widehat{M}'}^\Gamma}\right)/Z_{\widehat{H}}^{\Gamma} \subseteq Z_{\widehat{M}'}^\Gamma/Z_{\widehat{H}}^{\Gamma}.
        $$
        Every parameter outside of this subset contributes 0 recursively to the definition of the stable orbital integral.
        
        Arthur's definition for induced endoscopic parameters contains the additional condition that they define \emph{elliptic} endoscopic groups. In fact, our dimension condition $\dim \a_M^L = 1$ implies that \emph{every} parameter $t \in \cE_{M'}^\alpha(H)$ defines an elliptic endoscopic $H_t'$ for $H$, as follows. Since $\alpha^\vee|_{Z_{\widehat{M}'}^\Gamma}$ (indeed, even $\alpha^\vee|_{Z_{\widehat{M}}^\Gamma}$) is surjective,\footnote{Because $\alpha^\vee$ pairs nontrivially with a generic dual cocharacter defining $\widehat{M}$.} its kernel has dimension one less than $Z_{\widehat{M}'}^\Gamma$. But one less is exactly the dimension of $Z_{\widehat{H}}^\Gamma$, since $\dim \a_{M'}^H = 1$. In other words,
        $$
            \ker\left(\alpha^\vee|_{Z_{\widehat{M}'}^\Gamma}\right)^0 = \left(Z_{\widehat{H}}^{\Gamma}\right)^0.
        $$
        But the dual induced group $\widehat{H}_t'$ contains $\widehat{M}'$ as a Levi and also has $\alpha^\vee$ as a root, so its $\Gamma$-fixed center (the two actions of $\Gamma$ coincide on the center $Z_{\widehat{H}}$) satisfies
        $$
            Z_{\widehat{H}}^\Gamma \subseteq Z_{\widehat{H}_t'}^{\Gamma} \subseteq Z_{\widehat{M}'}^\Gamma \cap \ker(\alpha^\vee) = \ker\left(\alpha^\vee|_{Z_{\widehat{M}'}^\Gamma}\right).
        $$
        Taking identity components, these become equalities, which is characterization (d) of ellipticity in Proposition \ref{prop.ell.end.data}. 

        So we have shown that the number of nonzero endoscopic terms (including the trivial $t = 1$) appearing in the stabilization recursion is
        $$
            N_H := \left|\ker\left(\alpha^\vee|_{Z_{\widehat{M}'}^\Gamma}\right)/Z_{\widehat{H}}^{\Gamma}\right|
        $$
        Using this observation, we can recursively compute the stable orbital integral. We claim that for any $H$ appearing in the recursion, $S^H_{M', \{v\}}(Y_{b'})$ is equal to
        \begin{equation}
            N_H^{-1} C = \left|\ker\left(\alpha^\vee|_{Z_{\widehat{M}'}^\Gamma}\right)/Z_{\widehat{H}}^{\Gamma}\right|^{-1} C.\label{eqn.stab.rk1.recurs.formula}
        \end{equation}
        Inducting on $\dim H$, every nonidentity parameter $t$ defines a proper endoscopic group $H_t'$, for which the inductive hypothesis applies, giving a recursive contribution
        $$
            \left| Z_{\widehat{H}_t'}^{\Gamma}/Z_{\widehat{H}}^\Gamma\right|^{-1} S^{H'_t}_{M',\{v\}}(Y_{b'}) = \left| Z_{\widehat{H}_t'}^{\Gamma}/Z_{\widehat{H}}^\Gamma\right|^{-1} N_{H_t}^{-1} C = N_H^{-1} C
        $$
        Equation (\ref{eqn.stable.def.recurse}) therefore becomes
        $$
            S^H_{M', \{v\}}(Y_{b'}) = J^{H, 1}_{M', \{v\}}(Y_{b'}) - \sum_{\substack{t \in \cE_{M'}(H)\\t \neq 1}} \left| Z_{\widehat{H}_t'}^{\Gamma}/Z_{\widehat{H}}^\Gamma\right|^{-1} S^{H'_t}_{M',\{v\}}(Y_{b'}) = C - \left(N_H - 1\right) N_H^{-1} C = N_H^{-1} C,
        $$
        as claimed. 

        Finally, it remains to consider the combinatorics appearing in the actual WFL, (\ref{eqn.V.wfl}) with $V = \{v\}$. Applying the same argument from $M' \subset H$ to $M \subset L$, the only induced endoscopic groups $L_s'$ for $s = s_M z$ which contribute to the endoscopic side are those containing $\alpha$, which occurs precisely when
        $$
            \alpha^\vee(s) = 1.
        $$
        There exists at least one such $z$ for which the equation holds because $\alpha^\vee|_{Z_{\widehat{M}}^\Gamma}$ is surjective. Its solution set is therefore a trivial torsor for
        $$
            \ker\left(\alpha^\vee|_{Z_{\widehat{M}}^\Gamma}\right)/Z_{\widehat{L}}^{\Gamma},
        $$
        and all such parameters are elliptic by the same $\dim \a_M^L = 1$ argument as above. We write $\cE^{\alpha}_{\fM'}(L)$ for the set of such parameters, whose cardinality is
        $$
            \left|\ker\left(\alpha^\vee|_{Z_{\widehat{M}}^\Gamma}\right)/Z_{\widehat{L}}^{\Gamma}\right|.
        $$
        The right side of (\ref{eqn.V.wfl}) then becomes, since stable endoscopic orbital integrals for $s \not\in \cE^{\alpha}_{\fM'}(L)$ vanish,
        \begin{align*}
            |Z^\Gamma_{\widehat{M}'}/ Z_{\widehat{M}}^\Gamma| \sum_{s \in \cE_{\fM'}^{\alpha}(L)} |Z_{\widehat{L'_s}}^\Gamma/ Z_{\widehat{L}}^\Gamma |^{-1} S_{M', \{v\}}^{L'_s}(Y_{b'}) &= |Z^\Gamma_{\widehat{M}'}/ Z_{\widehat{M}}^\Gamma| \sum_{s \in \cE_{\fM'}^{\alpha}(L)} |Z_{\widehat{L'_s}}^\Gamma/ Z_{\widehat{L}}^\Gamma |^{-1} N_{L_s'}^{-1} C \\
            &= |Z^\Gamma_{\widehat{M}'}/ Z_{\widehat{M}}^\Gamma| \sum_{s \in \cE_{\fM'}^{\alpha}(L)} \left|\ker\left( \alpha^\vee|_{Z_{\widehat{M}'}^\Gamma}\right)/Z^\Gamma_{\widehat{L}}\right|^{-1} C \\
            &= |Z^\Gamma_{\widehat{M}'}/ Z_{\widehat{M}}^\Gamma| \left|\ker\left(\alpha^\vee|_{Z_{\widehat{M}'}^\Gamma}\right)/\ker\left(\alpha^\vee|_{Z_{\widehat{M}}^\Gamma}\right)\right|^{-1} C.
        \end{align*}
        Note that the index of the kernel terms is finite, since both have identity component $(Z_{\widehat{L}}^\Gamma)^0$. Since $\alpha^\vee$ is surjective on $Z_{\widehat{M}}^\Gamma$, the second isomorphism theorem applied to the kernel of $\chi := \alpha^\vee|_{Z_{\widehat{M}'}^\Gamma}$ gives
        $$
            \ker\left(\alpha^\vee|_{Z_{\widehat{M}'}^\Gamma}\right)/\ker\left(\alpha^\vee|_{Z_{\widehat{M}}^\Gamma}\right) = \ker(\chi)/\left(\ker(\chi) \cap Z_{\widehat{M}}^\Gamma\right) \cong \ker(\chi) Z_{\widehat{M}}^\Gamma/ Z_{\widehat{M}}^\Gamma = Z_{\widehat{M}'}^\Gamma/Z_{\widehat{M}}^\Gamma,
        $$
        so the scalars cancel and we have only $C = J_{M,\{v\}}^{L,s_M}(Y_b)$. This proves the $\fM'$-WFL on $v$ when $L \supsetneq M$, so the last requirement for peelability of $v$ is established.
    \end{proof}
\end{lemma}

\section{Local constancy of weighted orbital integrals.}

We prove a very weak version of local constancy of local $s$-orbital integrals over finite residue field which is sufficient for our purposes.
\begin{proposition}[{cf.\ \cite[\S 18.7]{CLII}}]\label{prop.local.constancy}
    Let $b \in \fc_M(\O_{0,v}) \cap \fc_M^{\Greg}(F_{0,v})$, and write $X = \epsilon_M(b) \in \fm(\O_{0,v})$. There exists an integer $n$ such that for every $b' \in \fc_M(\O_{0,v})$ satisfying
    $$
        b' \equiv b \mod z_v^n,
    $$
    and the corresponding lift $X' \in \epsilon_M(b')$, there exists
    $$
        h \in M(\O_{0,v}), \qquad h \equiv 1 \mod z_v
    $$
    such that
    \begin{enumerate}
        \item conjugation by $h$ identifies the regular centralizers $J_b \cong J_{b'}$,
        \item for every $g \in G(F_{0,v})$,
        $$
            \Ad(g^{-1}) X \in \fg(\O_{0,v}) \Longleftrightarrow \Ad(h^{-1}g^{-1}) X' \in \fg(\O_{0,v}),
        $$
        \item for every $F_{0,v}$-parabolic $P$ with Levi $M$,
        $$
            H_{P,v}^{\loc}(hg) = H_{P,v}^{\loc}(g).
        $$
    \end{enumerate}
    In particular, the local $s$-weighted orbital integrals $J_{M,v}^{G,s}$ associated to $b$ and $b'$ are equal.
    \begin{proof}
        Consider the action map over $F_{0,v}$
        $$
            \Phi : M \times \fj \to \fm, \quad (m, Y) \mapsto \Ad(m)Y,
        $$
        where $\fj = \Lie(J_b)$. A differential calculation and semisimplicity of $X$ shows this is smooth at $(1, X)$. Fix arbitrarily small neighborhoods 
        $$
            U_h := \ker(M(\O_{0,v}) \to M(\O_{0,v}/z_v^N)) \subset M(F)
        $$
        of 1 and 
        $$
            U_Z := X + z_v^{N_1} \Lie(J_b)(\O_{0,v}) \subset \fj(F)
        $$
        of $X$. Hensel's lemma as a $p$-adic inverse function theorem says that $\Phi(U_h, U_Z)$ contains a sufficiently small ball around $\Ad(1)X = X$. In other words, for any $X' \equiv X \mod z_v^n$, there exist $h \approx 1 \in U_h, Z\approx X \in U_Z$ such that 
        $$
            X' = \Ad(h)Z.
        $$
        Since $\epsilon_M$ is defined over $\O_{0,v}$, this applies to our nearby Hensel lift $X'$. By definition, $Z \in \fj(F)$, so that $[X,Z] = 0$, and $Z - X \in z_v^{N_1} \Lie(J_b)(\O_{0,v})$.

        Since $Z \in \fj(F)$, we have at least 
        $$
            J_b \subseteq Z_M(Z).
        $$
        But $Z$ is sufficiently close to regular semisimple $X$, so it is also regular semisimple and its centralizer is a maximal torus. We conclude 
        $$
            J_b = Z_M(Z),
        $$
        So $\Ad(h)$ sends $J_b$ to $Z_M(X') = J_{b'}$ over $F_{0,v}$. But $h$ is integral, so the conjugation extends to an $\O_{0,v}$-morphism
        $$
            \Ad(h) : J_b \to M_{\O_{0,v}}.
        $$
        Generically it lands in $J_{b'}$, and $J_b$ is flat, so the morphism factors through $J_{b'}$. The same argument applied to the inverse shows $J_{b} \cong J_{b'}$.

        The second claim is the heart of the proof, and proceeds as in \cite[Proposition 18.7.1]{CLII}. We will show only the forward claim; the converse is almost a perfectly symmetric argument. We want to show that $\Ad(g^{-1})(X) \in \fg(\O_{0,v})$ if and only if $\Ad(g^{-1})(Z) \in \fg(\O_{0,v})$. Fix an $F_{0,v}$-parabolic $P$ with Levi $M$ and consider the iwasawa decomposition $P = MN$ of $g$
        $$
            g = mnk, \qquad k \in G(\O_{0,v}).
        $$
        Suppose that $\Ad(g^{-1})(X) \in \fg(\O_{0,v})$. We can discard $k$. By assumption, writing $Y = \Ad(m^{-1})(X)$ we have
        $$
            \Ad((mn)^{-1})(X) = (\Ad(n^{-1})Y - Y) + Y.
        $$
        The difference lies in $N(F_{0,v})$, and $Y$ lies in $M(F_{0,v})$. By the Levi decomposition there is no intersection, so since the sum of the two terms is integral each term must be likewise. The $\O_{0,v}$-morphism $n \mapsto \Ad(n^{-1}) Y - Y$ is an $F_{0,v}$-isomorphism whose inverse has denominators of valuation bounded by a fixed power of $\det(\Ad(Y)|_{\fn})$. This boundedness can be checked after a finite unramified extension which splits the group, where we use \cite[Lemme 18.7.2]{CLII}. Note that $Y$ and $X$ have the same determinant valuation on $\fn$ because conjugation by $m$ intertwines the endomorphisms. So we have shown that
        $$
            \Ad(g^{-1})(X) \in \fg(\O_{0,v})
        $$
        implies $n$ lies in a bounded subset of $B \subset N(F_{0,v})$, independent of $g$. 

        Now we need to bound $\delta_m =\Ad(m^{-1})\delta$, where $\delta = Z - X$. The universal property of the regular centralizer and integrality of $Y$ implies that the $F_{0,v}$-map
        $$
            \Ad(m^{-1}) : J_{b} \to Z_M(Y)
        $$
        extends over $\O_{0,v}$. It follows that 
        $$
            \Ad(m^{-1})\Lie(J_b)(\O_{0,v}) \subset \fm(\O_{0,v}).
        $$
        Since $\delta \in z_v^{N_1} \Lie(J_b)(\O_{0,v})$, we have $\delta_m \in z_v^{N_1} \fm(\O_{0,v})$. Choosing $N_1$ large enough relative to the bounded subset $B \subset N(F_{0,v})$ that $n$ lives in, we have
        $$
            \Ad(n^{-1})(Y + \delta_m) = Y + (\Ad(n^{-1} Y - Y) + \delta_m + (\Ad(n^{-1} \delta_m - \delta_m)
        $$
        with all four pieces integral, i.e., in $\fg(\O_{0,v})$. But 
        $$
            \Ad(g^{-1})Z = \Ad(n^{-1}) \Ad(m^{-1})(X + \delta) = \Ad(n^{-1}) (Y + \delta_m),
        $$
        so this gives the second claim. 

        For the third claim, since $h \equiv 1 \mod z_v$ is an integral element of $M$, every rational character takes unit value on $h$, and so $H_{M,v}^{\loc}(h) = 0$. But
        $$
            H_{P,v}^{\loc}(hg) = H_{M,v}^{\loc}(h) + H_{P,v}^{\loc}(g),
        $$
        so $v_M^{G,\loc}(hg) = v_M^{G, \loc}(g)$. When global data are fixed, projection and multiplication by $\deg(k_v/k_0)$ give the same invariance for the global heights, and so likewise for both the global weight function $v_M^G$ of Arthur and the weight function $w_M^\xi$ of Chaudouard--Laumon \ref{def.weights.CL}.

        It remains to show $J_{M,v}^{G,s}(Y_b) = J_{M,v}^{G,s}(Y_{b'})$. The isomorphism $\Ad(h) : J_b \to J_{b'}$ induces a bijection of the domains and identifies bounded subgroups, so preserves the normalizations of the measures. The discriminant valuation factor is continuous, i.e., locally constant. Finally, the exact sequence in $\tau$-cohomology and the Kottwitz isomorphism are natural along $\Ad(h)$, so 
        $$
            \inv_{b'}(hg) = \inv_b(g).
        $$
        If one prefers the presentation of the invariant pairing in terms of $\widehat{T}^{W_{b,v}}$, then one should be cautious about the Weyl group ambiguity implicit in the choice of geometric point $\wtilde{b} \mapsto b$ upstairs in the cameral cover. However, making a choice of lift for ${b}$ induces a corresponding such choice for every $b'$ sufficiently close to $b$ since the Weyl group is finite---choose a neighborhood of $b$ small enough that its preimage is a disjoint union of lifts, and then the component containing $\wtilde{b}$ is the distinguished one. We have shown that all requisite parts of the $s$-orbital integral satisfy the desired local constancy.
    \end{proof}
\end{proposition}

\section{Globalization and the proof of the WFL.}\label{subsect.globalize.param}

We are at last prepared to prove Theorem \ref{thm.wfl}, the Lie algebra version of Arthur's WFL for Lie algebras. The final stage is that of globalizing our arbitrary local parameter in a way that satisfies all our global hypotheses, and is mostly analogous to \cite[Theorem 19.1.1]{CLII}, \cite[\S 8.6]{ngofl}, except that we make use of Poonen's formalism for Bertini theorems over finite fields in order to work directly over our given finite field. We begin from our local setting as in \S \ref{subsect.local.setup}. In particular, $\O_{0,v}$ is the ring of integers of equal-characteristic local field $F_{0,v}$ with residue field $k_v$. As in the statement of Theorem \ref{thm.wfl}, let $\fM_v' = (s_M, \rho_{v, s_M})$ be an unramified local endoscopic datum with endoscopic group $M'/\O_{0,v}$. Let $\Gamma_v$ be an unramified (hence cyclic) Galois group large enough to split $G$ and $M'$, assumed prime to $p$. 

\subsection{Globalization of the group and endoscopic datum.}
We recall the requisite globalization statement. It is well-known but we provide a proof for completeness.
\begin{lemma}[{cf.\ \cite[\S 8.6]{ngofl}}]\label{lem.globalize.group}
    Let $k_v$ finite of characteristic $p$, and $\Gamma_v = \langle \gamma \rangle \cong \Z/n\Z$ for $(n, p) = 1$. Then there exists $X_0/k_v$ smooth, projective, geometrically connected, two distinct points $v, \infty \in X_0(k_v)$, and a geometrically connected finite \'etale $\Gamma_v$-torsor
    $$
        X_0' \to X_0
    $$
    such that $\tau_v \mapsto \gamma$, $\tau_\infty \mapsto 1$. Equivalently, $D_v = \Gamma_v$ and $D_\infty = 1$.
    \begin{proof}
        Choose any elliptic curve $E/k_v$. For some finite extension $k_m/k_v$ we have $E[n](k_m) \cong (\Z/n)^2$. So $A = \Res_{k_m/k_v}(E_{k_m})$ is an abelian variety such that there exists $\psi : A(k_v) \twoheadrightarrow \Gamma_v$. Choose $a \in A(k_v)$ mapping to generator $\gamma$, and consider the Lang map
        $$
            L_A = \on{Frob}_q - 1 : A \to A.
        $$
        Quotienting the Lang cover by $H = \ker(\psi) \subset A(k_v)$ gives an \'etale $\Gamma_v$-cover 
        $$
            A/H \to A.
        $$
        Frobenius acts on the fiber over $x \in A(k_v)$ by $\psi(x)$, so at the identity $x = e$ it acts trivially and at $x = a$ it acts by $\gamma$. By \cite[Corollary 3.4]{poonenbertini} there exists a smooth projective geometrically integral curve $X_0 \subset A$ containing both $e$ and $a$. Writing $\infty := e$ and $v := a$ and pulling back gives our desired \'etale cover
        $$
            X_0' = X_0 \times_A (A/H) \to X_0,
        $$
        The arithmetic Frobenius has trivial image (forced by the split point over $\infty$) but the arithmetic monodromy is surjective $\pi_1(X_0, \infty) \twoheadrightarrow \Gamma_v$. The juxtaposition ensures a surjection from the fundamental group of $X = (X_0)_{\ol{k_v}}$:
        $$
            \pi_1(X, \ol{\infty}) \twoheadrightarrow \Gamma_v
        $$
        so the cover is geometrically connected.
    \end{proof}
\end{lemma}

Our local datum prescribes an action of $\Gamma_v$ on a split pinned model whose corresponding twist along the unramified cover of the disk $\Spec \O_{0,v}$ defines $(G_v, M_v, \fM'_v)$. Instead twisting the constant objects by the same action over the cover $X'_0 \to X_0$ constructs
\begin{enumerate}
    \item $G/X_0$ quasisplit reductive group with $X_0$-maximal torus $T$ satisfying the global characteristic hypotheses of \S \ref{subsect.char.assumptions}, split at $\infty$, and such that $G|_{\Spec \O_{0,v}} \cong G_v$;
    \item $M/X_0$ semistandard Levi such that $M|_{\Spec \O_{0,v}} \cong M_v$;
    \item $\fM'$ an elliptic endoscopic datum for $M$ on $X_0$ with endoscopic group $M'/X_0$ such that $\fM'|_{\Spec \O_{0,v}} \cong \fM'_v$ and $M'$ is split at $\infty$.
\end{enumerate}
We fix a finite set $S \subset |X_0| \setminus \{v, \infty\}$ of points. For convenience, by Chebotarev we are free to choose them to split in $X'$ and to have sufficiently high residue field such that they admit regular elements in the split torus. We will define $D$ to be a large even effective divisor with $\on{supp}(D) = S$, say by fixing such a $D_1$ to begin with and then defining $D$ to be a sufficiently positive even multiple of $D_1$.

\subsection{Specifying finitely many jets.}

Now fix a local parameter $b' \in \fc_{M'}(\O_{0,v}) \cap \fc_{M'}^{\Greg}(F_{0,v})$, write $b \in \fc_M(\O_{0,v})$ for its image, and write 
$$
    Y_{b'} = \epsilon_{M'}(b') \in \fm'(\O_{0,v}), \qquad Y_b = \epsilon_M(b) \in \fm(\O_{0,v})
$$
for the images under Kostant sections. 

Fix $n$ a positive integer greater than $v(\fD_G(b))$ and large enough to satisfy the conclusion of Proposition \ref{prop.local.constancy} for not only $M \subset G$, but for every endoscopic pair which appears in the (finite) recursion defining the local stable orbital integrals of Definition \ref{def.stable.integral.local} and in the statement of the WFL \ref{eqn.local.wfl}. In particular, this ensures that if a global parameter $a_{M'}$ agrees with the local parameter $b'$ up to order $n$, then it and all its images under any Chevalley maps will satisfy a sufficiently high congruence to identify the corresponding orbital integral values.

For any finite effective divisor $Z$ on $X_0$, the restriction map
$$
    \res_Z : \cA_{M'}(k_v) = H^0(X_0, \fc_{M',D}) \to H^0(Z, \fc_{M',D}|_Z)
$$
is surjective for $D > 2g - 2 + \deg Z$. In particular, we can prescribe finite-order jets for the restriction to finitely many points on the curve, and for $D \gg 0$ we will be able to find a global section of $\cA_{M'}(k_v)$ mapping to it. In particular, we prescribe the following conditions:
\begin{itemize}
    \item At $v$, we prescribe $$a_{M'}|_{nv} = b'|_{nv},$$
    where $n$ is the local constancy parameter. In other words, we prescribe the $n$-th order jet $a_{M'}|_{\O_{0,v}} \equiv b' \mod z_v^n$. Since $n > v(\fD_G(b))$, all such $a$ have the same discriminant valuation as $b$.
    \item At $\infty$, we prescribe $a_{M'}(\infty) = \chi_{M'}(t_\infty)$ for a fixed $t_\infty \in \ft_{M',D,\infty}^{\Greg}(k_v)$.
    \item At each $x \in S = \on{supp}(D)$, we prescribe $a_{M'}(x)$ to match a given arbitrary regular section of $\ft_{D}(x)$, and in particular can certainly enforce $G$-regularity at each such point with a 0-jet.
    \item At one other auxiliary $x_{\mathrm{ell}}$, chosen to split completely in $X' \to X$ and have sufficiently large residue degree, we enforce ellipticity in the special fiber as follows: choose an elliptic Weyl element $w \in W_{M'}(\ol{x}_{\el}) = W_{M'}(x_{\mathrm{ell}})$, and consider the $k_{x_{\mathrm{ell}}}$-rational torus $T_w \subset M'|_{x_{\mathrm{ell}}}$ (such that Frobenius acts on the character lattice by $w$). This torus is elliptic, and for the residue field chosen sufficiently large there exists a $G$-regular
    $$
        Y \in \Lie(T_w)(k_{x_{\mathrm{ell}}}).
    $$
    We prescribe the 0-jet $a_{M'}(x_{\mathrm{ell}}) = \chi_{M'}(Y)$.
\end{itemize}
That is, $Z := n[v] + \infty + S + x_{\mathrm{ell}}$. The last condition ensures $M'$-ellipticity of the global section $a_{M'}$: the decomposition group at $x_{\mathrm{ell}}$ contributes an elliptic element conjugate to $(w,1)$ to the global cameral monodromy. This shows that $a_{M'}$ does not have monodromy subgroup lying within $W_{L'}(\infty) \rtimes \im(\rho_G)$ for any proper Levi, since $w \not\in W_{L'}(\infty)$, so by Proposition \ref{prop.ngo.comparison.wa}, the element is $M'$-elliptic with \emph{any} choice of rigidifying datum $t$. In particular, it lies in the compatible locus of Definition \ref{def.comp.loc} for $M'$ vacuously, and its image is $M$-elliptic by Proposition \ref{prop.ell.end.data} (whether considered arithmetically or geometrically, equivalent as in Remark \ref{rmk.arith.geo.data.correspond}).

\subsection{Bertini--Poonen transversality.}

So far we have only enforced finitely many conditions, so for $D = mD'$ sufficiently positive (supported on $S$) the map $\res_Z$ is surjective. Moreover, the fiber $V_m \subset H^0(X_0, \fc_{M',mD'})$ is an affine space, a torsor for $H^0(X_0, \fc_{M', mD'}(-Z))$. It remains to handle the transversality condition, which we choose to do directly on $\F_q$ using a Bertini--Poonen sieve as in \cite{poonenbertini}. Consider the reduced divisor
$$
    \Delta = (\fD_G|_{\fc_{M',D}})_{\red} = \left(\fD_{M'} + \fR_{M'}^M + (\chi_M^{M'})^\star \fR_M^G \right)_{\red}.
$$
At every closed point $x \not\in \on{supp}(Z)$, we wish to impose that $a_{M'}$ meets $\Delta$ with multiplicity at most 1. This is a first-jet condition.

\begin{lemma}[{\cite[Proposition 5.1.1]{ngoetal}, cf.\ \cite[Theorem 1.2]{poonenbertini}, \cite[Theorem 8.1]{poonensquarefree}}]\label{lem.bert.glob}
    Write $D = mD'$ for $D'$ finite effective and supported on $S$. For $m$ sufficiently positive, there exists
    $$
        a_{M'} \in H^0(X_0, \fc_{M',mD'}) = \cA_{M'}(k_v)
    $$
    such that $a_{M'}$ has the prescribed jets above on $Z$ and $a^\star_{M'} \Delta$ is reduced on $X_0 \setminus |Z|$.
    \begin{proof}
        This is essentially the first result cited. We outline how to adapt the arguments of Poonen to give the Lemma. For each $x \not\in |Z|$, we wish to bound the size of
        $$
            B_x(m) = \{a \in V_m : v_x(a^\star \Delta) \ge 2\}.
        $$
        Let $r_m$ denote the rank of $\fc_{M',D}$ as a vector bundle, and write $q_x$ for the size of the residue field at $x$. The first-order jet space for $a$ at $x$ has $q_x^{2r_m}$ points. If $y = a_{M'}(x)$ lies in a smooth point of $\Delta_x$, then there are $O(q_x^{r_m-1})$ points $y$ can equal, and tangency imposes a single linear condition on the derivative, leaving $q_x^{r_m-1}$ values of the derivative. In total, there are 
        $$
            O(q_x^{2r_m - 2}) 
        $$
        bad jets at smooth points. On the other hand, since $\Delta$ is defined inside $\fc_{M',D}$ by twisting from the constant divisor, inside each fiber we have
        $$
            \codim_{(\fc_{M',D})_x}(\Delta_x^{\mathrm{sing}}) \ge 2;
        $$
        the divisor is reduced so its singular locus is of higher codimension. So the set of singular values $y$ can take is of order $O(q_x^{r_m - 2})$, and the derivative can be chosen freely. So we have again 
        $$
            O(q_x^{2r_m - 2})
        $$
        bad jets. In particular, if the space of jets is equidistributed, then the probability of a bad jet at the point $x$ is
        $$
            \mathbb{P}(B_x(m)) = q_x^{-2} = q^{-2\deg(x)}.
        $$

        We divide the closed points of $X$ into 3 ranges, as in \textit{loc.\ cit.} \begin{enumerate}
            \item (Low-degree points). Fix $R$, independent of $m$, and consider the finite collection of points satisfying $\deg(x) < R$. Just as done above, for $m \gg 0$ the first-order jets at such points can be imposed directly, and for $m >> 0$ the restriction map
            $$
                V_m \to \prod_{\substack{x \not\in Z \\ \deg(x) < R}} H^0(2x, \fc_{M',mD'})
            $$
            is surjective. The map is affine-linear, so its fibers are equal in size and points $a_{M'} \in V_m$ have images which equidistribute over the codomain. The probability of avoiding bad jets at all low-degree points is then
            $$
                \mathbb{P}(\text{good at all $x$, $\deg(x) < R$}) = \prod_{\deg(x) < R} (1 - q^{-2 \deg(x)}) > 0.
            $$
            \item (Medium-degree points). There is some constant $c > 0$ such that for points with $\deg(x) \le cm$, the restriction map $V_m \to H^0(2x, \fc_{M', mD'})$ is still surjective: the surjectivity holds whenever
            $$
                H^1(X_0, \fc_{M', mD'}(-Z - 2x)) = 0,
            $$
            Indeed,
            $$
                \deg(D) > 2g - 2 + \deg Z + 2\deg x
            $$
            suffices. For each degree $R \le e \le cm$, the number of $\F_{q^e}$-points on the curve is bounded by $O(q^{e})$ (cf. \cite[Lemma 2.4]{poonenbertini}), and each closed point of degree $e$ contributes $e$ such. Summing, we see by the union bound that the probability that there exists a point with a bad jet in the range $R \le \deg(x) \le cm$ is
            $$
                \mathbb{P}(\exists x, R \le \deg(x) \le cm, B_x(m)) \le \sum_{e = R}^{cm} O(q^e/e) q^{-2e} \le \sum_{e \ge R} q^{-e} = O(q^{-R}).
            $$
            In particular, we can choose $R$ originally such that 
            $$
                \limsup_{m \to \infty} \mathbb{P}(\exists x, R \le \deg(x) \le cm, B_x(m)) < \epsilon,
            $$
            so taking $m$ sufficiently large the bad probability will be bounded arbitrarily close to this limit, say $< 2\epsilon$.

            \item (High-degree points). Consider now the set of points satisfying $\deg(x) > cm$, for which jet-surjectivity fails. Here we use a different estimate of Poonen for the density of squarefree values of a squarefree multivariable function. Deleting finitely many additional points, we can assume $\fc_{M',D}$ is trivial on $U := X_0 \setminus (S \cup |Z| \cup S_1)$. Since $D$ was supported on $S$, the resulting family of sections is independent of $m$. Choosing a trivialization as a family of fiber coordinates $y_1, \ldots, y_r$, the reduced discriminant divisor has on $U$ an equation
            $$
                \Delta|_U = V(F(y_1,\ldots, y_r))
            $$
            with $F \in H^0(U, \O_U)[y_1,\ldots, y_r]$ squarefree. Then \cite[Theorem 8.1]{poonensquarefree} implies that a positive proportion of sections over $U$ are squarefree. Moreover, since the allowed degrees of poles on $S$ go to infinity, the proof shows that the sections can be chosen to globalize to $H^0(X_0, \fc_{M',D})$. See the proof of \cite[Proposition 5.1.1]{ngoetal} for more details. 
        \end{enumerate}
        Combining these three, we find that a positive proportion of sections in $V_m$ are transverse to $\Delta$ on $X_0 \setminus |Z|$, and in particular such a globalization exists.
    \end{proof}
\end{lemma}

\subsection{Goodness.}

Finally, it remains to check that when $\deg(D)$ is chosen sufficiently large relative to the complexity of the parameter $b'$, our globalizations $a_{M'}$ necessarily land inside the good locus $\cA_G^{\good}$ of Definition \ref{def.good.locus}. This follows from Theorem \ref{thm.suff.criterion.for.goodness}, taken with $S = \{v\}$, as long as 
$$
    \deg(D) > 2g - 2 + \frac{|W_G|}{2} v(\fD_G(a)),
$$
which is determined by the local parameter $b'$.

\subsection{Proof of the WFL}

\begin{proof}[Proof of Theorem \ref{thm.wfl}]
    Choose $D = 2 m D'$ sufficiently positive for the Bertini lemma and such that
    $$
        \deg(D) > \max \left\{ 2g, 2g - 2 + \deg Z, 2g-2 + \frac{|W_G|}{2} v(\fD_G(b))\right\}.
    $$
    Any globalized parameter $a_{M'}$ as in the conclusion of Lemma \ref{lem.bert.glob} is $M'$-elliptic by the condition at $x_{\mathrm{ell}}$ and lies in the good locus by Theorem \ref{thm.suff.criterion.for.goodness}. Thus, choosing any $t$ to give a parameter $(a_{M'}, t) \in \bA_{M'}$, the image $(a_M, t) \in \bA_M$ is $M$-elliptic by Proposition \ref{prop.ell.end.data}. The parameter (as well as the endoscopic datum $\fM'$) therefore satisfy all the hypotheses of Theorem \ref{thm.global.wfl}, which implies that the $\fM'$-WFL holds on $V = |X_0|$. Peelability of the regular places then tells us that the $\fM'$-WFL holds on the finite set $\{v\} \cup \on{supp}(a_{M'}^\star \fD_{M'} + a_{M'}^\star \fR_{M'}^M + a_M^\star \fR_M^G) \subset |X_0|$. We have engineered our global parameter to intersect this divisor transversely at smooth points, so Lemmas \ref{lem.transv.disc}, \ref{lem.tranv.result.1}, and \ref{lem.tranv.result.2} imply that all places other than $v$ are peelable. So the Peeling Lemma \ref{lem.peeling} implies that the $\fM'$-WFL for $M \subset G$ holds on $\{v\}$. Since we have engineered our global cover to satisfy $\Gamma = D_v = \Gamma_v$ and $k_v = k_0$, as in Remark \ref{rmk.monodr.local.global.discrep.wfl} this recovers the local statement of the WFL. 
\end{proof}

\chapter{The nonstandard WFL}\label{sect.nonstd}

In \cite{Chaudouardnonstandardexample}, the nonstandard WFL is discussed for the nonstandard pair $\Sp_{2n}$ and $\on{SO}_{2n+1}$. There, it is proposed (at least in the split case) that the local nonstandard WFL should follow from the global cohomological Theorem \ref{thm.no.new.supports} and our stable point count, Theorem \ref{thm.stable.point.count}, followed by a splitting/peeling formalism like that performed in the last chapter or in \cite[\S 16-18]{CLII}. This chapter realizes that proposal in full generality. In \S \ref{subsect.nonstd.setup}, we recall the formalism of nonstandard endoscopy as in \cite{WaldTordu} and \cite{ngofl}. Then we deduce from our stable point count Theorem \ref{thm.stable.point.count} a global version of the nonstandard WFL.

\section{Groups \textit{appari\'es} and the statement of the nonstandard WFL.}\label{subsect.nonstd.setup}

We recall the formalism of the nonstandard WFL. Let $\G_1, \G_2$ be two split reductive groups over algebraically closed field $k$ equipped with pinnings. Associated to each we have its based root datum
$$
    (X^\star(\T_i), X_\star(\T_i), \ul{\Phi}_i, \ul{\Phi}_i^\vee, \Delta_i).
$$
\begin{definition}[{cf.\ \cite[D\'efinition 1.12.1]{ngofl}, \cite[\S 3.7]{WaldTordu}}]
    An \emph{isogeny} of (based) root data $\psi$ from $\G_1$ to $\G_2$ consists of a pair of transpose isomorphisms
    $$
        \psi^\star : X^\star(\T_2)_{\Q} \to X^\star(\T_1)_{\Q}
    $$
    and
    $$
        \psi_\star : X_\star(\T_1)_{\Q} \to X_\star(\T_2)_{\Q}
    $$
    which identify the roots up rational scalars in the following sense
    \begin{enumerate}
        \item $\psi^\star$ identifies the set of $\Q$-spans of the roots of $\ul{\Phi}_2$ (resp.\ the simple roots $\Delta_2$) bijectively onto the set of $\Q$-spans of the roots of $\ul{\Phi}_1$ (resp.\  $\Delta_1$), and
        \item $\psi_\star$ identifies the set of $\Q$-spans of the coroots $\ul{\Phi}_1^\vee$ (resp.\ the simple coroots $\Delta_1^\vee$) onto the set of $\Q$-spans of $\ul{\Phi}_2^\vee$ (resp.\ $\Delta_2^\vee$).
    \end{enumerate}
    We say that $\G_1, \G_2$ are \emph{paired} if there exists an isogeny of based data between them.
\end{definition}

\begin{example}
    If $\G_1, \G_2$ are semisimple and agree up to isogeny (i.e., have the same adjoint form), then they have isogeneous root data. More interestingly, if $\G$ is a simple group and $\widehat{\G}$ is its split dual group over $k$, then the two have isogeneous root data.
\end{example}

\begin{definition}
    Let $\G_1, \G_2$ be paired split groups over $k$. Define $\Out_{12}$ as the set of automorphisms of $X^\star(\T_1)_{\Q}$ which preserve $X^\star(\T_1), \ul{\Phi}_1, \Delta_1$ but also preserve $X^\star(\T_2), \ul{\Phi}_2, \Delta_2$ under transport along $\psi^\star$. Necessarily then, for $\sigma \in \Out_{12}$ the transpose $\sigma^\vee$ preserves the cocharacter/coroot data for both groups under $\psi_\star$.

    We say that two (pinned) quasisplit reductive groups $G_1, G_2$ over a base $S$ are \emph{paired} if they are defined by Galois-twisting paired split groups $\G_1, \G_2$ along a common $\Out_{12}$-torsor $\rho_{12}$. Note that a local \'etale cover which splits $T_1$ also splits $T_2$, and vice versa.
\end{definition}

\begin{definition}\label{def.good.char.pair}
    We say that the characteristic of $S$ is \emph{good with respect to $\psi$} if at every point $s \in S$ with residue characteristic $p$, 
    $$
        \psi^\star : X^\star(\T_2) \otimes \Z_p \xrightarrow{\sim} X^\star(\T_1) \otimes \Z_p.
    $$
    We say that the characteristic of $S$ is \emph{good with respect to $(G_1, G_2, \psi)$} if it is good with respect to $\psi$ and also good for each group $G_1, G_2$ separately, in the sense of Definition \ref{def.good.characteristic}, with respect to the corresponding tori. 
\end{definition}
    In particular, good characteristic implies that for every root $\alpha_1 \in \ul{\Phi}_1$, the rational scalar $b(\alpha_1)$ such that $\psi^\star(\alpha_2) = b(\alpha_2) \alpha_1$ for a root $\alpha_1 \in \ul{\Phi}_1$ should be a \emph{$p$-adic unit}. This implies induced isomorphisms at the Lie algebra level, so paired quasisplit groups in good characteristic have a natural transfer of their invariant-theoretic data. 
\begin{proposition}[{\cite[Lemme 1.12.6, Lemme 2.4.7]{ngofl}}]\label{prop.paired.gps.invariants.correspond}
    Let $(G_1, G_2, \psi)$ be paired groups over a base $S$ in good characteristic. Let $T_1, T_2$ the maximal $S$-tori corresponding to the pinnings. Then there exist isomorphisms over $S$
    $$
        \ft_1 \xrightarrow{\sim} \ft_2, \qquad W_{G_1} \xrightarrow{\sim} W_{G_2}
    $$
    and an induced $S$-isomorphism of the characteristic spaces
    $$
        \nu : \fc_{G_1} \xrightarrow{\sim} \fc_{G_2}
    $$
    compatible with the Chevalley maps for each group.
\end{proposition}
 Until the end of the chapter, we adopt the convention that when objects are canonically identified along the pairing, we drop the subscript $i$, i.e., $\t$, $W_G$, $\fc_G$. In particular, since the discriminant divisors are the ramification divisors for the Chevalley maps, we have $\fD_G := \fD_{G_1} = \fD_{G_2}$.

Since our goal is to compare stable weighted orbital integrals between the two groups, we require compatibly chosen semistandard Levis. The isogeny $\psi$ induces a bijection from semistandard Levis in $\G_1$ to those in $\G_2$, and by twisting we get a corresponding bijection from semistandard $S$-Levis $M_1 \subset G_1$ to $M_2 \subset G_2$, characterized by the bijection of their geometric roots (resp.\ coroots) under $\psi^\star$ (resp.\ $\psi_\star$). Indeed, such Levis are paired as reductive groups in their own right. For any such pair, the Proposition above implies the resultant divisors coincide: $\fR_{M_1}^{G_1} = \fR_{M_2}^{G_2}$.
\begin{definition}\label{def.nonstd.datum}
    If $G_1, G_2$ over base $S$ are paired by isogeny $\psi$ and $M_i \subset G_i$ are semistandard $S$-Levis whose geometric roots biject under $\psi^\star$, then we call the quintuple $(M_1, G_1, M_2, G_2, \psi)$ a \emph{nonstandard endoscopic datum}.
\end{definition}

\begin{remark}[Compatible measures]\label{rmk.nonstd.measures}
   The isomorphism
$$
    \psi_\star : \ul{X}_\star(T_1)(\infty)_{\Q} \to \ul{X}_\star(T_2)(\infty)_{\Q}
$$
is naturally equivariant for $W_G(\infty) \rtimes \Gamma$, where $\Gamma$ is a prime-to-$p$ finite monodromy quotient splitting both $T_i$. Further extending scalars to $\R$ it gives isomorphisms between all relevant relative and geometric cocharacter spaces, and in particular
$$
    \psi_\star: \ga_{T_1} \cong \ga_{T_2}. 
$$
Recall that we have equipped these spaces with $W_G(\infty) \rtimes \Gamma$-equivariant measures, cf.\ Definition \ref{def.cochar.space.measure.choice}, and equipped all induced geometric and relative spaces, i.e., $\ga_{L_i}$, $\ga_{L_i}^{M_i}$, $\a_{L_i}$, $\a_{L_i}^{M_i}$, with the induced measures. Until the end of the chapter, assume that whenever paired groups appear in good characteristic, that the measures chosen are compatible with $\psi_\star$. The weight function $v_M^G$ (Definition \ref{def.arthur.wt.ftn}) depends upon this choice, as does $\vol(\a_M/\cX_\star(M))$ and Arthur's coefficient $d$ (Definition \ref{def.Arthur.d}). 
   
   In good characteristic, we choose Haar measures for corresponding regular centralizers $J_1, J_2$ which correspond in the sense of \cite[\S 1.12]{ngofl}: $\psi$ twists to give an isogeny between $J_1$ and $J_2$, in particular an isomorphism of their Lie algebras. The inverse function theorem makes this an isomorphism between sufficiently small compact open subgroups, which we can translate to equip $J_1(F_v)$ and $J_2(F_v)$ with compatible measures.

   In this chapter, all local or adelic orbital integrals and their stabilizations (Definitions \ref{def.s.orbital.integral.adelic}, \ref{def.local.s.integral}, and \ref{def.s.integral.V}) are defined relative to these compatible measure choices on $J_i(F_{0,v})$ and their induced products on $J_i(\A_0)$ and $J_i(\A_{0,V})$, rather than normalizations which equip the maximal compacts with volume 1 for both centralizers. The ambient groups retain their hyperspecial-volume-one normalizations.
\end{remark}

The remainder of the chapter proves the following conjecture of Waldspurger, the ``nonstandard WFL.''

\begin{notation}\label{not.commens.rel.index}
Fix the convention that for two commensurate lattices $A,B \subset V$ inside a rational or real vector space, $[A : B]$ denotes the relative index, i.e.,
$$
    [A : B] = \frac{[A : A \cap B]}{[B : A \cap B]}.
$$
\end{notation}

\begin{theorem}[{The nonstandard WFL, \cite[Conjecture 3.7]{WaldTordu}}]\label{thm.nonstd.WFL}
    Let $F_{0,v}$ be an equal-characteristic local field with ring of integers $\O_{0,v}$. Let $(M_1, G_1, M_2, G_2, \psi)$ be a nonstandard endoscopic datum over $\Spec \O_{0,v}$, assumed in good characteristic. Let $$b \in \fc_{M_1}(\O_{0,v}) \cap \fc_{M_1}^{G_1\text{-reg}}(F_{0,v}) = \fc_{M_2}(\O_{0,v}) \cap \fc_{M_2}^{G_2\text{-reg}}(F_{0,v})$$ be a generically regular parameter, and write
    $$
        Y_{1,b} = \epsilon_{M_1}(b), \qquad Y_{2,b} = \epsilon_{M_2}(b)
    $$
    for the Kostant lifts. Then the local stable orbital integrals (defined as in Definition \ref{def.stable.integral.local} with compatible choices of measure on $J_i(F_v)$ as above) satisfy
    $$
        S_{M_1, v}^{G_1}(Y_{1,b}) = c_{M_1, M_2}^{G_1, G_2} S_{M_2, v}^{G_2}(Y_{2,b}),
    $$
    where 
    $$
        c_{M_1, M_2}^{G_1, G_2} = \frac{[(\Z \langle \Phi_{G_2}^\vee \rangle)^{\Gamma} : \psi_\star((\Z \langle \Phi_{G_1}^\vee \rangle)^{\Gamma})]}{[(\Z \langle \Phi_{M_2}^\vee \rangle)^{\Gamma} : \psi_\star((\Z \langle \Phi_{M_1}^\vee \rangle)^{\Gamma})]}
    $$
    is the coefficient of Waldspurger.
\end{theorem}

\section{The global nonstandard WFL.}\label{subsect.nonstd.global}

We begin in the global setting described at the beginning of \S \ref{sect.global.point.counts}. Let $(M_1, G_1, M_2, G_2, \psi)$ be a nonstandard endoscopic datum over a curve $X_0/k_0 =\F_q$ split at $\infty$, in good characteristic.

\begin{proposition}[{\cite[\S 4.18, Proposition 4.18.1]{ngofl}}]\label{prop.nonstd.same.chevalley.lie.torus}
    If $G_1, G_2/X$ are paired and $D$ is any divisor, then the isomorphism $\nu : \fc_{G_1} \to \fc_{G_2}$ induces by twisting an isomorphism
    $$
        \fc_{G_1,D} \xrightarrow{\sim} \fc_{G_2,D}
    $$
    and on global sections an identification of the Hitchin bases
    $$
        \cA_G := \cA_{G_1} = \cA_{G_2}.
    $$
    Moreover, we have isogeneous identity components $J_{G_1}^0$ and $J_{G_2}^0$; this induces a homomorphism of the Picard stacks
    $$
        \cP_{G_1} \longrightarrow \cP_{G_2}
    $$
    over this common base which is an isogeny on the neutral components. 
\end{proposition}

\begin{remark}
    In the proof of \cite[Proposition 4.18.1]{ngofl}, the author asserts an isomorphism $J_1^0 \cong J_2^0$. This is not true in general. For example, when $G_1 = \SL_3$ and $G_2 = \PGL_3$ in good characteristic, and a parameter is regular semisimple with spectral curve corresponding on generic points to a non-Galois cubic extension $E/F$, the generic regular centralizers are
    $$
        J_1 = \ker(N_{E/F}), \qquad J_2 = \Res_{E/F} \G_m/\G_m,
    $$
    connected and isogeneous tori, but not isomorphic. The isogeny is sufficient to construct the homomorphism $\cP_{G_1} \to \cP_{G_2}$, so no downstream consequences are affected.
\end{remark}

Identifying the split tori with the fiber at infinity and thereby $\Lie(\T_i)$ with $\ft_{i,\infty}$, we have canonical identifications of the pullbacks $\ctA_{G_1} = \ctA_{G_2}$. The poset of semistandard Levis and their various Chevalley maps likewise biject, so we deduce an identification of the compatible loci $\bA_{G_1} = \bA_{G_2}$ of Definition \ref{def.comp.loc} such that $\bA_{M_1\el} = \bA_{M_2\el}$ for each corresponding pair of Levis. We abusively write without a subscript these common loci, e.g. $\cA_{G\el}, \bA_{M\el}$. The homomorphism of Picard stacks defined above pulls back to $\bP_{G_1} \to \bP_{G_2}$ over this base. The monodromy homomorphisms $\mu_{a_1}$ and $\mu_{a_2}$ of Definition \ref{def.monodr.inv.Out} for corresponding Hitchin base points induce actions on $\ul{X}_\star(T_i)(\infty)_{\Q}$ which correspond under $\psi_\star$. In particular, there is canonically a common monodromy group $W_a = W_{a_1} = W_{a_2}$.

\begin{remark}\label{rmk.point.count.indep.of.meas.choice}
    Since our choice of $\psi$-compatible measures breaks some of our previous normalization conventions, at this point we should remark on the dependence of our orbital integrals and prior theorems on the choices of measure. If the measure on $\ga_T$ is inflated by a factor of $\lambda$, then our adelic orbital integrals $J_M^{G,s}$ will inflate by the same factor. 
    The same scalar factor is then inherited by $S_M^{G}$. On the other hand, the term $\vol(\a_{M_i}/\cX_\star(M))$ appearing in the denominator of Theorems \ref{thm.s.summand.orbital.integral} and \ref{thm.stable.point.count} will likewise scale by $\lambda$. All other terms are unaffected by the choice. On the other hand, the measure on $J(\A_0)$ is inflated by a factor $\lambda$, then the induced measure on $(J_{a_M}(\A) \backslash G(\A))^{\tau}$ and the value of the orbital integral $J_{M}^{G,s}$ will scale by a factor of $\lambda^{-1}$, and $S_{M}^{G}$ likewise. The factor $\vol(J(F_0) \backslash J(\A_0)^1)$, and thus $c(J)$, scales by a corresponding factor of $\lambda$. It is now clear that Theorems \ref{thm.s.summand.orbital.integral} and \ref{thm.stable.point.count} are fully independent of either choice, and we can freely assume that the theorems hold with measures chosen compatible with $\psi$ in the sense of Remark \ref{rmk.nonstd.measures}. 
\end{remark} 

It is easy to see that the degree maps and $\xi$-truncation correspond across the pairing $\psi$, and one expects the $s$-decompositions to be closely related. In the spirit of perverse continuation, we will never need any such identifications. Instead, the trace of Frobenius computation of Theorem \ref{thm.stable.point.count} will allow us to compare the stable orbital integrals without knowing any identification beyond the elliptic locus, \cite[Th\'eor\`eme 8.8.2]{ngofl}.

First, we will require the following fundamental lattice-theoretic computation.

\begin{lemma}\label{lem.covol.coroots}
    Let $H \supseteq T$ be any quasisplit reductive $X$-group, with $\Gamma$ a finite monodromy group acting on its geometric based root datum $(X^\star(\T), \ul{\Phi}_H, X_\star(\T), \ul{\Phi}_H^\vee, \ul{\Delta}_H)$. Then, with respect to any compatible family of measures on the various cocharacter spaces $\a_{L}^M$ as in Definition \ref{def.cochar.space.measure.choice}, we have
    $$
        \vol(\a_{T}^H/\Z \langle \ul{\Phi}^\vee_H\rangle_\Gamma ) = \frac{|\left(X_\star(\T)/\Z\langle \ul{\Phi}^\vee_{H} \rangle\right)_{\Gamma, \tors}|}{|X_\star(\T)_{\Gamma, \tors}|}\times \frac{\vol(\a_T/\cX_\star(T))}{\vol(\a_H/\cX_\star(H))}.
    $$
    \begin{proof}
        On the other hand, the torsion-free part $(X_\star(\T)/\Z\langle \ul{\Phi}^\vee_{H} \rangle)_{\Gamma, \tf}$ is our global dual lattice $\cX_\star(H)$ as in Proposition \ref{prop.dual.lattice.global.invs.coinvs.mod.tors}:
        $$
            \Hom((X_\star(\T)/\Z\langle \ul{\Phi}^\vee_{H} \rangle)_{\Gamma, \tf}, \Z) = \Hom(X_\star(\T)/\Z\langle \ul{\Phi}^\vee_{H} \rangle, \Z)^{\Gamma} = X^\star(H).
        $$
        
        We now have the lattices of coroots for each $H$; our goal is to compare these lattices to the standard lattice dual to the global characters $\cX_{\star}(T) = X_\star(\T)_{\Gamma, \tf} \subset \a_{T}$. Consider the embedding
        $$
            \Z\langle \ul{\Phi}^\vee_{H} \rangle \hookrightarrow X_\star(\T).
        $$
        The action of $\Gamma$ is compatible with our fixed pinning, so permutes the simple coroots. In particular, the coroot lattice is a permutation lattice and its coinvariants are torsion-free, $\Z \langle \ul{\Phi}^\vee_{H} \rangle_{\Gamma, \tf} = \Z \langle \ul{\Phi}^\vee_{H} \rangle_{\Gamma}$. We deduce a short exact sequence on coinvariants 
        $$
            \begin{tikzcd}[ampersand replacement=\&]
            0 \ar[r] \& \Z\langle \ul{\Phi}^\vee_{H} \rangle_{\Gamma, \tf} \ar[r] \& X_\star(\T)_{\Gamma} \ar[r] \& \left(X_\star(\T)/\Z\langle \ul{\Phi}^\vee_{H} \rangle\right)_{\Gamma} \ar[r] \& 0
            \end{tikzcd}.
        $$
        Note the left-exactness -- the kernel of $\Z \langle \ul{\Phi}^\vee_{H} \rangle_{\Gamma, \tf} \to X_\star(\T)_{\Gamma}$ would be torsion, but the domain is torsion-free. Passing to the quotient by torsion, the sequence 
        $$
            \begin{tikzcd}[ampersand replacement=\&]
             0 \ar[r] \& \Z\langle \ul{\Phi}^\vee_{H} \rangle_{\Gamma, \tf} \ar[r] \& \cX_\star(T) \ar[r] \& \cX_\star(H) \ar[r] \& 0
            \end{tikzcd}
        $$
        need not be exact in the middle because the inclusion of $\Z\langle \ul{\Phi}^\vee_{H} \rangle_{\Gamma, \tf}$ may not be saturated, but writing $S(H) \supseteq \Z\langle \ul{\Phi}^\vee_{H} \rangle_{\Gamma, \tf}$ for its saturation inside the torsion-free quotient $\cX_\star(T)$, then we have a short exact sequence
        $$
            \begin{tikzcd} 0 \ar[r] & S(H) \ar[r] & \cX_\star(T) \ar[r] & \cX_\star(H) \ar[r] & 0
            \end{tikzcd}. 
        $$
        This decomposition occurs across the ambient orthogonal decomposition $\a_T = \a_T^H \oplus \a_H$, so the covolume with respect to the induced measures is multiplicative:
        $$
            \vol(\a_T/\cX_\star(T)) = \vol(\a_T^H/S(H)) \vol(\a_H/\cX_\star(H)).
        $$
        It remains to compute the index of $$\Z\langle \ul{\Phi}^\vee_{H} \rangle_{\Gamma, \tf} = \Z\langle \ul{\Phi}^\vee_{H} \rangle_{\Gamma} \subseteq S(H).$$
        Consider the preimage $S'(H)$ of $S(H)$ under the quotient by torsion $q : X_\star(\T)_{\Gamma} \to X_\star(\T)_{\Gamma, \tf} = \cX_\star(T)$. The third isomorphism theorem gives
        $$
            X_\star(\T)_{\Gamma}/S'(H) \cong X_\star(\T)_{\Gamma, \tf}/S(H),
        $$
        so a simple application of the snake lemma identifies $\ker(S'(H) \to S(H)) \cong X_\star(\T)_{\Gamma, \tors}$. Thus, 
        \begin{equation}
            [S(H) : \Z\langle \ul{\Phi}^\vee_{H} \rangle_{\Gamma, \tf}] \cong \frac{[S'(H) : \Z\langle \ul{\Phi}^\vee_{H} \rangle_{\Gamma, \tf}]}{|X_\star(\T)_{\Gamma, \tors}|}.\label{eqn.saturation.index.form.partial}
        \end{equation}
        But on the other hand, the isomorphism 
        $$
            X_\star(\T)_{\Gamma}/S'(H) \cong \cX_\star(H) = \left(X_\star(\T)/\Z\langle \ul{\Phi}^\vee_{H} \rangle\right)_{\Gamma, \tf}
        $$
        implies by another simple snake lemma the short exact sequence
        $$
            \begin{tikzcd}
                0 \ar[r] & \Z\langle \ul{\Phi}^\vee_{H} \rangle_{\Gamma} \ar[r] & S'(H) \ar[r] & \left(X_\star(\T)/\Z\langle \ul{\Phi}^\vee_{H} \rangle\right)_{\Gamma, \tors} \ar[r] & 0
            \end{tikzcd}.
        $$
        Substituting the resulting index formula into (\ref{eqn.saturation.index.form.partial}) gives the Lemma.
    \end{proof}
\end{lemma}

\begin{corollary}[The global nonstandard WFL]\label{cor.global.nonstd.WFL}
    Let $a_M \in \bA_{M_1\el}(k_0) = \bA_{M_2\el}(k_0)$ a parameter with image $a \in \bA_G^{\good}$. Then, for measures chosen compatible with $\psi$ as above, we have an equality of adelic stable orbital integrals 
    $$
        S_{M_1}^{G_1}(a_M) = c_{M_1, M_2}^{G_1, G_2} S_{M_2}^{G_2}(a_M),
    $$
    where $c_{M_1, M_2}^{G_1, G_2}$ is the coefficient of Waldspurger.
    \begin{proof}
        Assume first that $G_1$ has anisotropic center. Since $\a_{G_1} = \a_{G_2}$, this is equivalent to $G_2$ having anisotropic center. Thus both are admissible parameters for Theorem \ref{thm.stable.point.count}. 
        We deduce equations for $i = 1,2$:
        \begin{equation}
            \Tr(\tau^{-1}, ((j_{ell})_{!\star}(\pH^\bullet(Rf_{G_i,\star}^{G_i\el} \ol{\Q}_l)_1))_a) = \frac{c(J_{i, a_{M_i}}) \cdot q^{\deg(D) |\Phi^{G_i}|/2} |\pi_0(Z^\Gamma_{\widehat{M}_i})|}{\vol(\a_{M_i}/\cX_\star(M_i)) \cdot |\pi_0(\widehat{T}_i^{W_a})| \cdot |\pi_0(Z^\Gamma_{\widehat{G}_i})|} S_{M_i}^{G_i}(a_M),\label{eqn.nonstd.global.stable.trace}
        \end{equation}
        where $J_{i, a_{M_i}}$ is the $M_i$-regular centralizer of Ng\^{o} and $c(J) = \vol(J(F_0) \backslash J(\A_0)^1)\cdot  |\ker(J(F)_{\tau} \to J(\A)_{\tau})|$.
        
        By Ng\^{o}'s cohomological nonstandard Fundamental Lemma \cite[Th\'eor\`eme 8.8.2]{ngofl}, we have an isomorphism (over $k_0$) of the fine stable components up to semisimplification on the common elliptic locus $\cA_{G\el}$:
        $$
            \pHi(Rf_{G_1,\star}^{G_1\el} \ol{\Q}_l)_1^{\ss} \cong \pHi(Rf_{G_2,\star}^{G_2\el} \ol{\Q}_l)^{\ss}_1.
        $$
        In general, intermediate extension along an open may not be exact on perverse sheaves, so need not commute with semisimplification. In general, one has a decomposition $(j_{!\star} K)^{\ss} = j_{!\star}(K^{\ss}) \oplus B$, where $B$ is supported on the boundary. But we know by Corollary \ref{cor.coh.intermed.ext.from.ell.locus} that for both groups every semisimple summand has support which intersects the elliptic locus, so no new summands can appear and we have 
        $$
            (j_{ell})_{!\star} \left(\pHi(Rf_{G_i,\star}^{G_i\el} \ol{\Q}_l)_1^{\ss}\right) = \left((j_{ell})_{!\star} \pHi(Rf_{G_i,\star}^{G_i\el} \ol{\Q}_l)_1\right)^{\ss}.
        $$
        Taking the stalk at $a$ and trace of Frobenius, which is insensitive to the semisimplification, identifies the left sides of (\ref{eqn.nonstd.global.stable.trace}), so their right sides are equal. 
        
        The factor of $q$ cancels immediately. Since our parameter $a_M$ is generically $M$-regular semisimple the regular centralizer is generically a torus, and the $F_0$-tori $J_i = J_{i, a_M}$ are isogenous with compatible as in Remark \ref{rmk.nonstd.measures}. We claim that 
        \begin{equation}
            \frac{c(J_1)}{|\pi_0(\widehat{T}_1^{W_a})|} = \frac{c(J_2)}{|\pi_0(\widehat{T}_2^{W_a})|} \label{eqn.cpi0.scaling}
        \end{equation}
        by an argument passing through the Tamagawa measures on these tori.
        
        First, choosing corresponding nonzero invariant top forms $\omega_i$, our choice of compatible local measures on $J_i(F_{0,v})$ satisfy
        $$
            dj_{i,v} = b_v |\omega_i|_v
        $$
        for the same positive scalar $b_v$. We equip $J_i(F_0)$ with the counting measure and  the adelic points with the associated product. The isomorphism $X^\star(T_1)_{\Q} \cong X^\star(T_2)_{\Q}$ induces equality on their local  local and global Artin $L$-factors, so the induced compatible measures $dj_{i}$ on the adelic quotient differ from their Tamagawa measures by the same scalar $B > 0$.
        
        By \cite[Lemma 12.2.5]{CLII}, the kernel term in the definition (\ref{eqn.def.of.cJ}) of $c(J)$ is equal to $\Sha(J_i) := \ker\left(H^1(F_0, J_i) \to \prod_v H^1(F_{0,v}, J_i)\right)$, i.e.,
        $$
            c(J_i) = \vol_{dj_i}(J(F_0)\backslash J(\A_0)^1) |\Sha(J_i)| = B \tau(J_i) |\Sha(J_i)|.
        $$
        The Ono--Oesterl\'e formula (\cite{OnoTamagawa}) gives
        $$
            \tau(J_i) = \frac{|H^1(F_0, X^\star(J_i))|}{|\Sha(J_i)|} = \frac{|\pi_0(\widehat{J}_i^{\Gamma_{F_0}})|}{|\Sha(J_i)|}
        $$
        and the rational cameral marking at $\infty$ identifies
        $$
            |\pi_0(\widehat{J}_i^{\Gamma_{F_0}})| = |\pi_0(\widehat{T}_i^{W_a})|
        $$
        by taking the dual of 
        $$
            X_\star(J_{i})_{\pi_1(U_a, \infty), \tors} \cong X_\star(T_i)_{W_a}
        $$
        and identifying the image of arithmetic and geometric monodromy, 
        as in the discussion after Kottwitz's lemma in \S \ref{subsect.weighted.s.orbital.global}. So (\ref{eqn.cpi0.scaling}) holds.

        It remains to identify Waldspurger's constant 
        $c_{M_1, M_2}^{G_1, G_2}$ with
        \begin{equation}
            \frac{|\pi_0(Z^\Gamma_{\widehat{M}_2})|}{|\pi_0(Z^\Gamma_{\widehat{M}_1})|} \times \frac{|\pi_0(Z^\Gamma_{\widehat{G}_1})|}{|\pi_0(Z^\Gamma_{\widehat{G}_2})|} \times \frac{\vol(\a_{M_1}/\cX_\star(M_1))}{\vol_{\psi}(\a_{M_2}/\cX_\star(M_2))},\label{eqn.lattice.constant.Walds}
        \end{equation}
        where we have written $\psi$ to emphasize the compatible choice of measures.
        
        Let us translate the statement of Lemma \ref{lem.covol.coroots} to the dual side. We have for each $H_i = T_i, M_i, G_i$
        $$
            X^\star(Z_{\widehat{H}_i}) = X_\star(\T_i)/\Z\langle \ul{\Phi}^\vee_{H_i} \rangle,
        $$
        since $Z_{\widehat{H}_i} \subseteq \widehat{T}_i$ is cut out by the roots of $\widehat{H}_i$ (coroots of $H_i$). Taking invariants on the dual side corresponds to $\Gamma$-coinvariants, and taking the component group corresponds to the torsion part, so
        $$
            |\pi_0(Z_{\widehat{H}_i}^{\Gamma})| = \left|(X_\star(\T_i)/\Z\langle \ul{\Phi}^\vee_{H_i} \rangle)_{\Gamma, \tors}\right|.
        $$
        In particular, we could specialize $H_i = T_i$ to get 
        $$
            |\pi_0(Z_{\widehat{T}_i}^{\Gamma})| = \left|X_\star(\T_i)_{\Gamma, \tors}\right|.
        $$
        So Lemma \ref{lem.covol.coroots} becomes the statement
        $$
            \vol(\a_{T_i}^{H_i}/\Z \langle \ul{\Phi}_{H_i}^\vee \rangle_\Gamma) = \frac{|\pi_0(Z_{\widehat{H}_i}^{\Gamma})|}{|\pi_0(Z_{\widehat{T}_i}^{\Gamma})|} \times \frac{\vol(\a_{T_i}/\cX_\star(T_i))}{\vol(\a_{H_i}/\cX_\star(H_i))}
        $$
        where as above we note that $\Z \langle \ul{\Phi}_{H_i}^\vee \rangle_\Gamma$ is torsion-free to identify itself with its image in the torsion-free quotient space.

        We apply this Lemma four times, to each of the four groups in our datum, using the $\psi_\star$-induced measures on spaces for $G_2, M_2$. When applied to the $G_i$, our anisotropic center hypothesis makes the covolume term in $\a_{G_i} = 0$ identically 1. All terms involving just the tori cancel, and the expression (\ref{eqn.lattice.constant.Walds}) becomes
        $$
            \frac{\vol_{\psi}(\a_{T_2}^{M_2}/ \Z\langle \ul{\Phi}_{M_2}^\vee \rangle_\Gamma)}{\vol(\a_{T_1}^{M_1}/ \Z\langle \ul{\Phi}_{M_1}^\vee \rangle_\Gamma)} \times \frac{\vol(\a_{T_1}^{G_1}/ \Z\langle \ul{\Phi}_{G_1}^\vee \rangle_\Gamma)}{\vol_{\psi}(\a_{T_2}^{G_2}/ \Z\langle \ul{\Phi}_{G_2}^\vee \rangle_\Gamma)}.
        $$
        By definition, this can be written in terms of the relative indices
        $$
            \frac{[\Z \langle \ul{\Phi}_{G_2}^\vee \rangle_{\Gamma} : \psi_\star(\Z\langle \ul{\Phi}_{G_1}^\vee \rangle_{\Gamma})]}{[\Z \langle \ul{\Phi}_{M_2}^\vee \rangle_{\Gamma} : \psi_\star(\Z\langle \ul{\Phi}_{M_1}^\vee \rangle_{\Gamma})]}.
        $$
        It remains only to see that these coinvariants can be replaced with invariants. This follows from the $\Gamma$-equivariance of the map $\psi_\star$: a basis for the invariant space $\Z \langle \ul{\Phi}_{H}^\vee \rangle^{\Gamma}$ is given by sums over $\Gamma$-orbits of simple coroots $\sum_{\alpha^{\vee}} \alpha^\vee$ for $H$, while the corresponding basis element in $\Z \langle \ul{\Phi}_{H}^\vee \rangle_{\Gamma}$ is the class of any $\alpha^\vee$. But whatever scalar $\psi_\star$ acts by on $\alpha^\vee_1$
        $$
            \psi_\star(\alpha_1^\vee) = b(\alpha_1^\vee) \alpha_2^\vee
        $$
        must by $\Gamma$-equivariance be constant on orbits, so $\psi_\star$ acts by the same diagonal matrix on the invariant and coinvariant spaces. We conclude that (\ref{eqn.lattice.constant.Walds}) is equal to    
        $$
            \frac{[\Z \langle \ul{\Phi}_{G_2}^\vee \rangle^{\Gamma} : \psi_\star(\Z\langle \ul{\Phi}_{G_1}^\vee \rangle^{\Gamma})]}{[\Z \langle \ul{\Phi}_{M_2}^\vee \rangle^{\Gamma} : \psi_\star(\Z\langle \ul{\Phi}_{M_1}^\vee \rangle^{\Gamma})]} = c_{M_1, M_2}^{G_1, G_2}
        $$
        and the Corollary holds.

        Now let $G_i$ general paired quasisplit reductive groups. As in the proof of the global endoscopic WFL, Theorem \ref{thm.global.wfl}, passing to the anisotropic quotients $\ol{G}_i = G_i/A_{G_i}^{\spl}$, $\ol{M}_i = M_i/A_{G_i}^{\spl}$ does not affect ordinary weighted $s$-orbital integrals. The pairing $\psi$ already identified the central cocharacter spaces
        $$
            Z_i = \{x \in X_\star(\T_i)_{\Q} : \langle \alpha, x \rangle = 0 \;\forall \alpha\}
        $$
        and was $\Gamma$-equivariant, so identifies the subspaces $Z_i^\Gamma$ corresponding to the $A_{G_i}^{\spl}$. Quotienting by those subspaces deduces a pairing $\ol{\psi}$ from $\ol{G}_1$ to $\ol{G}_2$. 
         
        All endoscopic groups contain the split center, so descend, and in fact we have a bijection for any ambient $H \supset A_{G_i}^\spl$
        $$
            \cE_{M}(H) \xrightarrow{\sim} \cE_{\ol{M}}(\ol{H}),
        $$
        by the following. For any Levi $L \subset H$, one has a dual sequence of $\Gamma$-fixed centers
        $$
            \begin{tikzcd}
                1 \ar[r] & Z_{\widehat{\ol{L}}}^\Gamma \ar[r] & Z_{\widehat{L}}^\Gamma \ar[r] & \widehat{A_{G_i}^\spl}^\Gamma \ar[r] & 1
            \end{tikzcd};
        $$
        the right-exactness follows from the splitness $A_{G_i}^\spl$. It follows that we have an isomorphism
        $$
            Z_{\widehat{M}}^\Gamma/Z_{\widehat{H}}^\Gamma \xrightarrow{\sim} Z_{\widehat{\ol{M}}}^\Gamma/Z_{\widehat{\ol{H}}}^\Gamma
        $$
        Ellipticity of an elliptic parameter is easily seen to be preserved since the dual split central torus is a common quotient, so this gives the bijection of endoscopic parameters. The same bijection shows that the scalars
        $$
            |Z^\Gamma_{\widehat{M}}/Z^\Gamma_{\widehat{H}}|
        $$
        remain unchanged upon passage to the quotient, so the entire recursion defining the stable orbital integrals passes unchanged. 

        Finally, quotienting by the $Z_i^\Gamma$ induces a canonical isomorphism between the geometric coroot lattices for $G_i$ and $\ol{G}_i$ (resp.\ $M_i$ and $\ol{M}_i$), and $\ol{\psi}$ acts by the same scalars. So both terms in Waldspurger's coefficient are unaffected. We conclude that the statement of the nonstandard WFL descends term-by-term to the statement for $(\ol{M}_1, \ol{G}_1, \ol{M}_2, \ol{G}_2, \ol{\psi})$, already proven.

                                                                            \end{proof}
\end{corollary}

\section{Splitting formalism for the nonstandard WFL.}\label{subsect.nonstd.splitting}

Retain the global setting of the previous section, and fix a parameter $a_M \in \bA_{M_1\el}(k_0) = \bA_{M_2\el}(k_0)$. Now we proceed to deduce the local nonstandard WFL from the global result, Corollary \ref{cor.global.nonstd.WFL}, above. The strategy is essentially the same as the previous chapter. Recall that we have defined (Definition \ref{def.s.integral.V}) partial stable orbital integrals $S_{M,V}^G(a_M)$ for any $V \subset |X_0|$, and that the stable orbital integrals satisfy a splitting formula, (\ref{eqn.splitting.formula.stable}) from Lemma \ref{lem.splitting.lemma}: for any $V = V_1 \bigsqcup V_2 \subseteq |X_0|$,
$$
    S_{M,V}^{G}(a_M) = \sum_{L_1, L_2 \supseteq M} e_M^G(L_1, L_2) S_{M,V_1}^{L_1}(a_M) S_{M,V_2}^{L_2}(a_M),
$$
where
$$e_M^G(L_1,L_2) := d_M^G(L_1, L_2) \left|\frac{Z_{\widehat{L_1}}^\Gamma \cap Z_{\widehat{L_2}}^\Gamma}{Z_{\widehat{G}}^\Gamma} \right|^{-1}.$$

\begin{definition}\label{def.nonstd.WFL.V}
    Given nonstandard endoscopic datum $(M_1, G_1, M_2, G_2, \psi)$ on $X_0$ and a set of places $V \subseteq |X_0|$, say that \emph{the nonstandard WFL holds on $V$} if for every corresponding pair of intermediate Levis $L_1 \longleftrightarrow L_2$,
    $$
        S_{M_1,V}^{L_1}(a_M) = c_{M_1, M_2}^{L_1, L_2} S_{M_2, V}^{L_2}(a_M).
    $$
    Say that $V$ is \emph{peelable} if the nonstandard WFL holds on $V$ and if in addition
    $$
        S_{M_1,V}^{M_1}(a_M) = S_{M_2,V}^{M_2}(a_M) \neq 0.
    $$
\end{definition}

\begin{remark}\label{rmk.nonstd.local.global.discrep}
    As in the endoscopic case, this does not necessarily reproduce the local statement of the nonstandard WFL when $V = \{v\}$, because the weights appearing in the orbital integrals are global and the endoscopic recursion which defines the stable orbital integrals is determined by global endoscopic parameters. However, as observed in Remark \ref{rmk.monodr.local.global.discrep.wfl}, these discrepancies vanish up to a common scalar if we engineer $D_v = \Gamma$ at a given place $v$, and vanish entirely if $k_v = k_0$. Our globalization will ensure these, so when we prove that the nonstandard WFL holds on $V = \{v\}$, this will imply Theorem \ref{thm.nonstd.WFL}.
\end{remark}

We have a peeling lemma.
\begin{lemma}[nonstandard Peeling Lemma]\label{lem.peeling.nonstd}
    With notations as above, let $V = V_1 \bigsqcup V_2 \subseteq |X_0|$ be a nontrivial decomposition of sets of places. Suppose that the nonstandard WFL holds on $V$ and that $V_2$ is peelable. Then the nonstandard WFL holds on $V_1$.
    \begin{proof}
        The induction proceeds exactly as in Lemma \ref{lem.peeling}, where, as we explain below, Lemma \ref{lem.e.compat.c} is used to identify the coefficients that appear in the inductive step. The base case $L_i = M_i$ decomposes by the product formula on each side
        $$
            S_{M_1,V_1}^{M_1}(a_M) S_{M_1,V_2}^{M_1}(a_M) = S_{M_1,V}^{M_1}(a_M) = c_{M_1, M_2}^{M_1, M_2} S_{M_2,V}^{M_2}(a_M) = c_{M_1, M_2}^{M_1, M_2}S_{M_2,V_1}^{M_2}(a_M)S_{M_2,V_2}^{M_2}(a_M)
        $$
        and the nonvanishing implies the identity on $V_1$. For the inductive step, fix $M_i \subsetneq L_i \subseteq G_i$. Applying the splitting formulae to the nonstandard WFL on $V$ gives
        $$
            \sum_{A_1,B_1 \supseteq M_1} e_{M_1}^{L_1}(A_1, B_1) S_{M_1, V_1}^{A_1}(a_M) S_{M_1, V_2}^{B_1}(a_M) = c_{M_1, M_2}^{L_1, L_2} \sum_{A_2, B_2 \supseteq M_2} e_{M_2}^{L_2}(A_2, B_2) S_{M_2, V_1}^{A_2}(a_M) S_{M_2, V_2}^{B_2}(a_M),
        $$
        where we can bijectively pair intermediate Levis $A_1 \longleftrightarrow A_2, B_1 \longleftrightarrow B_2$. Using the nonstandard WFL on $V_2$ then gives
        $$
            \sum_{A_1,B_1 \supseteq M_1} e_{M_1}^{L_1}(A_1, B_1) c_{M_1, M_2}^{B_1, B_2} S_{M_1, V_1}^{A_1}(a_M) S_{M_2, V_2}^{B_2}(a_M) = \sum_{A_2, B_2 \supseteq M_2} c_{M_1, M_2}^{L_1, L_2} e_{M_2}^{L_2}(A_2, B_2) S_{M_2, V_1}^{A_2}(a_M) S_{M_2, V_2}^{B_2}(a_M).
        $$
        For all terms except $(A,B) = (L,M)$, the inductive hypothesis together with the identification of constants in Lemma \ref{lem.e.compat.c} give equality. For that term, the $S_{M_i, V_2}^{M_i}(a_M)$ are nonzero and equal, so cancel; the remaining constants equate by Lemma \ref{lem.e.compat.c} (and are nonzero), cancelling to give the desired
        $$
            S_{M_1, V_1}^{L_1}(a_M) = c_{M_1, M_2}^{L_1, L_2} S_{M_2, V_1}^{L_2}(a_M).
        $$
    \end{proof}
\end{lemma}

\begin{lemma}\label{lem.e.compat.c}
    For any diamond of paired semistandard Levis
    $$
        \begin{tikzcd}
            & L_1 & \; \ar[rr, leftrightarrow] & & \; & L_2 & \\
            A_1 \ar[ur, hook] & & B_1 \ar[ul, hook] \ar[rr, leftrightarrow, "\psi_\star"] & & A_2\ar[ur, hook] & & B_2\ar[ul, hook] \\
            & M_1\ar[ur, hook] \ar[ul, hook] & \;\ar[rr, leftrightarrow] & & \; & M_2\ar[ur, hook]\ar[ul, hook] &
        \end{tikzcd}
    $$
    we have an equality of constants $$
            e_{M_1}^{L_1}(A_1, B_1) c_{M_1, M_2}^{A_1, A_2} c_{M_1, M_2}^{B_1, B_2}  = c_{M_1, M_2}^{L_1, L_2} e_{M_2}^{L_2}(A_2, B_2)
        $$
    \begin{proof}
        The identity is vacuous if $d_{M_1}^{L_1}(A_1, B_1) = d_{M_2}^{L_2}(A_2, B_2)$ vanishes, so assume otherwise. The essential lattice-theoretic identity is that for either index
        \begin{equation}
            \left|\frac{Z_{\widehat{A}}^\Gamma \cap Z_{\widehat{B}}^\Gamma}{Z_{\widehat{L}}^\Gamma} \right| = d_M^L(A,B) \frac{|\pi_0(Z^\Gamma_{\widehat{M}})| \cdot |\pi_0(Z^\Gamma_{\widehat{L}})|}{|\pi_0(Z^\Gamma_{\widehat{A}})| \cdot |\pi_0(Z^\Gamma_{\widehat{B}})|} \frac{\vol(\a_{A}/\cX_\star(A)) \vol(\a_{B}/\cX_\star(B))}{\vol(\a_{M}/\cX_\star(M)) \vol(\a_{L}/\cX_\star(L))}.\label{eqn.lattice.comp.ABML}
        \end{equation}
        To help organize this, notate for any diagonalizable group $C$ over $\C = \ol{\Q}_l$ the total index
        $$
            \nu(C) := \frac{\on{covol}(X^\star(C)_{\tf})}{|X^\star(C)_{\tors}|},
        $$
        the covolume taken with respect to some established choice of measure on $X^\star(C)_{\R}$. When $C = Z_{\widehat{H}}^\Gamma$, we obtain by the same computations as in Corollary \ref{cor.global.nonstd.WFL}
        $$
            \nu(Z_{\widehat{H}}^\Gamma) = \frac{\vol(\a_H/\cX_\star(H))}{|\pi_0(Z_{\widehat{H}}^\Gamma)|}.
        $$
        The proof will now compute $\nu(Z_{\widehat{A}}^\Gamma \cap Z_{\widehat{B}}^\Gamma)$ two ways, once in terms of $A$ and $B$ and once in terms of $M$ and $L$. Both require the following covolume computation.

        \begin{lemma}\label{lem.index.mult.ses.fg.gps}
            Suppose 
            $$
                \begin{tikzcd}
                    0 \ar[r] & Y_0 \ar[r] & Y_1 \ar[r] & Y_2 \ar[r] & 0
                \end{tikzcd}
            $$
            is an exact sequence of finitely generated abelian groups. Tensoring to $\R$ gives a short exact sequence of real vector spaces
            $$
                \begin{tikzcd}
                    0 \ar[r] & V_0 \ar[r] & V_1 \ar[r] & V_2 \ar[r] & 0
                \end{tikzcd},
            $$
            and assume we have equipped these spaces with measures. Let $\Delta$ denote the absolute determinant of the short exact sequence, i.e. $\Delta = |\det(i \oplus s : V_0 \oplus V_2 \to V_1)|$ for any choice of splitting $s: V_2 \to V_1$. Then
            \begin{equation}
                \frac{\vol(V_1/(Y_1)_{\tf})}{|(Y_1)_{\tors}|} = \Delta \frac{\vol(V_0/(Y_0)_{\tf})}{|(Y_0)_{\tors}|} \frac{\vol(V_2/(Y_2)_{\tf})}{|(Y_2)_{\tors}|}.\label{eqn.index.mult.ses.fg.gps}
            \end{equation}
            \begin{proof}
                Choose bases of $(Y_0)_{\tf}, (Y_2)_{\tf}$, and lift the latter to $V_1$. They generate a sub-lattice $\Lambda \subseteq (Y_1)_{\tf}$. A snake-lemma computation similar to the one in the proof of Lemma \ref{lem.covol.coroots} gives
                $$
                     [(Y_1)_{\tf} : \Lambda] = \frac{|(Y_0)_{\tors}| \cdot |(Y_2)_{\tors}|}{|(Y_1)_{\tors}|}.
                $$
                Details left to the reader.
            \end{proof}
        \end{lemma}

        First consider the short-exact sequence of fixed centers
        $$
            \begin{tikzcd}
                1 \ar[r] & Z_{\widehat{A}}^\Gamma \cap Z_{\widehat{B}}^\Gamma \ar[r, "{(1, -1)}"] & Z_{\widehat{A}}^\Gamma \times Z_{\widehat{B}}^\Gamma \ar[r] & Z_{\widehat{M}}^\Gamma \ar[r] & 1
            \end{tikzcd}.
        $$
        Surjectivity follows from \cite[Lemma 1.1]{arthurtransferwoi}; this argument already appeared in the proof of Lemma \ref{lem.tranv.result.2}, Step 2. Taking characters gives
        $$
            \begin{tikzcd}
                0 \ar[r] & X^\star(Z_{\widehat{M}}^\Gamma) \ar[r] & X^\star(Z_{\widehat{A}}^\Gamma) \oplus X^\star(Z_{\widehat{B}}^\Gamma) \ar[r] & X^\star(Z_{\widehat{A}}^\Gamma \cap Z_{\widehat{B}}^\Gamma) \ar[r] & 0
            \end{tikzcd}.
        $$
        We apply (\ref{eqn.index.mult.ses.fg.gps}) to this sequence. Note that tensoring with $\R$ recovers the short-exact sequence of cocharacter spaces
        $$
            \begin{tikzcd}
                0 \ar[r] & \a_M \ar[r] & \a_A \oplus \a_B \ar[r] & \a_L \ar[r] & 1
            \end{tikzcd}.
        $$
        Splitting off the common orthogonal summand $\a_L$, the absolute determinant is $\Delta = d_M^L(A,B)$. So Lemma \ref{lem.index.mult.ses.fg.gps} gives
    \begin{equation}
            \nu(Z_{\widehat{A}}^\Gamma) \nu(Z_{\widehat{B}}^\Gamma) = d_M^L(A,B) \nu(Z_{\widehat{M}}^\Gamma) \nu(Z_{\widehat{A}}^\Gamma \cap Z_{\widehat{B}}^\Gamma).\label{eqn.first.nu(int).comp}
        \end{equation}

        The second computation of $\nu(Z_{\widehat{A}}^\Gamma \cap Z_{\widehat{B}}^\Gamma)$ passes instead through the finite-index containment of $Z_{\widehat{L}}^\Gamma$. Dualizing gives
        $$
            \begin{tikzcd}
                0 \ar[r] & X^\star((Z_{\widehat{A}}^\Gamma \cap Z_{\widehat{B}}^\Gamma)/Z_{\widehat{L}}^\Gamma) \ar[r] & X^\star(Z_{\widehat{A}}^\Gamma \cap Z_{\widehat{B}}^\Gamma) \ar[r] & X^\star(Z_{\widehat{L}}^\Gamma) \ar[r] & 0
            \end{tikzcd}.
        $$
        Since the first term is finite (this follows from nonvanishing of the $d$) the induced map on real vector spaces is just $\a_L \cong \a_L$ and has determinant 1. Lemma \ref{lem.index.mult.ses.fg.gps} therefore becomes the elementary statement
        $$
            \nu(Z_{\widehat{A}}^\Gamma \cap Z_{\widehat{B}}^\Gamma) = \left|\frac{Z_{\widehat{A}}^\Gamma \cap Z_{\widehat{B}}^\Gamma}{Z_{\widehat{L}}^\Gamma}\right|^{-1} \nu(Z_{\widehat{L}}^\Gamma).
        $$
        Substituting this into (\ref{eqn.first.nu(int).comp}) gives the equality (\ref{eqn.lattice.comp.ABML}) in terms of the indices $\nu$. 
        
        Solving (\ref{eqn.lattice.comp.ABML}), applied to each index $i$, for the coefficient $e$ gives
        $$
            e_{M_i}^{L_i}(A_i, B_i) = \frac{\nu(Z_{\widehat{A}_i}^\Gamma) \nu(Z_{\widehat{B}_i}^\Gamma)}{\nu(Z_{\widehat{L}_i}^\Gamma) \nu(Z_{\widehat{M}_i}^\Gamma)}.
        $$
        On the other hand, Lemma \ref{lem.covol.coroots} gives, in terms of our index $\nu$, the description
        $$
        c_{M_1, M_2}^{H_1, H_2} = \frac{\nu(Z^\Gamma_{\widehat{M}_1}) \nu(Z^\Gamma_{\widehat{H}_2})}{\nu(Z^\Gamma_{\widehat{M}_2}) \nu(Z^\Gamma_{\widehat{H}_1})}.
        $$
        These formulae verify the equality of Lemma \ref{lem.e.compat.c}.    
    \end{proof}
\end{lemma}

\section{Peeling the regular places.}\label{subsect.nonstd.reg.peeling}

We now need to demonstrate that sufficiently regular places are peelable, as we did for the endoscopic WFL in \S \ref{subsect.peel.regular}-\ref{subsect.transverse.places}. We retain our global setup $(M_1, G_1, M_2, G_2, \psi)$ and our parameter $a_M \in \bAe{M_i}(k_0) \cap \bA^{\good}(k_0)$. 

In the endoscopic case we performed a basic computation of the $J_{M,v}^{L,s}$ at regular places in Lemma \ref{lem.unit.disc.result.vanish}, then packaged these over a set $V$ of regular places to deduce (Corollary \ref{cor.reg.v.peelable}) that regular places are peelable; the stable orbital integrals have a corresponding result. All notation at a fixed place $v$ is the standard local notation for local fields and rings of integers as in \S \ref{subsect.local.setup}.
\begin{lemma}\label{lem.nonstd.peel.regular}
    Let $M \subset G$ a containment of semistandard Levis.
    \begin{enumerate}
        \item (Unit-discriminant normalization) Suppose $V \subseteq |X_0|$ is a set of places such that $a_M$ avoids $\fD_M$ on $V$. Then with respect to the unramified volume-one normalized measure,
        $$
            S_{M,V}^{M}(a_M) = J_{M,V}^{M, 1}(a_M) = 1.
        $$
        \item (Unit-resultant vanishing) Let $M \subsetneq L$ an intermediate Levi, so $b$ is generically $L$-regular semisimple. Suppose $V \subseteq |X_0|$ is a set of places such that $a_M$ avoids $\fR_M^G$ on $V$. Then
        $$
            S_{M,V}^L(a_M) = 0.
        $$
    \end{enumerate}
    In particular, if $a_M$ avoids $\fD_G$ on $V$ then $V$ is peelable.
    \begin{proof}
        The corresponding statements for the $J^s$ are \textit{loc.\ cit.} The stable orbital integrals inherit the corresponding result by their defining recursion. The renormalization trick of Remark \ref{rmk.local.renorm.D} is used to handle places on the support of $D$. Then the nonstandard WFL
        $$
            S_{M_1, V}^{L_1}(a_M) = c_{M_1, M_2}^{L_1, L_2} S_{M_2, V}^{L_2}(a_M)
        $$
        reduces to $1 = 1$ when $M_i = L_i$ and $0 = 0$ otherwise.
    \end{proof}
\end{lemma}

\section{Transverse places.}\label{subsect.nonstd.transverse}

Since our nonstandard endoscopic datum depends only on $M_i \subseteq G_i$, there are only two interesting transverse cases to check. As for the endoscopic WFL, the crucial new computation is the transverse-resultant case, which will reduce to reusing the rank 1 computation of the stable orbital integrals from Lemma \ref{lem.tranv.result.2}.

\begin{lemma}\label{lem.nonstd.transv.disc}
    Let $v \in |X_0|$ be a place disjoint from $\on{supp}(D)$, and write $b$ for the restriction of the global parameter $a_M$ to the completed disk around $v$. Suppose
    \begin{equation}
        v(\fD_{M}(b)) = 1, \qquad v(\fR_{M}^G(b)) = 0.\label{eqn.nonstd.transv.disc}
    \end{equation}
    Then $v$ is peelable.
    \begin{proof}
        Fix a compatible pair of intermediate $M_i \subseteq L_i \subseteq G_i$. If $M_i = L_i$, then the identity is essentially the (unweighted) nonstandard WFL, \cite[Th\'eor\`eme 1.12.7]{ngofl}. The only discrepancy is the discriminant factor $|D^{M_i}(Y_{i, a_M})|_v^{1/2}$ included in our weighted orbital integrals. But indeed
        $$
            |D^{M_1}(Y_{1, a_M})|_v^{1/2} = |D^{M_2}(Y_{2, a_M})|_v^{1/2};
        $$
        this can be checked after splitting \'etale cover, and then the condition that for corresponding roots 
        $$
            \alpha_2 \circ \psi_\star = b(\alpha_2) \alpha_1
        $$
        with $b(\alpha_2)$ a $p$-adic unit implies that the corresponding $Y_i$ intersect each branch of the discriminant with the same multiplicity at $v$.
        
        For $M_i \neq L_i$, the unit-resultant vanishing from Lemma \ref{lem.nonstd.peel.regular} renders the equality vacuous. It remains to check that the stable orbital integral is nonzero; we checked this in the proof of Lemma \ref{lem.transv.disc}.
    \end{proof}
\end{lemma}

\begin{lemma}\label{lem.nonstd.transv.result}
    Let $v \in |X_0|$ be a place disjoint from $\on{supp}(D)$, and write $b$ for the restriction of the global parameter $a_M$ to the completed disk around $v$. Suppose
    \begin{equation}
        v(\fD_{M}(b)) = 0, \qquad v(\fR_{M}^G(b)) = 1.\label{eqn.nonstd.transv.result}
    \end{equation}
    Then $v$ is peelable.
    \begin{proof}
        The unit-discriminant result from Lemma \ref{lem.nonstd.peel.regular} implies the result for $M_i = L_i$, as well as the nonvanishing condition of peelability, so assume $M_i \subsetneq L_i$. If $v(\fR_{M_i}^{L_i}(b)) = 0$, then the nonstandard WFL is the identity $0 = 0$ from the unit-resultant result from the same Lemma, so assume $v(\fR_{M_i}^{L_i}(b)) = 1$. Since $\psi$ scales corresponding roots by $p$-adic units, it preserves the valuation of each root. So the ``bad'' geometric roots $\alpha_i$ for the two groups correspond, and correspondingly the geometric $R_{i,\alpha}^{\on{geom}}$ correspond under $\psi$. 

        We apply the computation performed in \ref{lem.tranv.result.2} to the two groups separately. If $\dim \a_{M_i}^{L_i} > 1$ then
        $$
            S_{M_i,\{v\}}^{L_i}(Y_{i,b}) = 0
        $$
        and the nonstandard WFL is certainly true, so assume $\dim \a_{M_i}^{L_i} = 1$. Then (\ref{eqn.stab.rk1.recurs.formula}) gives
        $$
            S_{M_i, \{v\}}^{L_i}(Y_{i,b}) = (1 - q_v^{-1}) \deg(v) \frac{\|h_{i, \alpha}\|}{N_i}, 
        $$
        where $h_{i,\alpha} \in \a_{M_i}^{L_i}$ is the image of $\alpha_i^\vee$ and 
        $$
            N_i := \left|\ker\left( \alpha_i^\vee |_{Z_{\widehat{M}_i}^\Gamma}\right) / Z_{\widehat{L}_i}^\Gamma \right|.
        $$
        Indeed, these formulae hold even with our choice of $\psi$-compatible normalizations for the measures: the calculation requires normalized volume for integral points on the centralizers, here unramified tori by vanishing of the discriminant valuation. Transporting Haar measure across the Lie-algebra identification from $\psi$ preserves these normalizations, cf. \cite[\S 1.12]{ngofl}. It remains only to interpret $N_i$ in a way which allows us to compare them. 

        Write $$\Lambda_i := \on{proj}_{\a_{M_i}^{L_i}} \left(\Z \langle \ul{\Phi}_{L_i}^\vee \rangle_\Gamma\right).$$ This is a lattice in the one-dimensional space $\a_{M_i}^{L_i}$, containing but not necessarily equal to $\Z h_{i, \alpha}$. The crucial observation is that $N_i$ measures the discrepancy:
        \begin{equation}
            N_i = [\Lambda_i : \Z h_{i, \alpha}].\label{eqn.Ni.const.nonstd.interp}
        \end{equation} 
        To see why, choose a $\Gamma$-stable positive system for which $M_i$ is standard. The simple coroots of $M_i$ are a $\Gamma$-stable subset of those on $L_i$, and passing to orbit classes gives the inclusion $\Z \langle \ul{\Phi}_{M_i}^\vee \rangle_\Gamma \hookrightarrow \Z \langle \ul{\Phi}_{L_i}^\vee \rangle_\Gamma$. Passing to $\a_{M_i}^{L_i}$ is exactly the quotient by this subset, so
        \begin{equation}
            \Lambda_i \cong \Z \langle \ul{\Phi}_{L_i}^\vee \rangle_\Gamma/\Z \langle \ul{\Phi}_{M_i}^\vee \rangle_\Gamma.\label{eqn.Lambdai.coker}
        \end{equation}
        This lattice can also be measured relative to $X_\star(\T_i)$:
        $$
            \begin{tikzcd}
                0 \ar[r] & \Z \langle \ul{\Phi}_{L_i}^\vee \rangle_\Gamma/ \Z \langle \ul{\Phi}_{M_i}^\vee \rangle_\Gamma \ar[r] & \left(X_\star(\T_i)/\Z \langle \ul{\Phi}_{M_i}^\vee \rangle\right)_\Gamma \ar[r] & \left(X_\star(\T_i)/\Z \langle \ul{\Phi}_{L_i}^\vee \rangle\right)_\Gamma \ar[r] & 0
            \end{tikzcd}.
        $$
        Left-exactness follows because, as in Lemma \ref{lem.covol.coroots}, $\Z \langle \ul{\Phi}_{L_i}^\vee \rangle_\Gamma \to X_\star(\T_i)_\Gamma$ is injective: its kernel is torsion but the domain is torsion-free. Quotienting by the common subgroup $\Z \langle \ul{\Phi}_{M_i}^\vee \rangle_\Gamma$ and applying right exactness of coinvariants gives left-exactness.

        Recall that for any of $M_i, L_i$,
        $$
            X^\star(Z_{\widehat{H}_i}^\Gamma) = \left(X_\star(\T_i)/\Z \langle \ul{\Phi}_{H_i}^\vee \rangle\right)_\Gamma.
        $$
        The above sequence therefore is dual to the short exact sequence for the inclusion $Z_{\widehat{L}_i}^\Gamma \hookrightarrow Z_{\widehat{M}_i}^\Gamma$. 
        
        The cocharacter $\alpha_i^\vee$ is trivial on $Z_{\widehat{L}_i}^\Gamma$, by construction, so its restriction is a character
        $$
            \alpha_i^\vee |_{Z_{\widehat{M}_i}^\Gamma} \in X^\star\left(Z_{\widehat{M}_i}^\Gamma/Z_{\widehat{L}_i}^\Gamma\right).
        $$
        As in the proof of Lemma \ref{lem.tranv.result.2}, it is surjective. Restricting to its kernel gives a short exact sequence
        $$
        \begin{tikzcd}
            0 \ar[r] & \ker\left( \alpha_i^\vee |_{Z_{\widehat{M}_i}^\Gamma}\right)/Z_{\widehat{L}_i}^\Gamma \ar[r] & Z_{\widehat{M}_i}^\Gamma/Z_{\widehat{L}_i}^\Gamma \ar[r] & Z_{\widehat{M}_i}^\Gamma/\ker\left( \alpha_i^\vee |_{Z_{\widehat{M}_i}^\Gamma}\right) \ar[r] & 0
        \end{tikzcd}.
        $$
        Taking $X^\star$, we have shown that the middle term gives $\Lambda_i$. On the other hand, the quotient is $\G_m$, and the dual inclusion $\Z \hookrightarrow \Lambda_i$ sends the generator $1 \mapsto h_{i, \alpha}$, the image of $\alpha_i^\vee$ inside $\a_{M_i}^{L_i}$.
        
        Thus
        $$
            X^\star\left(\ker\left( \alpha_i^\vee |_{Z_{\widehat{M}_i}^\Gamma}\right) / Z_{\widehat{L}_i}^\Gamma \right) \cong \Lambda_i/\Z h_{i, \alpha},
        $$
        which proves by taking cardinalities (\ref{eqn.Ni.const.nonstd.interp}). In the line $\a_{M_i}^{L_i}$, 
        $$
            \|h_{i, \alpha}\| = [\Lambda_i : \Z h_{i,\alpha}] \on{covol}_{\a_{M_i}^{L_i}}(\Lambda_i),
        $$
        so we can rewrite our stable orbital integrals in a way amenable to comparison:
        \begin{equation}
            S_{M_i, \{v\}}^{L_i}(Y_{i,b}) = (1 - q_v^{-1}) \deg(v) \on{covol}_{\a_{M_i}^{L_i}}(\Lambda_i).\label{eqn.nonstd.stable.covol}
        \end{equation}

        As in the proof of the global nonstd WFL, we can rewrite Waldspurger's constant equally in terms of $\Gamma$-invariants or coinvariants. Choosing the latter and passing through the isomorphism (\ref{eqn.Lambdai.coker}), we can rewrite the constant simply as
        $$
            c_{M_1, M_2}^{L_1, L_2} = \frac{[(\Z \langle \Phi_{L_2}^\vee \rangle)_{\Gamma} : \psi_\star((\Z \langle \Phi_{L_1}^\vee \rangle)_{\Gamma})]}{[(\Z \langle \Phi_{M_2}^\vee \rangle)_{\Gamma} : \psi_\star((\Z \langle \Phi_{M_1}^\vee \rangle)_{\Gamma})]} = [\Lambda_2 : \psi_\star(\Lambda_1)].
        $$
        Using the measures taken compatibly along $\psi$, we can rewrite this in terms of lattice covolumes and then pass one across $\psi$:
        $$
            c_{M_1, M_2}^{L_1, L_2} = \frac{\on{covol}_{\a_{M_2}^{L_2}}(\psi_\star\Lambda_1)}{\on{covol}_{\a_{M_2}^{L_2}}(\Lambda_2)} = \frac{\on{covol}_{\a_{M_1}^{L_1}}(\Lambda_1)}{\on{covol}_{\a_{M_2}^{L_2}}(\Lambda_2)}.
        $$
        Combining with (\ref{eqn.nonstd.stable.covol}) gives immediately the nonstandard WFL.

    \end{proof}
\end{lemma}

\section{Globalization and the proof of the nonstandard WFL.}\label{subsect.nonstd.globalization}

We begin now in the local setting, with notation as in \S \ref{subsect.local.setup}. In particular, $F_{0,v}$ is a local field with residue field $k_v$ and $\O_{0,v}$ is its ring of integers. Let $(M_{1}, G_{1}, M_2, G_2, \psi)$ be a nonstandard endoscopic datum over $\Spec \O_{0,v}$ in good characteristic. In particular, $G_1, G_2/\O_{0,v}$ are both defined by twisting along a common $\Out_{12}$-torsor. Writing $\Gamma$ for the finite image of the monodromy, by Lemma \ref{lem.globalize.group} we can engineer a global nonstandard datum, which we also denote $(M_{1}, G_{1}, M_2, G_2, \psi)$, over some smooth, projective, geometrically connected curve $X_0/k_v$ which restricts to the original datum at the place $v$, and is split at some other place $\infty$. 

The construction moreover engineers $D_v = \Gamma$, so as in Remark \ref{rmk.nonstd.local.global.discrep}, the statement of the nonstandard WFL at $v$ does actually reproduce Waldspurger's local nonstandard WFL of Theorem \ref{thm.nonstd.WFL}.  
Fix our local parameter $b \in \fc_{M}(\O_{0,v}) \cap \fc_M^{\Greg}(F_{0,v})$ and write $Y_{i,b} \in \fm_i(\O_{0,v})$ for its Kostant lifts. By Proposition \ref{prop.local.constancy}, we can find an $n$ larger than the discriminant valuation and large enough that the values of the stable orbital integrals $S_{M_i,v}^{G_i}(Y_{i,b})$ are determined by the $n$-th order jets of $b$. By Lemma \ref{lem.bert.glob} and the discussion prior, applied when $M'_i = M_i$, we can find $a_M \in \cA_{M_1}(k_v) = \cA_{M_2}(k_v)$ such that
\begin{itemize}
    \item $a_M|_{nv} \equiv b|_{nv}$,
    \item $a_M(\infty) = \chi_M(t_\infty)$ for a fixed $t_\infty \in \ft_{M,D,\infty}^{\Greg}(k_v)$ (everything appearing here is canonically identified between $M_1 \subseteq G_1$ and $M_2 \subseteq G_2$, by Proposition \ref{prop.nonstd.same.chevalley.lie.torus}),
    \item $a_M(x)$ is $G_i$-regular at each $x \in \on{supp}(D)$,
    \item $a_M$ is $M_i$-elliptic, and
    \item $a_M$ intersects $\fD_{M} + \fR_M^G$ transversally away from $\{v, \infty\} \cup \on{supp}(D)$.
\end{itemize}
Write $a \in \cA_G(k_v)$ for the image of $a_M$. We require $D$ to have been sufficiently positive that
$$
    \deg(D) > 2g - 2 + \frac{|W_{G}|}{2} v(\fD_G(a)),
$$
(since $n$ is greater than the discriminant valuation of $b$, the right side is determined by the local parameter $b$ at $v$). Then Theorem \ref{thm.suff.criterion.for.goodness} implies that $a \in \cA_{G_i}^{\good}$ for both $i = 1,2$.\footnote{It is good that we have such a numerical condition; note that we have no reason to know the good loci for $G_1, G_2$ coincide!}

\begin{proof}[Proof of Theorem \ref{thm.nonstd.WFL}]
We have arranged our global parameter $a_M$ to be $M_i$-elliptic. Thus, choosing any $t$ gives a parameter $(a_M, t) \in \bAe{M}$. The parameter satisfies all hypotheses of Corollary \ref{cor.global.nonstd.WFL}, so the nonstandard WFL for $a_M$ and $(M_1, G_1, M_2, G_2, \psi)$ holds on $V = |X_0|$. By Lemmas \ref{lem.nonstd.peel.regular}, \ref{lem.nonstd.transv.disc}, \ref{lem.nonstd.transv.result}, all places except $v$ are peelable, so the nonstandard Peeling Lemma \ref{lem.peeling.nonstd} implies that the nonstandard WFL holds on $\{v\}$. By Remark \ref{rmk.nonstd.local.global.discrep}, this recovers exactly Waldspurger's statement in terms of local endoscopic data. This completes the proof.
\end{proof}

\printbibliography
\end{document}